\documentclass[12pt]{amsart}

\usepackage{amscd,amssymb,amsthm,verbatim}

\newcommand{\al}{\alpha}
\newcommand{\area}{\operatorname{area}}

\newcommand{\be}{\beta}

\newcommand{\C}{{\mathbb C}}

\newcommand{\const}{\operatorname{const.}}
\newcommand{\core}{\operatorname{core}}

\newcommand{\D}{\partial}
\newcommand{\De}{\Delta}

\newcommand{\diam}{\operatorname{diam}}
\newcommand{\Diff}{\operatorname{Diff}}

\newcommand{\dist}{\operatorname{dist}}
\newcommand{\Div}{\operatorname{div}}

\newcommand{\dvol}{\operatorname{dvol}}
\newcommand{\End}{\operatorname{End}}

\newcommand{\f}{{\mathcal F}}

\newcommand{\Ga}{\Gamma}

\newcommand{\gk}{\operatorname{GK}}

\newcommand{\grad}{\operatorname{grad}}
\newcommand{\h}{{\mathcal H}}

\newcommand{\Hess}{\operatorname{Hess}}

\newcommand{\HH}{\operatorname{H}}

\newcommand{\Id}{\operatorname{Id}}

\newcommand{\im}{\operatorname{im}}
\newcommand{\Image}{\operatorname{Im}}

\newcommand{\inj}{\operatorname{inj}}
\newcommand{\injrad}{\operatorname{InjRad}}
\newcommand{\interior}{\operatorname{Interior}}
\newcommand{\Int}{\operatorname{int}}

\newcommand{\Isom}{\operatorname{Isom}}
\newcommand{\Ker}{\operatorname{Ker}}

\newcommand{\length}{\operatorname{length}}
\renewcommand{\L}{{\mathcal L}}
\newcommand{\La}{\Lambda}

\newcommand{\lra}{\longrightarrow}
\newcommand{\M}{{\mathcal M}}

\newcommand{\Nil}{\operatorname{Nil}}

\newcommand{\Om}{\Omega}

\newcommand{\pt}{\operatorname{pt}}
\newcommand{\Q}{{\mathbb Q}}
\newcommand{\R}{{\mathbb R}}
\newcommand{\I}{{\mathbb I}}

\newcommand{\reg}{\operatorname{reg}}

\newcommand{\restr}{\mbox{\Large \(|\)\normalsize}}

\newcommand{\Ric}{\operatorname{Ric}}
\newcommand{\Rm}{\operatorname{Rm}}

\newcommand{\si}{\sigma}

\newcommand{\sing}{\operatorname{sing}}
\newcommand{\s}{{\mathcal S}}
\newcommand{\SL}{\operatorname{SL}}
\newcommand{\SO}{\operatorname{SO}}
\newcommand{\Sol}{\operatorname{Sol}}

\newcommand{\spec}{\operatorname{spec}}

\newcommand{\supp}{\operatorname{supp}}

\renewcommand{\th}{\theta}

\newcommand{\Tr}{\operatorname{Tr}}

\newcommand{\vol}{\operatorname{vol}}
\newcommand{\Z}{{\mathbb Z}}
\newcommand{\cangle}{\widetilde{\angle}}
\newcommand{\ctits}{C_T}
\newcommand{\de}{\delta}

\newcommand{\eps}{\epsilon}
\newcommand{\ga}{\gamma}
\newcommand{\la}{\lambda}
\newcommand{\ol}{\overline}
\renewcommand{\r}{{\mathcal R}}
\newcommand{\ra}{\rightarrow}
\newcommand{\Si}{\Sigma}
\newcommand{\tangle}{\partial_T}
\newcommand{\tits}{\partial_T}

\numberwithin{equation}{section}
\theoremstyle{plain}

\newtheorem{lemma}[equation]{Lemma}
\newtheorem{theorem}[equation]{Theorem}
\newtheorem{proposition}[equation]{Proposition}
\newtheorem{corollary}[equation]{Corollary}

\newtheorem{sublemma}[equation]{Sublemma}
\newtheorem{claim}[equation]{Claim}

\theoremstyle{definition}
\newtheorem{definition}[equation]{Definition}

\theoremstyle{remark}
\newtheorem{remark}[equation]{Remark}
\begin{document}
\title{Notes on Perelman's papers}
\author{Bruce Kleiner}
\address{Department of Mathematics\\
Yale University\\
New Haven, CT  06520-8283\\
USA}
\email{bruce.kleiner@yale.edu}
\author{John Lott}
\address{Department of Mathematics\\
University of Michigan\\
Ann Arbor, MI  48109-1109\\
USA}
\email{lott@umich.edu}
\thanks{Research supported by NSF grants DMS-0306242 and DMS-0204506,
and the Clay Mathematics Institute}
\date{February 13, 2013}
\maketitle

\section{Introduction}

These are notes on Perelman's papers 
``The Entropy Formula for the Ricci Flow and its Geometric
Applications'' \cite{Perelman} and ``Ricci Flow with Surgery 
on Three-Manifolds' \cite{Perelman2}.  
In these two remarkable preprints, which were posted on the ArXiv 
in  2002 and 2003, Grisha Perelman
announced a proof of the Poincar\'e Conjecture, and more generally Thurston's
Geometrization Conjecture, using the Ricci flow approach of Hamilton.  
Perelman's proofs are concise and,
at times, sketchy. The purpose of these notes is to provide the details
that are missing in
\cite{Perelman} and \cite{Perelman2}, which contain Perelman's arguments for 
the Geometrization Conjecture.  

Among other things, 
we cover the construction of the Ricci flow with surgery 
of \cite{Perelman2}.
We also discuss the long-time behavior of the 
Ricci flow with surgery, which is needed for the full Geometrization
Conjecture.
The papers \cite{Colding-Minicozzi,Perelman3},
which are not covered in these notes, 
each provide a shortcut in the case of the Poincar\'e Conjecture.
Namely, these papers show
that if the initial manifold is simply-connected then the 
Ricci flow with surgery becomes extinct in a
finite time, thereby removing the issue of the long-time
behavior. Combining this claim with the proof of existence of
Ricci flow with surgery gives 
the shortened proof in the simply-connected case.

These notes are intended for readers with a solid background in geometric analysis.
Good sources for background material on Ricci flow are 
\cite{Chow-Knopf,Chow-Lu-Ni,Hamilton,Topping}.
The notes are self-contained but are designed to be read along with
\cite{Perelman,Perelman2}.
For the most part we follow the format of 
\cite{Perelman,Perelman2} and use the section numbers
of \cite{Perelman,Perelman2} to label our
sections. We have done this in order to respect the
structure of \cite{Perelman,Perelman2} and to facilitate the use of
the present notes as a companion to \cite{Perelman,Perelman2}.
 In some places we have rearranged Perelman's arguments
or provided alternative arguments,
but we have refrained from an overall reorganization.  

Besides providing details for Perelman's proofs, we have included some 
expository material in the form of overviews and appendices.
Section \ref{overview0} contains an overview of the Ricci flow approach to
geometrization of $3$-manifolds. Sections \ref{overview1} and \ref{overview2} contain
overviews of \cite{Perelman} and \cite{Perelman2}, respectively.
The appendices discuss some background material and
techniques that are used throughout the notes.

Regarding the proofs, the papers \cite{Perelman,Perelman2}  contain some incorrect statements 
and incomplete arguments, which we have attempted to point out to the reader.
(Some of the mistakes in \cite{Perelman} were corrected in \cite{Perelman2}.)
We did not find any
serious problems, meaning problems that cannot be corrected using the methods
introduced by Perelman.

We will refer to Section X.Y of \cite{Perelman} as I.X.Y, and
Section X.Y of \cite{Perelman2} as II.X.Y.
A reader may wish to start with the overviews, which 
explain the logical structure of the arguments and the
interrelations between the sections.
It may also be helpful to browse through the appendices before delving into the 
main body of the material.

These notes have gone through various versions, which were
posted at \cite{KleinerLottwebsite}.
An initial version with notes on \cite{Perelman}
was posted in June 2003.
A version covering \cite{Perelman,Perelman2} was posted in September
2004. 
After the May 2006 version of these notes
 was  posted on the ArXiv, expositions of Perelman's
work appeared in
\cite{Cao-Zhu} and \cite{Morgan-Tian}.
 \\ \\
\noindent
{\em Acknowledgements.}
In the preparation of the September 2004 version of our notes,
we benefited from a 
workshop on Perelman's surgery procedure
that was held in August 2004, at Princeton.
We thank the participants of that meeting, as well as 
the Clay Mathematics Institute for supporting it.

We are grateful for comments and corrections that we have
received regarding earlier versions of these
notes.  
We especially thank G\'erard Besson for comments on Theorem
\ref{noncollthm} and Lemma \ref{bdrysmoothness}, Peng Lu for
comments on Lemma \ref{II.4.5}, Bernhard Leeb for
discussions on the issue of uniqueness for the standard
solution and Richard Bamler for comments on 
Proposition \ref{propII.6.3} and Proposition \ref{propII.6.4}.
Along with these people we thank
Mike Anderson, Albert Chau,
Ben Chow, Xianzhe Dai, Jonathan Dinkelbach,
Sylvain Maillot, John Morgan, 
Joan Porti, Gang Tian, Peter Topping, Bing Wang, 
Guofang Wei, Hartmut Weiss, 
Burkhard Wilking, Jon Wolfson, Rugang Ye and Maciej Zworski.

We thank the referees for their detailed and helpful comments on
this paper.

We thank the Clay Mathematics Institute for supporting the writing
of the notes on \cite{Perelman2}.

\tableofcontents

\section{A reading guide}

Perelman's papers contain a wealth of results about Ricci flow.
We cover all of these results, whether or not they
are directly relevant to the Poincar\'e and Geometrization Conjectures.

Some readers may wish to take an abbreviated route that
focuses on the proof of the 
Poincar\'e Conjecture or the Geometrization Conjecture.
Such readers can try the following itinerary.  

Begin with the overviews in Sections \ref{overview0} and \ref{overview1}.
Then review Hamilton's compactness theorem and its variants, as described in
Appendix \ref{subsequence}; an exposition is in 
\cite[Chapter 7]{Topping}. 
Next, read I.7 (Sections \ref{I.7}-\ref{I.7.3}),
followed by I.8.3(b) (Section \ref{secI.8.3}).  
After reviewing the theory of Riemannian
manifolds and Alexandrov spaces of nonnegative sectional curvature 
(Appendix \ref{Alexandrov} and references therein),
proceed to I.11 (Sections \ref{secI.11.1}-\ref{secI.11.9}), followed by
II.1.2 and I.12.1 (Sections \ref{II.1.2}-\ref{I.12.1}).

At this point, the reader should be ready for 
the overview of Perelman's second paper in Section \ref{overview2},  and 
can proceed with II.1-II.5 (Sections
\ref{notationterminology}-\ref{surgeryflow}).  In conjunction with
one of the finite extinction time results 
\cite{Colding-Minicozzi,Colding-Minicozzi2,Perelman3},
this completes the proof of the Poincar\'e Conjecture.

To proceed with the rest of the proof of the Geometrization Conjecture, 
the reader can begin with the large-time estimates for nonsingular
Ricci flows, which appear in I.12.2-I.12.4 
(Sections \ref{I.12.2}-\ref{I.12.4}).
The reader can then go to II.6 and II.7 (Sections \ref{II.6}-\ref{II.7.4}).

The main topics that are missed by such an abbreviated route are
the ${\mathcal F}$ and ${\mathcal W}$ functionals
(Sections \ref{I.1.1}-\ref{I.5}), Perelman's differential Harnack inequality
(Section \ref{I.9}),  pseudolocality (Sections
\ref{statement of I.10.1}-\ref{I.10.5}) and Perelman's alternative proof
of cusp incompressibility (Section \ref{II.8}).

\section{An overview of the Ricci flow approach to 
$3$-manifold geometrization}
\label{overview0}

This section is an overview of the Ricci flow approach to
$3$-manifold geometrization.  We make no attempt to present the
history of the ideas that go into the argument.
We caution the reader that for the sake of readability,
in many places we have suppressed technical points and deliberately
oversimplified the story.
 The overview will introduce the argument in 
three passes, with successively greater precision and detail : we start with a 
very crude sketch, then expand this to a step-by-step outline of the 
strategy, and then move on to more detailed commentary on specific
points.  

Other overviews may be found in \cite{Cao-Chow,Morgan}.  The 
primary objective of our exposition is 
to prepare the reader for a more detailed study of Perelman's work.

We refer the reader to Appendix \ref{geom} for the statement 
of the geometrization conjecture.

By convention, all manifolds and Riemannian metrics in this section 
will be smooth.
We follow the notation of
\cite{Perelman,Perelman2} for pointwise quantities:  $R$ 
denotes the scalar curvature, 
$\Ric$ the Ricci curvature,  and $|\Rm|$ the largest absolute value of 
the sectional curvatures.   
An inequality such as $\Rm\geq C$ means that all of the 
sectional curvatures at a
point or in a region, depending on the context, are bounded below by $C$.
In this section, we will specialize to three dimensions.

\subsection{The definition of Ricci flow, and some basic properties}

Let $M$ be a compact  $3$-manifold and let $\{g(t)\}_{t\in [a,b]}$
be a smoothly varying family of  Riemannian metrics on $M$.
Then $g(\cdot)$ satisfies the {\em Ricci flow equation} if 
\begin{equation}
\frac{\D g}{\D t}(t) \: = \:
- \: 2 \: \Ric(g(t))
\end{equation}
holds for every $t\in[a,b]$.
Hamilton showed in \cite{Hamiltonnnnn} that
for any Riemannian metric $g_0$ on $M$, there is a $T\in (0,\infty]$
with the property that there is a (unique) solution $g(\cdot)$
to the Ricci flow equation defined on the time interval $[0,T)$ with
$g(0) = g_0$, so
that if $T<\infty$ then the curvature of $g(t)$ becomes unbounded
as $t\ra T$.  We refer to this maximal solution as the Ricci flow with
initial condition $g_0$. If $T<\infty$ then we call $T$ the blow-up
time.  A basic example 
is the shrinking  round $3$-sphere, with $g_0 \: = \: r_0^2 \: g_{S^3}$ and
$g(t) \: = \: (r_0^2 - 4t) \:  g_{S^3}$, in which case $T \: = \: \frac{r_0^2}{4}$.

Suppose that $M$ is simply-connected.
Based on the round $3$-sphere example, one could hope that every Ricci flow  
on $M$ blows up in finite time and becomes round while shrinking to a  point,
as $t$ approaches the blow-up time $T$. If so, then
by rescaling and taking a limit as $t \rightarrow T$, one would show
that $M$ admits a metric of constant positive sectional curvature and 
therefore,
by a classical theorem, is diffeomorphic to $S^3$. The analogous argument does work
in two dimensions \cite[Chapter 5]{Chow-Knopf}. Furthermore, if the initial 
metric $g_0$
 has positive Ricci curvature then Hamilton
showed in \cite{Hamilton} that this is the correct picture:
the manifold shrinks to a point in finite time and becomes round
as it shrinks.

One is then led to ask what can happen if $M$ is simply-connected
but $g_0$ does not have positive Ricci curvature.  Here a new 
phenomenon can occur --- the Ricci flow solution may become singular
before it has time to shrink to a point, due to a possible neckpinch. 
A neckpinch is modeled by a product region $(- c, c) \times S^2$ in which
one or many $S^2$-fibers separately shrink to a point at time $T$, due to the
positive curvature of $S^2$.   The formation of  neckpinch (and other)
singularities prevents one from continuing the Ricci flow.  
In order to continue the evolution some intervention is required, and this
is the role of surgery.  Roughly speaking, the
idea of surgery is to remove a neighborhood diffeomorphic to 
$(- c^\prime, c^\prime) \times S^2$ containing the shrinking $2$-spheres,
and cap off the resulting boundary components by 
gluing in $3$-balls.   Of course the
topology of the manifold changes during surgery -- for instance it may become
disconnected -- but it changes in a controlled way. 
The postsurgery Riemannian manifold is smooth,
so one may  restart the Ricci flow using 
it as an initial condition.  
When continuing the flow one may encounter further neckpinches,
which give rise to further surgeries, etc.
One hopes that
eventually all of the connected components
shrink to points while becoming round, i.e. that the Ricci flow solution 
has a {\em finite extinction time}.

\subsection{A rough outline of the Ricci flow proof of the Poincar\'e
Conjecture}
We now give a step-by-step glimpse of the proof, 
stating the needed steps as claims.

One starts with a compact orientable $3$-manifold $M$ with an 
arbitrary metric $g_0$. For the moment we do not assume that
$M$ is simply-connected. Let $g(\cdot)$ be the Ricci flow
with initial condition $g_0$, defined on $[0,T)$.
Suppose that $T <\infty$.   Let $\Omega \subset M$ be 
the set of points $x \in M$ for which  
$\lim_{t \rightarrow T^-} R(x,t) $ exists 
and is finite.
Then $M - \Omega$ is the part of $M$ that is going singular.
(For example, in the case of a single standard neckpinch, $M - \Omega$
is a $2$-sphere.) The first claim says what $M$ looks like
near this singularity set.

\begin{claim} \cite{Perelman2} \label{introclaim1}
The set $\Omega$ is  open and as $t\ra T$, the evolving metric
$g(\cdot)$ converges smoothly 
on compact subsets of $\Omega$ to
a Riemannian metric $\overline{g}$.  
There is a geometrically defined neighborhood
$U$ of $M - \Omega$ such that each connected component of $U$
is either \\ \\
\noindent
A. Compact and diffeomorphic to $S^1 \times S^2$,
$S^1 \times_{\Z_2} S^2$ or $S^3/\Gamma$, where $\Gamma$ is a finite
subgroup of $\SO(4)$ that acts freely and isometrically on the round $S^3$.
(In writing $S^1 \times_{\Z_2} S^2$, the generator of $\Z_2$ acts on
$S^1$ by complex conjugation and on $S^2$ by the antipodal map.
Then $S^1 \times_{\Z_2} S^2$ is diffeomorphic to $\R P^3 \# \R P^3$.)

or \\ \\
\noindent
B. Noncompact and diffeomorphic to $\R\times S^2$,
$\R^3$
or the twisted line
bundle $\R \times_{\Z_2} S^2$ over $\R P^2$.  \\

In Case B, the connected component meets $\Omega$ in 
geometrically controlled collar regions diffeomorphic to $\R\times S^2$. 
\end{claim}

Thus Claim \ref{introclaim1} provides a topological description of 
a neighborhood $U$ of the region
$M-\Omega$ where the Ricci flow is going singular, along with some
geometric control on $U$.

\begin{claim} \cite{Perelman2} \label{introclaim2}
There is a well-defined way 
to perform surgery on $M$, which yields a smooth  post-surgery  
manifold $M'$ with a Riemannian metric $g'$.
\end{claim}

Claim \ref{introclaim2} means that there is a well-defined 
procedure for  specifying the part of $M$ that will be removed, and
for gluing caps on the resulting manifold with boundary.  
The discarded part corresponds to the neighborhood $U$ in Claim \ref{introclaim1}.
The procedure is required to satisfy a number
of additional conditions which we do not mention here.

Undoing the surgery, i.e. going from the postsurgery manifold to the
presurgery manifold, amounts to restoring some discarded 
components (as in Case A of Claim \ref{introclaim1})
and performing connected sums of some
of the components of the postsurgery manifold, along with some
possible connected sums with a finite number of new $S^1 \times S^2$ and
$\R P^3$ factors.  The $S^1 \times S^2$ comes from the case when a
surgery does not disconnect the connected component where it is
performed.  The $\R P^3$ factors arise from the twisted line bundle
components in Case B of Claim \ref{introclaim1}.

After performing a surgery one lets the new manifold evolve under
the Ricci flow until one encounters the next blowup time (if 
there is one).  One then performs further surgery, lets the new manifold
evolve, and repeats the process.

\begin{claim} \cite{Perelman2} \label{introclaim3}
One can arrange the surgery procedure so that the surgery times
do not accumulate.
\end{claim}

If the surgery times were to accumulate, then one would have trouble continuing
the flow further, effectively killing the whole program.
Claim \ref{introclaim3} implies that by alternating Ricci flow and surgery,
one  obtains an evolutionary process that is defined for all time  
(though the manifold may become the empty set from some time
onward).  We call this   {\em Ricci flow with surgery}.

\begin{claim} \cite{Colding-Minicozzi,Colding-Minicozzi2,Perelman3} 
\label{introclaim4}
If the original manifold $M$ is simply-connected then any Ricci flow with
surgery on $M$ becomes extinct in finite time.
\end{claim}

Having a finite extinction time means that from some time onwards, the
manifold is the empty set. 
More generally, the same proof shows that
if the prime decomposition of the original manifold $M$ 
has no aspherical factors, then every Ricci flow with surgery on $M$
becomes extinct in finite time.
(Recall that a connected manifold $X$ is aspherical if
$\pi_k(X) = 0$ for all $k > 1$ or, equivalently, if its
universal cover is contractible.)

The Poincar\'e Conjecture follows immediately from the above claims.
From Claim \ref{introclaim4}, 
after some finite time the manifold is the empty set. 
From Claims \ref{introclaim1}, \ref{introclaim2} and \ref{introclaim3},
the original manifold $M$ is diffeomorphic to a connected sum of
factors that are each $S^1 \times S^2$ or a standard quotient
$S^3/\Gamma$ of $S^3$. As we are assuming that $M$ is simply-connected,
van Kampen's theorem implies that $M$ is diffeomorphic to a connected
sum of $S^3$'s, and hence is diffeomorphic to $S^3$.

\subsection{Outline of the proof of the Geometrization Conjecture}
We now drop the assumption that $M$ is
simply-connected. The main difference is that
Claim \ref{introclaim4} no longer applies, so the Ricci flow with surgery may
go on forever in a nontrivial way.  (We remark that Claim \ref{introclaim4}
is needed only for a shortened proof of the Poincar\'e Conjecture; the proof
in the general case is logically independent of Claim \ref{introclaim4} and
also implies the Poincar\'e Conjecture.)
The possibility that there are infinitely many surgery times is
not excluded, although it is not
known whether this can actually happen.  

A simple example of a Ricci flow that does not become extinct is when
$M = H^3/\Gamma$, where $\Gamma$ is a freely-acting cocompact discrete subgroup
of the orientation-preserving isometries of hyperbolic space $H^3$. If
$g_{hyp}$ denotes the metric on $M$ of constant sectional curvature
$-1$ and $g_0 \: = \: r_0^2 \: g_{hyp}$ then
$g(t) \: = \: (r_0^2 \: + \: 4t) \: g_{hyp}$. Putting 
$\widehat{g}(t) \: = \: \frac{1}{t} \: g(t)$, one finds that
$\lim_{t \rightarrow \infty} \widehat{g}(t)$ is the metric on $M$
of constant sectional curvature $- \: \frac14$, independent of
$r_0$.

Returning to the general case,
let $M_t$ denote the time-$t$ manifold in a Ricci flow with surgery.
(If $t$ is a surgery time
then we consider the postsurgery manifold.)
If for some $t$
a component of
$M_t$ admits a metric with
nonnegative scalar curvature then one can show  that
the component becomes extinct or admits a flat metric; either possiblity
is good enough when we are trying to prove the Geometrization
Conjecture for the initial manifold $M$.
So we will assume
that for every  $t$, each component of $M_t$ has a point with strictly
negative scalar curvature.

Motivated by the hyperbolic
example, we consider the metric $\widehat{g}(t) \: = \: \frac{1}{t} \: g(t)$
on $M_t$. Given $x \in M_t$, define the {\em intrinsic scale} $\rho(x,t)$ to be 
the radius $\rho$ such that $\inf_{B(x,\rho)} \Rm \: = \: - \: \rho^{-2}$,
where $\Rm$ denotes the sectional curvature of $\widehat{g}(t)$; this
is well-defined because the scalar curvature is negative somewhere in 
the connected component of $M_t$ containing $x$.
Given $w > 0$, define the
{\em $w$-thick part} of $M_t$ by 
\begin{equation}
M^+(w,t) \: = \:
\{ x \in M_t \: : \: \vol(B(x, \rho(x,t))) \: > \: w \: \rho(x,t)^3 \}.
\end{equation}
It is not excluded that $M^+(w,t) \: = \: M_t$ or $M^+(w,t) \: = \:
\emptyset$.
The next claim says that for any $w > 0$, as time goes on,
$M^+(w,t)$ approaches the $w$-thick part of a manifold of
constant sectional curvature $- \frac14$.

\begin{claim} \cite{Perelman2} \label{introclaim5}
There is a finite collection
$\{ (H_i, x_i) \}_{i=1}^k$ of complete pointed finite-volume $3$-manifolds
with constant sectional curvature $- \frac14$ and, for large $t$,
a decreasing function $\alpha(t)$ tending to zero and a family
of maps 
\begin{equation}
f_t \: : \bigsqcup_{i=1}^k H_i\:\supset\: \bigsqcup_{i=1}^k 
\:\:B \left( x_i, \frac{1}{\alpha(t)} \right)
 \rightarrow M_t\,, 
\end{equation}
such that\\
1. $f_t$ is $\alpha(t)$-close to being an isometry. \\
2. The image of $f_t$ contains $M^+(\alpha(t), t)$.\\
3. The image under $f_t$ of a cuspidal torus of
$\{H_i\}_{i=1}^k$ is incompressible in $M_t$.
\end{claim} 

The proof of Claim \ref{introclaim5} uses earlier work by
Hamilton \cite{Hamiltonn}.

\begin{claim} \cite{Perelman2,Shioya-Yamaguchi} \label{introclaim6}
Let $Y_t$ be the truncation of $\bigcup_{i=1}^k H_i$ obtained by
removing horoballs at distance approximately $\frac{1}{2 \alpha(t)}$
from the basepoints $x_i$. Then for large $t$, $M_t - f_t(Y_t)$ is
a graph manifold.
\end{claim}

Claim \ref{introclaim6} reduces to a statement in Riemannian geometry
about $3$-manifolds that are locally volume-collapsed with a lower
bound on sectional curvature.

Claims \ref{introclaim5} and \ref{introclaim6}, along with Claims
\ref{introclaim1}-\ref{introclaim3}, imply the geometrization conjecture,
cf. Appendix \ref{geom}.

In the remainder of this section, we will discuss some of the claims in 
more detail.

\subsection{Claim \ref{introclaim1} and the structure of singularities}
Claim \ref{introclaim1} is derived from a more localized statement, which
says that near points of large scalar curvature, the Ricci flow looks
very special : it is well-approximated by a special kind of
model Ricci flow, called a $\kappa$-solution.

\begin{claim}
\label{introclaim7}
Suppose that we have a given Ricci flow solution on a finite time interval.
If $x\in M$ and the scalar curvature $R(x,t)$ is large then 
in the time-$t$ slice,
there is a ball  centered at $x$
of radius comparable to $R(x,t)^{-\frac12}$ in which the geometry
of the Ricci flow 
is close to that of a ball in a $\kappa$-solution.
\end{claim}
The quantity $R(x,t)^{-\frac12}$ is sometimes called the {\em curvature scale}
at $(x,t)$, because it scales like a distance.  We will  define $\kappa$-solutions
below, but mention here that they are Ricci flows with nonnegative sectional curvature,
and they are {\em ancient}, i.e. defined on a time interval of the form
$(-\infty,a)$. 

The strength of Claim \ref{introclaim7} comes from the fact that there
is a good description of  $\kappa$-solutions.
 
\begin{claim} \cite{Perelman2} \label{introclaim8}
Any three-dimensional oriented $\kappa$-solution $(M_\infty, g_\infty(\cdot))$ falls 
into one of the
following types : \\
(a) A finite isometric quotient of the round shrinking $3$-sphere. \\
(b) A Ricci flow on a  manifold  diffeomorphic to $S^3$ or $\R P^3$. \\
(c) A standard shrinking round neck on $\R \times S^2$ \\
(d) A Ricci flow on a manifold diffeomorphic to $\R^3$, each time slice of
which is  asymptotically necklike at infinity.\\
(e)  The  $\Z_2$-quotient $\R \times_{\Z_2} S^2$ of a  shrinking round neck. 
\end{claim}
Together, Claims \ref{introclaim7} and \ref{introclaim8} say that where the 
scalar curvature is large, there is a region of diameter comparable to the
curvature scale  where one  sees either a closed manifold of known topology
(cases (a) and (b)),  a neck region
(case (c)),  a neck region capped off by a $3$-ball (case (d)), or a neck
region capped off by a twisted line bundle over $\R P^2$ (case (e)).
Applying this statement to every point of large scalar curvature at a time $t$
just prior to the blow-up time $T$,
one obtains a cover of $M$ by regions with special geometry and topology.
Any  overlaps occur in neck-like regions, 
permitting one to splice them together to form the  connected  components
with known topology  whose existence is asserted in Claim \ref{introclaim1}.  

Claim \ref{introclaim7} is proved using  a rescaling (or blow-up) argument.
This is a standard technique in geometric analysis and PDE's for treating
scale-invariant equations, such as the Ricci flow equation.  
The claim is equivalent to the statement
that if $\{(x_i, t_i)\}_{i=1}^\infty$ is a sequence of spacetime points
for which  $\lim_{i \rightarrow \infty} R(x_i, t_i) = \infty$, then 
by rescaling  the Ricci flow and passing to a subsequence, one obtains
a new sequence of Ricci flows which converges to a $\kappa$-solution.
More precisely, 
view $(x_i, t_i)$ as a new
spacetime basepoint and spatially expand the solution around $(x_i, t_i)$
by $R(x_i,t_i)^{\frac12}$.  For dimensional reasons,
in order for rescaling to produce a new Ricci flow solution
one must also expand the time factor by $R(x_i,t_i)$. The new 
Ricci flow solution, with time parameter $s$, is given by
\begin{equation}
\overline{g}_i(s) \: = \: R(x_i,t_i) \: 
g \left( R(x_i,t_i)^{-1}  \: s \: + \: t_i \right).
\end{equation}
The new time interval for $s$ is
$\left[ - \:  R(x_i,t_i) \: t_i,0 \right]$.
One would then hope to take an appropriate limit $(M_\infty, 
\overline{g}_\infty)$ of a subsequence of these rescaled solutions 
$\{(M,  \overline{g}_i(\cdot))\}_{i=1}^\infty$.
(Technically speaking,
one uses smooth convergence of sequences of Ricci flows with basepoints;
this notion of convergence allows us to focus on what happens near
the spacetime points $(x_i, t_i)$.)
Any such  limit solution $(M_\infty,  \overline{g}_\infty)$ will be an ancient
solution,  since
$\lim_{i \rightarrow \infty} -R(x_i,t_i)\,t_i = -\infty$.
Furthermore, from a $3$-dimensional
result of Hamilton and Ivey (see Appendix \ref{phiappendix}),
any limit solution will have nonnegative sectional
curvature.

Although this sounds promising,
a major  problem was to show that a limit solution
actually exists.
To prove this, one would like to invoke Hamilton's compactness
theorem \cite{Hamiltonnn}.   In the present situation, the
compactness theorem says
that the sequence  of rescaled Ricci flows  
$\{(M,  \overline{g}_i(\cdot))\}_{i=1}^\infty$
has a smoothly convergent subsequence provided two conditions are met: \\ \\
\noindent
A. For every $r>0$ and sufficiently large $i$,
the sectional curvature of $\overline{g_i}$ is  bounded uniformly independent of $i$
at each point $(x,s)$ in  spacetime such that $x$ lies in the $\overline{g_i}$-ball
$B_0(x_i,r)$ and $s\in [-r^2,0]$ and \\ \\
\noindent
B.  The  injectivity radii
$\inj(x_i, 0)$ in the time-$0$ slices of the $\overline{g_i}$'s have a 
uniform positive lower bound.

For the moment, we ignore the issue of verifying condition A,
and simply assume that it holds for the sequence 
$\{(M, \overline{g}_i(\cdot))\}_{i=1}^\infty$.
 In the presence
of the sectional curvature bounds in condition A,  a lower 
bound on the injectivity
radius is known to be equivalent to a lower bound on the volume of
metric balls.  In terms of the original Ricci flow solution, this 
becomes the condition that
\begin{equation} \label{bound}
r^{-3} \vol(B_t(x, r)) \ge \kappa \: > \: 0,
\end{equation} 
where  $B_t(x, r)$ is 
an arbitrary
metric $r$-ball in a time-$t$ slice, and the   curvature
bound   $|\Rm| \: \le \: \frac{1}{r^2}$ holds in $B_t(x,r)$.
The number $\kappa$ could depend on the
given Ricci flow solution, but the bound (\ref{bound})
should hold for all $t \in [0, T)$ and all $r < \rho$, where $\rho$ is a relevant
scale.

One of the outstanding achievements of \cite{Perelman} is to prove that
for an arbitrary Ricci flow defined on a
finite time interval, equation
(\ref{bound}) does hold with appropriate values of $\kappa$ and $\rho$.
In fact, the proof works in arbitrary dimension.
This result is called
a ``no local collapsing
theorem'' because it excludes the phenomenon of Cheeger-Gromov
collapse, in which a sequence of Riemannian manifolds has uniformly
bounded curvature, but fails to converge because the injectivity 
radii tend to zero.

One can then
apply the no local collapsing theorem to the preceding rescaling argument,
provided that one has the needed sectional curvature bounds, in order to construct the 
ancient solution $(M_\infty,\overline g_\infty)$.  
In the blowup limit the condition that $r < \rho$ goes away, and so we
can say that $(M_\infty,\overline g_\infty)$ is $\kappa$-noncollapsed
(i.e. satisfies (\ref{bound})) at all scales. In addition, in the three-dimensional
case one can show that $(M_\infty,\overline g_\infty)$ has bounded sectional
curvature. To summarize, $(M_\infty,\overline g_\infty)$ is a 
{\em $\kappa$-solution},
meaning that it is an ancient Ricci flow solution with nonnegative curvature operator on
each time slice and bounded sectional curvature
on compact time intervals,
which is $\kappa$-noncollapsed
at all scales.

With  the no local collapsing theorem in place, most of the proof of  
Claim \ref{introclaim7} is concerned with showing that in the rescaling
argument, we effectively have the needed curvature
bounds of condition A.  The argument is a tour-de-force with many ingredients,
including earlier work
of Hamilton and the theory of 
Riemannian manifolds with nonnegative sectional curvature.

\subsection{The proof of Claims \ref{introclaim2} and \ref{introclaim3}}
Claim \ref{introclaim1} allows one to take the limit of the evolving metric
$g(\cdot)$ as $t\ra T$, on the open set $\Omega$ where the metric is not
becoming singular.  It also provides geometrically defined regions -- the 
connected components of the open set $U$ -- which one removes during surgery.
Each boundary component of the resulting manifold is
a nearly round $2$-sphere with a nearly cylindrical collar,
because the collar regions in Case B of Claim \ref{introclaim1} have a neck-like
geometry.  This enables one to glue in $3$-balls with a standard 
metric, using a partition of unity  construction.

The Ricci flow starting with the postsurgery metric may also go singular after
a finite time. If so, one can appeal to Claim \ref{introclaim1} again to perform
surgery.  
However, the elapsed time between successive surgeries will depend on 
the scales at which surgeries are performed. Unless one performs
the surgeries very carefully, the surgery times may accumulate.

The way to rule out an accumulation of surgery times is to
arrange the surgery procedure so that a surgery at time $t$ removes a definite amount
of volume $v(t)$.  That is, a surgery at time $t$ should be performed at a
definite scale $h(t)$. In order to guarantee that this is possible, one needs to 
establish a
quantitative version of Claim \ref{introclaim1} for a Ricci flow with surgery,
which applies not just at the first surgery time $T$ but also at a later
surgery
time $T^\prime$. The output of this quantitative version can depend on 
the surgery time $T^\prime$
and the time-zero
metric, but it should be independent of whether or when surgeries occur 
before time $T^\prime$.

The 
general idea of the proof is similar to that of Claim \ref{introclaim1},
except that one has to carefully prescribe the surgery procedure in order
to control the effect of the earlier
surgeries. In particular, one of Perelman's
remarkable achievements is  a version of the no local collapsing theorem
for Ricci flows with surgery.

We refer the  reader to Section \ref{overview2} for a 
more detailed overview of the proof of Claim \ref{introclaim3}, and for 
further discussion of Claims \ref{introclaim5} and \ref{introclaim6}.

\section{Overview of
 {\em The Entropy Formula for the Ricci 
Flow and its Geometric
Applications} 
\cite{Perelman}}
\label{overview1}

The paper \cite{Perelman} deals with nonsingular Ricci flows;
the surgery process is considered in \cite{Perelman2}.
In particular, the final conclusion of \cite{Perelman} concerns
Ricci flows that are singularity-free and exist for all positive time. 
It does not apply
to compact $3$-manifolds with finite fundamental group or
positive scalar curvature.

The purpose of the present overview is not to give a comprehensive
summary of the results of \cite{Perelman}. Rather we indicate its
organization and the interdependence of its sections, for the
convenience of the reader.  
Some of the remarks in the overview
may only become clear once the reader has absorbed 
a portion of the detailed notes. 

Sections I.1-I.10, along with the first part of I.11, deal with Ricci flow on
$n$-dimensional manifolds.  The second part of I.11, and Sections I.12-I.13, deal
more specifically with Ricci flow on $3$-dimensional manifolds.
The main result is that geometrization holds if a compact $3$-manifold
admits a Riemannian metric which is the initial metric of a smooth
Ricci flow.  This was previously shown in \cite{Hamiltonn} under the
additional assumption that the sectional curvatures in the Ricci flow
are $O(t^{-1})$ as $t \rightarrow \infty$.

The paper \cite{Perelman} can be divided into four main parts.  

Sections I.1-I.6 construct
certain entropy-type functionals ${\mathcal F}$ and ${\mathcal W}$
that are monotonic under Ricci flow. The functional ${\mathcal W}$ is
used to prove a no local collapsing theorem.

Sections I.7-I.10 introduce and apply another monotonic quantity,
the reduced volume $\tilde{V}$. It is also used to prove a
no local collapsing theorem.  The construction of $\tilde{V}$ uses
a modified notion of a geodesic, called an ${\mathcal L}$-geodesic.
For technical reasons the reduced volume $\tilde{V}$ seems to be
easier to work with than the ${\mathcal W}$-functional, and is used
in most of the sequel.  
A reader who wants to focus on the Poincar\'e Conjecture or the
Geometrization Conjecture could
in principle start with I.7.

Section I.11 is concerned with $\kappa$-solutions, meaning nonflat
ancient solutions that are $\kappa$-noncollapsed at all scales
(for some $\kappa > 0$) and have bounded nonnegative curvature
operator on each time slice. In three dimensions, a blowup limit of a
finite-time singularity will be a $\kappa$-solution. 

Sections I.12-I.13 are about three-dimensional Ricci flow solutions.
It is shown that high-scalar-curvature regions are modeled by
rescalings of $\kappa$-solutions. A decomposition of the
time-$t$ manifold into ``thick" and ``thin" pieces is described. 
It is stated that as $t \rightarrow \infty$, the thick piece becomes
more and more hyperbolic, with incompressible cuspidal tori, and the thin
piece is a graph manifold. More details of these assertions
appear in \cite{Perelman2}, which also
deals with the necessary modifications if the solution has singularities.

We now describe each of these four parts in a bit more detail.

\subsection{I.1-I.6}

In these sections $M$ is assumed to be a closed $n$-dimensional
manifold.

A functional $F(g)$ of Riemannian metrics $g$ is said to be monotonic
under Ricci flow if $F(g(t))$ is nondecreasing in $t$ whenever $g(\cdot)$ is
a Ricci flow solution.
Monotonic quantities are an important tool for understanding Ricci flow. One wants
to have useful monotonic quantities, in particular with a characterization
of the Ricci flows for which $F(g(t))$ is constant in $t$.

Formally thinking of Ricci flow as a flow on the space of 
metrics, one way to get a monotonic quantity  would be if the Ricci flow
were the gradient flow of a functional $F$. In Sections I.1-I.2, a functional ${\mathcal F}$
is introduced whose gradient flow is not quite Ricci flow, but only
differs from the Ricci flow by the action of diffeomorphisms.  
(If one formally considers the Ricci flow as a flow on the space of
metrics modulo diffeomorphisms then it turns out to be the gradient flow of a functional
$\lambda_1$.) The functional ${\mathcal F}$ actually depends on a Riemannian metric
$g$ and a function $f$. If $g(\cdot)$ satisfies the Ricci flow equation
and $e^{-f(\cdot)}$ satisfies a conjugate or ``backward'' 
heat equation, in terms of $g(\cdot)$, then
${\mathcal F}(g(t),f(t))$ is nondecreasing in $t$. Furthermore, it is
constant in $t$ if and only if $g(\cdot)$ is a gradient steady soliton with
associated function $f(\cdot)$. Minimizing ${\mathcal F}(g,f)$ over all
functions $f$ with $\int_M e^{-f} \: dV \: = \: 1$ gives the monotonic
quantity $\lambda_1(g)$, which turns out to be the
lowest eigenvalue of $- \: 4 \: \triangle \: + \: R$. 

In Section I.3 a modified ``entropy'' functional ${\mathcal W}(g,f,\tau)$ is
introduced.  It is nondecreasing in $t$ provided that $g(\cdot)$ is a Ricci flow,
$\tau \: = \: t_0 \: - \: t$ and $(4 \pi \tau)^{- \: \frac{n}{2}} \:
e^{-f}$ satisfies the conjugate heat equation. The functional
${\mathcal W}$ is constant on a Ricci flow if and only if the flow is a gradient
shrinking soliton that terminates at time $t_0$. Because of this,
${\mathcal W}$ is more suitable than ${\mathcal F}$ when one wants
information that is localized in spacetime.

In Section I.4 the entropy functional ${\mathcal W}$ is used to prove
a no local collapsing theorem. The statement is that if $g(\cdot)$ is a given
Ricci flow on a finite time interval $[0, T)$ then for any (scale) $\rho > 0$,
there is a number $\kappa > 0$ so that if $B_t(x,r)$ is a time-$t$ ball
with radius $r$ less than $\rho$, on which $|\Rm| \: \le \: \frac{1}{r^2}$, 
then $\vol(B_t(x,r)) \: \ge \:
\kappa \: r^n$. The method of proof is to show that if 
$r^{-n} \: \vol(B_t(x,r))$ is very small then the evaluation of ${\mathcal W}$
at time $t$ is very negative, which contradicts the monotonicity of
${\mathcal W}$.

The significance of a no local collapsing theorem is that it allows one to 
use Hamilton's compactness theorem to construct blowup limits of finite time
singularities, and more generally to understand high curvature regions.

Section I.5 and I.6 are not needed in the sequel.  Section I.5 gives some
thermodynamic-like equations in which ${\mathcal W}$
appears as an entropy.  Section I.6 motivates the construction of the
reduced volume of Section I.7. 

\subsection{I.7-I.10}

A new monotonic quantity, the reduced volume $\tilde{V}$, 
is introduced in I.7.  It is defined in terms of so-called
${\mathcal L}$-geodesics. Let $(p, t_0)$ be a fixed spacetime point.
Define backward time by $\tau = t_0 -t$. Given a curve
$\gamma(\tau)$ in $M$ defined for $0 \le \tau \le \overline{\tau}$ (i.e.
going backward in real time) with $\gamma(0) \: = \: p$, its
${\mathcal L}$-length is 
\begin{equation}
{\mathcal L}(\gamma) = \int_0^{\overline{\tau}} \sqrt{\tau} \:
\left( |\dot{\gamma}(\tau)|^2 \: + \: R(\gamma(\tau),t_0 - \tau) \right)
\: d\tau.
\end{equation}
One can compute the first and second variations of ${\mathcal L}$, in
analogy to what is done in Riemannian geometry.  

Let $L(q, \overline{\tau})$
be the infimum of ${\mathcal L}(\gamma)$ over curves $\gamma$ with
$\gamma(0) = p$ and
$\gamma(\overline{\tau}) = q$. Put $l(q, \overline{\tau}) \: = \:
\frac{L(q, \overline{\tau})}{2\sqrt{\overline{\tau}}}$. The reduced volume
is defined by $\tilde{V}(\overline{\tau}) = 
\int_M \overline{\tau}^{- \frac{n}{2}} \: e^{-l(q, \overline{\tau})} \:
\dvol(q)$. The remarkable fact is that if $g(\cdot)$ is a Ricci flow
solution then $\tilde{V}$ is nonincreasing in $\overline{\tau}$, i.e.
nondecreasing in real time $t$. Furthermore, it is constant if and only if
$g(\cdot)$ is a gradient shrinking soliton that terminates at time
$t_0$. The proof of monotonicity uses a subtle cancellation between
the $\overline{\tau}$-derivative of $l(\gamma(\overline{\tau}), \overline{\tau})$
along an ${\mathcal L}$-geodesic and the Jacobian of the
so-called ${\mathcal L}$-exponential map.

Using a differential inequality, it is shown that for each $\overline{\tau}$ there
is some point $q(\overline{\tau}) \in M$ so that
$l(q(\overline{\tau}), \overline{\tau}) \: \le \: \frac{n}{2}$. This is then
used to prove a no local collapsing theorem :  Given a Ricci flow
solution $g(\cdot)$ 
defined on a finite time interval $[0, t_0]$
and a scale $\rho > 0$ there
is a number $\kappa > 0$
with the following property. For $r < \rho$, suppose that
$|\Rm| \le \frac{1}{r^2}$ on the
``parabolic'' ball $\{(x,t) \: : \: \dist_{t_0}(x, p) \: \le \: r, \: \:
t_0 - r^2 \: \le \: t \: \le t_0\}$. Then
$\vol(B_{t_0}(p, r)) \: \ge \: \kappa \: r^n$. The number $\kappa$ can be chosen
to depend on $\rho$, $n$, $t_0$ and bounds on the geometry of the initial metric
$g(0)$.

The hypotheses of the no local collapsing theorem proved using $\tilde{V}$ are
more stringent than those of the no local collapsing theorem proved using $\mathcal{W}$,
but the consequences turn out to be the same.  The no local collapsing theorem is
used extensively in Sections I.11 and I.12 when extracting a convergent
subsequence from a sequence of Ricci flow solutions.

Theorem I.8.2 of Section I.8 says that under appropriate
assumptions, local $\kappa$-noncollapsing extends forwards in
time to larger distances.  
This will be used in I.12 to analyze long-time behavior.
The statement is that for any $A < \infty$ there is a
$\kappa = \kappa(A) > 0$ with the following property. Given $r_0 > 0$ and 
$x_0 \in M$, suppose that
$g(\cdot)$ is defined for $t \in [0, r_0^2]$, with
$|\Rm(x,t)| \: \le \: r_0^{-2}$ for all $(x,t)$ satisfying
$\dist_0(x, x_0) < r_0$, and the volume of the time-zero ball
$B_0(x_0, r_0)$ is at least $A^{-1} r_0^n$. Then the metric
$g(t)$ cannot be $\kappa$-collapsed on scales less than $r_0$
at a point $(x, r_0^2)$ with  $\dist_{r_0^2}(x, x_0) \le Ar_0$.

A localized version of the ${\mathcal W}$-functional appears in I.9.
Section I.10, 
which is not needed for the sequel
but is of independent
interest, shows if a point in a time slice 
lies in a ball with quantitatively bounded geometry then 
at nearby later times, the curvature at the point is quantitatively
bounded.  That is, there is a damping effect with regard to exterior high-curvature 
regions. The ``bounded geometry'' assumptions on the initial
ball are a lower bound on its scalar curvature and an assumption that
the isoperimetric constants of subballs are close to the Euclidean
value.

\subsection{I.11}

Section I.11 contains an analysis of $\kappa$-solutions. As mentioned before, in
three dimensions $\kappa$-solutions arise as blowup limits of finite-time singularities and,
more generally, as limits of rescalings of high-scalar-curvature regions.

In addition to the no local collapsing theorem, some of the tools used to analyze
$\kappa$-solutions are Hamilton's Harnack inequality for Ricci flows with
nonnegative curvature operator, and the comparison geometry of
nonnegatively curved manifolds.

The first result is that any time slice of a $\kappa$-solution
has vanishing asymptotic volume ratio $\lim_{r \rightarrow \infty} r^{-n} \vol(B_t(p, r))$.  
This apparently technical result is used to show that if a
$\kappa$-solution $(M, g(\cdot))$ has scalar curvature normalized to equal
one at some spacetime point $(p,t)$ then there is an {\it a priori} upper bound on
the scalar curvature $R(q,t)$ at other points $q$ in terms of $\dist_t(p,q)$.
Using the curvature
bound, it is shown that a sequence $\{(M_i, p_i, g_i(\cdot)) \}_{i=1}^\infty$ of pointed 
$n$-dimensional $\kappa$-solutions, normalized so that $R(p_i, 0) =1$ for each $i$, has
a convergent subsequence whose limit satisfies all of the requirements to be a
$\kappa$-solution, except possibly the condition of having bounded
sectional curvature on each time slice.

In three dimensions this statement is improved by showing that
the sectional curvature will be bounded on
each 
compact time interval,
so the space of pointed $3$-dimensional $\kappa$-solutions
$(M, p, g(\cdot))$ with $R(p,0) = 1$ is in fact compact. 
This is used to draw conclusions about the global geometry of
$3$-dimensional $\kappa$-solutions.

If $M$ is a compact $3$-dimensional $\kappa$-solution
then Hamilton's theorem about
compact $3$-manifolds with nonnegative Ricci curvature implies that $M$
is diffeomorphic to a spherical space form. If $M$ is noncompact then assuming
that $M$ is oriented, it follows easily that $M$ is
diffeomorphic to $\R^3$, or isometric to the round shrinking cylinder $\R \times S^2$
or its $\Z_2$-quotient $\R \times_{\Z_2} S^2$.

In the noncompact case it is shown that after rescaling, each time slice is neck-like at
infinity.  More precisely, considering a given time-$t$ slice, for each $\epsilon > 0$ 
there is a compact subset $M_\epsilon \subset M$ so that if $y \notin M_\epsilon$
then the pointed manifold $(M, y, R(y,t) g(t))$ is $\epsilon$-close
to the standard cylinder $\left[ - \: \frac{1}{\epsilon}, \frac{1}{\epsilon} \right] \times
S^2$ of scalar curvature one.

\subsection{I.12-I.13}

Sections I.12 and I.13 deal with three-dimensional Ricci flows.

Theorem I.12.1 uses the results of Section I.11 to model the
high-scalar-curvature regions of a Ricci flow. 
Let us assume a
pinching condition of the form $\Rm \: \ge \: - \: \Phi(R)$ for an
appropriate function $\Phi$ with $\lim_{s \rightarrow \infty}
\frac{\Phi(s)}{s} = 0$.
(This will eventually follow from Hamilton-Ivey pinching, 
 cf. Appendix \ref{phiappendix}.)
Theorem I.12.1 says that given
numbers $\epsilon, \kappa > 0$, one can find $r_0 > 0$ with the
following property. Suppose that $g(\cdot)$ is a Ricci flow solution
defined on some time interval $[0, T]$ that satisfies the
pinching condition and is $\kappa$-noncollapsed at scales less than
one.   Then for any point $(x_0, t_0)$ with $t_0 \ge 1$ and
$Q = R(x_0, t_0) \ge r_0^{-2}$, after scaling by the factor $Q$, the
solution in the region $\{(x,t) \: : \: \dist_{t_0}^2(x, x_0) \le (\epsilon Q)^{-1},
t_0 - (\epsilon Q)^{-1} \le t \le t_0\}$ is $\epsilon$-close to the 
corresponding subset of a $\kappa$-solution.

Theorem I.12.1 says in particular that near a first singularity, the geometry is modeled
by a $\kappa$-solution, for some $\kappa$. This fact is used in
\cite{Perelman2}. Although Theorem I.12.1 is not
used directly in \cite{Perelman}, its method of proof is used in Theorem I.12.2.

The method of proof of Theorem I.12.1 is by contradiction. If it were not
true then there would be a sequence $r_0^{(i)} \rightarrow 0$ and 
a sequence 
$(M_i, g_i(\cdot))$ of
Ricci flow solutions that satisfy the assumptions, each with a spacetime point
$\left( x^{(i)}_0, t^{(i)}_0 \right)$ that does not satisfy the conclusion. To consider
first a special
case,
suppose that each point $\left( x^{(i)}_0, t^{(i)}_0 \right)$ 
is the first point at which a certain curvature threshold $R_i$ is
achieved, i.e.
$R(y, t) \: \le \: R \left( x^{(i)}_0, t^{(i)}_0 \right)$ for each
$y \in M_i$ and $t \in \left[ 0,  t^{(i)}_0 \right]$.
Then after rescaling the Ricci flow
$g_i(\cdot)$ by $Q_i \: = \: R \left( x^{(i)}_0, t^{(i)}_0 \right)$ and shifting the
time parameter, one has the curvature bounds on the time interval
$[- Q_i t^{(i)}_0,0]$ that form part of  the
hypotheses of Hamilton's compactness theorem.  Furthermore, the
no local collapsing theorem gives the lower injectivity radius bound
needed to apply Hamilton's theorem and take a convergent subsequence of the pointed 
rescaled solutions.
The limit will be a $\kappa$-solution, giving the contradiction.

In the general case, one effectively proceeds by induction on the size of the
scalar curvature.  By modifying the choice of points $\left( x^{(i)}_0, t^{(i)}_0 \right)$,
one can assume that the conclusion of the theorem holds for all of the points $(y,t)$ in a
large spacetime neighborhood of $\left( x^{(i)}_0, t^{(i)}_0 \right)$
that have $R(y,t) > 2Q_i$. One then shows that one has the curvature bounds
needed to form the time-zero slice of the putative $\kappa$-solution. One shows
that this ``time-zero'' metric can be extended backward in time to form a
$\kappa$-solution, thereby giving the contradiction.

The rest of Section I.12 begins the analysis of the long-time behaviour of
a nonsingular $3$-dimensional Ricci flow. There are two main results,
Theorems I.12.2 and I.12.3. They extend curvature bounds forward and
backward in time, respectively.  

Theorem I.12.2 roughly says that if
one has $|\Rm| \le r_0^{-2}$ on a spacetime region of spatial
size $r_0$ and temporal size $r_0^2$, and if one has a lower bound on the
volume of the initial time face of the region, then one gets scalar curvature
bounds on much larger spatial balls at the final time. More precisely,
for any $A < \infty$
there are numbers $K = K(A) < \infty$ and
$\rho = \rho(A) < \infty$ so that with the
hypotheses and notation of Theorem I.8.2, if in addition 
$r_0^2 \: \Phi(r_0^{-2}) \: < \: \rho$
then
$R(x, r_0^2) \le K r_0^{-2}$ for points $x$ lying in the ball of 
radius $Ar_0$ 
around $x_0$
at time $r_0^2$.

Theorem I.12.3 says that if one has a lower bound on volume and sectional
curvature on a ball at a certain time then one obtains an upper scalar
curvature bound on a smaller ball at an earlier time.  More precisely,
given $w > 0$ there exist $\tau = \tau(w) > 0$, $\rho = \rho(w) > 0$ and
$K = K(w) < \infty$ with the following property.  Suppose that 
a ball
$B(x_0, r_0)$ at time $t_0$ has volume bounded below by
$w r_0^3$ and sectional curvature bounded below by $- \: r_0^{-2}$. Then
$R(x,t) < K r_0^{-2}$ for
$t \in [t_0 - \tau r_0^2, t_0]$ and $\dist_t(x,x_0) \: < \: \frac14 r_0$,
provided that $\phi(r_0^{-2}) < \rho$. 

Applying a back-and-forth argument using Theorems I.12.2 and I.12.3,
along with the pinching condition, one concludes, roughly speaking,
that if a metric ball of small radius $r$ has infimal sectional curvature
exactly equal to $- \: r^{-2}$ then the ball has a small volume
compared to $r^3$. Such a ball can be said to be locally
volume-collapsed with respect to a lower sectional curvature bound.

Section I.13 defines the thick-thin decomposition of a large-time slice
of 
a nonsingular
Ricci flow and shows the geometrization. 
Rescaling the metric to $\widehat{g}(t) \: = \: t^{-1} \: g(t)$,
there is a universal function $\Phi$ so that for large $t$, the metric
$\widehat{g}(t)$ satisfies the $\Phi$-pinching condition. In terms of the
original unscaled metric, given $x \in M$ let $\widehat{r}(x,t) > 0$ be the unique
number such that $\inf \Rm \big|_{B_t(x,\widehat{r})} \: = \: - \:
\widehat{r}^{-2}$. 

Given $w > 0$, define the $w$-thin part 
$M_{thin}(w,t)$ of the time-$t$ slice to be the points $x \in M$ so that
$\vol(B_t(x, \widehat{r}(x,t))) \: < \: w \: \widehat{r}(x,t)^3$. 
That is, a point
in $M_{thin}(w,t)$ lies in a ball that is locally volume-collapsed
with respect to a lower sectional curvature bound. Put
$M_{thick}(w,t) = M - M_{thin}(w,t)$. One shows that for large $t$, the subset
$M_{thick}(w,t)$ has bounded geometry in the sense that there
are numbers $\overline{\rho} = \overline{\rho}(w) > 0$ and 
$K = K(w) < \infty$ so that $|\Rm| \: \le \: K t^{-1}$ on 
$B(x, \overline{\rho} \sqrt{t})$ and 
$\vol(B(x, \overline{\rho} \sqrt{t})) \: \ge \: \frac{1}{10} \: w \: (\overline{\rho} \sqrt{t}))^3$,
whenever $x \in M_{thick}(w,t)$. 

Invoking arguments of Hamilton (that are written out in more detail in
\cite{Perelman2}) one can take a sequence $t \rightarrow \infty$ and
$w \rightarrow 0$ so that $M_{thick}(w,t)$ converges to a complete finite-volume
manifold with constant sectional curvature $- \: \frac14$, whose cuspidal
tori are incompressible in $M$. On the other hand, a result from Riemannian
geometry implies that for large $t$ and small $w$, $M_{thin}(w,t)$ is homeomorphic to
a graph manifold; again a more precise statement appears in \cite{Perelman2}.
The conclusion is that $M$ satisfies the geometrization conjecture. Again, one is
assuming in 
I.13
that the Ricci flow is nonsingular
for all times.

\section{I.1.1. The ${\mathcal F}$-functional and its
monotonicity} \label{I.1.1}
\label{Derivation of first equation of Section 1.1}

The goal of this section is to show that in an appropriate sense, Ricci
flow is a gradient flow on the space of metrics. 
We introduce the entropy functional ${\mathcal F}$.
We compute its formal variation and show that the corresponding
gradient flow is a modified Ricci flow.

In Sections \ref{I.1.1} through \ref{I.5} of these notes,  $M$ 
is a closed manifold.  
We will use the Einstein summation convention freely. We also follow
Perelman's convention that a condition like $a > 0$ means that
$a$ should be considered to be a small parameter, while a condition
like $A < \infty$ means that $A$ should be considered to be a large
parameter. This convention is only for pedagogical purposes and
may be ignored by a logically minded reader.

Let ${\mathcal M}$ denote the space of smooth Riemannian metrics $g$ on
$M$. We think of ${\mathcal M}$ formally as an infinite-dimensional
manifold.  The tangent space $T_g {\mathcal M}$ consists of the
symmetric covariant $2$-tensors $v_{ij}$ on $M$. Similarly, $C^\infty(M)$
is an infinite-dimensional manifold with $T_f C^\infty(M) = 
C^\infty(M)$. 
The diffeomorphism group
$\Diff(M)$ acts on ${\mathcal M}$ and $C^\infty(M)$ by pullback.

Let $dV$ denote the Riemannian volume density associated to a metric $g$.
We use the convention that
$\triangle \: = \: \Div \: \grad$.

\begin{definition}
The ${\mathcal F}$-functional ${\mathcal F} \: : \: {\mathcal M} \times C^\infty(M)
\rightarrow \R$ is given by
\begin{equation}
{\mathcal F}(g, f) \: = \:
\int_M \left( R \: + \: |\nabla f|^2 \right) \: e^{-f} \: dV .
\end{equation}
\end{definition}

Given $v_{ij} \in T_g{\mathcal M}$ and $h \in T_f C^\infty(M)$, the
evaluation of the differential $d{\mathcal F}$ on $(v_{ij}, h)$ is
written as $\delta {\mathcal F}(v_{ij}, h)$.
Put $v \: = \: g^{ij} \: v_{ij}$.
\begin{proposition} (cf. I.1.1)
We have 
\begin{equation}
\delta {\mathcal F}(v_{ij}, h) \: = \:
\int_M e^{-f} \: \left[ - \: v_{ij} 
(R_{ij} \: + \: \nabla_i \nabla_j f) \: + \: 
\left( \frac{v}{2} \: - \: h \right) \: (2 \triangle f \: - \:
|\nabla f|^2 \: + \: R) \right] \: dV.
\end{equation}
\end{proposition}
\begin{proof}
From a standard formula,
\begin{equation}
\delta R \: = \: - \: \triangle v \: + \: \nabla_i \nabla_j v_{ij}
\: - \: R_{ij} v_{ij}.
\end{equation}
As 
\begin{equation}
|\nabla f|^2 \: = \: g^{ij} \nabla_i f \nabla_j f,
\end{equation}
we have
\begin{equation}
\delta |\nabla f|^2 \: = \: - \: v^{ij} \nabla_i f \nabla_j f
\: + \: 2 \: \langle \nabla f, \nabla h \rangle.
\end{equation}
As $dV \: = \: \sqrt{\det(g)} \: dx_1 \ldots dx_n$, we have
$\delta (dV) \: = \: \frac{v}{2} \: dV$, so
\begin{equation} \label{measvar}
\delta \left( e^{-f} \: dV \right) \: = \: \left( \frac{v}{2} \: - \: h
\right) \: e^{-f} \: dV.
\end{equation}
Putting this together gives
\begin{align} \label{putting}
\delta {\mathcal F} \: = \: 
\int_M \:  e^{-f} \: & \left[  - \: \triangle v \: + \:
\nabla_i \nabla_j v_{ij} \: - \: R_{ij} v_{ij} \: - \:
v_{ij} \nabla_i f \nabla_j f \: + \right. \\
& \left.  2 \: \langle \nabla f, \nabla h \rangle
\: + \: (R \: + \: |\nabla f|^2 ) \: \left( \frac{v}{2} \: - \: h \right)
\right] \:  dV. \notag
\end{align}

The goal now is to rewrite the right-hand side of
(\ref{putting}) so that $v_{ij}$ and $h$ appear
algebraically, i.e. without derivatives.
As 
\begin{equation} \label{id}
\triangle e^{-f} \: = \: (|\nabla f|^2 \: - \: \triangle f) \: e^{-f},
\end{equation}
we have
\begin{equation}
\int_M e^{-f} \: [ - \: \triangle v ] \: dV \: = \:
- \int_M (\triangle e^{-f} )\: v \: dV \: = \:
\int_M e^{-f} \: (\triangle f \: - \: |\nabla f|^2 ) \: v \: dV.
\end{equation}
Next,
\begin{align}
\int_M e^{-f} \nabla_i \nabla_j v_{ij} \: dV \: & = \:
\int_M (\nabla_i \nabla_j e^{-f}) \: v_{ij} \: dV \: = \:
- \: \int_M \nabla_i ( e^{-f} \nabla_j f) \: v_{ij} \: dV \\
& = \:
\int_M e^{-f} (\nabla_i f \nabla_j f \: - \: \nabla_i \nabla_j f) 
\: v_{ij} \: dV. \notag
\end{align}
Finally,
\begin{align}
2 \: \int_M e^{-f} \langle \nabla f, \nabla h \rangle \: dV \: & = \:
- \: 2 \: \int_M  \langle \nabla e^{-f}, \nabla h \rangle \: dV \: = \:
2 \: \int_M  (\triangle e^{-f}) \: h \: dV \\
& = \:
2 \: \int_M e^{-f} \:
(|\nabla f|^2 \: - \: \triangle f) \: h \: dV. \notag
\end{align}
Then 
\begin{align} \label{Fvariation}
\delta {\mathcal F} \: = \:
\int_M e^{-f} \: & \left[ 
\left( \frac{v}{2} \: - \: h \right) (2 \triangle f \: - \:
2 |\nabla f|^2 ) \: - \: v_{ij} (R_{ij} \: + \: \nabla_i \nabla_j f)
\right. \\
& 
\left.  + \: \left( \frac{v}{2} \: - \: h \right) \: 
(R \: + \: |\nabla f|^2 )
\right] \: dV \notag \\
=  \: \int_M e^{-f} \: & \left[ - \: v_{ij} 
(R_{ij} \: + \: \nabla_i \nabla_j f) \: + \: 
\left( \frac{v}{2} \: - \: h \right) \: (2 \triangle f \: - \:
|\nabla f|^2 \: + \: R) \right] \: dV. \notag
\end{align}
This proves the proposition.
\end{proof}

We would like to get rid of the $\left( \frac{v}{2} \: - \: h \right) \: (2 \triangle f \: - \:
|\nabla f|^2 \: + \: R)$ term in (\ref{Fvariation}).  We can do this by restricting
our variations so that $\frac{v}{2} \: - \: h \: = \: 0$. From
(\ref{measvar}), this amounts to assuming that assuming $e^{-f} \: dV$ is fixed.
We now fix a smooth measure $dm$ on $M$ and relate $f$ to $g$ by
requiring that $e^{-f} \: dV \: = \: dm$. Equivalently, we define
a section $s \: : \: {\mathcal M} \rightarrow {\mathcal M} \times C^\infty(M)$
by $s(g) \: = \: \left( g, \ln \left( \frac{dV}{dm}\right) \right)$.
Then the composition ${\mathcal F}^m \: = \: 
{\mathcal F}\circ s$ is a function on ${\mathcal M}$ and its differential
is given by 
\begin{equation}
d{\mathcal F}^m (v_{ij}) \: = \: \int_M e^{-f} \: \left[ - \: v_{ij} 
(R_{ij} \: + \: \nabla_i \nabla_j f) \right] \: dV.
\end{equation}

Defining a formal Riemannian metric on ${\mathcal M}$ by
\begin{equation}
\langle v_{ij}, v_{ij} \rangle_{g} \: = \: \frac12
\int_M v^{ij} \: v_{ij} \: dm, 
\end{equation}
the gradient flow of ${\mathcal F}^m$ on ${\mathcal M}$
is given by
\begin{equation} \label{Fflow}
(g_{ij})_t \: = \: -2 \: (R_{ij} \: + \: \nabla_i \nabla_j f).
\end{equation}
The induced flow equation for $f$ is
\begin{equation}
f_t \: = \: \frac{\partial}{\partial t} \: \ln \left( \frac{dV}{dm}\right) 
\: = \: \frac12 \: g^{ij} \: (g_{ij})_t \: = \: - \: \triangle f \: - \: R.
\end{equation}
As with any gradient flow, the function ${\mathcal F}^m$ is nondecreasing
along the flow line with its derivative being given by the
length squared of the gradient, i.e.
\begin{equation} \label{Fchange}
{\mathcal F}^m_t \: = \: 2 \int_M |R_{ij} \: + \: \nabla_i \nabla_j f|^2 \: dm,
\end{equation}
as follows from (\ref{Fvariation}) and (\ref{Fflow})

We now perform time-dependent 
diffeomorphisms to transform (\ref{Fflow}) into the
Ricci flow equation. If $V(t)$ is the time-dependent generating vector field
of the diffeomorphisms then the new equations for $g$ and $f$ become
\begin{align}
(g_{ij})_t \: & = \: -2 \: (R_{ij} \: + \: \nabla_i \nabla_j f) \: + \:
{\mathcal L}_V g, \\
f_t & \: = \: - \: \triangle f \: - \: R \: + \: {\mathcal L}_V f. \notag
\end{align}
Taking $V \: = \: \nabla f$ gives
\begin{align} \label{newevolution}
(g_{ij})_t \: & = \: -2 \: R_{ij}, \\
f_t & \: = \: - \: \triangle f \: - \: R \: + \: |\nabla f|^2. \notag
\end{align}
As ${\mathcal F}(g, f)$ is unchanged by a simultaneous pullback of $g$ and $f$,
and the right-hand side of (\ref{Fchange}) is also unchanged under a
simultaneous pullback, it follows that 
under the new evolution equations (\ref{newevolution})
we still have
\begin{equation} \label{Fder}
\frac{d}{dt} \: {\mathcal F}(g(t), f(t)) \: = \:
2 \int_M |R_{ij} \: + \: \nabla_i \nabla_j f|^2 \: e^{-f} \: dV
\end{equation}
(This can also be checked directly).

Because of the diffeomorphisms that we applied,
$g$ and $f$ are no longer related by
$e^{-f} \: dV \: = \: dm$. We do have that
$\int_M e^{-f} \: dV$ is constant in $t$, as 
$e^{-f} \: dV$ is related to $dm$ by a diffeomorphism.

The relation between $g$ and $f$ is as follows: 
we solve (or assume that we have a solution for) 
the first equation in (\ref{newevolution}), with some initial
metric.
Then given the solution $g(t)$, we require that $f$ satisfy 
the second equation in (\ref{newevolution}) (which is in
terms of $g(t)$). 

The second equation in
(\ref{newevolution}) can be written as
\begin{equation} \label{backwards}
\frac{\partial}{\partial t} \: e^{-f} \: = \: - \: \triangle e^{-f}
\: + \: R \: e^{-f}.
\end{equation}
As this is a backward heat equation, we cannot solve for $f$
forward in time starting with an arbitrary smooth function.
Instead, (\ref{newevolution}) will be applied by
starting with a solution for $(g_{ij})_t \: = \: -2 \: R_{ij}$ on
some time interval $[t_1, t_2]$ and then solving
(\ref{backwards}) backwards in time on $[t_1, t_2]$
(which can always be done) starting with some initial
$f(t_2)$. Having done this, the solution $(g(t), f(t))$ on $[t_1, t_2]$
will satisfy (\ref{Fder}).

\section{Basic example for I.1}

In this section we compute ${\mathcal F}$ in a Euclidean example.

Consider $\R^n$ with the standard metric, constant in time.
Fix $t_0 \: > \: 0$. Put $\tau \: = \: t_0 \: - \: t$ and
\begin{equation}
f(t, x) \: = \: \frac{|x|^2}{4\tau} \: + \: \frac{n}{2} \ln(4 \pi \tau),
\end{equation}
so
\begin{equation}
e^{-f} \: = \: (4 \pi \tau)^{-n/2} \: e^{- \frac{|x|^2}{4\tau}}.
\end{equation}
This is the standard heat kernel when considered for $\tau$ going from
$0$ to $t_0$, i.e. for $t$ going from $t_0$ to $0$.
One can check that $(g,f)$ solves (\ref{newevolution}).
As
\begin{equation} \label{Gaussianintegral}
\int_{\R^n} e^{- \frac{|x|^2}{4\tau}} \: dV \: = \: (4 \pi \tau)^{n/2}, 
\end{equation}
$f$ is properly normalized.
Then $\nabla f \: = \: \frac{{\bf x}}{2 \tau}$ and 
$|\nabla f|^2 \: = \: \frac{|x|^2}{4 \tau^2}$.
Differentiating (\ref{Gaussianintegral}) with respect to $\tau$ gives
\begin{equation} \label{Gaussianintegral2}
\int_{\R^n} \frac{|x|^2}{4\tau^2} \: 
e^{- \frac{|x|^2}{4\tau}} \: dV \: = \: (4 \pi \tau)^{n/2} \: 
\frac{n}{2\tau}, 
\end{equation}
so
\begin{equation}
\int_{\R^n} |\nabla f|^2 \: e^{- f} \: dV \: = \: \frac{n}{2\tau}.
\end{equation}
Then ${\mathcal F}(t) \: = \: \frac{n}{2\tau} \: = \: \frac{n}{2(t_0-t)}$.
In particular, this is nondecreasing as a function of 
$t \in [0, t_0)$.

\section{I.2.2. The $\lambda$-invariant and its applications} 
\label{I.2.2}

In this section we define $\lambda(g)$ and show that it is nondecreasing under
Ricci flow.  We use this to show that a steady breather on a compact manifold
is a gradient steady soliton.

\begin{proposition}
Given a metric $g$, there is a unique minimizer 
$\overline{f}$ of ${\mathcal F}(g, f)$ under the constraint
$\int_M e^{-f} \: dV \: = \: 1$. 
\end{proposition}
\begin{proof}
Write
\begin{equation}
{\mathcal F} \: = \: \int_M \left( R e^{-f} \: + \: 4 |\nabla e^{-f/2}|^2
\right) \: dV.
\end{equation}
Putting $\Phi \: = \: e^{- f/2}$,
\begin{equation}
{\mathcal F} \: = \: \int_M \left( 4 |\nabla \Phi|^2 \: + \: R \: \Phi^2
\right) \: dV \: = \:
\int_M \Phi \left( - \: 4 \triangle \Phi \: + \: R \Phi
\right) \: dV. 
\end{equation}
The constraint equation becomes $\int_M \Phi^2 \: dV \: = \: 1$.
Then $\lambda$ is the smallest eigenvalue of
$- \: 4 \triangle \: + \: R$ and $e^{-\overline{f}/2}$ 
is a corresponding normalized
eigenvector.  As the operator is a Schr\"odinger operator, there is a
unique normalized positive eigenvector \cite[Chapter XIII.12]{Reed-Simon}.
\end{proof}

\begin{definition}
The $\lambda$-functional is given by
$\lambda(g) \: =\: {\mathcal F}(g, \overline{f})$.
\end{definition}

If $g(t)$ is a smooth family of metrics then it follows from eigenvalue
perturbation theory that $\lambda(g(t))$ and $\overline{f}(t)$ are
smooth in $t$ \cite[Chapter XII]{Reed-Simon}.

\begin{proposition} (cf. I.2.2) \label{mono}
If $g(\cdot)$ is a Ricci flow solution then $\lambda(g(t))$ is
nondecreasing in $t$.
\end{proposition}
\begin{proof}
Consider a time
interval $[t_1, t_2]$, and the minimizer $\overline{f}(t_2)$. In
particular, $\lambda(t_2) \: = \: {\mathcal F}(g(t_2), 
\overline{f}(t_2))$. Put $u(t_2) \: = \: e^{-\overline{f}(t_2)}$. Solve the
backward heat equation 
\begin{equation}
\frac{\partial u}{\partial t} \: = \: - \: \triangle u \: + \: Ru
\end{equation}
backward on $[t_1, t_2]$. 

We claim that
$u(x^\prime,t^\prime) > 0$ 
for all $x^\prime \in M$ and $t^\prime \in [t_1, t_2]$. To see this,
we take $t^\prime \in [t_1, t_2)$, 
and let $h$ be the solution to the forward heat equation
$\frac{\partial h}{\partial t} \: = \: \triangle h$ on $(t^\prime, t_2]$
with $\lim_{t \rightarrow t^\prime} h(t) \: = \: \delta_{x^\prime}$.
We have
\begin{equation}
\frac{d}{dt} \: \int_M u(t) \: h(t) \: dV \: = \:
\int_M \left[ (\partial_t u \: + \: \triangle u \: - \: Ru) \: v \: + \: u \:
(\partial_t h \: - \: \triangle h) \right] \: dV \: = \: 0.
\end{equation}
One knows that $h(t) > 0$ for all $t \in (t^\prime, t_2]$. Then
\begin{equation}
u(x^\prime, t^\prime) \: = \: \int_M u(x, t^\prime) \:\delta_{x^\prime}(x) \:
dV(x) \: = \: \lim_{t \rightarrow t^\prime} 
\int_M u(t) \: h(t) \:
dV \: = \:
\int_M u(t_2) \: h(t_2) \: dV \: > \: 0.
\end{equation} 

For $t \in [t_1, t_2]$, define $f(t)$ by $u(t) \: = \: e^{-f(t)}$.
By (\ref{Fder}),
$
{\mathcal F}(g(t_1), f(t_1)) \: \le \: 
{\mathcal F}(g(t_2), f(t_2))$. By the definition of
$\lambda$, $\lambda(t_1) \: = \: {\mathcal F}(g(t_1), \overline{f}(t_1)) 
\: \le \: {\mathcal F}(g(t_1), f(t_1))$. (We are using the fact that
$\int_M e^{-f(t_1)} \: dV(t_1) \: = \:
\int_M e^{-f(t_2)} \: dV(t_2) \: = \: 1$.)
Thus $\lambda(t_1) \: \le \: \lambda(t_2)$.
\end{proof}

\begin{definition}
A {\em steady breather} is a Ricci flow solution on an interval
$[t_1, t_2]$ that satisfies the equation
$g(t_2) \: = \: \phi^* g(t_1)$ for some $\phi \in \Diff(M)$.
\end{definition}

Steady soliton solutions are steady breathers. 

Again, we are assuming that $M$ is compact. The next result is not
essential for the sequel, but gives a good illustration of how a
monotonicity formula is used.

\begin{proposition}  (cf. I.2.2) \label{steadybreather}
A steady breather is a gradient steady soliton.
\end{proposition}
\begin{proof}
We have $\lambda(g(t_2)) \: = \: \lambda(\phi^* g(t_1)) \: = \:
\lambda(g(t_1))$. Thus we have equality in Proposition
\ref{mono}. Tracing through the proof, 
${\mathcal F}(g(t), f(t))$ must be constant in $t$. From
(\ref{Fder}), $R_{ij} \: + \: \nabla_i \nabla_j f \: = \: 0$.
Then $R \: + \: \triangle f \: = \: 0$ and so
(\ref{newevolution}) becomes (\ref{steady}). 
\end{proof}

One can sharpen Proposition \ref{mono}.

\begin{lemma} (cf. Proposition I.1.2) \label{2/3}
\begin{equation}
\frac{d\lambda}{dt} \: \ge \: \frac{2}{n} \: \lambda^2(t).
\end{equation}
\end{lemma}
\begin{proof}
Given a time interval $[t_1, t_2]$, with the notation of the proof of Proposition \ref{mono} we have
\begin{align}
\lambda(t_1) \: & \le \: {\mathcal F}(g(t_1), f(t_1)) \: = \:
{\mathcal F}(g(t_2), f(t_2)) \: - \: 2 \: \int_{t_1}^{t_2} \int_M |R_{ij} + \nabla_i \nabla_j f|^2 \:
e^{-f} \: dV \: dt \\
& = \: \lambda(t_2)  \: - \: 2 \: \int_{t_1}^{t_2} \int_M |R_{ij} + \nabla_i \nabla_j f|^2 \:
e^{-f} \: dV \: dt. \notag
\end{align}
Then
\begin{equation}
\frac{d\lambda}{dt} \Big|_{t=t_2}
\: = \: \lim_{t_1 \rightarrow t_2^-} \frac{\lambda(t_2) - \lambda(t_1)}{t_2 - t_1} \: \ge \: 
2 \int_M |R_{ij} \: + \: \nabla_i \nabla_j {f}|^2 \: e^{- {f}}
\: dV,
\end{equation}
where the right-hand side is evaluated at time $t_2$.

Hence for all $t$,
\begin{equation} \label{dlambdadt}
\frac{d\lambda}{dt} \: \ge \:
2 \int_M |R_{ij} \: + \: \nabla_i \nabla_j \overline{f}|^2 \: e^{-\overline{f}}
\: dV
\end{equation}
and so
\begin{equation}
\frac{d\lambda}{dt} \: \ge \: \frac{2}{n} \int_M (R + \triangle \overline{f})^2 \: 
e^{- \overline{f}} \: dV.
\end{equation}
From the Cauchy-Schwarz inequality and the fact that 
$\int_M e^{- \overline{f}} \: dV \: = \: 1$, 
\begin{equation}
\left( \int_M (R + \triangle \overline{f}) \: e^{- \overline{f}} \: dV \right)^2 
\: \le \: \int_M (R + \triangle \overline{f})^2 \: 
e^{- \overline{f}} \: dV.
\end{equation}
Finally, (\ref{id}) gives
\begin{equation} \label{using}
\int_M (R + \triangle \overline{f}) \: e^{- \overline{f}} \: dV \: = \:
\int_M (R + | \nabla \overline{f}|^2 ) \: e^{- \overline{f}} \: dV \: = \:
{\mathcal F}(g(t), \overline{f}(t)) \: = \: \lambda(t).
\end{equation}
This proves the lemma.
\end{proof}
\section{I.2.3. The rescaled $\lambda$-invariant}

In this section we show the monotonicity of a scale-invariant
version of $\lambda$.  This will be used in Section \ref{II.8}.
We then show that an expanding breather on a compact manifold
is a gradient expanding soliton.

Put $\overline{\lambda}(g) \: = \: \lambda(g) \: V(g)^{\frac{2}{n}}$.
As $\overline{\lambda}$ is scale-invariant, it is constant in $t$
along a steady, shrinking or expanding soliton solution.

\begin{proposition} \label{expandingprop}  (cf. Claim of I.2.3)
If $g(\cdot)$ is a Ricci flow solution and 
$\overline{\lambda}(g(t)) \: \le \: 0$ for some $t$ then
$\frac{d}{dt} \overline{\lambda}(g(t)) \: \ge \: 0$.
\end{proposition}
\begin{proof}
We have
\begin{equation}
\frac{d\overline{\lambda}}{dt} \: = \:
\frac{d\lambda}{dt} \: V(t)^{\frac{2}{n}} \: - \:
\frac{2}{n} \: V(t)^{\frac{2-n}{n}} \: \lambda(t) \: \int_M R \: dV.
\end{equation}
From (\ref{dlambdadt}),
\begin{equation}
\frac{d\overline{\lambda}}{dt} \: \ge \:
V(t)^{\frac{2}{n}} \: \left[
2 \int_M |R_{ij} \: + \:
\nabla_i \nabla_j \overline{f}|^2 \: e^{- \overline{f}} \:  dV 
 \: - \:
\frac{2}{n} \: V(t)^{-1} \: \lambda(t) \: \int_M R \: dV \right].
\end{equation}
Using the spatially-constant function $\ln V(t)$ as a test function for ${\mathcal F}$
gives
\begin{equation} \label{ttest}
\lambda(t) \: \le \: V(t)^{-1} \: \int_M R \: dV.
\end{equation}
The assumption that $\lambda(t) \: \le \: 0$ gives
\begin{equation}
- \: \lambda(t)^2 \: \le \: - \: V(t)^{-1} \: \lambda(t) \: \int_M R \: dV
\end{equation}
and so
\begin{equation}
\frac{d\overline{\lambda}}{dt} \: \ge \:
V(t)^{\frac{2}{n}} \: \left[
2 \int_M |R_{ij} \: + \:
\nabla_i \nabla_j \overline{f}|^2 \: e^{- \overline{f}} \: dV 
 \: - \:
\frac{2}{n} \: \lambda(t)^2 \right].
\end{equation}
Next,
\begin{equation}
|R_{ij} \: + \:
\nabla_i \nabla_j \overline{f}|^2 \: \ge \:
\left| R_{ij} \: + \:
\nabla_i \nabla_j \overline{f} \: - \: \frac{1}{n} \: (R + \triangle \overline{f}) \: g_{ij}
\right|^2 \: + \: \frac{1}{n} \: (R + \triangle \overline{f})^2.  
\end{equation}
Using (\ref{using}), one obtains
\begin{align} \label{cs}
\frac{d\overline{\lambda}}{dt} \: \ge \: 2 \:
V(t)^{\frac{2}{n}} \: & \left[
\int_M \left| R_{ij} \: + \:
\nabla_i \nabla_j \overline{f} \: - \: \frac{1}{n} \: (R + \triangle \overline{f}) \: 
g_{ij} \right|^2 \: e^{- \overline{f}} \: dV
\: + \: \frac{1}{n} \: \int_M (R + \triangle \overline{f})^2 \: e^{- \overline{f}} \: dV \right. \\
& \left. - \: \frac{1}{n} \: \left( \int_M (R + \triangle \overline{f}) \: e^{- \overline{f}} \: dV
\right)^2 \right]. \notag
\end{align}
As $\int_M e^{- \overline{f}} \: dV \: = \: 1$, the Cauchy-Schwarz inequality implies
that the right-hand side of (\ref{cs}) is nonnegative.
\end{proof}

\begin{corollary} \label{expp}
If $\overline{\lambda}$ is a
constant nonpositive number on an interval $[t_1, t_2]$ then the
Ricci flow solution is a 
gradient soliton.
\end{corollary}
\begin{proof}
From equation (\ref{cs}), we obtain that $R + \triangle \overline{f} \: = \: \alpha(t)$ for some
function $\alpha$ that is spatially constant, and
$R_{ij} \: + \: \nabla_i \nabla_j \overline{f} \: = \: \frac{\alpha(t)}{n} \: g_{ij}$.
Thus $g$ evolves by diffeomorphisms and dilations.
After a shift of the time parameter,
$\alpha(t)$ is proportionate to $t$
cf. \cite[Lemma 2.4]{Chow-Knopf}. 
This is the gradient soliton equation of Appendix \ref{Riccisolitons}.
\end{proof}

\begin{definition}
An {\em expanding breather} is a Ricci flow solution on $[t_1, t_2]$ that
satisfies $g(t_2) \: = \: c \: \phi^* g(t_1)$ for some
$c > 1$ and $\phi \in \Diff(M)$. 
\end{definition}

Expanding soliton solutions are
expanding breathers.

Again, we are assuming that $M$ is
compact.

\begin{proposition}
An expanding breather is a gradient expanding soliton.
\end{proposition}
\begin{proof}
First, $\overline{\lambda}(t_2) \: = \: \overline{\lambda}(t_1)$.
As $V(t_2) \: > \: V(t_1)$, we must have
$\frac{dV}{dt} \: > \: 0$ for some $t \in [t_1, t_2]$.
From (\ref{ttest}), $\frac{dV}{dt} \: = \:
-  \int_M R \: dV \: \le \: - \: \lambda(t) \: V(t)$, so
$\overline{\lambda}(t)$ must be negative for some 
$t \in [t_1, t_2]$.
Proposition \ref{expandingprop} implies that
$\overline{\lambda}(t_1) \: < \: 0$. Then as
$\overline{\lambda}(t_2) \: = \: \overline{\lambda}(t_1)$, 
it follows that $\overline{\lambda}$ is a negative
constant on $[t_1, t_2]$. From Corollary \ref{expp}, the solution is a gradient expanding
soliton.
\end{proof}

\section{I.2.4. Gradient steady solitons on compact manifolds}

It was shown in Section \ref{I.2.2} that a steady breather on a compact manifold
is a gradient steady soliton.  We now show that it is in fact Ricci flat.
This was previously shown in \cite[Theorem 20.1]{Hamilton}.

\begin{proposition}
A gradient steady solution on a compact manifold is Ricci flat.
\end{proposition}
\begin{proof}
As we are in the equality case of Proposition \ref{mono},
the function $f(t)$ must be the minimizer of ${\mathcal F}(g(t), \cdot)$
for all $t$. That is,
\begin{equation}
(- \: 4 \triangle \: + \: R) \: e^{-\frac{f}{2}} \: = \: \lambda \:
e^{-\frac{f}{2}}
\end{equation}
for all $t$, where $\lambda$ is constant in $t$.
Equivalently, $2 \triangle f \: - \: |\nabla f|^2 \: + \: R \: = \:
\lambda$. As $R \: + \: \triangle f \: = \: 0$, we have
$\triangle f \: - \: |\nabla f|^2 \: = \: \lambda$. Then
$\triangle e^{-f} \: = \: - \lambda \: e^{-f}$. Integrating gives
$0 \: = \: \int_M \triangle e^{-f} \: dV \: = \: - \: \lambda \:
\int_M e^{-f} \: dV$, so $\lambda \: = \: 0$. Then
$0 \: = \: - \: \int_M e^{-f} \: \triangle e^{-f} \: dV \: = \:
\int_M \left| \nabla e^{-f} \right|^2 \: dV$, 
so $f$ is constant and $g$ is Ricci flat.  
\end{proof}

A similar argument shows that a gradient expanding soliton on a
compact manifold comes from an Einstein metric with negative
Ricci curvature.

\section{Ricci flow as a gradient flow}

We have shown in Section \ref{I.1.1} that the modified Ricci flow is the gradient flow
for the functional ${\mathcal F}^m$ on the space of metrics ${\mathcal M}$.
One can ask if the unmodified Ricci flow is a gradient flow. This turns out to be true
provided that one
considers it as a flow on the space ${\mathcal M}/\Diff(M)$.

As mentioned in II.8, Ricci flow is the gradient flow
for the function $\lambda$. More precisely, this statement  is valid on ${\mathcal M}/\Diff(M)$,
with the latter being equipped with an appropriate metric.
To see this, we first
consider $\lambda$ as a function on the space of metrics ${\mathcal M}$.
Here the formal Riemannian metric on ${\mathcal M}$ comes from saying
that for $v_{ij} \in T_g {\mathcal M}$,
\begin{equation} \label{riemmetric}
\langle v_{ij}, v_{ij} \rangle \: = \: \frac12 \: 
\int_M v^{ij} \: v_{ij} \: \Phi^2 \: dV(g),
\end{equation}
where $\Phi \: = \: \Phi(g)$ is the unique
normalized positive eigenvector corresponding to $\lambda(g)$.

\begin{lemma}
The formal gradient flow of $\lambda$ 
is
\begin{equation}
\frac{\partial g_{ij}}{\partial t} \: = \:
- \: 2 \: \left( R_{ij} \: - \: 2 \: \nabla_i \nabla_j \ln \Phi 
\right).
\end{equation}
\end{lemma}
\begin{proof}
We set
\begin{equation}
\label{eqnlambdadef}
\lambda(g) \: = \: \inf_{f \in C^\infty(M) 
\: : \: \int_M e^{-f} \: dV \: = \: 1} 
{\mathcal F}(g, f).
\end{equation}
To calculate the variation in $\lambda$
due to a variation $\delta g_{ij} \: = \: v_{ij}$, 
we let $h=\delta f$ be the variation induced by letting $f$ be the 
minimizer in (\ref{eqnlambdadef}). Then 
\begin{equation}
0=\delta \left(\int_M e^{-f}\: dV\right) = \int_M \left(\frac{v}{2}-h\right)\:e^{-f}\:dV.
\end{equation}
Now  equation
(\ref{Fvariation}) gives
\begin{align} \label{lambdavar}
\delta \lambda(v_{ij}) 
\: = \: 
\int_M e^{-f} \: & \left[  - \: v_{ij} 
(R_{ij} \: + \: \nabla_i \nabla_j f) \: + \right. \\
& \left. \left( \frac{v}{2} \: - \: h \right) \: (2 \triangle f \: - \:
|\nabla f|^2 \: + \: R) \right] \: dV. \notag
\end{align}
As $\Phi \: = \: e^{-\frac{f}{2}}$ satisfies
\begin{equation}
- \: 4 \: \triangle \Phi \: + \: R \Phi \: = \: \lambda \Phi,
\end{equation}
it follows that
\begin{equation} \label{feqn}
2 \: \triangle f \: - \: | \nabla f|^2 \: + \: R \: = \: \lambda.
\end{equation}
Hence the last term in (\ref{lambdavar}) vanishes, and 
(\ref{lambdavar}) becomes
\begin{equation}
\delta \lambda(v_{ij}) 
\: = \: 
- \: \int_M e^{-f} \: v_{ij} 
(R_{ij} \: + \: \nabla_i \nabla_j f) \: dV.
\end{equation}
The corresponding gradient flow is
\begin{equation} \label{gradflow}
\frac{\partial g_{ij}}{\partial t} \: = \:
- \: 2 \: \left( R_{ij} \: + \: \nabla_i \nabla_j f
\right)\: = \:
- \: 2 \: \left( R_{ij} \: - \: 2 \: \nabla_i \nabla_j \ln \Phi 
\right). 
\end{equation}
\end{proof}

We note that it follows from (\ref{feqn}) that
\begin{equation}
\nabla_j \left( (R_{ij} \: + \: \nabla_i \nabla_j f) \: e^{-f} \right)
\: = \: 0.
\end{equation}
This implies that the gradient vector field of $\lambda$ is perpendicular
to the infinitesimal diffeomorphisms at $g$, as one would expect.

In the sense of  \cite{Bourguignon}, the quotient space ${\mathcal M}/\Diff(M)$ 
is a stratified
infinite-dimensional Riemannian manifold, with the strata corresponding
to the possible isometry groups $\Isom(M, g)$.
We give it the quotient Riemannian metric coming from
(\ref{riemmetric}). 
The modified Ricci flow (\ref{gradflow}) on 
${\mathcal M}$ projects to a flow on ${\mathcal M}/\Diff(M)$
that coincides with the projection of the unmodified Ricci flow
$\frac{dg}{dt} \: = \: - \: 2 \Ric$. The upshot is that  the
Ricci flow, as a flow on ${\mathcal M}/\Diff(M)$, is the gradient flow of
$\lambda$, the latter now being considered as a function on 
${\mathcal M}/\Diff(M)$.

One sees an intuitive explanation for Proposition \ref{steadybreather}.
If a gradient flow on a finite-dimensional manifold has a periodic
orbit then it must be a fixed-point.  Applying this principle formally
to the Ricci flow on ${\mathcal M}/\Diff(M)$, one infers that a steady
breather only evolves by diffeomorphisms.

\section{The ${\mathcal W}$-functional}
\label{Basic Example for Section 3.1}

\begin{definition}
The {\em ${\mathcal W}$-functional} ${\mathcal W} \: : \: {\mathcal M} \times C^\infty(M) \times
\R^+ \rightarrow \R$ is given by
\begin{equation}
{\mathcal W}(g, f, \tau) \: = \: \int_M \left[ \tau \left(
|\nabla f|^2 + R \right) \: + \: f \: - \: n \right] \:
(4\pi \tau)^{-\frac{n}{2}} \: e^{-f} \: dV.
\end{equation}
\end{definition}

The ${\mathcal W}$-functional
is a scale-invariant variant of ${\mathcal F}$. It
has the symmetries ${\mathcal W}(\phi^* g, \phi^* f, \tau) \: = \:
{\mathcal W}(g, f, \tau)$ for $\phi \in \Diff(M)$, and 
${\mathcal W}(cg, f, c\tau) \: = \: {\mathcal W}(g, f, \tau)$
for $c > 0$. Hence it is constant in $t \: = \: - \: \tau$ along
a gradient shrinking soliton defined for $t \in (-\infty, 0)$,
as in Appendix \ref{Riccisolitons}. In this sense, ${\mathcal W}$ is
constant on gradient shrinking solitons just as ${\mathcal F}$ is
constant on gradient steady solitons.

As an example of a gradient shrinking soliton,
consider $\R^n$ with the flat metric, constant in time
$t \in (-\infty, 0)$. Put $\tau \: = \: - \: t$ and
\begin{equation}
f(t, x) \: = \: \frac{|x|^2}{4\tau},
\end{equation}
so
\begin{equation}
e^{-f} \: = \: e^{- \frac{|x|^2}{4\tau}}.
\end{equation}
One can check that $(g(t), f(t), \tau(t))$ satisfies 
(\ref{neweqn}) and (\ref{Wnorm}).
Now
\begin{equation}
\tau ( |\nabla f|^2 \: + \: R) \: + \: f \: - \: n \: = \:
\tau \cdot \frac{|x|^2}{4 \tau^2} \: + \: \frac{|x|^2}{4\tau} \: - \: n
\: = \:  \frac{|x|^2}{2\tau} \: - \: n.
\end{equation}
It follows from (\ref{Gaussianintegral}) and (\ref{Gaussianintegral2}) 
that ${\mathcal W}(t) \: = \: 0$ for all $t$.

\section{I.3.1. Monotonicity of the ${\mathcal W}$-functional}

In this section
we compute the variation of ${\mathcal W}$, in analogy with
the computation in Section \ref{I.1.1} of the variation of ${\mathcal F}$.
We then show that a shrinking breather on a compact manifold
is a gradient shrinking soliton.

As in Section \ref{I.1.1}, we write $\delta g_{ij} \: = \: v_{ij}$ and $\delta f \: = \: h$.
Put $\sigma \: = \: \delta \tau$. 

\begin{proposition}
We have
\begin{align}
\delta {\mathcal W}(v_{ij}, h, \sigma) \: = \:
\int_M & \left[ \sigma  (R \: + \: |\nabla f|^2) \:
- \: \tau v_{ij} (R_{ij} \: + \: \nabla_i \nabla_j f) \: + \: h \: +
\right.  \\
& \left.
\left[ \tau (2 \triangle f \: - \: |\nabla f|^2 \: + \: R) \: + \: f
\: - \: n \right] \:  
\left( \frac{v}{2} \: - \: h \: - \: \frac{n\sigma}{2\tau} \right) \right] \:
(4 \pi \tau)^{-n/2} \: e^{-f} \: dV. \notag
\end{align}
\end{proposition}
\begin{proof}
One finds
\begin{equation}
\delta \left( (4 \pi \tau)^{-n/2} \: e^{-f} \: dV \right) \: = \:
\left( \frac{v}{2} \: - \: h \: - \: \frac{n\sigma}{2\tau} \right) \:
(4 \pi \tau)^{-n/2} \: e^{-f} \: dV.
\end{equation}
Writing
\begin{equation}
{\mathcal W} \: = \: \int_M \left[ \tau (R \: + \: |\nabla f|^2) \: + \: f
\: - \: n \right] \: (4 \pi \tau)^{-n/2} \: e^{-f} \: dV
\end{equation}
we can use 
(\ref{Fvariation}) to obtain
\begin{align}
\delta {\mathcal W} \: = \:
\int_M & \left[ \sigma  (R \: + \: |\nabla f|^2) \: + \:
\tau  \left( \frac{v}{2} \:  - \: h \: \right) (2 \triangle f \: - \:
2|\nabla f|^2) \: - \: \tau v_{ij} (R_{ij} \: + \: \nabla_i \nabla_j f) \: +
\right. \\
& \left. h \: + \:
\left[ \tau (R \: + \: |\nabla f|^2) \: + \: f
\: - \: n \right] \:  
\left( \frac{v}{2} \: - \: h \: - \: \frac{n\sigma}{2\tau} \right) \right] \:
(4 \pi \tau)^{-n/2} \: e^{-f} \: dV.
\end{align}
Then (\ref{id}) gives
\begin{align}
\delta {\mathcal W} \:
= \:
\int_M & \left[ \sigma  (R \: + \: |\nabla f|^2) \:
- \: \tau v_{ij} (R_{ij} \: + \: \nabla_i \nabla_j f) \: + \: h \: +
\right. \\
& \left.
\left[ \tau (2 \triangle f \: - \: |\nabla f|^2 \: + \: R) \: + \: f
\: - \: n \right] \:  
\left( \frac{v}{2} \: - \: h \: - \: \frac{n\sigma}{2\tau} \right) \right] \:
(4 \pi \tau)^{-n/2} \: e^{-f} \: dV. \notag
\end{align}
This proves the proposition.
\end{proof}

We now fix a smooth measure $dm$ on $M$ with mass $1$ and
relate $f$ to $g$ and $\tau$ by
requiring that $(4 \pi \tau)^{-n/2} \: e^{-f} \: dV \: =\: dm$. Then
$\frac{v}{2} \: - \: h \: - \: \frac{n\sigma}{2\tau} \: = \: 0$ and
\begin{equation} \label{Wvar}
\delta {\mathcal W} \: = \: \int_M  \left[ \sigma  (R \: + \: |\nabla f|^2) \:
- \: \tau v_{ij} (R_{ij} \: + \: \nabla_i \nabla_j f) \: + \: h \right] \:
(4 \pi \tau)^{-n/2} \: e^{-f} \: dV.
\end{equation}

We now consider $\frac{d{\mathcal W}}{dt}$ when
\begin{align} \label{oldeqn}
(g_{ij})_t \: & = \: - \: 2 \: (R_{ij} \: + \: \nabla_i \nabla_j f ), \\ 
f_t \: & = \: - \: \triangle f \: - \: R \: + \: \frac{n}{2\tau}, \notag \\
\tau_t \: & = \: -1. \notag
\end{align}
To apply (\ref{Wvar}), we put
\begin{align}
v_{ij} \: & = \: - \: 2 \: (R_{ij} \: + \: \nabla_i \nabla_j f ), \\ 
h \: & = \: - \: \triangle f \: - \: R \: + \: \frac{n}{2\tau}, \notag \\
\sigma \: & = \: -1. \notag
\end{align}
We do have
$\frac{v}{2} \: - \: h \: - \: \frac{n\sigma}{2\tau} \: = \: 0$.
Then from (\ref{Wvar}),
\begin{align} \label{Wtimeder}
\frac{d{\mathcal W}}{dt} \: = \:
\int_M & \left[- \:  (R \: + \: |\nabla f|^2) 
\: + \: 2 \: \tau |R_{ij} \: + \: \nabla_i \nabla_j f|^2 \: 
\right. \\
& -\left. \triangle f \: - \: R \: + \: \frac{n}{2\tau}
 \right] \:
(4 \pi \tau)^{-n/2} \: e^{-f} \: dV \notag \\
= \:
\int_M & \left[- \:  2(R \: + \: \triangle f) 
\: + \: 2 \: \tau |R_{ij} \: + \: \nabla_i \nabla_j f|^2 
\: + \: \frac{n}{2\tau}
 \right] \:
(4 \pi \tau)^{-n/2} \: e^{-f} \: dV \notag \\
= \: \int_M &
2 \: \tau |R_{ij} \: + \: \nabla_i \nabla_j f
\: - \: \frac{1}{2\tau} \: g_{ij}|^2 \:
(4 \pi \tau)^{-n/2} \: e^{-f} \: dV. \notag
\end{align}
Adding a Lie derivative to the right-hand side of
(\ref{oldeqn}) gives the new flow equations
\begin{align} \label{neweqn}
(g_{ij})_t \: & = \: - \: 2 \: R_{ij}, \\ 
f_t \: & = \: - \: \triangle f \: + \: |\nabla f|^2
\: - \: R \: + \: \frac{n}{2\tau}, \notag \\
\tau_t \: & = \: -1, \notag
\end{align}
with (\ref{Wtimeder}) still holding. We no longer have
$(4 \pi \tau)^{-n/2} \: e^{-f} \: dV \: = \: dm$, but we do have
\begin{equation} \label{Wnorm}
\int_M (4 \pi \tau)^{-n/2} \: e^{-f} \: dV \: = \: 1.
\end{equation}

We now want to look at the variational problem of 
minimizing ${\mathcal W}(g,f,\tau)$ under the constraint
that $\int_M (4 \pi \tau)^{-n/2} \: e^{-f} \: dV \: = \: 1$.
We write
\begin{equation}
\mu(g, \tau) \: = \: \inf_f \{ {\mathcal W}(g,f,\tau) \: : \: 
\int_M (4 \pi \tau)^{-n/2} \: e^{-f} \: dV \: = \: 1\}.
\end{equation}
Making the change of variable $\Phi \: = \: e^{- \frac{f}{2}}$, we are
minimizing
\begin{equation}
(4 \pi \tau)^{-n/2} \: \int_M \left[ \tau (4 |\nabla \Phi|^2 \: + \: R \Phi^2) \: - \:
2 \Phi^2 \log \Phi \: - \: n \: \Phi^2 \right] \: dV
\end{equation}
under the constraint $(4 \pi \tau)^{-n/2} \: \int_M \Phi^2 \: dV \: = \: 1$.
From  \cite[Section 1]{Rothaus (1981)} the infimum is finite and
there is a positive
continuous minimizer $\Phi$. It will be a weak solution of the
variational equation
\begin{equation}
\tau (- 4 \triangle + R) \Phi \: = \: 2\Phi \log \Phi \: + \: 
(\mu(g, \tau) + n) \Phi.
\end{equation}
From elliptic theory, $\Phi$ is smooth. Then $f \: = \: - \: 2 \log \Phi$
is also smooth.

As in Section \ref{I.2.2}, it follows that
$\mu(g(t), t_0 - t)$ is nondecreasing in $t$ for a Ricci flow
solution, where $t_0$ is any fixed number and $t < t_0$. 
If it is constant in $t$ then the solution must be
a gradient shrinking soliton that goes singular at time $t_0$.

\begin{definition}
A {\em shrinking breather} is a Ricci flow solution on $[t_1, t_2]$ that
satisfies $g(t_2) \: = \: c \: \phi^* g(t_1)$ for some
$c < 1$ and $\phi \in \Diff(M)$. 
\end{definition}

Gradient shrinking soliton solutions
are shrinking breathers.

Again, we are assuming that $M$ is
compact.

\begin{proposition}
A shrinking breather is a gradient shrinking soliton.
\end{proposition}
\begin{proof}
Put $t_0 \: = \: \frac{t_2 - c t_1}{1-c}$. Then if
$\tau_1 \: = \: t_0 - t_1$ and $\tau_2 \: = \: t_0 - t_2$, we
have $\tau_2 \: = \: c\tau_1$. Hence
\begin{equation}
\mu(g(t_2), \tau_2) \: = \: \mu \left( \frac{\tau_2}{\tau_1} \:
\phi^* g(t_1), \tau_2 \right) \: = \: 
\mu \left(\phi^* g(t_1), \tau_1 \right)
 \: = \: 
\mu \left(g(t_1), \tau_1 \right).
\end{equation}
It follows that the solution is a gradient shrinking soliton.
\end{proof}

\section{I.4. The no local collapsing theorem I}

In this section we prove the no local collapsing theorem.

\begin{definition}
A smooth Ricci flow solution $g(\cdot)$ on a time interval
$[0, T)$ is said to be {\em locally collapsing} at $T$ if there
is a sequence of times $t_k \rightarrow T$ and a sequence 
of metric balls $B_k = B(p_k, r_k)$ at times $t_k$ such that
$r_k^2/t_k$ is bounded, $|\Rm|(g(t_k)) \: \le \: r_k^{-2}$
in $B_k$ and $\lim_{k \rightarrow \infty}
r_k^{-n} \: \vol(B_k) \: = \: 0$.
\end{definition}

\begin{remark}
In the definition of noncollapsing, $T$ could be infinite.
This is why it is written that $r_k^2/t_k$ stays bounded, while
if $T \: < \: \infty$ then this is obviously the same as saying that
$r_k$ stays bounded.
\end{remark}

\begin{theorem} \label{noncollthm} (cf. Theorem I.4.1)
If $M$ is closed and $T < \infty$ then $g(\cdot)$ is not locally
collapsing at $T$.
\end{theorem}
\begin{proof}
We first sketch the idea of the proof.
In Section \ref{Basic Example for Section 3.1} we showed that in the
case of flat $\R^n$, taking $e^{-f}(x) \: = \:
e^{- \frac{|x|^2}{4\tau}}$, we get ${\mathcal W}(g, f, \tau) \: = \: 0$.
So putting $\tau \: = \: r_k^2$ and $e^{-f_k}(x) \: = \:
e^{- \frac{|x|^2}{4r_k^2}}$, we have ${\mathcal W}(g, f_k, r_k^2) \: = \: 0$.
In the collapsing case, the idea is to use a test function $f_k$ so that
\begin{equation}
e^{-f_k}(x) \: \sim \:
e^{- c_k} \: e^{- \frac{\dist_{t_k}(x, p_k)^2}{4r_k^2}},
\end{equation}
where $c_k$ is determined by the normalization condition
\begin{equation} \label{normalization}
\int_M (4 \pi r_k^2)^{- n/2} \: e^{-f_k} \: dV \: = \: 1.
\end{equation} 
The main difference between computing (\ref{normalization}) in $M$ and
in $\R^n$ comes from the difference in volumes, which means that
$e^{-c_k} \sim \frac{1}{r_k^{-n} \: vol(B_k)}$. 
In particular, as $k \rightarrow
\infty$, we have $c_k \rightarrow - \infty$. 

Now that $f_k$ is normalized correctly, the main difference between
computing ${\mathcal W}(g(t_k), f_k, r_k^2)$ in $M$, and the
analogous computation for the
Gaussian in $\R^n$,
comes from the $f$ term in the integrand of ${\mathcal W}$.  Since
$f_k \sim c_k$, this will drive ${\mathcal W}(g(t_k), f_k, r_k^2)$ to
$- \infty$ as $k \rightarrow \infty$, so $\mu(g(t_k), r_k^2)
\rightarrow - \infty$; by the monotonicity of $\mu(g(t),t_0-t)$
it follows that $\mu(g(0),t_k+r_k^2)\ra -\infty$ as $k\ra \infty$.
This contradicts the fact that $\mu(g(0),\tau)$ is a continuous function 
of $\tau$.

To write this out precisely, let us put $\Phi \: = \: e^{- f/2}$, so
that 
\begin{equation}
{\mathcal W}(g, f, \tau) \: = \: (4 \pi \tau)^{-n/2} \:
\int_M \left[ 4 \tau \: |\nabla \Phi|^2 \: + \:
( \tau R \: - \: 2 \ln \Phi
\: - \: n ) \: \Phi^2 \right] \: dV.
\end{equation}
For the argument, it is enough to obtain small values of ${\mathcal W}$ for 
positive $\Phi$.   Since $\lim_{s\ra 0}(-2\ln s)\:s^2\:=\:0$, 
by an approximation it is enough to obtain small values of ${\mathcal W}$ for 
nonnegative $\Phi$, where the integrand is declared to be 
$4 \tau \: |\nabla \Phi|^2$ at points where $\Phi$ vanishes. Take 
\begin{equation}
\Phi_k(x) \: = \: e^{- c_k/2} \: \phi(\dist_{t_k}(x, p_k)/r_k),
\end{equation}
where $\phi \: : \: [0, \infty) \rightarrow [0, 1]$ is 
a monotonically nonincreasing function such that
$\phi(s) \: = \: 1$ if $s \in [0, 1/2]$, $\phi(s) \: = \: 0$ if $s \: \ge \: 1$
and $|\phi^\prime(s)| \: \le \: 10$ for $s \in [1/2, 1]$.
The function $\Phi_k$ is {\it a priori} only Lipschitz, but by smoothing
it slightly we can use $\Phi_k$ in the variational formula 
to bound ${\mathcal W}$ from above.

The constant $c_k$ is determined by
\begin{equation}
e^{c_k} \: = \: \int_M (4 \pi r_k^2)^{-n/2} \: 
\phi^2(\dist_{t_k}(x, p_k)/r_k) \: dV \: \le \:
(4 \pi r_k^2)^{-n/2} \: 
\vol(B_k).
\end{equation}
Thus $c_k \rightarrow - \infty$.
Next,
\begin{equation} \label{mess}
{\mathcal W}(g(t_k), f_k, r_k^2) \: = \: (4 \pi r_k^2)^{-n/2} \:
\int_M \left[ 4 r_k^2 \: |\nabla \Phi_k|^2 \: + \:
( r_k^2 R \: - \: 2 \ln \Phi_k
\: - \: n ) \: \Phi_k^2 \right] \: dV.
\end{equation}
Let $A_k(s)$ be the mass of the distance sphere $S(p_k, r_k s)$ around $p_k$.
Put 
\begin{equation}
\overline{R}_k(s) \: = \: r_k^{2} \: A_k(s)^{-1} 
\: \int_{S(p_k, r_k s)} R \: d\area. 
\end{equation}
We can compute the integral in (\ref{mess}) radially to get
\begin{equation}
{\mathcal W}(g(t_k), f_k, r_k^2) \: 
= \:  \frac{\int_0^1 \left[ 4 (\phi^\prime(s))^2 \: + \:
(\overline{R}_k(s) \: + \: c_k \: - \: 2 \: \ln \phi(s) \: 
\: - \: n ) \: \phi^2(s) \right] \: A_k(s) \: 
ds}{\int_0^1 \phi^2(s) \: A_k(s) \: ds}.
\end{equation}
The expression 
$4 (\phi^\prime(s))^2 \: - \: 2 \: \ln \phi(s) \: \phi^2(s)$ vanishes
if $s \notin [1/2, 1]$, and is bounded above by $400 \: + \: e^{-1}$ if
$s \in [1/2, 1]$. Then the lower
bound on the Ricci curvature and the Bishop-Gromov
inequality give
\begin{align}
\frac{\int_0^1 \left[ 4 (\phi^\prime(s))^2 \: - \: 2 \: \ln \phi(s) 
\: \phi^2(s) \right] \: A_k(s) \: 
ds}{\int_0^1 \phi^2(s) \: A_k(s) \: ds} \: & \le \: 401 \: 
\frac{\vol(B(p_k, r_k)) - \vol(B(p_k, r_k/2))}{\vol(B(p_k, r_k/2))} \\
& \le \: 401 \: \left( 
\frac{\int_0^1 \sinh^{n-1}(s) 
\: ds}{\int_0^{1/2} \sinh^{n-1}(s) \: ds} \: - \: 1 \right).
\notag
\end{align}
Next, from the upper bound on scalar curvature,
$\overline{R}_k(s) \: \le \: n(n-1)$ for $s \in [0,1]$.
Putting this together gives
${\mathcal W}(g(t_k), f_k, r_k^2) \: \le \: \const \: + \: c_k$ and 
so ${\mathcal W}(g(t_k), f_k, r_k^2) \rightarrow -
\infty$ as $k \rightarrow \infty$.

Thus $\mu(g(t_k), r_k^2) \rightarrow - \infty$. 
For any $t_0>t$, $\mu(g(t), t_0 - t)$ is
nondecreasing in $t$. Hence $\mu(g(0), t_k \: + \: r_k^2) \: \le \:
\mu(g(t_k), r_k^2)$, so $\mu(g(0), t_k \: + \: r_k^2) \rightarrow
- \infty$. Since $T$ is finite, $t_k$ and $r_k^2$ are 
uniformly bounded, 
and $t_k$ uniformly positive,
which
contradicts the fact that $\mu(g(0),\tau)$ is a continuous function of $\tau$.
\end{proof}

\begin{remark}
In the preceding argument we only used the upper bound on scalar
curvature and the lower bound on Ricci
curvature, i.e. in the definition
of local collapsing one could have assumed that
$R(g_{ij}(t_k)) \: \le \: n(n-1) \: r_k^{-2}$ in $B_k$ and
$\Ric(g_{ij}(t_k)) \: \ge \: - \: (n-1) \: r_k^{-2}$ in $B_k$.
In fact, one can also remove the lower bound on Ricci curvature
(observation of Perelman, communicated by Gang Tian).
The necessary ingredients of the preceding argument were that \\
1. $r_k^{-n} \: \vol(B(p_k, r_k)) \rightarrow 0$, \\
2. $r_k^2 \: R$ is uniformly bounded above on $B(p_k, r_k)$ and \\
3. $\frac{\vol(B(p_k, r_k))}{\vol(B(p_k, r_k/2))}$ is uniformly bounded above.

Suppose only that $r_k^{-n} \: \vol(B(p_k, r_k)) \rightarrow 0$ and 
for all $k$, $r_k^2 \: R \: \le \: n(n-1)$ on $B(p_k, r_k)$. If
$\frac{\vol(B(p_k, r_k))}{\vol(B(p_k, r_k/2))} \: < \: 3^n$ for all
$k$ then we are done.  If not, suppose that for a given $k$,
$\frac{\vol(B(p_k, r_k))}{\vol(B(p_k, r_k/2))} \: \ge \: 3^n$.
Putting $r_k^\prime \: = \: r_k/2$, we have that
$(r_k^\prime)^{-n} \: \vol(B(p_k, r_k^\prime)) \: \le \:
r_k^{-n} \: \vol(B(p_k, r_k))$ and 
$(r_k^\prime)^2 \: R \: \le \: n(n-1)$ on $B(p_k, r_k^\prime)$.
We replace $r_k$ by $r_k^\prime$. If now 
$\frac{\vol(B(p_k, r_k))}{\vol(B(p_k, r_k/2))} \: < \: 3^n$ then we
stop. If not then we repeat the process and replace $r_k$ by $r_k/2$.
Eventually we will achieve that 
$\frac{\vol(B(p_k, r_k))}{\vol(B(p_k, r_k/2))} \: < \: 3^n$.
Then we can apply the preceding argument to this new sequence of pairs
$\{(p_k, r_k)\}_{k=1}^\infty$.
\end{remark}

\begin{definition} \label{noncollapsed1} (cf. Definition I.4.2) 
We say that a metric $g$ is {\em $\kappa$-noncollapsed on the scale
$\rho$} if every metric ball $B$ of radius $r < \rho$, which satisfies
$|\Rm(x)| \: \le \: r^{-2}$ for every $x \in B$, has volume at least
$\kappa r^n$.
\end{definition}

\begin{remark}
We caution the reader that this definition differs slightly from the 
definition of noncollapsing that is used from section I.7 onwards.
\end{remark}

Note that except for the overall scale $\rho$, the $\kappa$-noncollapsed
condition is scale-invariant.
From the proof of Theorem \ref{noncollthm} we extract the following
statement. Given a Ricci flow defined on
an interval $[0, T)$, with $T < \infty$, and a scale $\rho$, there
is some number $\kappa = \kappa(g(0), T, \rho)$ so that the solution is
$\kappa$-noncollapsed on the scale $\rho$ for all $t \in [0, T)$.
We note that the estimate on $\kappa$ deteriorates as $T \rightarrow
\infty$, as there are Ricci flow solutions that collapse at 
long time.

\section{I.5. The ${\mathcal W}$-functional as a time derivative} \label{I.5}

We will only discuss one formula from I.5, showing that along a Ricci flow,
${\mathcal W}$ is itself the time-derivative of an integral expression.

Again, we put $\tau \: = \: - \: t$. Consider the evolution equations
(\ref{oldeqn}), with
$(4 \pi \tau)^{-n/2} \: e^{-f} \: dV \: =\: dm$. Then
\begin{align}
\frac{d}{d\tau} \left( \tau \int_M \left( f \: - \: \frac{n}{2} \right)
\: dm \right) \: & = \:
\int_M \left( f \: - \: \frac{n}{2} \right)
\: dm \: + \: \tau \: \int_M  \left( \triangle f \: + \: R \: - \:
\frac{n}{2\tau} \right) \: dm \\
& = \:
\int_M \left( f \: - \: \frac{n}{2} \right)
\: dm \: + \: \tau \: \int_M  \left( |\nabla f|^2 \: + \: R \: - \:
\frac{n}{2\tau} \right) \: dm \notag \\
& = \: {\mathcal W}(g(t), f(t), \tau). \notag
\end{align}
With respect to the evolution equations
(\ref{neweqn}) obtained by performing diffeomorphisms, we get
\begin{equation}
\frac{d}{d\tau} \left( \tau \int_M \left( f \: - \: \frac{n}{2} \right)
\: (4 \pi \tau)^{-n/2} \: e^{-f} \: dV \right) \: 
 = \: {\mathcal W}(g(t), f(t), \tau). 
\end{equation}

Similarly, with respect to (\ref{newevolution}),
\begin{equation}
\frac{d}{dt} \left(- \int_M f
\: e^{-f} \: dV \right) \: 
 = \: {\mathcal F}(g(t), f(t)). 
\end{equation}

\section{I.7. Overview of reduced length and reduced volume} \label{I.7}

We first give a brief summary of I.7.
In I.7, the variable $\tau \: = \: t_0 \: - t$ is used and so the
corresponding Ricci flow equation is $(g_{ij})_\tau \: = \: 2 R_{ij}$.
The goal is to prove a no local collapsing theorem
by means of the ${\mathcal L}$-lengths of curves
$\gamma \: : \: [\tau_1, \tau_2] \rightarrow M$, defined by
\begin{equation} \label{defL}
{\mathcal L}(\gamma) \: = \:
\int_{\tau_1}^{\tau_2} \sqrt{\tau} \: \left( R(\gamma(\tau)) \: + \:
\big| \dot{\gamma}(\tau) \big|^2 \right) \: d\tau,
\end{equation}
where the scalar curvature $R(\ga(\tau))$ and the norm $|\dot{\gamma}(\tau)|$
are evaluated using the metric at time $t_0-\tau$.
Here $\tau_1 \ge 0$.
With $X \: = \: \frac{d\gamma}{d\tau}$, 
the corresponding ${\mathcal L}$-geodesic equation is
\begin{equation}
\nabla_X X \: - \: \frac12 \: \nabla R \: + \: \frac{1}{2\tau} X \: + \:
2 \Ric(X, \cdot) \: = \: 0,
\end{equation}
where again the connection and curvature are taken at the corresponding time,
and the $1$-form $\Ric(X,\cdot)$ has been identified with the
corresponding dual vector field.

Fix $p \in M$.
Taking $\tau_1 \: = \: 0$ and $\gamma(0) \: = \: p$,
the vector $v \: = \: \lim_{\tau \rightarrow 0} 
\sqrt{\tau} \: X(\tau)$ is well-defined in $T_pM$
and is called the initial vector of the geodesic. The 
${\mathcal L}$-exponential map ${\mathcal L}exp_\tau \: : \: T_p M
\rightarrow M$ sends $v$ to $\gamma(\tau)$.

The function $L(q, \overline{\tau})$ is the infimal ${\mathcal L}$-length of
curves $\gamma$ with $\gamma(0) \: = \: p$ and $\gamma(\overline{\tau}) 
\: = \: q$. Defining the reduced length by 
\begin{equation} \label{reducedl}
l(q, {\tau}) \: = \: \frac{L(q, \tau)}{2 \sqrt{\tau}}
\end{equation}
and the reduced volume by
\begin{equation} \label{reducedv}
\tilde{V}(\tau) \: = \: \int_M \tau^{-\frac{n}{2}} \: e^{-l(q, \tau)} \:
dq, 
\end{equation}
the goal is to show that $\tilde{V}(\tau)$ is
nonincreasing in $\tau$, i.e. nondecreasing in $t$.
To do this, one uses the ${\mathcal L}$-exponential map to
write $\tilde{V}(\tau)$ as an integral over $T_pM$:
\begin{equation}
\tilde{V}(\tau) \: = \: \int_{T_pM} \tau^{-\frac{n}{2}} \: 
e^{-l({\mathcal L}exp_\tau(v), \tau)} \: {\mathcal J}(v, \tau) \: 
\chi_{\tau}(v) \:
dv,
\end{equation} 
where ${\mathcal J}(v, \tau) \: = \:
\det \: d\left( {\mathcal L}exp_\tau \right)_v$ is the Jacobian factor
in the change of variable and 
$\chi_\tau$ is a cutoff function related to the
${\mathcal L}$-cut locus of $p$.
To show that $\tilde{V}(\tau)$ is
nonincreasing in $\tau$ it suffices to show that
$\tau^{-\frac{n}{2}} \: 
e^{-l({\mathcal L}exp_\tau(v), \tau)} 
\: {\mathcal J}(v, \tau)$ is nonincreasing
in $\tau$ or, equivalently, that $- \: \frac{n}{2} \: \ln(\tau) \:
- \: l({\mathcal L}exp_\tau(v), \tau) \: + \: \ln \:{\mathcal J}(v, \tau)$ is
nonincreasing in $\tau$. Hence it is necessary to compute
$\frac{dl({\mathcal L}exp_\tau(v), \tau)}{d\tau}$  and
$\frac{d{\mathcal J}(v, \tau)}{d\tau}$. The
computation of the latter will involve the ${\mathcal L}$-Jacobi fields.

The fact that $\tilde{V}(\tau)$ is
nonincreasing in $\tau$ is then used to show that the Ricci flow
solution cannot be  collapsed near $p$.

\section{Basic example for I.7} \label{7Euclidean}
 
 In this section we say what the various expressions of I.7 become in
 the model case of a flat Euclidean Ricci solution.
 
If $M$ is flat $\R^n$ and $p \: = \: \vec{0}$ then the unique
${\mathcal L}$-geodesic $\gamma$ with $\gamma(0) \: =\: \vec{0}$ and
$\gamma(\overline{\tau}) \: = \: \vec{q}$ is 
\begin{equation}
\gamma(\tau) \: = \: \left( \frac{\tau}{\overline{\tau}}
\right)^{\frac12} \: \vec{q} \: = \: 2 \: \tau^{\frac12} \: \vec{v}.
\end{equation}
The function $L$ is given by
\begin{equation}
L(q, \overline{\tau}) \: = \: \frac12 \: \overline{\tau}^{\: -
\frac12} \: |q|^2
\end{equation}
and the reduced length (\ref{reducedl}) is given by
\begin{equation}
l(q, \overline{\tau}) \: = \: \frac{|q|^2}{4 \overline{\tau}}.
\end{equation}
The function $\overline{L}(q, \overline{\tau}) \: = \:
2 \: \overline{\tau}^{\frac12} \: L(q, \overline{\tau})$
is
\begin{equation}
\overline{L}(q, \overline{\tau}) \: = \: |q|^2.
\end{equation}
Then
\begin{equation}
\tilde{V}(\tau) \: = \: 
\int_{\R^n} \tau^{-\frac{n}{2}} \: e^{-\frac{|q|^2}{4 {\tau}}} \:
d^nq \: = \: (4\pi)^{\frac{n}{2}} 
\end{equation}
is constant in $\tau$.

\section{Remarks about ${\mathcal L}$-Geodesics and 
${\mathcal L}\exp$} \label{Lremarks}

In this section we discuss the variational equation corresponding to
(\ref{defL}).

To derive the ${\mathcal L}$-geodesic equation, as in Riemannian geometry
we consider a
$1$-parameter family of curves 
$\gamma_s \: : \: [\tau_1, \tau_2] \rightarrow
M$, parametrized by $s \in (-\epsilon, \epsilon)$. Equivalently, we
have a map $\tilde\gamma(s, \tau)$ with $s \in (-\epsilon, \epsilon)$ and
$t \in [\tau_1, \tau_2]$. Putting $X \: = \: 
\frac{\partial \tilde\gamma}{\partial \tau}$ and 
$Y \: = \: 
\frac{\partial \tilde\gamma}{\partial s}$, we have
$[X, Y] \: = \: 0$. Then $\nabla_X Y \: = \: \nabla_Y X$.
Restricting to the curve $\gamma(\tau) \: = \: \tilde\gamma(0, \tau)$ and
writing $\delta_Y$ as shorthand for $\frac{d}{ds} \big|_{s=0}$, we
have $(\delta_Y \gamma)(\tau) \:  =\: Y(\tau)$ and $(\delta_Y X)(\tau) \: = \: 
(\nabla_X Y)(\tau)$.
Then
\begin{equation}
\delta_Y {\mathcal L} \: = \:
\int_{\tau_1}^{\tau_2} \sqrt{\tau} \: \left( \langle Y, \nabla R \rangle
\: + \: 2 \: \langle \nabla_X Y, X \rangle \right) \: d\tau.
\end{equation}
Using the fact that $\frac{dg_{ij}}{d\tau} \: = \: 2 R_{ij}$, we have
\begin{equation}
\frac{d\langle Y, X \rangle}{d\tau} \: = \:
\langle \nabla_X Y, X \rangle \: + \: 
\langle Y, \nabla_X X \rangle \: + \: 
2 \: \Ric(Y, X).
\end{equation}
Then
\begin{align}
\label{firstvariation}
& \int_{\tau_1}^{\tau_2} \sqrt{\tau} \: \left( \langle Y, \nabla R \rangle
\: + \: 2 \: \langle \nabla_X Y, X \rangle \right) \: d\tau \: = \\
& \int_{\tau_1}^{\tau_2} \sqrt{\tau} \: \left( \langle Y, \nabla R \rangle
\: + \: 2 \: \frac{d}{d\tau} \langle Y, X \rangle \: - \: 2
\langle Y, \nabla_X X \rangle \: - \: 4 \: \Ric(Y, X) 
 \right) \: d\tau \: = \notag \\
& 2 \: \sqrt{{\tau}} \: \langle X, Y \rangle \big|_{\tau_1}^{\tau_2}
\: + \: 
\int_{\tau_1}^{\tau_2} \sqrt{\tau} \: \left\langle Y, \nabla R
\: - \: 2 \nabla_X X  \: - \: 4 \: \Ric(X, \cdot) \: -\: \frac{1}{\tau}
\: X \right\rangle 
 \: d\tau. \notag
\end{align}
Hence the ${\mathcal L}$-geodesic equation is
\begin{equation} \label{Lgeod}
\nabla_X X \: - \: \frac12 \: \nabla R \: + \: \frac{1}{2\tau} X \: + \:
2 \Ric(X, \cdot) \: = \: 0.
\end{equation}

We now discuss some technical issues about ${\mathcal L}$-geodesics and
the ${\mathcal L}$-exponential map.
We are assuming that $(M,g(\cdot))$ is a Ricci flow, where the curvature 
operator
of $M$ is uniformly bounded on a $\tau$-interval $[\tau_1,\tau_2]$, and
each $\tau$-slice $(M,g(\tau))$ is complete for $\tau\in [\tau_1,\tau_2]$.
By Appendix \ref{applocalder},  for every $\tau'<\tau_2$ there 
is a constant $D< \infty$ such that
\begin{equation}
\label{dr}
|\nabla R(x,\tau)|< \frac{D}{\sqrt{\tau_2 - \tau}}
\end{equation}
for all $x\in M$, $\tau \in [\tau_1, \tau_2)$.

Making the change of variable $s \: = \: \sqrt{\tau}$ in the formula
for ${\mathcal L}$-length, we get
\begin{equation}
\label{changetos}
{\mathcal L}(\gamma) \: = \: 
2 \: \int_{s_1}^{s_2} \: \left( \frac14 
\Big| \frac{d\gamma}{ds} \Big|^2 \: + \:
s^2 \: R(\gamma(s)) \right) \: ds.
\end{equation}
The Euler-Lagrange equation becomes
\begin{equation}
\label{sE-L}
\nabla_{\hat X}\hat X-2s^2\nabla R+4s\Ric(\hat X,\cdot)=0,
\end{equation}
where $\hat X \: = \: \frac{d\gamma}{ds}=2sX$.
Putting $s_1 \: = \: \sqrt{\tau_1}$,
it follows from standard existence theory
for ODE's that for each $p \in M$ and $v\in T_pM$,
there is a unique solution $\gamma(s)$ to 
(\ref{sE-L}), defined on an interval $[s_1,s_1+\eps)$, with $\gamma(s_1) = p$
and
\begin{equation}
\frac{1}{2} \: \gamma^\prime(s_1) \: = \: 
\lim_{\tau\ra \tau_1}\sqrt{\tau}\frac{d\gamma}{d\tau} \: = \: v.
\end{equation}
If $\gamma(s)$ is defined for $s \in [s_1,s']$ then 
\begin{align}
\frac{d}{ds}|\hat X|^2=\frac{d}{ds} \: \langle \hat X,\hat X\rangle
&= 4s\Ric(\hat X,\hat X)+2\langle\nabla_{\hat X}\hat X,\hat X\rangle \\
&= - 4s\Ric(\hat X,\hat X) + 4s^2 \langle \nabla R, \hat X \rangle \notag
\end{align}
and so if $\hat X (s) \neq 0$ then
\begin{equation}
\frac{d}{ds}|\hat X|=
\frac{1}{2|\hat X|} \: \frac{d}{ds}|\hat X|^2 \: = \:
- 2s |\hat X|
 \Ric \left( \frac{\hat X}{|\hat X|}, \frac{\hat X}{|\hat X|} \right)
+ 2s^2 \left\langle \nabla R, \frac{\hat X}{|\hat X|} \right\rangle.
\end{equation}
By (\ref{dr}),
\begin{equation}
\frac{d}{ds}|\hat X| \: \le \: C_1 |\hat X| + 
\frac{C_2}{\sqrt{s_2 - s}}
\end{equation}
for appropriate constants $C_1$ and $C_2$, where
$s_2 \: = \: \sqrt{\tau_2}$.
Since the metrics $g(\tau)$
are uniformly comparable for $\tau\in [\tau_1,\tau_2]$, we conclude
(by a continuity argument 
in $s$
) that the ${\mathcal L}$-geodesic $\gamma_v$
with $\frac{1}{2} \gamma_v^\prime(s_1)=v$ is defined on the whole 
interval 
$[s_1,s_2]$.
In particular, in terms of the original variable $\tau$,
for each 
${\tau}\in [\tau_1,\tau_2]$
and each $p\in M$,
we get a globally defined and smooth ${\mathcal L}$-exponential map
${\mathcal L}\exp_\tau:T_pM\ra M$ which takes each $v\in T_pM$
to $\gamma({\tau})$, where $v \: = \: \lim_{\tau^\prime \rightarrow \tau_1}
\sqrt{\tau^\prime} \frac{d\gamma}{d\tau^\prime}$.
Note that unlike in the case of Riemannian geometry, 
${\mathcal L}\exp_\tau(0)$ may
not be $p$, because of the $\nabla R$ term in
(\ref{Lgeod}).

We now fix $p\in M$, take $\tau_1=0$, and let $L(q,\bar\tau)$ be the
minimizer function as in Section \ref{I.7}.
We can imitate the traditional Riemannian geometry proof that geodesics
minimize for a short time.
Using the change of variable $s = \sqrt{\tau}$
and the implicit function theorem,  
there is an $r=r(p)>0$ (which varies continuously with $p$) such that
for every $q\in M$ with $d(q,p)\leq 10r$ at $\tau=0$, and every 
$0<\bar\tau\leq r^2$,
there is a unique ${\mathcal L}$-geodesic 
$\gamma_{(q,\bar\tau)}:[0,\bar\tau]\ra M$,
starting at $p$ and ending at $q$, which remains within the ball $B(p,100r)$
(in the $\tau=0$ slice $(M,g(0))$),
and $\gamma_{(q,\bar\tau)}$ varies smoothly with $(q,\bar\tau)$.
Thus, the ${\mathcal L}$-length of $\gamma_{(q,\bar\tau)}$ varies smoothly
with $(q,\bar\tau)$, and defines a function $\hat L(q,\bar\tau)$ near $(p,0)$.
We claim that $\hat L = L$ near $(p, 0)$. 
Suppose that $q \in B(p,r)$ and let $\alpha : [0, \bar \tau]
\rightarrow M$ be a smooth curve whose ${\mathcal L}$-length is
close to $L(q,\bar\tau)$. If $r$ is small, relative to the
assumed curvature bound, then $\alpha$ must stay within
$B(p, 10r)$. Equations
(\ref{grad1}) and (\ref{timed}) below imply that
\begin{equation}
\frac{d}{d\tau} \hat{L}(\alpha(\tau), \tau) \: = \:
\langle 2 \sqrt{\tau} X, \frac{d\alpha}{d\tau} \rangle \: + \: 
\sqrt{\tau} (R - |X|^2) \: \le \:
\sqrt{\tau} \left( R + \Big| \frac{d\alpha}{d\tau} \Big|^2 \right) \: = \: 
\frac{d}{d\tau} \: \left(
{\mathcal L} \text{length}(\alpha \big|_{[0, \tau]}) \right).
\end{equation}
Thus
$\gamma_{(q, \bar \tau)}$ minimizes when $(q, \bar \tau)$ is close to
$(p, 0)$.

We can now deduce that for all $(q,\bar\tau)$, there is an 
${\mathcal L}$-geodesic
$\gamma:[0,\bar\tau]\ra M$ which has infimal ${\mathcal L}$-length
among all piecewise smooth curves starting at $p$ and ending at $q$
(with domain $[0,\bar\tau]$).  This can be done by  imitating the 
usual broken geodesic argument, using the fact that 
for $x,y$ in a given small ball of $M$ and
for sufficiently
small time intervals $[\tau^\prime, \tau^\prime+\epsilon] \subset 
[0, \overline{\tau}]$, there is a unique
minimizer $\gamma$ for $\int_{\tau^\prime}^{\tau^\prime+\epsilon}
\sqrt{\tau} \left( R(\gamma(\tau)) \: +\: \left| \dot{\gamma}(\tau)
\right|^2 \right) \: d\tau$ with 
$\gamma(\tau^\prime) \: = \: x$ and
$\gamma(\tau^\prime+\epsilon) \: = \: y$. Alternatively, using the
change of variable $s \: = \: \sqrt{\tau}$, one can
take a minimizer of ${\mathcal L}$ among $H^{1,2}$-regular curves.

Another technical issue is the justification of the change of variables
from $M$ to $T_pM$ in the proof of monotonicity of reduced volume.
Fix $p\in M$ and $\tau>0$, and let ${\mathcal L}\exp_\tau:T_pM\ra M$
be the map which takes $v\in T_pM$ to $\gamma_v(\tau)$, where
$\gamma_v:[0,\tau]\ra M$ is the unique ${\mathcal L}$-geodesic
with $\sqrt{\tau^\prime}\frac{d\gamma_v}{d\tau^\prime}\ra v$ as 
$\tau^\prime\ra 0$.
Let ${\mathcal B}_\tau \subset M$ be the set of points which are
either endpoints of more than one minimizing ${\mathcal L}$-geodesic,
or which are the endpoint of a minimizing geodesic $\gamma_v:[0,\tau]\ra M$
where $v\in T_pM$ is a critical point of ${\mathcal L}\exp_\tau$.
We will call ${\mathcal B}_\tau$ 
the time-$\tau$ ${\mathcal L}$-cut locus of $p$. It is a closed
subset of $M$.
Let ${\mathcal G}_\tau \subset M$ be the complement of ${\mathcal B}_\tau$ and
let $\Omega_\tau \subset T_pM$ be the corresponding set of initial
conditions for minimizing ${\mathcal L}$-geodesics.  
Then $\Omega_\tau$ is an open set, and ${\mathcal L}\exp_\tau$
maps it diffeomorphically onto ${\mathcal G}_\tau$.  We claim that
${\mathcal B}_\tau$ 
has measure zero.   By Sard's theorem, to prove this it suffices
to prove that the set ${\mathcal B}'_\tau$ of points 
$q\in {\mathcal B}_\tau$ which are 
regular values of ${\mathcal L}\exp_\tau$, has measure zero.
Pick $q\in {\mathcal B}'_\tau$, and distinct points $v_1,v_2\in T_pM$
such that $\gamma_{v_i}:[0,\tau]\ra M$ are both minimizing geodesics
ending at $q$.  Then ${\mathcal L}\exp_\tau$ is a local diffeomorphism
near each $v_i$. The first variation formula and the implicit function theorem
then show that there are neighborhoods
$U_i$ of $v_i$, and a smooth hypersurface $H$  passing through
$q$, such that if we have points $w_i\in U_i$ with 
\begin{equation}
q' \: = \: {\mathcal L}\exp_\tau(w_1)={\mathcal L}\exp_\tau(w_2)\quad 
\mbox{and}\quad{\mathcal L}\length(\gamma_{w_1})={\mathcal L}
\length(\gamma_{w_2}),
\end{equation}
then $q'$ lies on $H$.  Thus ${\mathcal B}'_\tau$ 
is contained in a countable union
of hypersurfaces, and hence has measure zero. 

Therefore one may compute the integral of any integrable function on $M$
by pulling it back to $\Omega_\tau \subset T_pM$ and using the change
of variables formula. 
Note that
if $\tau \: \le \: \tau^\prime$ then $\Omega_{\tau^\prime} \subset
\Omega_\tau$.

\section{I.(7.3)-(7.6). First derivatives of $L$}

In this section we do some preliminary calculations leading up to the
computation of the second variation of $L$.

A remark about the notation : $L$ is a function of
a point $q$ and a time $\overline{\tau}$. The notation
$L_{\overline{\tau}}$ refers to the partial derivative with
respect to $\overline{\tau}$, i.e. differentiation while keeping
$q$ fixed.  The notation $\frac{d}{d\tau}$ refers to differentiation
along an ${\mathcal L}$-geodesic, i.e. simultaneously
varying both the point and the time.

If $q$ is not in the time-$\overline{\tau}$ 
${\mathcal L}$-cut locus of $p$, let
$\gamma \: : \: [0, \overline{\tau}] \rightarrow M$ be the unique minimizing
${\mathcal L}$-geodesic from $p$ to $q$, with length $L(q, \overline{\tau})$.
If $c \: : \: (- \epsilon, \epsilon) \rightarrow M$ is a short curve
with $c(0) \: = \: q$, consider the $1$-parameter
family of minimizing ${\mathcal L}$-geodesics
$\tilde\gamma(s, \tau)$ 
with $\tilde\gamma(s, 0) \: = \: p$ and $\tilde\gamma(s, \overline{\tau})
\: =\: c(s)$. Putting $Y(\tau) \: =\: 
\frac{\partial \tilde\gamma(s, \tau)}{\partial s} \Big|_{s=0}$,
equation (\ref{firstvariation}) gives
\begin{equation} 
\langle \nabla L, c^\prime(0) \rangle \: = \: 
\frac{dL(c(s), \overline{\tau})}{ds} \Big|_{s=0} \: = \:
2 \: \sqrt{\overline{\tau}} \:
\langle X(\overline{\tau}), Y(\overline{\tau}) \rangle. 
\end{equation}
Hence 
\begin{equation} \label{grad1}
(\nabla L)(q, \overline{\tau}) \: = \: 
2 \: \sqrt{\overline{\tau}} \: X(\overline{\tau})
\end{equation}
and
\begin{equation} \label{gradsq}
|\nabla L|^2(q, \overline{\tau}) \: = \: 
4 \: {\overline{\tau}} \: |X(\overline{\tau})|^2 \: = \:
- \: 4 \: {\overline{\tau}} \: R(q) \: + \:
4 \: {\overline{\tau}} \: \left( R(q) \: + \: |X(\overline{\tau})|^2 \right).
\end{equation}

If we simply extend the ${\mathcal L}$-geodesic $\gamma$ 
in $\overline{\tau}$, we
obtain
\begin{equation} \label{grad2}
\frac{dL(\gamma(\overline{\tau}),\overline{\tau})}{d\overline{\tau}} \: = \: 
\sqrt{\overline{\tau}} \: \left( R(\gamma(\overline{\tau})) \: + \:
|X(\overline{\tau})|^2 \right).
\end{equation}
As
\begin{equation}
\frac{dL(\gamma(\overline{\tau}),\overline{\tau})}{d\overline{\tau}} \: = \:
L_{\overline{\tau}}(q, \overline{\tau}) \: + \:
\langle (\nabla L)(q, \overline{\tau}), X(\overline{\tau}) \rangle,
\end{equation}
equations (\ref{grad1}) and (\ref{grad2}) give
\begin{align} \label{timed}
L_{\overline{\tau}}(q, \overline{\tau}) \: & = \:
\sqrt{\overline{\tau}} \: \left( R(q) \: + \:
|X(\overline{\tau})|^2 \right) \: - \: 
\langle (\nabla L)(q, \overline{\tau}), X(\overline{\tau}) \rangle \\
&=\: 
2 \sqrt{\overline{\tau}} R(q) \:
- \: \sqrt{\overline{\tau}} \: \left( R(q) \: + \:
|X(\overline{\tau})|^2 \right). \notag
\end{align}

When computing $\frac{dl(\gamma(\tau), \tau)}{d\tau}$, it will be useful
to have a formula for $R(\gamma(\tau)) \: + \: \big| X(\tau) \big|^2$.
As in I.(7.3),
\begin{equation}
\frac{d}{d\tau} \left( R(\gamma(\tau)) \: + \: \big| X(\tau) \big|^2 \right)
\: = \:
R_\tau \: + \: \langle \nabla R, X \rangle \: + \: 2 \langle \nabla_X X, X
\rangle \: + \: 2 \Ric(X, X).
\end{equation}
Using the ${\mathcal L}$-geodesic equation (\ref{Lgeod}) gives
\begin{align} \label{ddt}
\frac{d}{d\tau} \left( R(\gamma(\tau)) \: + \: \big| X(\tau) \big|^2 \right)
\: & = \:
R_\tau \: + \: \frac{1}{\tau} \: R \: + \:
2 \langle \nabla R, X \rangle \: 
- \: 2 \Ric(X, X) \: - \: \frac{1}{\tau} \: (R + |X|^2) \\
& = \: - \: H(X) \: - \: \frac{1}{\tau} \: (R + |X|^2), \notag
\end{align}
where 
\begin{equation} \label{Hexp}
H(X) \: = \: - \: R_\tau \: - \: \frac{1}{\tau} \: R \: - \:
2 \langle \nabla R, X \rangle \: 
+ \: 2 \Ric(X, X)
\end{equation}
is the expression of (\ref{Heqn})
after the change $\tau \: = \:- t$ and $X \rightarrow -X$.
Multiplying (\ref{ddt}) by $\tau^{\frac32}$ and integrating gives
\begin{equation} \label{parts}
\int_0^{\overline{\tau}} \tau^{\frac32} \: 
\frac{d}{d\tau} \left( R(\gamma(\tau)) \: + \: \big| X(\tau) \big|^2 \right)
\: d\tau
\: = \: - \: K \: - \: L(q, \overline{\tau}),
\end{equation}
where
\begin{equation}
K \: = \: \int_0^{\overline{\tau}} \tau^{\frac32} \: H(X(\tau)) \: d\tau.
\end{equation}
Then integrating the left-hand side of (\ref{parts}) by parts gives
\begin{equation} \label{sub4}
\overline{\tau}^{\frac32} \: 
\left( R(\gamma(\overline{\tau})) \: + \: \big| X(\overline{\tau}) 
\big|^2 \right)
\: = \: - \: K \: + \: \frac12 \: L(q, \overline{\tau}).
\end{equation}
Plugging this back into  (\ref{timed}) and (\ref{gradsq}) gives
\begin{equation} \label{touse1}
L_{\overline{\tau}}(q, \overline{\tau}) \: = \:
2 \sqrt{\overline{\tau}} R(q) \:
- \: \frac{1}{2\overline{\tau}} \: L(q,\overline{\tau}) \: + \:
\frac{1}{\overline{\tau}} \: K
\end{equation}
and
\begin{equation} \label{ttouse1}
|\nabla L|^2(q, \overline{\tau}) \: = \: 
- \: 4 \: {\overline{\tau}} \: R(q) \: + \:
\frac{2}{\sqrt{\overline{\tau}}} \: L(q, \overline{\tau}) \: - \:
\frac{4}{\sqrt{\overline{\tau}}} \: K.
\end{equation}

\section{I.(7.7). Second variation of ${\mathcal L}$}

In this section we compute the second variation of ${\mathcal L}$.
We use it to compute the Hessian of $L$ on $M$.

To compute the second variation $\delta_Y^2 {\mathcal L}$, we start with
the first variation equation
\begin{equation}
\delta_Y {\mathcal L} \: = \:
\int_{0}^{\overline{\tau}} \sqrt{\tau} \: \left( \langle Y, \nabla R \rangle
\: + \: 2 \: \langle \nabla_Y X, X \rangle \right) \: d\tau.
\end{equation}
Recalling that $\delta_Y \gamma(\tau) \: =\: Y(\tau)$ and 
$\delta_Y X(\tau) \: = \: (\nabla_Y X)(\tau)$,
the second variation is
\begin{align}
\delta_Y^2 {\mathcal L} \: & = \:
\int_{0}^{\overline{\tau}} \sqrt{\tau} \: \left( Y \cdot Y \cdot R
\: + \: 2 \: \langle \nabla_Y \nabla_Y X, X \rangle 
\: + \: 2 \: \langle \nabla_Y X, \nabla_Y X \rangle \right) \: d\tau \\
& = \: 
\int_{0}^{\overline{\tau}} \sqrt{\tau} \: \left( Y \cdot Y \cdot R
\: + \: 2 \: \langle \nabla_Y \nabla_X Y, X \rangle 
\: + \: 2 \: \big| \nabla_X Y \big|^2 \right) \: d\tau \notag \\
& = \: 
\int_{0}^{\overline{\tau}} \sqrt{\tau} \: \left( Y \cdot Y \cdot R
\: + \: 2 \: \langle \nabla_X \nabla_Y Y, X \rangle
\: + \: 2 \: \langle R(Y, X) Y, X \rangle 
\: + \: 2 \: \big| \nabla_X Y \big|^2 \right) \: d\tau, \notag
\end{align}
where the notation $Z \cdot u$ refers to the directional derivative, i.e.
$Z\cdot u\;=\; i_Z\;du$.
In order to deal with the $\langle \nabla_X \nabla_Y Y, X \rangle$ term,
we have to compute $\frac{d}{d\tau} \langle \nabla_Y Y, X \rangle$.

From the general equation for the Levi-Civita
connection in terms of the metric
\cite[(1.29)]{Chavel}, if $g(\tau)$ is a $1$-parameter family of metrics,
with $\dot{g} \: = \: \frac{dg}{d\tau}$
and $\dot{\nabla} \: = \: \frac{d\nabla}{d\tau}$, then
\begin{equation}
2 \langle \dot{\nabla}_X Y, Z \rangle \: = \:
(\nabla_X \dot{g}) (Y, Z) \: + \:
(\nabla_Y \dot{g}) (Z, X) \: - \:
(\nabla_Z \dot{g}) (X, Y).
\end{equation}
In our case $\dot{g} \: = \: 2 \: \Ric$ and so
\begin{align} \label{xxx}
\frac{d}{d\tau} \langle \nabla_Y Y, X \rangle \: = \: &
\langle \nabla_X \nabla_Y Y, X \rangle \: + \:
\langle \nabla_Y Y, \nabla_X X \rangle \: + \:
2 \Ric ( \nabla_Y Y, X ) \: + \:
\langle \dot{\nabla}_Y Y, X \rangle \\
 = \: &
\langle \nabla_X \nabla_Y Y, X \rangle \: + \:
\langle \nabla_Y Y, \nabla_X X \rangle \: + \notag \\
& 2 \Ric ( \nabla_Y Y, X ) \: + \:
2 (\nabla_Y \Ric) (Y, X) \: - \:
(\nabla_X \Ric) (Y, Y). \notag
\end{align}

(Although we will not need it, we can write
\begin{align}
2 Y \cdot \Ric(Y, X) \: - \: X \cdot \Ric(Y, Y) \: = \: &
2 (\nabla_Y \Ric)(Y, X) \: + \:
2 \Ric(\nabla_Y Y, X) \: + 2 \Ric(Y, \nabla_Y X) \\
& - \: 
(\nabla_X \Ric) (Y, Y) \: - \: 
2 \Ric(\nabla_X Y, Y) \notag \\
 = \: &
2 \Ric ( \nabla_Y Y, X ) \: + \:
2 (\nabla_Y \Ric) (Y, X) \notag \\
& - \:
(\nabla_X \Ric) (Y, Y) \: - \: 
2 \Ric([X, Y], Y).
\notag
\end{align}
We are assuming that the variation field $Y$ satisfies $[X, Y] \: = \: 0$
(this was used in deriving the ${\mathcal L}$-geodesic equation).
Hence one obtains the formula
\begin{equation}
\frac{d}{d\tau} \langle \nabla_Y Y, X \rangle \: = \:
\langle \nabla_X \nabla_Y Y, X \rangle \: + \:
\langle \nabla_Y Y, \nabla_X X \rangle \: + 
2 Y \cdot \Ric(Y, X) \: - \: X \cdot \Ric(Y, Y)
\end{equation}
of I.7.)

Next, using (\ref{xxx}),
\begin{align}
2 \sqrt{\overline{\tau}} \langle \nabla_Y Y, X \rangle \: = \: &
2 \int_0^{\overline{\tau}} \frac{d}{d\tau} 
\left( \sqrt{\tau} \langle \nabla_Y Y, X \rangle  \right) \: d\tau \\
 = \: &
\int_0^{\overline{\tau}}
\sqrt{\tau} \left[ \frac{1}{\tau} \: 
\langle \nabla_Y Y, X \rangle \: + \:
2  \:  \frac{d}{d\tau} \: \langle \nabla_Y Y, X \rangle \right]
\: d\tau \notag \\
\: = \: & 
\int_0^{\overline{\tau}}
\sqrt{\tau} \left[ \frac{1}{\tau} \: 
\langle \nabla_Y Y, X \rangle \: + \:
2 \langle \nabla_X \nabla_Y Y, X \rangle \: + \:
2 \langle \nabla_Y Y, \nabla_X X \rangle \: + \: \right. \notag \\
& \left. 
\: \: \: \: \:\: \: \: \: \: \: \: \: \: \: 4 \Ric ( \nabla_Y Y, X ) \: + \:
4 (\nabla_Y \Ric) (Y, X) \: - \: 2
(\nabla_X \Ric) (Y, Y) \right] \: d\tau \notag \\
\: = \: & 
\int_0^{\overline{\tau}}
\sqrt{\tau} \left[
2 \langle \nabla_X \nabla_Y Y, X \rangle \: + \:
(\nabla_Y Y) R \: + \: \right. \notag \\
& \left. 
\: \: \: \: \: \: \: \: \: \: \: \: \: \: \: 
4 (\nabla_Y \Ric) (Y, X) \: - \: 2
(\nabla_X \Ric) (Y, Y) \right] \: d\tau \notag \\
& - \:  
\int_0^{\overline{\tau}}
\sqrt{\tau} \left\langle \nabla_Y Y,
\nabla R
\: - \: 2 \nabla_X X  \: - \: 4 \: \Ric(X, \cdot) \: -\: \frac{1}{\tau}
\: X
\right\rangle \: d\tau. \notag
\end{align}
(Of course the last term vanishes if $\gamma$ is an ${\mathcal L}$-geodesic,
but we do not need to assume this here.)

The quadratic form $Q$ representing the Hessian of ${\mathcal L}$ on the
path space is given by
\begin{align}
Q(Y, Y) \:  = \: & \delta_Y^2 {\mathcal L} \: - \: \delta_{\nabla_Y Y} 
{\mathcal L} \\
= \:
& \delta_Y^2 {\mathcal L} \: - \: 
2 \sqrt{\overline{\tau}} \langle \nabla_Y Y, X \rangle \: 
-  \notag \\
& \int_0^{\overline{\tau}}
\sqrt{\tau} \left\langle \nabla_Y Y,
\nabla R
\: - \: 2 \nabla_X X  \: - \: 4 \: \Ric(X, \cdot) \: -\: \frac{1}{\tau}
\: X
\right\rangle. \notag
\end{align}
It follows that
\begin{align} \label{Hess}
Q(Y, Y) \: = \: &
\int_0^{\overline{\tau}}
\sqrt{\tau} \left[ Y \cdot Y \cdot R \: - \: (\nabla_Y Y) R \: + \: 
2 \langle R(Y, X) Y, X \rangle \: + \right. \\
& \left. 2 \: |\nabla_X Y|^2 
\: - \: 4 (\nabla_Y \Ric) (Y, X) \: + \: 2
(\nabla_X \Ric) (Y, Y) \right] \: d\tau \notag \\
= \: & \int_0^{\overline{\tau}}
\sqrt{\tau} \left[\Hess_R (Y,Y) \: + \: 
2 \langle R(Y, X) Y, X \rangle \: + \: 2 \: |\nabla_X Y|^2 \right. \notag \\
&
\: \: \: \: \: \: \: \: \: \: \: \: \: \: \:
 \left. - \: 4 (\nabla_Y \Ric) (Y, X) \: + \: 2
(\nabla_X \Ric) (Y, Y) \right] \: d\tau. \notag 
\end{align}
There is an associated second-order differential 
operator $T$ on vector fields $Y$ given by saying that
$2 \sqrt{\tau} \langle Y, TY\rangle$ equals the integrand of (\ref{Hess}) minus
$2 \frac{d}{d\tau} \left( \sqrt{\tau} \: \langle \nabla_X Y, Y \rangle 
\right)$.
Explicitly,
\begin{equation}
TY \: = \: - \: \nabla_X \nabla_X Y \: - \: \frac{1}{2\tau} \: \nabla_X Y
\: + \: \frac12 \: \Hess_R(Y, \cdot) \: - \: 2 \: (\nabla_Y \Ric)(X, \cdot)
\: - \: 2 \: \Ric(\nabla_Y X, \cdot).
\end{equation}
Then
\begin{equation} \label{Qeqn}
Q(Y, Y) \: = \: 2 \int_0^{\overline{\tau}}
\sqrt{\tau} \langle Y, TY \rangle \: d\tau \: + \:
2 \: \sqrt{\overline{\tau}} \: \langle \nabla_X Y(\overline{\tau}),
Y(\overline{\tau}) \rangle.
\end{equation}
An ${\mathcal L}$-Jacobi field along an ${\mathcal L}$-geodesic
is a field $Y(\tau)$ that is annihilated by $T$. 

The Hessian of the function
$L(\cdot, \overline{\tau})$ can be computed as follows.
Assume that $q \in M$ is outside of the time-$\overline{\tau}$
${\mathcal L}$-cut locus. Let $\gamma \: : \: [0, \overline{\tau}]
\rightarrow M$ be the minimizing ${\mathcal L}$-geodesic with
$\gamma(0) \: = \: p$ and $\gamma(\overline{\tau}) \: = \: q$.
Given $w \in T_qM$, take a short geodesic $c \: : \: (-\epsilon, 
\epsilon) \rightarrow M$ with $c(0) \: = \: q$ and $c^\prime(0) \: = \: w$.
Form the $1$-parameter
family of ${\mathcal L}$-geodesics $\tilde\gamma(s, \tau)$ with
$\tilde\gamma(s, 0) \: = \: p$ and 
$\tilde\gamma(s, \overline{\tau}) \: = \: c(s)$.
Then $Y(\tau) \: = \: \frac{\partial \tilde\gamma(s, \tau)}{\partial s} 
\Big|_{s=0}$ is an ${\mathcal L}$-Jacobi field $Y$ along $\gamma$ with
$Y(0) \: = \: 0$ and $Y(\overline{\tau}) \: = \: w$. 
We have
\begin{equation} \label{hesseq}
\Hess_L(w, w) \: = \: \frac{d^2 L(c(s), 
\overline{\tau})}{ds^2} \Big|_{s=0} \: = \:
Q(Y, Y) \: = \:
2 \: \sqrt{\overline{\tau}} \: \langle \nabla_X Y(\overline{\tau}),
Y(\overline{\tau}) \rangle.
\end{equation}

From (\ref{Qeqn}),
a minimizer of
$Q(Y, Y)$, among fields $Y$ with given values at the endpoints,
is an ${\mathcal L}$-Jacobi field.
It follows that 
$\Hess_L(w, w) \: \le \: Q(Y, Y)$ for {\em any} field $Y$ along $\gamma$
satisfying
$Y(0) \: = \: 0$ and $Y(\overline{\tau}) \: = \: w$.

\section{I.(7.8)-(7.9). Hessian bound for $L$} \label{7.9}

In this section we use a test variation field in order to estimate
the Hessian of $L$.

If $Y(\overline{\tau})$ 
is a unit vector at $\gamma(\overline{\tau})$, solve for
$\widetilde{Y}(\tau)$ in the equation
\begin{equation} \label{diffeq}
\nabla_X \widetilde{Y} \: = \: - \: \Ric(\widetilde{Y}, \cdot) 
\: + \: \frac{1}{2\tau} \:\widetilde{Y},
\end{equation}
on the interval $0<\tau\leq \overline{\tau}$ with the endpoint condition
$\widetilde{Y}(\overline{\tau}) \: = \: 
Y(\overline{\tau})$.
(For this section, we change notation from the
$Y$ in I.7 to $\widetilde{Y}$.)
Then
\begin{equation}
\frac{d}{d\tau} \langle \widetilde{Y}, \widetilde{Y} \rangle \: = \: 
2 \Ric(\widetilde{Y}, \widetilde{Y})
\: + \: 2 \langle \nabla_X \widetilde{Y}, \widetilde{Y} \rangle \: = \:
\frac{1}{\tau} \langle \widetilde{Y}, \widetilde{Y} \rangle,
\end{equation}
so
$\langle \widetilde{Y}(\tau), \widetilde{Y}(\tau) \rangle \: = \: 
\frac{\tau}{\overline{\tau}}$.
Thus we can extend $\widetilde{Y}$ continuously to the interval $[0,\overline{\tau}]$
by putting $\widetilde{Y}(0) \: = \: 0$.
Substituting into
(\ref{Hess}) gives
\begin{align} \label{sub1}
Q(\widetilde{Y}, \widetilde{Y}) \: = \: & \int_0^{\overline{\tau}}
\sqrt{\tau} \left[\Hess_R (\widetilde{Y},\widetilde{Y}) \: + \: 
2 \langle R(\widetilde{Y}, X) \widetilde{Y}, X \rangle 
\: + \: 2 \: \left| - \Ric(\widetilde{Y}, \cdot) \: + \:
\frac{1}{2\tau} \widetilde{Y} \right|^2 \right. \\
&
\: \: \: \: \: \: \: \: \: \: \: \: \: \: \:
 \left. - \: 4 (\nabla_{\widetilde{Y}} \Ric) (\widetilde{Y}, X) \: + \: 2
(\nabla_X \Ric) (\widetilde{Y}, \widetilde{Y}) \right] \: d\tau \notag \\
= \: & \int_0^{\overline{\tau}}
\sqrt{\tau} \left[\Hess_R (\widetilde{Y},\widetilde{Y}) \: + \: 
2 \langle R(\widetilde{Y}, X) \widetilde{Y}, X \rangle  \: + \: 2
(\nabla_X \Ric) (\widetilde{Y}, \widetilde{Y}) \: \right. \notag \\
&
\: \: \: \: \: \: \: \: \: \: \: \: \: \: \:
\left. - \: 4 (\nabla_{\widetilde{Y}} \Ric) (\widetilde{Y}, X) \:
+ \: 2 \: |\Ric(\widetilde{Y}, \cdot)|^2 \: - \: \frac{2}{\tau} 
\Ric(\widetilde{Y}, \widetilde{Y})
\: + \:
\frac{1}{2\tau \overline{\tau}}
\right] \: d\tau. \notag 
\end{align}
From
\begin{align}
\frac{d}{d\tau} \Ric(\widetilde{Y}(\tau), \widetilde{Y}(\tau)) \: & = \:
\Ric_\tau(\widetilde{Y}, \widetilde{Y}) \: + \: 
(\nabla_X \Ric)(\widetilde{Y}, \widetilde{Y}) \: + \: 2 
\Ric(\nabla_X \widetilde{Y}, \widetilde{Y}) \\
& = \: \Ric_\tau(\widetilde{Y}, \widetilde{Y}) \: + \: 
(\nabla_X \Ric)(\widetilde{Y}, \widetilde{Y}) \: + \:
\frac{1}{\tau} \: \Ric(\widetilde{Y}, \widetilde{Y}) \: - \: 
2 |\Ric(\widetilde{Y}, \cdot)|^2, \notag
\end{align}
one obtains
\begin{align} \label{sub2}
& - \: 2 \sqrt{\overline{\tau}} \Ric({Y}({\overline{\tau}}), 
{Y}({\overline{\tau}})) \: = \: - \: 2 \:
\int_0^{\overline{\tau}} \frac{d}{d\tau} \left(
\sqrt{\tau} \: \Ric(\widetilde{Y}, \widetilde{Y}) \right) \: d\tau \: =  \\
& - \: \int_0^{\overline{\tau}} \sqrt{\tau} 
\left[ \frac{1}{\tau} \: \Ric(\widetilde{Y}, \widetilde{Y}) 
\: + \: 2 \Ric_\tau (\widetilde{Y}, \widetilde{Y})
\: + \: 2 (\nabla_X \Ric)(\widetilde{Y}, \widetilde{Y}) \: + \:
\frac{2}{\tau} \: \Ric(\widetilde{Y}, \widetilde{Y}) 
\: - \: 4 |\Ric(\widetilde{Y}, \cdot)|^2 \right]. \notag
\end{align}
Combining (\ref{sub1}) and (\ref{sub2}) gives
\begin{equation} \label{sub3}
\Hess_L(Y({\overline{\tau}}), Y({\overline{\tau}})) \: \le \: 
Q(\widetilde{Y}, \widetilde{Y}) \: = \: 
\frac{1}{\sqrt{\overline{\tau}}}
 - \: 2 \sqrt{\overline{\tau}} \Ric({Y}({\overline{\tau}}), 
{Y}({\overline{\tau}})) \: - \: \int_0^{\overline{\tau}}
\sqrt{\tau} \: H(X, \widetilde{Y}) \: d\tau,
\end{equation}
where
\begin{align} \label{Htilde}
H(X, \widetilde{Y}) \: = \: & - \: \Hess_R 
(\widetilde{Y},\widetilde{Y}) \: - \: 2 \langle 
R(\widetilde{Y}, X)\widetilde{Y}, X \rangle
\: -\: 4 \left( \nabla_X \Ric(\widetilde{Y}, \widetilde{Y}) 
\: - \: \nabla_{\widetilde{Y}} \Ric(\widetilde{Y}, X) \right) \\
&  - \: 2 \Ric_\tau(\widetilde{Y}, \widetilde{Y}) \: + \: 
2 \big| \Ric(\widetilde{Y}, \cdot) \big|^2
\: - \: \frac{1}{\tau} \: \Ric(\widetilde{Y}, \widetilde{Y}). \notag
\end{align}
is the expression appearing in 
(\ref{Hexpression}), after the change $\tau = - t$ and $X \rightarrow -X$.
Note that $H(X, \widetilde{Y})$ is a quadratic form in
$\widetilde{Y}$. For its relation to the expression
$H(X)$ from (\ref{Hexp}), see Appendix \ref{appharnack}.

\section{I.(7.10). The Laplacian of $L$} \label{7.10}

In this section we estimate $\triangle L$.

Let $\{Y_i(\overline{\tau})\}_{i=1}^n$ be an orthonormal basis of
$T_{\gamma(\overline{\tau})}M$.
Solve for $\widetilde{Y}_i(\tau)$ from (\ref{diffeq}).
Putting
$\widetilde{Y}_i(\tau) 
\: = \: \left( \frac{\tau}{\overline{\tau}} \right)^{1/2}
e_i(\tau)$,
the vectors $\{e_i(\tau)\}_{i=1}^n$ form an orthonormal basis of
$T_{\gamma(\tau)}M$.
Substituting into 
(\ref{sub3}) and summing over $i$ gives
\begin{equation}
\triangle L \: \le \: 
\frac{n}{\sqrt{\overline{\tau}}} \: - \: 2 \sqrt{\overline{\tau}} R
\: - \: \frac{1}{\overline{\tau}} \: \int_0^{\overline{\tau}} 
\tau^{3/2} \: \sum_i H(X, e_i) \: d\tau.
\end{equation}
Then from (\ref{together}),
\begin{align} \label{touse2}
\triangle L \: \le \: &
\frac{n}{\sqrt{\overline{\tau}}} \: - \: 2 \sqrt{\overline{\tau}} R
\: - \: \frac{1}{\overline{\tau}} \: \int_0^{\overline{\tau}} 
\tau^{3/2} \: H(X) \: d\tau \\
= \: & \frac{n}{\sqrt{\overline{\tau}}} \: - \: 2 \sqrt{\overline{\tau}} R
\: - \: \frac{1}{\overline{\tau}} \: K. \notag
\end{align}

\section{I.(7.11). Estimates on ${\mathcal L}$-Jacobi fields}

In this section we estimate the growth rate of an ${\mathcal L}$-Jacobi field.

Given an ${\mathcal L}$-Jacobi field $Y$, we have
\begin{equation}
\frac{d}{d\tau} |Y|^2 \: =\: 2 \Ric(Y, Y) \: + \: 2 \langle \nabla_X Y, Y
\rangle \: = \: 2 \Ric(Y, Y) \: + \: 2 \langle \nabla_Y X, Y \rangle.
\end{equation}
Thus
\begin{equation}
\frac{d |Y|^2}{d\tau} \: \Big|_{\tau = \overline{\tau}}
= \: 2 \Ric(Y(\overline{\tau}), Y(\overline{\tau})) 
\: + \: \frac{1}{\sqrt{\overline{\tau}}} \:
\Hess_L(Y(\overline{\tau}), Y(\overline{\tau})), 
\end{equation}
where we have used (\ref{hesseq}).

Let $\tilde{Y}$ be a field along $\gamma$ as in Section \ref{7.9},
satisfying (\ref{diffeq}) with $\tilde{Y}(\overline{\tau}) \: = \:
{Y}(\overline{\tau})$ and $|{Y}(\overline{\tau})| \: = \: 1$.
Then from (\ref{sub3}),
\begin{equation}
\frac{1}{\sqrt{\overline{\tau}}} \:
\Hess_L(Y(\overline{\tau}), Y(\overline{\tau})) \: \le \:
\frac{1}{\overline{\tau}} \: - \:
2 \Ric(Y(\overline{\tau}),Y(\overline{\tau})) \: - \: 
\frac{1}{\sqrt{\overline{\tau}}} \int_0^{\overline{\tau}}
\sqrt{\tau} \: H(X, \tilde{Y}) \: d\tau.
\end{equation}
Thus
\begin{equation} \label{dY}
\frac{d |Y|^2}{d\tau} \: \Big|_{\tau = \overline{\tau}} \: \le \: 
\frac{1}{\overline{\tau}} \: - \:
\frac{1}{\sqrt{\overline{\tau}}} \int_0^{\overline{\tau}}
\sqrt{\tau} \: H(X, \tilde{Y}) \: d\tau.
\end{equation}
The inequality is sharp if and only if the first inequality in
(\ref{sub3}) is an equality.
This is the case if and only if $\tilde{Y}$ is actually the
${\mathcal L}$-Jacobi field $Y$, in which case
\begin{equation} \label{sharp}
\frac{1}{\overline{\tau}} \: = \: 
\frac{d |\tilde{Y}|^2}{d\tau} \: \Big|_{\tau = \overline{\tau}} \: = \:
\frac{d |Y|^2}{d\tau} \: \Big|_{\tau = \overline{\tau}}
\: = \:
2 \Ric(Y(\overline{\tau}),Y(\overline{\tau})) \: + \:
\frac{1}{\sqrt{\overline{\tau}}} \:
\Hess_L(Y(\overline{\tau}), Y(\overline{\tau})).
\end{equation}

\section{Monotonicity of 
the reduced volume $\widetilde{V}$} \label{monoV}

In this section we show that the reduced volume $\widetilde{V}(\tau)$ is
monotonically nonincreasing in $\tau$.

Fix $p \in M$. Define $l(q, \tau)$ 
as in (\ref{reducedl}). In order to show that
$\tilde{V}(\tau)$ is well-defined in the noncompact case,
we will need a lower bound on $l(q, \tau)$. For later
use, we prove something slightly more general. Recall
that we are assuming that we have bounded curvature on
compact time intervals,
and that time slices are complete.

\begin{lemma} \label{lowerbound}
Given $0 < \overline{\tau}_1 \le \overline{\tau}_2$,
there constants $C_1, C_2 > 0$ so that for all $\tau \in 
[\overline{\tau}_1, \overline{\tau}_2]$ and all
$q \in M$, we have
\begin{equation} \label{lbound}
l(q, \tau) \: \ge \:
C_1 \: d(p,q)^2 \: 
- \: C_2.
\end{equation}
\end{lemma}
\begin{proof}
We write ${\mathcal L}$ in the form (\ref{changetos}).
Given an ${\mathcal L}$-geodesic $\gamma$ with
$\gamma(0) \: = \: p$ and $\gamma({\tau}) \: = \: q$, we obtain
\begin{equation}
{\mathcal L}(\gamma) \: \ge \:
\frac12 \int_0^{\sqrt{\tau}} \left| \frac{d\gamma}{ds} \right|^2 \: ds \:
- \: \const
\end{equation}
As the multiplicative change in the metric between times
$\overline{\tau}_1$ and $\overline{\tau}_2$ is bounded by a factor
$e^{\const \: (\overline{\tau}_2 - \overline{\tau}_1)}$, it follows that
$L(q, \tau) \: \ge \:
\const \: d(p,q)^2 \:
- \: \const$, where the distance is measured at time $\tau$.
The lemma follows.
\end{proof}

Define 
$\tilde{V}(\tau)$
as in (\ref{reducedv}).
As the volume of time-$\tau$ balls in $M$ increases at most
exponentially fast in the radius, it follows that
$\tilde{V}(\tau)$ is well-defined.
From the discussion in Section \ref{Lremarks}, we can write
\begin{equation}
\tilde{V}(\tau) \: = \: \int_{T_pM} \tau^{-\frac{n}{2}} \: 
e^{-l({\mathcal L}exp_\tau(v),\tau)} \: {\mathcal J}(v, \tau) \: 
\chi_{\tau}(v) \:
dv,
\end{equation} 
where ${\mathcal J}(v, \tau) \: = \:
\det \: d\left( {\mathcal L}exp_\tau \right)_v$ is the Jacobian factor
in the change of variable and 
$\chi_\tau$ is the characteristic function of the
time-$\tau$ domain $\Omega_\tau$ of Section \ref{Lremarks}. 

We first show
that for each $v$, the expression
$- \: \frac{n}{2} \ln(\tau) \: 
- \: l({\mathcal L}exp_\tau(v), \tau) \: + \: \ln {\mathcal J}(v, \tau)$ is
nonincreasing in $\tau$. Let $\gamma$ be the ${\mathcal L}$-geodesic
with initial vector $v \in T_pM$.
From (\ref{grad2}) and (\ref{sub4}),
\begin{equation} \label{bbound1}
\frac{dl(\gamma(\tau), \tau)}{d\tau} \: \Big|_{\tau = \overline{\tau}} \: = \:
- \: \frac{1}{2\overline{\tau}} \: l(\gamma(\overline{\tau})) \: + \:
\frac12 \: \left( R(\gamma(\overline{\tau})) \: + \:
|X(\overline{\tau})|^2 \right) \: =\:
- \: \frac12 \: \overline{\tau}^{-\frac32} \: K.
\end{equation}

Next, let 
$\{Y_i\}_{i=1}^n$ be a basis for the Jacobi fields along $\gamma$
that vanish at $\tau \: =\: 0$.
We can write
\begin{equation}
\ln {\mathcal J}(v, \tau)^2 \: = \:
\ln \det \: \left( \left( d\left( {\mathcal L}exp_\tau \right)_v \right)^*
d\left( {\mathcal L}exp_\tau \right)_v \right) \: = \: \ln \det(S(\tau))
\: + \: \const, 
\end{equation}
where $S$ is the matrix
\begin{equation}
S_{ij}(\tau) \: = \: \langle Y_i(\tau), Y_j(\tau) \rangle.
\end{equation}
Then
\begin{equation}
\frac{d\ln {\mathcal J}(v, \tau)}{d\tau} \: = \: \frac12 \: 
\Tr \left( S^{-1} \frac{dS}{d\tau} \right).
\end{equation}
To compute the derivative at $\tau \: = \: \overline{\tau}$, we can choose
a basis so that $S(\overline{\tau}) \: = \: I_n$, i.e.
$\langle Y_i(\overline{\tau}), Y_j(\overline{\tau}) \rangle
\: = \: \delta_{ij}$.
Then using (\ref{dY}) and computing as in Section \ref{7.10},
\begin{equation} \label{bbound2}
\frac{d\ln {\mathcal J}(v, \tau)}{d\tau} 
\: \Big|_{\tau = \overline{\tau}}
\: = \: \frac12 \: \sum_{i=1}^n
\frac{d|Y_i|^2}{d\tau} \: \Big|_{\tau = \overline{\tau}}
\: \le \: \frac{n}{2\overline{\tau}} \: - \:
\frac{1}{2} \: \overline{\tau}^{-\frac32} K.
\end{equation}
If we have equality then (\ref{sharp}) holds for each
$Y_i$, i.e. $2 \Ric \: + \: \frac{1}{\sqrt{\overline{\tau}}} \:
\Hess_L \: =\: \frac{g}{\overline{\tau}}$ at $\gamma(\overline{\tau})$.

From (\ref{bbound1}) and (\ref{bbound2}), we deduce that
$\tau^{-\frac{n}{2}} \: 
e^{-l({\mathcal L}exp_\tau(v),\tau)} \: {\mathcal J}(v, \tau)$ is
nonincreasing in $\tau$. 
Finally, recall that if
$\tau \: \le \: \tau^\prime$ then $\Omega_{\tau^\prime} \subset
\Omega_\tau$, so $\chi_\tau(v)$ is nonincreasing in $\tau$. Hence
$\tilde{V}(\tau)$ is nonincreasing in $\tau$. If it is not
strictly decreasing then we must have
\begin{equation}
2 \Ric(\tau) \: + \: \frac{1}{\sqrt{\tau}} \: \Hess_{L(\tau)}
\: = \: \frac{g(\tau)}{\tau}.
\end{equation}
on all of $M$. Hence we have a gradient shrinking soliton solution.

\section{I.(7.15). A differential inequality for $L$} \label{I.(7.15)}

In this section we discuss an important differential inequality
concerning the reduced length $l$. We use the differential
inequality to estimate  $\min l(\cdot, \tau)$ from above.
We then give a lower bound on $l$.

With $\overline{L}(q, \tau) \: = \: 2 \sqrt{\tau} \: L(q, \tau)$, 
equations (\ref{touse1}) and (\ref{touse2}) imply that
\begin{equation} \label{barrier}
\overline{L}_{\overline{\tau}} \: + \: \triangle \overline{L}
\: \le \: 2n
\end{equation}
away from the time-$\overline{\tau}$ 
${\mathcal L}$-cut locus of $p$. 
We will eventually apply the maximum principle to this
differential inequality.  However to do so, we must 
discuss several  senses in which
the inequality can be made global on $M$, i.e. how it can be interpreted
on the cut locus.

The first sense is that of a barrier differential inequality.
Given $f \in C(M)$ and a function $g$ on $M$, one says that
$\triangle f \: \le \: g$ {\em in the sense of barriers}
if for all $q \in M$ and $\epsilon > 0$,
there is a neighborhood $V$ of $q$ and some $u_\epsilon \in C^2(V)$ so that 
 $u_\epsilon(q) = f(q)$, $u_\epsilon \ge f$ on $V$ and 
$\triangle u_\epsilon \: \le \: g + \epsilon$ on $V$
\cite{Calabi}. 
There is a similar spacetime definition for
$\triangle f \: -\: \frac{\partial f}{\partial t} \: \le \: g$
\cite{Dodziuk}. The point of barrier differential inequalities is
that they allow one to apply the maximum principle just as with
smooth solutions.

We illustrate this by constructing a barrier function for
$\overline{L}$ in (\ref{barrier}). Given the spacetime point
$(q, \overline{\tau})$, let $\gamma \: : \: [0, \overline{\tau}]
\rightarrow M$ be a minimizing ${\mathcal L}$-geodesic with
$\gamma(0) = p$ and $\gamma(\overline{\tau}) = q$. Given a small
$\epsilon > 0$, let 
$u_\epsilon(q^\prime, \overline{\tau}^\prime)$ be the minimum of
\begin{equation}
\int_{\epsilon}^{\overline{\tau}^\prime} \sqrt{\tau} \: \left( 
R(\gamma_2(\tau)) \: + \:
\big| \dot{\gamma_2}(\tau) \big|^2 \right) \: d\tau \: + \:
\int_{0}^{\epsilon} \sqrt{\tau} \: \left( 
R(\gamma(\tau)) \: + \:
\big| \dot{\gamma}(\tau) \big|^2 \right) \: d\tau 
\end{equation}
among curves $\gamma_2 \: : \: [\epsilon, \overline{\tau}^\prime] 
\rightarrow M$
with $\gamma_2(\epsilon) = \gamma(\epsilon)$ and
$\gamma_2(\overline{\tau}^\prime) \: = \: q^\prime$.
Because the
new basepoint $(\gamma(\epsilon), \epsilon)$ 
is moved in along $\gamma$
from $p$, the minimizer $\ga_2$  will be unique and will vary smoothly with $q'$, when
$q'$ is close to $q$; otherwise a second minimizer or a ``conjugate 
point'' would imply that $\ga$ was not minimizing. Thus the function $u_\epsilon$ 
is smooth in a spacetime 
neighborhood $V$ of
$(q, \overline{\tau})$. Put 
$U_\epsilon(q^\prime, \overline{\tau}^\prime) \: = \:
2 \sqrt{\overline{\tau}^\prime} \: 
u_\epsilon(q^\prime, \overline{\tau}^\prime)$.
By construction, $U_\epsilon \: \ge \: \overline{L}$ in $V$ and
$U_\epsilon(q, \overline{\tau}) \: = \: \overline{L}(q, \overline{\tau})$.
For small $\epsilon$, 
$(U_\epsilon)_{\overline{\tau}^\prime} \: + \: \triangle U_\epsilon$ 
will be bounded above
on $V$ by something close to $2n$. Hence
$\overline{L}$ satisfies (\ref{barrier}) globally on $M$ in the barrier sense.

As we are assuming bounded curvature on compact time intervals,
we can now apply the maximum principle of Appendix \ref{maxprin}
to conclude that the minimum of $\overline{L}(\cdot, \overline{\tau})
\: - \: 2n \overline{\tau}$ is nonincreasing in $\overline{\tau}$.  
(Note that from Lemma \ref{lowerbound}, the minimum of
$\overline{L}(\cdot, \overline{\tau})
\: - \: 2n \overline{\tau}$ exists.)

\begin{lemma} \label{neg}
For small positive $\overline{\tau}$, we have 
$\min \overline{L}(\cdot, \overline{\tau})
\: - \: 2n \overline{\tau} < 0$.
\end{lemma}
\begin{proof}
Consider the
static curve at the point $p$. Then 
for small $\overline{\tau}$, we have 
$\overline{L}(\cdot, \overline{\tau}) \le \const \overline{\tau}^2$,
from which the claim follows.
\end{proof}

(Being a bit more careful with the
estimates in the proof of Lemma \ref{lowerbound}, one sees that
$\lim_{\overline{\tau} \rightarrow 0} 
\min \overline{L}(\cdot, \overline{\tau}) \: = \: 0$.)
Then for $\overline{\tau} > 0$, we must have
$\min \overline{L}(\cdot, \overline{\tau}) \: \le \: 
2n \overline{\tau}$, so $\min l(\cdot, \overline{\tau}) \: \le \: 
\frac{n}{2}$.

The other sense of a differential inequality is the distributional
sense, i.e. $\triangle f \: \le \: g$ if for every nonnegative
compactly-supported smooth function $\phi$ on $M$,
\begin{equation}
\int_M (\triangle \phi) \: f \: dV \: \le \: \int_M \phi \: g \: dV.
\end{equation}
A general fact is that a barrier differential inequality implies
a distributional differential inequality
\cite{Ishii,Ye1}.

We illustrate this by giving an alternative proof that 
$\tilde{V}(\tau)$ is nonincreasing in $\tau$.
From (\ref{touse1}), (\ref{ttouse1}) and
(\ref{touse2}),
one finds that in the barrier sense (and hence in the distributional sense
as well)
\begin{equation}
l_{\overline{\tau}} \: - \: \triangle l \: + \: |\nabla l|^2 \: - \: R
\: + \frac{n}{2\overline{\tau}} \: \ge \: 0
\end{equation}
or, equivalently, that
\begin{equation} \label{illustrate}
\left( \partial_{\overline{\tau}} \: -\: \triangle \right) \:
\left( \tau^{-\frac{n}{2}} \: e^{-l}  \: dV \right) \: \le \: 0.
\end{equation}
Then for all nonnegative $\phi \in C^\infty_c(M)$ and
$0 \: < \: \overline{\tau}_1 \: \le \overline{\tau}_2$, one obtains
\begin{align} \label{pphieqn}
& \int_M \phi \: \overline{\tau}_2^{-\frac{n}{2}} \:
e^{-l(\cdot,\overline{\tau}_2)} \: dV(\overline{\tau}_2) \: - \:
\int_M \phi \: \overline{\tau}_1^{-\frac{n}{2}} \:
e^{-l(\cdot,\overline{\tau}_1)} \: dV(\overline{\tau}_1) \: 
\: = \: \\
& \int_{\overline{\tau}_1}^{\overline{\tau}_2} \int_M  \phi \: 
\left( \partial_{\overline{\tau}} \: -\: \triangle \right) \: \left(
\tau^{-\frac{n}{2}} \:
e^{-l(\cdot,\tau)} \: dV(\tau) \right) \: d\tau \: + \:
\int_{\overline{\tau}_1}^{\overline{\tau}_2} \int_M (\triangle \phi) \:
\tau^{-\frac{n}{2}} \:
e^{-l(\cdot,\tau)} \: dV(\tau) \: d\tau \: \le \notag \\
& 
\int_{\overline{\tau}_1}^{\overline{\tau}_2} \int_M (\triangle \phi) \:
\tau^{-\frac{n}{2}} \:
e^{-l(\cdot,\tau)} \: dV(\tau) \: d\tau.   \notag
\end{align} 
We can find a sequence $\{\phi_i\}_{i=1}^\infty$ of such functions $\phi$ 
with range $[0,1]$ so
that $\phi_i$ is one on $B(p, i)$, vanishes outside of $B(p, i^2)$, and
$\sup_M |\triangle \phi_i| \: \le \: i^{-1}$,
uniformly in $\tau \in [\overline{\tau}_1, \overline{\tau}_2]$.
Then to finish the argument it suffices to have a good upper bound on
$e^{-l(\cdot,\tau)}$ in terms of $d(p, \cdot)$, uniformly in 
$\tau \in [\overline{\tau}_1, \overline{\tau}_2]$.
This is given by Lemma \ref{lowerbound}.
The monotonicity of $\tilde{V}$ follows.

We also note the equation
\begin{equation} \label{illustrate2}
2 \triangle l \: - \: |\nabla l|^2  \: + \: R \:  + \:
\frac{l-n}{\overline{\tau}} \: \le \: 0,
\end{equation}
which follows from (\ref{ttouse1}) and
(\ref{touse2}).

Finally, suppose that the Ricci flow exists on a time interval
$\tau \in [0, \tau_0]$. From the maximum principle of
Appendix \ref{maxprin}, $R(\cdot, \tau) \: \ge \: - \: \frac{n}{2(\tau_0 - \tau)}$.
Then one obtains a lower bound on $l$ as before, using the
better lower bound for $R$.

\section{I.7.2. Estimates on the reduced length}

In this section we
suppose that our solution has nonnegative curvature operator.
We use this to derive estimates on the reduced length $l$.

We refer to Appendix \ref{appharnack} for
Hamilton's differential Harnack inequality.  We consider a Ricci flow defined on
a  time interval $t \in [0, \tau_0]$, with bounded nonnegative curvature operator,
and put $\tau \: = \: \tau_0 - t$. The differential Harnack inequality gives the
nonnegativity of the expression in (\ref{Hexpression}). 
Comparing this with the formula for $H(X,Y)$ in (\ref{Htilde}), we can write
the nonnegativity as
\begin{equation}
\left( H(X, Y) \: + \: \frac{\Ric(Y,Y)}{\tau} \right) \: + \: 
\frac{\Ric(Y,Y)}{\tau_0 - \tau} \: \ge \: 0,
\end{equation}  
or
\begin{equation} \label{ineqtoshow}
H(X, Y) \: \ge \: - \: \Ric(Y, Y) \left( \frac{1}{\tau} \: + \: 
\frac{1}{\tau_0 - \tau} \right).
\end{equation}
Then
\begin{equation}
H(X) \: \ge \: - \: R \: \left( \frac{1}{\tau} \: + \:
\frac{1}{\tau_0 - \tau} \right).
\end{equation}
As long as $\tau \: \le \: (1-c) \: \tau_0$, equations (\ref{gradsq}) and (\ref{sub4})
give
\begin{align}
4 \tau |\nabla l|^2 \: & = \: - \: 4 \tau R \: + \: 4l \: - \:
\frac{4}{\sqrt{\tau}} \int_0^\tau \tilde{\tau}^{3/2} H(X) \: d\tilde{\tau} \\
& \le \: - \: 4 \tau R \: + \: 4l \: + \:
\frac{4}{\sqrt{\tau}} \int_0^\tau \tilde{\tau}^{3/2} R \: 
\left( \frac{1}{\tilde{\tau}} \: + \:
\frac{1}{\tau_0 - \tilde{\tau}} \right)
  \: d\tilde{\tau} \notag \\
& = \: - \: 4 \tau R \: + \: 4l \: + \:
\frac{4}{\sqrt{\tau}} \int_0^\tau \sqrt{\tilde\tau} R \: 
\frac{\tau_0}{\tau_0 - \tilde{\tau}}
\: d\tilde{\tau} \notag \\
& \le \: - \: 4 \tau R \: + \: 4l \: + \:
\frac{4}{c \sqrt{\tau}} \int_0^\tau \sqrt{\tilde\tau} R 
\: d\tilde{\tau} \notag \\
& \le \: - \: 4 \tau R \: + \: 4l \: + \:
\frac{8l}{c}, \notag
\end{align}
where the last line uses (\ref{defL}).
Thus 
\begin{equation} \label{I.(7.16)}
|\nabla l|^2 \: + \: R \: \le \: \frac{Cl}{\tau}
\end{equation}
for a constant $C = C(c)$.
One shows similarly that
\begin{equation}
\frac{d}{d\tau} \ln |Y|^2 \: \le \: \frac{1}{\tau} \: (Cl+1)
\end{equation}
for an ${\mathcal L}$-Jacobi field $Y$, using (\ref{dY}).

\section{I.7.3. The no local collapsing theorem II} \label{I.7.3}

In this section we use the reduced volume to prove a no-local-collapsing
theorem for a Ricci flow on a finite time interval.

\begin{definition} \label{noncollapsed2}
We now say that a Ricci flow solution $g(\cdot)$ defined on a time
interval $[0, T)$ is {\em $\kappa$-noncollapsed on the scale
$\rho$} if for each
$r < \rho$ and all $(x_0, t_0) \in M \times [0,T)$ 
with $t_0 \ge r^2$, 
whenever it is true that
$|\Rm(x,t)| \: \le \: r^{-2}$ for every $x \in B_{t_0}(x_0, r)$
and $t \in [t_0 - r^2, t_0]$, then we also have 
$\vol(B_{t_0}(x_0, r)) \: \ge \: \kappa r^n$.
\end{definition}

Definition \ref{noncollapsed2} differs from
Definition \ref{noncollapsed1} by the requirement that the
curvature bound holds in the entire parabolic region 
$B_{t_0}(x_0, r) \times [t_0 - r^2, t_0]$ instead of just on
the ball $B_{t_0}(x_0, r)$ in the final time slice.  Therefore a Ricci
flow which is $\kappa$-noncollapsed in the sense of 
Definition \ref{noncollapsed1} is also 
$\kappa$-noncollapsed in the sense of 
Definition \ref{noncollapsed2}.

\begin{theorem} \label{nolocalcollapse}
Given numbers $n \in \Z^+$,  $T < \infty$
and $\rho, K, c > 0$, there is a number $\kappa = \kappa(n,K,c,\rho, T)>0$
with the following property.
Let $(M^n, g(\cdot))$ be a Ricci flow solution 
defined on a
time interval $[0, T)$ with $T < \infty$, such that the curvature $|\Rm|$
is bounded on every compact subinterval $[0,T']\subset [0,T)$. 
Suppose that $(M, g(0))$ is a complete
Riemannian manifold with $|\Rm| \: \le \: K$ and
$\inj(M, g(0)) \: \ge \: c \: > \: 0$. Then
the Ricci flow solution is $\kappa$-noncollapsed on the scale $\rho$,
in the sense of Definition \ref{noncollapsed2}.
Furthermore, with the other constants fixed, we can take
$\kappa$ to be nonincreasing in $T$.
\end{theorem}
\begin{proof}
We first observe that the existence of $\L$-geodesics and the 
monotonicity of the reduced volume are valid in this setting; see Section \ref{Lremarks}. 

Suppose that the theorem were false.  Then for given $T < \infty$
and $\rho, K, c > 0$, there are : \\
1. A sequence
$\{(M_k, g_k(\cdot))\}_{k=1}^\infty$ of Ricci flow solutions,
each defined on the
time interval $[0, T)$, with $|\Rm| \: \le \: K$ on $(M_k, g_k(0))$ and
$\inj(M_k, g_k(0)) \: \ge \: c$, \\
2. Spacetime points
$(p_k, t_k) \in M_k \times [0, T)$ and \\
3. Numbers $r_k \in (0, \rho)$ \\
having the following property :
$t_k \ge r_k^2$ and
if we put $B_k \: = \: B_{t_k}(p_k, r_k) \subset M_k$ then
$|\Rm|(x, t) \: \le \: r_k^{-2}$ whenever $x \in B_k$ and
$t \in [t_k - r_k^2, t_k]$, but
$\epsilon_k \: = \: r_k^{-1} \: \vol(B_k)^{\frac1n} \rightarrow 0$
as $k \rightarrow \infty$.
From short-time curvature estimates along with the assumed bounded geometry
at time zero, there is some $\overline{t} > 0$ so that
we have uniformly bounded geometry on the time interval
$[0, \overline{t}]$. In particular, we may assume that each
$t_k$ is greater than $\overline{t}$.

We define $\tilde{V}_k$ using curves
going backward in real time from the
basepoint $(p_k, t_k)$, i.e. forward in $\tau$-time from
$\tau = 0$.
The first step is to show that
$\tilde{V}_k(\epsilon_k r_k^2)$ is small. Note that 
$\tau \:  = \: \epsilon_k r_k^2$ corresponds to a real time of
$t_k - \epsilon_k r_k^2$, which is very close to $t_k$.

Given an ${\mathcal L}$-geodesic $\gamma(\tau)$ with
$\gamma(0) \: = \: p_k$ and
velocity vector
$X(\tau) \: = \: \frac{d\gamma}{d\tau}$, its initial vector is
$v \: = \: \lim_{\tau \rightarrow 0} \sqrt{\tau} \: X(\tau)
\in T_{p_k} M_k$.
We first want
to show that if $|v| \: \le \: .1 \: \epsilon_k^{- 1/2}$ then
$\gamma$ does not escape from $B_k$ in time
$\epsilon_k r_k^2$.

We have
\begin{align}
\frac{d}{d\tau} \langle X(\tau), X(\tau) \rangle \: = \:
& 2 \: \Ric(X, X) \: + \: 2 \langle X, \nabla_X X \rangle \\
\: = \: &
2 \: \Ric(X, X) \: + \: \langle X, 
\nabla R \: - \: \frac{1}{\tau} X \: - \: 4 \Ric(X, \cdot)
 \rangle \notag \\
\: = \: &
 - \: \frac{|X|^2}{\tau} \: - \: 2 \: \Ric(X, X) \: + \: \langle X, 
\nabla R \rangle, \notag
\end{align}
so
\begin{equation}
\frac{d}{d\tau} \left( \tau |X|^2 \right) \: = \:
- \: 2 \: \tau \: \Ric(X, X) \: + \: \tau \: \langle X, 
\nabla R \rangle.
\end{equation}
Letting $C$ denote a generic $n$-dependent constant,
for $x \in B(p_k, r_k/2)$ and $t \in [t_k \: - \: r_k^2/2,
t_k]$, the fact that $g_k$ satisfies the Ricci flow gives
an estimate $|\nabla R|(x, t) \: \le \: C r_k^{-3}$, as follows
from the case $l=0$, $m=1$ of
Appendix \ref{applocalder}.
Then in terms of dimensionless variables,
\begin{equation}
\Big| \frac{d}{d(\tau/r_k^2)} \left( \tau |X|^2 \right) \Big| \: \le \:
C \: \tau \: |X|^2 \: + \: C \: (\tau/r_k^2)^{1/2} \: 
(\tau |X|^2)^{1/2}.
\end{equation}
Equivalently,
\begin{equation}
\Big| \frac{d}{d(\tau/r_k^2)} \left( \sqrt{\tau} \: 
|X| \right) \Big| \: \le \:
C \: \sqrt{\tau} \: |X| \: + \: C \: \left(
\tau/r_k^2 \right)^{1/2}.
\end{equation}

Let us rewrite this as
\begin{equation} \label{diffmess}
\Big| \frac{d}{d( \frac{\tau}{\epsilon_k r_k^2})} \left( 
\epsilon_k^{\frac12} \: \sqrt{\tau} \: 
|X| \right) \Big| \: \le \:
C \: \epsilon_k \: \left( \epsilon_k^{\frac12} \: \sqrt{\tau} \: |X| \right)
\: + \: C \: \epsilon_k^2 \: \left(
\frac{\tau}{\epsilon_k r_k^2} \right)^{1/2}.
\end{equation}
We are interested in the time range when $\frac{\tau}{\epsilon_k r_k^2}
\in [0,1]$ and the initial condition satisfies
$\lim_{\tau \rightarrow 0} \epsilon_k^{\frac12} \: \sqrt{\tau} \: 
|X|(\tau) \: \le \: .1$. Then because of the $\epsilon_k$-factors
on the right-hand side, it follows from (\ref{diffmess})
that for large $k$, we will have $\epsilon_k^{\frac12} \: \sqrt{\tau} \: 
|X|(\tau) \: \le \: .11$ for all $\tau \in [0, \epsilon_k r_k^2]$.
Next,
\begin{equation}
\int_0^{\epsilon_k r_k^2} |X(\tau)| \: d\tau \: = \:
\epsilon_k^{- \: \frac12} \: \int_0^{\epsilon_k r_k^2} 
\epsilon_k^{\frac12} \: \sqrt{\tau} \: 
|X|(\tau) \: \frac{d\tau}{\sqrt{\tau}} \: \le \:
.11 \:
\epsilon_k^{- \: \frac12} \: \int_0^{\epsilon_k r_k^2} 
\frac{d\tau}{\sqrt{\tau}} \: = \: .22 \: r_k.
\end{equation}

From the Ricci flow equation $g_\tau \: = \: 2 \: \Ric$, it follows 
that the metrics $g(\tau)$ between $\tau \: = \: 0$ and 
$\overline{\tau} \: = \: \epsilon_k r_k^2$
are $e^{C \epsilon_k}$-biLipschitz close to each other.  Then for
$\epsilon_k$ small,  the length
of $\gamma$, as measured with the metric at time $t_k$, will be at most
$.3 \: r_k$. This shows that $\gamma$ does not leave $B_k$ within
time $\epsilon_k r_k^2$.

Hence the contribution to $\tilde{V}_k(\epsilon_k r_k^2)$ coming from vectors
$v \in T_{p_k} M_k$ 
with $|v| \: \le \: .1 \: \epsilon_k^{- \: \frac12}$ is at most
$\int_{B_k} (\epsilon_k r_k^2)^{- \frac{n}{2}} \:
e^{- \: l(q, \epsilon_k r_k^2)} \: dq$. 
We now want to give a lower
bound on $l(q, \epsilon_k r_k^2)$ for $q \in B_k$. Given the
${\mathcal L}$-geodesic $\gamma \: : \: [0, \epsilon_k r_k^2] \rightarrow M_k$
with $\gamma(0) \: = \: p_k$ and $\gamma(\epsilon_k r_k^2) \: = \: q$,
we have
\begin{equation}
{\mathcal L}(\gamma) \: \ge \: \int_0^{\epsilon_k r_k^2} \sqrt{\tau} \:
R(\gamma(\tau)) \: d\tau \: \ge \: - \: 
\int_0^{\epsilon_k r_k^2} \sqrt{\tau} \:
n(n-1) \: r_k^{-2} \: d\tau \: = \: - \: \frac23 \: n(n-1) \: 
\epsilon_k^{\frac32} \: r_k. 
\end{equation}
Then
\begin{equation}
l(q, \epsilon_k r_k^2) \: \ge \: - \: \frac13 \: n(n-1) \: \epsilon_k.
\end{equation}
Thus the contribution to $\tilde{V}_k(\epsilon_k r_k^2)$ coming from vectors
$v \in T_{p_k}M_k$ 
with $|v| \: \le \: .1 \: \epsilon_k^{- \: \frac12}$ is at most
\begin{equation}
e^{\frac13 \: n(n-1) \: \epsilon_k} \: 
(\epsilon_k r_k^2)^{- \frac{n}{2}} \: \vol_{t_k - 
\epsilon_k r_k^2}(B_k) \: \le \:
e^{\frac13 \: n(n-1) \: \epsilon_k} \: 
e^{\const \: \frac{1}{r_k^2} \: \epsilon_k r_k^2} \:
\epsilon_k^{\frac{n}{2}},
\end{equation}
which is less than $2 \epsilon_k^{\frac{n}{2}}$ for large $k$.

To estimate
the contribution to $\tilde{V}_k(\epsilon_k r_k^2)$ coming from vectors
$v \in T_{p_k}M_k$ with $|v| \: > \: .1 \: \epsilon_k^{- \: \frac12}$, we
can use the previously-shown monotonicity of the integrand in $\tau$.
As $\tau \rightarrow 0$, the Euclidean calculation of Section \ref{7Euclidean}
shows that 
$\tau^{- n/2} \: e^{- \: 
l({\mathcal L} exp_\tau(v), \tau)} \:
{\mathcal J}(v, \tau) \rightarrow 2^n \: e^{- \: |v|^2}$.
Then for all $\tau > 0$ and all $v \in \Omega_\tau$,
\begin{equation}
\tau^{- n/2} \: e^{- \: 
l({\mathcal L}_\tau(v), \tau)} \:
{\mathcal J}(v, \tau) \: \le \: 2^n \: e^{- \: |v|^2},
\end{equation}
giving
\begin{equation}
\int_{T_{p_k}M_k - B(0, .1 \: \epsilon_k^{-1/2})}
\tau^{- n/2} \: e^{- \: l({\mathcal L}_\tau(v), \tau)} \:
J(v, \tau) \: \chi_\tau \: d^nv \: \le \:
2^n \int_{T_{p_k}M_k - B(0, .1 \: \epsilon_k^{-1/2})}
\: e^{- \: |v|^2} \: d^nv \: \le \:
 e^{- \:  \frac{1}{10 \epsilon_k}}
\end{equation}
for $k$ large.

The conclusion is that $\lim_{k \rightarrow \infty} 
\tilde{V}_k(\epsilon_k r_k^2) \: = \: 0$. We now claim that there is a
uniform positive lower bound on $\tilde{V}_k(t_k)$.

To estimate $\tilde{V}_k(t_k)$ (where $\tau \: = \: t_k$ corresponds to
$t = 0$), we choose a point $q_k$ at time $t = \frac{\overline{t}}{2}$, i.e. at
$\tau \: = \: t_k - \frac{\overline{t}}{2}$, for which 
$l(q_k, t_k - \overline{t}/2) 
\: \le \: \frac{n}{2}$; see Section \ref{I.(7.15)}. Then we consider the concatenation of a
fixed curve $\gamma^{(k)}_1 \: : \: 
[0,   t_k - \overline{t}/2] \rightarrow M_k$, having
$\gamma^{(k)}_1(0) \: = \: p_k$ and 
$\gamma^{(k)}_1(t_k - \overline{t}/2) \: = \: q_k$, 
with a fan of curves 
$\gamma^{(k)}_2 \: : \: [t_k - \overline{t}/2, t_k] \rightarrow M_k$
having $\gamma^{(k)}_2(t_k - \overline{t}/2) \: = \: q_k$. 
Because of the uniformly bounded geometry
in the spacetime region with $t \in [0, \overline{t}/2]$, 
we can get an upper bound in
this way
for $l(\cdot, t_k)$ in a region around $q_k$. Integrating 
$e^{-l(\cdot, t_k)}$, we get a positive lower bound on $\tilde{V}_k(t_k)$
that is uniform in $k$.

As
$\epsilon_k r_k^2 \rightarrow 0$, the monotonicity of $\tilde{V}$
implies that $\tilde{V}_k(t_k) \le \tilde{V}_k(\epsilon_k r_k^2)$
for large $k$, which is a contradiction.
\end{proof}

\section{I.8.3. Length distortion estimates}
\label{secI.8.3}

The distortion of distances under Ricci flow can be estimated in terms
of the Ricci tensor.
We first mention a crude estimate.

\begin{lemma} \label{crudelemma}
If $\Ric \: \le \: (n-1) \: K$ then for $t_1 > t_0$,
\begin{equation}
\frac{\dist_{t_1}(x_0, x_1)}{\dist_{t_0}(x_0, x_1)} \: \ge \:
e^{-(n-1)K(t_1-t_0)}.
\end{equation}
\end{lemma}
\begin{proof}
For any curve $\gamma \: : \: [0,a] \rightarrow M$, we have
\begin{equation}
\frac{d}{dt} L(\gamma) \: =\:
\frac{d}{dt} \int_0^a 
\sqrt{ \left\langle \frac{d\gamma}{ds}, \frac{d\gamma}{ds}
\right\rangle} \: ds \: = \:
- \: \int_0^a \Ric \left( \frac{d\gamma}{ds}, \frac{d\gamma}{ds}
\right) \: \frac{ds}{\left| \frac{d\gamma}{ds} \right|} \: \ge \:
- \: (n-1) \: K \: L(\gamma).
\end{equation}
Integrating gives
\begin{equation}
\frac{L(\gamma) \big|_{t_1}}{L(\gamma) \big|_{t_0}} \: \ge \:
e^{-(n-1)K(t_1-t_0)}.
\end{equation}
The lemma follows by taking 
$\gamma$ to be a minimal geodesic at time $t_1$ between
$x_0$ and $x_1$.
\end{proof}

\begin{remark} \label{similar}
By a similar argument, if
$\Ric \: \ge \: - \: (n-1) \: K$ then for $t_1 > t_0$,
\begin{equation}
\frac{\dist_{t_1}(x_0, x_1)}{\dist_{t_0}(x_0, x_1)} \: \le \:
e^{(n-1)K(t_1-t_0)}.
\end{equation}
\end{remark}

We can write the conclusion of Lemma \ref{crudelemma} as
\begin{equation} \label{conclusion}
\frac{d}{dt} \dist_t(x_0, x_1) \: \ge \: - \: (n-1) K \:
\dist_t(x_0, x_1), 
\end{equation}
where the derivative is interpreted in the sense of forward difference quotients.

The estimate in Lemma \ref{crudelemma} is multiplicative. 
We now give an estimate that is additive in the distance.

\begin{lemma} \label{distancedist} (cf. Lemma I.8.3(b))
Suppose $\dist_{t_0}(x_0, x_1) \: \ge \: 2r_0$, and
$\Ric(x, t_0) \: \le \: (n-1) \: K$ for all $x\in B_{t_0}(x_0, r_0) \cup
B_{t_0}(x_1, r_0)$. Then
\begin{equation}
\frac{d}{dt} \dist_t(x_0, x_1) \: \ge \:
- \: 2 \: (n-1) \: \left( \frac23 \: K \: r_0 \: + \: r_0^{-1} \right)
\end{equation}
at time $t = t_0$.
\end{lemma}
\begin{proof}
If $\gamma$ is a normalized minimal geodesic from $x_0$ to $x_1$ with
velocity field $X(s) = \frac{d\gamma}{ds}$ then
for any piecewise-smooth 
normal
vector field $V$ along $\gamma$ that vanishes at the
endpoints, the second variation formula gives
\begin{equation}
\int_0^{d(x_0, x_1)} \left( \Big| \nabla_X V \Big|^2 \: + \:
\langle R(V, X)V, X \rangle \right) \: ds \: \ge \: 0.
\end{equation}
Let $\{e_i(s)\}_{i=1}^{n-1}$ be a parallel orthonormal frame along
$\gamma$ that is perpendicular to $X$.
Put $V_i(s) \: = \: f(s) \: e_i(s)$, where
\begin{equation}
f(s) \: = \: 
\begin{cases}
\frac{s}{r_0} & \text{ if } 0 \le s \le r_0, \\
1 & \text{ if } r_0 \le s \le d(x_0, x_1)-r_0, \\
\frac{d(x_0, x_1)-s}{r_0} & \text{ if } d(x_0, x_1)-r_0 \le s \le d(x_0, x_1).
\end{cases}
\end{equation}
Then $\Big| \nabla_X V_i \Big| \: = \: |f^\prime(s)|$ and
\begin{equation}
\int_0^{d(x_0, x_1)}  \Big| \nabla_X V_i \Big|^2 \: ds \: = \:
2 \int_0^{r_0} \frac{1}{r_0^2} \: ds \: = \: \frac{2}{r_0}.
\end{equation}
Next,
\begin{align}
\int_0^{d(x_0, x_1)} 
\langle R(V_i, X)V_i, X \rangle \: ds \: = &
\int_0^{r_0} \frac{s^2}{r_0^2} \langle R(e_i, X)e_i, X \rangle \: ds \: +  \\
& \int_{r_0}^{d(x_0,x_1)-r_0} \langle R(e_i, X)e_i, X \rangle \: ds \: + \notag \\
& \int_{d(x_0,x_1)-r_0}^{d(x_0,x_1)} \frac{(d(x_0,x_1) - s)^2}{r_0^2} \langle R(e_i, X)e_i, X \rangle \: ds.
\notag
\end{align}
Then
\begin{align}
0 \: & \le \: \sum_{i=1}^{n-1} \int_0^{d(x_0, x_1)} \left( \Big| \nabla_X V_i \Big|^2 \: + \:
\langle R(V_i, X)V_i, X \rangle \right) \: ds \\ 
& = \: 
\frac{2(n-1)}{r_0} \: - \: \int_0^{d(x_0, x_1)} \Ric(X,X) \: ds \: + \:
\int_0^{r_0} \left( 1 - \frac{s^2}{r_0^2} \right) \: \Ric(X, X) \: ds \: + \notag \\
& \: \: \: \: \: \:
\int_{d(x_0,x_1)-r_0}^{d(x_0,x_1)}  \left( 1 - \frac{(d(x_0,x_1) - s)^2}{r_0^2} \right) \: \Ric(X, X) \: ds. \notag
\end{align}
This gives
\begin{align}
\frac{d}{dt} \: \dist_t(x_0,x_1) \: & = \: - \int_0^{d(x_0,x_1)} \Ric(X, X) \: ds \\
& \ge \:
- \: \frac{2(n-1)}{r_0} \: - \: \int_0^{r_0} \left( 1 - \frac{s^2}{r_0^2} \right) \: \Ric(X, X) \: ds \notag \\
& \: \: \: \: \: \: - \:
\int_{d(x_0,x_1)-r_0}^{d(x_0,x_1)}  \left( 1 - \frac{(d(x_0,x_1) - s)^2}{r_0^2} \right) \: \Ric(X, X) \: ds
\notag \\
& \ge \: - \: \frac{2(n-1)}{r_0} \: - \: 2(n-1) K \cdot \frac23 \: r_0, \notag
\end{align}
which proves the lemma.
\end{proof}

We now give an additive version of Lemma \ref{crudelemma}.
\begin{corollary} \label{bound1}
\cite[Theorem 17.2]{Hamilton}
If $\Ric \: \le \: K$ with $K > 0$ 
then for all $x_0, x_1 \in M$,
\begin{equation}
\frac{d}{dt} \: \dist_t(x_0, x_1) \: \ge \: - \: \const(n) \: K^{1/2}.
\end{equation} 
\end{corollary}
\begin{proof}
Put $r_0 \: = \: K^{-1/2}$. 
If $\dist_t(x_0, x_1) \: \le 2r_0$ then the corollary follows from
(\ref{conclusion}).
If $\dist_t(x_0, x_1) \: > 2r_0$ then
it follows from Lemma \ref{distancedist}.
\end{proof}

The proof of the next lemma is similar to that of Lemma
\ref{distancedist} and is given in I.8.

\begin{lemma} (cf. Lemma I.8.3(a)) \label{I.8.3a}
Suppose that
$\Ric(x, t_0) \: \le \: (n-1) \: K$ on $B_{t_0}(x_0, r_0)$. Then
the distance function $d(x,t) = \dist_t(x,x_0)$ satisfies 
\begin{equation}
d_t - \triangle d \: \ge \: 
- \: (n-1) \: \left( \frac23 \: K \: r_0 \: + \: r_0^{-1} \right)
\end{equation}
at time $t = t_0$, outside of $B_{t_0}(x_0, r_0)$. The inequality must
be understood in the barrier sense (see Section \ref{I.(7.15)}) if
necessary.
\end{lemma}

\section{I.8.2. No local collapsing propagates forward in time and
to larger scales}

This section is concerned with a localized version of the
no-local-collapsing theorem.  The main result, Theorem \ref{8.2thm},
says that noncollapsing propagates forward in time and to a
larger distance scale.
 
We first give a local version of Definition \ref{noncollapsed2}.

\begin{definition} (cf. Definition of I.8.1)
A Ricci flow solution is said to be {\em $\kappa$-collapsed at
$(x_0, t_0)$, on the scale $r>0$}, if $|\Rm|(x,t) \: \le \: r^{-2}$
for all $(x,t)\in B_{t_0}(x_0,r)\times [t_0-r^2,t_0]$, 
but $\vol(B_{t_0}(x_0,r^2)) \: \le \:
\kappa r^n$.
\end{definition}

\begin{theorem} (cf. Theorem I.8.2) \label{8.2thm}
For any $0<A<\infty$, there is some $\kappa = \kappa(A)>0$ with
the following property.   Let 
$g(\cdot)$ be a Ricci flow solution defined for
$t \in [0, r_0^2]$, having complete time slices and uniformly bounded
sectional curvature.
Suppose that $\vol(B_0(x_0, r_0)) \geq A^{-1} r_0^n$
and that $|\Rm|(x,t) \le \frac{1}{n r_0^2}$ for all
$(x,t)\in B_0(x_0,r_0)\times [0,r_0^2]$. Then the 
solution cannot
be $\kappa$-collapsed on a scale less than $r_0$ at any point
$(x, r_0^2)$ with $x \in B_{r_0^2}(x_0,Ar_0)$.
\end{theorem}

\begin{remark} \label{8.2fix}
In \cite[Theorem I.8.2]{Perelman} the assumption is that
$|\Rm|(x,t) \le r_0^{-2}$.
We make the slightly stronger assumption 
that $|\Rm|(x,t) \le \frac{1}{n r_0^2}$.
The extra factor of $n$
is needed in order to assert in the proof that the region
$\{(y, t) \: : \:
\dist_{\frac12}(y,x_0)
\le \frac{1}{10}, \: t \in [0, \frac12 ]\}$ has
bounded geometry; see below.  Clearly this change of
hypothesis does not make any substantial difference in the
sequel.
\end{remark}

\begin{proof}
We follow the lines of the proof of Theorem \ref{nolocalcollapse}.
By scaling, we can take $r_0 = 1$. Choose $x \in M$ with
$\dist_{1}(x,x_0) < A$. Define the reduced volume 
$\tilde{V}(\tau)$ by 
means of curves starting at $(x, 1)$. An effective lower bound
on $\tilde{V}(1)$ would imply that the solution is not $\kappa$-collapsed
at $(x, 1)$, on a scale less than $1$, for an appropriate $\kappa > 0$.

We first note that the geometry of the region $\{(y, t) \: : \:
\dist_{\frac12}(y,x_0)
\le \frac{1}{10}, \: t \in [0, \frac12 ]\}$ is uniformly bounded.
To see this, the upper sectional curvature bound implies that
$\Ric \: \le \: 1$, so the 
distance distortion estimate of
Section \ref{secI.8.3} implies that 
$B_{\frac12}(x_0, \frac{1}{10}) \subset B_0(x_0, 1)$.
In particular, $|\Rm|(y,t) \: \le \: \frac{1}{n}$ 
on the region. By Remark \ref{similar},
if $\dist_{\frac12}(y,x_0)
\le \frac{1}{10}$ and $t \in [0, \frac12 ]$ then
$B_0(y, \frac{1}{1000}) \subset B_t(y, \frac{1}{100})$.
The Bishop-Gromov inequality gives a lower bound for the
time-zero volume of $B_0(y, \frac{1}{1000})$, of the form
$\vol_0(B_0(y, \frac{1}{1000})) \: \ge \: C_1(n,A)$.
The Ricci flow equation then gives a lower bound for
the time-$t$ volume of $B_0(y, \frac{1}{1000})$, of the form
$\vol_t(B_0(y, \frac{1}{1000})) \: \ge \: C_2(n,A)$.
Thus the time-$t$ volume of $B_t(y, \frac{1}{100})$ satisfies
$\vol(B_t(y, \frac{1}{100})) \: \ge \: C_2(n,A)$. 
This, along with the uniform sectional curvature bound,
implies that the region has uniformly bounded geometry.

If we have an effective upper bound on $\min_y l(y, \frac12)$,
where $y$ ranges over points that satisfy $\dist_{\frac12}(y,x_0)
\le \frac{1}{10}$, then we obtain a lower bound on 
$\tilde{V} \left( 1 \right)$.
Thus it suffices to obtain an effective upper bound on $\min_y
l(y, \frac{1}{2})$ or,
equivalently, on $\min_y \overline{L}(y, \frac{1}{2})$ 
(as defined using ${\mathcal L}$-geodesics
from $(x,1)$) for $y$ satisfying $\dist_{\frac{1}{2}}(y, x_0) \: \le \:
\frac{1}{10}$.
Applying the maximum principle to (\ref{barrier}) gave an upper bound
on $\inf_M \overline{L}$. The idea is to spatially
localize this estimate near $x_0$,
by means of a radial function $\phi$.

Let $\phi = \phi(u)$ be a smooth function that equals $1$ on
$(- \infty, \frac{1}{20})$, equals infinity on 
$(\frac{1}{10}, \infty)$ and is increasing on 
$(\frac{1}{20}, \frac{1}{10})$, with
\begin{equation} \label{phieqn}
2 (\phi^\prime)^2/\phi - \phi^{\prime \prime} \: \ge \:
(2A+100n) \phi^\prime - C(A) \phi
\end{equation}
for some constant $C(A) < \infty$.
To satisfy (\ref{phieqn}),
it suffices to take $\phi(u) \: = \: \frac{1}{e^{(2A+100n) 
(\frac{1}{10} - u)} - 1}$ for $u$ near $\frac{1}{10}$.

We claim that $\overline{L} \: + \: 2n \: + \: 1 \: \ge 1$ for
$t \ge \frac12$. To see this, from the end of Section \ref{I.(7.15)},
\begin{equation}
R(\cdot, \tau) \: \ge \: - \: \frac{n}{2(1-\tau)}.
\end{equation}
Then for $\overline{\tau} \in [0, \frac12]$,
\begin{equation}
L(q, \overline{\tau}) \: \ge \: - \: 
\int_0^{\overline{\tau}} \sqrt{\tau} \: \frac{n}{2(1-\tau)}
d\tau \: \ge \: - \: n \:
\int_0^{\overline{\tau}} \sqrt{\tau} \:
d\tau \: = \: - \: \frac{2n}{3} \:
{\overline{\tau}}^{3/2}.
\end{equation}
Hence
\begin{equation}
\overline{L}(q, \overline{\tau}) \: = \: 2 \: \sqrt{\overline{\tau}} \:
L(q, \overline{\tau}) \: \ge \: - \frac{4n}{3} \: 
\overline{\tau}^2 \: \ge \:
- \: \frac{n}{3},
\end{equation}
which proves the claim.

Now put
\begin{equation}
h(y,t) \: = \: \phi(d(y,t) - A(2t-1)) \:
(\overline{L}(y,1-t) +2n+1),
\end{equation}
where $d(y,t) = \dist_t(y,x_0)$. 
It follows from the above claim that $h(y,t) \: \ge \: 0$
if $t \ge \frac12$.
Also,
\begin{equation}
\min_y h(y,1) \: \le \: h(x,1) \: = \: 
\phi(\dist_1(x,x_0) - A) \cdot (2n+1) \: = \: 2n+1.
\end{equation}
As $\phi$ is infinite on $(\frac{1}{10}, \infty)$ and
$\overline{L}(\cdot, \frac12) +2n+1 \ge 1$, the
minimum of $h(\cdot, \frac12)$ is achieved at some
$y$ satisfying $d(y,\frac12) \: \le \: \frac{1}{10}$.

The calculations in I.8 give
\begin{equation}
\square h \: \ge \: -(2n+C(A)) h
\end{equation}
at a minimum point of $h$,
where $\square \: = \: \partial_t - \triangle$. 
Then $\frac{d}{dt} h_{min}(t) \: \ge \: -(2n+C(A)) \: h_{min}(t)$,
so 
\begin{equation}
h_{min} \left( \frac12 \right) 
\: \le \: e^{n+\frac{C(A)}{2}} \: h_{min}(1) 
\: \le \: (2n+1) \: e^{n+\frac{C(A)}{2}}.
\end{equation}
It follows that 
\begin{equation}
\min_{y \: : \: d(y, \frac12) \le \frac{1}{10}} 
\overline{L}(y, \frac12) +2n+1 \: \le 
\: (2n+1) \: e^{n+\frac{C(A)}{2}}.
\end{equation}
This implies the theorem.
\end{proof}

\section{I.9. Perelman's differential Harnack inequality}
\label{I.9}

This section is concerned with a localized version of the ${\mathcal W}$-functional.
It is mainly used in I.10.  

Let $g(\cdot)$ be a Ricci flow solution on a manifold $M$, defined for
$t \in (a,b)$. Put
$\square \: = \: \partial_t - \triangle$. For $f_1, f_2 \in C^\infty_c((a,b) \times M)$, we have
\begin{align}
0 \: & = \: \int_a^b \frac{d}{dt} \int_M f_1(t,x) f_2(t,x) \: dV \: dt \\
&  = \:
\int_a^b \int_M ((\partial_t - \triangle) f_1) \: f_2 \: dV \: + \:
\int_a^b \int_M f_1 \: (\partial_t + \triangle \: - \: R) f_2 \: dV \notag \\
&  = \:
\int_a^b \int_M (\square f_1) \: f_2 \: dV \: - \:
\int_a^b \int_M f_1 \:  \square^* f_2 \: dV, \notag
\end{align}
where $\square^* \: = \: - \: \partial_t - \triangle \: + \: R$. In this sense,
$\square^*$ is the formal adjoint to $\square$.

Now suppose that the Ricci flow is defined for $t \in [0,T)$.
Suppose that 
\begin{equation}
u \: = \: (4\pi(T-t))^{- \: \frac{n}{2}} \: e^{-f}
\end{equation}
satisfies $\square^* u \: = \: 0$. Put
\begin{equation} \label{veqn}
v \: = \: [(T-t)(2\triangle f - |\nabla f|^2 + R) + f - n ] \: u.
\end{equation}
If $M$ is compact then using (\ref{id}),
\begin{equation}
{\mathcal W}(g_{ij},f,T-t) \: = \: \int_M v \: dV.
\end{equation}

\begin{proposition} (cf. Proposition I.9.1) \label{localized}
\begin{equation} \label{I.(9.1)}
\square^*v \: = \: - \: 2(T-t) \: \Big| R_{ij} + \nabla_i \nabla_j f - \frac{g_{ij}}{2(T-t)}
\Big|^2 \: u.
\end{equation}
\end{proposition}
\begin{proof}
We note that the right-hand side of I.(9.1) should be multiplied by $u$.

To prove the proposition, we first claim that 
\begin{equation} \label{der}
\frac{d\triangle}{dt} \: = \: 2 \: R_{ij} \: \nabla_i \nabla_j.
\end{equation}
To see this, for $f_1, f_2 \in C^\infty_c(M)$, we have
\begin{equation}
\int_M f_1 \: \triangle f_2 \: dV \: = \: - \:
\int_M \langle df_1, df_2 \rangle \: dV.
\end{equation}
Differentiating with respect to $t$ gives
\begin{equation}
\int_M f_1 \: \frac{d\triangle}{dt} f_2 \: dV \: - \:
\int_M f_1 \: \triangle f_2 \: R \: dV
\: = \: - \: 2 \:
\int_M \Ric(df_1, df_2) \: dV
\: + \: \int_M \langle df_1, df_2 \rangle \: R \: dV,
\end{equation}
so
\begin{equation}
\frac{d\triangle}{dt} f_2 \: - \:
R \: \triangle f_2
\: = \: 2 \: \nabla_i (R_{ij} \: \nabla_j f_2)
\: - \: \nabla_i (R \:  \nabla_i f_2).
\end{equation}
Then (\ref{der}) follows from the traced second Bianchi identity.

Next, one can check that $\square^* u \: = \: 0$ is equivalent to
\begin{equation}
(\partial_t \: + \: \triangle) \: f \: = \: \frac{n}{2} \: \frac{1}{T-t}
\: + \: |\nabla f|^2 \: - \: R.
\end{equation}
Then one obtains
\begin{align}
u^{-1} \: \square^* v \: = \: &
- \: (\partial_t \: + \: \triangle) \: \left[ (T - t) \:
(2 \triangle f \: - \: |\nabla f|^2 \: + \: R) \: + \: f \right] 
\: - \\
& 2 \langle \nabla \left[ (T - t) \:
(2 \triangle f \: - \: |\nabla f|^2 \: + \: R) \: + \: f \right], 
u^{-1} \nabla u \rangle \notag \\
\: = \: & 2 \triangle f \: - \: |\nabla f|^2 \: + \: R \: - \:
(T-t) \: (\partial_t \: + \: \triangle) \: 
(2 \triangle f \: - \: |\nabla f|^2 \: + \: R) \notag \\
&- \: 
(\partial_t \: + \: \triangle) \: f  \: + \:
2 \: (T-t) \: \langle \nabla
(2 \triangle f \: - \: |\nabla f|^2 \: + \: R), 
\nabla f \rangle \: + \: 2 \: |\nabla f|^2. \notag
\end{align}
Now
\begin{align}
(\partial_t \: + \: \triangle) \: 
(2 \triangle f \: - \: |\nabla f|^2 \: + \: R) \: = \:
& 2 (\partial_t \triangle) f \: + \: 2 \triangle \:
(\partial_t \: + \: \triangle) \: f \\
& - \: (\partial_t \: + \: \triangle) \: |\nabla f|^2 \: + \: 
(\partial_t \: + \: \triangle) \: R \notag \\
 \: = \: & 4 \: R_{ij} \: \nabla_i \nabla_j f \: \: + \:
2 \: \triangle \: (|\nabla f|^2 \: - \: R) \: - \: 2 \:
\Ric(df, df) \notag \\
& - \: 2 \langle \nabla f_t, \nabla f \rangle \: - \:
\triangle |\nabla f|^2 
\: + \: \triangle R \: + \: 2 \: |\Ric|^2 \: + \: \triangle R \notag \\
 \: = \: & 4 \: R_{ij} \: \nabla_i \nabla_j f \: \: + \:
2 \: \triangle \: |\nabla f|^2 \: - \: 2 \:
\Ric(df, df) \notag \\
& - \: 2 \langle \nabla 
(- \: \triangle f \: + \: |\nabla f|^2 \: - \: R),
\nabla f \rangle \: - \:
\triangle |\nabla f|^2  \: + \: 2 \: |\Ric|^2. \notag
\end{align}
Hence the term in $u^{-1} \: \square^* v$ proportionate to
$(T-t)^{-1}$ is
\begin{equation}
- \: \frac{n}{2} \: \frac{1}{T-t}.
\end{equation}
The term proportionate to $(T-t)^0$ is
\begin{equation}
2 \triangle f \: - \: |\nabla f|^2 \: + \: R \: - \: |\nabla f|^2
\: + \: R \: + \: 2 |\nabla f|^2 \: = \: 2 (\triangle f \: + \: R).
\end{equation}
The term proportionate to $(T-t)$ is $(T-t)$ times
\begin{align}
& - \: 4 \: R_{ij} \: \nabla_i \nabla_j f \: \: - \:
2 \: \triangle \: |\nabla f|^2 \: + \: 2 \:
\Ric(df, df) \: + \\
&  2 \langle \nabla 
(- \: \triangle f \: + \: |\nabla f|^2 \: - \: R),
\nabla f \rangle \: + \:
\triangle |\nabla f|^2  \: - \: 2 \: |\Ric|^2 \: + \notag \\
& 2 \: \langle \nabla
(2 \triangle f \: - \: |\nabla f|^2 \: + \: R), 
\nabla f \rangle \: = \notag \\
&- \: 4 \: R_{ij} \: \nabla_i \nabla_j f \: \: - \:
\triangle \: |\nabla f|^2 \: + \: 2 \:
\Ric(df, df) \: + \:  2 \langle \nabla \triangle f,
\nabla f \rangle \: - \: 2 \: |\Ric|^2 \: = \notag \\
&- \: 4 \: R_{ij} \: \nabla_i \nabla_j f \: \: - \:
2 \: |\Hess(f)|^2 \: - \: 2 \: |\Ric|^2. \notag
\end{align}
Putting this together gives
\begin{equation}
\square^* v \: = \: - \: 2 \: (T-t) \:
\Big|R_{ij} \: + \: \nabla_i \nabla_j f \: - \: \frac{1}{2(T-t)} \: g_{ij}
\Big|^2 \: u.
\end{equation}
This proves the proposition.
\end{proof}

As a consequence of Proposition \ref{localized},
\begin{align}
\frac{d}{dt} {\mathcal W}(g_{ij}, f, T-t) \: & = \:
\frac{d}{dt} \int_M v \: dV \: = \:
\int_M (\partial_t \: + \: \triangle \: - \: R)v \: dV \\ 
& = \: 2 \: (T-t) \: \int_M
\Big|R_{ij} \: + \: \nabla_i \nabla_j f \: - \: \frac{1}{2(T-t)} \: g_{ij}
\Big|^2 \: u \: dV. \notag
\end{align}
In this sense, Proposition \ref{localized} is a local version of
the monotonicity of ${\mathcal W}$.

\begin{corollary} (cf. Corollary I.9.2) \label{whenever}
If $M$ is closed, or whenever the maximum principle holds, then
$\max v/u$ is nondecreasing in $t$.
\end{corollary}
\begin{proof}
We note that the statement of Corollary I.9.2
should have $\max v/u$ instead of
$\min v/u$.

To prove the corollary, we have
\begin{equation}
(\partial_t \: + \: \triangle) \: \frac{v}{u} \: = \:
\frac{v \square^* u \: - \: u \square^* v}{u^2} \: - \:
\frac{2}{u} \: \left\langle \nabla u, \nabla \frac{v}{u} 
\right\rangle.
\end{equation}
As $\square^* u = 0$ and $\square^* v \le 0$,
the corollary now follows from the maximum principle.
\end{proof}

We now assume that the Ricci flow solution is defined on 
the closed interval $[0, T]$.

\begin{corollary} (cf. Corollary I.9.3) \label{corI.9.3}
Under the same assumptions, if the solution is defined for
$t \in [0,T]$
and $u$ tends to a $\delta$-function as
$t \rightarrow T$ then $v \le 0$ for all $t < T$.
\end{corollary}
\begin{proof}
Suppose that $h$ is a positive solution of $\square h = 0$. Then
\begin{equation}
\frac{d}{dt} \: \int_M h v \: dV \: = \:
\int_M \left( (\square h) \: v \: - \:
h \square^* v \right) \:  dV \: = \:
- \: \int_M h \: \square^* v \: dV \: \ge \: 0.
\end{equation}
As $t \rightarrow T$, the computation of $\int_M hv$ approaches
the flat-space calculation, which one finds to be zero; see
\cite{Ni} for details. 
(Strictly speaking, the paper \cite{Ni} deals with the case when 
$M$ is closed. It is indicated that the proof should extend to the
noncompact setting.)
Thus
$\int_M h(t_0) v(t_0) \: dV$ is nonpositive for all $t_0<T$. As $h(t_0)$ can be
taken to be an arbitrary positive function, and then flowed forward to
a positive solution of $\square h = 0$, it follows that
$v(t_0) \le 0$ for all $t_0 < T$.
\end{proof}

The next result compares the function $f$ used in the ${\mathcal  W}$-functional
and the function $l$ used in the reduced volume.

\begin{corollary} (cf. Corollary I.9.5) \label{compare}
Under the assumptions of the previous corollary, let $p \in M$ be the
point where the limit $\delta$-function is concentrated. Then
$f(q,t) \le l(q, T-t)$, where $l$ is the reduced distance defined using
curves starting from $(p, T)$. 
\end{corollary}
\begin{proof}
Equation (\ref{illustrate}) implies that
$\square^* \left( (4 \pi \tau)^{-n/2} \: e^{-l} \right) \: \le \: 0$.
(This corrects the statement at the top of page 23 of I.) 
From this and the fact that 
$\square^* \left( (4 \pi \tau)^{-n/2} \: e^{-f} \right) \: = \: 0$,
the argument of the proof of Corollary \ref{whenever}  gives that
$\max e^{f-l}$ is nondecreasing in $t$, so
$\max (f-l)$ is nondecreasing in $t$. As $t \rightarrow T$ one
obtains the flat-space result, namely that $f \: - \: l$
vanishes.  Thus $f(t) \: \le \: l(T - t)$ for all $t \in [0, T)$.
\end{proof}

\begin{remark}
To give an alternative proof of Corollary \ref{compare}, putting $\tau \: = \: T \: - \: t$,
Corollary I.9.4 of \cite{Perelman} says that for any smooth curve $\gamma$,
\begin{equation}
\frac{d}{d\tau} \: f(\gamma(\tau), \tau) \: \le \:
\frac{1}{2} \left( R(\gamma(\tau), \tau) \: + \: \big| \dot{\gamma}(t)
\big|^2 \right) \: - \: \frac{1}{2\tau} \: f(\gamma(\tau), \tau),
\end{equation}
or 
\begin{equation}
\frac{d}{d\tau} \: \left( \tau^{1/2} f(\gamma(\tau), \tau) \right) \: \le \:
\frac{1}{2} \: \tau^{1/2} \:
\left( R(\gamma(\tau), \tau) \: + \: \big| \dot{\gamma}(t)
\big|^2 \right).
\end{equation}
Take $\gamma$ to be a curve emanating from
$(p, T)$.
For small $\tau$,
\begin{equation}
f(\gamma(\tau), \tau) \sim d(p, \gamma(\tau))^2/4\tau =
O(\tau^0).
\end{equation}
Then
integration gives $\tau^{1/2} f \: \le \: \frac12 \: L$, or $f \: \le \: l$.
\end{remark}

\section{The statement of the pseudolocality theorem}
\label{statement of I.10.1}

The next theorem says that, in a localized sense, if the initial
data of a Ricci flow solution
has a lower bound on the scalar curvature and satisfies
an isoperimetric inequality close to that of Euclidean space
then there is a sectional curvature bound in a forward region.
The result is not used in the sequel.

\begin{theorem} (cf. Theorem I.10.1) \label{thmI.10.1}
For every $\alpha > 0$ there exist $\delta,\epsilon > 0$ with the
following property.  Suppose that we have a 
smooth 
pointed Ricci flow solution
$(M, (x_0, 0), g(\cdot))$  
defined for $t \in [0, (\epsilon r_0)^2]$,
such that each time slice
is complete.
Suppose that
for any $x \in B_0(x_0, r_0)$ and $\Omega \subset B_0(x_0, r_0)$,
we have $R(x,0) \ge - r_0^{-2}$ and $\vol(\partial \Omega)^n
\ge (1-\delta) \: c_n \: \vol(\Omega)^{n-1}$, where $c_n$ is the
Euclidean isoperimetric constant.  Then 
$|\Rm|(x,t) <
\alpha t^{-1} \: + \: (\epsilon r_0)^{-2}$
 whenever
$0 < t \le (\epsilon r_0)^2$ and 
$d(x,t) = \dist_t(x, x_0) \le
\epsilon r_0$.
\end{theorem}

The sectional curvature bound $|\Rm|(x,t) <
\alpha t^{-1} \: + \: (\epsilon r_0)^{-2}$ necessarily blows up
as $t \rightarrow 0$, as nothing was assumed about
the sectional curvature at $t = 0$.

We first sketch the idea of the proof of Theorem
\ref{thmI.10.1}.
It is an argument by contradiction.
One takes a Ricci flow solution that satisfies the assumptions and
picks a point $(\overline{x}, \overline{t})$ where the
desired curvature bound does not hold.  One can assume, roughly speaking,
that $(\overline{x}, \overline{t})$ is the first point in the given
solution where the bound does not hold. (This will give the curvature
bound needed for taking
a limit in a sequence of counterexamples.)  One now
considers the solution $u$ to the conjugate heat equation, starting as
a $\delta$-function at $(\overline{x}, \overline{t})$, and the 
corresponding function $v$. We know that $v \: \le \: 0$. The first
goal is to get a negative upper bound for the integral of $v$ over
an appropriate ball $B$ at
a time $\widetilde{t}$ near $\overline{t}$; see Section \ref{Claim3}. 
The argument to get such
a bound is by contradiction.  If there were not such a bound then one
could consider a rescaled sequence of counterexamples with
$\int_B v \: dV \: \rightarrow 0$, and try to take a 
limit. If one has the injectivity radius bounds needed to take a limit
then one obtains a limit solution with
$\int_B v \: dV \: = \: 0$, which implies that the limit solution is a
gradient shrinking soliton, which violates curvature assumptions. If one
doesn't have the injectivity radius bounds
then one can do a further rescaling to see that
in fact $\int_B v \: dV \: \rightarrow -\infty$ for some subsequence, 
which is a contradiction.

If $M$ is compact then
$\int_M v \: dV$ is monotonically nondecreasing in $t$.
As (\ref{I.(9.1)}) is a localized version of this statement,
whether $M$ is compact or noncompact
we can use a cutoff function $h$ and equation (\ref{I.(9.1)})
to get a negative upper bound on $\int_M h v \: dV$ at time $t = 0$.
Finally, $\int_M v \: dV$ is the expression that
appears in the
logarithmic Sobolev inequality. If the isoperimetric constant
is sufficiently close to the Euclidean value $c_n$ then one concludes
that $\int_M hv \: dV$ must be bounded below by a constant close to zero,
which contradicts the negative upper bound on $\int_M hv \: dV$.

\section{Claim 1 of I.10.1. A 
point selection 
argument}

In Theorem \ref{thmI.10.1}, 
we can assume that $r_0 = 1$ and $\alpha < \frac{1}{100n}$. Fix $\alpha$
and put $M_\alpha = \{(x,t) \: : \: |\Rm(x,t)| \: \ge \: \alpha t^{-1} \}$.

The next lemma says that if we have a point $(x,t)$ where the conclusion of
Theorem \ref{thmI.10.1} does not hold then there is another point
$(\overline{x}, \overline{t})$ with $|\Rm(\overline{x}, \overline{t})|$
large (relative to $\overline{t}^{-1}$) so that any other such point
$(x^\prime, t^\prime)$ either has 
$t^\prime > \overline{t}$
or is much farther from $x_0$ than $\overline{x}$ is.

\begin{lemma} (cf. Claim 1 of I.10.1) \label{Claim 1}
For any $A > 0$, if $g(\cdot)$ is a Ricci flow solution for
$t \in [0, \epsilon^2]$, with $A \epsilon < \frac{1}{100n}$, and
$|\Rm|(x,t) \ge \alpha t^{-1} + \epsilon^{-2}$ 
for some
$(x,t)$ satisfying 
$t \in (0, \epsilon^2]$
and $d(x,t) \le \epsilon$, then
one can find $(\overline{x}, \overline{t}) \in M_\alpha$ with
$\overline{t} \in (0, \epsilon^2]$ and $d(\overline{x}, \overline{t}) <
(2A+1) \epsilon$, such that
\begin{equation} \label{I.(10.1)}
|\Rm(x^\prime,t^\prime)| \: \le \: 4 \: |\Rm(\overline{x}, \overline{t})|
\end{equation}
whenever
\begin{equation} \label{I.(10.2)}
(x^\prime, t^\prime) \in M_\alpha, \: \: \: \: t^\prime \in (0, \overline{t}],
 \: \: \: \: d(x^\prime, t^\prime) \le d(\overline{x}, \overline{t}) +
A |\Rm|^{- \: \frac12}(\overline{x}, \overline{t}).
\end{equation}
\end{lemma}
\begin{proof}
The proof is by a point selection argument as in
Appendix \ref{pointpicking}. By assumption, there is a point $(x,t)$
satisfying $t \in (0, \epsilon^2]$, 
$d(x, t) \le \epsilon$ and
$|\Rm(x, t)| \ge \alpha t^{-1} + \epsilon^{-2}$. 
Clearly $(x,t) \in M_\alpha$. Define points
$(x_k, t_k)$ inductively as follows. First, $(x_1, t_1) = (x, t)$. Next,
suppose that $(x_k, t_k)$ is constructed but cannot be taken for
$(\overline{x}, \overline{t})$. Then there is some point
$(x_{k+1}, t_{k+1}) \in M_\alpha$ such that $0 < t_{k+1} \le t_k$,
$d(x_{k+1}, t_{k+1}) \le d(x_k, t_k) + A |\Rm|^{- \: \frac12}(x_k,t_k)$ and
$|\Rm|(x_{k+1},t_{k+1}) > 4|\Rm|(x_k, t_k)$. Continuing in this way,
the point $(x_k, t_k)$ constructed has
$|\Rm|(x_k, t_k) \ge 4^{k-1} |\Rm|(x_1, t_1) \ge 4^{k-1} \epsilon^{-2}$.
Then 
\begin{align}
d(x_k, t_k) & \le d(x_1, t_1) + A |\Rm|^{- \: \frac12}(x_1, t_1) + \ldots
+ A |\Rm|^{- \: \frac12}(x_{k-1}, t_{k-1}) \\
& \le \epsilon + 2A
|\Rm|^{- \: \frac12}(x_1, t_1) \le (2A+1) \epsilon. \notag
\end{align}
As the solution is smooth, the induction process must terminate after a 
finite number of steps and the last value $(x_k, t_k)$ can be taken
for $(\overline{x}, \overline{t})$.
\end{proof}

\section{Claim 2 of I.10.1. Getting parabolic regions}

In Lemma \ref{Claim 1}, we know that (\ref{I.(10.1)}) is satisfied under the
condition (\ref{I.(10.2)}). The spacetime region described in (\ref{I.(10.2)})
is
not a product region, due to the fact that $d(x, t)$ is time-dependent.
The next goal is to obtain the estimate (\ref{I.(10.1)}) on a product region in
spacetime; this will be necessary when taking limits of Ricci flow solutions.
To get the estimate on a product region, one needs to bound 
how fast distances are changing with respect to $t$.

\begin{lemma} (cf. Claim 2 of I.10.1) \label{Claim 2}
For the point $(\overline{x}, \overline{t})$ constructed in Lemma \ref{Claim 1},
\begin{equation} \label{I.(10.1a)}
|\Rm(x^\prime,t^\prime)| \: \le \: 4 \: |\Rm(\overline{x}, \overline{t})|
\end{equation}
holds whenever
\begin{equation} \label{I.(10.3)}
\overline{t} \: - \: \frac12 \: \alpha Q^{-1} \: \le \:
t^\prime \: \le \: \overline{t}, 
\: \: \: \:
\dist_{\overline{t}}(x^\prime, \overline{x}) \: \le \: \frac{1}{10} \: A Q^{- \: \frac12},
\end{equation}
where $Q = |\Rm(\overline{x}, \overline{t})|$.
\end{lemma}
\begin{proof}
We first claim that if $(x^\prime, t^\prime)$ satisfies 
$\overline{t} \: - \: \frac12 \: \alpha Q^{-1} \: \le \: t^\prime \: \le \:
\overline{t}$ and $d(x^\prime, t^\prime) \: \le \: d(\overline{x}, \overline{t})
\: + \: A Q^{-1/2}$ then $|\Rm|(x^\prime, t^\prime) \: \le \: 4Q$. To see this,
if $(x^\prime, t^\prime) \in M_\alpha$ then it is true by Lemma \ref{Claim 1}.  If
$(x^\prime, t^\prime) \notin M_\alpha$ then 
$|\Rm|(x^\prime, t^\prime) \: < \: \alpha (t^\prime)^{-1}$.
As $(\overline{x}, \overline{t}) \in M_\alpha$, we know that
$Q \: \ge \: \alpha \overline{t}^{-1}$. Then $t^\prime \: \ge \:
\overline{t} \: - \: \frac12 \: \alpha Q^{-1} \: \ge \: 
\frac12 \: \overline{t}$ and so $|\Rm|(x^\prime, t^\prime) \: < \: 2 \: 
\alpha \overline{t}^{-1} \: \le \: 2Q$.

Thus we have a uniform curvature bound on the time-$t^\prime$ distance ball
$B(x_0, d(\overline{x}, \overline{t})
\: + \: A Q^{-1/2})$, provided that
$\overline{t} \: - \: \frac12 \: \alpha Q^{-1} \: \le \: t^\prime \: \le \:
\overline{t}$. We now claim that the time-$\overline{t}$ ball
$B(x_0, d(\overline{x}, \overline{t})
\: + \: \frac{1}{10} \: A Q^{-1/2})$ lies in the time-$t^\prime$ distance ball
$B(x_0, d(\overline{x}, \overline{t})
\: + \: A Q^{-1/2})$. To see this,
applying Lemma \ref{distancedist} with $r_0 \: = \:
\frac{1}{100} A Q^{-1/2}$ and the above curvature bound, if
$x^\prime$ is in the time-$\overline{t}$ ball
$B(x_0, d(\overline{x}, \overline{t})
\: + \: \frac{1}{10} \: A Q^{-1/2})$ then
\begin{equation} \label{32.4}
\dist_{t^\prime}(x_0, x^\prime) \: - \: 
\dist_{\overline{t}}(x_0, x^\prime) \: \le \:
\frac12 \: \alpha Q^{-1} \: \cdot \: 2(n-1) \:
\left( \frac23 \cdot 4Q ( \frac{1}{100} A Q^{-1/2}) \: + \:
100 \: A^{-1} Q^{1/2} \right). 
\end{equation}
Assuming that $A$ is sufficiently large
(we'll take $A \rightarrow \infty$ later)
and using the fact that $\alpha \: < \: \frac{1}{100n}$, it follows that
$d(x^\prime, t^\prime) \: \le \: 
d(x^\prime, \overline{t}) \: + \: 
\frac{1}{2} A Q^{-1/2} \: \le \: 
d(\overline{x},\overline{t}) + A Q^{- \: \frac12}$, 
which is what we want to show.  We note that
the argument also shows that is indeed self-consistent to use the curvature
bounds in the application of Lemma \ref{distancedist}.

Now suppose that  $(x^\prime, t^\prime)$ satisfies 
(\ref{I.(10.3)}). By the triangle inequality,
$x^\prime$ lies in the time-$\overline{t}$ distance ball
$B(x_0, d(\overline{x}, \overline{t}) \: + \: \frac{1}{10} \: A Q^{-1/2})$.
Then $x^\prime$ is in the time-$t^\prime$ distance ball
$B(x_0, d(\overline{x}, \overline{t})
\: + \: A Q^{-1/2})$ and so $|\Rm(x^\prime, t^\prime)| \: \le \: 4Q$,
which proves the lemma.
\end{proof}

\section{Claim 3 of I.10.1. An upper bound on the integral of $v$} \label{Claim3}

We first make some remarks about the fundamental solution to the
backward heat equation. Let $(M, (\overline{x}, b), g(\cdot))$ be a smooth 
one-parameter family of complete pointed Riemannian manifolds,
parametrized by $t \in (a,b]$. The fundamental solution 
$u$ of the backward heat equation is a positive solution of
$\square^* u = 0$ on $M \times  (a,b)$ 
such that 
$ u(\cdot, t)$ converges  to $\delta_{\overline{x}}$ 
in the distributional sense,  as $t\ra b^-$.
It is constructed as follows (cf. \cite[Section 3]{Dodziuk}).
Let $\{ D_i \}_{i=1}^\infty$ be an exhaustion of $M$ by
an increasing sequence of smooth compact codimension-zero
submanifolds-with-boundary containing $\overline{x}$ in
the interior. Let $u^{(i)}$ be the unique solution of
$\square^* u^{(i)} = 0$ on $D_i \times (a,b)$ with 
$\lim_{t \rightarrow b^-} u^{(i)}(x, t) = \delta_{\overline{x}}(x)$,
as constructed using Dirichlet boundary conditions on $D_i$.
If $D_i \subset D_j$ then $u^{(i)} \le u^{(j)}$ on $D_i$,
using the maximum principle as in \cite[Lemma 3.1]{Dodziuk}.
Then the fundamental solution  is defined to be the limit 
$u=\lim_{i \rightarrow \infty} u^{(i)}$, with smooth  convergence
 on compact subsets of $M \times  (a,b)$.
 The function $u$ is
independent of the choice of exhaustion sequence $\{D_i\}_{i=1}^\infty$.
For any $t \in (a,b)$, we have $\int_M u(x,t) \: dV(x) \: \le \: 1$.
If $\int_M u(x,t) \: dV(x) \: = \: 1$ for all $t$ then 
we say that $(M, (\overline{x}, b), g(\cdot))$ is stochastically complete
for $\square^*$.
This will be the case if one has
bounded curvature on compact time intervals, but need not be the case in
general. 

\begin{lemma} \label{stoccomplete}
Let $\{(M_k, (\overline{x}_k, b), g_k(\cdot)) \}_{k=1}^\infty$ be a 
sequence of 
manifolds as above, 
each defined on the time interval
$(a,b]$. 
Suppose that $\lim_{k \rightarrow \infty}
(M_k, (\overline{x}_k, b), g_k(\cdot)) \: = \: 
(M_\infty, (\overline{x}_\infty, b), g_\infty(\cdot))$ 
in the pointed smooth topology, and that 
$(M_\infty, (\overline{x}_\infty, b), g_\infty(\cdot))$
is stochastically complete for $\square^*$.  
Then after passing to a subsequence,
the fundamental solutions $\{ u_k \}_{k=1}^\infty$ converge smoothly
on compact subsets of $M_\infty \times (a,b)$ to the fundamental
solution $u_\infty$.
(Of course, we use the
pointed diffeomorphisms inherent in the statement of pointed convergence
in order to compare the $u_k$'s with $u_\infty$.)
\end{lemma}
\begin{proof}
From the uniform upper $L^1$-bound on $\{ u_k(\cdot, t) \}_{k=1}^\infty$ 
and parabolic
regularity, after passing to a subsequence we can assume that
$\{ u_k\}_{k=1}^\infty$ converges smoothly on compact subsets of
$M_\infty \times (a,b)$ to some function $U$. From the construction of
$u_\infty$, it follows easily that $u_\infty \le U$. For any
$t \in (a,b)$, we have
\begin{align}
\int_{M_\infty} ( U(x,t) \: - \: u_\infty(x,t) ) \: \dvol(x) \: & = \:
\int_{M_\infty} \liminf_k ( u_k(x, t) \: - \: u_\infty(x,t) ) \: \dvol(x) \\
& \le \: \liminf_k \int_{M_\infty} ( u_k(x, t) \: - \: u_\infty(x,t) ) \: 
\dvol(x) \: \le \: 0, \notag
\end{align}
so $U \: = \: u_\infty$.
\end{proof}

Starting the proof of Theorem \ref{thmI.10.1}, we suppose that the
theorem is not true. Then there are
sequences $\epsilon_k \rightarrow 0$ and
$\delta_k \rightarrow 0$, and
pointed Ricci flow solutions $(M_k, (x_{0,k}, 0), g_k(\cdot))$
which
satisfy the hypotheses of the theorem but for which there is a point
$(x_k, t_k)$ with $0 < t_k \le \epsilon_k^2$, 
$d(x_k, t_k) \le \epsilon_k$
and 
$|\Rm|(x_k, t_k) \geq \alpha t_k^{-1} + \epsilon_k^{-2}$.
Given the flow $(M_k, g_k(\cdot))$, we reduce $\epsilon_k$ as much as
possible so that there is still such a point $(x_k, t_k)$. Then
\begin{equation} \label{largebound}
|\Rm|(x, t) < \alpha t_k^{-1} + 2\epsilon_k^{-2}
\end{equation}
whenever
$0 < t \le \epsilon_k^2$ and 
$d(x, t) \le \epsilon_k$.
Put $A_k \: = \: \frac{1}{100n\epsilon_k}$. Construct points
$(\overline{x}_k, \overline{t}_k)$ as in Lemma \ref{Claim 1}.
Consider fundamental
solutions $u_k = (4\pi(\overline{t}_k - t))^{- \: \frac{n}{2}} \: 
e^{-f_k}$ of $\square^* u_k = 0$ satisfying 
$\lim_{t \rightarrow \overline{t}_k^-} u(x, t) = \delta_{\overline{x}_k}(x)$.
Construct the corresponding functions
$v_k$ from (\ref{veqn}).

\begin{lemma} (cf. Claim 3 of I.10.1) \label{Claim 3}
There is some $\beta > 0$ so that for 
all sufficiently large $k$,
there is some
$\widetilde{t}_k \in [\overline{t}_k \: - \: \frac12 \: \alpha Q_k^{-1}, \overline{t}_k]$ with
$\int_{B_k} v_k \: dV_k \: \le \: - \beta$, where $Q_k \: = \: |\Rm|(\overline{x}_k,
\overline{t}_k)$ and $B_k$ is the time-$\widetilde{t}_k$ ball of radius
$\sqrt{\overline{t}_k - \widetilde{t}_k}$ centered at $\overline{x}_k$.
\end{lemma}
\begin{proof} 
Suppose that the claim is not true.  After passing to a subsequence, we
can assume that for any choice of $\widetilde{t}_k$, 
$\liminf_{k \rightarrow \infty}
\int_{B_k} v_k \: dV_k \: \ge \: 0$.

Consider the pointed solution $(M_k, (\overline{x}_k, \overline{t}_k), g_k(\cdot))$ 
parabolically rescaled by $Q_k$. 
Suppose first that there is a subsequence so that
the injectivity radii of the scaled metrics at
$(\overline{x}_k, \overline{t}_k)$ are bounded away from zero.
Since $A_k \rightarrow \infty$, we can use 
Lemma \ref{Claim 2} and
Appendix \ref{subsequence} to 
take a subsequence that converges to a complete Ricci flow solution 
$(M_\infty, (\overline{x}_\infty, \overline{t}_\infty), g_\infty(\cdot))$ on a
time interval $(\overline{t}_\infty - \frac12 \alpha, \overline{t}_\infty]$, with
$|\Rm| \le 4$ and $|\Rm|(\overline{x}_\infty, \overline{t}_\infty) = 1$.
Consider the fundamental solution $u_\infty$ of $\square^*$ on $M_\infty$
with $\lim_{t \rightarrow \overline{t}_\infty^-} u_\infty(x_\infty, t) \: = \:
\delta_{\overline{x}_\infty}(x_\infty)$.
As before,
let $u_k$ be the fundamental solution 
of $\square^*$ on $M_k$
with $\lim_{t \rightarrow \overline{t}_k^-} u_k(x_k, t) \: = \:
\delta_{\overline{x}_k}(x_k)$.
In view of the pointed convergence of the rescalings of 
$(M_k, (\overline{x}_k, \overline{t}_k), g_k(\cdot))$ to
$(M_\infty, (\overline{x}_\infty, \overline{t}_\infty), 
g_\infty(\cdot))$,
Lemma \ref{stoccomplete} implies that
after passing to a further subsequence we can ensure that
$\lim_{k \rightarrow \infty} u_k \: = \: u_\infty$, with 
smooth convergence on
compact subsets of $M_\infty \times (\overline{t}_\infty
 - \frac12 \alpha, \overline{t}_\infty)$.
(The curvature bounds on 
$(M_\infty, (\overline{x}_\infty, \overline{t}_\infty), 
g_\infty(\cdot))$ ensure that it is stochastically complete
for $\square^*$.)
From Corollary \ref{corI.9.3}, $v_\infty \: \le \: 0$.
Note that we are applying Corollary \ref{corI.9.3} on
$M_\infty \times (\overline{t}_\infty 
- \frac12 \alpha, \overline{t}_\infty)$,
where we have the curvature bounds needed to use the maximum
principle.

Given $\widetilde{t}_\infty \in (\overline{t}_\infty - 
\frac12 \alpha, \overline{t}_\infty)$,
let $B_\infty$ be the time-$\widetilde{t}_\infty$ ball of radius
$\sqrt{\overline{t}_\infty - \widetilde{t}_\infty}$ 
centered at $\overline{x}_\infty$.
In view of the smooth convergence 
$\lim_{k \rightarrow \infty} u_k \: = \: u_\infty$ on
compact subsets of 
$M_\infty \times (\overline{t}_\infty - \frac12 \alpha, 
\overline{t}_\infty)$,
it follows that
$\int_{B_\infty} v_\infty \: dV_\infty \: = \: 0$ at time $\widetilde{t}_\infty$, so
$v_\infty$ vanishes on $B_\infty$ at time $\widetilde{t}_\infty$. Let $h$ be a solution
to $\square h = 0$ on $M_\infty \times [\widetilde{t}_\infty, \overline{t}_\infty)$ with
$h(\cdot, \widetilde{t}_\infty)$ a nonnegative nonzero function supported in 
$B_\infty$. As in the proof of Corollary \ref{corI.9.3},
$\int_{M_\infty} h v_\infty \: dV_\infty$ is nondecreasing in $t$ and vanishes for
$t \: = \: \widetilde{t}_\infty$ and  $t \: \rightarrow \: \overline{t}_\infty$.
Thus $\int_{M_\infty} h v_\infty \: dV_\infty$ vanishes for all 
$t \in [\widetilde{t}_\infty, \overline{t}_\infty)$. However, for $t \in 
(\widetilde{t}_\infty, \overline{t}_\infty)$, $h$ is strictly positive and
$v_\infty$ is nonpositive.  Thus $v_\infty$ vanishes on $M_\infty$ for all
$t \in (\widetilde{t}_\infty, \overline{t}_\infty)$, and so
\begin{equation}
\Ric(g_\infty) \: + \: \Hess f_\infty \: - \:
\frac{1}{2(\overline{t} - t)} \: g_{\infty} \: = \: 0.
\end{equation}
on this interval.
We know that $|\Rm| \le 4$ on
$M_\infty \times (\overline{t}_\infty - 
\frac12 \alpha, \overline{t}_\infty]$. 
From the evolution equation,
\begin{equation}
\frac{dg_{\infty}}{dt} \: = \: - \: 2 \Ric(g_\infty) \: = \:
2 \Hess f_\infty \: - \:
\frac{1}{\overline{t} - t} \: g_{\infty}.
\end{equation}
It follows that the supremal and infimal
sectional curvatures of $g_\infty(\cdot,t)$ go like
$(\overline{t}_\infty \: - \: t)^{-1}$.
Hence $g_\infty$ is flat,
which contradicts the fact that 
$|\Rm|(\overline{x}_\infty, \overline{t}_\infty) \: = \: 1$.

Suppose now that there is a subsequence so that
the injectivity radii of the scaled metrics 
at $(\overline{x}_k, \overline{t}_k)$ tend to zero.
Parabolically rescale $(M_k,(\overline{x}_k, \overline{t}_k), g_k(\cdot))$ 
further so that the injectivity radius becomes one. After passing to a
subsequence we will have convergence to a flat Ricci flow solution
$(-\infty, 0] \times L$. The complete 
flat manifold $L$ can be described as the total space of a flat 
orthogonal
$\R^m$-bundle over a flat compact manifold $C$. After separating
variables, the
fundamental solution $u_\infty$ on $L$ will be Gaussian in the fiber directions and will
decay exponentially fast to a constant in the base directions, i.e.
$u(x, \tau) \: \sim \: (4\pi \tau)^{-m/2} \: e^{- \frac{|x|^2}{4\tau}} \:
\frac{1}{\vol(C)}$, where $|x|$ is the fiber norm. With this for $u_\infty$,
one finds that
$v_\infty \: = \: (m - n) \: \left( 1 + \frac12 \ln(4 \pi \tau) \right) \: u_\infty$.
With $B_\tau$ the ball around a basepoint 
of radius $\sqrt{\tau}$, the integral of $u$ over
$B_\tau$ has a positive limit as $\tau \rightarrow \infty$, and so 
$\lim_{\tau \rightarrow \infty} \int_{B_\tau} v_\infty \: dV_\infty \: = \: - \: \infty$.
Then there are times $\widetilde{t}_k
\in [\overline{t}_k \: - \: \frac12 \: \alpha Q_k^{-1}, \overline{t}_k]$ so that
$\lim_{k \rightarrow \infty}
\int_{B_k} v_k \: dV_k \: = \: - \infty$, which is a contradiction.
\end{proof}

\section{Theorem I.10.1. Proof of the pseudolocality theorem} \label{pfI.10.1}

Continuing with the proof of Theorem \ref{thmI.10.1}, we now use
Lemma \ref{Claim 3} to get a contradiction to a log Sobolev inequality. 
For simplicity of notation, we drop the subscript $k$ and deal with a 
particular $(M_k, (\overline{x}_k, \overline{t}_k), g_k(\cdot))$ for
$k$ large.
Define a smooth function $\phi$ on $\R$ which is
one on $(- \infty, 1]$, decreasing on $[1,2]$ and zero on $[2, \infty]$,
with $\phi^{\prime \prime} \ge  - 10 \phi^\prime$ and
$(\phi^\prime)^2 \le 10 \phi$. To construct $\phi$ we can 
take the function which is $1$ on $(-\infty, 1]$,
$1 \: - \: 2(x-1)^2$ on $[1, 3/2]$, $2(x-2)^2$ on $[3/2, 2]$ and
$0$ on $[2, \infty)$, and
smooth it slightly.

Put $\widetilde{d}(y,t) = d(y, t) + 200 n \sqrt{t}$. We claim that if
$10A \epsilon \le \widetilde{d}(y,t) \le 20A\epsilon$ then
$d_t(y,t) \: - \: \triangle d(y,t) \: + \: \frac{100n}{\sqrt{t}} \: \ge \: 0$.
To see this, recalling that $t \in [0, \epsilon^2]$, 
if $10A \epsilon \le \widetilde{d}(y,t) \le 20A\epsilon$ and $A$ is
sufficiently large then
$9A\epsilon \: \le \: d(y, t) \: \le \: 
21A\epsilon$. We apply
Lemma \ref{I.8.3a} with the parameter $r_0$ of
Lemma \ref{I.8.3a} equal to $\sqrt{t}$. As
$r_0 \le \epsilon$, we have $y \notin B(x_0, r_0)$. 
From (\ref{largebound}), on $B(x_0, r_0)$ we have
$|\Rm|(\cdot,t) \: \le \: \alpha \: t^{-1} \: + \: 2 \: \epsilon^{-2}$.
Then from Lemma \ref{I.8.3a}, at $(y, t)$ we have
\begin{align}
d_t \: - \: \triangle d \: & \ge \: - \: (n-1) \:
\left( \frac23 (\alpha \: t^{-1} \: + \: 2 \: \epsilon^{-2} ) t^{1/2}
\: + \: t^{-1/2} \right) \\
& = \: - \: (n-1) \:
\left( 1 + \frac23 \alpha + \frac43 \: 
\epsilon^{-2} t \right) \: t^{-1/2}. \notag
\end{align}
It follows that
$d_t \: - \: \triangle d \: + \: \frac{100n}{\sqrt{t}} \: \ge \: 0$.

Now put $h(y,t) = \phi \left( \frac{\widetilde{d}(y,t)}{10A\epsilon} \right)$.
Then $\square h = \frac{1}{10A\epsilon} \left( 
d_t \: - \: \triangle d \: + \: \frac{100n}{\sqrt{t}} \right) \phi^\prime \: - \:
\frac{1}{(10A\epsilon)^2} \phi^{\prime \prime}$, where the arguments of
$\phi^\prime$ and $\phi^{\prime \prime}$ are 
$\frac{\widetilde{d}(y,t)}{10A\epsilon}$. Where
$\phi^\prime \neq 0$, 
we have $d_t \: - \: \triangle d \: + \: \frac{100n}{\sqrt{t}} \: \ge \: 0$.
The fundamental solution $u(x,t) = (4\pi(\overline{t}-t))^{- \: \frac{n}{2}} \:
e^{-f(x,t)}$ of $\square^*$ is positive for $t \in [0, \overline{t})$ 
and we have
$\int_M u \: dV \: \le \: 1$ for all $t$.
(Recall that we are not assuming stochastic completeness.)
Then
\begin{align}
\left( \int_M hu \: dV \right)_t \: & = \: \int_M ((\square h) u \: - \: h
\square^* u) \: dV \:  = \: \int_M (\square h) u \: dV \: \le \: - \:
\frac{1}{(10A\epsilon)^2} \int_M \phi^{\prime \prime} u \: dV \\
& \le \frac{10}{(10A\epsilon)^2} \int_M \phi u \: dV \:
 \le \:
\frac{10}{(10A\epsilon)^2} \int_M u \: dV \: \le \: \frac{10}{(10A\epsilon)^2}. \notag
\end{align}
Hence 
\begin{equation} \label{ut}
\int_M hu \: dV \Big|_{t=0} \ge \int_M hu \: dV \Big|_{t=\overline{t}} - 
\frac{\overline{t}}{(A\epsilon)^2} \ge 1 - A^{-2}.
\end{equation}

Similarly, using Proposition \ref{localized} and Corollary \ref{corI.9.3},
\begin{align}
\left(- \int_M hv \: dV \right)_t \: & = \: - \: \int_M 
((\square h) v - h \square^* v) \: dV \: \le \: - \int_M (\square h) v 
\: dV \\
& \le \:
\frac{1}{(10A\epsilon)^2} \: \int_M \phi^{\prime \prime} v \: dV \: \le \: - \:
\frac{10}{(10A\epsilon)^2} \: \int_M \phi v \: dV \: = \: - \:
\frac{10}{(10A\epsilon)^2} \: \int_M h v \: dV. \notag
\end{align}
Consider the time $\widetilde{t}$ of Lemma \ref{Claim 3}. As 
$(\overline{x}, \overline{t}) \in M_\alpha$, 
$\widetilde{t} \in
[\overline{t}/2, \overline{t}]$. Then $\sqrt{\overline{t} - 
\widetilde{t}} \: \le \: 2^{-1/2} \: \epsilon$ and so for large $A$,
$h$ will be one on the ball $B$
at time $\widetilde{t}$ of radius $\sqrt{\overline{t} - 
\widetilde{t}}$ centered at $\overline{x}$, using
(\ref{32.4}). Then at time $\widetilde{t}$,
\begin{equation}
- \int_M hv \: dV \: \ge \: - \int_B v \: dV \: \ge \: \beta.
\end{equation}
Thus
\begin{equation}
- \int_M hv \: dV \Big|_{t=0} \: \ge \: \beta \: e^{- 
\frac{\widetilde{t}}{(A\epsilon)^2}} \: \ge \:
\beta \: e^{- \frac{\overline{t}}{(A\epsilon)^2}} \: \ge \: \beta 
\left( 1 - \frac{\overline{t}}{(A \epsilon)^2} \right)
\: \ge \: \beta (1 - A^{-2}).  
\end{equation}

Working at time $0$, put $\widetilde{u} = hu$ and
$\widetilde{f} = f - \log h$.
In what follows we implicitly integrate over $\supp(h)$.
We have
\begin{equation}
\beta (1 - A^{-2}) \le - \int_M hv \: dV \: = \:
\int_M [(-2\triangle f + |\nabla f|^2 - R) \overline{t}
- f + n]hu \: dV.
\end{equation}

We claim that
\begin{equation}
\int_M \left( - 2 \triangle f \: + \: |\nabla f|^2 \right) \: h e^{-f} \: 
dV \: = \: \int_M \left( - \: |\nabla \widetilde{f}|^2 \: + \:
\frac{|\nabla h|^2}{h^2} \right) \: he^{-f} \: dV.
\end{equation}
This follows from
\begin{align}
\int_M \left( - 2 \triangle f \: + \: |\nabla f|^2 \right) \: h e^{-f} \: 
dV \: & = \: 
\int_M \left( 2 \langle \nabla f, \nabla(he^{-f}) \rangle \: + \: 
|\nabla f|^2 \: h e^{-f} \right) \: 
dV \\
& = \: \int_M \left( 2 \langle \nabla f, \frac{\nabla h}{h} -
\nabla f \rangle \: + \: 
|\nabla f|^2  \right) \: h e^{-f} \: 
dV \notag \\
& = \: \int_M \langle \nabla f, 2 \frac{\nabla h}{h} -
\nabla f \rangle \: h e^{-f} \: dV \notag \\
& = \:
\int_M \langle \nabla \widetilde{f} + \frac{\nabla h}{h}, \frac{\nabla h}{h} -
\nabla \widetilde{f} \rangle \: h e^{-f} \: dV \notag \\
& = \:
\int_M \left( - \: |\nabla \widetilde{f}|^2 \: + \:
\frac{|\nabla h|^2}{h^2} \right) \: he^{-f} \: dV. \notag 
\end{align}
Then
\begin{align}
& \int_M [(-2\triangle f + |\nabla f|^2 - R) \overline{t}
- f + n]hu \: dV \: = \\
& \int_M [- \overline{t} |\nabla \widetilde{f}|^2 - \widetilde{f} + n]
\widetilde{u} \: dV \: + \:
\int_M [\overline{t}(|\nabla h|^2/h - Rh) - h \log h ] u \: dV. \notag
\end{align}

Next, $\frac{|\nabla h|^2}{h} \: \le \: \frac{10}{(10A\epsilon)^2}$ and
$- \: Rh \: \le \: 1$ (from the assumed lower bound on $R$ at time zero). Then
\begin{equation}
\int_M \overline{t} \: \left( \frac{|\nabla h|^2}{h} \: - \: R h \right)
\: u \: dV \: \le \: \epsilon^2 \: \left( \frac{10}{(10A\epsilon)^2} \: + \:
1 \right) \: \le \: A^{-2} \: + \: \epsilon^2.
\end{equation}
Also,
\begin{align}
- \int_M uh \log{h} \: dV \: & = \: -
\int_{B(x_0, 20 A \epsilon) - B(x_0, 10 A \epsilon)} uh \log{h} \: dV \: 
\le \: \int_{M - B(x_0, 10 A \epsilon)} u \: dV \\
&\le \: 1 \: - \: \int_{B(x_0, 10 A \epsilon)} u \: dV. \notag
\end{align}
Putting $\overline{h}(y) = \phi \left( \frac{d(y)}{5A\epsilon} \right)$, a result
similar to (\ref{ut}) shows that 
\begin{equation}
\int_{B(x_0, 10 A \epsilon)} u \: dV \: \ge \:
\int_{M} \overline{h} u \: dV \: \ge \:
 1 \: - \: c A^{-2}
 \end{equation}
for an appropriate constant $c$. 
Putting this together gives
\begin{equation} \label{somme}
\beta(1 - A^{-2}) \: \le \: \int_M \left( - \overline{t} |\nabla \widetilde{f}|^2 - \widetilde{f} + n
\right)
\widetilde{u} \: dV \: + \: (1+c) A^{-2} + \epsilon^2.
\end{equation}

Put $\widehat{g} = \frac{1}{2\overline{t}} g$,
$\widehat{u} = (2 \overline{t})^{\frac{n}{2}} \widetilde{u}$ and define
$\widehat{f}$ by $\widehat{u} \: = \: (2\pi)^{-\frac{n}{2}} \: e^{-\widehat{f}}$.
From (\ref{ut}) and (\ref{somme}), if we restore the subscript $k$ then
$\lim_{k \rightarrow \infty} \int_M \widehat{u}_k \: d\widehat{V}_k \: = \: 1$ and for large $k$,
\begin{equation} 
\frac12 \: \beta \: \le \: \int_{M_k} 
\left( - \: \frac12 \: |\nabla \widehat{f}_k|^2 - \widehat{f}_k + n
\right)
\widehat{u}_k \: d\widehat{V}_k.
\end{equation}
If we normalize $\widehat{u}_k$ by putting
$U_k \: = \: \frac{\widehat{u}_k}{\int_{M_k} \widehat{u}_k \: d\widehat{V}_k}$, 
and
define $F_k$ by $U_k \: = \: (2\pi)^{- \: \frac{n}{2}} \: e^{- \: F_k}$, then
for large $k$, we also have
\begin{equation} \label{contra1}
\frac12 \: \beta \: \le \: \int_{M_k} 
\left( - \frac12 |\nabla {F}_k|^2 - {F}_k + n
\right)
{U}_k \: d\widehat{V}_k.
\end{equation}

On the other hand, the
logarithmic Sobolev inequality for $\R^n$ \cite[I.(8)]{Beckner} says that
\begin{equation} \label{Beckner}
\int_{\R^n} \left( - \: \frac12 \: |\nabla F|^2 \: - \: F \: + \: n \right) U \: dV
\: \le \: 0,
\end{equation}
provided that the compactly-supported function
$U \: = \: (2 \pi)^{-n/2} \: e^{-F}$ satisfies
$\int_{\R^n} U \: dV \: = \: 1$.
As was mentioned to us by Peter Topping, one can get a sharper inequality
by applying (\ref{Beckner}) to the rescaled function
$U_c(x) \: = \: c^n \: U(cx)$ and optimizing with respect to $c$.
The result is
\begin{equation}
\int_{\R^n} |\nabla F|^2 \: U \: dV \: \ge \: n \: e^{1 \: - \: \frac{2}{n} \:
\int_{\R^n} F \: U \: dV}.
\end{equation}
Given this inequality on $\R^n$,
one can use a symmetrization argument to prove the same
inequality for a compactly-supported function on any complete Riemannian
manifold, provided that the Euclidean isoperimetric
inequality holds for domains in the support of $U$.
See, for example, \cite[Proposition 4.1]{Ni3} which gives 
the symmetrization argument for
(\ref{Beckner}), attributing it to Perelman.
Again using the inequality for $\R^n$, if instead we have
$\vol(\partial \Omega)^n
\ge (1-\delta_k) \: c_n \: \vol(\Omega)^{n-1}$ for domains
$\Omega \subset \supp(U_k)$ then the symmetrization argument gives
\begin{equation} \label{contra2}
\int_{M_k} |\nabla {F}_k|^2 \: {U}_k \: d\widehat{V}_k \: \ge \:
(1 - \delta_k)^{\frac{2}{n}} \: 
n \: e^{1 \: - \: \frac{2}{n} \:
\int_{M_k} {F}_k \: U_k \: d\widehat{V}_k}.
\end{equation}

Equations (\ref{contra1}) and (\ref{contra2}) imply that
\begin{equation}
\frac{n}{2} \: \left( (1 - \delta_k)^{\frac{2}{n}} \: e^{1 - \frac{2}{n} \int_{M_k}
F_k \: U_k \: d\widehat{V}_k} \: - \: 1 \: - \:
\left( 1 - \frac{2}{n} \int_{M_k}
F_k \: U_k \: d\widehat{V}_k \right) \right) \: \le \: - \: \frac{\beta}{2}.
\end{equation}
However,
\begin{equation}
\lim_{k \rightarrow \infty} \inf_{x \in \R} 
\left(
(1 - \delta_k)^{\frac{2}{n}} \: e^{x} \: - \: 1 \: - \: x \right) \: = \: 0.
\end{equation}
This is a contradiction.

\section{I.10.2. The volumes of future balls}

The next result gives a lower bound on the volumes of future
balls.

\begin{corollary} (cf. Corollary I.10.2) \label{corI.10.2}
Under the assumptions of Theorem \ref{thmI.10.1}, for
$0 < t \le (\epsilon r_0)^2$ we have
$\vol(B_t(x, \sqrt{t})) \: \ge \: c t^{\frac{n}{2}}$ for $x \in B_0(x_0, \epsilon r_0)$,
where $c = c(n)$ is a universal constant. 
\end{corollary}
\begin{proof} (Sketch)
If the corollary were not true then taking a sequence of counterexamples,
we can center ourselves around the
collapsing balls $B(x, \sqrt{t})$ to obtain functions $f$ as in Section \ref{pfI.10.1}.
As in the proof of Theorem \ref{noncollthm}, the volume condition
along with the fact that $\int_M (2 \pi)^{-n/2} \: e^{-f} \: dV \rightarrow 1$
means that $f \rightarrow - \infty$, which implies that
$\int_M \left( - \: \frac12 \: |\nabla f|^2 \: - \: f \: + \: n \right) u \: dV 
\rightarrow \infty$. This contradicts the logarithmic Sobolev inequality.
\end{proof}

\section{I.10.4. $\kappa$-noncollapsing at future times}

The next result gives $\kappa$-noncollapsing at future times.

\begin{corollary} (cf. Corollary I.10.4)
There are $\delta,\epsilon > 0$ such that for any $A > 0$ there exists
$\kappa = \kappa(A) > 0$ with the following property.  
Suppose that we have a Ricci flow solution
$g(\cdot)$ defined for $t \in [0, (\epsilon r_0)^2]$ which
has bounded $|\Rm|$ and complete time slices. Suppose that
for any $x \in B(x_0, r_0)$ and $\Omega \subset B(x_0, r_0)$,
we have $R(x,0) \ge - r_0^{-2}$ and $\vol(\partial \Omega)^n
\ge (1-\delta) \: c_n \: \vol(\Omega)^{n-1}$, where $c_n$ is the
Euclidean isoperimetric constant. If $(x,t)$ satisfies
$A^{-1} (\epsilon r_0)^2 \le t \le (\epsilon r_0)^2$ and
$\dist_t(x, x_0) \le A r_0$ then $g(\cdot)$ is not
$\kappa$-collapsed at $(x,t)$ on scales less than $\sqrt{t}$.
\end{corollary}
\begin{proof}
Using Theorem \ref{thmI.10.1} and
Corollary \ref{corI.10.2}, we can apply
Theorem \ref{8.2thm} starting at time $A^{-1} (\epsilon r_0)^2$.
\end{proof}

\section{I.10.5. Diffeomorphism finiteness}
\label{I.10.5}

In this section we prove the diffeomorphism finiteness of Riemannian
manifolds with local isoperimetric inequalities, a lower bound on 
scalar curvature and an upper bound on volume.

\begin{theorem}
Given $n \in \Z^+$, there is a $\delta > 0$ with the following property.  For any
$r_0, V > 0$, there are finitely many diffeomorphism types
of compact $n$-dimensional Riemannian manifolds $(M, g_0)$ satisfying \\
1. $R \: \ge \: - r_0^{-2}$. \\
2. $\vol(M, g_0) \: \le \: V$. \\
3. Any domain $\Omega \subset M$ contained in a metric
$r_0$-ball satisfies $\vol(\partial \Omega)^n
\ge (1-\delta) \: c_n \: \vol(\Omega)^{n-1}$, where $c_n$ is the
Euclidean isoperimetric constant.
\end{theorem}
\begin{proof}
Choose $\alpha > 0$. Let $\delta$ and $\epsilon$ be the parameters
of Theorem \ref{thmI.10.1}.
Consider Ricci flow $g(\cdot)$ starting from $(M, g_0)$. 
Let $T > 0$ be the maximal number so that a smooth flow exists
for $t \in [0, T)$.  If $T < \infty$ then
$\lim_{t \rightarrow T^-} \sup_{x \in M} |\Rm(x, t)| \: = \: \infty$. It
follows from Theorem \ref{thmI.10.1} that $T > (\epsilon r_0)^2$. 
Put $\widehat{g} = g((\epsilon r_0)^2)$. Theorem \ref{thmI.10.1} gives a 
uniform double-sided sectional curvature bound on $(M, \widehat{g})$. 
Corollary \ref{corI.10.2} gives a uniform lower bound on the volumes
of $(\epsilon r_0)$-balls in $(M, \widehat{g})$. Let $\{x_i\}_{i=1}^N$ be a
maximal $(2\epsilon r_0)$-separated net in $(M, \widehat{g})$. 

From the lower bound $R \: \ge \: - r_0^{-2}$ on $(M, g_0)$ and the
maximum principle, we have 
$R(x,t)  \: \ge \: - r_0^{-2}$ for $t \in [0, (\epsilon r_0)^2]$. 
Then the Ricci flow equation gives
a uniform upper bound on $\vol(M, \widehat{g})$. This implies a uniform
upper bound on $N$ or, equivalently, a uniform upper bound on 
$\diam(M, \widehat{g})$. The theorem now follows from the diffeomorphism
finiteness of $n$-dimensional Riemannian
manifolds with double-sided sectional curvature bounds,
upper bounds on diameter and lower bounds on volume.
\end{proof}

\section{I.11.1. $\kappa$-solutions}
\label{secI.11.1}

\begin{definition}
Given $\kappa > 0$, a {\em $\kappa$-solution} is a Ricci flow
solution $(M,g(\cdot))$ that is defined on a time interval of the 
form $(-\infty,C)$
(or $(-\infty,C]$)
such that 

$\bullet$ The curvature $|\Rm|$ is bounded on each compact time interval
$[t_1,t_2]\subset (-\infty,C)$
(or $(-\infty,C]$), and each time slice $(M,g(t))$ is complete.

$\bullet$ The curvature operator is nonnegative and the scalar curvature
is everywhere positive.

$\bullet$ The Ricci flow is $\kappa$-noncollapsed at all scales.
\end{definition}

By abuse of terminology, we may sometimes write that 
``$(M, g(\cdot))$ is a $\kappa$-solution'' if it is a $\kappa$-solution
for some $\kappa > 0$.

From Appendix \ref{appharnack},
$R_t  \:  \ge \: 0$ for an ancient solution.
This implies the essential equivalence of the notions of
$\kappa$-noncollapsing in
Definitions \ref{noncollapsed1} and \ref{noncollapsed2} , when
restricted to ancient solutions.
Namely, if a solution is $\kappa$-collapsed in the sense of 
Definition \ref{noncollapsed2} then it is
automatically $\kappa$-collapsed in the sense of 
Definition \ref{noncollapsed1}.  Conversely,
if a time-$t_0$ slice of an ancient solution is collapsed in the sense
of Definition \ref{noncollapsed1} 
then the fact that $R_t  \:  \ge \: 0$, together with bounds on
distance distortion, implies that it is collapsed in the sense of
Definition \ref{noncollapsed2} (possibly for a different value of $\kappa$).

The relevance of $\kappa$-solutions is that a blowup limit of a
finite-time singularity on a compact manifold will be a $\kappa$-solution.

For examples of $\kappa$-solutions, if $n \ge 3$ then 
there is a $\kappa$-solution
on the cylinder $\R \times S^{n-1}(r)$,
where the radius satisfies $r^2(t) \: = \: r_0^2 \: - \:
2(n-2)t$. There is also a $\kappa$-solution on the
$\Z_2$-quotient $\R \times_{\Z_2} S^{n-1}(r)$, where
the generator of $\Z_2$ acts by reflection on $\R$ and by
the antipodal map on $S^{n-1}$.
On the other hand,  the quotient solution on
$S^1 \times S^{n-1}(r)$ is not $\kappa$-noncollapsed for any 
$\kappa > 0$,
as can be seen by looking at large negative time.

Bryant's gradient steady soliton is a 
three-dimensional $\kappa$-solution given by
$g(t) \: = \: \phi_t^* g_0$, where $g_0 \: = \: dr^2 \: + \: \mu(r) \: 
d\Theta^2$ is a certain rotationally symmetric metric on $\R^3$. It has
sectional curvatures that go like $r^{-1}$, and $\mu(r) \sim r$.
The gradient function $f$ satisfies $R_{ij} \: + \: \nabla_i \nabla_j f
\: = \: 0$, with $f(r) \: \sim \: -2r$. Then for $r$ and $r \: - \: 2t$
large, $\phi_t(r, \Theta) \: \sim \: (r-2t, \Theta)$. In particular,
if $R_0 \in C^\infty(\R^3)$ is the scalar curvature function of
$g_0$ then $R(t, r, \Theta) \sim R_0(r-2t, \Theta)$. 

To check the conclusion of Corollary \ref{corI.11.8}
in this case,
given a point $(r_0, \Theta) \in \R^3$ at time $0$, the scalar
curvature goes like $r_0^{-1}$. Multiplying the soliton metric by
$r_0^{-1}$ and sending $t \rightarrow r_0 t$ gives the asymptotic metric
\begin{equation}
d(r/\sqrt{r_0})^2 \:
+ \: \frac{r-2r_0t}{r_0} \: d\Theta^2.
\end{equation}
Putting $u \: = \: (r - r_0)/\sqrt{r_0}$,
the rescaled metric is approximately 
\begin{equation}
du^2 \: + \: 
\left( 1 \: + \: \frac{u}{\sqrt{r_0}} \: - \: {2t} \right)
\: d\Theta^2.
\end{equation}
Given $\epsilon \: > 0$, this will be $\epsilon$-biLipschitz close to 
the evolving cylinder
$du^2 \: + \: 
\left( 1 \: - \: {2t} \right)
\: d\Theta^2$ provided that $|u| \: \le \: \epsilon \sqrt{r_0}$,
i.e. $|r-r_0| \: \le \: \epsilon r_0$. To have an $\epsilon$-neck,
we want this to hold whenever $|r \: - \: r_0|^2 \: \le \: 
(\epsilon r_0^{-1})^{-1}$. This will be the case if $r_0 \: \ge \:
\epsilon^{-3}$. Thus $M_\epsilon$ is approximately
\begin{equation}
\{ (r, \Theta) \in \R^3 \: : \: r \: \le \: \epsilon^{-3} \}  
\end{equation}
and $Q \: = \: R(x_0,0) \sim \epsilon^3$. 
Then $\diam(M_\epsilon) \sim \epsilon^{-3}$
and at the origin $0 \in M_\epsilon$,
$R(0, 0) \sim \epsilon^0$. It follows that for the value of
$\kappa$ corresponding to this solution, $C(\epsilon, \kappa)$ must
grow at least as fast as $\epsilon^{-3}$ as $\epsilon \rightarrow 0$.

\section{I.11.2. Asymptotic solitons} \label{I.11.2}

This section shows that every $\kappa$-solution has a gradient shrinking
soliton buried inside of it, in an asymptotic sense as $t \rightarrow -\infty$.
Such a soliton will be called an {\em asymptotic soliton}.

Heuristically, the existence of an asymptotic soliton is
a consequence of the compactness results and 
the monotonicity of the reduced volume.  Taking an appropriate 
sequence of
spacetime points going backward in time, one constructs a
limiting rescaled solution. As the limit reduced volume
is constant in time, the monotonicity formula implies that
this limit solution is a gradient shrinking soliton. This is
the basic idea but
the rigorous argument is a bit more subtle.

Pick an arbitrary point $(p, t_0)$ in the $\kappa$-solution $(M, g(\cdot))$. Define
the reduced volume $\widetilde{V}(\tau)$ and the reduced length
$l(q, \tau)$ as in Section \ref{I.7}, by means of curves starting from
$(p, t_0)$, with $\tau = t_0 - t$. From 
Section \ref{I.(7.15)}, for each $\tau> 0$ there is some
$q(\tau)\in M$ such that $l(q(\tau), \tau) \le \frac{n}{2}$. (Note that
$l \ge 0$ from the curvature assumption.)

\begin{proposition} (cf. Proposition I.11.2) \label{asympsoliton}
There is a sequence $\overline{\tau}_i \rightarrow \infty$ so that if we
consider the solution $g(\cdot)$ on the time interval
$[t_0 - \overline{\tau}_i, t_0 - \frac12 \overline{\tau}_i]$  and 
parabolically rescale it
at the point $(q(\overline{\tau}_i), t_0 - 
\overline{\tau}_i)$ by the factor $\overline{\tau_i}^{-1}$ then as
$i \rightarrow \infty$, the
rescaled solutions
converge to a nonflat gradient shrinking soliton (restricted to $[-1, - \frac12]$).
\end{proposition}
\begin{proof}

Equation (\ref{I.(7.16)}) implies that
$|\nabla l^{1/2}|^2 \: \le \: \frac{C}{4\tau}$, and so
\begin{equation}
|l^{1/2}(q,\tau) \: - \: l^{1/2}(q(\tau),\tau)| \: 
\le \: \sqrt{\frac{C}{4\tau}}
\: \dist_{t_0-\tau}(q, q(\tau)).
\end{equation}
We apply this estimate initially at some fixed time $\tau = \overline{\tau}$, to
obtain
\begin{equation}
l(q, \overline{\tau}) \: \le \: \left( \sqrt{\frac{C}{4\overline{\tau}}}
\: \dist_{t_0- \overline{\tau}}(q, q(\overline{\tau})) \: + \: 
\sqrt{\frac{n}{2}} \right)^2.
\end{equation}
From (\ref{touse1}), (\ref{ttouse1}) and (\ref{I.(7.16)}), 
\begin{equation}
\partial_\tau l \: = \: \frac{R}{2} \: - \: \frac{|\nabla l|^2}{2}
\: - \: \frac{l}{2\tau} \: \ge \: - \: \frac{(1+C)l}{2\tau}.
\end{equation}
This implies that for $\tau \in \left[ \frac12 \overline{\tau},
\overline{\tau} \right]$,
\begin{equation} \label{upperl}
l(q, {\tau}) \: \le \: \left( \frac{\overline{\tau}}{\tau}
\right)^{\frac{1+C}{2}} \:
\left( \sqrt{\frac{C}{4\overline{\tau}}}
\: \dist_{t_0- \overline{\tau}}(q, q(\overline{\tau})) \: + \: 
\sqrt{\frac{n}{2}} \right)^2.
\end{equation}
Also from (\ref{I.(7.16)}), we have
$\tau R \: \le \: Cl$.
Then we can plug in the previous bound on 
$l$ to get an upper bound on $\tau R$ for
$\tau \in \left[ \frac12 \overline{\tau}, \overline{\tau} \right]$.
The upshot is that for any $\epsilon > 0$, one can find $\delta > 0$
so that both $l(q,\tau)$ and $\tau R(q, t_0-\tau)$ do not exceed
$\delta^{-1}$ whenever $\frac12 \overline{\tau} \le \tau \le \overline{\tau}$
and $\dist^2_{t_0-\overline{\tau}} (q, q(\overline{\tau})) \: \le \:
\epsilon^{-1} \overline{\tau}$.

Varying $\overline{\tau}$, as the rescaled solutions
(with basepoints at $(q(\overline{\tau}), t_0 - \overline{\tau})$) 
are uniformly noncollapsing and have
uniform curvature bounds on balls, Appendix \ref{subsequence} 
implies that
we can take
a sequence $\overline{\tau}_i \rightarrow \infty$
to get a pointed limit 
$(\overline{M}, \overline{q}, \overline{g}(\cdot))$ that is
a complete Ricci flow solution (in the backward parameter $\tau$) for
$\frac12 \: < \: \tau \: < \: 1$. 
We may assume that we have locally Lipschitz convergence of $l$ to
a limit function $\overline{l}$.

We define the reduced volume $\tilde{V}(\tau)$ for the limit solution
using the limit function $\overline{l}$. 
We claim that for any 
$\tau \in \left( \frac12, 1 \right)$, if we put $\tau_i = \tau \overline{\tau}_i$
then the number $\tilde{V}(\tau)$ for
the limit solution is the limit of numbers $\tilde{V}(\tau_i)$ for the
original solution. One wishes to apply dominated convergence to the
integrals $\int_M e^{-l(q, \tau_i)} \: 
\tau_i^{-n/2} \: \dvol(q,t_0 - \tau_i)$. 
(Note that $\tau_i^{-n/2} \: \dvol(q,t_0 - \tau_i) \: = \:
\tau^{-n/2} \: \overline{\tau}_i^{-n/2} \: \dvol(q,t_0 - \tau_i)$ and
$\overline{\tau}_i^{-n/2} \: \dvol(q,t_0 - \tau_i)$ is the volume form for the
rescaled metric $\overline{\tau}_i^{-1} \: g(t_0 - \tau_i)$.)
However, to do so one needs uniform lower bounds on
$l(q, \tau^\prime)$ for the original solution in terms of 
$d_{t_0-\tau^\prime}(q, q(\tau^\prime))$, for $\tau^\prime \in (-\infty, 0)$. 
By an argument of Perelman,
written in detail in \cite{Ye1}, one does indeed have a lower bound of the
form
\begin{equation} \label{Perelmanbound}
l(q, \tau^\prime) \: \ge \: - \: l(q(\tau^\prime)) \: - \: 1 \: + \: C(n) \:
\frac{d_{t_0-\tau^\prime}(q, q(\tau^\prime))^2}{\tau^\prime}.
\end{equation}
The nonnegative curvature gives polynomial volume growth for 
distance balls,
so using (\ref{Perelmanbound}) 
one can apply dominated convergence to the integrals
$\int_M e^{-l(q, \tau_i)} \: 
\tau_i^{-n/2} \: \dvol(q,t_0 - \tau_i)$.
Thus $\lim_{i \rightarrow \infty} \tilde{V}(\tau_i) \: = \: 
\tilde{V}(\tau)$.

As (\ref{upperl}) gives
a uniform upper bound on $l$ on an appropriate 
ball around $q(\tau_i)$, and there is a lower 
volume bound on the ball, it follows that as $i \rightarrow \infty$,
$\tilde{V}(\tau_i)$ is uniformly bounded away from zero. From this
argument and the monotonicity of $\tilde{V}$, 
$\tilde{V}(\tau)$ is a positive constant $c$ as a function of
$\tau$, namely the limit of the reduced volume of the original
solution as real time goes to $- \infty$. As the original solution is
nonflat, the constant $c$ is strictly less than the limit of the
reduced volume of the original solution
as real time goes to zero, which is $(4 \pi)^{\frac{n}{2}}$. 

Next, we will apply (\ref{illustrate}) and (\ref{illustrate2}).
As (\ref{illustrate}) holds distributionally
for each rescaled solution, it follows that
it holds distributionally for $\overline{l}$. In particular, the
nonpositivity implies that the left-hand side of (\ref{illustrate}), when
computed for the limit solution, is actually a nonpositive measure.
If the left-hand side of
(\ref{illustrate}) (for the limit solution)
were not strictly zero then using (\ref{pphieqn})  we would conclude
that $\frac{d\tilde{V}}{d\tau}$ 
is somewhere negative, which is a contradiction.
(We use the fact that (\ref{Perelmanbound}) passes to the limit
to give a similar lower bound on $\overline{l}$.)
Thus we must have equality in (\ref{illustrate}) for the limit solution.
This implies equality in (\ref{touse2}),
which implies equality in  (\ref{illustrate2}).
Writing
 (\ref{illustrate2}) as
\begin{equation}
(4 \triangle \: - \: R) \: e^{-\frac{l}{2}} \: = \:
\frac{l-n}{\tau} \: e^{-\frac{l}{2}},
\end{equation}
elliptic theory gives smoothness of $\overline{l}$. 

In I.11.2 it is said that equality in  (\ref{illustrate2}) implies equality
in (\ref{bbound2}), which implies that one has a gradient shrinking soliton.
There is a problem with this argument, as the use of (\ref{bbound2})
implicitly
assumes that the solution is defined for all $\tau \ge 0$, which we do not
know. (The function $\overline{l}$ is only defined by a limiting
procedure, and not in terms of ${\mathcal L}$-geodesics on some
Ricci flow solution.) However, one can instead use Proposition \ref{localized}, with 
$f = \overline{l}$. 
Equality in (\ref{illustrate2}) implies that $v = 0$, so (\ref{I.(9.1)})
directly gives the gradient shrinking soliton equation.
(The problem with the argument using (\ref{bbound2}), and its resolution using
(\ref{I.(9.1)}), were pointed out by the UCSB group.)

If the gradient shrinking soliton $\overline{g}(\cdot)$ is
flat then, as it will be $\kappa$-noncollapsed
at all scales, it must be $\R^n$. From the soliton equation,
$\partial_i \partial_j \overline{l} \: = \:
\frac{\overline{g}_{ij}}{2\tau}$ and 
$\triangle \overline{l} \: = \: \frac{n}{2\tau}$. Putting this
into the equality (\ref{illustrate2}) gives
$|\nabla \overline{l}|^2 \: = \: \frac{\overline{l}}{\tau}$. 
It follows that the level sets of $\overline{l}$ are distance spheres.
Then (\ref{illustrate}) implies
that with an appropriate choice of
origin,  
$\overline{l} \: = \: \frac{|x|^2}{4\tau}$.  The reduced volume
$\tilde{V}(\tau)$ for the limit solution is now computed to be
$(4 \pi)^{\frac{n}{2}}$, which is a contradiction.  Therefore the
gradient shrinking soliton is not flat.
\end{proof}

We remark that the gradient soliton constructed here does
not, {\it a priori}, have bounded curvature on 
compact time intervals,
i.e.
it may not be a $\kappa$-solution.  In the $2$ and
$3$-dimensional cases
one can prove this using additional
reasoning.  See Section \ref{altpf11.3} where
it is shown that $2$-dimensional $\kappa$-solutions are
round spheres, and Section \ref{pf11.7} where the $3$-dimensional
case is discussed.

\section{I.11.3. Two dimensional $\kappa$-solutions} \label{I.11.3}

The next result is a classification of two-dimensional
$\kappa$-solutions.  It is important when doing dimensional
reduction.

\begin{corollary} (cf. Corollary I.11.3) \label{2Dsoliton}
The only oriented two-dimensional $\kappa$-solution is
the shrinking round $2$-sphere.
\end{corollary}
\begin{proof}
First, the
only nonflat oriented 
nonnegatively curved gradient shrinking 2-D soliton is the
round $S^2$. The reference \cite{Hamsurface} given in I.11.3 for this
fact does not actually cover it, as the reference
only deals with compact solitons. A proof using 
Proposition \ref{asympsoliton} to rule
out the noncompact case appears in \cite{Ye2}.

Given this, the limit solution in Proposition \ref{asympsoliton}
is a shrinking round $2$-sphere.  Thus the rescalings 
$\overline{\tau}_i^{-1} g(t_0 - \overline{\tau}_i)$ converge to
a round $2$-sphere as $i \rightarrow \infty$.  However, by \cite{Hamsurface}
the Ricci flow makes an
almost-round $2$-sphere become more round. Thus any given time slice
of the original $\kappa$-solution must be a round $2$-sphere. 
\end{proof}

\begin{remark}
One can employ a somewhat different
line of reasoning to prove Corollary \ref{2Dsoliton}; see Section
\ref{altpf11.3}.
\end{remark}

\section{I.11.4. Asymptotic scalar curvature and asymptotic volume ratio} 
\label{11.4proof}

In this section we first show that the asymptotic scalar curvature ratio of
a $\kappa$-solution is infinite.  We then show that the asymptotic
volume ratio vanishes. The proofs are somewhat rearranged from
those in I.11.4.  They are logically independent of
Section \ref{I.11.3}, i.e. also cover the case $n=2$. 
We will use results from
Appendices \ref{appharnack} and \ref{Alexandrov}, 
in particular (\ref{harnackinequality}).

\begin{definition}
If $M$ is a complete connected Riemannian manifold then its {\em asymptotic
scalar curvature ratio} is ${\mathcal R} \: = \:
\limsup_{x \rightarrow \infty} \: R(x) \: d(x, p)^2$. It is independent
of the choice of basepoint $p$.
\end{definition}

\begin{theorem} \label{asympscalar}
Let $(M,g(\cdot))$ be a noncompact $\kappa$-solution.
Then the asymptotic scalar curvature ratio ${\mathcal R}$
is infinite for each time slice.
\end{theorem}
\begin{proof}
Suppose that $M$ is $n$-dimensional,
with $n\geq 2$.   Pick $p\in M$ and consider a time-$t_0$
slice $(M,g(t_0))$.  
We deal with the cases ${\mathcal R} \in (0, \infty)$ and
${\mathcal R} = 0$ separately and show that they lead to contradictions. \\ \\
{\bf Case 1: $0<{\mathcal R}<\infty$.}  We choose a sequence 
$x_k\in M$ such that $d_{t_0}(x_k,p)\ra \infty$ and 
$R(x_k,t_0)d^2(x_k,p)\ra{\mathcal R}$.
Consider the rescaled pointed solution
$(M,x_k,g_k(t))$ with
$g_k(t) =R(x_k,t_0)g(t_0 + \frac{t}{R(x_k,t_0)})$ and $t\in (-\infty,0]$.
We have $R_k(x_k,0)=1$, and for all $b>0$, for sufficiently large $k$,
we have $R_k(x,t)\leq R_k(x,0)\leq \frac{2{\mathcal R}}{d_k^2(x,p)}$
for all $x$ such that $d_k(x,p)>b$.  
Fix numbers $b, B > 0$ so that  $b<\sqrt{{\mathcal R}}<B$.
The $\kappa$-noncollapsing assumption gives a uniform positive
lower bound on the
injectivity radius of $g_k(0)$ at $x_k$,
and so by Appendix \ref{subsequence}
we may extract a pointed limit solution 
$(M_\infty,x_\infty,g_\infty(\cdot))$, defined on a time interval
$(-\infty,0]$ from the sequence
$(M_k,x_k,g_k(\cdot))$ where $M_k =\{x\in M\mid b<d_k(x,p)<B\}$.
Note that $g_\infty$ has nonnegative curvature operator and 
the time slice $(M_\infty,g_\infty(0))$ is locally
isometric to an annular
portion of a nonflat metric cone, since $(M_k,p,g_k(0))$
Gromov-Hausdorff converges to the Tits cone $\ctits(M,g(t_0))$.
(We use the word ``locally'' because the annulus in $\ctits(M,g(t_0))$
need not be geodesically convex in $\ctits(M,g(t_0))$, so we are only 
saying that the distance functions in small balls match up.)
When $n=2$ this contradicts the fact that $R_\infty(x_\infty,0)=1$.
When $n\geq 3$, we will derive a contradiction from 
Hamilton's curvature evolution equation
\begin{equation}
\label{curvevol}
\Rm_t=\Delta \Rm+Q(\Rm).
\end{equation}  

Let $d_v : \ctits(M,g(t_0)) \rightarrow \R$ be the distance function
from the vertex and let 
$\rho:M_\infty\ra \R$ be the pullback of $d_v$ under the inclusion of
the annulus $M_\infty$ in $\ctits(M,g(t_0))$.

\begin{lemma}
The metric cone structure on $(M_\infty,g_\infty(0))$  is
smooth, i.e. $\rho$ is a smooth function.
\end{lemma}
\begin{proof}
Consider a unit speed geodesic segment $\ga$ in the Tits cone 
$\ctits(M,g(t_0))$,
such that $\ga$ is disjoint from the vertex $v\in \ctits(M,g(t_0))$.
Note that since $\ctits(M,g(t_0))$ is a Euclidean cone over
the Tits boundary $\tits (M,g(t_0))$, the
geodesic $\ga$ lies in the cone over a geodesic segment
$\hat\ga\subset\tits (M,g(t_0))$.
Thus $\ga$ lies in a 
$2$-dimensional
locally convex flat subspace of $\ctits(M,g(t_0))$.
Also, as in a flat $2$-dimensional cone,
the second derivative of the composite function 
$d_v^2\circ\ga$ is identically $2$. 

Since $\rho$ is obtained from $d_v$ by composition with a locally isometric
embedding $(M_\infty,g_\infty(0))\ra \ctits(M,g(t_0))$, the composition
of $\rho^2$ with any unit speed geodesic segment in $M_\infty$ 
also has second derivative identically equal to $2$.

Since $\rho^2$ is Lipschitz, Rademacher's theorem implies that $\rho^2$
is differentiable almost everywhere.   
Let $y\in M_\infty$ be a point of differentiability of $\rho^2$.
 If the injectivity radius of $g_\infty(0)$ at 
$y$ is $>a$, then the function $h_y$ given by  the composition 
$$
T_yM_\infty\supset B(0,a)\stackrel{\exp_y}{\lra}\;M_\infty
\stackrel{\rho^2}{\lra}\R
$$
has the property that its second radial derivative is identically $2$,
and it is differentiable at the origin $0\in T_yM_\infty$.
Therefore 
$h_y$ is a second order polynomial with Hessian identically $2$,
and is smooth.  As the injectivity radius is a continuous function, this
implies that $\rho^2$ is smooth everywhere in $M_\infty$.

Since $\rho^2$ is strictly positive on $M_\infty$, 
it follows that $\rho=\sqrt{\rho^2}$ 
is smooth as well.
\end{proof}

By the lemma, we may choose a smooth local orthonormal
frame $e_1,\ldots,e_n$ for $(M_\infty,g_\infty(0))$
near $x_\infty$ such that $e_1$ points radially outward (with respect to the 
cone structure), and $e_2,e_3$ span a $2$-plane at $x_\infty$ with strictly
positive curvature; 
such a $2$-plane exists because $R_\infty(x_\infty, 1) = 1$.
Put $P = e_1 \wedge e_2$.
In terms of the curvature operator, the fact that
$\Rm_\infty(e_1,e_2,e_2,e_1)=0$ is equivalent to 
$\langle P, \Rm_\infty P \rangle \: = \: 0$.
As the curvature operator is nonnegative, it follows that
$\Rm_\infty P \: = \: 0$. (In fact, this is true for any metric cone.)
Differentiating gives
\begin{equation} \label{diff1}
\left( \nabla_{e_i} \Rm_\infty \right) P \: + \: \Rm_\infty (\nabla_{e_i} P) \: = \: 0
\end{equation}
and
\begin{equation} \label{diff2}
\left( \triangle \Rm_\infty \right) P \: + \: 2 \sum_i \left( \nabla_{e_i} \Rm_\infty \right)
\nabla_{e_i} P \: + \: \Rm_\infty \left( \triangle P \right) \: = \: 0.
\end{equation}
Taking the inner product of (\ref{diff2}) with $P$ gives
\begin{align}
0 \: & = \: \langle P, \left( \triangle \Rm_\infty \right) P \rangle  \: + \: 2 \sum_i
\langle P, \left( \nabla_{e_i} \Rm_\infty \right)
\nabla_{e_i} P \rangle \\ 
& = \: \langle P, \left( \triangle \Rm_\infty \right) P \rangle  \: + \: 2 \sum_i
\langle \nabla_{e_i} P, \left( \nabla_{e_i} \Rm_\infty \right)
P \rangle. \notag
\end{align}
Then  (\ref{diff1}) gives
\begin{equation}
\langle P, (\triangle \Rm_\infty) P \rangle \: = \:  2 \:
\sum_i \langle \nabla_{e_i} P,  \Rm_\infty \left( \nabla_{e_i} P \right)  \rangle.
\end{equation}

As the sphere of distance $r$ from the vertex in a metric cone has 
principal curvatures $\frac{1}{r}$, we have
$\nabla_{e_3} e_1 \: = \: - \: \frac{1}{r} \: e_3$.
Then
\begin{equation}
\nabla_{e_3} (e_1 \wedge e_2) \: = \:
\left( \nabla_{e_3} e_1 \right) \wedge e_2 \: + \:
e_1 \wedge \nabla_{e_3} e_2 \: = \: \frac{1}{r} \left( e_2 \wedge e_3 \right) \: + \: 
e_1 \wedge \nabla_{e_3} e_2.
\end{equation}
This shows that $\nabla_{e_3} P$ has a nonradial component
$\frac{1}{r} \: e_2 \wedge e_3$. Thus
$(\Delta \Rm_\infty)(e_1,e_2,e_2,e_1)>0$. 
The zeroth order quadratic term $Q(\Rm)$ appearing in (\ref{curvevol})
is nonnegative when $\Rm$ is nonnegative, so we conclude that 
$\partial_t\Rm_\infty(e_1,e_2,e_2,e_1)>0$ at $t=0$.  This means that
$\Rm_\infty(-\eps)(e_1,e_2,e_2,e_1)<0$ for $\eps>0$ sufficiently small,
which is impossible. \\ \\
{\bf Case 2: ${\mathcal R}=0$.}   
Let us take any sequence $x_k\in M$ with 
$d_{t_0}(x_k,p)\ra \infty$.  Set $r_k = d_{t_0}(x_k,p)$, put
\begin{equation}
g_k(t) =  r_k^{-2}g(t_0 + r_k^2 t)
\end{equation}
for $t\in (-\infty,0]$, and let $d_k(\cdot, \cdot)$ be the distance function
associated to $g_k(0)$.
For any $0<b<B$, put
\begin{equation}
M_k(b,B) \: = \: \{x\in M\mid 0<b<d_k(x,p)<B\}.
\end{equation}
Since ${\mathcal R}=0$, we get that $\sup_{x\in M_k(b,B)}|\Rm_k(x,0)|\ra 0$
as $k\ra\infty$.  Invoking the $\kappa$-noncollapsed assumption as in
the previous case, we
may assume that $(M,p,g_k(0))$ Gromov-Hausdorff converges to a metric
cone $(M_\infty,p_\infty,g_\infty)$ (the Tits cone $\ctits(M,g(t_0))$)
which is flat and smooth away from the
vertex $p_\infty$, and the convergence is smooth away from $p_\infty$.

The ``unit sphere'' in $\ctits(M,g(t_0))$
defines a compact smooth hypersurface $S_\infty$
in $(M_\infty - \{p_\infty\},g_\infty(0))$ whose principal curvatures
are identically $1$. If $n \ge 3$ then $S_\infty$ must be a quotient of
the standard $(n-1)$-sphere by the free action of a finite group of
isometries. We have a sequence $S_k\subset M_k$
of approximating smooth hypersurfaces whose principal curvatures
(with respect to $g_k(0)$) go to $1$ as $k\ra\infty$.  In view of the
convergence to $(M_\infty,p_\infty,g_\infty)$, for 
sufficiently large $k$,
the inward principal curvatures of $S_k$ with respect to $g_k(0)$
are close to $1$. As $M$ has nonnegative curvature,
$S_k$ is diffeomorphic to a sphere \cite[Theorem A]{Eschenburg}.
Thus $S_\infty$ is isometric to the standard $(n-1)$-sphere, 
and so 
$\ctits(M,g(t_0))$ is isometric to $n$-dimensional Euclidean space.
Then  $(M,g(t_0))$ is isometric to $\R^n$, which 
contradicts the definition of a $\kappa$-solution.

In the case $n=2$ we know that $S_\infty$ is diffeomorphic to a circle but
we do not know {\it a priori} that it has length $2\pi$.
To handle the case $n = 2$, we use the fact that $g_k(t)$ is
a Ricci flow solution, to
extract a limiting smooth incomplete
time-independent Ricci flow solution 
$(M_\infty \backslash p_\infty, g_\infty(t))$ for $t \in [-1, 0]$.
Note that this solution is unpointed.
In view of the
convergence to the limiting solution,
for sufficiently large $k$,
the inward principal curvatures of $S_k$ with respect to $g_k(t)$
are close to $1$ for all $t\in [-1,0]$. 
This implies that
$S_k$ bounds a domain $B_k\subset M$ whose diameter with
respect to $g_k(t)$ is uniformly bounded above, say by $10$ 
(see Appendix \ref{Alexandrov}).

Applying the Harnack
inequality (\ref{harnackinequality}) 
with $y_k\in S_k$ (at time $0$) and $x\in B_k$
(at time $-1$), we see that 
$\sup_{x\in B_k}|\Rm_k(x,-1)|\ra 0$ as $k\ra \infty$.  
Thus $(B_k,p,g_k(-1))$ Gromov-Hausdorff converges to a flat manifold 
$(B_\infty,\bar p_\infty,g_\infty(-1))$
with convex boundary. As all of the principal curvatures of
$\partial B_\infty$ are $1$,
$B_\infty$ must be isometric to a Euclidean unit ball.
This implies that $S_\infty$ is isometric to the standard $S^1$ of
length $2\pi$, 
and we obtain a contradiction as before.
\end{proof}

\begin{definition}
If $M$ is a complete $n$-dimensional
Riemannian manifold with nonnegative Ricci curvature
then its {\em asymptotic
volume ratio} is ${\mathcal V} \: = \:
\lim_{r \rightarrow \infty} \: r^{-n} \: \vol(B(p, r))$. It is independent
of the choice of basepoint $p$.
\end{definition}

\begin{proposition} (cf. Proposition I.11.4) \label{propI.11.4}
Let $(M,g(\cdot))$ be a noncompact $\kappa$-solution.
Then the asymptotic volume ratio ${\mathcal V}$
vanishes for each time slice $(M,g(t_0))$.  Moreover, there is a 
sequence of points $x_k\in M$ going to infinity such that the
pointed sequence $\{ (M,(x_k,t_0),g(\cdot)) \}_{k=1}^\infty$ converges, 
modulo rescaling
by $R(x_k,t_0)$, to a $\kappa$-solution which isometrically splits off 
an $\R$-factor.
\end{proposition}
\begin{proof}
Consider the time-$t_0$ slice.
Suppose that ${\mathcal V}>0$.
As ${\mathcal R}=\infty$,  there
are sequences $x_k\in M$ and  $s_k>0$  
such that $d_{t_0}(x_k,p)\ra \infty$, $\frac{s_k}{d_{t_0}(x_k,p)}\ra 0$,
$R(x_k, t_0 )s_k^2\ra\infty$, and $R(x, t_0)\leq 2R(x_k, t_0)$ for all
$x\in B_{t_0}(x_k,s_k)$
\cite[Lemma 22.2]{Hamilton}.
Consider the rescaled pointed solution
$(M,x_k,g_k(t))$ with
$g_k(t) =R(x_k,t_0) \: g(t_0 + \frac{t}{R(x_k,t_0)})$ and $t\in (-\infty,0]$.
As $R_t \ge 0$, we have
$R_k(x,t)\leq 2$ whenever $t\leq 0$ and
$d_k(x,x_k)\leq R(x_k, t_0)^{1/2} s_k$, 
where $d_k$ is the distance function  for $g_k(0)$
and $R_k(\cdot,\cdot)$ is the scalar curvature of $g_k(\cdot)$. 
The $\kappa$-noncollapsing assumption gives a uniform positive
lower bound on the
injectivity radius of $g_k(0)$ at $x_k$,
so by Appendix \ref{subsequence}
we may extract a complete pointed limit solution 
$(M_\infty,x_\infty,g_\infty(t))$, $t\in (-\infty,0]$, of a subsequence of
the sequence of pointed Ricci flows.
By relative volume comparison, $(M_\infty,x_\infty,g_\infty(0))$
has positive asymptotic volume ratio.  By Appendix \ref{Alexandrov},
the Riemannian manifold
$(M_\infty,x_\infty,g_\infty(0))$ is isometric to an Alexandrov
space which splits off a line, which means that it is a Riemannian 
product $\R\times N$. This implies a product structure for
earlier times; 
see Appendix \ref{maxprin}.
Now when $n=2$, we have a contradiction,
since $R(x_\infty,0)=1$ but $(M_\infty,g_\infty(0))$ is a product
surface, and must therefore be flat.  When $n>2$ we obtain
a $\kappa$-solution on an $(n-1)$-manifold with positive asymptotic
volume ratio at time zero, and by induction this is impossible.
\end{proof}

\section{In a $\kappa$-solution, the curvature and the normalized
volume control each other}

In this section we show that, roughly speaking, in a 
$\kappa$-solution
the curvature and the normalized volume control each other.

\begin{corollary} \label{11.4cors}
1.  If $B(x_0,r_0)$ is a ball in a time slice of a $\kappa$-solution,
then the normalized volume $r_0^{-n}\vol(B(x_0,r_0))$ is controlled
(i.e. bounded away from zero) $\Longleftrightarrow$
the normalized scalar curvature $r_0^2R(x_0)$ is controlled
(i.e. bounded above).

2. If $B(x_0,r_0)$ is a ball in a time slice of a $\kappa$-solution,
then the normalized volume $r_0^{-n}\vol(B(x_0,r_0))$ is almost
maximal $\Longleftrightarrow$
the normalized scalar curvature $r_0^2R(x_0)$ is almost zero.

3. (Precompactness) 
If $(M_k,(x_k,t_k),g_k(\cdot))$ is a sequence of pointed $\kappa$-solutions
(without the assumption that $R(x_k, t_k) = 1$)
and for some $r>0$, the $r$-balls $B(x_k,r)\subset (M_k,g_k(t_k))$ 
have controlled normalized volume, then a subsequence converges
to an ancient solution $(M_\infty,(x_\infty,0),g_\infty(\cdot))$
which has nonnegative curvature operator, and is $\kappa$-noncollapsed
(though {\it a priori} the curvature may be unbounded on a given time slice).

4. 
There is a constant 
$\eta=\eta(n,\kappa)$ such that for every $n$-dimensional
$\kappa$-solution $(M,g(\cdot))$, and all $x\in M$, we have
$|\nabla R|(x,t)\leq \eta R^{\frac{3}{2}}(x,t)$ and
$|R_t|(x,t) \leq \eta R^2(x,t)$.
More generally, there are scale invariant bounds on 
all derivatives of the curvature tensor, that only depend
on 
$n$ and
$\kappa$. That is, for 
each $\rho,k,l<\infty$ there is a constant 
$C=C(n,\rho,k,l,\kappa) < \infty$
such that 
$\left|\frac{\D^k}{\D t^k}\; 
\nabla^l\Rm\right|(y,t) \: \le \: C\;R(x,t)^{(k+\frac{l}{2}+1)}$ for any
$y \in B_t(x,\rho R(x,t)^{-\frac12})$.

5.  There is a function $\al:[0,\infty)\ra[0,\infty)$ depending only
on $\kappa$ such that $\lim_{s\ra\infty}\al(s)=\infty$,
and  for every  $\kappa$-solution $(M,g(\cdot))$ and
$x,y\in M$, we have $R(y)d^2(x,y)\geq \al(R(x)d^2(x,y))$.
\end{corollary}
\begin{proof}
Assertion 1, $\Longrightarrow$.  
Suppose we have a sequence of $\kappa$-solutions
$(M_k,g_k(\cdot))$, and sequences 
$t_k\in (-\infty, 0]$, $x_k\in M_k$, $r_k>0$, such that at time $t_k$,
the normalized volume of $B(x_k,r_k)$ is 
$\geq c>0$, and $R(x_k,t_k)r_k^2\ra \infty$.
By Appendix \ref{pointpicking}, 
for each $k$, we can find $y_k\in B(x_k,5r_k)$, 
$\bar r_k\leq r_k$, such that $R(y_k,t_k)\bar r_k^2\geq R(x_k,t_k)r_k^2$,
and $R(z,t_k)\leq 2R(y_k,t_k)$ for all $z\in B(y_k,\bar r_k)$. 
Note that by relative volume comparison,
whenever $\widetilde{r}_k \le \overline{r}_k$ we have
\begin{equation}
\label{volbig}
\frac{\vol(B(y_k,{\widetilde r}_k))}{\widetilde{r}_k^n}\geq
\frac{\vol(B(y_k,\bar r_k))}{\bar r_k^n}\geq 
\frac{\vol(B(y_k,10r_k))}{(10r_k)^n}
\geq \frac{\vol(B(x_k,r_k))}{(10r_k)^n}\geq\frac{c}{10^n}.
\end{equation}
Rescaling the sequence of pointed solutions $(M_k,(y_k,t_k),g_k(\cdot))$
by $R(y_k,t_k)$, we get a sequence satisfying the hypotheses
of Appendix \ref{subsequence} (we use here the fact that $R_t \ge 0$
for an ancient solution), so it accumulates on a limit
flow $(M_\infty,(y_\infty,0),g_\infty(\cdot))$ which is 
a $\kappa$-solution.   
By (\ref{volbig}), the asymptotic volume ratio of $(M,g_\infty(0))$ is 
$\geq \frac{c}{10^n}>0$.  This contradicts Proposition \ref{propI.11.4}.

\bigskip
Assertion 3.  By relative volume comparison, it follows that  every 
$r$-ball in $(M_k,g_k(t_k))$ has normalized volume bounded below
by a ($k$-independent) function of its distance to $x_k$.
By 1, this implies that  the curvature of $(M_k,g_k(t_k))$ is bounded
by a $k$-independent function of the distance to $x_k$, 
and hence we can apply Appendix \ref{subsequence} to extract
a smoothly converging subsequence.

\bigskip
Assertion 1, $\Longleftarrow$.  Suppose we have a sequence
$(M_k,g_k(\cdot))$ of $\kappa$-solutions, and sequences
$x_k\in M_k$, $r_k>0$, such that $R(x_k,t_k)r_k^2<c$ for all
$k$, but $r_k^{-n}\vol(B(x_k,r_k))\ra 0$.
For large $k$, we can choose $\bar r_k \in (0, r_k)$ such that
$\bar r_k^{-n}\vol(B(x_k,\bar r_k))=\frac{1}{2}c_n$
where $c_n$ is the volume of the unit Euclidean $n$-ball.
By relative volume comparison, $\frac{\bar r_k}{r_k}\ra 0$.
Applying 3, we see that the pointed sequence $(M_k,(x_k,t_k),g_k(\cdot))$,
rescaled by the factor $\bar r_k^{-2}$, accumulates
on a pointed ancient solution $(M_\infty,(x_\infty,0),g_\infty(\cdot))$,
such that the ball $B(x_\infty,1)\subset (M_\infty,g_\infty)$ has 
normalized volume $\frac{1}{2}c_n$ at $t=0$.

Suppose the ball $B(x_\infty,1)\subset (M_\infty,g_\infty(0))$ were flat.
Then by the Harnack inequality (\ref{harnackinequality}) (applied to the approximators)
we would have $R_\infty(x,t)=0$ for all $x\in M_\infty$, $t\leq 0$,
i.e. $(M_\infty,g_\infty(t))$ would be a time-independent flat manifold.
It cannot be $\R^n$ since $\vol(B(x_\infty,1)) \: = \: \frac12 \: c_n$.
But flat manifolds other than Euclidean space have zero asymptotic volume ratio
(as follows from the Bieberbach theorem that if $N=\R^n/\Gamma$ is
a flat manifold and $\Ga$ is nontrivial then there is a $\Gamma$-invariant
affine subspace $A\subset \R^n$ of dimension at least $1$ 
on which $\Gamma$ acts cocompactly).
This contradicts the assumption that the sequence $(M_k,g_k(\cdot))$
is $\kappa$-noncollapsed.  Thus
$B(x_\infty,1)\subset (M_\infty,g_\infty(0))$ is not flat, which means,
by the Harnack inequality, that the 
scalar curvature of $g_\infty(0)$ is strictly positive
everywhere.  Therefore, with respect to $g_k$, we have
\begin{equation}
\liminf_{k\ra \infty}R(x_k,t_k)r_k^2=
\liminf_{k\ra\infty}(R(x_k,t_k))\bar 
r_k^2)\left(\frac{r_k}{\bar r_k}\right)^2
\geq \const
\liminf_{k\ra\infty}\left(\frac{r_k}{\bar r_k}\right)^2=\infty,
\end{equation}
which is a contradiction.

\bigskip
Assertion 2, $\Longrightarrow$.  Apply 1,  
the precompactness criterion, and the fact
that a nonnegatively-curved manifold whose balls
have normalized volume $c_n$ must be flat.

\bigskip
Assertion 2, $\Longleftarrow$. Apply 1, the precompactness criterion,
and the Harnack inequality (\ref{harnackinequality}) (to the approximators).

\bigskip
Assertion 4.  This follows by rescaling $g$ so that $R(x,t)=1$, and
applying 1 and 3.

\bigskip
Assertion 5.  The quantity $R(z)d^2(u,v)$ is scale invariant.
If the assertion failed then  we would have sequences $(M_k,g_k(\cdot))$,
$x_k,y_k\in M_k$, such that $R(y_k)=1$ and
$d(x_k,y_k)$ remains bounded, but the curvature at $x_k$ blows up.
This contradicts $1$ and $3$.

\section{An alternate proof of Corollary \ref{2Dsoliton}
using Proposition \ref{propI.11.4} and Corollary \ref{11.4cors}}
\label{altpf11.3}

In this section we give an alternate proof of Corollary \ref{2Dsoliton}.
It uses Proposition \ref{propI.11.4} and Corollary \ref{11.4cors}
To clarify the chain of logical dependence, we remark that
this section is concerned with $2$-dimensional $\kappa$-solutions,
and does not use anything from Sections \ref{I.11.2} or \ref{I.11.3}.   It does use
Proposition \ref{propI.11.4}.  However, we avoid circularity here because the proof
of Proposition \ref{propI.11.4} given in Section \ref{11.4proof}, unlike the
proof in \cite{Perelman}, does not use Corollary \ref{2Dsoliton}.

\begin{lemma}
\label{2dvol}
There is a constant $v=v(\kappa)>0$ such that if
$(M,g(\cdot))$ is a 2-dimensional $\kappa$-solution
({\it a priori} either compact or noncompact), $x,y\in M$ and 
$r = d(x,y)$ then 
\begin{equation}
\label{volnotsmall}
\vol(B_t(x,r))\geq vr^2.
\end{equation}
\end{lemma}
\begin{proof}
If the lemma were not true then there would be a sequence $(M_k,g_k(\cdot))$
of 2-dimensional $\kappa$-solutions, and sequences $x_k,y_k\in M_k$, $t_k\in\R$
such that $r_k^{-2}\vol(B_{t_k}(x_k,r_k))\ra 0$, where
$r_k = d(x_k,y_k)$.  Let $z_k$ be the midpoint of a shortest segment
from $x_k$ to $y_k$ in the $t_k$-time slice $(M_k,g_k(t_k))$. For large $k$,
choose
$\bar r_k \in (0, r_k/2)$ such that 
\begin{equation}
\label{halfvol}
\bar r_k^{-2}\vol(B_{t_k}(z_k,\bar r_k))=\frac{\pi}{2},
\end{equation}
 i.e.
half the area of the  unit disk in $\R^2$.   
As
\begin{equation}
\frac{\pi}{2} = \bar r_k^{-2}\vol(B_{t_k}(z_k,\bar r_k)) \le
\bar r_k^{-2}\vol(B_{t_k}(x_k, r_k)) =
(\bar r_k / r_k)^{-2} \:  r_k^{-2} \vol(B_{t_k}(x_k, r_k)), 
\end{equation}
it follows that
$\lim_{k\ra\infty}\frac{\bar r_k}{ r_k}=0$.  
Then by part 3 of Corollary \ref{11.4cors},
the sequence of pointed Ricci flows $(M_k,(z_k,t_k),g_k(\cdot))$,
when rescaled by $\bar r_k^{-2}$, accumulates on a complete Ricci flow
$(M_\infty,(z_\infty,0),g_\infty(\cdot))$.  The segments from $z_k$
to $x_k$ and $y_k$ accumulate on a line in $(M_\infty,g_\infty(0))$,
and hence $(M_\infty,g_\infty(0))$ splits off a line.  By (\ref{halfvol}),
$(M_\infty,g_\infty(0))$ cannot be isometric to $\R^2$, and hence
must be a cylinder.  Considering the approximating
Ricci flows, we get a contradiction to the $\kappa$-noncollapsing assumption.
\end{proof}

\bigskip
Lemma \ref{2dvol} implies that the asymptotic volume ratio of any
noncompact 2-dimensional $\kappa$-solution is at least $v>0$.
By Proposition \ref{propI.11.4}
we therefore conclude that every 2-dimensional $\kappa$-solution
is compact.  (This was implicitly assumed in the
proof of Corollary I.11.3 in \cite{Perelman}, as its reference
\cite{Hamsurface} is about compact surfaces.)

Consider the family $\f$ 
of 2-dimensional $\kappa$-solutions $(M,(x,0),g(\cdot))$
with $\diam(M,g(0))=1$.  By Lemma \ref{2dvol}, there is uniform
lower bound on the volume of the $t=0$ time slices of $\kappa$-solutions in 
$\f$. Thus $\f$  is compact in the smooth topology by part 3 of 
Corollary \ref{11.4cors}
(the precompactness leads to compactness in view of the diameter bound).
This implies 
(recall that $R > 0$) that there is a constant $K\geq 1$ such that 
every time slice of every 2-dimensional $\kappa$-solution has
$K$-pinched curvature.    

Hamilton has shown that volume-normalized Ricci 
flow on compact surfaces with positively pinched initial data
converges exponentially fast to a constant curvature metric
\cite{Hamsurface}.  
His argument shows that there is a small $\eps>0$, depending 
continuously on the initial data, so that when the volume of the 
(unnormalized) solution
has  gone down by a factor of 
at least $\eps^{-1}$,
the pinching is at most the square 
root of the initial pinching.  By the compactness of the family $\f$, 
this $\eps$ can be chosen
uniformly when we take the initial data to be the $t=0$ time slice
of a  $\kappa$-solution in $\f$.
  
Now let $K_0$ be the worst pinching of
a 2-dimensional $\kappa$-solution, and let
 $(M,g(\cdot))$ be a $\kappa$-solution where the curvature pinching of 
$(M,g(0))$ is $K_0$.  Choosing $t<0$ such that 
$\eps\vol(M,g(t))=\vol(M,g(0))$,
the previous paragraph implies the curvature pinching of $(M,g(t))$ is at
least $K_0^2$. This would
 contradict the fact that $K_0$ is the upper bound on the
pinching for all $\kappa$-solutions, unless $K_0=1$.

\section{I.11.5. A volume bound}

In this section we give a consequence of Proposition \ref{propI.11.4}
concerning the volumes of metric balls in Ricci flow solutions
with nonnegative curvature operator.

\begin{corollary} (cf. Corollary I.11.5) \label{corI.11.5}
For every $\epsilon > 0$, there is an $A < \infty$ with the following
property.  Suppose that we have a sequence of (not necessarily
complete) Ricci flow solutions $g_k(\cdot)$ with nonnegative
curvature operator, defined on $M_k \times [t_k, 0]$, such that \\
1.  For each $k$, the time-zero ball $B(x_k, r_k)$ has compact closure
in $M_k$. \\
2. For all $(x, t) \in B(x_k, r_k) \times
[t_k, 0]$, $\frac12 R(x,t) \le R(x_k, 0) = Q_k$. \\
3. $\lim_{k \rightarrow \infty} t_k Q_k = - \infty$. \\
4. $\lim_{k \rightarrow \infty} r_k^2 Q_k = \infty$. \\ 
Then for large $k$, $\vol(B(x_k, A Q_k^{- \: \frac12})) \le 
\epsilon (A Q_k^{- \: \frac12})^n$ at time zero.
\end{corollary}
\begin{proof}
Given $\epsilon > 0$, suppose that the corollary is not true.
Then there is a sequence of such Ricci flow solutions with
$\vol(B(x_k, A_k Q_k^{- \: \frac12})) \: > \: \epsilon 
(A_k Q_k^{- \: \frac12})^n$ at time zero, where
 $A_k \rightarrow \infty$.  
 By Bishop-Gromov, $\vol(B(x_k, Q_k^{- \: \frac12})) \: > \: \epsilon 
Q_k^{- \: \frac{n}{2}}$ at time zero, so we can parabolically rescale by
$Q_k$ and take a convergent subsequence. The limit $(M_\infty,
g_\infty(\cdot))$ will be a nonflat complete ancient solution with
nonnegative curvature operator, bounded curvature and
 ${\mathcal V}(0) \: > \: 0$. 
 By Proposition \ref{propI.11.4}, it cannot be $\kappa$-noncollapsed
 for any $\kappa$. 
Thus for each $\kappa \: > \: 0$,
there are a point $(x_\kappa, t_\kappa) \in M_\infty \times (-\infty,0]$ and a radius
$r_\kappa$ so that $|\Rm(x_\kappa, t_\kappa)| \: \le \:
r_\kappa^{-2}$ on the time-$t_\kappa$ ball $B(x_\kappa, r_\kappa)$,
but $\vol(B(x_\kappa, r_\kappa)) \: < \: \kappa \:
r_\kappa^n$. From the Bishop-Gromov inequality,
${\mathcal V}(t_\kappa) \: < \: \kappa$ for the limit solution.

We claim that
${\mathcal V}(t)$ is nonincreasing in $t$. To see this, we have
$\frac{\dvol(U)}{dt} \: = \: \int_U R \: dV  \ge 0$ for any domain $U
\subset M_\infty$.
Also, as $R \le 2$ on $M_\infty \times (-\infty, 0]$, 
Corollary \ref{bound1}
gives that distances on $M_\infty$
decrease at most linearly in $t$,
which implies the claim.

Thus ${\mathcal V}(0) = 0$ for the limit solution, which is a contradiction.
\end{proof}

\section{I.11.6. Curvature bounds for Ricci flow solutions with
nonnegative curvature operator, assuming a lower volume bound}

In this section we show that for a Ricci flow solution with nonnegative
curvature operator, a lower bound on the volume of a ball
implies an earlier upper curvature bound on a slightly smaller ball.
This will be used in Section \ref{I.12.3}.

\begin{corollary} (cf. Corollary I.11.6) \label{corI.11.6}
For every $w > 0$, there are $B = B(w) < \infty$, $C = C(w) < \infty$ and
$\tau_0 = \tau_0(w) > 0$ with the following properties. \\
(a) Take  $t_0 \in  [- r_0^2, 0)$. 
Suppose that we have a (not necessarily complete) Ricci flow
solution $(M, g(\cdot))$, defined for $t \in [t_0, 0]$, so that at time zero
the metric ball $B(x_0, r_0)$ has compact closure. Suppose that 
for each $t \in [t_0, 0]$,
$g(t)$ has nonnegative curvature operator and $\vol(B_t(x_0, r_0))
\ge w r_0^n$. Then 
\begin{equation} \label{Rbound3}
R(x,t) \le C r_0^{-2} + B (t - t_0)^{-1}
\end{equation}
whenever
$\dist_t(x, x_0) \le \frac14 r_0$. \\
(b)
Suppose that we have a (not necessarily complete) Ricci flow
solution $(M, g(\cdot))$, defined for $t \in [- \tau_0 r_0^2, 0]$, 
so that at time zero
the metric ball $B(x_0, r_0)$ has compact closure. Suppose that 
for each $t \in [- \tau_0 r_0^2, 0]$,
$g(t)$ has nonnegative curvature operator. 
If we assume a time-zero volume bound $\vol(B_0(x_0, r_0))
\ge w r_0^n$ then
\begin{equation} 
R(x,t) \le C r_0^{-2} + B (t + \tau_0 r_0^2)^{-1}
\end{equation}
whenever $t \in [- \tau_0 r_0^2,0]$ and $\dist_t(x, x_0) \le \frac14 r_0$.
\end{corollary}
\begin{remark}
The statement in \cite[Corollary 11.6(a)]{Perelman}
does not have any constraint on $t_0$. 
In our proof we seem to need that $-t_0 \le c r_0^2$ for some arbitrary
but fixed
constant $c < \infty$.
(The statement $R(\overline{x}, \overline{t}) >
C + B (\overline{t} - t_0)^{-1}$ in 
\cite[Proof of Corollary 11.6(a)]{Perelman} is the issue.) 
For simplicity we take
$- t_0 \le r_0^2$. This point does not affect the proof of 
Corollary \ref{corI.11.6}(b), which is what ends up getting used.
\end{remark}
\begin{proof}
For part (a), we can assume that $r_0 \: = \: 1$. Given $B,C > 0$, 
suppose that $g(\cdot)$ is a Ricci flow solution for
$t \in [t_0, 0]$ that satisfies the hypotheses of the corollary, 
with $R(x, t) \: > \: C \: + \: B(t - t_0)^{-1}$
for some $(x, t)$ satisfying $\dist_t(x, x_0) \: \le \: \frac{1}{4}$.
Following the notation of the proof of Theorem \ref{thmI.10.1},
except changing the $A$ of Theorem \ref{thmI.10.1} to
$\widehat{A}$,
put $\widehat{A} \: = \: \lambda C^{\frac12}$ and $\alpha \: = \:
\min(\lambda^2 C^{\frac12}, B)$, where we will take $\lambda$ to be a 
sufficiently small number that only depends on $n$. Put
\begin{equation}
M_\alpha \: = \: \{(x^\prime, t^\prime) \: : \: R(x^\prime, t^\prime) \: \ge \: 
\alpha (t^\prime-t_0)^{-1}\}.
\end{equation}
Clearly $(x, t) \in M_\alpha$.
 
We first go through the analog of the proof of Lemma \ref{Claim 1}.
We claim that there is some $(\overline{x}, \overline{t}) \in
M_\alpha$, with $\overline{t} \in (t_0, 0]$ and 
$\dist_{\overline{t}}(\overline{x}, x_0) \: \le \: \frac13$, such that
$R(x^\prime, t^\prime) \: \le \: 2 Q \: = 
\: 2R(\overline{x}, \overline{t})$ whenever
$(x^\prime, t^\prime) \in M_\alpha$, $t^\prime \in (t_0, \overline{t}]$
and $\dist_{t^\prime}(x^\prime, x_0) \le 
\dist_{\overline{t}}(\overline{x}, x_0) \: + \: \widehat{A}Q^{- \: \frac12}$.
Put $(x_1, t_1) = (x, t)$. Inductively, if we cannot take
$(x_k, t_k)$ for $(\overline{x}, \overline{t})$ then there is some
$(x_{k+1}, t_{k+1}) \in M_\alpha$ with 
$t_{k+1} \in (t_0, t_k]$,
$R(x_{k+1}, t_{k+1}) > 2
R(x_k, t_k)$ and $\dist_{t_{k+1}}(x_{k+1}, x_0) \le 
\dist_{t_k}(x_k, x_0) \: + \: \widehat{A}R(x_k, t_k)^{- \: \frac12}$. As the
process must terminate, we end up with $(\overline{x}, \overline{t})$
satisfying 
\begin{equation}
\dist_{\overline{t}}(\overline{x}, x_0) \: \le \: \frac14
\: + \: \frac{1}{1-\sqrt{1/2}} \: \widehat{A} \: R(x,t)^{- \: \frac12}
\: \le \: \frac13
\end{equation}
if $\lambda$ is sufficiently small.

Next, we go through the analog of the proof of Lemma \ref{Claim 2}.
As in the proof of Lemma \ref{Claim 2}, $R(x^\prime, t^\prime) \le 2
R(\overline{x}, \overline{t})$ whenever
$\overline{t}-\frac12 \alpha Q^{-1} \le t^\prime \le \overline{t}$
and $\dist_{t^\prime}(x^\prime, x_0) \le 
\dist_{\overline{t}}(\overline{x}, x_0) + \widehat{A} Q^{- \: \frac12}$. We
claim that the time-$\overline{t}$ ball
$B(x_0, \dist_{\overline{t}}(\overline{x}, x_0) + \frac{1}{10} \widehat{A}
Q^{-1/2})$ is contained in the time-$t^\prime$ ball
$B(x_0, \dist_{\overline{t}}(\overline{x}, x_0) + \widehat{A} Q^{- \: \frac12})$.
To see this, we apply Lemma \ref{distancedist} with
$r_0 \: = \: \frac{1}{2} Q^{-1/2}$ to give
\begin{equation}
\dist_t(x_0, \overline{x}) \: - \: 
\dist_{\overline{t}}(x_0, \overline{x}) \: \le \: \const(n) \:
\alpha Q^{-1/2} \: \le \: \lambda \const(n) \widehat{A} Q^{-1/2}.
\end{equation}
If $\lambda$ is sufficiently small then the claim follows. The argument
also shows that it is consistent to use the curvature bound when
applying Lemma \ref{distancedist}.

Hence $R(x^\prime, t^\prime) \le 2
R(\overline{x}, \overline{t})$ whenever
$\overline{t} \: - \: \frac{1}{2}  \alpha 
Q^{-1} \: \le \: t^\prime \: \le \: \overline{t}$ and
$\dist_{\overline{t}}(x^\prime, x_0) \: \le \:
\dist_{\overline{t}}(\overline{x}, x_0) + \frac{1}{10} \widehat{A}
Q^{-1/2}$. It follows that 
$R(x^\prime, t^\prime) \le 2
R(\overline{x}, \overline{t})$ whenever
$\overline{t} \: - \: \frac{1}{2}  \alpha 
Q^{-1} \: \le \: t^\prime \: \le \: \overline{t}$ and
$\dist_{\overline{t}}(x^\prime, \overline{x}) \: \le \: \frac{1}{10} 
\widehat{A} Q^{-1/2}$. This shows that there is an $A^\prime = A^\prime
(B, C)$, which goes
to infinity as $B, C \rightarrow \infty$, so that
$R(x^\prime, t^\prime) \le 2
R(\overline{x}, \overline{t})$ whenever
$\overline{t} \: - \: A^\prime
Q^{-1} \: \le \: t^\prime \: \le \: \overline{t}$ and
$\dist_{\overline{t}}(x^\prime, \overline{x}) \: \le \: A^\prime
Q^{-1/2}$.

Now suppose that Corollary \ref{corI.11.6}(a) is not true. 
Fixing $w > 0$, 
for any sequences $\{B_k\}_{k=1}^\infty$ and $\{C_k\}_{k=1}^\infty$
going to infinity and
for each $k$, there is a Ricci flow solution $g_k(\cdot)$
which satisfies the hypotheses of
the corollary but for which $R(x_k,t_k) \ge C_k \: + \: B_k (t_k - t_{0,k})^{-1}$
for some point $(x_k, t_k)$ satisfying $\dist_{t_k}(x_k, x_{0,k}) \le \frac14$.
We can assume that $\lambda^2 \: C_k^{\frac12} \: \ge \: B_k$.
From the preceding discussion, there is a sequence $A^\prime_k \rightarrow
\infty$ and points $(\overline{x}_k, \overline{t}_k)$ 
with $\dist_{\overline{t}_k}(\overline{x}_k, {x}_{0,k}) \le \frac13$ so that
$R(x_k^\prime, t_k^\prime) \le 2
R(\overline{x}_k, \overline{t}_k)$ whenever
$\overline{t}_k \: - \: A_k^\prime
Q_k^{-1} \: \le \: t_k^\prime \: \le \: \overline{t}_k$ and
$\dist_{\overline{t}_k}(x_k^\prime, \overline{x}_k) \: \le \: A_k^\prime
Q_k^{-1/2}$, where 
\begin{equation}
Q_k \: = \: R(\overline{x}_k, \overline{t}_k) \: \ge \:
B_k \: (\overline{t}_k - t_{0,k})^{-1} \: \ge \: B_k.
\end{equation}
By Corollary \ref{corI.11.5}, for any $\epsilon > 0$ there is some 
$A = A(\epsilon)< \infty$
so that for large $k$,
\begin{equation}
\vol(B(\overline{x}_k, A/\sqrt{Q_k})) \: \le \: \epsilon
(A/\sqrt{Q_k})^n
\end{equation}
at time zero. By the Bishop-Gromov inequality, 
$\vol(B(\overline{x}_k, 1)) \: \le \: \epsilon$ for large $k$, since
$Q_k \rightarrow \infty$.
If we took $\epsilon$ sufficiently small from the beginning then we 
would get a contradiction to the fact that
\begin{equation}
\vol(B(\overline{x}_k, 1)) \: \ge \: \vol \left( B \left( x_{0,k}, \frac23 \right) \right) \: \ge 
\left( \frac23 \right)^n \vol(B(x_{0,k}, 1)) \: \ge \left( \frac23 \right)^n w.
\end{equation}

For part (b), the idea is to choose the parameter 
$\tau_0$ sufficiently small so that
we will still have the estimate 
$\vol(B(x_0, r_0)) \: \ge \: 5^{-n} \: w \: r_0^n$ for the time-$t$ ball
$B(x_0, r_0)$ when $t \in [-\tau_0 r_0^2, 0]$,
and so we can apply part (a) with $w$ replaced by $\frac{w}{5}$.
The value of $\tau_0$ will emerge from the proof.
More precisely, putting $r_0 \: = \: 1$ and with a given 
$\tau_0$, let $\tau$ be the
largest number in $[0, \tau_0]$ so that the time-$t$ ball
$B(x_0, 1)$ satisfies
$\vol(B(x_0, 1)) \: \ge \: 5^{-n} \: w$ whenever $t \in [-\tau, 0]$.
If $\tau \: < \: \tau_0$ then at time $- \tau$, we have
$\vol(B(x_0, 1)) \: = \: 5^{-n} \: w$.
The conclusion of part (a) holds in the sense that
\begin{equation} \label{cc}
R(x, t) \: \le \: C(5^{-n} w) \: + \: B(5^{-n} w) ( t + \tau)^{-1}
\end{equation}
whenever $t \in [- \tau, 0]$ and $\dist_t(x, x_0) \: \le \: \frac14$.
Lemma \ref{distancedist}, along with (\ref{cc}), implies that the 
time-$(- \tau)$ ball $B(x_0, \frac14)$ contains the
time-$0$ ball $B(x_0, \frac14 \: - \: 10(n-1)(\tau \sqrt{C} \: + \:
2 \sqrt{B \tau}))$. From the nonnegative curvature, the 
time-$(- \tau)$ volume of the first ball is at least as large as the
time-$0$ volume of the second ball.  Then
\begin{align} \label{lll}
5^{-n} \: w \: & = \: \vol(B(x_0, 1)) 
\: \ge \: \vol(B(x_0, \frac14)) \\
& \ge \: 
\vol(B(x_0, \frac14 \: - \: 10(n-1)(\tau \sqrt{C} \: + \:
2 \sqrt{B \tau}))) \notag \\
& \ge 
(\frac14 \: - \: 10(n-1)(\tau \sqrt{C} \: + \:
2 \sqrt{B \tau}))^n \: \vol(B(x_0, 1)) \notag \\
& \ge \:
(\frac14 \: - \: 10(n-1)(\tau \sqrt{C} \: + \:
2 \sqrt{B \tau}))^n \: w, \notag
\end{align}
where the balls on the top line of (\ref{lll}) are at time-($- \tau$),
and the other balls are at time-$0$. 
Thus $\frac14 \: - \: 10(n-1)(\tau \sqrt{C} \: + \:
2 \sqrt{B \tau}) \: \le \: \frac15$.
This contradicts our assumption
that $\tau \: < \: \tau_0$ provided that
$\frac14 \: - \: 10(n-1)(\tau_0 \sqrt{C} \: + \:
2 \sqrt{B \tau_0}) \: = \: \frac15$.
\end{proof}

Finally, we give a version of Corollary \ref{corI.11.6}(b)  where instead of assuming a
nonnegative curvature operator, we assume that the 
curvature operator
in the time-dependent ball of radius $r_0$ around $x_0$ is bounded below by
$- \: r_0^{-2}$.

\begin{corollary} (cf. end of Section I.11.6)\label{I.11.6end}
For every $w > 0$, there are $B = B(w) < \infty$, $C = C(w) < \infty$ and
$\tau_0 = \tau_0(w) > 0$ with the following property.
Suppose that we have a (not necessarily complete) Ricci flow
solution $(M, g(\cdot))$, defined for $t \in [- \tau_0 r_0^2, 0]$, so that at time zero
the metric ball $B(x_0, r_0)$ has compact closure. Suppose that 
for each $t \in [- \tau_0 r_0^2, 0]$, the 
curvature operator in the
time-$t$ ball $B(x_0, r_0)$ is bounded below by $- \: r_0^{-2}$. 
If we assume a time-zero volume bound $\vol(B_0(x_0, r_0))
\ge w r_0^n$ then
\begin{equation} 
R(x,t) \le C r_0^{-2} + B (t + \tau_0 r_0^2)^{-1}
\end{equation}
whenever $t \in [- \tau_0 r_0^2,0]$ and $\dist_t(x, x_0) \le \frac14 r_0$.
\end{corollary}
\begin{proof}
The blowup argument goes through as before.
The only real difference is that the volume of the
time-$(- \tau)$ ball $B(x_0, \frac14)$ will be at least
$e^{-\const \tau r_0^{-2}}$ times the volume of the
time-$0$ ball $B(x_0, \frac14 \: - \: 10(n-1)(\tau \sqrt{C} \: + \:
2 \sqrt{B \tau}))$.
\end{proof}

\section{I.11.7. Compactness of the space of three-dimensional
$\kappa$-solutions}
\label{pf11.7}

In this section we prove a compactness result for the space of three-dimensional
$\kappa$-solutions.  The three-dimensionality assumption is used to
show that
the limit solution has bounded curvature.

If a three-dimensional $\kappa$-solution $M$ is compact then it is
diffeomorphic to a quotient of $S^3$ or $\R \times S^2$, as it has 
nonnegative curvature and is nonflat.
If its asymptotic 
soliton (see Section \ref{I.11.2}) is also closed then 
$M$ is a quotient of the round $S^3$ or $\R \times S^2$. 
There are $\kappa$-solutions on $S^3$ and $\R P^3$ with noncompact
asymptotic soliton; see \cite[Section 1.4]{Perelman2}.
They are not isometric to the round metric; this corrects the statement in
the first paragraph of \cite[Section 11.7]{Perelman}.

\begin{theorem} (cf. Theorem I.11.7) \label{thmI.11.7}
Given $\kappa > 0$, the set of oriented
three-dimensional $\kappa$-solutions is compact modulo scaling.  That is, from any
sequence of such solutions and points $(x_k, 0)$, 
after appropriate dilations 
we can extract a
smoothly converging subsequence that satisfies the same conditions.
\end{theorem}
\begin{proof}
If $(M_k, (x_k, 0), g_k(\cdot))$ is a sequence of such $\kappa$-solutions with
$R(x_k, 0) = 1$ then
parts 1 and 3 of Corollary \ref{11.4cors} imply there is a subsequence that
converges to an ancient solution $(M_\infty, (x_\infty, 0), g_\infty(\cdot))$
which has nonnegative curvature operator and is $\kappa$-noncollapsed.
The remaining issue is to show that it has bounded curvature. 
Note that $R_t\geq 0$ since $g_\infty(\cdot)$  is  a limit of a sequence
of Ricci flows satisfying $R_t\geq 0$.  
Hence it is
enough to show that $(M_\infty, g(0))$ has bounded scalar curvature.

If not, there is a sequence of points $y_i$ going to infinity in $M_\infty$
such that $R(y_i, 0) \rightarrow \infty$ and
$R(y,0) \le 2 R(y_i,0)$ for $y \in B(y_i, A_i R(y_i,0)^{- \: \frac12})$, where
$A_i \rightarrow \infty$; compare \cite[Lemma 22.2]{Hamilton}. 
Using the $\kappa$-noncollapsing, a subsequence of the rescalings
$(M_\infty, y_i, R(y_i, 0) g_\infty)$ will converge to a limit manifold
$N_\infty$.
As in the proof of Proposition \ref{propI.11.4}
from Appendix \ref{Alexandrov},
$N_\infty$ will split off a line.
By Corollary \ref{2Dsoliton} or Section \ref{altpf11.3}, 
$N_\infty$ must be the standard solution
on $\R \times S^2$. Thus $(M_\infty, g(0))$ contains a sequence
$D_i$ of neck regions, with their 
cross-sectional radii
tending to zero as
$i \rightarrow \infty$.  

Note that $M_\infty$
has to be 1-ended.  Otherwise, it would contain a line,
and would therefore have to split off a line isometrically
\cite[Theorem 8.17]{Cheeger-Ebin}.
But then $M_\infty$, the product of a line and a surface,
could not have neck regions with cross-sections tending to zero.

From the theory of nonnegatively curved manifolds
\cite[Chapter 8.5]{Cheeger-Ebin},  there is an
exhaustion $M_\infty \: = \: \bigcup_{t\ge 0}C_t$ 
by nonempty totally convex compact sets $C_t$ so that
$(t_1 \le t_2) \Rightarrow (C_{t_1} \subset C_{t_2})$, and
\begin{equation} \label{C}
C_{t_1} \: = \: \{ q \in C_{t_2} \: : \: \dist(q, \partial C_{t_2})
\: \ge \: t_2 - t_1\}. 
\end{equation}
Now consider a 
neck region $D$ which is  close to a cylinder.  Note by
triangle comparison -- or simply because the distance function
in $D$ is close to that of a product metric --  any 
minimizing geodesic  segment $\ga\subset D$
of length large compared to cross-sectional radius of $D$ must
be nearly orthogonal to the cross-section.   
It follows from this and (\ref{C}) that if $t>0$ and 
$\partial C_t$ contains a point $p\in D$ such that $d(p,\D D)$
is large compared to the cross-section of $D$, then 
 $\partial C_t\cap D$ is an approximate $2$-sphere 
cross-section of $D$. 
Fix such a neck region $D_0$ and let $C_{t_0}$ be the
corresponding convex set.
As $M_\infty$ 
has one end, $\partial C_{t_0}$ has only one connected 
component,
namely the approximate $2$-sphere cross-section.

For all $t \: > \: t_0$, there is a distance-nonincreasing retraction
$r \: : \: C_t \rightarrow C_{t_0}$ which maps $C_t \: - \: C_{t_0}$ onto
$\partial C_{t_0}$ \cite{Sharafutdinov}. Let $D$ be a neck region
with a very small cross-section and let $C_t$ be a convex set so that
$\partial C_t$ intersects $D$ in an approximate $2$-sphere cross-section.
Then $\partial C_t$ consists entirely of this approximate cross-section.
The restriction of $r$ to $\partial C_t$ is distance-nonincreasing, but
will map the $2$-sphere $\partial C_t$ onto
the $2$-sphere $\partial C_{t_0}$. This is a contradiction. 
\end{proof}

\begin{remark} The statement of \cite[Theorem 11.7]{Perelman} is about
noncompact $\kappa$-solutions but the proof works whether the
solutions are compact or noncompact.
\end{remark}

\begin{remark}
 One may wonder where we have used the fact that we have a
Ricci flow solution, i.e. whether the curvature is bounded for
any $\kappa$-noncollapsed Riemannian $3$-manifold with nonnegative
sectional curvature.  Following the above argument, we could again split off
a line in a rescaling around high-curvature points.  However, we would
not necessarily know that the ensuing nonnegatively-curved
surface is compact. ({\it A priori}, it could be a smoothed-out cone, for
example.)  In the case of a Ricci flow, the compactness comes from 
Corollary \ref{2Dsoliton} or Section \ref{altpf11.3}.
\end{remark}

\begin{corollary}
Let $(M,g(\cdot))$ be a $3$-dimensional $\kappa$-solution.
Then any asymptotic soliton constructed as in Section \ref{I.11.2} is
also a $\kappa$-solution.
\end{corollary}

\section{I.11.8. Necklike behavior at infinity of a three-dimensional
$\kappa$-solution - weak version}

The next corollary says that outside of a compact region,  any oriented noncompact
three-dimensional $\kappa$-solution looks necklike (after rescaling). 
In this section we give a simple argument to prove the corollary, except for
a diameter bound on the compact region.  In the next section we give an
argument that also proves the diameter bound.

More information on three-dimensional $\kappa$-solutions is in 
Section \ref{II.1}. 

\begin{definition}
Fix $\epsilon > 0$. Let $(M, g(\cdot))$ be an oriented three-dimensional $\kappa$-solution.
We say that a point $x_0 \in M$ is the {\em center of an $\epsilon$-neck} if the
solution $g(\cdot)$ in the set 
$\{(x,t) \: : \: - (\epsilon Q)^{-1} < t \le 0, \dist_0(x, x_0)^2 < (\epsilon Q)^{-1}\}$,
where $Q = R(x_0, 0)$, is, after scaling with the factor $Q$, $\epsilon$-close
in some fixed smooth topology to the corresponding subset of the
evolving round cylinder (having scalar curvature one at time zero).
(See Definition \ref{closenessdef} below for a more precise statement.)

We let $M_\epsilon$ denote the points in $M$ that are not centers of
$\epsilon$-necks.
\end{definition}

\begin{corollary} (cf. Corollary I.11.8) \label{corI.11.8}
For any $\epsilon > 0$, there exists $C = C(\epsilon, \kappa) > 0$ such that
if $(M, g(\cdot))$ is an oriented noncompact three-dimensional $\kappa$-solution then \\
1.  $M_\epsilon$ is compact with $\diam(M_\epsilon) \le C Q^{- \: \frac12}$ and \\
2. $C^{-1} Q \le R(x,0) \le CQ$ whenever $x \in M_\epsilon$, \\
where $Q = R(x_0, 0)$ for some $x_0 \in \partial M_\epsilon$.
\end{corollary}
\begin{proof}
We  prove here the claims of 
Corollary \ref{corI.11.8}, except for the
diameter bound. In the next section we give another argument which
also proves the diameter bound.

We claim first that $M_\epsilon$ is compact.  Suppose not.
Then there is a sequence of points $x_k \in M_\epsilon$ going to
infinity. Fix a basepoint $x_0 \in M$. Then
$R(x_0) \: \dist_0^2(x_0, x_k) \rightarrow \infty$. By
part 5 of Corollary \ref{11.4cors},
$R(x_k) \: \dist_0^2(x_0, x_k) \rightarrow \infty$.  Rescaling
around $(x_k, 0)$ to make its scalar curvature one, we can use 
Theorem \ref{thmI.11.7} to extract a convergent subsequence
$(M_\infty, x_\infty)$.
As in the proof of Proposition \ref{propI.11.4}, we can say that 
$(M_\infty, x_\infty)$ splits off a line.  Hence for large $k$,
$x_k$ is the center of an $\epsilon$-neck, which is a contradiction.

Next we claim that for any $\epsilon$, there exists
$C = C(\epsilon, \kappa) > 0$ such that if $g_{ij}(t)$ is a $\kappa$-solution 
then for
any point $x \in M_\epsilon$, there is a point $x_0 \in
\partial M_\epsilon$ such that $\dist_0(x, x_0) \: \le \: CQ^{-1/2}$
and $C^{-1} Q \: \le \: R(x, 0) \: \le \: CQ$, where
$Q \: = \: R(x_0, 0)$.

If not then there is a sequence $\{M_i\}_{i=1}^\infty$
of $\kappa$-solutions along with points
$x_i \in M_{i,\epsilon}$ such that for each $y_i \in 
\partial M_{i,\epsilon}$, we have \\
1. $\dist_0^2(x_i, y_i) \: R(y_i, 0) \: \ge \: i$ or \\
2. $R(y_i, 0) \: \ge \: i \: R(x_i, 0)$ or \\
3. $R(x_i, 0) \: \ge \: i \: R(y_i, 0)$. 

Rescale the metric on $M_i$ so that $R(x_i, 0) \: = \: 1$. From
Theorem \ref{thmI.11.7}, a subsequence of the pointed spaces
$(M_i, x_i)$ will converge smoothly to a $\kappa$-solution
$(M_\infty, x_\infty)$. Also, $x_\infty \in M_{\infty, \epsilon}$.

Taking a subsequence, we
can assume that 1. occurs for each $i$, or 2. occurs for each $i$,
or 3. occurs for each $i$. If $M_\infty \neq M_{\infty, \epsilon}$, 
choose $y_\infty \in
\partial M_{\infty,\epsilon}$. Then $y_\infty$ is the limit of
a subsequence of points $y_i \in \partial M_{i,\epsilon}$. 

If 1. occurs for each $i$ then $\dist_0^2(x_\infty, y_\infty) \: 
R(y_\infty, 0) \: = \: \infty$, which is impossible. 
If 2. occurs for each $i$ then $R(y_\infty, 0) \: = \: \infty$, which
is impossible.
If 3. occurs for each $i$ then $R(y_\infty, 0) \: = \: 0$.
It follows from (\ref{harnackinequality}) that $M_\infty$ is flat,
which is impossible, as $R(x_\infty, 0) \: = \: 1$.

Hence $M_\infty \: = \: M_{\infty,\epsilon}$, i.e. no point in the
noncompact ancient solution $M_\infty$ is the center of an $\epsilon$-neck.
This contradicts the previous conclusion that 
$M_{\infty,\epsilon}$ is compact.
\end{proof}

\section{I.11.8. Necklike behavior at infinity of a three-dimensional
$\kappa$-solution - strong version}
\label{strongversion}

The following corollary is an application of the compactness result
Theorem \ref{thmI.11.7}.  
it is a refinement of 
\cite[Cor. I.11.8]{Perelman}.

\begin{corollary}
\label{neckstructure}
For all $\kappa>0$, there exists an $\eps_0>0$ such that 
for all $0<\eps<\eps_0$ there exists an $\al=\al(\eps,\kappa)$
with the property that for  any $\kappa$-solution $(M,g(\cdot))$,
and at any time $t$, precisely one of the following holds
($M_\eps$ denotes the set of points which are not centers of $\eps$-necks
at time $t$):

A.  $(M,g(\cdot))$ is round cylindrical flow, and 
so every point at every time is the center of an $\eps$-neck for all
$\eps>0$.

B.  $M$ is noncompact, $M_\eps\neq\emptyset$, and for
all $x,y\in M_\eps$, we have $R(x)d^2(x,y)<\al$.

C.  $M$ is compact, and there is a pair of points
$x,y\in M_\eps$ such that $R(x)d^2(x,y)>\al$,
\begin{equation}
M_\eps\subset B(x,\al R(x)^{-\frac{1}{2}})\cup 
B(y,\al R(y)^{-\frac{1}{2}}),
\end{equation}
and there is a minimizing geodesic $\ol{xy}$ such that
every $z\in M - M_\eps$ satisfies $R(z)d^2(z,\ol{xy})<\al$.

D. $M$ is compact and there exists a point $x\in M_\eps$
such that $R(x)d^2(x,z)<\al$ for all $z\in M$.
\end{corollary}

\begin{lemma}
\label{annoy}
For all $\eps>0$, $\kappa>0$,  there exists 
$\al=\al(\eps,\kappa)$  with the following 
property.  Suppose $(M,g(\cdot))$ is any 
$\kappa$-solution, $x,y,z\in M$, and 
at time $t$ we have $x, y \in M_\epsilon$
and $R(x)d^2(x,y)>\al$.  Then at time $t$ either 
$R(x)d^2(z,x)<\al$ or $R(y)d^2(z,y)<\al$ or ($R(z)d^2(z,\ol{xy})<\al$ and 
$z \notin M_\epsilon$).
\end{lemma}
\begin{proof}
Pick $\eps>0,\,\kappa>0$, and suppose no such $\al$ exists.  Then
there is a sequence $\al_k\ra\infty$, a  
sequence of $\kappa$-solutions $(M_k,g_k(\cdot))$,
and sequences $x_k,y_k,z_k\in M_k$, $t_k\in \R$ violating the $\al_k$-version
of the statement
for all $k$.  In particular, $x_k, y_k \in (M_k)_\epsilon$ and
\begin{equation}
R(x_k,t_k)d^2_{t_k}(x_k,y_k)\ra\infty,\; 
R(x_k,t_k)d^2_{t_k}(z_k,x_k)\ra\infty,\;
\mbox{and}\,R(y_k,t_k)d^2_{t_k}(z_k,y_k)\ra\infty.
\end{equation}

 Let $z_k'\in\ol{x_ky_k}$ be a point in $\ol{x_ky_k}$ 
nearest $z_k$ in $(M_k,g_k(t_k))$.

We first show that $R(x_k,t_k)d_{t_k}^2(z_k',x_k)\ra\infty$.
If not, we may pass to a subsequence on which $R(x_k,t_k)d_{t_k}^2(z_k',x_k)$
remains bounded.   Applying Theorem \ref{thmI.11.7}, we may
pass to a subsequence and 
rescale by $R(x_k,t_k)$,  
to make the sequence $(M_k,(x_k,t_k),g_k(\cdot))$
converge to a $\kappa$-solution $(M_\infty,(x_\infty,0),g_\infty(\cdot))$,
the segments $\ol{x_ky_k}\subset (M_k,g_k(t_k))$
converge to a ray $\ol{x_\infty\xi}\subset (M_\infty,g_\infty(0))$,
and the segments $\ol{z_k'z_k}$ converge to
a ray $\ol{z_\infty'\eta}$.  Recall that the comparison angle
$\cangle_{z^\prime_\infty}(u,v)$ tends to the Tits angle
$\tangle(\xi,\eta)$ as $u\in\ol{z_\infty'\xi}$, $v\in \ol{z_\infty'\eta}$
tend to infinity.  Since $d(z_k,z_k')=d(z_k,\ol{x_ky_k})$ we must
have $\tangle(\xi,\eta)\geq\frac{\pi}{2}$.  Now consider a sequence
$u_k\in \ol{z_\infty'\xi}$ tending to infinity.  By 
Theorem \ref{thmI.11.7}, part 5 of Corollary \ref{11.4cors},
and the remarks about Alexandrov spaces in Appendix \ref{Alexandrov},
if we rescale $(M_\infty,(u_k,0),g_\infty(\cdot))$ by $R(u_k,0)$,
we get round cylindrical flow as a limit.   When $k$ is sufficiently
large,
we may find an almost product region $D\subset (M_\infty,g_\infty(\cdot))$
containing $u_k$ which is disjoint from $\ol{z_\infty'\eta}$,
and whose cross-section $\Si\times\{0\}\subset\Si\times (-1,1)\simeq D$
intersects the ray $\ol{z_\infty'\xi}$ transversely at a single point.
This implies that $\Si\times \{0\}$ separates the two ends
of $\ol{z_\infty'\xi}\cup\ol{z_\infty'\eta}$ from each other; hence
$M_\infty$ is two-ended, and $(M_\infty,g_\infty(\cdot))$ is round
cylindrical flow.   This contradicts the assumption that
$x_k$ is not the center of
an $\eps$-neck.  Hence $R(x_k,t_k)d_{t_k}^2(z_k',x_k)\ra\infty$,
and similar reasoning shows that $R(y_k,t_k)d_{t_k}^2(z_k',y_k)\ra\infty$.

By part 5 of Corollary \ref{11.4cors}, we therefore have
$R(z_k',t_k)d_{t_k}^2(z_k',x_k)\ra\infty$ and 
$R(z_k',t_k)d_{t_k}^2(z_k',y_k)\ra\infty$.  Rescaling the sequence
$(M_k,(z_k',t_k),g_k(\cdot))$ by
$R(z_k',t_k)$, we get convergence to round cylindrical flow
(since any limit flow contains a line), and $\ol{z_k'z_k}$
subconverges to a segment orthogonal to the $\R$-factor, which
implies that $R(z_k',t_k)d_{t_k}^2(z_k,z_k')$ is bounded
and $z_k$ is the center of an $\eps$-neck for large $k$.  This contradicts
our assumption that the $\al_k$-version of the lemma is 
violated for each $k$.
\end{proof}

\bigskip
\noindent
{\em Proof of Corollary \ref{neckstructure}.}
Let $(M,g(\cdot))$
be a $\kappa$-solution, and $\eps>0$.

\medskip
{\em Case 1: Every $x\in (M,g(t))$
is the center of an $\eps$-neck.}  In this case, if $\eps>0$ is sufficiently
small, $M$ fibers over a 
$1$-manifold with fiber $S^2$.  If the 
$1$-manifold is homeomorphic to $\R$, then
$M$ has two ends, which implies that the flow
$(M,g(\cdot))$ is an evolving round cylinder.  If the 
base of the fibration were a circle, then the universal
cover $(\tilde M,\tilde g(t))$ would split off
a line, which would imply that the universal covering
flow would be a round cylindrical flow; but this
would violate the $\kappa$-noncollapsed assumption
at very negative times.  Thus A holds in this
case.

{\em Case 2: There exist $x,y\in M_\eps$ such that
$R(x)d^2(x,y)>\al$.}  By Lemma \ref{annoy} and Corollary \ref{11.4cors}
part 5, for all $z \in M - \left( B(x,\al R(x)^{-\frac{1}{2}})\cup 
B(y,\al R(y)^{-\frac{1}{2}}) \right)
$, we have $R(z)d^2(z,\ol{xy})<\alpha$ and $z \notin M_\epsilon$.
This implies (again by
Corollary \ref{11.4cors} part 5) that there exists a
$\gamma=\gamma(\eps,\kappa)$ such that for every $z\in M$
there is a $z'\in\ol{xy}$ for which  $R(z')d^2(z',z)<\gamma$,
which means that $M$ must be compact, and C holds.

{\em Case 3: $M_\eps\neq \emptyset$, and for all $x,y\in M_\eps$, 
we have $R(x)d^2(x,y)<\al$.}
If $M$ is noncompact then
we are in case B and are done, so assume that $M$ is compact.
Pick $x\in M_\eps$, and suppose  $z\in M$ maximizes
$R(x)d^2(\cdot,x)$.  
If $R(x)d^2(z,x)\geq\al$,  then $z$ is the center of an $\eps$-neck, and we 
may look at the cross-section $\Si$ of the neck region.
If $\Si$ separates $M$, then when $\eps>0$ is sufficiently
small, we  get a contradiction to the assumption
that $z$ maximizes $R(x)d^2(x,\cdot)$.  Hence
$\Si$ cannot separate $M$, and there is a loop
passing through $x$ which intersects $\Si$ transversely
at one point.  It follows that the universal covering flow
$(\tilde M,\tilde g(\cdot))$ is cylindrical flow, a contradiction.
Hence $R(x)d^2(x,z)<\al$ for all $z\in M$, so D holds.
\end{proof}

\section{More properties of $\kappa$-solutions}

In  this section we prove some additional properties of
$\kappa$-solutions. In particular, Corollary \ref{longsegment} implies
that if 
$z$
lies in a geodesic segment $\eta$ in a $\kappa$-solution $M$
and if the endpoints of $\eta$ are sufficiently far from 
$z$ (relative to $R(z)^{- \: \frac12}$) then $z \notin M_\epsilon$.
The results of this section will be used in the proof of Theorem 
\ref{thmI.12.1}.

\begin{proposition}
\label{neckcrit}
For all $\kappa>0$, $\alpha>0$, $\theta>0$, 
there exists a $\beta(\kappa,\alpha,\theta)<\infty$
such that if $(M,g(t))$ is a time slice of a $\kappa$-solution,
$x,y_1,y_2\in M$, $R(x)d^2(x,y_i)>\beta$ for $i=1,2$, and $\cangle_x(y_1,y_2)\geq \theta$,
then (a) $x$ is the center of
 an $\alpha$-neck, and (b) $\cangle_x(y_1,y_2)\geq \pi-
\alpha$.
\end{proposition}
\begin{proof}
The proof of this is similar to the first part of the proof of
Lemma \ref{annoy}.  Note that when $\alpha$ is small, then after enlarging
$\beta$ if necessary, the neck region around $x$ will separate 
$y_1$ from $y_2$; this implies (b).
\end{proof}

\begin{corollary} \label{longsegment}
For all $\kappa>0$, $\eps>0$, there exists a $\rho=\rho(\kappa,\eps)$ such that
if $(M,g(t))$ is a time slice of a $\kappa$-solution, $\eta\subset (M,g(t))$
is a minimizing geodesic segment with endpoints $y_1,y_2$, $z\in M$, $z'\in \eta$
is a point in $\eta$ nearest $z$, and $R(z')d^2(z',y_i)>\rho$
for $i=1,2$, then
$z,z'$ are centers of 
$\eps$-necks, and $\max(R(z)d^2(z,z'),R(z')d^2(z,z'))<4\pi^2$.
\end{corollary}
\begin{proof}
Pick $\eps'>0$.  Under the assumptions, if 
\begin{equation}
\min(R(z')d^2(z',y_1),R(z')d^2(z',y_2))
\end{equation}
is sufficiently large, we can apply the preceding proposition
to the triple $z',y_1,y_2$, to conclude that $z'$ is the center of an 
$\eps'$-neck.  Since the shortest segment from $z$ to $z'$
is orthogonal to $\eta$, when $\eps'$ is small enough 
the segment $\ol{zz'}$ will lie close to an $S^2$ cross-section
in the approximating round cylinder, which 
gives $R(z')d^2(z,z')\lesssim 2\pi^2$.
\end{proof}

\section{I.11.9. Getting a uniform value of $\kappa$}
\label{secI.11.9}

\begin{proposition} \label{propI.11.9}
There is a $\kappa_0 > 0$ so that if $(M, g(\cdot))$ is an oriented
three-dimensional $\kappa$-solution, for some
$\kappa > 0$, then it is a $\kappa_0$-solution or
it is a quotient of the round shrinking $S^3$.
\end{proposition}
\begin{proof}
Let $(M, g(\cdot))$ be a $\kappa$-solution.
Suppose that for some $\kappa^\prime>0$, 
the solution is $\kappa^\prime$-collapsed at some scale.  
After
rescaling, we can assume that there is a point $(x_0, 0)$ so that
$|\Rm(x, t)| \: \le \: 1$ for all $(x, t)$ satisfying
$\dist_0(x, x_0) \: < \: 1$ and $t \in [-1, 0]$, with
$\vol(B_0(x_0, 1)) \: < \: \kappa^\prime$.  Let $\widetilde{V}(t)$ denote
the reduced volume as a function of $t \in (- \infty, 0]$, as
defined using curves from $(x_0, 0)$. It is nondecreasing in $t$.
As in the proof of Theorem \ref{nolocalcollapse}, there is an estimate
$\widetilde{V}(-\kappa^\prime) \: \le \: 3 (\kappa^\prime)^{3/2}$. Take a sequence
of times $t_i \rightarrow - \infty$. For each $t_i$, choose $q_i \in M$
so that $l(q_i, t_i) \: \le \: \frac{3}{2}$. From the proof of Proposition
\ref{asympsoliton}, for all $\epsilon > 0$
there is a $\delta > 0$ such that $l(q, t)$
does not exceed $\delta^{-1}$ whenever $t \in [t_i, t_i/2]$ and
$\dist_{t_i}^2(q, q_i) \: \le \: \epsilon^{-1} t_i$. Given the
monotonicity of $\widetilde{V}$ and the upper bound on $l(q,t)$,
we obtain an upper bound on the volume of the time-$t_i$ ball
$B(q_i, \sqrt{t_i/\epsilon})$ of the form $\const \: t_i^{3/2} \:
e^{\delta^{-1}} \: (\kappa^\prime)^{3/2}$. 

On the other hand, from
Proposition \ref{asympsoliton}, a subsequence of the rescalings
of the ancient solution around $(q_i, t_i)$ converges to a nonflat
gradient shrinking soliton.  If the gradient shrinking soliton is compact
then it must be a quotient of the round shrinking $S^3$
\cite{Hamiltonnnnn}. Otherwise,
Corollary \ref{cornosoliton} says that if the gradient shrinking soliton is
noncompact then it must be an evolving cylinder or its
$\Z_2$-quotient.  Fixing $\epsilon$, this gives a lower bound on
$\vol(B_{t_i} (q_i, \sqrt{t_i/\epsilon}))$ in terms of the noncollapsing
constants of the evolving cylinder and its $\Z_2$-quotient. Hence there
is a universal constant $\kappa_0$ so that if
$\kappa^\prime<\kappa_0$ then we obtain a contradiction to the
assumption of $\kappa^\prime$-collapsing.
\end{proof}
\begin{remark} The hypotheses of Corollary \ref{cornosoliton} 
assume a global upper bound on the sectional curvature of any time slice, 
which in the $n$-dimensional case
is not {\it a priori} true for the asymptotic soliton of
Proposition \ref{asympsoliton}.  However, in our $3$-dimensional case, 
the argument of Theorem \ref{thmI.11.7} shows that there is such an upper bound.
\end{remark}

\section{II.1.2. Three-dimensional noncompact $\kappa$-noncollapsed
gradient shrinkers are standard}
\label{II.1.2}

In this section we show that any complete oriented 
$3$-dimensional noncompact 
$\kappa$-noncollapsed 
gradient shrinking soliton with
bounded nonnegative curvature
is either the evolving round cylinder $\R \times S^2$ or its
$\Z_2$-quotient.

The basic example of a gradient shrinking soliton is the
metric on $\R \times S^2$ which gives the $2$-sphere
a radius of $\sqrt{-2t}$ at time $t \in (-\infty, 0)$. With
coordinates $(s, \theta)$ on $\R \times S^2$, the function $f$ is given by
$f(t, s, \theta) \: = \: - \: \frac{s^2}{4t}$.

\begin{lemma} \label{nosoliton} (cf. Lemma of II.1.2)
There is no complete oriented $3$-dimensional noncompact
$\kappa$-noncollapsed gradient shrinking soliton with 
bounded positive sectional curvature.
\end{lemma}
\begin{proof}
The idea of the proof is to show that the soliton has the qualitative
features of a shrinking cylinder, and then to get a contradiction to
the assumption of positive sectional curvature.

Applying $\nabla_i$ to the gradient shrinker equation
\begin{equation} \label{IIeqn}
\nabla_i \nabla_j f \: + \: R_{ij} \: + \: \frac{1}{2t} \: g_{ij}
\: = \: 0
\end{equation}
gives
\begin{equation}
\triangle \nabla_j f \: + \: \nabla_i R_{ij} \: = \: 0.
\end{equation}
As $\nabla_i R_{ij} \: = \: \frac{1}{2} \: \nabla_j R$ and
$\triangle \nabla_j f \: = \: \nabla_j \triangle f \: + \:
R_{jk} \nabla_k f \: = \: \nabla_j \left( -R-\frac{n}{2t} \right) \: + \: 
R_{jk} \nabla_k f$, 
we obtain
\begin{equation} \label{gradR}
\nabla_i R \: = \: 2 \: R_{ij} \: \nabla_j f.
\end{equation}

Fix a basepoint $x_0 \in M$ and consider a normalized minimal geodesic
$\gamma \: : \: [0, \overline{s}] \rightarrow M$ in the time $-1$ slice
with $\gamma(0) \: = \: x_0$. Put $X(s) \: = \: \frac{d\gamma}{ds}$.
As in the proof of Lemma \ref{distancedist},
$\int_0^{\overline{s}} \Ric(X, X) \: ds \: \le \:
\const$ for some constant independent of $\overline{s}$.
If $\{Y_i\}_{i=1}^3$ are orthonormal parallel vector fields along $\gamma$
then
\begin{equation}
\left( \int_0^{\overline{s}} |\Ric(X, Y_1)| \: ds \right)^2 \: \le \:
\overline{s} \: \int_0^{\overline{s}} |\Ric(X, Y_1)|^2 \: ds \: \le \:
\overline{s} \: \sum_{i=1}^3 \int_0^{\overline{s}} |\Ric(X, Y_i)|^2 \: ds.
\end{equation}
Thinking of $\Ric$ as a self-adjoint linear operator on $TM$, 
$\sum_{i=1}^3 |\Ric(X, Y_i)|^2 \: = \: \langle X, \Ric^2 X \rangle$.
In terms of a pointwise orthonormal frame  $\{e_i\}$
of eigenvectors of $\Ric$,
with eigenvalues $\lambda_i$, write $X \: = \: \sum_{i=1}^3 X_i e_i$.
Then 
\begin{equation}
\langle X, \Ric^2 X \rangle \: = \: \sum_{i=1}^3 \lambda_i^2 \: X_i^2 
\: \le \: (\sum_{i=1}^3 \lambda_i) \: (\sum_{i=1}^3 \lambda_i \: X_i^2) 
\: = \: R \cdot \Ric(X, X). 
\end{equation}
Hence
\begin{equation}
\left( \int_0^{\overline{s}} |\Ric(X, Y_1)| \: ds \right)^2 \: \le \: 
(\sup_M R) \: \overline{s} \: 
\int_0^{\overline{s}} \Ric(X, X) \: ds \: \le \:
\const \: \overline{s}.
\end{equation}

Multiplying (\ref{IIeqn}) by $X^i X^j$ and summing gives
$\frac{d^2f(\gamma(s))}{ds^2} \: + \: \Ric(X, X) \: - \: \frac{1}{2}
\: =\: 0$. Then 
\begin{equation}
\frac{d f(\gamma(s))}{ds} \Big|_{s = \overline{s}} \: = \:
\frac{d f(\gamma(s))}{ds} \Big|_{s = 0} \: + \: 
\frac{1}{2} \: \overline{s} \: - \: \int_0^{\overline{s}} \Ric(X, X) \: ds
\: \ge \: \frac{1}{2} \: \overline{s} \: - \: \const
\end{equation}
This implies that there is a compact subset of $M$ outside of which
$f$ has no critical points.

If $Y$ is a unit vector field perpendicular to $X$ then
multiplying (\ref{IIeqn}) by $X^i Y^j$ and summing gives
$\frac{d}{ds} (Y \cdot f)(\gamma(s)) \: + \: \Ric(X, Y)
\: =\: 0$. Then 
\begin{equation}
(Y \cdot f)(\gamma(\overline{s})) \: = \: 
(Y \cdot f)(\gamma(0)) \: - \: 
\int_0^{\overline{s}} \: \Ric(X, Y) \: ds
\end{equation}
and 
\begin{equation}
|(Y \cdot f)(\gamma(\overline{s}))| \: \le \:
\const (\sqrt{\overline{s}} \: +\: 1).
\end{equation}
For large $\overline{s}$, $|(Y \cdot f)(\gamma(\overline{s}))|$
is small compared to $(X \cdot f)(\gamma(\overline{s}))$. This
means that as one approaches infinity, the gradient of $f$ becomes
more and more parallel to the gradient of the distance function
from $x_0$, where by the latter we mean the vectors $X$ that are
tangent to minimal geodesics.

The gradient flow of $f$ is given by the equation
\begin{equation}
\frac{dx}{du} \: =\: (\nabla f)(x). 
\end{equation}
Then along a flowline, equation (\ref{gradR}) implies that
\begin{equation}
\frac{dR(x)}{du} \: = \: \left\langle \nabla R, \frac{dx}{du} \right\rangle
\: = \: 2 \Ric(\nabla f, \nabla f).
\end{equation}
In particular, outside of a compact set, $R$ is strictly increasing
along the flowlines. Put $\overline{R} \: = \: \limsup_{x \rightarrow \infty} 
R$.
Take points $x_\alpha$ tending toward infinity, with
$R(x_\alpha) \rightarrow \overline{R}$. Putting
${r}_\alpha \: = \: \sqrt{\dist_{-1}(x_0, x_\alpha)}$, we have
$\frac{{r}_\alpha}{\dist_{-1}(x_0, x_\alpha)} \rightarrow 0$ and
$R(x_\alpha) \: {r}_\alpha^2 \rightarrow \infty$.
Then the argument of the proof of Proposition \ref{propI.11.4}
shows that any convergent
subsequence of the rescalings around $(x_\alpha, -1)$ splits off
a line.  Hence the limit is a shrinking round cylinder with
scalar curvature $\overline{R}$ at time $-1$. Because our original
solution exists up to time zero, we must have 
$\overline{R} \le 1$.
Now equation
(\ref{shrinkingflow}) says that the Ricci flow
is given by $g(-t) \: = \: - \: t \: \phi_t^* g(-1)$, where
$\phi_t$ is the flow generated by $\nabla f$. 
It follows that $\inf_{x \in M} R(x,t) \: = \: Ct^{-1}$ for some $C > 0$.
That is, the curvature blows up uniformly as $t \rightarrow 0$.
Comparing this with the
singularity time of the shrinking round cylinder implies that
$\overline{R} = 1$. Performing a similar argument with any
sequence of $x_\alpha$'s tending toward infinity, with the property
that $R(x_\alpha)$ has a limit, shows that 
$\lim_{x \rightarrow \infty} R(x) \: = 1$.

Let $N$ denote a (connected component of a) 
level surface of $f$. At a point of $N$,
choose an orthonormal frame
$\{e_1, e_2, e_3\}$ with $e_3 \: = \: X$ normal to $N$.
From the Gauss-Codazzi equation, 
\begin{equation}
R^N \: = \: 2 \: K^N(e_1, e_2)
\: = \: 2 (K^M(e_1, e_2) \: + \: \det(S)),
\end{equation}
where $S$ is the shape operator. As
$R \: = \: 2 (K^M(e_1, e_2) \: + \: K^M(e_1, e_3) \: + \:
K^M(e_2, e_3))$ and
$\Ric(X, X) \: = \: K^M(e_1, e_3) \: + \:
K^M(e_2, e_3)$, we obtain
\begin{equation}
R^N \: = \: R \: - \: 2 \Ric(X, X) \: + \: 2 \det(S).
\end{equation}
The shape operator is given by 
$S \: = \: \frac{Hess f |_{TN}}{|\nabla f|}$. 
From (\ref{IIeqn}),
$\Hess f \: = \: \frac{1}{2} \: - \: \Ric$. We can
diagonalize $\Ric \Big|_{TN}$ to write
$\Ric \: = \:
\begin{pmatrix}
r_1 & 0 & c_1 \\
0 & r_2 & c_2 \\
c_1 & c_2 & r_3
\end{pmatrix}$, where $r_3 \: = \: \Ric(X, X)$. Then
\begin{align}
\det \left( \Hess f \Big|_{TN} \right) \: & = \: 
\left( \frac{1}{2} \: - \: r_1 \right) \left( \frac{1}{2} \: - \: r_2 \right)
\: = \: 
\frac{1}{4} \: \left( (1 - r_1 - r_2)^2 \: - \: (r_1 - r_2)^2 \right)
\\
& \le \: \frac{1}{4} \: (1 - r_1 - r_2)^2
\:  = \:
\frac{1}{4} \: (1 - R + \Ric(X, X))^2. \notag
\end{align}
This shows that the scalar curvature of $N$ is bounded above by
\begin{equation}
R \: - \: 2 \Ric(X, X) \: + \: \frac{(1-R+\Ric(X, X))^2}{2|\nabla f|^2}.
\end{equation}

If $|\nabla f|$ is large then
$1 \: - \: R \: + \: \Ric(X, X) \: < \: 2 \: |\nabla f|^2$. 
As $1 \: - \: R \: + \: \Ric(X, X)$ is positive when the
distance from $x$ to $x_0$ is large enough,
\begin{align}
\left( 1 \: - \: R \: + \: \Ric(X, X) \right)^2 \: & < \:
2 (1 \: - \: R \: + \: \Ric(X, X)) |\nabla f|^2 \\
& \le \:
2 (1 \: - \: R \: + \: \Ric(X, X)) |\nabla f|^2 \: + \:
2 |\nabla f|^2 \Ric(X, X) \notag
\end{align}
and so
\begin{equation}
\frac{\left( 1 \: - \: R \: + \: \Ric(X, X) \right)^2}{2|\nabla f|^2}
 \: < \:
1 \: - \: R \: + \: 2 \Ric(X, X).
\end{equation}
Hence
\begin{equation}
R \: - \: 2 \Ric(X, X) \: + \: \frac{(1-R+\Ric(X, X))^2}{2|\nabla f|^2}
\: < \: 1.
\end{equation}
This shows that $R^N < 1$ if $N$ is sufficiently far from $x_0$.

If $Y$ is a unit vector that is tangential to $N$ then from (\ref{IIeqn}),
\begin{equation}
\nabla_Y \nabla_Y f \: = \: \frac12 \: - \: \Ric(Y,Y).
\end{equation}
If $\{Y, Z, W\}$ is an orthonormal basis then
\begin{equation}
\Ric(Y, Y) \: = \: K^M(Y,Z) \: + \: K^M(Y, W) \: \le \:
K^M(Y,Z) \: + \: K^M(Y, W) \: + \: K^M(Z, W) \: = \: \frac12 \: R.
\end{equation}
Hence $\nabla_Y \nabla_Y f \: \ge \: \frac12 \: (1-R)$, which is
positive if $N$ is sufficiently far from $x_0$. Thus $N$ is convex
and so 
the area of the level set increases as the level increases.
On the other hand, 
we can take points $x_\alpha$ on the level sets going to
infinity, apply the previous splitting argument 
and use the fact that $\grad f$ becomes almost parallel to
$\grad d(\cdot, x_0)$. Within one of the approximate cylinders coming from the
splitting argument, there is a projection $\pi$ to its base $S^2$. As
a tangent plane of $N$ is almost perpendicular to $\Ker(d\pi)$,
the restriction of $\pi$ to $N$ is an almost-isometry from
$N$ to $S^2$.
By the monotonicity of $\area(N)$, we conclude that
$\area(N) \: \le \: 8\pi$ if $N$ is sufficiently far from $x_0$.
However as $N$ is a topologically a $2$-sphere, 
the Gauss-Bonnet theorem says that
$\int_N R^N \: dA \: = \: 8\pi$. This contradicts the facts that
$R^N \: < 1$ and $\area(N) \: \le \: 8\pi$.
\end{proof}

\begin{corollary} \label{cornosoliton}
The only complete oriented $3$-dimensional noncompact
$\kappa$-noncollapsed gradient shrinking solitons with 
bounded nonnegative sectional curvature are the round
evolving $\R \times S^2$ and its $\Z_2$-quotient
$\R \times_{\Z_2} S^2$.
\end{corollary}
\begin{proof}
Let $(M, g(\cdot))$ be a complete oriented $3$-dimensional noncompact
$\kappa$-noncollapsed gradient shrinking soliton with 
bounded nonnegative sectional curvature. By Lemma \ref{nosoliton}, $M$ cannot
have positive sectional curvature.  From Theorem \ref{splitting},
$M$ must locally split off an $\R$-factor.  Then the universal cover splits
off an $\R$-factor and so, by Corollary \ref{2Dsoliton} or Section
\ref{altpf11.3}, must be the standard $\R \times S^2$. From the 
$\kappa$-noncollapsing, $M$ must be $\R \times S^2$ or $\R \times_{\Z_2} S^2$.
\end{proof}

\section{I.12.1. Canonical neighborhood theorem}
\label{I.12.1}

In this section we show that a high-curvature region of a three-dimensional
Ricci flow is modeled by part of a $\kappa$-solution.

We first define the notion of $\Phi$-almost nonnegative curvature.

\begin{definition} (cf. I.12) \label{pinchingdef}
Let $\Phi \in C^\infty(\R)$ be a positive nondecreasing function 
such that
for positive $s$, $\frac{\Phi(s)}{s}$ is a decreasing function which
tends to zero as $s \rightarrow \infty$.
A Ricci flow solution is said to have {\em $\Phi$-almost nonnegative
curvature} if for all $(x,t)$, we have
\begin{equation} \label{phieqnn}
\Rm(x,t) \: \ge \: - \: \Phi(R(x,t)).
\end{equation}
\end{definition}

\begin{remark}
Note that $\Phi$-almost nonnegative curvature implies that the scalar curvature is uniformly
bounded below by $- \: 6 \: \Phi(0)$. 
The formulation of the pinching condition in \cite[Section 12]{Perelman}
is that there is a decreasing function $\phi$, tending to zero at infinity,
so that $\Rm(x,t) \: \ge \: - \: \phi(R(x,t)) \: R(x,t)$ for each
$(x,t)$. This formulation has a problem when $R(x,t) < 0$, if one
takes $\phi$ to be defined on all of $\R$. The condition in
Definition \ref{pinchingdef} is what comes out of the
three-dimensional
Hamilton-Ivey pinching result (applied to the rescaled metric
$\tilde{g}(t) = \frac{g(t)}{t}$) if we assume normalized initial
conditions; see Appendix \ref{phiappendix}.
\end{remark}

We note that since the sectional curvatures have to add up to
$R$, the lower bound (\ref{phieqnn}) implies a double-sided bound on the
sectional curvatures.  Namely,
\begin{equation} \label{double}
- \Phi(R) \: \le \:
\Rm \: \le \: \frac{R}{2} \: + \: \left( \frac{n(n-1)}{2} - 1 \right)
\: \Phi(R).
\end{equation}

The main use of the pinching condition is to show that
blowup limits have nonnegative sectional curvature.

\begin{lemma} \label{almost nonneg}
Let $\{(M_k, p_k, {g}_k)\}_{k=1}^\infty$ be a sequence of
complete pointed Riemannian manifolds with $\Phi$-almost
nonnegative curvature. 
Given a sequence $Q_k \rightarrow \infty$, put
$\overline{g}_k = Q_k g_k$ and
suppose that there
is a pointed smooth limit
$(M_\infty, p_\infty, \overline{g}_\infty) = \lim_{k \rightarrow \infty}
(M_k, p_k, \overline{g}_k)$. Then
$M_\infty$ has nonnegative sectional curvature.
\end{lemma}
\begin{proof}
First, the $\Phi$-almost
nonnegative curvature condition implies that the scalar
curvature of $M_k$ is bounded below uniformly in $k$.
For $m \in M_\infty$, let $m_k \in (M_k,\bar g_k)$ be
a sequence of approximants to $m$.  Then
$\lim_{k \rightarrow \infty} \overline{\Rm}(m_k) \: = \:  
\overline{\Rm}_\infty(m)$, where $\overline{\Rm}(m_k) \: = \:
Q_k^{-1} \: {\Rm}(m_k)$.
There are two possibilities : either the numbers $R(m_k)$ are
uniformly bounded above or they are not.  If they are uniformly
bounded above then (\ref{double}) implies that
$\Rm(m_k)$ is uniformly bounded above and below,
so $\overline{\Rm}_\infty(m)
= 0$. Suppose on the other hand that a subsequence of the numbers 
$R(m_k)$ tends to infinity.  We pass to this subsequence. Now
$\overline{R}_\infty(m) = \lim_{k \rightarrow \infty} \overline{R}(m_k)$
exists by assumption and is nonnegative.
Applying (\ref{phieqnn})
gives that $\overline{\Rm}_\infty(m)$, the limit of 
\begin{equation}
Q_k^{-1} \: {\Rm}(m_k) \: = \: 
\overline{R}(m_k) \: \frac{{\Rm}(m_k)}{R(m_k)},
\end{equation}
is nonnegative.
\end{proof}

We now prove the first version of the ``canonical neighborhood''
theorem.

\begin{theorem} (cf. Theorem I.12.1) \label{thmI.12.1}
Given $\epsilon, \kappa, \sigma > 0$ and a function $\Phi$ as above, one can
find $r_0 > 0$ with the following property.  Let $g(\cdot)$ be a Ricci flow
solution on a three-manifold $M$, defined for $0 \le t \le T$ with
$T \ge 1$. We suppose that for each $t$, $g(t)$ is complete, 
and the sectional curvature is bounded on compact time intervals.
Suppose that the Ricci flow has
$\Phi$-almost nonnegative curvature and is $\kappa$-noncollapsed
on scales less than $\sigma$. Then for any point $(x_0, t_0)$ with $t_0 \ge 1$
and $Q = R(x_0, t_0) \ge r_0^{-2}$, the solution in
$\{(x,t) \: : \: \dist_{t_0}^2(x, x_0) < (\epsilon Q)^{-1}, t_0 - (\epsilon Q)^{-1}
\le t \le t_0\}$ is, after scaling by the factor $Q$, $\epsilon$-close to the
corresponding subset of a $\kappa$-solution.
\end{theorem}

\begin{remark} \label{wangg}
Our statement of Theorem \ref{thmI.12.1} differs slightly from that
in \cite[Theorem 12.1]{Perelman}. 
First, we allow $M$ to be noncompact, provided that 
there is bounded sectional curvature on compact time intervals.
This generalization will be useful for later work.
More importantly, the statement in \cite[Theorem 12.1]{Perelman}
has noncollapsing at scales less than $r_0$, whereas we require
noncollapsing at scales less than $\sigma$. See Remark \ref{wang} for further
comment.  

In the phrase ``$t_0 \ge 1$'' there is an implied scale which comes
from the $\Phi$-almost nonnegativity assumption, and similarly for the
statement  ``scales less than $\sigma$''.
\end{remark}
\begin{proof}
We give a proof which differs in some points from the proof in 
\cite{Perelman} but which has the same ingredients. We first outline
the argument.

Suppose that the theorem is false.  Then for some $\epsilon, \kappa,\sigma>0$,
we have a sequence of such $\Phi$-nonnegatively curved 
$3$-dimensional Ricci flows $(M_k,g_k(\cdot))$
defined on intervals $[0,T_k]$, 
and sequences $r_k \rightarrow 0$, $\hat{x}_k\in M_k$, $\hat{t}_k\geq 1$ such
that $M_k$ is $\kappa$-noncollapsed on scales $< \sigma$ and
$Q_k = R(\hat{x}_k,\hat{t}_k) \ge r_k^{-2}$, but the $Q_k$-rescaled solution
in $B_{(\epsilon Q_k)^{-1/2}}(\hat{x}_k) \times
[\hat{t}_k - (\epsilon Q_k)^{-1}, \hat{t}_k]$ is not $\epsilon$-close to the
corresponding subset of any $\kappa$-solution.

We note that if the statement were false for $\epsilon$ then it would
also be false for any smaller $\epsilon$. Because of
this, somewhat paradoxically, we will begin the argument with a
given $\epsilon$ but will allow ourselves to make $\epsilon$ small
enough later so that the argument works. 
To be clear, we will eventually get a contradiction using a fixed
(small) value of $\epsilon$, but as the proof goes along we will
impose some upper bounds on this value in order for the proof to work. 
(If we tried to list all of
the constraints at the beginning of the argument then they would
look unmotivated.)

The goal is to get a contradiction based on the ``bad'' points
$(\widehat{x}_k, \widehat{t}_k)$.
In a sense, the method of proof of Theorem \ref{thmI.12.1} is an induction on
the curvature scale. For example, if we were to make the additional
assumption in the theorem 
 that $R(x, t) \le R(x_0, t_0)$ for all $x \in M$ and $t \le t_0$
then the theorem would be very easy to prove. We would just take
a convergent subsequence of the rescaled solutions, based at 
$(\hat{x}_k, \hat{t}_k)$, to get a $\kappa$-solution; this would
give a contradiction. This simple argument
can be considered to be the first step
in a proof by induction on curvature scale. In the proof of
Theorem \ref{thmI.12.1} one effectively proves the result at a
given curvature scale inductively by assuming that the result
is true at higher curvature scales.

The actual proof consists of four steps. 
Step 1 consists of replacing the sequence $(\widehat{x}_k, \widehat{t}_k)$
by another sequence of ``bad'' points $(x_k, t_k)$ which have the
property that points near
$(x_k, t_k)$ with distinctly higher scalar curvature are ``good'' points.
It then suffices to get a contradiction based on the existence of the
sequence $(x_k, t_k)$. 

In steps 2-4 one uses the points $(x_k, t_k)$
to build up a $\kappa$-solution, whose existence then contradicts the
``badness'' of the points $(x_k, t_k)$. 
More precisely, let
$(M_k,(x_k,t_k),\bar g_k(\cdot))$ be the result of rescaling $g_k(\cdot)$
by $R(x_k,t_k)$.
We will show that the sequence of pointed
flows $(M_k,(x_k,t_k),\bar g_k(\cdot))$ accumulates on a $\kappa$-solution
$(M_\infty, (x_\infty, t_0), \overline{g}_\infty(\cdot))$,
thereby obtaining a contradiction. 

In step 2 one takes a pointed limit of the manifolds
$(M_k, x_k,\bar g_k(t_k))$ in order to construct what
will become the final time slice of the $\kappa$-solution,
$(M_\infty, x_\infty, \overline{g}_\infty(t_0))$. In order to take this
limit, it is necessary to show that the manifolds 
$(M_k, x_k,\bar g_k(t_k))$ have uniformly bounded curvature on distance
balls of a fixed radius.  If this were not true then for some radius, 
a subsequence of the
manifolds $(M_k, x_k,\bar g_k(t_k))$ would have curvatures that
asymptotically blowup on the ball of that radius.  One shows that
geometrically, the curvature blowup is due to the asymptotic formation of 
a cone-like point at the blowup radius.  Doing a further
rescaling at this cone-like point, one obtains a Ricci flow solution
that ends on a part of a nonflat metric cone. This gives a contradiction
as in the case $0 < {\mathcal R} < \infty$ of Theorem \ref{asympscalar}.

Thus one can construct the pointed limit
$(M_\infty, x_\infty, \overline{g}_\infty(t_0))$. The goal now is
to show that $(M_\infty, x_\infty, \overline{g}_\infty(t_0))$ is the
final time slice of a $\kappa$-solution 
$(M_\infty, (x_\infty, t_0), \overline{g}_\infty(\cdot))$. In step 3
one shows that $(M_\infty, x_\infty, \overline{g}_\infty(t_0))$ extends
backward
to a Ricci flow solution on some time interval $[t_0 - \Delta, t_0]$, 
and that
the time slices have bounded nonnegative curvature.  In step 4 one
shows that the Ricci flow solution can be extended all the way to time 
$(-\infty, t_0]$, thereby constructing a $\kappa$-solution

\bigskip
{\em Step 1: Adjusting the choice of basepoints.}

We first modify the points 
$(\widehat{x}_k, \widehat{t}_k)$ slightly in time
to points $(x_k, t_k)$
so that in the given Ricci flow solution, there are no other ``bad'' points with
much larger scalar curvature in
a earlier time interval whose length is large compared to
$R(x_k, t_k)^{-1}$. 
(The phrase ``nearly the smallest curvature $Q$'' in 
\cite[Proof of Theorem 12.1]{Perelman} should read
``nearly the largest curvature $Q$''. This is clear from
the sentence in parentheses that follows.) The proof of
the next lemma is by pointpicking, as in Appendix \ref{pointpicking}. 

\begin{lemma}  \label{ptpicking}
We can find $H_k\ra\infty$, $x_k\in M_k$, and $t_k\geq\frac{1}{2}$
such that $Q_k =  R(x_k,t_k)\ra\infty$, and for all $k$
the conclusion of Theorem \ref{thmI.12.1} fails at $(x_k,t_k)$, 
but holds for
any $(y,t)\in M_k\times[t_k-H_kQ_k^{-1},t_k]$
for which $R(y,t)\geq 2R(x_k,t_k)$.
\end{lemma}
\begin{proof}
Choose $H_k\ra\infty$ such that 
$H_k(R(\widehat{x}_k,\widehat{t}_k))^{-1}\leq \frac{1}{10}$
for all $k$.  For each $k$, initially set $(x_k,t_k) = (\hat x_k,\hat t_k)$.
Put $Q_k = R(x_k, t_k)$ and look for a point in 
$M_k\times[t_k-H_kQ_k^{-1},t_k]$ at which 
Theorem \ref{thmI.12.1} fails,
and the scalar curvature is at least $2R(x_k,t_k)$. 
If such a point
exists, replace $(x_k,t_k)$  by this point; otherwise do nothing.
Repeat this until the second alternative occurs.
This process must terminate with a new choice of $(x_k,t_k)$
satisfying the lemma.
\end{proof}

Hereafter we use this modified sequence $(x_k, t_k)$. 
Let $(M_k,(x_k,t_k),\bar g_k(\cdot))$ be the result of rescaling $g_k(\cdot)$
by $Q_k=R(x_k,t_k)$. We use $\bar R_k$ to denote its scalar curvature;
in particular, $\bar R_k (x_k, t_k) = 1$.
Note that the rescaled time interval of Lemma \ref{ptpicking}
has duration $H_k \rightarrow \infty$;
this is what we want in order to try to extract an ancient solution.

\bigskip
{\em Step 2: For every $\rho<\infty$, the scalar curvature $\bar R_k$ 
is uniformly bounded 
on the $\rho$-balls $B(x_k,\rho)\subset (M_k,\bar g_k(t_k))$
(the argument for this is essentially equivalent to
\cite[Pf. of Claim 2 of Theorem I.12.1]{Perelman}).}
Before proceeding, we need some bounds which come from our choice
of basepoints, and the derivative bounds inherited (by approximation)
from $\kappa$-solutions.

\begin{lemma} \label{grad}
There is a constant $C = C(\kappa)$ so that for any $(x, t)$ in a
Ricci flow solution, if $R(x, t) > 0$ and the solution in
$B_t(x, (\epsilon R(x, t))^{-1/2}) \times [t-(\epsilon R(x, t))^{-1}, t]$
is $\epsilon$-close to a corresponding subset of a $\kappa$-solution
then $|\nabla R^{-1/2}|(x, t) \: \le \: C$ and
$|\partial_t R^{-1}|(x,t) \: \le \: C$.
\end{lemma}
\begin{proof}
This follows from the compactness in Theorem \ref{thmI.11.7}.
\end{proof}

Note that the same value of $C$ in Lemma \ref{grad}
also works for smaller $\epsilon$.
 
\begin{lemma} (cf. Claim 1 of I.12.1)
\label{claim1} For each $(\overline{x}, \overline{t})$ with
$t_k \: - \: \frac12 \: H_k \: Q_k^{-1} \le \overline{t} \le t_k$, we have
$R_k(x, t) \le 4 \overline{Q}_k$ whenever 
$\overline{t} \: - \: c \: \overline{Q}_k^{-1} \le t \le \overline{t}$ and
$\dist_{\overline{t}}(x, \overline{x}) \: \le \: c \: \overline{Q}_k^{-1/2}$,
where $\overline{Q}_k = Q_k + |R_k(\overline{x}, \overline{t})|$ and
$c = c(\kappa) > 0$ is a small constant. 
\end{lemma}
\begin{proof}
If $R_k(x, t) \le 2Q_k$ then there is nothing to show.  If
$R_k(x, t) > 2Q_k$, consider a spacetime curve $\gamma$ that goes linearly from
$(x,t)$ to $(x, \overline{t})$, and then goes from
$(x, \overline{t})$ to $(\overline{x}, \overline{t})$ along a
minimizing geodesic.  If there is a point on $\gamma$ with curvature
$2Q_k$, let $p$ be the nearest such point to $(x, t)$. If not, put
$p = (\overline{x}, \overline{t})$. From the conclusion of Lemma
\ref{ptpicking}, we can apply 
Lemma \ref{grad} along $\gamma$ from $(x, t)$ to $p$.
The claim follows from integrating the ensuing 
derivative bounds along $\gamma$.
\end{proof}

\begin{lemma} \label{rescaled}
In terms of the rescaled solution $\overline{g}_k(\cdot)$,
for each $(\overline{x}, \overline{t})$ with
$t_k \: - \: \frac12 \: H_k \le \overline{t} \le t_k$, we have
$\overline{R}_k(x, t) \le 4 \widetilde{Q}_k$ whenever 
$\overline{t} \: - \: c \: \widetilde{Q}_k^{-1} \le t \le \overline{t}$ and
$\dist_{\overline{t}}(x, \overline{x}) \: \le \: c \: \widetilde{Q}_k^{-1/2}$,
where $\widetilde{Q}_k = 1 + |\overline{R}_k(\overline{x}, \overline{t})|$.
\end{lemma}
\begin{proof}
This is just the rescaled version of Lemma \ref{claim1}.
\end{proof}

For all $\rho\geq 0$, put 
\begin{equation}
D(\rho) = \sup\{\bar R_k(x,t_k)\mid k\geq 1,\;x\in B(x_k,\rho)\subset
(M_k,\bar g_k(t_k))\},
\end{equation}
and let $\rho_0$ be the supremum of the $\rho$'s for which
$D(\rho)<\infty$.  Note that $\rho_0>0$, in view of Lemma \ref{rescaled}
(taking $(\bar x,\bar t)=(x_k,t_k)$).
Suppose that $\rho_0<\infty$.         
After passing to a subsequence if necessary, we can find
a sequence $y_k\in M_k$ with $\dist_{t_k}(x_k,y_k) \rightarrow \rho_0$
and $\bar R(y_k,t_k)\ra\infty$.    Let $\eta_k\subset (M_k,\bar g_k(t_k))$
be a minimizing geodesic segment from $x_k$ to $y_k$. 
Let $z_k\in \eta_k$ be the point on $\eta_k$ closest to $y_k$
at which $\bar R(z_k,t_k)=2$, and let $\ga_k$
be the subsegment of $\eta_k$ running from $y_k$ to $z_k$.  By
Lemma \ref{rescaled} the length of $\ga_k$ is bounded away from zero
independent of $k$.  Due to the $\Phi$-pinching
(see (\ref{double})),
for all $\rho<\rho_0$, we have a uniform bound on $|\Rm|$
on the balls $B(x_k,\rho)\subset (M_k,\bar g_k(t_k))$.  The 
injectivity radius is also controlled in $B(x_k,\rho)$, 
in view of the curvature bounds and the $\kappa$-noncollapsing.
Therefore after passing
to a subsequence, we can assume that the pointed sequence 
$(B(x_k,\rho_0),\bar g_k(t_k),x_k)$ converges in the pointed 
Gromov-Hausdorff topology (i.e. for all $\rho<\rho_0$ we have the 
usual Gromov-Hausdorff convergence) to a pointed $C^1$-Riemannian manifold 
$(Z,\bar g_\infty,x_\infty)$,  the segments $\eta_k$ converge to a  segment (missing
an endpoint)
$\eta_\infty\subset Z$ emanating from $x_\infty$, and $\gamma_k$
converges to  $\ga_\infty\subset\eta_\infty$.  Let $\bar Z$
denote the completion of $(Z,\bar g_\infty)$, and $y_\infty\in \bar Z$
the limit point of $\eta_\infty$.    Note that by
Lemma \ref{ptpicking} and
part 4 of Corollary \ref{11.4cors}, the Riemannian structure near
$\ga_\infty$ may be chosen to be many times differentiable.
(Alternatively, this follows from Lemma \ref{rescaled} and the Shi
estimates of Appendix \ref{applocalder}.)
In particular
the scalar curvature $\bar R_\infty$ is defined, differentiable,
and satisfies the bound in Lemma \ref{grad} near $\ga_\infty$.

\begin{lemma}
\label{nearlyround}
1. There is a function $c:(0,\infty)\ra \R$ depending only
on $\kappa$, with
$\lim_{t\ra 0}c(t)=\infty$, such that if $w\in\ga_\infty$ then
$\bar R_\infty(w) \: d(y_\infty,w)^2 \: > \: c(\eps)$.

2. There is a function $\eps':(0,\infty)\ra \R\cup \{\infty\}$ depending only
on $\kappa$, with $\lim_{t\ra 0}\eps'(t)=0$, such that if 
$w\in \ga_\infty$ and $d(y_\infty,w)$ is sufficiently small
then the pointed manifold $(Z,w, \bar R_\infty(w)\bar g_\infty)$ is 
$2\eps^\prime(\epsilon)$-close to a round cylinder in the $C^2$ topology.
\end{lemma}
\begin{proof}
It follows from Lemma \ref{ptpicking} that 
for all $w\in \ga_\infty$, the pointed Riemannian 
manifold
$(Z,w,\bar R_\infty(w)\bar g_\infty)$ is $2\eps$-close
to (a time slice of) a pointed $\kappa$-solution.
From the definition of pointed closeness, there is an
embedded region around $w$, large on the scale defined by 
$\bar R_\infty(w)$, which is close to the corresponding subset of
a pointed $\kappa$-solution.  This gives a lower bound on 
the distance $\rho_0-d(w,x_\infty)$ to the point of curvature
blowup, thereby proving part 1 of the lemma. 

We know that $(Z,w, \bar R_\infty(w) \bar g_\infty)$ is $2\eps$-close to
a pointed $\kappa$-solution $(N,\star,h(t))$ in the pointed
$C^2$-topology.   By Lemma \ref{grad}, 
we know  $\bar R_\infty(w)$ tends to $\infty$
as $d(w,x_\infty)\ra\rho_0$, for $w\in\ga_\infty$.  
In particular, $\bar R_\infty(w)d^2(w,x_\infty)\ra \infty$.
From part 1 above, we can choose $\epsilon$ small enough in order to
make
$\bar R_\infty(w)d^2(w,y_\infty)$ large enough to apply Proposition
\ref{neckcrit}. Hence the pointed
manifold $(N,\star,R(\star)h(t))$ is $\eps^{\prime \prime}(\eps)$-close to
a round cylinder, where $\eps^{\prime \prime}:(0,\infty)\ra\R$ is a function with
$\lim_{t\ra 0}\eps^{\prime \prime}(t)=0$. The lemma follows.
\end{proof}

Note that the function $\epsilon^\prime$ in Lemma
\ref{nearlyround} is independent of the particular manifold
$Z$ that arises in our proof.

From part 2 of Lemma \ref{nearlyround}, if
$\eps$ is small and $ w \in \ga_\infty$ is
sufficiently close to $y_\infty$ then
$(Z, w, \overline{R}_\infty(w) \overline{g}_\infty)$ is 
$C^2$-close to a round cylinder.  The cross-section
of the cylinder has diameter approximately
$\pi \: (\overline{R}_\infty(w)/2)^{- \: \frac12}$. 
If we form the
union of the balls $B(w,2\pi(R_\infty(w)/2)^{- \: \frac{1}{2}})$, as $w$
ranges over such points in $\ga_\infty$, then we obtain a 
connected Riemannian manifold $W$. 
By adding in the point $y_\infty$, we get a
metric space $\bar W$ which is locally complete, and geodesic
near $y_\infty$. As the original manifolds $M_k$ had
$\Phi$-almost nonnegative curvature, it follows
from Lemma \ref{almost nonneg}
that $W$ is
nonnegatively curved.  Furthermore, 
$y_\infty$ cannot be an interior point of any geodesic segment in
$\bar W$, since such a geodesic would have to pass through a cylindrical 
region near $y_\infty$
twice.  The usual
proof of the Toponogov triangle comparison inequality now applies near
$y_\infty$ since
minimizers remain in the smooth nonnegatively curved part of $\bar W$.
Then $\overline{W}$ has nonnegative curvature in the Alexandrov sense.

This implies that blowups of $(\bar W,y_\infty)$ converge
to the tangent cone $C_{y_\infty}\bar W$.
As $\overline{W}$ is three-dimensional, so is
$C_{y_\infty}\bar W$
\cite[Corollary 7.11]{Burago-Gromov-Perelman}.
It will be
$C^2$-smooth away from the vertex and nowhere flat, by 
part 2 of Lemma \ref{nearlyround}. 
Pick $z\in C_{y_\infty}\bar W $ such that $d(z,y_\infty)=1$.
Then 
for any $\delta < \frac12$, 
the ball $B(z,\delta)\subset C_{y_\infty}\bar W $
is the Gromov-Hausdorff limit of a sequence of 
rescaled balls 
$B(\widehat{z}_k,\widehat{r}_k)\subset (M_k,\bar g_k(t_k))$
where $\widehat{r}_k\ra 0$, whose center points $(\widehat{z}_k,t_k)$
satisfy the conclusions of Theorem 
\ref{thmI.12.1}.  Applying Lemma \ref{rescaled}
and Appendix \ref{phiappendix},
if $\delta = \delta(\kappa)$ is taken small enough then
we get the curvature bounds needed to extract a limiting
Ricci flow solution whose time zero slice is isometric
to $B(z,\delta)$.
Now we can apply
the reasoning from the $0<\r<\infty$ case of Theorem
\ref{asympscalar} to get a contradiction.  This completes step 2. 

\bigskip
{\em Step 3: The  sequence of pointed
flows $(M_k,\bar g_k(\cdot),(x_k,t_k))$ accumulates on a 
pointed Ricci flow $(M_\infty,\bar g_\infty(\cdot),(x_\infty,t_0))$
which is defined on a time interval $[t',t_0]$ with $t'< t_0$.}
By step 2, we know that the scalar curvature   of $(M_k,\bar g_k(t_k))$
at $y\in M_k$ is bounded by a function of  the distance from 
$y$ to $x_k$.
Lemma \ref{rescaled} extends this curvature control to a
backward parabolic neighborhood centered at $y$ whose radius depends only
on $d(y,x_k)$.  Thus we can conclude,
using $\Phi$-pinching (\ref{double}) and 
Shi's estimates (Appendix \ref{applocalder}),
that all derivatives of the
curvature $(M_k,\bar g_k(t_k))$ are controlled as a function
of the distance from $x_k$, which means that the 
sequence of pointed manifolds $(M_k,\bar g_k(t_k),x_k)$ accumulates
to a smooth manifold $(M_\infty,\bar g_\infty)$. 

From Lemma \ref{almost nonneg},
$M_\infty$ has nonnegative sectional curvature.
We claim that $M_\infty$ has bounded curvature.
If not then there is a sequence of points
$q_k \in M_\infty$ so that
$\lim_{k \rightarrow \infty} \overline{R}(q_k) \: = \: \infty$
and $\overline{R}(q) \le 2 \overline{R}(q_k)$ for
$q \in B(y_k, A_k \overline{R}(q_k)^{- \frac12})$, where
$A_k \rightarrow \infty$; compare \cite[Lemma 22.2]{Hamilton}.
Lemma \ref{ptpicking} implies that for large $k$, a rescaled neighborhood of
$(M_\infty, q_k)$ is $\epsilon$-close to the corresponding subset
of a time slice of a $\kappa$-solution. As in the proof of
Theorem \ref{thmI.11.7}, we obtain a sequence of neck-like regions 
in $M_\infty$ with smaller
and smaller cross-sections, which contradicts the
existence of the Sharafutdinov retraction.

By Lemma \ref{rescaled} again,
we now get curvature control
 on $(M_k,\bar g_k(\cdot))$ for a time
interval $[t_k-\De,t_k]$ for some $\De>0$, and hence
we can extract a  subsequence which converges to a pointed Ricci flow
$(M_\infty,(x_\infty,t_0),\bar g_\infty(\cdot))$ defined
for $t\in[t_0 - \De,t_0]$, which has nonnegative curvature and bounded
curvature on 
compact time intervals.

\bigskip
{\em Step 4: Getting an ancient solution.}
Let $(t',t_0]$ be the maximal time interval on which 
we can extract a limiting solution $(M_\infty,\bar g_\infty(\cdot))$
with bounded curvature on compact time intervals.
Suppose that $t' > -\infty$.
By Lemma \ref{rescaled} the maximum of the scalar curvature
on the time slice $(M_\infty,\bar g_\infty(t))$ must tend
to infinity as $t\ra t'$.  From the trace Harnack inequality,
$R_t \: + \: \frac{R}{t - t^\prime} \: \ge \: 0$, and so
\begin{equation}
\bar R_\infty(x,t)\leq Q \frac{t_0 - t^\prime}{t-t'},
\end{equation}  
where $Q$ is the maximum of the scalar curvature
on $(M_\infty,\bar g_\infty(t_0))$.  Combining this
with Corollary \ref{bound1},
we get 
\begin{equation}
\frac{d}{dt}d_t(x,y)\geq \const \: \sqrt{Q \frac{t_0 - t^\prime}{t-t'}}.
\end{equation}
Since the right hand side is integrable on $(t',t_0]$, and using
the fact that distances are nonincreasing in time (since $\Rm \geq 0$),
it follows that there is a constant $C$ such that
\begin{equation}
\label{addbd}
|d_t(x,y)-d_{t_0}(x,y)|<C
\end{equation}
for all $x,y\in M_\infty$, $t\in(t',t_0]$.

If $M_\infty$ is compact then by (\ref{addbd}) the diameter of
$(M_\infty,\bar g_\infty(t))$ is bounded independent of
$t\in (t',t_0]$.  Since the minimum of the scalar
curvature is increasing in 
time, it is also bounded independent of $t$.  Now the argument
in Step 2 shows that the curvature is bounded everywhere
independent of $t$.  (We can apply the argument of Step 2 to the
time-$t$ slice because the main ingredient was
Lemma \ref{ptpicking}, which holds for rescaled time $t$.)

We may therefore assume $M_\infty$ is noncompact.
To be consistent with the notation of I.12.1, we now relabel
the basepoint $(x_\infty, t_0)$ as $(x_0, t_0)$.  Since
nonnegatively curved manifolds are asymptotically conical
(see Appendix \ref{Alexandrov}),
there is a constant $D$ such that if $y\in M_\infty$,
and $d_{t_0}(y,x_0)>D$, then there is a point $x\in M_\infty$
such that
\begin{equation}
\label{dists}
d_{t_0}(x,y)=d_{t_0}(y,x_0)\quad \mbox{and} 
\quad d_{t_0}(x,x_0)\geq \frac{3}{2}d_{t_0}(y,x_0);
\end{equation} 
by (\ref{addbd})
the same conditions hold at all times $t\in(t',t_0]$, up to error
$C$.   If for some such $y$, and some $t\in (t',t_0]$ the scalar
curvature were large, then $(M_\infty,(y,t),\bar R_\infty(y,t)\bar 
g_\infty(t))$
would be $2\eps$-close to a $\kappa$-solution $(N,h(\cdot),(z,t_0))$.  When
$\eps$ is small we could  use Proposition \ref{neckcrit} to 
see that $y$ lies in a neck region $U$ in $(M_\infty,\bar g_\infty(t))$ 
of diameter $\approx R(y,t)^{-\frac{1}{2}}\ll 1$.

We claim that $U$ separates
$x_0$ from $x$ in the sense that $x_0$ and $x$ belong
to disjoint components of $M_\infty - U$, where $x_0, y$, and $x$
satisfy (\ref{dists}).
To see this, we
recall that if $M_\infty$ has more than one end then it splits
isometrically, in which case the claim is clear. If $M_\infty$ has
one end then we consider an exhaustion of $M_\infty$ by totally
convex compact sets $C_t$ as in Section \ref{pf11.7}. One of the
sets $C_t$ will have boundary consisting of an approximate $2$-sphere
cross-section in the neck region $U$, giving the separation.

(For another argument, suppose that $U$ does not separate $x_0$ from 
$x$. Since $\cangle_y(x_0, x) > \frac{\pi}{2}$ (from Proposition \ref{neckcrit}),
the segments $\overline{yx_0}$ and $\overline{yx}$ must exit $U$ by
different ends.  If $x_0$ and $x$ can be joined by a curve avoiding $U$
then there is a nonzero element of
$\HH^1(M_\infty, \Z)$. The corresponding infinite cyclic cover of $M_\infty$ will then
isometrically split off a line and its quotient $M_\infty$ will be
compact, which is a contradiction.  We thank one of the referees for
this argument.)

Obviously, at time $t_0$ the set
$U$ still separates $x_0$ from $x$. 
Since $\bar g_\infty$ has nonnegative curvature,
we have $\diam_{t_0}(U)\leq \diam_t(U)\ll 1$. 
Since $(M_\infty,\bar g_\infty({t_0}))$ has bounded geometry,
there cannot be topologically separating subsets of arbitrarily small
diameter. 
Thus there must be a uniform upper bound on $R(y, t)$
and the curvature of $(M_\infty,\bar g_\infty)$ is uniformly bounded
(in space and time)
outside a set of uniformly bounded diameter.  Repeating the reasoning
from Step 2, we get uniform bounds everywhere.  This contradicts
our assumption that the curvature blows up as $t\ra t'$.

It remains to show that the ancient solution is a $\kappa$-solution.
The only remaining point is to show that it is $\kappa$-noncollapsed
at all scales. This follows from the fact that the original Ricci
flow solutions $(M_k, g_k(\cdot))$ were $\kappa$-noncollapsed on 
scales less than the fixed number $\sigma$.  
\end{proof}

\begin{remark} \label{wang}
As mentioned in Remark \ref{wangg}, the statement of 
\cite[Theorem 12.1]{Perelman} instead assumes noncollapsing at all
scales less than $r_0$. Bing Wang pointed out that with this
assumption, after constructing the ancient solution in Step 4
of the proof, one only gets that it is $\kappa$-collapsed at
all scales less than one. Hence it may not be a $\kappa$-solution.
The literal statement of \cite[Theorem 12.1]{Perelman} is not
used in the rest of \cite{Perelman,Perelman2}, but rather its
method of proof.  Because of this, the change of hypotheses 
does not seem to lead to any problems. The method of proof
of Theorem \ref{thmI.12.1} is used in two different ways.
The first way is to construct the Ricci flow with surgery on a
fixed finite time interval, as in Section \ref{II.5}.
In this case the
noncollapsing at a given scale $\sigma$ comes from 
Theorem \ref{nolocalcollapse}, and its extension when surgeries
are allowed.  The second way is to analyze the large-time
behavior of the Ricci flow, as in the next few sections.
\end{remark}
 
\section{I.12.2. Later scalar curvature bounds on bigger balls
from curvature and volume bounds}
\label{I.12.2}

The next theorem roughly says that if one has a sectional curvature bound on a
ball, for a certain time interval, and
a lower bound on the volume of the ball at the initial time, then one obtains
an upper scalar curvature bound on a larger ball at the final time.

We first write out the corrected version of the theorem
(see II.6.2).
\begin{theorem} \label{thmI.12.2}
For any $A < \infty$, there exist $K = K(A) < \infty$ and
$\rho = \rho(A) >  0$ with the following property. Suppose in
dimension three we have a Ricci flow solution with
$\Phi$-almost nonnegative curvature.
Given $x_0 \in M$ and $r_0 > 0$, suppose that 
$r_0^2 \: \Phi(r_0^{-2}) \: < \: \rho$,
the solution is defined for
$0 \le t \le r_0^2$ and it has 
$|\Rm|(x, t) \le \frac{1}{3r_0^{2}}$
for all $(x, t)$ satisfying $\dist_0(x, x_0) < r_0$. Suppose in
addition that the volume of the metric ball $B(x_0, r_0)$ at
time zero is at least $A^{-1} r_0^3$.
Then $R(x, r_0^2) \le K r_0^{-2}$ whenever $\dist_{r_0^2}(x, x_0) <
A r_0$.
\end{theorem}

\begin{remark} \label{rmkI.12.2}
The added restriction that $r_0^2 \: \Phi(r_0^{-2}) \: < \: \rho$ (see II.6.2)
imposes an upper bound on $r_0$.  This is necessary, as
otherwise the conclusion would imply that neck pinches cannot occur.

There is an apparent gap in the proof of \cite[Theorem 12.2]{Perelman}, 
in the sentence
(There is a little subtlety...).  We instead follow the proof
of II.6.3(b,c) (see Proposition \ref{propII.6.3}(b,c)), 
which proves the same statement in the
presence of surgeries.

The volume assumption in the theorem is used to guarantee
noncollapsing, by means of Theorem \ref{8.2thm}. The reason
for the ``$3$'' in the hypothesis $|\Rm|(x, t) \le \frac{1}{3r_0^{2}}$
comes from Remark \ref{8.2fix}.
\end{remark}
\begin{proof}
The proof is in two steps.  In the first step one shows that
if $R(x, r_0^2)$ is large then a parabolic neighborhood of
$(x, r_0^2)$ is close to the corresponding subset of a
$\kappa$-solution.  In the second part one uses this to
prove the theorem.

The first step is the following lemma.

\begin{lemma} \label{12.2lemma}
For any $\epsilon > 0$ there exists 
$K = K(A, \epsilon) < \infty$
so that for any $r_0$,
whenever we have a solution as in the statement of the theorem 
and $\dist_{r_0^2}(x, x_0) <
A r_0$ then  \\
(a) $R(x, r_0^2) < K r_0^{-2}$ or \\
(b) The solution in
$\{(x^\prime, t^\prime) \: : \: \dist_{t}(x^\prime, x) <
(\epsilon Q)^{-1}, t - (\epsilon Q)^{-1} \le t^\prime \le t \}$ is,
after scaling by the factor $Q$, $\epsilon$-close to the
corresponding subset of a $\kappa$-solution. 

Here $t = r_0^2$ and
$Q = R(x,t)$.
\end{lemma}
\begin{remark}
One can think of this lemma as a
localized analog of Theorem \ref{thmI.12.1}, where ``localized'' refers to
the fact that both the hypotheses and the conclusion involve the
point $x_0$.
\end{remark}

\begin{proof}
To prove the lemma, 
suppose that there is
a sequence of such pointed solutions
$(M_k, x_{0,k}, g_k(\cdot))$, along with points
$\hat{x}_k \in M_k$, so that $\dist_{r_0^2}(\hat{x}_k, x_{0,k}) <
A r_0$ and $r_0^2 \: R(\hat{x}_k, r_0^2) \rightarrow \infty$, but 
$(\hat{x}_k, r_0^2)$ does not satisfy conclusion (b) of the
lemma. As in the proof of Theorem \ref{thmI.12.1}, we will allow ourselves to make
$\epsilon$ smaller during the course of the proof.

We first show that there is a sequence $D_k \rightarrow \infty$
and modified points $(\overline{x}_k, \overline{t}_k)$
with $\frac{3}{4} r_0^{2} \le \overline{t}_k \le r_0^{2}$,
$\dist_{\overline{t}_k}(\overline{x}_k, x_{0,k}) < (A+1) r_0$ and
$Q_k = R(\overline{x}_k, \overline{t}_k) \rightarrow \infty$, 
so that any point
$(x^\prime_k, t^\prime_k)$ with $R(x^\prime_k, t^\prime_k) > 2Q_k$, 
$\overline{t}_k - D_k^2 Q_k^{-1} \le
t^\prime_k \le \overline{t}_k$ and 
$\dist_{t^\prime_k}(x^\prime_k, {x}_{0,k}) 
< \dist_{\overline{t}_k}(\overline{x}_k, 
x_{0,k}) +
D_k Q_k^{-1/2}$ satisfies conclusion (b) of the lemma, 
but $(\overline{x}_k, \overline{t}_k)$ does not satisfy conclusion (b) of
the lemma. (Of course, in saying ``$(x^\prime_k, t^\prime_k)$
satisfies conclusion (b)'' or
``$(\overline{x}_k, \overline{t}_k)$ does not satisfy conclusion (b)'', we mean that
the $(x,t)$ in conclusion (b) is replaced by $(x^\prime_k, t^\prime_k)$ or
$(\overline{x}_k, \overline{t}_k)$, respectively.)

The construction of $(\overline{x}_k, \overline{t}_k)$ is by a pointpicking
argument.
Put $D_k = \frac{r_0 R(\hat{x}_k, r_0^2)^{1/2}}{10}$.
Start with $(x_k, t_k) = (\hat{x}_k, r_0^2)$
and look if there is a point
$(x^\prime_k, t^\prime_k)$ with $R(x^\prime_k, t^\prime_k) > 2
R(x_k, t_k)$, 
$\overline{t}_k - D_k^2 R(x_k, t_k)^{-1} \le
t^\prime_k \le \overline{t}_k$ and 
$\dist_{t^\prime_k}(x^\prime_k, {x}_{0,k}) 
< \dist_{{t}_k}({x}_k, 
x_{0,k}) +
D_k R(x_k, t_k)^{-1/2}$, but 
which does not have a neighborhood that is
$\epsilon$-close to the corresponding subset of a $\kappa$-solution.
If there is such a point, we replace $(x_k, t_k)$ by $(x^\prime_k, 
t^\prime_k)$ and repeat the process.  The process must terminate after
a finite number of steps to give a point
$(\overline{x}_k, \overline{t}_k)$ with the desired property.

(Note that the condition
$\dist_{t^\prime_k}(x^\prime_k, {x}_{0,k}) 
< \dist_{\overline{t}_k}(\overline{x}_k, 
x_{0,k}) +
D_k Q_k^{-1/2}$ involves the metric at time $t^\prime_k$. In order
to construct an ancient solution, one of the issues will be to
replace this by a condition that only involves the metric at time
$\overline{t}_k$, i.e. that involves a parabolic neighborhood
around $(\overline{x}_k, \overline{t}_k)$.)

Let $\overline{g}_k(\cdot)$ denote
the rescaling of the solution $g_k(\cdot)$ by $Q_k$. We
normalize the time interval of the rescaled solution
by fixing a number ${t}_\infty$
and saying that for all $k$,
the time-$\overline{t}_k$ slice of $(M_k, g_k)$ corresponds to the
time-${t}_\infty$ slice of $(M_k, \overline{g}_k)$.
Then the scalar
curvature $\overline{R}_k$ of $\overline{g}_k$ satisfies 
$\overline{R}_k(\overline{x}_k, {t}_\infty) \: = \: 1$.

By the argument of Step 2 of the proof of Theorem \ref{thmI.12.1}, 
a subsequence of the
pointed spaces $(M_k, \overline{x}_k, \overline{g}_k({t}_\infty))$ 
will smoothly converge to a nonnegatively-curved pointed space $(M_\infty,
\overline{x}_\infty,\overline{g}_\infty)$. 
By the pointpicking, if $m \in M_\infty$ has $\overline{R}(m) \ge 3$ 
then a parabolic neighborhood of $m$ is $\epsilon$-close 
to the corresponding region in a $\kappa$-solution.
It follows,
as in Step 3 of the proof of Theorem
\ref{thmI.12.1}, that the sectional curvature of $M_\infty$ will be bounded above by
some $C < \infty$. Using Lemma \ref{rescaled}, the metric on 
$M_\infty$ is the time-$t_\infty$ slice of
a nonnegatively-curved Ricci flow solution defined on some time interval 
$[t_\infty - c , t_\infty]$, with $c > 0$, 
and one has convergence of a subsequence $\overline{g}_k(t)
\rightarrow \overline{g}_\infty(t)$ for $t \in [t_\infty - c , t_\infty]$.
As $\overline{R}_t \ge 0$, the scalar curvature on this time interval
will be uniformly bounded above by $6C$ and so from  the
$\Phi$-almost nonnegative curvature (see (\ref{double})), the
sectional curvature will be uniformly bounded above on the time interval.
Hence we can apply Lemma \ref{distancedist} to get a uniform additive bound on the 
length
distortion between times $t_\infty - c$ and $t_\infty$ (see Step 4 of the
proof of Theorem \ref{thmI.12.1}). More precisely, in applying Lemma 
\ref{distancedist}, we
use the curvature bound coming from the hypotheses of the theorem near 
$x_0$, and the just-derived upper curvature bound near $\overline{x}_k$.

It follows that for a given $A^\prime > 0$, for large
$k$, if $t_k^\prime \in [\overline{t}_k - c Q_k^{-1}/2, \overline{t}_k]$ and 
$\dist_{\overline{t}_k}(x_k^\prime, \overline{x}_k) < A^\prime Q_k^{-1/2}$ 
then $\dist_{t^\prime_k}(x^\prime_k, {x}_{0,k}) 
< \dist_{\overline{t}_k}(\overline{x}_k, 
x_{0,k}) +
D_k Q_k^{-1/2}$. In particular, if a point
$(x_k^\prime, t_k^\prime)$ lies in the
parabolic neighborhood given by
$t_k^\prime \in [\overline{t}_k - c Q_k^{-1}/2, \overline{t}_k]$ and 
$\dist_{\overline{t}_k}(x_k^\prime, \overline{x}_k) < A^\prime Q_k^{-1/2}$, 
and  has 
$R(x_k^\prime, t_k^\prime) > 2Q_k$, then it has a neighborhood that
is $\epsilon$-close to the corresponding subset of a $\kappa$-solution.

As in Step 4 of the proof of Theorem \ref{thmI.12.1}, we now extend 
$(M_\infty,
\overline{g}_\infty, \overline{x}_\infty)$ backward to an ancient solution 
$\overline{g}_\infty(\cdot)$, defined for $t \in (-\infty, t_\infty]$.
To do so, we use the fact that if the solution
is defined backward to a time-$t$ slice then 
the length distortion bound, along with the pointpicking,
implies that a point $m$ in a time-$t$ slice with $\overline{R}_\infty(m)
> 3$ has a neighborhood that is $\epsilon$-close to the corresponding 
subset of a $\kappa$-solution. The ancient solution is $\kappa$-noncollapsed
at all scales since the original solution was $\kappa$-noncollapsed at
some scale, by Theorem \ref{8.2thm}.
Then we obtain smooth convergence
of parabolic regions of the points $(\overline{x}_k, \overline{t}_k)$ to
the $\kappa$-solution, which is a contradiction to the choice of the
$(\overline{x}_k, \overline{t}_k)$'s.
\end{proof}

We now know that regions of high scalar curvature are modeled by
corresponding
regions in $\kappa$-solutions. To continue with the proof of the theorem,
fix $A$ large. 
Suppose that the theorem is not true.  Then there are \\
1. Numbers $\rho_k \rightarrow 0$, \\
2. Numbers $r_{0,k}$ with 
$r_{0,k}^{2} \: \Phi(r_{0,k}^{-2}) \: \le \: \rho_k$, 
\\
3. Solutions
$(M_k, g_k(\cdot))$ defined for $0 \le t \le r_{0,k}^2$, \\
4. Points $x_{0,k} \in M_k$ and \\
5. Points $x_k \in M_k$\\
so that \\
a. $|\Rm|(x, t) \le r_{0,k}^{-2}$ for
all $(x, t) \in M_k \times [0, r_{0,k}^2]$ 
satisfying $\dist_0(x, x_{0,k}) < r_{0,k}$, \\
b. The volume of the metric ball $B(x_{0,k}, r_{0,k})$ at
time zero is at least $A^{-1} r_{0,k}^3$ and \\
c.  $\dist_{r_{0,k}^2}(x_k, x_{0,k}) <
A r_{0,k}$, but \\
d. $r_{0,k}^2 \: R(x_k, r_{0,k}^2) \rightarrow \infty$.

We now apply Step 2 of the proof of Theorem \ref{thmI.12.1} to obtain a contradiction.
That is, we take a subsequence of
$\{(M_k, x_{0,k}, r_{0,k}^{-2} \: g_k(r_{0,k}^2))\}_{k=1}^\infty$ that converges on a maximal ball.
The only difference is that in Theorem \ref{thmI.12.1}, the nonnegative curvature
of $W$ came from the $\Phi$-almost nonnegative curvature assumption on the
original manifolds $M_k$ along with the fact (with the notation of the proof of
Theorem \ref{thmI.12.1}) 
that the numbers $Q_k = R(x_k, t_k)$, which
we used to rescale, go to infinity.
In the present case 
the rescaled scalar curvatures $r_{0,k}^2 \: R(x_{0,k}, r_{0,k}^2)$ at the 
basepoints $x_{0,k}$
stay bounded.  
However, if a point $y \in W$ is a limit of points
$\widetilde{x}_k \in M_k$ then the equations
\begin{equation}
\Rm(\widetilde{x}_k, r_{0,k}^2) \: \ge \: - \: \Phi(R(\widetilde{x}_k, r_{0,k}^2))\end{equation} in the form
\begin{equation}
r_{0,k}^2 \: \Rm(\widetilde{x}_k, r_{0,k}^2) \: \ge \: - \: 
\frac{\Phi \left(  r_{0,k}^{2}
R(\widetilde{x}_k, r_{0,k}^2) \cdot r_{0,k}^{-2}
\right)}{ R(\widetilde{x}_k, r_{0,k}^2)} \: 
r_{0,k}^2 \: R(\widetilde{x}_k, r_{0,k}^2)
\end{equation}
pass to the limit to give $\Rm(y) \: \ge 0$
(using that $y \in W$, so $r_{0,k}^2 \: \Rm(\widetilde{x}_k, r_{0,k}^2) \rightarrow
R(y) > 0$).
This is enough to carry out the argument.
\end{proof}

\section{I.12.3. Earlier scalar curvature bounds on smaller balls from 
lower curvature bounds and volume bounds}
\label{I.12.3}

The main result of this section
says that if one has a lower bound on volume and sectional curvature
on a ball at a certain time then one obtains an upper scalar curvature bound on
a smaller ball at an earlier time.

We first prove a result in Riemannian geometry saying that under
certain hypotheses, metric balls have subballs of a controlled size
with almost-Euclidean volume.

\begin{lemma} \label{Alexlemma}
Given $w^\prime > 0$ and $n \in \Z^+$, there is a number $c = c(w^\prime,n) > 0$ with the
following property.
Let $B$ be a radius-$r$ ball with compact closure in an $n$-dimensional Riemannian 
manifold. Suppose that the sectional curvatures of $B$ are bounded below by
$- \: r^{-2}$.
Suppose that $\vol(B) \ge w^\prime r^n$. Then there is a subball $B^\prime \subset B$ of
radius $r^\prime \ge cr$ so that $\vol(B^\prime) \: \ge \: \frac12 \: \omega_n \:
(r^\prime)^n$, where $\omega_n$ is the volume of the unit ball in $\R^n$.
\end{lemma}
\begin{proof}
Suppose that the lemma is not true. Rescale
so that $r \: = \: 1$. Then there is a sequence of Riemannian
manifolds $\{M_i\}_{i=1}^\infty$ with balls $B(x_i, 1) \subset M_i$
having compact closure
so that $\Rm \Big|_{B(x_i, 1)} \: \ge \: - \: 1$ and
$\vol(B(x_i, 1)) \: \ge \: w^\prime$, but with the property that all balls
$B(x_i^\prime, r^\prime) \subset B(x_i, 1)$ with $r^\prime \: \ge \:
i^{-1}$ satisfy
$\vol(B(x_i^\prime, r^\prime)) \: < \: \frac{1}{2} \: \omega_n 
(r^\prime)^n$. After taking a subsequence, we can assume that
$\lim_{i \rightarrow \infty} (B(x_i, 1), x_i) \: = \: (X, x_\infty)$ in the
pointed Gromov-Hausdorff topology. From
\cite[Theorem 10.8]{Burago-Gromov-Perelman}, the Riemannian volume
forms $\dvol_{M_i}$ converge weakly to the three-dimensional Hausdorff
measure $\mu$ of $X$. From \cite[Corollary 6.7 and 
Section 9]{Burago-Gromov-Perelman}, for any $\epsilon > 0$, 
there are small balls $B(x_\infty^\prime,
r^\prime) \subset X$ with compact closure in $X$ such that $\mu(B(x_\infty^\prime,
r^\prime)) \: \ge \: (1-\epsilon) \: \omega_n (r^\prime)^n$.
This gives a contradiction.
\end{proof}

\begin{theorem} (cf. Theorem I.12.3) \label{thmI.12.3}
For any $w > 0$ there exist $\tau = \tau(w) > 0$, $K = K(w) < \infty$ and
$\rho = \rho(w) > 0$ with the following property.  Suppose that
$g(\cdot)$ is a Ricci flow on a closed three-manifold $M$, defined for
$t \in [0, T)$, with $\Phi$-almost nonnegative curvature. Let $(x_0, t_0)$
be a spacetime point and let $r_0 > 0$ be a radius with $t_0 \ge 4 \tau r_0^2$
and 
$r_0^2 \: \Phi(r_0^{-2}) \: < \: \rho$.
Suppose that $\vol_{t_0}(B(x_0,r_0)) \ge w r_0^3$ and
the time-$t_0$ sectional curvatures on $B(x_0, r_0)$ are bounded below by
$- \: r_0^{-2}$. Then $R(x,t) \le K r_0^{-2}$ whenever
$t \in [t_0 - \tau r_0^2, t_0]$ and
$\dist_t(x, x_0) \: \le \: \frac14 \: r_0$.
\end{theorem}
\begin{proof}
Let $\tau_0(w)$, $B(w)$ and $C(w)$ be the constants of Corollary
\ref{I.11.6end}. Put $\tau(w) \: = \: \frac12 \: \tau_0(w)$ and
$K(w) = C(w) + 2 \frac{B(w)}{\tau_0(w)}$. The function $\rho(w)$ 
will be specified in the course of the proof.

Suppose that the theorem is not true. Take a counterexample with
a point $(x_0, t_0)$ and a radius $r_0 > 0$ such that the
time-$t_0$ ball $B(x_0, r_0)$ satisfies
the assumptions of the theorem, but the conclusion of the theorem
fails.  We claim there is a counterexample coming from a point 
$(\widehat{x}_0, \widehat{t}_0)$
and a radius $\widehat{r}_0 > 0$, with the additional property that for any
$(x^\prime_0, t^\prime_0)$ and $r^\prime_0$ having
$t^\prime_0 \in [\widehat{t}_0 - 2 \tau \widehat{r}_0^2, \widehat{t}_0]$ and $r^\prime_0 \le
\frac{1}{2} \widehat{r}_0$, if $\vol_{t_0^\prime}(B(x_0^\prime,r_0^\prime)) \ge w 
(r_0^\prime)^3$ and
the time-$t_0^\prime$ sectional curvatures on $B(x_0^\prime, r_0^\prime)$ 
are bounded below by
$- \: (r_0^\prime)^{-2}$ then $R(x,t) \le K (r_0^\prime)^{-2}$ whenever
$t \in [t_0^\prime - \tau (r_0^\prime)^2, t_0^\prime]$ and
$\dist_t(x, x_0^\prime) \: \le \: \frac14 \: r_0^\prime$.
This follows from a
pointpicking argument - suppose that it is not
true for the original $x_0, t_0, r_0$.  Then there are
$(x^\prime_0, t^\prime_0)$ and $r^\prime_0$ with
$t^\prime_0 \in [t_0 - 2 \tau r_0^2, t_0]$ and $r^\prime_0 \le
\frac{1}{2} r_0$, for which 
the assumptions of the theorem hold but the conclusion
does not. If the triple 
$(x^\prime_0, t^\prime_0, r^\prime_0)$ satisfies the claim then we
stop, and otherwise we iterate the procedure.  The iteration must
terminate, which provides the desired triple 
$(\widehat{x}_0, \widehat{t}_0, \widehat{r}_0)$. Note that 
$\widehat{t}_0 \: > \: t_0 - 4 \tau r_0^2 \: \ge \: 0$.

We relabel $(\widehat{x}_0, \widehat{t}_0, \widehat{r}_0)$ as
$({x}_0, {t}_0, {r}_0)$.
For simplicity, let us assume that
the time-$t_0$ sectional curvatures on $B(x_0, r_0)$ are strictly
greater than
$- \: r_0^{-2}$; the general case will follow from continuity.
Let $\tau^\prime > 0$ be the largest number such that
$\Rm(x,t) \: \ge \: -\: r_0^{-2}$ whenever
$t \in [t_0 - \tau^\prime r_0^2, t_0]$ and $\dist_t(x, x_0) \le r_0$.
If $\tau^\prime \ge 2 \tau = \tau_0(w)$ then Corollary
\ref{I.11.6end} implies that
$R(x,t) \le C r_0^{-2} + B (t-t_0+ 2 \tau r_0^2)^{-1}$ whenever
$t \in [t_0 - 2 \tau r_0^2, t_0]$ and 
$\dist_t(x, x_0) \: \le \: \frac14 \: r_0$.
In particular, $R(x,t) \le K$ whenever 
$t \in [t_0 - \tau r_0^2, t_0]$ and 
$\dist_t(x, x_0) \: \le \: \frac14 \: r_0$, which contradicts our
assumption that the conclusion of the theorem fails.

Now suppose that $\tau^\prime < 2 \tau$. Put
$t^\prime = t_0 - \tau^\prime r_0^2$.
From estimates on
the length and volume distortion under the Ricci flow, we know that
there are numbers $\alpha = \alpha(w) > 0$ and
$w^\prime = w^\prime(w) > 0$ so that the time-$t^\prime$ ball
$B(x_0, \alpha r_0)$ has volume at least $w^\prime (\alpha r_0)^3$.
From Lemma \ref{Alexlemma}, there is a subball
$B(x^\prime, r^\prime) \subset B(x_0, \alpha r_0)$ with
$r^\prime \ge c \alpha r_0$ and $\vol(B(x^\prime, r^\prime)) \: \ge \:
\frac12 \: \omega_3 \: (r^\prime)^3$. From the preceding
pointpicking argument, we have the estimate $R(x,t) \: \le \: K
(r^\prime)^{-2}$ whenever $t \in [t^\prime - \tau (r^\prime)^2,
t^\prime]$ and $\dist_t(x, x^\prime) \: \le \: \frac14 r^\prime$.
From the $\Phi$-almost nonnegative curvature, 
we have a bound 
$|\Rm|(x,t) \: \le \:
\const \: K \: (r^\prime)^{-2} \: + \: \const \: \Phi( K (r^\prime)^{-2})$
at
such a point (see (\ref{double})).  If $\rho(w)$ is taken sufficiently
small then we can ensure that $r_0$ is small enough, and
hence $r^\prime$ is small enough, to make
$K \: (r^\prime)^{-2} \: + \: \Phi( K (r^\prime)^{-2}) \: \le \:
2 \: K \: (r^\prime)^{-2}$.
Then we can apply
Theorem \ref{thmI.12.2} to a time interval ending at time
$t^\prime$, after a redefinition of its
constants, to obtain a bound of the form
$R(x, t^\prime) \: \le \: K^\prime (r^\prime)^{-2}$ 
whenever $\dist_{t}(x, x^\prime) \: \le \: 10 r_0$, where
$K^\prime$ is related to the constant 
$K$ of Theorem \ref{thmI.12.2}. (We also obtain a similar estimate at
times slightly less than $t^\prime$.)  Thus at such a point,
$\Rm(x, t^\prime) \: \ge \: - \: 
\Phi(K^\prime (r^\prime)^{-2})$.  If we choose $\rho(w)$ to be
sufficiently small to force $r^\prime$ to be sufficiently small to force
$- \: 
\Phi(K^\prime (r^\prime)^{-2})\: > \: - \: r_0^{-2}$
then we have
$\Rm > - r_0^{-2}$ on $\overline{B_{t^\prime}(x_0, r_0)} \subset
B_{t^\prime}(x^\prime, 10r_0)$, which 
contradicts the assumed maximality of $\tau^\prime$.

We note that in the application of Theorem \ref{thmI.12.2}
at the end of the proof, we must take into account the
extra hypothesis, in the notation of
Theorem \ref{thmI.12.2}, that
$r_0^2 \: \Phi(r_0^{-2}) \: < \: \rho$
(see Remark \ref{rmkI.12.2}).  This
will be satisfied if the $r_0$ in Theorem \ref{thmI.12.2} is small
enough, which is ensured
by taking the $\rho$ of Theorem \ref{thmI.12.3} small enough.
\end{proof}

\section{I.12.4. Small balls with strongly negative curvature are
volume-collapsed}
\label{I.12.4}

In this section we show that under certain hypotheses, if the infimal sectional 
curvature on an $r$-ball is exactly $- \: r^{-2}$ then the volume of the ball
is small compared to $r^3$.

\begin{corollary} (cf. Corollary I.12.4) \label{corI.12.4}
For any $w > 0$, one can find $\rho > 0$ with the following property.
Suppose that $g(\cdot)$ is a $\Phi$-almost nonnegatively curved
Ricci flow solution on a closed three-manifold $M$, defined for
$t \in [0, T)$ with $T \ge 1$. If $B(x_0, r_0)$ is a metric ball
at time $t_0 \ge 1$ with $r_0 < \rho$ and if
$\inf_{x \in B(x_0,r_0)}  \Rm(x, t_0)  = - r_0^{-2}$ then
$\vol(B(x_0, r_0)) \le w r_0^3$.
\end{corollary}
\begin{proof}
Fix $w > 0$. 
The number $\rho$ will be specified in the course of the proof.
Suppose that the corollary is not true, i.e. 
there is a Ricci flow solution as in the statement of the corollary along
with a metric ball $B(x_0, r_0)$ at a time $t_0 \ge 1$
so that  $\inf_{x \in B(x_0,r_0)}  \Rm(x, t_0)  = - r_0^{-2}$ 
and $\vol(B(x_0, r_0)) \: > \: w r_0^n$.
The idea is to use
Theorem \ref{thmI.12.3}, along with the $\Phi$-almost nonnegative
curvature, to get a double-sided sectional
curvature bound on a smaller
ball at an earlier time. Then one goes forward in time using 
Theorem \ref{thmI.12.2}, along with the $\Phi$-almost nonnegative
curvature, to get a lower sectional curvature bound on the
original ball, thereby obtaining a contradiction.

Looking at the hypotheses of Theorem \ref{thmI.12.3}, if
we require $r_0 < (4\tau)^{- \: \frac12}$ then
$4 \tau r_0^2 < 1 \le t_0$.
From Theorem \ref{thmI.12.3}, $R(x, t) \: \le \: K r_0^{-2}$
whenever  $t \in [t_0 - \tau r_0^2, t_0]$
and $\dist_t(x, x_0) \: \le \: \frac{1}{4} r_0$, provided that $r_0$ is
small enough that
$r_0^2 \: \Phi(r_0^{-2})$
is less than the $\rho$ of Theorem \ref{thmI.12.3}.
If in addition $r_0$ is sufficiently small then it follows that
$|\Rm(x,t)| \: \le \: \const \: \Phi(K r_0^{-2}) \: \le \: r_0^{-2}$. 

From the Bishop-Gromov inequality and the bounds on length and
volume distortion under Ricci flow, there is a small number $c$ so
that we are ensured that $|\Rm(x,t)| \le (cr_0)^{-2}$ for all $(x,t)$
satisfying $\dist_{t_0 - (cr_0)^2}(x, x_0) < cr_0$ and
$t \in [t_0 - (cr_0)^2, t_0]$, and in addition the volume of
$B(x_0, cr_0)$ at time $t_0 - (cr_0)^2$ is at least
$c(cr_0)^3$. Choosing the constant $A$ of Theorem \ref{thmI.12.2}
appropriately in terms of $c$, we can apply Theorem \ref{thmI.12.2}
to the ball $B(x_0, cr_0)$ and the time interval $[t_0 - (cr_0)^2, t_0]$
to conclude that at time $t_0$,
$R(\cdot, t_0) \Big|_{B(x_0, r_0)} \: \le \: K(A) \:
(cr_0)^{-2}$, where $K(A)$ is as in the statement of
Theorem \ref{thmI.12.2}.
From the $\Phi$-almost nonnegative curvature condition,
\begin{equation}
\Rm \Big|_{B(x_0, r_0)} \: \ge \: - \:
\Phi \left(K(A) \:
(cr_0)^{-2} \right).
\end{equation}
If $r_0$ is sufficiently small then we contradict the assumption that
$\inf \Rm(x, t_0) 
\Big|_{B(x_0,r_0)} \: = \: - \:
\frac{1}{r_0^2}$.
\end{proof}

\section{I.13.1. Thick-thin decomposition for nonsingular flows} \label{I.13.1}

The main result of this section says that if $g(\cdot)$ is a Ricci flow
solution on a closed oriented three-dimensional manifold $M$ that
exists for $t \in [0, \infty)$ then for large $t$, $(M, g(t))$ has a
thick-thin decomposition. 
A fuller description is in Sections \ref{II.7.1}-\ref{II.7.4}.

We assume that at time zero, the sectional curvatures are bounded below
by $-1$. This can always be achieved by rescaling the initial
metric.  Then we have the $\Phi$-almost nonnegative curvature result of
(\ref{phi}). 

If the metric $g(t)$ has nonnegative sectional curvature then it must be flat, as we are
assuming that the Ricci flow exists for all time. Let us assume that $g(t)$ is not flat, so 
it has some negative sectional curvature.
Given $x \in M$, consider the time-$t$ ball $B_t(x,r)$. Clearly if $r$ is 
sufficiently small then
$\Rm \Big|_{B_t(x,r)} \: > \: - \: r^{-2}$, while if $r$ is sufficiently large
(maybe greater than the diameter of $M$) then it is not true that 
$\Rm \Big|_{B_t(x,r)} \: > \: - \: r^{-2}$.
Let $\widehat{r}(x,t) > 0$ be the unique number such that
$\inf \Rm \Big|_{B_t(x,\widehat{r})} \: = \: - \: \widehat{r}^{-2}$.
Let $M_{thin}(w,t)$ be the set of points $x \in M$ for which
\begin{equation} \label{thinpart}
\vol(B_t(x, \widehat{r}(x,t))) \: < \: w \: \widehat{r}(x,t)^3.
\end{equation}
Put $M_{thick}(w,t) = M - M_{thin}(w,t)$. 

As the statement of (\ref{phi}) is invariant under parabolic rescaling
(although we must take $t \ge t_0$ for (\ref{phi}) to apply), if $t \ge t_0$ 
and we are interested in the Ricci flow at time $t$ then we can apply
Theorem \ref{thmI.12.2}, Theorem \ref{thmI.12.3} and Corollary \ref{corI.12.4}
to the rescaled flow $g(t^\prime) \: = \: t^{-1} g(tt^\prime)$.
From Corollary \ref{corI.12.4}, 
for any $w > 0$ we can find $\widehat{\rho} = \widehat{\rho}(w) > 0$ 
so that if $\widehat{r}(x,t) < \widehat{\rho} \sqrt{t}$ then $x \in M_{thin}(w,t)$,
provided that $t$ is
sufficiently large (depending on $w$). Equivalently, if 
$t$ is sufficiently large (depending on $w$) and $x \in M_{thick}(w,t)$ then
$\widehat{r}(x,t) \ge \widehat{\rho} \sqrt{t}$.

\begin{theorem} (cf. I.13.1) \label{thmI.13.1}
There are numbers $T = T(w) > 0$, $\overline{\rho} = \overline{\rho}(w) > 0$
and $K = K(w) < \infty$ so that if $t \ge T$ and
$x \in M_{thick}(w,t)$ then $|\Rm| \le K t^{-1}$ on $B_t(x, \overline{\rho} \sqrt{t})$,
and $\vol(B_t(x, \overline{\rho} \sqrt{t})) \: \ge \: \frac{1}{10} \: 
w \left( \overline{\rho} \sqrt{t} \right)^3$.
\end{theorem}
\begin{proof}
The method of proof is the same as in Corollary \ref{corI.12.4}.  By assumption,
$\Rm \Big|_{B_t(x, \widehat{r}(x,t))} \ge - \: \widehat{r}(x,t)^{-2}$ and
$\vol(B_t(x, \widehat{r}(x,t))) \: \ge \: w \: \widehat{r}(x,t)^3$.   As
$\widehat{r}(x,t) \ge \widehat{\rho} \sqrt{t}$, for any $c \in (0,1)$ we have
$\Rm \Big|_{B_t(x, c\widehat{\rho} \sqrt{t})} \ge - \: (c\widehat{\rho})^{-2} t^{-1}$. By
the Bishop-Gromov inequality,
\begin{align}
\vol(B_t(x, c\widehat{\rho} \sqrt{t})) \: & \ge \: 
\frac{
\int_0^{ \frac{c \widehat{\rho} \sqrt{t}}{\widehat{r}(x,t)}} \sinh^2(u) \: du}{\int_0^1 \sinh^2(u) \: du}
\: w \: \widehat{r}(x,t)^3
\: \ge \: \frac{1}{3\int_0^1 \sinh^2(u) \: du} \: w \: 
(c\widehat{\rho})^3 \: t^{\frac32} \\
& \ge \:  \frac{1}{10} \: w \: 
(c\widehat{\rho})^3 \: t^{\frac32}. \notag
\end{align}   
Considering Theorem \ref{thmI.12.3} with its $w$ replaced by $\frac{w}{10}$,     
if $c = c(w)$ is taken sufficiently small 
(to ensure $t \ge 4 \tau (c \widehat{\rho} \sqrt{t})^2$)
and $t$ is larger than a certain
$w$-dependent constant 
(to ensure 
$\frac{\Phi((c \widehat{\rho} \sqrt{t})^{-2})}{(c 
\widehat{\rho} \sqrt{t})^{-2}} <\rho$)
then we can apply Theorem \ref{thmI.12.3} with
$r_0 =  c \widehat{\rho} \sqrt{t}$ to obtain
$R(x^\prime,t^\prime) \le K^\prime(w) c^{-2} \widehat{\rho}^{-2} t^{-1}$  whenever
$t^\prime \in [t - \tau c^2 \widehat{\rho}^2 t, t]$ and
$\dist_{t^\prime}(x^\prime, x) \: \le \: \frac14 \: c \widehat{\rho} \sqrt{t}$. From the
$\Phi$-almost nonnegative curvature (see (\ref{double})), 
\begin{equation}
|\Rm|(x^\prime, t^\prime) \: \le \: \const \:
K^\prime c^{-2} \widehat{\rho}^{-2} t^{-1} \: + \: \const \: 
\Phi(K^\prime c^{-2} \widehat{\rho}^{-2} t^{-1}),
\end{equation}
which is bounded above by $2 \: \const \: K^\prime c^{-2} \widehat{\rho}^{-2} t^{-1}$
if $t$ is larger than a certain
$w$-dependent constant. Then from length and volume distortion estimates
for the Ricci flow,
we obtain a lower volume bound 
$\vol(B_{t^\prime}(x, c^\prime \widehat{\rho} \sqrt{t}))   \: \ge \:
w^{\prime} ( c^\prime \widehat{\rho} \sqrt{t})^3$ on a smaller ball of
controlled radius, for
some $c^\prime = c^\prime(w)$.  Using Theorem \ref{thmI.12.2}, we finally
obtain an upper bound $R \le K^{\prime \prime}(w) (c\widehat{\rho} \sqrt{t})^{-2}$
on  $B_t(x, c\widehat{\rho} \sqrt{t})$ and hence,
by the $\Phi$-almost nonnegative curvature, an upper bound of the form
$|\Rm| \le K(w) \: t^{-1}$ on $B_t(x, c\widehat{\rho} \sqrt{t})$,
provided that $t \ge T$ for an appropriate $T = T(w)$.
Taking $\overline{\rho} = c \widehat{\rho}$, the theorem follows.                                                                  
\end{proof}

\begin{remark}
The use of Theorem \ref{thmI.12.2} in the proof of Theorem \ref{thmI.13.1}
also gives an upper bound $|\Rm|(x,t) \le K(A,w) \: t^{-1}$ on
$B_t(x, A \overline{\rho} \sqrt{t})$ if $x \in M_{thick}(w,t)$, for any $A > 0$.
\end{remark}

We now take $w$ sufficiently small. Then for large $t$, $M_{thick}(w,t)$ has a boundary
consisting of tori that are incompressible in $M$ and the interior of
$M_{thick}(w,t)$ admits a complete Riemannian metric with constant
sectional curvature $- \: \frac14$ and finite volume; see
Sections 
\ref{hyprig} and \ref{inctori}. In addition,
$M_{thin}(w,t)$ is a graph manifold; see 
Section \ref{II.7.4}.

\section{Overview of 
{\em Ricci Flow with Surgery on 
Three-Manifolds}  \cite{Perelman2}} \label{overview2}

The paper \cite{Perelman2} is concerned with the Ricci flow on compact oriented
$3$-manifolds.  The main difference with respect to \cite{Perelman}
is that singularity formation is allowed, so the paper deals with 
a ``Ricci flow with surgery''.

The main part of the paper is concerned with
setting up the surgery procedure and showing that it is
well-defined, in the sense that surgery times do not
accumulate.  In addition, the long-time behavior of a
Ricci flow with surgery is analyzed.

The paper can be divided into three main parts.
Sections II.1-II.3 contain preparatory material about ancient
solutions, the so-called standard solution and the geometry
at the first singular time.  Sections II.4-II.5 set up the
surgery procedure and prove that it is well-defined.
Sections II.6-II.8 analyze the long-time behavior.

\subsection{II.1-II.3}

Section II.1 continues the analysis of three-dimensional
$\kappa$-solutions from I.11.  From I.11, any $\kappa$-solution
contains an ``asymptotic soliton'', a gradient shrinking
soliton that arises as a rescaled limit of the $\kappa$-solution
as $t \rightarrow - \infty$. It is shown that any such
gradient shrinking soliton must be a shrinking round cylinder $\R \times S^2$,
its $\Z_2$-quotient $\R \times_{\Z_2} S^2$
or a finite quotient of the
round shrinking $S^3$. Using this, one
obtains a finer description of the $\kappa$-solutions. In 
particular, any compact $\kappa$-solution must be isometric
to a finite quotient of the round shrinking $S^3$, or
diffeomorphic to $S^3$ or $\R P^3$. It is shown that there is
a universal number $\kappa_0 > 0$ so that any $\kappa$-solution
is a finite quotient of the round shrinking $S^3$ or is a
$\kappa_0$-solution. This implies universal
derivative bounds on the scalar curvature of a $\kappa$-solution.

Section II.2 defines and analyzes the Ricci flow of the
so-called standard solution.  This is a Ricci flow on $\R^3$ whose
initial metric is a capped-off half cylinder.  The surgery
procedure will amount to gluing in a truncated copy of the
time-zero slice of the standard solution. Hence one
needs to understand the Ricci flow on the standard solution
itself.  It is shown that the Ricci flow of the standard
solution exists on a maximal time interval $[0,1)$, and the
solution goes singular everywhere as $t \rightarrow 1$. 

The geometry of the solution at the first singular time $T$
(assuming that there is one) is considered in II.3.  Put
$\Omega \: = \: \{x \in M \: : \: \limsup_{t \rightarrow T^-} 
|\Rm(x,t)| < \infty \}$. Then $\Omega$ is an open subset of $M$,
and $x \in M-\Omega$ if and only if $\lim_{t \rightarrow T^-}
R(x,t) = \infty$. If $\Omega = \emptyset$ then for $t$ slightly less than $T$,
the manifold $(M, g(t))$ consists of nothing but high-scalar-curvature
regions. Using Theorem I.12.1, one shows that 
$M$ is diffeomorphic to 
$S^1 \times S^2$, $\R P^3 \# \R P^3$ or
a finite isometric quotient of $S^3$.

If $\Omega \neq \emptyset$ then there is a well-defined limit
metric $\overline{g}$ on $\Omega$, with scalar curvature
function $\overline{R}$. The set $\Omega$
could {\it a priori} have an infinite number of connected components,
for example if an infinite number of distinct $2$-spheres simultaneously
shrink to points at time $T$. For small $\rho > 0$,
put $\Omega_\rho \: = \: \{x \in \Omega \: : \: \overline{R}(x) \le \rho^{-2} \}$,
a compact subset of $M$.
The connected components of $\Omega$ can be divided into those
that intersect $\Omega_\rho$ and those that do not.  If a connected
component does not intersect $\Omega_\rho$ then it is a
``capped $\epsilon$-horn'' (consisting of a hornlike end capped off by a ball
or a copy of $\R P^3 - B^3$) or a ``double $\epsilon$-horn'' (with two
hornlike ends). If a connected component of $\Omega$ does intersect
$\Omega_\rho$ then it has a finite number of ends, each being
an $\epsilon$-horn. 

Topologically, the surgery procedure of II.4 will amount to taking each
connected component of $\Omega$ that intersects $\Omega_\rho$,
truncating each of its $\epsilon$-horns and gluing a $3$-ball onto
each truncated horn. 
The connected components of $\Omega$ that do not intersect
$\Omega_\rho$ are thrown away.
Call the new manifold $M^\prime$. At a time
$t$ slightly less than $T$, the region $M - \Omega_\rho$ consists
of high-scalar-curvature regions. Using the characterization of such
regions in I.12.1, one shows that $M$ can be reconstructed from
$M^\prime$ by taking the connected sum of its connected components,
along possibly with a finite number of $S^1 \times S^2$ and
$\R P^3$ factors.

\subsection{II.4-II.5}

Section II.4 defines the surgery procedure. 
A Ricci flow with surgery consists of a sequence of smooth $3$-dimensional Ricci flows on
adjacent time intervals with the property that for any two adjacent intervals, there is a
a compact $3$-dimensional submanifold-with-boundary that is common to the
final slice of the first time interval and the initial slice of the second time interval.

There are two {\it a priori} assumptions on a Ricci flow with surgery, the
pinching assumption and the
canonical neighborhood assumption.
The pinching assumption is a form of Hamilton-Ivey pinching. The canonical
neighborhood assumption says that every spacetime
point $(x,t)$ with $R(x,t) \ge r(t)^{-2}$ has a neighborhood which,
after rescaling, is $\epsilon$-close to one of the neighborhoods that occur
in a $\kappa$-solution or in a time slice of the standard solution. Here 
$\epsilon$ is a small but universal constant and $r(\cdot)$ is a decreasing
function, which is to be specified.

One wishes to define a Ricci flow with
surgery starting from any compact oriented $3$-manifold, say with a
normalized initial metric.  There are
various parameters that will enter into the definition : the above canonical neighborhood scale
$r(\cdot)$,
a nonincreasing function $\delta(\cdot)$ that decays to zero, the truncation
scale $\rho(t) = \delta(t) r(t)$ and the surgery scale $h$. In order
to show that one can construct the Ricci flow with surgery, it turns out that one wants to perform
the surgery only on necks with a radius that is very small compared to the
canonical neighborhood scale; this is the role of the parameter $\delta(\cdot)$.

Suppose that the Ricci flow with surgery is defined at times less than $T$, with
the {\it a priori} assumptions satisfied, and goes
singular at time $T$. Define the open subset $\Omega \subset M$ as before and 
construct the compact subset $\Omega_\rho \subset M$ 
using $\rho = \rho(T)$. Any connected component $N$ of $\Omega$ that
intersects $\Omega_\rho$ has a finite number of ends, each of which is an
$\epsilon$-horn.  This means that each point in the horn is in the center of
an $\epsilon$-neck, i.e. has a neighborhood
that, after rescaling, is $\epsilon$-close to a cylinder $\left[ - \frac{1}{\epsilon},
\frac{1}{\epsilon} \right] \times S^2$. In II.4.3 it is shown that as one goes down the
end of the horn, there is a self-improvement phenomenon; for any
$\delta > 0$, one can find $h < \delta \rho$ so that if a point $x$ in the horn has
$\overline{R}(x) \ge h^{-2}$ then it is actually in the center of a $\delta$-neck. 

With $\delta = \delta(T)$, let
$h$ be the corresponding number. One then
cuts off the $\epsilon$-horn at a $2$-sphere in the center of such a
$\delta$-neck and glues in a copy of a rescaled truncated standard solution.
One does this for each $\epsilon$-horn in $N$ and each connected component
$N$ that intersects $\Omega_\rho$, and throws away the connected components
of $\Omega$ that do not intersect $\Omega_\rho$. One lets the new manifold evolve 
under the Ricci flow. If one encounters another singularity then one again
performs surgery.  Based on an estimate on the volume change under a surgery,
one concludes that a finite number of surgeries occur in any finite time
interval. (However, 
one is not able to conclude from volume arguments that there is a finite
number of surgeries altogether.)

The preceding discussion was predicated on the condition that the
{\it a priori} assumptions hold for all times.  For the Ricci flow
before the first surgery time, the pinching condition
follows from the Hamilton-Ivey result.  One shows
that surgery can be performed so that it does not make the pinching any worse.
Then the pinching condition will hold up to the second surgery time, etc.
The main issue is to show that one can choose the parameters
$r(\cdot)$ and $\delta(\cdot)$ so that one knows {\it a priori} that the canonical 
neighborhood assumption, with parameter $r(\cdot)$, will hold for the
Ricci flow with surgery. (For any singularity time $T$, one needs to know that the
canonical neighborhood assumption holds for $t \in [0, T)$ in order to do the
surgery at time $T$.)

As a preliminary step, in Lemma II.4.5 it is shown that after one glues in a standard
solution, the result will still look similar to a standard solution, for as long of
a time interval as one could expect, unless the entire region gets removed by some
exterior surgery.

The result of II.5 is that the time-dependent
parameters $r(\cdot)$ and $\delta(\cdot)$ can be
chosen so as to ensure that the {\it a priori} assumptions hold.
The normalization of the initial metric implies that there is a time
interval $[0, C]$, for a universal constant $C$, on which the
Ricci flow is smooth and has explicitly bounded curvature.
On this time interval the canonical neighborhood assumption
holds vacuously, if $r(\cdot) \big|_{[0,C]}$ is sufficiently small.
To handle later times, the strategy is to divide $[\epsilon, \infty)$ into 
a countable sequence of finite time intervals and proceed by induction.
In II.5 the intervals $\{[2^{j-1} \epsilon, 2^j \epsilon]\}_{j=1}^\infty$ are used, 
although the precise choice of intervals is immaterial.  

We recall from I.12.1 that
in the case of smooth flows, the proof of the canonical neighborhood assumption
used the fact that one has $\kappa$-noncollapsing. It is not immediate that the
method of proof of I.12.1 extends to a Ricci flow with surgery. 
(It is exactly for this reason that one takes $\delta(\cdot)$ to be a 
time-dependent function which can be forced to be very small.)

Hence one needs to prove $\kappa$-noncollapsing and
the canonical neighborhoood assumption together. The main proposition of
II.5 says that there are decreasing sequences $r_j$, $\kappa_j$ and
$\overline{\delta}_j$ so that if $\delta(\cdot)$ is a function with
$\delta(\cdot) \big|_{[2^{j-1} \epsilon, 2^j \epsilon]} \le \overline{\delta}_j$ 
for each $j > 0$ then any Ricci flow with surgery, defined with the
parameters $r(\cdot)$ and $\delta(\cdot)$, is $\kappa_j$-noncollapsed on
the time interval $[2^{j-1} \epsilon, 2^j \epsilon]$ at scales less than
$\epsilon$ and satisfies the
canonical neighborhood assumption there. Here we take
$r(\cdot) \big|_{[2^{j-1} \epsilon, 2^j \epsilon)} = r_j$.

The proof of the proposition is by induction.  Suppose that it is true for
$1 \le j \le i$. In the induction step, 
besides defining the parameters $r_{i+1}$, $\kappa_{i+1}$ and
$\overline{\delta}_{i+1}$, one redefines $\overline{\delta}_i$. As one
only redefines $\overline{\delta}$ in the previous interval, there is no
circularity. 

The first step of the proof, Lemma II.5.2, consists of showing
that there is some $\kappa > 0$ so that for any $r$, one can
find $\overline{\delta} = \overline{\delta}(r) > 0$ with the following
property. Suppose that $g(\cdot)$ is a Ricci flow
with surgery defined on $[0, T)$, with $T \in [2^i \epsilon, 2^{i+1} \epsilon]$,
that satisfies the proposition on $[0, 2^i \epsilon]$.
Suppose that it also satisfies the canonical neighborhood assumption with
parameter $r$ on $[2^i \epsilon, T)$, and is constructed using a function
$\delta(\cdot)$ that satisfies $\delta(t) \: \le \: \overline{\delta}$ on
$[2^{i-1} \epsilon, T)$. Then it is 
$\kappa$-noncollapsed at all scales less than $\epsilon$. 

The proof of this
lemma is along the lines of the $\kappa$-noncollapsing result of I.7, with
some important modifications.  One again considers the ${\mathcal L}$-length of curves
$\gamma(\tau)$ starting from the point at which one wishes to prove the
noncollapsing. One wants to find a spacetime point
$(\overline{x}, \overline{t})$, with
$\overline{t} \in [2^{i-1} \epsilon, 2^i \epsilon]$, at which one has
an explicit upper bound on $l$. In I.7, the analogous statement came from
a differential inequality for $l$. In order to use this differential
equality in the present case, one needs to know that any curve $\gamma(\tau)$ 
that is competitive to be a minimizer for $L(\overline{x}, \overline{t})$ 
will avoid the surgery regions.
Choosing $\overline{\delta}$ small enough, one can ensure that the
surgeries in the time interval $[2^{i-1} \epsilon, T)$ are done on very
long thin necks. Using Lemma II.4.5, one shows that a curve $\gamma(\tau)$
passing near such a surgery region obtains a large value of ${\mathcal L}$,
thereby making it noncompetitive as a minimizer for 
$L(\overline{x}, \overline{t})$. (This is the underlying
reason that the surgery parameter
$\delta(\cdot)$ is chosen in a time-dependent way.)
One also chooses the point
$(\overline{x}, \overline{t})$ so that there is a small parabolic neighborhood
around it with a bound on its geometry. One can then run the argument of
I.7 to prove $\kappa$-noncollapsing in the time interval
$\left[ 2^i \epsilon, T \right)$.

The proof of the main proposition of II.5 is now by contradiction.  Suppose
that it is not true.  Then for some sequences $r^\al \rightarrow 0$ and
$\overline{\delta}^\al \rightarrow 0$, for each $\al$ there is a counterexample to the
proposition with $r_{i+1} \le r^\al$ and $\delta_i, \delta_{i+1} \le 
\overline{\delta}^\al$. That is, there is some spacetime point in the interval
$\left[ 2^i \epsilon, 2^{i+1} \epsilon \right]$ at which the
canonical neighborhood assumption fails.  Take a first such point 
$(x^\al, t^\al)$.
By Lemma II.5.2, one has $\kappa$-noncollapsing up to the time of this
first counterexample. Using this noncollapsing, one can consider
taking rescaled limits.  If there are no surgeries in an
appropriate-sized backward spacetime region around $(x^\al, t^\al)$
then one can extract a convergent subsequence as $\al \rightarrow \infty$
and construct, as in the proof of Theorem I.12.1, a limit $\kappa$-solution,
thereby giving a contradiction. If there are nearby interfering surgeries
then one argues, using Lemma II.4.5, that the point $(x^\al, t^\al)$ is in fact in a
canonical neighborhood, again giving a contradiction.

Having constructed the Ricci flow with surgery, if the initial
manifold is simply-connected then according to
\cite{Colding-Minicozzi,Colding-Minicozzi2,Perelman3}, 
there is a finite extinction time.
One then concludes that the Poincar\'e Conjecture holds.

\subsection{II.6-II.8}

Sections II.6 and II.8 analyze the large-time behavior of a Ricci flow with
surgery.

Section II.6 establishes back-and-forth curvature estimates. Proposition
II.6.3 is an analog of Theorem I.12.2 and Proposition II.6.4 is an analog
of Theorem I.12.3.  The proofs are along the lines of the proofs of
Theorems I.12.2 and I.12.3, but are complicated by the possible presence
of surgeries.

The thick-thin decomposition for large-time slices is considered in
Section II.7. Using monotonicity arguments of Hamilton, it is
shown that as $t \rightarrow \infty$ the metric on the $w$-thick part $M^+(w,t)$
becomes closer and closer to having constant negative sectional
curvature. Using a hyperbolic rigidity argument of Hamilton, it is stated
that the hyperbolic pieces stabilize in the sense that
there is a finite collection $\{ (H_i, x_i) \}_{i=1}^k$ of
pointed finite-volume $3$-manifolds of constant sectional
curvature $- \: \frac14$ so that for large $t$, the metric
$\widehat{g}(t) \: = \: \frac{1}{t} \: g(t)$ on the $w$-thick part
$M^+(w,t)$ approaches the metric on the $w$-thick part
of $\bigcup_{i=1}^k H_i$. It is stated that the cuspidal tori (if any) of
the hyperbolic pieces  are incompressible in $M$. To show this
(following Hamilton), if there is a compressing $3$-disk then one
takes a minimal such $3$-disk, say of area $A(t)$, and shows from a
differential inequality for $A(\cdot)$ that for large $t$ the function
$A(t)$ is negative, which is a contradiction.

Theorem II.7.4, a statement in Riemannian geometry, 
characterizes the thin part $M^-(w,t)$, for small $w$ and
large $t$, as a graph manifold.  The main hypothesis of the theorem
is that for each point $x$, there is a radius $\rho = \rho(x)$
so that the ball $B(x,\rho)$ has volume at most $w \rho^3$ and
sectional curvatures bounded below by $- \: \rho^{-2}$. In this
sense the manifold is locally volume collapsed with respect to a lower
sectional curvature bound. 

Section II.8 contains an alternative proof of the incompressibility
of cuspidal tori, using the functional $\lambda_1(g) = 
\lambda_1(-4 \triangle + R)$. 
(At the beginning of Section \ref{II.8}, we give a simpler argument
using the functional $R_{min}(g) \: \vol(M,g)^{\frac23}$.)
More generally, the functional $\lambda_1(g)$ is
used to define a topological invariant that determines the nature
of the geometric decomposition. First, the manifold $M$ admits a
Riemannian metric $g$ with $\lambda_1(g) > 0$ if and only if it
admits a Riemannian metric with positive scalar curvature, which
in turn is equivalent to saying that $M$ is diffeomorphic to a
connected sum of $S^1 \times S^2$'s and round quotients of
$S^3$. If $M$ does not admit a Riemannian metric with
$\lambda_1 > 0$, let $\overline{\lambda}$ be the supremum of
$\lambda_1(g) \cdot \vol(M, g)^{\frac23}$ over all Riemannian metrics $g$
on $M$. If $\overline{\lambda} = 0$ then $M$ is a graph manifold.
If $\overline{\lambda} < 0$ then the geometric decomposition of $M$
contains a nonempty hyperbolic piece, with total volume
$\left( - \: \frac23 \: \overline{\lambda} \right)^{\frac32}$. The proofs
of these statements use the monotonicity of
$\lambda(g(t)) \cdot \vol(M, g(t))^{\frac23}$, when it is nonpositive,
under a smooth Ricci flow. The main work is to show
that in the case of a Ricci flow with surgery, one can choose $\delta(\cdot)$
so that $\lambda(g(t)) \cdot \vol(M, g(t))^{\frac23}$ is
arbitrarily close to being nondecreasing in $t$.

\section{II. Notation and terminology}
\label{notationterminology}

$B(x,t,r)$ denotes the open metric ball of radius $r$, with respect to the metric
at time $t$, centered at $x$.

$P(x,t,r, \Delta t)$ denotes a parabolic neighborhood, that is the set of
all points $(x^\prime, t^\prime)$ with $x^\prime \in B(x,t,r)$ and
$t^\prime \in [t, t+\Delta t]$ or $t^\prime \in [t+\Delta t, t]$, depending on
the sign of $\Delta t$.

\begin{definition} \label{closenessdef}
We say that a Riemannian manifold
$(M_1, g_1)$ has distance $\le \: \epsilon$ in the $C^N$-topology 
to another Riemannian manifold $(M_2, g_2)$ if
there is a diffeomorphism $\phi \: : \: M_2 \rightarrow M_1$ so
that $\sum_{|I| \: \le \: N} \frac{1}{|I|!} \:
\parallel \nabla^I (\phi^* g_1 - g_2) \parallel_\infty \: \le 
\: \epsilon$.
An open set $U$ in a Riemannian $3$-manifold $M$ is an 
{\em $\eps$-neck}
if modulo rescaling, it has 
distance less than $\eps$,
in the $C^{[1/\epsilon]+1}$-topology, to the product of the 
round $2$-sphere of scalar curvature $1$ (and therefore
Gaussian curvature $\frac12$)
with an interval $I$ of length 
greater than
$2 \eps^{-1}$.
If a point $x \in M$ and a neighborhood $U$ of $x$ are specified then we 
will understand that ``distance''
refers to the pointed topology,
where the basepoint in $S^2 \times I$ projects to the center of $I$.

We make a similar definition of $\epsilon$-closeness in the
spacetime case, where $\nabla^I$ now includes time derivatives.
A subset of the form $U\times [a,b]\subset M\times [a,b]$,
where $U\subset M$ is open, sitting
in the spacetime of a Ricci flow is a {\em strong $\eps$-neck}
if after parabolic rescaling and time shifting, it has distance less
than $\eps$ to the
product Ricci flow defined on the time
interval $[-1,0]$ which,  at its final time,
is isometric to the  product of a  round $2$-sphere of scalar
curvature $1$ 
with an interval of length greater than $2\eps^{-1}$. (Evidently,
the time-$0$ slice of the product has $3$-dimensional scalar curvature
equal to $1$.) 
\end{definition}

Our definition of an $\epsilon$-neck 
differs in an insubstantial way from that on p. 1 of II.
In the definition of \cite{Perelman2}, 
a ball $B(x, t, \epsilon^{-1} r)$ is called
an $\epsilon$-neck if, after rescaling the metric with a factor
$r^{-2}$, it is $\epsilon$-close,
i.e. has distance less than $\epsilon$,
to {\em the corresponding subset of}
the standard neck 
$S^2 \times I$... 
(italicized words added by us). 
(The issue is that a large metric ball in the cylinder $\R \times S^2$
does not have a smooth boundary.)
Clearly after a slight change
of the constants, an $\epsilon$-neck in our sense is contained in an
$\epsilon$-neck in the sense of \cite{Perelman2}, 
and vice versa. An important fact
is that the notion of $(x, t)$ being contained in an $\epsilon$-neck
is an open condition with respect to the pointed $C^{[1/\epsilon]+1}$-topology
on Ricci flow solutions.

With an $\epsilon$-approximation 
$f \: : \: S^2 \times I \rightarrow U$ 
being understood, a cross-sectional sphere in $U$ will mean the
image of $S^2 \times \{ \lambda \}$ under $f$, for some
$\lambda \in (- \epsilon^{-1}, \epsilon^{-1})$. 
Any curve $\gamma$ in $U$ that intersects both 
$f(S^2 \times \{ \epsilon^{-1} \})$ and  $f(S^2 \times \{ -\epsilon^{-1} \})$ must
intersect each cross-sectional
sphere. If $\gamma$ is a minimizing  geodesic and $\epsilon$ is small
enough then $\gamma$ will intersect each cross-sectional sphere
exactly once.

There is a typo in the definition of a strong $\eps$-neck in \cite{Perelman2} :
the parabolic neighborhood should be $P(x,t,\eps^{-1}r,-r^2)$,
i.e. it should go backward in time rather than forward. We note that the
time interval involved in the
definition of strong $\eps$-neck, i.e. $1$ after rescaling, is different than the 
rescaled time interval $\epsilon^{-1}$ in Theorem \ref{thmI.12.1}.

In the next definition, $\I$ is an open interval and  $B^3$ is an open ball.  

\begin{definition}
A metric on 
$S^2 \times \I$ such that each point is contained in some $\epsilon$-neck
is called an {\em $\epsilon$-tube}, or an {\em $\epsilon$-horn}, or a {\em double
$\epsilon$-horn}, if the scalar curvature stays bounded on both ends,
stays bounded on one end and tends to infinity on the other, or
tends to infinity on both ends, respectively.

A metric on $B^3$ or $\R P^3 - \overline{B^3}$, such that each point
outside some compact subset is contained in an $\epsilon$-neck, is
called an $\epsilon$-cap or a capped $\epsilon$-horn, if the
scalar curvature stays bounded or tends to infinity on the end,
respectively.
\end{definition}

An example of an $\epsilon$-tube is $S^2 \times (- \epsilon^{-1}, \epsilon^{-1})$ with the
product metric. 
For a relevant example of an $\epsilon$-horn, consider the metric
\begin{equation}
g \: = \: dr^2 \: + \: \frac{1}{8\ln \frac{1}{r}} \: r^2 \: d\theta^2
\end{equation}
on $(0, R) \times S^2$, where $d\theta^2$ is the metric on $S^2$ with $R=1$.
From \cite{Angenent-Knopf}, the metric $g$ models a rotationally symmetric
neckpinch. Rescaling around $r_0$, we put
$s \: = \: \sqrt{8 \ln \frac{1}{r_0}} \: \left( \frac{r}{r_0} - 1 \right)$ and find
\begin{equation}
\frac{8 \ln \frac{1}{r_0} }{r_0^2} \: g \: = \:
ds^2 \: + \: \left( 1 + 
{ \frac{1}{\sqrt{8 \ln \frac{1}{r_0} }}  }
\left( 2 + 
\frac{1}{\ln \frac{1}{r_0}} 
\right) s + O(s^2) \right) \: d\theta^2.
\end{equation}
For small $r_0$, if we take
$\epsilon \sim \left( \ln \frac{1}{r_0} \right)^{- \: \frac14}$ then the region with 
$s \in \left( - \: \epsilon^{-1}, \epsilon^{-1} \right)$
will be $\epsilon$-biLipschitz close
to the standard cylinder. Note that as $r_0 \rightarrow 0$,
the constant $\epsilon$ improves; this is related to
Lemma \ref{lemmaII.4.3}.

An $\epsilon$-cap is the result of capping off an $\epsilon$-tube by a
$3$-ball or $\R P^3 - \overline{B^3}$
with an arbitrary metric.  A capped $\epsilon$-horn
is the result of capping off an $\epsilon$-horn by a
$3$-ball or $\R P^3 - \overline{B^3}$ with an arbitrary metric.

\begin{remark}
\label{remepsilon}
Throughout the rest of these notes, $\epsilon$ denotes a small
positive constant that is meant to be universal.  The precise value
of $\epsilon$ is unspecified.  If the statement of a lemma or
theorem invokes $\epsilon$ then the statement is meant to be true
uniformly with respect to the other variables, provided $\eps$ is
sufficiently small. When going through
the proofs one is allowed to make $\epsilon$ small enough
so that the arguments work, but one is only allowed to make a
finite number of such reductions.  
\end{remark}

\begin{lemma}
Let $U$ be
an $\epsilon$-neck in an $\epsilon$-tube (or horn) and let 
$S$ be a cross-sectional
sphere in $U$. Then $S$ separates the two ends of the tube (or horn).
\end{lemma}
\begin{proof}
Let $W$ denote the tube (or horn).
As any point $m \in W$ lies in some $\epsilon$-neck,
there is a unique lowest eigenvalue of the Ricci operator
$\Ric  \in \End(T_m W)$ at $m$. Let $\xi_m \subset T_mW$ be the
corresponding eigenspace.  As $m$ varies, the 
$\xi_m$'s form a smooth
line field $\xi$ on $W$, to which $S$ is transverse.

Suppose that $S$ does not separate the two ends of $W$. 
Then $S$ represents a trivial element of
$\HH_2(S^2 \times I) \cong \pi_2(S^2 \times I)$ and
there is an embedded $3$-disk $D \subset W$ for which
$\partial D = S$. This contradicts the fact that
the line field $\xi$ is transverse to $S$ and extends over $D$.
\end{proof}

\section{II.1. Three-dimensional $\kappa$-solutions} \label{II.1}

This section is concerned with properties of three-dimensional oriented
$\kappa$-solutions. For brevity, in the rest of these notes we will generally omit 
the phrases ``three-dimensional'' and ``oriented''.

If $(M, g(\cdot))$ is a $\kappa$-solution then its
topology is easy to describe.  By definition, $(M, g(t))$ has nonnegative
sectional curvature.  If it does not have strictly positive curvature
then the universal cover splits off a line (see Theorem \ref{splitting}), 
from which
it follows (using Corollary \ref{2Dsoliton} and the $\kappa$-noncollapsing)
that $(M, g(\cdot))$ is a standard shrinking cylinder
$\R \times S^2$ or its $\Z_2$ quotient $\R \times_{\Z_2} S^2$.
If $(M, g(t))$ has strictly positive curvature and $M$ is compact
then it is diffeomorphic to a spherical space form
\cite{Hamiltonnnnn}. 
If $(M, g(t))$ has strictly positive curvature and $M$ is noncompact
then it is diffeomorphic to $\R^3$ \cite{Cheeger-Gromoll}. The lemmas
in this section give more precise geometric information.
Recall that $M_\epsilon$ consists of the points in a $\kappa$-solution
which are not the center of an $\epsilon$-neck.

\begin{lemma} \label{noncompactsoliton}
If $(M, g(t))$ is a time slice of a
noncompact $\kappa$-solution and
$M_\epsilon \neq \emptyset$ then there is a compact
submanifold-with-boundary $X \subset M$ so that
$M_\epsilon \subset X$,  $X$ is diffeomorphic to $\overline{B^3}$ or
$\R P^3 - B^3$, and $M- \Int(X)$ is
diffeomorphic to $[0, \infty) \times S^2$.
\end{lemma}
\begin{proof}
If $(M,g(t))$ does not have positive sectional curvature and
$M_\epsilon \neq \emptyset$ then $M$ must be isometric to
$\R \times_{\Z_2} S^2$, in which case the lemma is easily seen to
be true with $X$ diffeomorphic to $\R P^3 - B^3$. 
Suppose that $(M, g(t))$ has positive sectional curvature.
Choose $x \in M_\epsilon$. Let $\gamma \: : \: [0, \infty) 
\rightarrow M$ be a ray with $\gamma(0) = x$. As $M_\epsilon$ is 
compact, there is some $a > 0$ so that if $t > a$ then
$\gamma(t) \notin M_\epsilon$. 
We can cover $(a, \infty)$ by open
intervals $V_j$ so that $\gamma \Big|_{V_j}$ is a geodesic segment
in an $\epsilon$-neck of 
rescaled length approximately
$2 \epsilon^{-1}$.  
Then we can find a cover of $(a, \infty)$ by 
linearly ordered open intervals
$U_i$, refining the previous cover, so that \\
1. The rescaled length of $\gamma \Big|_{U_i}$ is approximately
$\frac{1}{10} \: \epsilon^{-1}$. \\
2. Choosing some $x_i \in U_i \cap U_{i+1}$, the
rescaled length (with rescaling at $x_i$) of
$\gamma \Big|_{U_i \cap U_{i+1}}$ is approximately
$\frac{1}{40} \: \epsilon^{-1}$ and 
$\gamma \Big|_{U_i \cap U_{i+1}}$ lies in an $\epsilon$-neck $W_i$
centered at $x_i$.

Let $\phi_i$
be projection on the first factor in the assumed diffeomorphism
$W_i \cong  (- \epsilon^{-1}, \epsilon^{-1}) \times S^2$.
If $\epsilon$ is sufficiently small then
the composition $\phi_i \circ \gamma \big|_{U_i}  \: : \: U_i \rightarrow
(- \epsilon^{-1}, \epsilon^{-1})$ is a diffeomorphism onto its image.
Put $N_i = \phi_i^{-1}(\Image(\phi_i \circ \gamma \big|_{U_i})) \subset
W_i$ and
$p_i \: = \: (\phi_i \circ \gamma \big|_{U_i})^{-1} \circ \phi_i \: : \:
N_i \rightarrow U_i \subset (a, \infty)$. Then $p_i$ is $\epsilon$-close
to being a Riemannian submersion and
on overlaps $N_i \cap N_{i+1}$, the maps $p_i$ and $p_{i+1}$ are
$C^K$-close. Choosing an appropriate
partition of unity $\{b_i\}$ subordinate to the $U_i$'s, if $\epsilon$ is
small then the function
$f = \sum_i b_i p_i$ is a submersion from $\bigcup_i N_i$ to 
$(a, \infty)$. The fiber is seen to be $S^2$. Given $t \in (a, \infty)$,
put $X = M - f^{-1}(t, \infty)$. Then $M - \Int(X) = f^{-1}([t, \infty))$ 
is diffeomorphic to $[0, \infty) \times S^2$. 

Recall from Section \ref{pf11.7} that we have an exhaustion of
$M$ by certain convex compact subsets.  As $M$ is one-ended, the
subsets have connected boundary.
As in Section \ref{pf11.7}, if the boundary of such a subset intersects an
$\epsilon$-neck then the intersection will be a nearly
cross-sectional $2$-sphere in the $\epsilon$-neck. Hence with
an appropriate choice of $t$, the set $X$ will be isotopic to one
of our convex subsets and so diffeomorphic to a closed $3$-ball.
\end{proof}

\begin{lemma} \label{onlyneck}
If $(M, g(t))$ is a time slice of a $\kappa$-solution with
$M_\epsilon = \emptyset$ then the
Ricci flow is the evolving round cylinder $\R \times S^2$.
\end{lemma}
\begin{proof}
By assumption, each point $(x,t)$ lies in an $\epsilon$-neck. 
If $\epsilon$ is sufficiently small then piecing
the necks together, we conclude that $M$ must be diffeomorphic to
$S^1 \times S^2$ or $\R \times S^2$; see the proof of
Lemma \ref{noncompactsoliton} for a similar argument. Then the universal cover 
$\widetilde{M}$ is
$\R \times S^2$. As it has nonnegative sectional curvature and
two ends, Toponogov's theorem implies that $(\widetilde{M}, 
\widetilde{g}(t))$ splits off
an $\R$-factor. Using the strong maximum principle, the Ricci flow
on $\widetilde{M}$
splits off an $\R$-factor; see Theorem \ref{splitting}.
Using Corollary \ref{2Dsoliton}, it follows that $(\widetilde{M}, 
\widetilde{g}(t))$ is the evolving round cylinder $\R \times S^2$.
From the $\kappa$-noncollapsing, the quotient $M$ cannot be
$S^1 \times S^2$.
\end{proof}

A $\kappa$-solution has an asymptotic soliton
(Section \ref{I.11.2}) that is either compact or noncompact.
If the asymptotic soliton of a compact $\kappa$-solution
$(M, g(\cdot))$ is also compact then it must be a shrinking 
quotient of the round $S^3$ \cite{Hamiltonnnnn}, so the same is true of
$M$.  

\begin{lemma} \label{S3RP3}
If a $\kappa$-solution $(M, g(\cdot))$
is compact and has a noncompact asymptotic soliton 
then $M$ is diffeomorphic to $S^3$ or $\R P^3$.
\end{lemma}
\begin{proof}
We use Corollary \ref{neckstructure} in Section \ref{strongversion}.
First, we claim that the time slices of the type-D $\kappa$-solutions of
Corollary \ref{neckstructure} have a
universal upper bound on $\max_M R \cdot \diam(M)^2$.
To see this, we can rescale at the point $x \in M_\epsilon$ by
$R(x)$, after which the diameter is bounded above by
$2 \sqrt{\alpha}$. We then use Theorem \ref{thmI.11.7}
to get an upper bound on the rescaled
scalar curvature, which proves the claim. Given an upper bound on 
$\max_M R \cdot \diam(M)^2$, the asymptotic soliton
cannot be noncompact.

Thus we are in case C of Corollary \ref{neckstructure}.
Take a sequence $t_i \rightarrow - \infty$ and choose points $x_i, y_i
\in M_\epsilon(t_i)$ as in Corollary \ref{neckstructure}.C.
Rescale by
$R(x_i, t_i)$ and take a subsequence that converges to a pointed
Ricci flow solution
$(M_\infty, (x_\infty, t_\infty))$. The limit $M_\infty$ cannot be
compact, as otherwise we would have a uniform upper bound on 
$R \: \cdot \: \diam^2$ for $(M, g(t_i))$, which would
contradict the existence of the noncompact asymptotic soliton. 
Thus $M_\infty$ is a noncompact $\kappa$-solution. 
We can find compact sets $X_i \subset M$ 
containing $B(x_i, \alpha R(x_i, t_i)^{- \: \frac12})$ so that
$\{X_i\}$ converges to a set $X_\infty \subset M_\infty$ as in
Lemma \ref{noncompactsoliton}. 
Taking a further subsequence, we find similar compact 
sets $Y_i \subset M$ containing
$B(y_i, \alpha R(y_i, t_i)^{- \: \frac12})$ so that $\{Y_i\}$ converges to
a set $Y_\infty \subset M_\infty$ as in Lemma
\ref{noncompactsoliton}. In particular, for large $i$, $X_i$ and $Y_i$ are
each diffeomorphic to either $\overline{B^3}$ or $\R P^3 - B^3$.
Considering a minimizing geodesic
segment $\overline{x_i y_i}$ as in the statement of
Corollary \ref{neckstructure}.C, we
can use an argument as in the proof of Lemma 
\ref{noncompactsoliton} to construct a submersion from 
$M - (X_i \cup Y_i)$ to an interval, with fiber $S^2$. Hence $M$
is diffeomorphic to the result of gluing $X_i$ and $Y_i$ along
a $2$-sphere.
As $M$ has finite fundamental group, no more
than one of $X_i$ and $Y_i$ can be diffeomorphic to 
$\R P^3 \: - \: {B^3}$. Thus $M$ is diffeomorphic to
$S^3$ or $\R P^3$.
\end{proof}

From Lemma \ref{propI.11.9}, every ancient solution which is a 
$\kappa$-solution {\em for some $\kappa$} is either a $\kappa_0$-solution
or a metric quotient of the round $S^3$. 

\begin{lemma} \label{derivativeestimates}
There is a universal constant
$\eta$ such that at each point of every ancient solution that is
a $\kappa$-solution for some $\kappa$, we have estimates
\begin{equation} \label{derestimates}
|\nabla R| < \eta R^{\frac32}, \: \: \: \: \: \: |R_t| < \eta R^2.
\end{equation}
\end{lemma}
\begin{proof}
This is obviously true for metric quotients of the round $S^3$. For
$\kappa_0$-solutions it follows from the compactness result in
Theorem \ref{thmI.11.7}, after rescaling the scalar curvature at
the given point to be $1$.
\end{proof}

It is sometimes useful to rewrite (\ref{derestimates}) as a pair of estimates
on the spacetime derivatives of the quantity $R^{-1}$ at points where $R\neq 0$:
\begin{equation}
\label{eqnderestimatesrinverse}
|\nabla (R^{-\frac12})| < \frac{\eta}{2}, \: \: \: \: \: \: |(R^{-1})_t| < \eta .
\end{equation}

\begin{lemma} \label{nbhd}
For every sufficiently small $\epsilon > 0$ one can
find $C_1 = C_1(\epsilon)$ and $C_2 = C_2(\epsilon)$
such that for each point $(x, t)$ in every
$\kappa$-solution there is a radius 
$r \in [ R(x,t)^{-1/2}, C_1 R(x,t)^{-1/2}]$ and a neighborhood
$B$, 
$\overline{B(x,t,r)} \subset B \subset B(x,t,2r)$, 
which falls into one
of the four categories: \\
(a) $B$ is a strong $\epsilon$-neck 
(more precisely, $B$ is the slice of a
strong $\epsilon$-neck at its maximal time, and an appropriate parabolic
neighborhood of $B$
satisfies the condition to be a strong $\epsilon$-neck),
or \\
(b) $B$ is an $\epsilon$-cap, or \\
(c) $B$ is a closed manifold, diffeomorphic to $S^3$ or $\R P^3$, or\\
(d) $B$ is a closed manifold of constant positive sectional curvature.

Furthermore:

\begin{itemize}
\item The scalar curvature in $B$ at time $t$ is between
$C_2^{-1} R(x,t)$ and $C_2 R(x,t)$.

\item The volume of $B$ in cases (a), (b) and (c) is
greater than $C_2^{-1} R(x,t)^{- \: \frac32}$.

\item In case (b), there is an $\eps$-neck $U\subset B$ with compact complement in $B$
(i.e. the end of $B$ is entirely contained in the $\eps$-neck) such that the 
distance from $x$ to $U$ is at least $10000R(x,t)^{-1/2}$.

\item In case (c) the 
sectional curvature in $B$ at time $t$ is greater than $C_2^{-1} R(x,t)$.
\end{itemize}
\end{lemma}

\begin{remark}
The statement of the lemma is slightly stronger than the corresponding
statement in II.1.5, in that
we have $r \ge  R(x,t)^{-1/2}$ as opposed to $r > 0$.
\end{remark}

\begin{proof}
We may assume that we are talking about a 
$\kappa_0$-solution, as if $M$ is a metric quotient of a round sphere then
it falls into category (d) for any $r > \pi (R(x_i, t_i)/6)^{-1/2}$
(since then 
$M \: = \: \overline{B(x_i, t_i, r)} \: = \: B(x_i, t_i, 2r)$). 

Fix a small $\epsilon$ and suppose that the claim
is not true.  Then there is a sequence of $\kappa_0$-solutions $M_i$ that
together provide a counterexample.  That is, there is a sequence
$C_i \rightarrow \infty$ and a sequence of points $(x_i, t_i) \in
M_i \times (-\infty, 0]$ so that for any $r \in
[R(x_i, t_i)^{-1/2}, C_i R(x_i, t_i)^{-1/2} ]$ one cannot
find a $B$ between 
$\overline{B(x_i, t_i, r)}$
and $B(x_i, t_i, 2r)$ falling into
one of the four categories
and satisfying the subsidiary conditions with parameter $C_2 = C_i$.
Rescale the metric by $R(x_i, t_i)$ and take 
a convergent subsequence of 
$(M_i, (x_i, t_i))$ to obtain a limit $\kappa_0$-solution
$(M_\infty, (x_\infty, t_\infty))$. Then for any $r > 1$,
one cannot
find a $B_\infty$ between 
$\overline{B(x_\infty, t_\infty, r)}$
and 
$B(x_\infty, t_\infty, 2r)$ falling into
one of the four categories and satisfying the subsidiary conditions
for any parameter $C_2$.

If $M_\infty$ is compact then
for any $r$ greater than the diameter of the 
time-$t_\infty$ slice of $M_\infty$, 
$\overline{B(x_\infty, t_\infty, r)} = M_\infty = B(x_\infty, t_\infty, 2r)$
falls
into category (c) or (d). For the subsidiary conditions,
$M_\infty$ clearly has a lower volume bound, a positive lower
scalar curvature bound and an upper scalar
curvature bound. As a compact
$\kappa_0$-solution has positive sectional curvature, 
$M_\infty$ also has a lower sectional curvature bound. 
This is a contradiction.

If $M_\infty$ is noncompact then Lemma
\ref{noncompactsoliton} (or more precisely its proof) and 
Lemma \ref{onlyneck}
imply that for some $r>1$, there will be a $B$ between 
$\overline{B(x_\infty, t_\infty, r)}$
and $B(x_\infty, t_\infty, 2r)$
falling into category (a) or (b).
In case (b), by choosing the parameter $r$ sufficiently large,
the existence of the $\eps$-neck $U$
with the desired properties follows from the proof of 
Lemma \ref{noncompactsoliton}.
For the other subsidiary conditions,
$B$ clearly has a lower volume bound, a positive lower
scalar curvature bound and an upper scalar
curvature bound. 
This contradiction completes the proof of the lemma.
\end{proof}

\section{II.2. Standard solutions} \label{II.2}

The next few sections are concerned with the properties of  special
Ricci flow solutions on $M=\R^3$. 
We fix a smooth rotationally symmetric   metric $g_0$ which
is the result of gluing a hemispherical-type cap to a 
half-infinite cylinder of scalar curvature $1$. 
Among other properties, $g_0$ is complete and has
nonnegative  curvature operator.
We also assume that $g_0$ has scalar curvature
bounded below by $1$.

\begin{remark}
In Section \ref{surgery} we will further specialize the initial metric
$g_0$ of the standard solution, for technical convenience in
doing surgeries.
\end{remark}

\begin{definition}
\label{defstandardsolution}
A Ricci flow $(\R^3,g(\cdot))$  defined on a time interval $[0,a)$
is a {\em standard solution} if
it has  initial condition $g_0$, the curvature $|\Rm|$ is bounded 
on compact time intervals $[0,a']\subset [0,a)$, and it cannot
be extended to a Ricci flow with the same properties 
on a strictly longer time interval. 
\end{definition}

It will turn out that every standard solution 
is defined on the time interval $[0,1)$.

To motivate the next few sections, let us mention that the surgery
procedure will amount to gluing in a truncated copy of 
$(\R^3, g_0)$. The metric on this added region will then evolve as 
part of the Ricci flow that takes up after the surgery is performed. 
We will need to understand
the behavior of the Ricci flow after performing a surgery. Near
the added region, this will be modeled by a standard solution. 
Hence one first needs to understand the Ricci flow of a standard
solution.

The main results of II.2 concerning the Ricci flow on a standard solution
(Sections \ref{claim2}-\ref{standardlycompact}) are used in
II.4.5 (Lemma \ref{II.4.5}) to show that,
roughly speaking, the part of the manifold added by surgery 
acquires a large scalar curvature soon after the surgery time.  This
is used crucially in II.5 (Sections \ref{II.5.2section} and \ref{surgeryflow})
to adapt the noncollapsing argument of I.7 to Ricci flows
with surgery.

Our order of presentation of the material in 
II.2 is somewhat different than that of \cite{Perelman2}. 
In Sections \ref{claim2}-\ref{secclaim5}, we 
cover Claims 2, 4 and 5 of II.2. These are what's needed in the
sequel.  The other results of II.2, Claims 1 and 3, are concerned with proving the
uniqueness of the standard solution. Although it may seem
intuitively obvious that there should be a unique and
rotationally-symmetric standard solution, the argument
is not routine since the manifold is noncompact.

In fact, the uniqueness
is not really needed for the sequel. (For example, the method of proof of
Lemma \ref{II.4.5} produces {\em a} standard solution in a
limiting argument and it is enough to know certain properties
of this standard solution.) Because of this we will talk
about {\em a} standard solution rather than {\em the}
standard solution.

Consequently, we present the material so that we do not logically
need the uniqueness of the standard solution.
Having uniqueness does not
shorten the subsequent arguments any.
Of course, one can ask independently whether the standard
solution is unique. In Section \ref{claimII.2.1} we
show that a standard solution is rotationally symmetric.
In Section \ref{claim3} we sketch the argument for uniqueness.
Papers concerning the uniqueness of the standard solution are
\cite{Chen-Zhu,Lu-Tian}.

We end this section by collecting some basic facts about standard solutions.

\begin{lemma}
\label{lemstandard0}
Let $(\R^3,g(\cdot))$ be a standard solution.  Then 

(1) The curvature operator of $g$ is nonnegative.

(2) All derivatives of curvature are bounded for small time,
independent of the standard solution.

(3) The scalar curvature satisfies $\lim_{t \rightarrow a^-}
\sup_{x \in \R^3} R(x,t) \: = \: \infty$.

(4) $(\R^3,g(\cdot))$ is $\kappa$-noncollapsed at scales below $1$ on any time
interval contained in $[0,2]$, where
$\kappa$ depends only on the choice of the initial condition $g_0$.

(5) $(\R^3,g(\cdot))$ satisfies the conclusion of Theorem \ref{thmI.12.1},
in the sense that for any $\Delta t > 0$, there is an $r_0 > 0$ so that
for any point $(x_0, t_0)$ with $t_0 \ge \Delta t$ and
$Q = R(x_0, t_0) \ge r_0^{-2}$, the solution in
$\{(x,t) \: : \: \dist_{t_0}^2(x, x_0) < (\epsilon Q)^{-1}, t_0 - (\epsilon Q)^{-1}
\le t \le t_0\}$ is, after scaling by the factor $Q$, $\epsilon$-close to the
corresponding subset of a $\kappa$-solution. 

Moreover, any Ricci flow which satisfies all of the conditions of Definition 
\ref{defstandardsolution} except maximality of the time interval can be
extended to a standard solution. In particular, using short-time
existence \cite[Theorem 1.1]{Shi2 (1989)}, there is at least one standard solution.
\end{lemma}
\begin{proof}
(1) follows from \cite[Theorem 4.14]{Shi3 (1989)}.

(2) follows from Appendix \ref{applocalder}.

(3) 
In view of (1), this is equivalent to saying that
$\lim_{t \rightarrow a^-}
\sup_{x \in \R^3} |\Rm|(x,t) \: = \: \infty$. The argument for
this last assertion
is in \cite[Chapter 6.7.2]{Chow-Knopf}. The proof in 
\cite[Chapter 6.7.2]{Chow-Knopf}
is for the compact case but using the derivative estimates of
Appendix \ref{applocalder}, the same argument works in 
the present case.

(4) See Theorem \ref{nolocalcollapse}.

(5) See Theorem \ref{thmI.12.1}.

The final assertion of the lemma follows from the method of
proof of (3).
\end{proof}

\section{Claim 2 of II.2.
The blow-up time for a standard solution is $\leq 1$
} \label{claim2}

\begin{lemma} (cf. Claim 2 of II.2) \label{lemclaim2}
Let $[0,T_S)$ be the maximal time interval such that 
the curvature of all standard solutions is uniformly
bounded for every compact subinterval $[0,a]\subset [0,T_S)$.
Then on the time interval $[0,T_S)$, the family of standard solution converges 
uniformly at (spatial) infinity  to the standard Ricci flow 
on the round infinite cylinder $S^2\times \R$ of scalar curvature one.
In particular, $T_S$ is  at most $1$.
\end{lemma}
\begin{proof}
Let $\{\M_i\}_{i=1}^\infty$ be a sequence of standard solutions, and
let $\{x_i\}_{i=1}^\infty$ be a sequence tending to infinity
in the time-zero slice $M$.

By (2) of Lemma \ref{lemstandard0}, the gradient estimates in Appendix \ref{applocalder},
 and Appendix \ref{subsequence}, every subsequence
of $\{\M_i, (x_i,0)\}_{i=1}^\infty$ has a subsequence which
converges in the pointed smooth topology on the time interval $[0,T_S)$.
Therefore, it suffices to show that if $\{\M_i, (x_i,0)\}_{i=1}^\infty$
converges to some pointed Ricci flow
$(\M_\infty, (x_\infty, 0))$ then 
$\M_\infty$ is round cylindrical flow.

Since $g_i(0)=g_0$ for all $i$,  the sequence of pointed time-zero slices
$\{(M,x_i, g_i(0))\}_{i=1}^\infty$ converges in the pointed smooth topology
to the round cylinder, i.e. $(M_\infty,g_\infty(0))$ is a round cylinder
of scalar curvature $1$. Each time slice $(M_\infty,g_\infty(t))$
is biLipschitz equivalent to $(M_\infty,g_\infty(0))$. In particular,
it has two ends. As it also has nonnegative sectional curvature,
Toponogov's theorem implies that $(M_\infty,g_\infty(t))$ splits off
an $\R$-factor. Using the strong maximum principle, 
the Ricci flow $\M_\infty$ splits off an $\R$-factor;
see Theorem \ref{splitting}.
Then using the uniqueness of the Ricci flow on the round $S^2$, it
follows that $\M_\infty$ is a standard shrinking cylinder, which proves the lemma.

In particular, $T_S\leq 1$. 
\end{proof}

\section{Claim 4 of II.2. 
The blow-up time of a standard solution is $1$}
\label{secclaim4}

\begin{lemma} (cf. Claim 4 of II.2) \label{lemclaim4}
Let $T_S$ be as in  Lemma \ref{lemclaim2}. Then $T_S=1$.
In particular, every standard solution survives until time $1$.
\end{lemma}
\begin{proof}
First, there is an 
$\alpha > 0$ 
so that $T_S > \alpha$ \cite[Theorem 1.1]{Shi2 (1989)}.
In what follows we will apply Theorem \ref{thmI.12.1}. The 
hypothesis of Theorem \ref{thmI.12.1} says that the flow
should exist on a time interval of duration at least one,
but by 
rescaling we can apply Theorem \ref{thmI.12.1} just as well
with the alternative hypothesis that the
flow exists on a time interval of duration
at least $\alpha$.

Suppose that $T_S < 1$. 
Then there is a sequence of standard solutions
$\{\M_i\}_{i=1}^\infty$, 
times $t_i \rightarrow T_S$ and points $(x_i,t_i)\in \M_i$ so that
$\lim_{i \rightarrow \infty} R(x_i, t_i) = \infty$. 

We first
argue that no subsequence of the points $x_i$ can go to infinity (with respect
to the time-zero slice). Suppose, after relabeling the subsequence, that
$\{x_i\}_{i=1}^\infty$ goes to infinity.
From Lemma \ref{lemclaim2}, 
for any fixed $t^\prime < T_S$ the pointed solutions 
$(M,(x_i, 0), g_i(\cdot))$, defined for $t \in [0, t^\prime]$, 
approach that of the shrinking cylinder on the same time interval.
Lemma \ref{lemstandard0} and the characterization of high-curvature regions from
Theorem \ref{thmI.12.1} implies a
uniform bound on high-curvature regions of the time derivative of $R$, of the form
(\ref{derestimates}). Then taking $t^\prime$
sufficiently close to $T_S$, we get a contradiction.
We conclude that outside of a compact
region the curvature stays uniformly bounded as $t \rightarrow T_S$;
compare with the proof of Lemma \ref{claim1}.
(Alternatively, one could apply Theorem \ref{thmI.10.1} to compact
approximants, as is done in II.2.)

Thus we may assume that
the sequence $\{x_i\}_{i=1}^\infty$ stays in a compact region of the
time-zero slice. By Theorem \ref{thmI.12.1}, 
there is a sequence $\epsilon_i \rightarrow 0$
so that after rescaling the pointed solution
$(\M, (x_i, t_i))$ by $R(x_i, t_i)$, the result is $\epsilon_i$-close
to the corresponding subset of an ancient solution. By 
Proposition \ref{propI.11.4}, 
the ancient solutions have vanishing asymptotic volume ratio.
Hence for every $\beta > 0$, there is some $L < \infty$ so that
in the original unscaled solution, for large $i$ we have
$\vol \left( B(x_i, t_i, L \: R(x_i, t_i)^{- \: \frac12} ) \right)
\: \le \: \beta \left( L \: R(x_i, t_i)^{- \: \frac12} \right)^3$.
Applying the Bishop-Gromov inequality to the time-$t_i$ slices, we
conclude that for any $D >0$,
$\lim_{i \rightarrow \infty}
D^{-3} \vol(B(x_i, t_i, D)) \: = \: 0$.
However, this contradicts the
previously-shown fact that the solution extends smoothly to time $T_S < 1$
outside of a compact set.

Thus $T_S=1$. 
\end{proof}

\begin{lemma} \label{scalarblowup}
The infimal scalar curvature on the time-$t$ slice
tends to infinity as $t\ra 1^-$ uniformly for all standard solutions.
\end{lemma}
\begin{proof}
Suppose the lemma failed and 
let $\{(\M_i,(x_i,t_i))\}_{i=1}^\infty$ 
be a sequence of pointed standard solutions,
with 
$\{R(x_i, t_i)\}_{i=1}^\infty$ 
uniformly bounded and $\lim_{i \rightarrow \infty} t_i = 1$.

Suppose first that after passing to a subsequence, 
the points $x_i$ go to infinity in the time-zero
slice. From Lemma \ref{lemclaim2}, for any $t^\prime \in [0,1)$ we have
$\lim_{i \rightarrow \infty} R^{-1}(x_i, t^\prime) \: = \: 1-t^\prime$.
Combining this with the derivative estimate 
$\Big| \frac{\partial R^{-1}}{\partial t} \Big| \: \le \: \eta$ at high
curvature regions gives a contradiction; compare with the proof
of Lemma \ref{claim1}.
Thus the points $x_i$ stay in a compact region. 
We can now use
the bounded-curvature-at-bounded-distance argument in Step 2 of
the proof of Theorem \ref{thmI.12.1} to extract
a convergent subsequence of $\{(\M_i,(x_i,0))\}_{i=1}^\infty$
with a limit Ricci flow solution
$(\M_\infty,(x_\infty,0))$ that exists on
the time interval $[0,1]$.
(In this case, the nonnegative curvature of the blowup
region $W$ comes from the
fact that a standard solution has nonnegative curvature.)
As in Step 3 of the proof of Theorem
\ref{thmI.12.1}, $\M_\infty$ will have
bounded curvature for $t \in [0,1]$.
Note that $\M_\infty$ is a standard solution.  
This contradicts Lemma \ref{lemclaim2}.
\end{proof}

\section{Claim 5 of II.2. Canonical neighborhood property for
standard solutions}
\label{secclaim5}

Let $p$ be the center of the hemispherical region in the time-zero slice.

\begin{lemma} (cf. Claim 5 of II.2) \label{claim5}
Given  $\epsilon > 0$ sufficiently small, there are constants
$\eta = \eta(\epsilon)$, $C_1 = C_1(\epsilon)$ and
$C_2 = C_2(\epsilon)$ so that
every standard solution $\M$ satisfies the conclusions of
Lemmas \ref{derivativeestimates} and \ref{nbhd},
except that the $\epsilon$-neck neighborhood
need not be strong. (Here the constants do not depend on the
standard solution.) More precisely, 
any point $(x,t)$ is covered by one of the following cases : 

1.The time $t$ lies in $(\frac34, 1)$ and 
$(x,t)$ has an $\epsilon$-cap neighborhood
or a strong $\epsilon$-neck neighborhood as in 
Lemma \ref{nbhd}. 

2. $x \in B(p, 0, \epsilon^{-1})$, $t \in [0, \frac34]$ and
$(x,t)$ has an $\epsilon$-cap neighborhood as in Lemma \ref{nbhd}. 

3. $x \notin B(p, 0, \epsilon^{-1})$, $t \in [0, \frac34]$ and
there is an $\epsilon$-neck $B(x,t,\epsilon^{-1} r)$
such that the solution in $P(x,t,\epsilon^{-1}r, -t)$ is, after
scaling with the factor $r^{-2}$, $\epsilon$-close to the appropriate
piece of the evolving round infinite cylinder.

Moreover, we have an estimate $R_{min}(t) \: \ge \: \const \: (1-t)^{-1}$, where
the constant does not depend on the standard solution.
\end{lemma}
\begin{proof}
We first show that the conclusion of Lemma \ref{nbhd} is satisfied.

In view of Lemma \ref{scalarblowup}, there is a $\delta > 0$ so that
if $t \in (1-\delta, 1)$ then
we can apply Theorem \ref{thmI.12.1}
and Lemma \ref{nbhd} to a point $(x,t)$ to see that the conclusions
of Lemma \ref{nbhd} are satisfied in this case.
If $t \in [0, 1-\delta]$ and $x$ is sufficiently far from $p$ 
(i.e. $\dist_0(x,p) \ge D$ for an appropriate $D$) then Lemma \ref{lemclaim2} implies 
that $(x,t)$ has a strong $\epsilon$-neck neighborhood or
there is an $\epsilon$-neck $B(x,t,\epsilon^{-1} r)$
such that the solution in $P(x,t,\epsilon^{-1}r, -t)$ is, after
scaling with the factor $r^{-2}$, $\epsilon$-close to the appropriate
piece of the evolving round infinite cylinder. 

(To elaborate a bit on the last possibility, 
the issue here is that
there is no backward extension of the solution to $t<0$.
Because of this, if $t >0$ is close to $0$ then the backward
neighborhood 
$P(x,t,\epsilon^{-1}r, -t)$ will not exist for rescaled time one,
as required to have a strong $\epsilon$-neck neighborhood. 
Since $\inf_{x \in M} R(x,0) = 1$, we know from
(\ref{Rlowerbound}) that
$R(x,t) \: \ge \: \frac{1}{1- \frac23 t}$.
Then if $t > \frac35$, the time from the initial 
slice to $(x,t)$, after rescaled by the scalar curvature, is bounded below by
$t \: \frac{1}{1- \frac23 t} \: > \: 1$. In particular,
if $t \ge \frac34$ then $r^2 t$ is at least one and
we are ensured that
the backward neighborhood 
$P(x,t,\epsilon^{-1}r, -t)$ does contain a strong
$\epsilon$-neck neighborhood.)

If $t \in [0, 1-\delta]$ and 
$\dist_0(x,p) < D$ then,  provided that
$D$ and $\epsilon$ are chosen appropriately, we can say that
$(x,t)$ has an $\epsilon$-cap neighborhood.

We now show that the conclusion of 
Lemma \ref{derivativeestimates} is satisfied.
If $t \in [1-\delta, 1)$ then the conclusion follows from 
Theorem \ref{thmI.12.1} and Lemma \ref{derivativeestimates}.
If $\delta^\prime > 0$ is sufficiently
small and $t \in [0, \delta^\prime]$ then the conclusion follows from Appendix
\ref{applocalder}. If $t \in [\frac12 \: \delta^\prime, 
1 \: - \: \frac12 \: \delta]$
then we have an upper scalar curvature bound from Lemma \ref{lemclaim4}.
From Hamilton-Ivey pinching (see Appendix \ref{phiappendix}), this implies
a double-sided sectional curvature bound.  
The conclusion of Lemma \ref{derivativeestimates},
when $t \in [\frac12 \: \delta^\prime, 
1 \: - \: \frac12 \: \delta]$, now follows from
the Shi estimates of Appendix \ref{applocalder}. 

The last statement of the lemma follows from the estimate 
$\Big| \frac{\partial R^{-1}}{\partial t} \Big| \: \le \: \const$, which holds
for $t$ near $1$ 
(see Lemmas \ref{derivativeestimates}, 
\ref{lemstandard0}(5) and \ref{scalarblowup})
and then can be extended to all $t \in [0,1)$ (see Lemma \ref{claim2}).
From Lemma \ref{scalarblowup}, $\lim_{t \rightarrow 1} R^{-1}(x,t)=0$ for 
every $x$.  Thus $R^{-1}(x, t )\: \leq \: \const \: (1-t)$ for any $(x,t)$.
Equivalently,
\begin{equation}
R(x,t) \: \geq \: \const \: (1-t)^{-1}.
\end{equation}
\end{proof}

\section{Compactness of the space of standard solutions} 
\label{standardlycompact}

\begin{lemma} \label{standardcompactness}
The family ${\mathcal ST}$ of pointed standard solutions $\{(\M,(p,0))\}$ 
is compact with respect to
pointed smooth convergence. 
\end{lemma}
\begin{proof}
This follows immediately from  Appendix \ref{subsequence}
and the fact that the constant $T_S$ from Lemma \ref{lemclaim2}
is equal to $1$, by Lemma \ref{lemclaim4}.
\end{proof}

\section{Claim 1 of II.2. Rotational symmetry of standard solutions} \label{claimII.2.1}

Consider a standard solution $(M,g(\cdot))$.
Since the time-zero metric $g_0$ is rotationally symmetric, it is clear
by separation of variables that there is a rotationally symmetric
solution for some time interval $[0, T)$.  In this section we
show that every standard solution is rotationally
symmetric for each $t \in [0,1)$.
Of course this would follow from the uniqueness of the 
standard solution; see \cite{Chen-Zhu,Lu-Tian}. But the direct argument given
here is the first step toward a uniqueness proof as in \cite{Lu-Tian}.

\begin{lemma} (cf. Claim 1 of II.2)
Any Ricci flow solution in the space ${\mathcal ST}$ is
rotationally symmetric for all $t \in [0,1)$.
\end{lemma}
\begin{proof}
We first describe an evolution equation for vector fields which
turns out to send Killing vector fields to Killing vector fields.
Suppose that a vector field $u = \sum_m u^m \partial_m$ evolves by
\begin{equation} \label{Killingevolution}
u^m_t \: = \: 
u^{m \: \: k}_{\: \: \: \: ;\: \: \:  k} + R^m_{\: \: \: i} u^i.
\end{equation}
Then
\begin{align} \label{bigmess1}
\partial_t (u^m_{\: \: \: ; i}) \: & = \:  u^m_{t \: ; i} \: + \:
(\partial_t \Gamma^m_{\: \: ki}) \: u^k \\
& = (u^{m \: \: k}_{\: \:
\: \: ;\: \: \:  k} + R^m_{\: \: \: k} u^k)_{; \: i} \: + \:
(\partial_t \Gamma^m_{\: \: ki}) \: u^k \notag \\
& = u^{m \: \: k}_{\: \: \: \: ;\: \: \:  ki} + R^m_{\: \: \: k;i} u^k +
R^m_{\: \: \: k} u^k_{\: \: ;i} \: + \:
(\partial_t \Gamma^m_{\: \: ki}) \: u^k \notag \\
& = u^{m \: \: k}_{\: \: \: \: ;\: \: \:  ik} -
R^m_{\: \: \: lki} \: u^{l \: \: \: k}_{\: \: ;} - R^k_{\: \: \: lki} \: 
u^{m \: \: \: l}_{\: \: \: \: ;}
+ R^m_{\: \: \: k;i} u^k +
R^m_{\: \: \: k} u^k_{\: \: ;i} \: + \:
(\partial_t \Gamma^m_{\: \: ki}) \: u^k \notag \\
& = u^{m \: \: \: \: \: \: k}_{\: \: ;\:ki} -
R^m_{\: \: \: lki} \: u^{l \: \: \: k}_{\: \: ;} - R_{li} \: 
u^{m \: \: \: l}_{\: \: \: \: ;}
+ R^m_{\: \: \: k;i} u^k +
R^m_{\: \: \: k} u^k_{\: \: ;i} \: + \:
(\partial_t \Gamma^m_{\: \: ki}) \: u^k \notag \\
& = u^{m \: \: \: \: \: \: k}_{\: \: ;\:ik} -
(R^m_{\: \: \: lki} u^l)_;^{\: \: k} -
R^m_{\: \: \: lki} \: u^{l \: \: \: k}_{\: \: ;} - R_{ki} \: 
u^{m \: \: \: k}_{\: \: \: \: ;}
+ R^m_{\: \: \: k;i} u^k +
R^m_{\: \: \: k} u^k_{\: \: ;i} \: + \:
(\partial_t \Gamma^m_{\: \: ki}) \: u^k \notag \\
& = u^{m \: \: \: \: \: \: k}_{\: \: ;\:ik} -
R^{m \: \: \: \: \: \: \: k}_{\: \: \: lki;} u^l -
2 R^m_{\: \: \: lki} \: u^{l \: \: \: k}_{\: \: ;} - R_{ik} \: 
u^{m \: \: \: k}_{\: \: \: \: ;}
+ R^m_{\: \: \: k;i} u^k +
R^m_{\: \: \: k} u^k_{\: \: ;i} \: + \:
(\partial_t \Gamma^m_{\: \: ki}) \: u^k. \notag
\end{align}
Contracting the second Bianchi identity gives
\begin{equation} \label{bigmess2}
R_{mlki;}^{\: \: \: \: \: \: \: \: \: \: \: k} = R_{il;m} - R_{im;l}.
\end{equation}
Also,
\begin{align} \label{bigmess3}
\partial_t \Gamma^m_{\: \: ki} & = \partial_t ( g^{ml} \Gamma_{lki}) = 
2 R^{ml} \Gamma_{lki} - g^{ml} (R_{lk,i} + R_{li,k} - R_{ik,l}) \\
& = - R^m_{\: \: \: k;i} - R^m_{\: \: \: i;k} + 
R^{\: \: \: \: \: m}_{ik;}.  \notag
\end{align}
Substituting (\ref{bigmess2}) and (\ref{bigmess3}) in (\ref{bigmess1})
gives
\begin{equation}
\partial_t (u^m_{\: ; i}) \: = \:
u^{m \: \: \: \: \: \: k}_{\: \: ;\:ik} -
2 R^m_{\: \: \: lki} \: u^{l \: \: \: k}_{\: \: ;} - R_{ik} \: 
u^{m \: \: \: k}_{\: \: \: \: ;} +
R^m_{\: \: \: k} u^k_{\: \: ;i}.
\end{equation}
Then
\begin{align}
\partial_t (u_{j ; i}) \: & = \:
\partial_t (g_{jm} u^m_{\: \: \: ; i}) \: = \:
-2 R_{jm} u^m_{\: \: \: ; i} + g_{jm} \: \partial_t (u^m_{\: \: \: ; i}) \\ 
& = \: 
u^{\: \: \: \: \: \: \: \: \: k}_{j  ;\:ik} -
2 R_{jlki} \: u^{l \: \: \: k}_{\: \: ;} - R_{ik} \: 
u^{\: \: \: \: k}_{j  \: ;} -
R_{jk} u^k_{\: \: ;i} \notag \\ 
& = \: 
u^{\: \: \: \: \: \: \: \: \: k}_{j  ;\:ik} +
2 R_{ikjl} \: u^{l \: \: \: k}_{\: \: ;} - R_{ik} \: 
u^{\: \: \: \: k}_{j  \: ;} -
R_{kj} u^k_{\: \: ;i}. \notag
\end{align}
Equivalently, writing $v_{ij} = u_{j ; i}$ gives
\begin{equation}
\partial_t v_{ij} \:  = \:
v^{\: \: \: \: \: \: \: k}_{ij  ;k} +
2 R^{\: \: k \: \: l}_{i \: \: j} \: v_{kl} - R_{i}^{\: \: k}
v_{kj} -
R^k_{\: \: j} v_{ik}. \notag
\end{equation}
Then putting $L_{ij} = v_{ij} + v_{ji}$ gives
\begin{equation}
\partial_t L_{ij} \:  = \:
L_{ij;k}^{\: \: \: \: \: \: \: k} +
2 R^{\: \: k \: \: l}_{i \: \: j} \: L_{kl} - R_{i}^{\: \: k} 
L_{kj} -
R^k_{\: \: j} L_{ik}.  \notag
\end{equation}

For any $\lambda \in \R$, we have
\begin{equation}
\partial_t (e^{2\lambda t} \: L_{ij} L^{ij}) \:  = \:
2 \lambda \: (e^{2\lambda t} \: L_{ij} L^{ij}) \: + \:
(e^{2\lambda t} L_{ij} L^{ij})_{;k}^{\: \: \: \: k} -
2 \: e^{2\lambda t} \: L_{ij;k} \: L^{ij;k} + Q(\Rm, e^{\lambda t}L),
\end{equation}
where $Q(\Rm, L)$ is an algebraic expression that is 
linear in the curvature tensor $\Rm$ and
quadratic in $L$.
Putting $M_{ij} = e^{\lambda t}L_{ij}$ gives
\begin{equation} \label{M2}
\partial_t (M_{ij} M^{ij}) \:  = \:
2 \lambda \: M_{ij} M^{ij} \: + \:
(M_{ij} M^{ij})_{;k}^{\: \: \: \: k} -
2 \: M_{ij;k} \: M^{ij;k} + Q(\Rm, M).
\end{equation}

Suppose that we have a Ricci flow solution $g(t)$, $t \in [0, T]$,
with $g(0) \: = \: g_0$. Let $u(0)$ be a rotational Killing vector field
for $g_0$. Let $u_\infty(0)$ be its restriction to (any) $S^2$, which we will
think of as the $2$-sphere at spatial infinity.
Solve (\ref{Killingevolution}) for $t \in [0, T]$ with $u(t)$ bounded at
spatial infinity for each $t$;
due to the asymptotics coming from Lemma \ref{lemclaim2}
(which is independent of the
rotational symmetry question), there is no problem
in doing so. Arguing as in the proof of Lemma \ref{lemclaim2}, 
one can show that
for any $t \in [0, T]$, at spatial infinity $u(t)$ converges to
$u_\infty(0)$.
Construct $M_{ij}(t)$ from $u(t)$. As $u(0)$ is a Killing vector field,
$M_{ij}(0) = 0$. For any $t \in [0, T]$,
at spatial infinity the tensor
$M_{ij}(t)$ converges smoothly to zero.
Suppose that $\lambda$ is sufficiently negative,
relative to the $L^\infty$-norm of the sectional curvature on the 
time interval $[0, T]$. 
We can apply the maximum principle to (\ref{M2}) to conclude that
$M_{ij}(t) = 0$ for all $t \in [0, T]$. Thus $u(t)$ is a Killing
vector field for all $t \in [0, T]$.

To finish the argument,
as Ben Chow pointed out, any Killing vector field $u$
satisfies 
\begin{equation}
u^{m \: \: k}_{\: \: \: \: ;\: \: \:  k} + R^m_{\: \: \: i} u^i
\: = \: 0.
\end{equation}
To see this, we use the Killing field equation to write
\begin{align}
0 \: & = \: u_{m;k}^{\: \: \: \: \: \: \: \: k} \: + \:
u_{k;m}^{\: \: \: \: \: \: \: \: k} \: = \:
u_{m;k}^{\: \: \: \: \: \: \: \: k} \: + \:
u_{k;m}^{\: \: \: \: \: \: \: \: k} \: - \: u_{k; \: \: m}^{\: \: \: \: k}
\: = \:  u_{m; \: \: \: k}^{\: \: \: \: \: k} \: + \:
u_{\: \: ;mk}^{k} \: - \: u_{\: \: ;km}^{k} \\
& = \:
u_{m; \: \: k}^{\: \: \: \: \: \: k} \: - \: R^k_{\: \: imk} \: u^i \: = \:
u_{m; \: \: k}^{\: \: \: \: \: k} \: + \: R_{mi} \: u^i. \notag
\end{align}

Then from (\ref{Killingevolution}), $u_t^m \: = \: 0$ and the Killing vector
fields are not
changing at all.  This implies that $g(t)$ is rotationally
symmetric for all $t \in [0,T]$.
\end{proof}

\section{Claim 3 of II.2. Uniqueness of the standard solution} \label{claim3}

In this section, which is not needed for the sequel, 
we outline an argument for the  
uniqueness of the standard solution.
We do this for the convenience of the reader.
Papers on the uniqueness issue are
\cite{Chen-Zhu,Lu-Tian}.
Our argument is somewhat different than that of 
\cite[Proof of Claim 3 of Section 2]{Perelman2},
which seems to have some unjustified statements.

In general,
suppose that we have
two Ricci flow solutions $\M \equiv (M, g(\cdot))$ and 
$\widehat{\M} \equiv (M, \widehat{g}(\cdot))$ with bounded curvature
on each compact time interval and the same initial condition.
We want to show that they coincide.
As the set of times for which $g(t) = \widehat{g}(t)$
is closed, it suffices to show that it is relatively open.  Thus
it is enough to show that $g$ and $\widehat{g}$ agree on $[0, T)$ 
for some small $T$.

We will carry out the Deturck trick in this noncompact setting,
using a time-dependent background metric as in 
\cite[Section 2]{Anderson-Chow}. The idea is to
define a $1$-parameter family of metrics $\{h(t)\}_{t \in [0, T)}$ by
$h(t) \: = \: \phi^{-1}(t)^* g(t)$, where
$\{\phi(t)\}_{t \in [0, T)}$ is a $1$-parameter family of
diffeomorphisms of $M$ whose generator is the negative of the time-dependent
vector field
\begin{equation} \label{Weqn}
W^i(t) \: = \: h^{jk} \: \left( \Gamma(h)^i_{\: jk} \: - \: 
\Gamma(\widehat{g})^i_{\: jk} \right),
\end{equation}
with $\phi_0 \: = \: \Id$.
More geometrically, as in \cite[Section 6]{Hamilton}, we consider
the solution of the harmonic heat flow equation
$\frac{\partial F}{\partial t} \: = \: \triangle F$ for maps
$F : M \rightarrow M$ between the
manifolds $(M, g(t))$ and $(M, \widehat{g}(t))$, with $F(0) = \Id$. 

We now specialize to the case when 
$(M, g(\cdot))$ and 
$(M, \widehat{g}(\cdot))$ come from standard solutions.
The technical issue, which we do not address here, is to show
that a solution to the harmonic heat flow
will exist for some time interval $[0, T)$ with uniformly
bounded derivatives; see \cite{Chen-Zhu}.
One is allowed to use the asymptotics of 
Section \ref{claim2} here and
from Section \ref{claimII.2.1},
one can also assume that all of the metrics are
rotationally invariant. In the rest of this section we assume
the existence of such a solution $F$.
By further reducing the time interval if
necessary, we may assume that $F(t)$ is a diffeomorphism of $M$
for each $t \in [0, T)$. Then $h(t) \: = \: F^{-1}(t)^* g(t)$.
Clearly $h(0) \: = \: g(0) \: = \: 
\widehat{g}(0)$.

By Section \ref{claim2}, $g$ and $\widehat{g}$ have the
same spatial asymptotics, namely that of the shrinking cylinder.
We claim that this is also true for $h$. That is, we claim that
$(\M, h(\cdot))$ converges
smoothly to the shrinking cylinder solution on $[0, T)$.
It suffices to show that $F$ converges smoothly to the identity
on $[0, T)$.
Suppose not.  Let $\{x_i\}_{i=1}^\infty$ be a sequence of
points in
the time-zero slice so that no subsequence of
the pointed spacetime maps $(F, (x_i, 0))$ converges to the identity.
Using the derivative
bounds, we can extract a subsequence that converges 
to some $\widetilde{F} \: : \: [0, T) \times \R \times S^2 \rightarrow 
\R \times S^2$ in the pointed smooth topology.
However, $\widetilde{F}$ will satisfy the harmonic heat flow
equation from the shrinking cylinder $\R \times S^2$ to itself, with 
$\widetilde{F}(0)$ being the identity, and will have bounded derivatives.
The uniqueness of $\widetilde{F}$ follows by standard methods. 
Hence $\widetilde{F}(t)$ is the identity for all $t \in [0, T)$, 
which is a contradiction.

By construction, the family of metrics $\{h(t)\}_{t \in [0, T)}$ satisfies
the equation
\begin{equation} \label{heqn}
\frac{dh_{ij}}{dt} \: = \: - \: 2 \: R_{ij}(h) \: + \:
\nabla(h)_i W_j \: + \: \nabla(h)_j W_i.
\end{equation}
In local coordinates, 
the right-hand side of (\ref{heqn}) is a polynomial in 
$h_{ij}$, $h^{ij}$, $h_{ij,k}$ and $h_{ij,kl}$. The leading term in
(\ref{heqn}) is
\begin{equation} \label{hder}
\frac{dh_{ij}}{dt} \: = \: h^{kl} \: \partial_k \partial_l h_{ij}
\: + \: \ldots.
\end{equation}
A particular solution of (\ref{heqn}) is $h(t) \: = \:
\widehat{g}(t)$, since if we had $g \: = \: \widehat{g}$ then
we would have $W = 0$ and $\phi_t = \Id$.

Put $w(t) \: = \: h(t) - \widehat{g}(t)$. We claim that
$w$ satisfies an equation of the form
\begin{equation} \label{wder}
\frac{dw}{dt} \: = \: - \: \nabla(\widehat{g})^* \nabla(\widehat{g}) w
\: + \: Pw \: + \: Qw,
\end{equation}
where $P$ is a first-order operator and $Q$ is a zeroth-order
operator. To obtain the leading derivative terms in (\ref{wder}), 
using (\ref{hder}) we write
\begin{align}
\frac{dw_{ij}}{dt} \: & = \: h^{kl} \: \partial_k \partial_l h_{ij} \: - \:
\widehat{g}^{kl} \: \partial_k \partial_l \widehat{g}_{ij} \: + \: \ldots \\
& = \:
\widehat{g}^{kl} \: \partial_k \partial_l w_{ij} \: + \:
\left( h^{kl} \: - \: \widehat{g}^{kl} \right) \:
\partial_k \partial_l h_{ij} \: + \: \ldots \notag \\
& = \:
\widehat{g}^{kl} \: \partial_k \partial_l w_{ij} \: - \:
\widehat{g}^{ka} \: w_{ab} \: h^{bl} \:
\partial_k \partial_l h_{ij} \: + \: \ldots \notag \\
& = \:
\widehat{g}^{kl} \: \partial_k \partial_l w_{ij} \: - \:
\widehat{g}^{ka} \: h^{bl} \:
\partial_k \partial_l h_{ij} \: w_{ab} \: + \: \ldots \notag
\end{align}
A similar procedure can be carried out for the lower order terms,
leading to (\ref{wder}). By construction the operators $P$ and $Q$
have smooth coefficients which, when expressed in terms of orthonormal frames,
will be bounded on $M$. In fact, as
$h$ and $\widehat{g}$
have the same spatial asymptotics, it follows from
\cite[Proposition 4]{Anderson-Chow} that
the operator on the right-hand side of (\ref{wder}) converges at
spatial infinity to
the Lichnerowicz Laplacian $\triangle_L(\widehat{g})$.

By assumption, $w(0) \: = \: 0$. We now claim that $w(t) \: = \: 0$
for all $t \in [0, T)$.
Let $K \subset M$ be a codimension-zero compact submanifold-with-boundary.
For any $\lambda \in \R$, we have
\begin{align}
& e^{2\lambda t} \: \frac{d}{dt} \: \left( \frac{1}{2} \: e^{-2\lambda t} \: 
\int_K |w(t)|^2 \: \dvol_{\widehat{g}(t)} \right) \:  = \\
& \int_K \left[ \left( -\lambda - \frac{R}{2} \right) |w|^2 \: + \:
\langle w, \frac{dw}{dt} \rangle \right] \: \dvol_{\widehat{g}(t)} = \notag \\
&
\int_K \left[ \left( -\lambda - \frac{R}{2} \right) |w|^2 \: - \:
|\nabla(\widehat{g}) w|^2 \: + 
w_{ab} P^{abcdi} \nabla(\widehat{g})_i w_{cd} \: + \: 
\langle w,  Q w \rangle \right] \: \dvol_{\widehat{g}(t)} 
\pm \notag \\
& \int_{\partial K} \langle w, \nabla_n w \rangle \: 
\dvol_{\partial \widehat{g}(t)} = \notag \\
& \int_K \left[ \left( -\lambda - \frac{R}{2} \right) |w|^2 \: - \:
|\nabla(\widehat{g})_i w^{cd} - \frac{1}{2} P^{abcdi} w_{ab}|^2 \:  
+ \frac{1}{4} |P^{abcdi} w_{ab}|^2 \: + \: 
\langle w,  Q w \rangle \right] \: \dvol_{\widehat{g}(t)} \pm \notag \\
& \int_{\partial K} \langle w, \nabla_n w \rangle 
\: \dvol_{\partial \widehat{g}(t)}. \notag
\end{align}
Choose
\begin{equation}
\lambda \: > \: \sup_{v \neq 0} \frac{
\frac{1}{4} |P^{abcdi} v_{ab}|^2 \: + \: 
\langle v,  Q v \rangle
}{\langle v, v \rangle}.
\end{equation}
On any subinterval $[0, T^\prime] \subset [0, T)$, since
$w$ converges to zero at infinity
and $(M, \widehat{g}(t))$
is standard at infinity, by choosing $K$ appropriately we can make
$\int_{\partial K} \langle w, \nabla_n w \rangle 
\: \dvol_{\partial \widehat{g}(t)}$ small.
It follows that there is an exhaustion
$\{K_i\}_{i=1}^\infty$ of $M$ so that
\begin{equation}
e^{2\lambda t} \: \frac{d}{dt} \: \left( \frac{1}{2} \: e^{-2\lambda t} \: 
\int_{K_i} |w(t)|^2 \: \dvol_{\widehat{g}(t)} \right) \: \le \:
\frac{1}{i}
\end{equation}
for $t \in [0, T^\prime]$.
Then
\begin{equation}
\int_{K_i} |w(t)|^2 \: \dvol_{\widehat{g}(t)} \: \le \: 
\frac{e^{2\lambda t}-1}{\lambda i}
\end{equation}
for all $t \in [0, T^\prime]$.
Taking $i \rightarrow \infty$ gives $w(t) = 0$.

Thus $h \: = \: \widehat{g}$. From (\ref{Weqn}), $W = 0$ and so
$h \: = \: g$. This shows that if $\M, \widehat{\M} \in
{\mathcal ST}$ then $\M = \widehat{\M}$.

\section{II.3. Structure at the first singularity time} \label{II.3}

This section is concerned with the structure of the Ricci flow solution at the
first singular time, in the case when the solution does go singular.

Let $M$ be a connected closed oriented $3$-manifold.  Let $g(\cdot)$ be a Ricci flow on
$M$ defined on a maximal time interval $[0,T)$ with $T < \infty$. One knows
that $\lim_{t \rightarrow T^-} \max_{x \in M} |\Rm|(x,t) \: = \: \infty$.

From Theorem \ref{nolocalcollapse} and Theorem \ref{thmI.12.1}, given
$\epsilon > 0$ there are numbers $r = r(\epsilon) > 0$ and
$\kappa = \kappa(\epsilon) > 0$ so that for any
point $(x,t)$ with $Q = R(x,t)  \ge r^{-2}$, the solution in 
$P(x,t, (\epsilon Q)^{- \: \frac12}, (\epsilon Q)^{- 1})$ is (after rescaling
by the factor $Q$) $\epsilon$-close
to the corresponding subset of a $\kappa$-solution. By Lemma
\ref{derivativeestimates}, the estimate (\ref{derestimates})
holds at $(x,t)$,
provided $\epsilon$ is sufficiently small. In addition, 
there is a neighborhood $B$ of $(x,t)$ 
as described in Lemma  \ref{nbhd}.  In particular, $B$ is a strong $\epsilon$-neck,
an $\epsilon$-cap or a closed manifold with positive sectional curvature.

If $M$ has positive sectional curvature at some time $t$ then it
is diffeomorphic
to a finite quotient of the round $S^3$ and shrinks to a point at time
$T$ \cite{Hamiltonnnnn}. 
The topology of $M$ satisfies the conclusion of the geometrization
conjecture and $M$ goes extinct in a finite time. Therefore
for the remainder of this section we will assume that the sectional curvature
does not become everywhere positive.

We now look at the behavior of the Ricci flow as one approaches the
singular time $T$. 

\begin{definition}
Define a subset $\Omega$ of $M$ by
\begin{equation}
\Omega \: = \: \{x \in M \: : \: \sup_{t \in [0,T)} |\Rm|(x,t) \: < \: \infty \}.
\end{equation}
\end{definition}

Suppose that $x \in M - \Omega$, so there is a sequence of times $\{t_i\}$ in $[0,T)$ with
$\lim_{i \rightarrow \infty} t_i =T$  and
$\lim_{i \rightarrow \infty} |\Rm|(x,t_i) = \infty$. As $\min_M R(\cdot, t)$ in 
nondecreasing in $t$, the largest sectional curvature at $(x, t_i)$ 
goes to infinity as $i \rightarrow \infty$. Then by the $\Phi$-almost nonnegative sectional 
curvature result of Appendix \ref{phiappendix}, $\lim_{i \rightarrow \infty} R(x, t_i) = \infty$.
From the time-derivative estimate of (\ref{derestimates}), $\lim_{t \rightarrow T^-}
R(x,t) = \infty$. Thus $x \in M-\Omega$ if and only if $\lim_{t \rightarrow T^-}
R(x,t) = \infty$.

\begin{lemma} \label{omegaopen}
$\Omega$ is open in $M$. 
\end{lemma}
\begin{proof}
Given $x \in \Omega$, using the time-derivative estimate in
(\ref{eqnderestimatesrinverse})
gives a bound of the form $|R(x,t)| \le C$ for $t \in [0,T)$.
Then using the spatial-derivative estimate in 
(\ref{eqnderestimatesrinverse}) gives
a number $\widehat{r} > 0$ so that
so that $|R(\cdot,t)| \le 2C$ on $B(x,t,\widehat{r})$, for each
$t \in [0, T)$.
The $\Phi$-almost nonnegative sectional 
curvature implies a bound of the form $|\Rm(\cdot,t)| \le C^\prime$ on
$B(x,t,\widehat{r})$, for each
$t \in [0, T)$. 
Then the length-distortion estimate of Section \ref{secI.8.3}
implies that we can pick a neighborhood $N$ of $x$ so that
$|\Rm| \le C^\prime$ on $N \times [0,T)$.
Thus $N \subset \Omega$.
\end{proof}

\begin{lemma} \label{omeganoncompact}
Any connected component $C$ of $\Omega$ is noncompact.
\end{lemma}
\begin{proof}
Since $M$ is connected, if $C$ were
compact then it would  be all of $M$.  This contradicts the assumption that
there is a singularity at time $T$. 
\end{proof}

We remark that {\it a priori}, the structure of $M - \Omega$ can be quite
complicated.  For example, it is not ruled out that an accumulating
collection of
$2$-spheres in $M$ can simultaneously shrink to points. 
That is, $M - \Omega$ could have
a subset of the form $(\{ 0 \} \cup \{\frac{1}{i}\}_{i=1}^\infty )
\times S^2 \subset (-1, 1) \times S^2$, the picture being that
$\Omega$ contains a sequence of smaller and smaller adjacent double horns. 
One could even imagine a
Cantor set's worth of $2$-spheres simultaneously shrinking, although conceivably
there may be additional arguments to rule out both of these cases.

\begin{lemma} \label{standard}
If $\Omega = \emptyset$ then $M$ is diffeomorphic
to $S^3$, $\R P^3$, $S^1 \times S^2$ or $\R P^3 \# \R P^3$.
\end{lemma}
\begin{proof}
The time-derivative estimate in (\ref{eqnderestimatesrinverse}) implies that 
for $t$ slightly less than $T$, we have
$R(x,t) \ge r^{-2}$ for all
$x \in M$. Thus at that time, every $x \in M$ has a neighborhood that is in an
$\epsilon$-neck or an $\epsilon$-cap, as described in Lemma \ref{nbhd}.
(Recall that we have already excluded the positively-curved case of the lemma.) 

As in the proof of Lemma \ref{noncompactsoliton}, by splicing together the projection
maps associated with neck regions, one obtains an open subset  $U\subset M$ 
and a $2$-sphere
fibration $U\ra N$ where the fibers are nearly totally geodesic, and the complement
of $U$  is contained in a union of $\eps$-caps.   It follows that $U$
is connected. If
 there are any $\epsilon$-caps then there must be exactly two of them $U_1,\;U_2$,
and they may be chosen to intersect $U$ in connected  open sets $V_i=U_i\cap U$
which are isotopic to product regions in both $U$ and in the $U_i$'s. 
 The caps being diffeomorphic to $\overline{B^3}$ or
$\R P^3 - B^3$, it follows that $M$ is diffeomorphic to $S^3$, $\R P^3$ or $\R P^3
\# \R P^3$ if $U\neq M$; otherwise  $M$ is diffeomorphic to an $S^2$ bundle over a circle,
and the orientability assumption implies that this bundle is diffeomorphic
to $S^1\times S^2$.
\end{proof}

In the rest of this section we assume that $\Omega \neq \emptyset$.
From the local derivative estimates of Appendix \ref{applocalder}, there is a smooth
Riemannian metric $\overline{g} \: =\: \lim_{t \rightarrow T^-} g(t) \Big|_\Omega$
on $\Omega$. Let $\overline{R}$ denote its scalar curvature. 
Thus 
the scalar curvature function extends to a continuous function on 
the subset $(M\times [0,T))\cup (\Omega\times\{T\})\subset M\times [0,T]$.

\begin{lemma} \label{volumeomega}
$(\Omega, \overline{g})$ has finite volume.
\end{lemma}
\begin{proof}
From the lower scalar curvature bound of (\ref{Rlowerbound}) and the formula
$\frac{d}{dt} \vol(M, g(t)) \: = \: - \: \int_M R \: \dvol_M$,
we obtain an estimate of the form
$\vol(M, g(t)) \: \le \: \const + \const t^{\frac32}$, for
$t < T$. 
The lemma follows.
\end{proof}

\begin{lemma}
\label{lemrinverseextends}
There is a open neighborhood $V$ of $(M-\Omega)\times\{T\}$
in $M\times [0,T]$ such that $R^{-1}$  extends to a 
continuous function on $V$ which vanishes on $(M-\Omega)\times\{T\}$.
\end{lemma}
\begin{proof}
As observed above Lemma \ref{omegaopen}, $x\in M-\Omega$ if and only if
$\lim_{t\ra T^-}\;R^{-1}(x,t)=0$. The lemma follows by applying
(\ref{eqnderestimatesrinverse}) to suitable spacetime paths.
\end{proof}

\begin{definition}
For 
$\rho < \frac{r}{2}$,
put $\Omega_\rho \: = \: \{x \in \Omega \: : \:
\overline{R}(x) \le \rho^{-2} \}$.
\end{definition}

\begin{lemma}
\label{lemoverlinerisproper}
The function $\overline{R}:\Omega\ra \R$ is proper; equivalently,
if $\{x_i\}\subset \Omega$ is a sequence which leaves every compact
subset of $\Omega$, then $\lim_{i \rightarrow \infty} \overline{R}(x_i) = \infty$.
In particular, 
$\Omega_\rho$ is a compact subset of $M$ for every $\rho<r$.
\end{lemma}
\begin{proof}
Suppose $\{x_i\}\subset \Omega$ is a sequence such that
$\{\overline{R}(x_i)\}$ is bounded.  After passing to a subsequence,
we may assume that $\{x_i\}$ converges to some point $x_\infty\in M$.
But $R^{-1}$ is well-defined and continuous on $V$, and vanishes on 
$(M-\Omega)\times\{T\}$, so we must have $x_\infty\in \Omega$.
\end{proof}

We now consider the connected components of $\Omega$ according to whether
they intersect $\Omega_\rho$ or not. First, let $C$ be a connected component of
$\Omega$ that does not intersect $\Omega_\rho$.
Given $x \in C$, there is a neighborhood $B_x$ of $x$ which is $\epsilon$-close 
to a region as described in Lemma \ref{nbhd}.
From Lemma \ref{omeganoncompact}, the neighborhood $B_x$ cannot be of
type (c) or (d) in the terminology of Lemma \ref{nbhd}.

We now introduce some terminology.

If a manifold $Z$ is diffeomorphic to $\R^3$ or $\R P^3-\ol{B^3}$ 
then any embedded $2$-sphere
$\Si\subset Z$ separates $Z$ into two connected subsets, one
of which has compact closure and the other contains the end
of $Z$. We refer to the first component as the {\em compact
side} and the other component as the {\em noncompact side}.

An open subset $R$ of a Riemannian manifold is a {\em good cylinder} if:

\begin{itemize}
\item It is
$\eps$-close, modulo rescaling, to a segment of a round
cylinder of scalar curvature $1$.

\item The diameter of $R$ is approximately
$100$ times its cross-section.

\item Every point in $R$, lies
in an $\eps$-neck in the ambient Riemannian manifold.
\end{itemize}
From Lemma \ref{nbhd},
every $\eps$-cap neighborhood $B_x$
contains a good cylinder lying in the $\eps$-neck at the end
of $B_x$.

\begin{lemma} \label{doublehorn}
Suppose that for all $x \in C$, 
the neighborhood $B_x$ can be taken to be a strong
$\epsilon$-neck as in case (a) of Lemma \ref{nbhd}. 
Then $C$ is a double $\epsilon$-horn.
\end{lemma}
\begin{proof}
Each point $x$ has an $\epsilon$-neck neighborhood.
We can glue these $\epsilon$-necks together to form a
submersion from $C$ to a $1$-manifold, with fiber $S^2$;
cf. the proof of Lemma \ref{noncompactsoliton}. (We can 
do the gluing by successively adding on good
cylinders, where the intersections of successive cylinders 
have diameter approximately 10 times the diameter of the cross-sections.) 
In view of Lemma
\ref{lemrinverseextends}, it follows in this case that
$C$ is a double $\epsilon$-horn.
\end{proof}

\bigskip
\begin{lemma}
Suppose that there is some $x \in C$ whose neighborhood
$B_x$ is an $\epsilon$-cap as in case (a) of Lemma \ref{nbhd}.
Then $C$ is a capped $\epsilon$-horn.
\end{lemma}
\begin{proof}

Put $p_1=x$.  Let $R$ be a good cylinder in the $\eps$-neck
at the end of $B_{p_1}$.  Now glue on successive good cylinders
to $R$,
as in the proof of the preceding lemma,
going away from $p_1$. \\ \\
Case 1 :
Suppose this gluing process can be continued
indefinitely.  Then taking the union of $B_{p_1}$ with 
the good cylinders, we obtain an open subset $W$ of $C$ which
is diffeomorphic to $\R^3$ or $\R P^3-\ol{B^3}$. 
We claim that $W$ is a closed subset of $\Om$.
If not then there is a sequence $\{x_k\}_{k=1}^\infty
\subset W$ converging to some 
$x_\infty\in \Om-W$. This implies that $\{\ol{R}(x_k)\}_{k=1}^\infty$
remains bounded. In view of the overlap condition between successive
good cylinders, a subsequence of
$\{x_k\}_{k=1}^\infty$ lies in an infinite
number of mutually disjoint good cylinders, whose volumes have a positive lower
bound (because of the upper scalar curvature bound at the points $x_k$).
This contradicts Lemma \ref{volumeomega}.

Thus $W$ is open and closed in $\Om$. Hence $W=C$ and 
we are done. \\ \\
Case 2 :
Now suppose that the gluing process cannot be continued beyond some 
good cylinder $R_1$.  Then there must be a point $p_2\in R_1$
such that $B_{p_2}$ is an $\eps$-cap. Also, note that 
the union $W_1$ of $B_{p_1}$ with the good cylinders is diffeomorphic
to $\R^3$ or $\R P^3-\overline{B^3}$, and that $R_1$ has compact
complement in $W_1$.   Let $\Si\subset R_1$ be a cross-sectional $2$-sphere
passing through $p_2$.  

We first claim that if $V$ is the compact side of $\Si$ in $W_1$,
then $V$ coincides with the 
compact side $V'$ of $\Si$ in $B_{p_2}$.  To see this, note that $V$ and $V'$
are both connected open sets disjoint from $\Si$, with topological 
frontiers $\partial V=\partial V'=\Si$.  Then
$V - V^\prime \: = \: V \cap (C - (V^\prime \cup \Sigma))$
and we obtain two
open decompositions
\begin{equation}
\label{eqndecompv}
V=(V\cap V')\sqcup (V-V'),\quad V'=(V\cap V')\sqcup (V'-V).
\end{equation}
If $V\cap V'=\emptyset$, then $\ol{V}\cup\ol{V'}$ is a union of 
two compact manifolds with the same boundary $\Si$, and disjoint interiors.
Hence it is an open
and closed subset of the connected component $C$, which contradicts
Lemma \ref{omeganoncompact}.  Thus $V\cap V'$ is nonempty.
By (\ref{eqndecompv}) and the connectedness of $V$ and $V^\prime$,
we get $V\subset V'$ and $V'\subset V$,
so $V=V'$ as claimed.

Next, we claim that if $R_2\subset B_{p_2}$ is a good cylinder 
with compact complement in $B_{p_2}$, then $R_2$ is disjoint from $W_1$.  
To see this,  note that  $R_2$ is disjoint from $\Si$ because 
$p_2 \in \Si$ and the diameter
of $\Si$ is close to $\pi (\ol{R}(p_2)/6)^{-1/2}$, whereas by Lemma \ref{nbhd}
there is an $\eps$-neck $U\subset B_{p_2}$ with compact complement in 
$B_{p_2}$, at distance at least $9000 \ol{R}(p_2)^{-1/2}$ from $p_2$.
Thus $R_2$ must lie in the noncompact side of $\Si$ in $B_{p_2}$, and 
hence is disjoint from $V$. 
As the good cylinder
$R_1 \ni p_2$ lies within $B(p_2, 1000 \ol{R}(p_2)^{-1/2}) \subset \Omega$,
it follows that $R_2$ is also disjoint from $R_1$,
so $R_2$ is disjoint from $W_1=V\cup R_1\subset B_{p_2}$.

We continue adding good cylinders to $R_2$ as long as we can.
If we come
to another cap point $p_3$ then we jump to its cap
$B_{p_3}$ and continue the process.  
When so doing, we encounter 
successive cap points $p_1, p_2, \ldots$ with associated caps
$B_{p_1}\subset B_{p_2} \subset \ldots$ and disjoint good cylinders
$R_1, R_2, \ldots$.
Since the ratio 
$\frac{\sup_{B_{p_k}}\,\ol{R}}{\inf_{B_{p_k}}\,\ol{R}}$ has an 
{\it a priori}
bound by Lemma \ref{nbhd}, 
in view of the disjoint good cylinders in $B_{p_k}$
we get $\vol(B_{p_k})\geq \const  k\,\ol{R}(p_1)^{-3/2}$.
Then Lemma \ref{volumeomega} gives 
an upper bound on $k$.  Hence we encounter a finite
number of cap points. 
Arguing as in Case 1,
we conclude that $C$ is a capped $\eps$-horn.
\end{proof}

We note that there could be an infinite number of connected components of
$\Omega$ that do not intersect $\Omega_\rho$. 

Now suppose that $C$ is a connected component of $\Omega$ that intersects
$\Omega_\rho$. As $C$ is noncompact, there must be some point $x \in C$ that is
not in $\Omega_\rho$. Again, any such $x$ has a neighborhood $B$ as in Lemma
\ref{nbhd}. 
If one of the boundary components of $B$ intersects
$\Omega_{2\rho}$ then we terminate the process in that direction.
For the directions of the boundary components of $B$ that do not 
intersect $\Omega_{2\rho}$,
we perform the above algorithm of looking for an adjacent
$\epsilon$-neck, etc.
The only difference from before 
is that in at
least one direction any such
sequence of overlapping $\epsilon$-necks will be finite, as it must
eventually intersect 
$\Omega_{2\rho}$. 
(In the other direction it may terminate 
in $\Omega_{2\rho}$, in an $\epsilon$-cap, or not terminate at all.)
Once a cross-sectional $2$-sphere intersects $\Omega_{2\rho}$, if
$\epsilon$ is small then the entire $2$-sphere lies in $\Omega_\rho$.
Thus any connected component of $C - (C \cap \Omega_\rho)$
is contained in an $\epsilon$-tube with both boundary components in $\Omega_\rho$, 
an $\epsilon$-cap with boundary in $\Omega_\rho$ or an $\epsilon$-horn with boundary in 
$\Omega_\rho$. We note that $\Omega_\rho$ need not have a nice boundary.

There is an {\it a priori} $\rho$-dependent lower bound for the volume of any such
connected component of $C - (C \cap \Omega_\rho)$, in view of the fact that
it contains $\epsilon$-necks that adjoin $\Omega_\rho$.
From Lemma \ref{volumeomega},
there is a finite number of connected components of $\Omega$ that
intersect $\Omega_\rho$. Any such connected component has a finite
number of ends, each being an $\epsilon$-horn. Note that the
$\epsilon$-horns can be made disjoint, each with a quantitative lower
volume bound.

The surgery procedure, which will be described in detail in Section 
\ref{II.4.4}, is performed as follows.  First, one throws away all
connected components of $\Omega$ that do not intersect $\Omega_\rho$.
For each connected component $\Omega_j$ of $\Omega$ that intersects
$\Omega_\rho$ and for each $\epsilon$-horn of $\Omega_j$, take
a cross-sectional sphere that lies far in the $\epsilon$-horn.
Let $X$ be what's left after cutting the $\epsilon$-horns 
at these $2$-spheres and removing the tips. The 
(possibly-disconnected) postsurgery manifold $M^\prime$ is the 
result of capping off $\partial X$ by $3$-balls. 

We now discuss how to reconstruct the original manifold
$M$ from $M^\prime$.

\begin{lemma} \label{reconstruct}
$M$ is the result of taking connected sums of components of $M^\prime$ 
and possibly taking additional connected sums with a finite number of
$S^1 \times S^2$'s and $\R P^3$'s.
\end{lemma}
\begin{proof}

At a time shortly before $T$,
each point of $M - X$ has a neighborhood as in Lemma
\ref{nbhd}. The components of $M-X$ are $\epsilon$-tubes and
$\epsilon$-caps.  Writing $M^\prime = X \cup \bigcup B^3$ and
$M = X \cup (M-X)$, one builds $M$ from
$M^\prime$ as follows. 
If the boundary of an $\epsilon$-tube of $M-X$ lies in 
two disjoint components of $X$ then it gives rise to
a connected sum of two components of $M^\prime$.  If the
boundary of an $\epsilon$-tube lies in a single connected 
component of $X$ then it gives rise to the 
connected sum of the corresponding component of $M^\prime$ 
with a new copy of $S^1 \times S^2$.
If an $\epsilon$-cap 
in $M-X$
is a $3$-ball it does not have any effect
on $M^\prime$. If
an $\epsilon$-cap is $\R P^3 - B^3$ then it gives rise to the
connected sum of the corresponding component of $M^\prime$ 
with a new copy of $\R P^3$. The lemma follows.
\end{proof}

\begin{remark}
We do not assume that
the diameter of $(M, g(t))$ stays bounded as $t \rightarrow T$;
it is an open question whether this is the case.
\end{remark}

\section{Ricci flow with surgery:
the general setting}
\label{Ricci flow with surgery}

In this section we introduce some notation and terminology in order to
treat Ricci flows with surgery.

The principal purpose of sections II.4 and II.5 is to show that
one can prescribe the surgery procedure in such a way that Ricci flow
with surgery is well-defined for all time.  This involves showing
that 

$\bullet$ One can give a sufficiently precise description of 
the formation of singularities so that one can envisage defining a geometric
surgery.
In the case of the formation of the first singularity, 
such a description was given in Section \ref{II.3}. 

$\bullet$ The sequence of surgery times cannot accumulate.

\noindent
The argument in Section \ref{II.3} strongly uses both the $\kappa$-noncollapsing
result
of Theorem \ref{nolocalcollapse} and the characterization
of the geometry in a spacetime region around a point $(x_0, t_0)$ with
large scalar curvature, as given in Theorem \ref{thmI.12.1}.  
The proofs of both of these results use the 
smoothness of the solution at times before $t_0$. If surgeries occur
before $t_0$ then one must have strong control on the scales at which the
surgeries occur, in order to extend the arguments of Theorems
\ref{nolocalcollapse} and \ref{thmI.12.1}.
This forces one to consider time-dependent scales.
    
Section II.4 introduces Ricci flow with surgery, in varying degrees of 
generality. Our treatment of this material follows Perelman's.
We have added some terminology to help formalize the surgery process.
There is some arbitrariness in this formalization, but the version given
below seems adequate.
 
For later  use, we now summarize the relevant notation
that we introduce.  More precise  definitions will be given below.
We will avoid using new notation as much as possible. 

$\bullet$ $\M$ is a Ricci flow with surgery.

$\bullet$ $\M_t$ is the time-$t$ slice of $\M$.

$\bullet$ $\M_{\reg}$ is the set of regular points of $\M$.

$\bullet$ If $T$ is a singular time then $M_T^-$ is the limit of time slices
$\M_t$ as $t \rightarrow T^-$ (called $\Omega$ in II.4.1) and
$M_T^+$ is the outgoing time slice (for example,  the result of performing
surgery on $\Omega$). If $T$ is a nonsingular time then
$\M_T^- = \M_T^+ = \M_T$. 

The basic notion of a Ricci flow with surgery is simply a sequence of
Ricci flows which ``fit together'' in the sense that the final 
(possibly singular) time slice 
of each flow is isometric, modulo surgery, to the initial time slice
of the next one. 
\begin{definition}   \label{flowwithsurgery}
A {\em Ricci flow with surgery} is given by 

$\bullet$   A collection of Ricci flows
$\{(M_k\times[t_k^-,t_k^+),g_k(\cdot))\}_{1\leq k \leq N}$,
where $N\leq \infty$, $M_k$ is a compact (possibly empty) manifold,
$t_k^+=t_{k+1}^-$ for  all
$1\leq k<N$, and the flow $g_k$ 
goes singular at $t_k^+$ 
for each $k < N$. We allow $t_N^+$ to be $\infty$.

$\bullet$ A collection of limits
$\{(\Om_k,\bar g_k)\}_{1\leq k\leq N}$, in the sense of Section \ref{II.3},
at the respective final times $t_k^+$ that are
singular if $k < N$. (Recall that
$\Omega_k$ is an open subset of
$M_k$.)

$\bullet$ A collection of isometric embeddings 
$\{\psi_k:X_k^+\ra X_{k+1}^-\}_{1\leq k<N}$
where $X_k^+\subset \Om_k$ and $X_{k+1}^-\subset M_{k+1}$,
$1 \le k< N$, are compact $3$-dimensional submanifolds with 
boundary.   The $X_k^\pm$'s are the subsets which survive the 
transition from one flow to the next, and the $\psi_k$'s give
the identifications between them.   

We will say that $t$ is a {\em singular time} if $t=t_k^+=t_{k+1}^-$
for some $1\leq k<N$, or $t=t_N^+$ and the metric goes singular
at time $t_N^+$.
\end{definition}

A Ricci flow with surgery does not necessarily have to have any
real surgeries, i.e. it could be a smooth nonsingular flow. 
Our definition allows Ricci flows with surgery that are more
general than those appearing in the argument for geometrization,
where the transitions/surgeries have a very special form.  Before
turning to these more special flows in Section \ref{II.4.4}, 
we first discuss some 
basic features of Ricci flow with surgery.

It will be convenient to associate a (non-manifold) spacetime $\M$ to
the Ricci flow with surgery.  This is constructed by taking the disjoint union 
of  the smooth manifolds-with-boundary  
\begin{equation}
\left(M_k\times [t_k^-,t_k^+)\right)\cup \left(\Om_k\times\{t_k^+\}\right)
\subset M_k\times [t_k^-,t_k^+]
\end{equation}
 for $1\leq k\leq N$
and making identifications using the $\psi_k$'s as gluing
maps.   We denote the quotient space by $\M$ and the quotient map
by $\pi$.  We will sometimes also use $\M$ to refer to the
whole Ricci flow with surgery structure,
rather than just the associated spacetime.
The {\em time-$t$ slice $\M_t$ of $\M$} is the image 
of the time-$t$ slices of the constituent Ricci
flows under the quotient map.  
 
If $t=t_k^+$ is a singular time then we put 
$\M_t^-=\pi(\Omega_k\times\{t_k^+\})$;  if in addition $t\neq t_N^+$ then 
we put $\M_t^+=\pi(M_{k+1}\times\{t_{k+1}^-\})$. If $t$ is not a singular time
then we put  $\M_t^+=\M_t^-=\M_t$. We refer to 
$\M_t^+$ and $\M_t^-$ as the forward and backward time
slices, respectively.

Let us summarize the structure of $\M$ near a singular time $t = t_k^+ = t_{k+1}^-$.
The backward time slice $\M_t^-$ is a copy of $\Omega_k$. The forward
time slice $\M_t^+$ is a copy of $M_{k+1}$. The time slice
$\M_t$ is the result of gluing $\Omega_k$ and $M_{k+1}$ using $\psi_k$.
Thus it is the disjoint union of $\Omega_k - X_k^+$, $M_{k+1} - X_{k+1}^-$ and
$X_k^+ \cong X_{k+1}^-$.  If $s > 0$ is small then in going from $M_{t-s}$ to
$M_{t+s}$, the topological change is that we remove
$M_k - X_k^+$ from $M_k$ and add $M_{k+1} - X_{k+1}^-$.

We let $\M_{(t,t')} = \cup_{\bar t\in(t,t')} \M_{\bar t}$ denote
the  {\em time slab} between $t$ and $t'$, i.e. 
the union of the time slices between $t$ and $t'$.
The closed time
slab $\M_{[t,t']}$ is defined to be the closure of $\M_{(t,t')}$
in $\M$, so $\M_{[t,t']}=\M_t^+\cup\M_{(t,t')}\cup \M_{t'}^-$.
We  (ab)use  the notation
$(x,t)$ to denote a point $x\in\M$ lying in the time $t$ slice $\M_t$, 
even though
$\M$ may no longer be a product.

The spacetime $\M$
has three types of points:  \\
1. The  $4$-manifold
points, which include all points at 
nonsingular times in
$(t^-_1, t^+_N)$ and all
points in $\pi(\interior(X_k^+)\times \{t_k^+\})$ 
(or $\pi(\interior(X_{k}^-)\times \{t_{k}^-\})$) for $1 \le k < N$,
\\
2. The boundary points of $\M$, which are the images in $\M$ of
$M_1\times \{t_1^-\}$,
$\Omega_N \times \{t^+_N\}$,
$(\Om_k - X_k^+)\times \{t_k^+\}$ for 
$1\leq k<N$, and $(M_k - X_k^-)\times \{t_k^-\}$ for $1<k\leq N$, and \\
3. The ``splitting'' points, which
 are the images in $\M$ of $\D X_k^+\times\{t_k^+\}$
for $1\leq k<N$.  

Here the classification of points is according 
to the smooth structure, not  the topology. In fact, $\M$ is  a topological
manifold-with-boundary.
We say that $(x,t)$ is {\em regular} if it is either a $4$-manifold
point, or it lies in the initial time slice $\M_{t_1^-}$ or final time slice
$\M_{t_N^+}$.
Let  $\M_{\reg}$ denote the set of regular points. It has a natural
smooth structure since the gluing maps $\psi_k$, being isometries
between smooth Riemannian manifolds, are smooth maps.

Note that the Ricci
flows on the $M_k$'s define a Riemannian metric $g$ on the ``horizontal''
subbundle of the tangent bundle of $\M_{\reg}$.  It follows from the 
definition of the Ricci flow that $g$ is actually smooth on $\M_{\reg}$.

We metrize each time slice $\M_t$, and the forward and backward time
slices $\M_t^\pm$, by infimizing the path length of piecewise smooth paths.
We allow our distance functions to be infinite,
since the infimum will be infinite when points lie in different components.
If $(x,t)\in \M_t$ and $r>0$ then we let $B(x,t,r)$ denote the 
corresponding metric ball. Similarly,  $B^\pm(x,t,r)$ denotes
the ball in $\M_t^\pm$ centered at $(x,t)\in\M_t^\pm$.  
A ball $B(x,t,r)\subset\M_t$ is {\em proper} if the distance
function $d_{(x,t)}:B(x,t,r)\ra [0,r)$ is a proper function;
a proper ball ``avoids singularities'', except possibly at its
frontier.
Proper balls $B^\pm(x,t,r)\subset\M_t^\pm$ are defined likewise.

An {\em admissible curve} in $\M$ is a path $\ga:[c,d]\ra \M$,
with $\ga(t)\in \M_t$ for all $t\in [c,d]$, such that 
for each $k$, the part of $\ga$ landing in $\M_{[t_k^-t_k^+]}$
lifts to a smooth map into 
$M_k\times[t_k^-,t_k^+)\;\cup \;\Om_k\times\{t_k^+\}$.  We will
use $\dot{\ga}$ to denote the ``horizontal'' part of the velocity
of an admissible curve $\ga$.    
If $t < t_0$,
a point $(x,t)\in\M$ is 
{\em accessible from $(x_0,t_0)\in\M$} if there is an admissible curve
running from $(x,t)$ to $(x_0,t_0)$.  An 
admissible curve $\ga:[c,d]\ra\M$ is {\em static} if
its lifts to the product spaces have constant first component. That is,
the points in the image of a  static curve are ``the same'',
modulo the passage of time and identifications taking place
at surgery times.
 A {\em barely
admissible curve} is an admissible curve $\ga:[c,d]\ra\M$
such that the  image is not contained
in $\M_c^+ \cup \M_{\reg}\cup \M_d^-$.
If $\ga:[c,d]\ra\M$ is barely admissible then 
there is a surgery time $t=t_k^+=t_{k+1}^-\in(c,d)$ such that
$\ga(t)$ lies in 
\begin{equation}
\pi(\D X_k^+\times\{t_k^+\})
=\pi(\D X_{k+1}^-\times\{t_{k+1}^-\}).
\end{equation}

If $(x,t)\in\M_t^+$, $r>0$, and $\De t>0$ then we define the 
forward {\em parabolic region } 
$P(x,t,r,\De t)$ to be the union of 
(the images of) the static
admissible curves $\ga:[t,t^\prime] \ra \M$ starting in $B^+(x,t,r)$, where
$t^\prime \le t+\De t$. That is, we take the union of all the maximal
extensions of all static
curves, up to time $t + \De t$, starting from the
initial time slice $B^+(x,t,r)$.
When $\De t<0$, the parabolic region $P(x,t,r,\De t)$ is defined 
similarly using static admissible curves ending in $B^-(x,t,r)$.

If $Y\subset \M_t$, and $t\in[c,d]$ then we say that $Y$ 
is {\em unscathed in $[c,d]$} if every point $(x,t)\in Y$
lies on  a static curve defined on the time
interval $[c,d]$.  If, for instance, $d=t$ then this will force
$Y\subset \M_t^-$.  
The term ``unscathed'' is intended to capture the idea that 
the set is unaffected by singularities and surgery.
(Sometimes Perelman uses the phrase 
``the solution is defined in $P(x,t,r,\De t)$'' as synonymous with
``the solution is unscathed in $P(x,t,r,\De t)$'', for example in
the definition of canonical neighborhood in II.4.1.)
We  may use the 
notation  $Y\times [c,d]$ for the set of points lying on static
curves $\ga:[c,d]\ra\M$ which pass through $Y$, when
$Y$ is unscathed on $[c,d]$.
Note that if $Y$ is open and
unscathed on $[c,d]$ then we can think of the Ricci flow on $Y\times [c,d]$
as an ordinary (i.e. surgery-free) Ricci flow.

The definitions of $\eps$-neck, $\eps$-cap, $\eps$-tube and
(capped/double) $\eps$-horn 
from Section \ref{notationterminology}
do not require modification for 
a Ricci flow with surgery, since they are just special types of Riemannian 
manifolds; they will turn up as subsets of forward or backward 
time slices of a Ricci flow with surgery. A {\em strong $\eps$-neck} 
is a subset  of the form $U\times [c,d]\subset \M$, where $U\subset \M_d^-$
is an open set  that is unscathed on the interval $[c,d]$, 
which is a strong $\eps$-neck in the sense of 
Section \ref{notationterminology}.

\section{II.4.1. {\it A priori} assumptions}

This section introduces the notion of canonical neighborhood.

The following definition captures the geometric structure that
emerges by combining Theorem \ref{thmI.12.1} and its extension
to Ricci flows with surgery (Section \ref{II.5})
with the geometric description of $\kappa$-solutions.   The idea
is that blowups either yield $\kappa$-solutions, whose
structure is well understood from Section \ref{II.1}, or there are surgeries nearby
in the recent  past, in which case the local geometry resembles
that of the standard solution.   Both alternatives produce canonical
neighborhoods.

\begin{definition}(Canonical neighborhoods, cf. Definition in II.4.1)
\label{canonicalnbhddef}
Let $\eps>0$ be 
small enough so that Lemmas \ref{nbhd} and \ref{claim5} hold.
Let $C_1$ be the maximum of $30 \epsilon^{-1}$ and the $C_1(\eps)$'s of Lemmas
\ref{nbhd} and \ref{claim5}. Let $C_2$ be the maximum of the
$C_2(\epsilon)$'s of Lemmas
\ref{nbhd} and \ref{claim5}.
Let  $r:[a,b]\ra (0,\infty)$ be a positive nonincreasing function.
A Ricci flow with surgery $\M$ defined on the time interval $[a,b]$ satisfies
the {\em $r$-canonical neighborhood assumption} if every
$(x,t)\in \M_t^\pm$ with scalar curvature $R(x,t)\geq r(t)^{-2}$
has a canonical neighborhood in the corresponding (forward/backward)
time slice, as in Lemma \ref{nbhd}.  More precisely, there is
an  $\hat r\in (R(x,t)^{-\frac12},C_1R(x,t)^{-\frac12})$ 
and an open 
set $U\subset \M_t^\pm$ with 
$\ol{B^\pm(x,t,\hat r)}\subset U\subset B^\pm(x,t,2\hat r)$ that
falls into one of the
following categories : 

(a) $U\times [t-\De t,t]\subset\M$  is a strong $\eps$-neck
for some $\De t>0$.  (Note that after parabolic rescaling the scalar
curvature at $(x,t)$ becomes $1$, so the scale factor must be
$\approx R(x,t)$, which implies that $\De t\approx R(x,t)^{-1}$.)

(b) $U$ is an $\eps$-cap which, after rescaling, is $\epsilon$-close to the
corresponding piece of a $\kappa_0$-solution or a time slice of a standard
solution (cf. Section \ref{II.2}).

(c) $U$ is a closed manifold diffeomorphic to $S^3$ or $RP^3$.

(d) $U$ is  $\epsilon$-close to 
a closed manifold of constant positive sectional 
curvature.

\noindent
Moreover, the scalar curvature in $U$ lies between $C_2^{-1}R(x,t)$
and $C_2R(x,t)$.
In cases (a), (b), and  (c), the volume of $U$ is greater than 
$C_2^{-1}R(x,t)^{-\frac32}$. In case (c), the infimal
sectional curvature of $U$ is greater than $C_2^{-1}R(x,t)$.

Finally, we require that
\begin{equation}
\label{II.(1.3)estimates}
|\nabla R(x, t)|<\eta R(x, t)^{\frac32}, \: \: \:  \: \: \: \: 
\left|\frac{\D R}{\D t} (x, t)\right|<\eta R(x, t)^2,
\end{equation}
where $\eta$ is the constant from (\ref{derestimates}).  Here the time 
dervative $\frac{\D R}{\D t} (x, t)$ should be interpreted
as a one-sided derivative when the point $(x,t)$ is added or removed
during surgery at time $t$.
\end{definition}

\begin{remark}
Note that the smaller of the 
two balls in $\ol{B^\pm(x,t,\hat r)}\subset U\subset B^\pm(x,t,2\hat r)$ is closed, 
in order to 
make it easier to check the openness of the canonical neighborhood condition.
The requirement that $C_1$ be at least $30 \epsilon^{-1}$ will be used
in the proof of Lemma \ref{extendit}; see Remark \ref{30}.
\end{remark}

\begin{remark}
For convenience, in case (b) we have added the
extra condition that $U$ is $\epsilon$-close to the
corresponding piece of a $\kappa_0$-solution or 
a time slice of a standard solution.
One does not need this
extra condition, but it is consistent to add it.
We remark that when surgery is performed according to the recipe
of Section \ref{II.4.4}, if a point $(p,t)$ lies in $\M_t^+   -  \M_t^-$ (i.e.
it  is ``added'' by surgery) then
it will sit in an $\eps$-cap, because $\M_t^+$ will resemble
a standard solution from Section \ref{II.2} near $(p,t)$.   Points lying
somewhat further out on the capped neck will belong to a strong 
$\eps$-neck which extends backward in time prior to the surgery. 
 \end{remark} 

The next condition, which will ultimately be guaranteed by the
Hamilton-Ivey curvature pinching result and careful surgery, 
is also essential in blowup arguments \`a la Section
\ref{I.12.1}.
\begin{definition}($\Phi$-pinching)
\label{phipinchingdef}
Let  $\Phi \in C^\infty(\R)$ be a positive nondecreasing 
function such that for positive $s$,
$\frac{\Phi(s)}{s}$ is a decreasing function which tends to zero
as $s \rightarrow \infty$.  The Ricci flow with surgery $\M$
satisfies the {\em $\Phi$-pinching assumption} if for all
$(x,t) \in \M$, one has $\Rm(x,t) \geq -\Phi(R(x,t))$.
\end{definition}

We remark that the notion of $\Phi$-pinching here is somewhat
different from Perelman's $\phi$-pinching.  The purpose of 
this definition is to distill out the properties of the
Hamilton-Ivey pinching condition which are needed in the rest of the
proof.  

\begin{definition} \label{apriori}
A Ricci flow with surgery satisfies the {\em {\it a priori} assumptions}
if it satisfies the $\Phi$-pinching and $r$-canonical
neighborhood assumptions
on the time interval of the flow. Note that
the {\it a priori} assumptions depend on $\epsilon$, the function $r(t)$ of
Definition \ref{canonicalnbhddef} 
and the function $\Phi$ of Definition \ref{phipinchingdef}.
\end{definition}

\section{II.4.2. Curvature bounds from the {\it a priori} assumptions}

In this section we state some technical lemmas about Ricci flows with surgery
that satisfy the {\em {\it a priori}} assumptions of the previous section.

The first one is the surgery 
analog of  Lemma \ref{claim1}.

\begin{lemma}(cf. Claim 1 of II.4.2)
\label{claim1II.4.2}
Given $(x_0,t_0)\in\M$ put $Q = |R(x_0,t_0)|+r(t_0)^{-2}$.
Then $R(x,t)\leq 8Q$ for all 
$(x,t) \in P(x_0,t_0,\frac12\eta^{-1}Q^{-\frac12},
-\frac18\eta^{-1}Q^{-1})$, where $\eta$ is the constant from
(\ref{II.(1.3)estimates}).
\end{lemma}
\begin{proof}
The lemma follows from the estimates (\ref{II.(1.3)estimates}).
One integrates these derivative bounds along a subinterval of
a path that goes 
in $B(x_0,t_0,\frac12\eta^{-1}Q^{-\frac12})$ and then
backward in time along a static path. See the proof of
Lemma \ref{claim1}.
We also use the fact that if $t^\prime \le t_0$ and 
$R(x^\prime, t^\prime) \ge Q$ then the inequalities
(\ref{II.(1.3)estimates}) are valid at $(x^\prime, t^\prime)$,
since $r(\cdot)$ is nonincreasing.
\end{proof}

The next lemma expresses the main consequence of
Claim 2 of II.4.2. 

\begin{lemma}(cf. Claim 2 of II.4.2)
\label{claim2II.4.2}
If $\epsilon$ is small enough then the following holds.
Suppose that $\M$ is a Ricci flow with surgery that satisfies the
$\Phi$-pinching assumption.
Then for any $A<\infty$ and $\widehat{r} > 0$ there exist 
$\xi=\xi(A) > 0$ and $K=K(A,\widehat{r}) < \infty$
with the following property. 
Suppose that $\M$ also satisfies the $r$-canonical neighborhood assumption
for some function $r(\cdot)$.
Then for any time $t_0$, if
$(x_0,t_0)$ is a point so that $Q = R(x_0,t_0) > 0$ satisfies
$\frac{\Phi(Q)}{Q} < \xi$
and $(x, t_0)$ is a point so that
$\dist_{t_0}(x_0, x) \le A Q^{- \: \frac12}$ then $R(x, t_0) \le KQ$,
where $K = K(A, r(t_0))$.
\end{lemma}
\begin{proof}
The proof is similar to the  proof of Step 2 of Theorem
\ref{thmI.12.1}. 
(The canonical neighborhood assumption
replaces Step 1 of the proof of Theorem \ref{thmI.12.1}.)
Assuming that the lemma fails, one obtains a
piece of a nonflat metric cone as a blowup limit.
Using the canonical neighborhood assumption, one concludes
that the corresponding points in $\M$ 
have a neighborhood of type (a), i.e. a strong
$\epsilon$-neck, since
the neighborhoods of type (b), (c) and (d) of Definition 
\ref{canonicalnbhddef} are not close to a piece of metric
cone. A strong $\epsilon$-neck, has the time interval needed to apply
the strong maximum
principle as in Step 2 of the proof of Theorem \ref{thmI.12.1}, in
order to get a contradiction.
\end{proof}

\section{II.4.3. $\delta$-necks in $\epsilon$-horns}

In this section we show that an $\epsilon$-horn has a
self-improving property as one goes down the horn.  
For any $\delta > 0$, if the scalar curvature at a point
is sufficiently large then the point actually lies in a
$\delta$-neck.

In the statement of the 
next lemma we will write $\Omega$
synonymously with the $\M_T^-$ of Section
\ref{Ricci flow with surgery}.

\begin{lemma}(cf. Lemma II.4.3) \label{lemmaII.4.3}
\label{hexists}

Given the pinching function $\Phi$, a number 
$\widehat{T} \in (0, \infty)$, a positive nonincreasing
function $r \: : \: [0, \widehat{T}] \rightarrow \R$ and a number $\delta \in
\left( 0, \frac12 \right)$,
there is a 
nonincreasing 
function $h\::\:[0,\widehat{T}]\ra \R$ with
$0 < h(T) < \delta^2 r(T)$ 
so that the following property is satisfied.  
Let
$\M$ be a Ricci flow with surgery defined on $[0, T)$, with $T < \widehat{T}$,
which satisfies the {\it a priori}
assumptions (Definition \ref{apriori}) and which goes singular
at time $T$. Let $(\Omega, \overline{g})$ denote the time-$T$ limit, in the sense of
Section \ref{II.3}.
Put $\rho = \de\; r(T)$ and
\begin{equation}
\Om_\rho=\{(x,t)\in \Omega \mid \overline{R}(x,T)\leq \rho^{-2}\}.
\end{equation}
\noindent
Suppose that $(x,T)$
lies in an $\eps$-horn $\h\subset \Omega$ whose boundary is contained in 
$\Om_\rho$.  
Suppose also that  $\overline{R}(x,T)\geq h^{-2}(T)$. 
Then the parabolic
region $P(x,T,\de^{-1}\overline{R}(x,T)^{-\frac12},-\overline{R}(x,T)^{-1})$ is
contained in a strong $\de$-neck. (As usual, $\epsilon$ is 
a fixed constant that is small enough so that the result holds
uniformly  with respect to the other variables.)
\end{lemma}
\begin{proof}
Fix $\delta \in (0,1)$. Suppose that the claim is not true.
Then there is a sequence of
Ricci flows with surgery $\M^\al$  and points $(x^\alpha, T^\alpha) \in
\M^\al$ with $T^\alpha < \widehat{T}$ such that \\
1. $\M^\al$ satisfies 
the $\Phi$-pinching
and $r$-canonical neighborhood assumptions, \\
2. $\M^\al$ goes singular at time $T^\al$,\\
3. $(x^\al,T^\al)$ belongs to an $\eps$-horn $\h^\al\subset \Omega^\alpha$
whose boundary is contained in $\Om_\rho^\al$, and \\
4. $\overline{R}(x^\al,T^\al) \rightarrow \infty$, but \\
5. For each $\alpha$,
 $P(x^\al,T^\al,\de^{-1}\overline{R}(x^\al,T^\al)^{-\frac12},
 -\overline{R}(x^\al,T^\al)^{-1})$
is not contained in a strong $\delta$-neck.  

Recall that when $\eps$ is small enough, any cross-sectional $2$-sphere
sitting in an $\eps$-neck $V\subset \h^\al$ 
separates the ends of $\h^\al$;
see Section \ref{notationterminology}.  We may
find a properly embedded minimizing geodesic $\ga^\al\subset \h^\al$ 
which joins the two
ends of $\h^\al$.   As $\ga^\al$ must intersect a cross-sectional
$2$-sphere containing $(x^\al,T^\al)$, it must pass within distance 
$\leq 10\overline{R}(x^\al,T^\al)^{-\frac12}$ of $(x^\al,T^\al)$,
when $\eps$ is small.  
Let $y^\al$ be the endpoint of $\gamma^\al$ contained in
$\Omega^\al_\rho$ and let $\widehat{y}^\al$ be the first point,
moving along $\gamma^\al$ from the noncompact end of $\h^\al$
toward $y^\al$,  
where $\overline{R}(\widehat{y}^\al, T^\al) = \rho^{-2}$. As the gradient bound
$|\nabla \overline{R}^{- \: \frac12}| \le \frac12 \eta$ is valid along
$\gamma^\al$ starting from $\widehat{y}^\al$ and going out the noncompact end
(since such points on $\ga^\al$ have scalar curvature greater than
$r(T^\al)^{-2}$), we have $\dist_{T^\al}(x^\al, y^\al) \ge 
\dist_{T^\al}(x^\al, \widehat{y}^\al) \ge \frac{2}{\eta}
\left(\rho - \overline{R}(x^\al, T^\al)^{- \: \frac12} \right)$.
Let $L^\al$ denote the time-$T^\al$ distance from $x^\al$ to the other
end of $\h^\al$. Since $\overline{R}$ goes to infinity as one exits the end,
Lemma \ref{claim2II.4.2} implies that $\lim_{\al \rightarrow \infty}
\overline{R}(x^\al, T^\al)^{\frac12} L^\al \: = \: \infty$.
From the existence of $\ga^\al$, whose length in either direction from
$x^\al$ is large compared to $\overline{R}(x^\al, T^\al)^{-\frac12}$,
it is clear that for large $\al$,
the canonical neighborhood of $(x^\al, T^\al)$ must be
of type (a) or (b) in the terminology of Definition \ref{canonicalnbhddef}. 
By Lemmas \ref{lemoverlinerisproper} and \ref{claim2II.4.2},
we also know that for any fixed $\sigma < \infty$, for large
$\al$ the ball $B(x^\al, T^\al, \sigma \overline{R}(x^\al, T^\al)^{-\frac12})$
has compact closure in the time-$T^\al$ slice of $\M^\al$.
 
By Lemma \ref{claim2II.4.2}, after rescaling the metric on the
time-$T^\alpha$ slice by $\overline{R}(x^\al,T^\al)$ we have uniform curvature
bounds on distance balls.
We also have a uniform lower bound on the injectivity radius at
$(x^\al,T^\al)$ of the rescaled solution, in view of its
canonical neighborhood.  Hence after passing to
a subsequence, we may take a pointed 
smooth complete limit $(M^\infty, x^\infty, g_\infty)$
of the time-$T^\al$ slices, where the derivative bounds needed to
take a smooth limit come from the canonical neighborhood assumption.
By the $\Phi$-pinching
assumption, $M^\infty$ will have nonnegative curvature.

After passing to a subsequence, we can also assume that 
the $\ga^\al$'s converge to
a minimizing geodesic $\gamma$ in $M^\infty$ that passes within distance $10$
from $x^\infty$. 
The
rescaled length of $\gamma^\al$ from $x^\al$ to $y^\al$ is
bounded below by $\frac{2}{\eta} \left( \overline{R}(x^\al, T^\al)^{\frac12} \rho - 1 \right)$,
which tends to infinity as $\al \rightarrow \infty$.
We have shown that the rescaled length of $\gamma^\al$ from 
$x^\al$ to the other end of $\h^\al$ also tends to infinity as
$\al \rightarrow \infty$.
It follows that $\gamma$ is bi-infinite. 
Thus by Toponogov's
theorem, $M^\infty$ splits off an $\R$-factor. Then for large $\al$,
the canonical neighborhood of $(x^\al,T^\al)$
must be an $\epsilon$-neck, and $M^\infty = \R \times S^2$
for some positively curved metric on $S^2$.  In particular,
$M^\infty$ has scalar curvature uniformly bounded above.

Any point
$\widehat{x} \in M^\infty$ is a limit of points
$(\widehat{x}^\alpha, T^\alpha) \in \M^\alpha$. As 
$R_\infty(\widehat{x}) > 0$ and $R(x^\al,T^\al) \rightarrow \infty$, it
follows that $R(\widehat{x}^\alpha, T^\alpha) \rightarrow \infty$. Then
for large $\alpha$, $(\widehat{x}^\al,T^\al)$ is in a canonical
neighborhood which, in view of the $\R$-factor in $M^\infty$, must 
be a strong $\epsilon$-neck.  From the upper bound on the scalar
curvature of $M^\infty$, along with the time interval involved in
the definition of a strong $\epsilon$-neck, it follows that
we can parabolically rescale the pointed flows 
$(\M^\al,x^\al,T^\al)$ by $R(x^\al,T^\al)$, shift
time and
extract a smooth pointed limiting Ricci flow 
$(\M^\infty,x^\infty,0)$ which is defined on a time
interval $(\xi,0]$, for some $\xi<0$.  

In view of the strong $\epsilon$-necks around the
points $(\widehat{x}^\alpha, T^\alpha)$, if we take 
$\xi$ close to zero
then we are ensured that the Ricci flow
$(\M^\infty,x^\infty,0)$ has positive
scalar curvature $R_\infty$. Given $(\widehat{x}, t) \in \M^\infty$, 
as $R_\infty(\widehat{x}, t) > 0$ and $R(x^\al,T^\al) \rightarrow \infty$, 
the $\Phi$-pinching implies that
the time-$t$ slice $\M^\infty_t$ has nonnegative curvature at
$\widehat{x}$. Thus $\M^\infty$ has nonnegative curvature.
The time-$0$ slice $\M^\infty_0$
splits off an $\R$-factor,
which means that the same will be true of all time slices;
cf. the proof of Lemma \ref{lemclaim2}. Hence 
$\M^\infty$ is a product Ricci flow.

Let $\xi$ be the minimal negative number so that after parabolically rescaling
the pointed flows 
$(\M^\al,x^\al,T^\al)$ by $R(x^\al,T^\al)$, we can extract a limit
Ricci flow $(\M^\infty, (x^\infty,0), g_\infty(\cdot))$ which is the product of $\R$ with 
a positively curved Ricci flow on $S^2$, and is defined on the
time interval $(\xi, 0]$.  We claim that $\xi =  - \infty$. Suppose
not, i.e. $\xi > - \infty$. Given $(\widehat{x}, t) \in \M^\infty$, 
as $R_\infty(\widehat{x}, t) > 0$ and 
$R(x^\al,T^\al) \rightarrow \infty$, it follows that
$(x, t)$ is a limit of points $(\widehat{x}^\al, t^\al) \in \M^\alpha$ 
that lie in canonical
neighborhoods.  In view of the $\R$-factor in $\M^\infty$, 
for large $\al$ these
canonical neighborhoods must be strong $\epsilon$-necks. 
This implies in particular that 
$\left( R_\infty^{-2} \: \frac{\partial R_\infty}{\partial t}\right)
(\widehat{x},t) > 0$, so
$\frac{\partial R_\infty}{\partial t}(\widehat{x},t) > 0$.
Then there
is a uniform upper bound $Q$ for the scalar curvature on $\M^\infty$. 
Extending  backward from a time-$(\xi + \frac{1}{100Q})$ slice, 
we can construct a limit Ricci flow that exists on some time interval
$(\xi^\prime, 0]$ with $\xi^\prime < \xi$. As before, using the
strong $\epsilon$-neck condition and the $\Phi$-pinching, if
$\xi^\prime$ is sufficiently close to $\xi$ then we are ensured
that the Ricci flow on $(\xi^\prime, 0]$ is the product of $\R$
with a positively curved Ricci flow on $S^2$. This is a contradiction.

Thus we obtain an ancient solution $\M^\infty$ 
with the property that each point
$(x, t)$ lies in a strong $\epsilon$-neck. Removing the $\R$-factor gives an
ancient solution on $S^2$. In view of the fact that each time slice is
$\epsilon$-close to the round $S^2$, up to rescaling, it follows that
the ancient solution on $S^2$ must be the standard shrinking
solution (see Sections \ref{I.11.3} and \ref{altpf11.3}). 
Then $\M^\infty$ is the standard shrinking solution
on $\R \times S^2$. Hence for an infinite number of
$\alpha$,  $P(x^\al,T^\al,\de^{-1}R(x^\al,T^\al)^{-\frac12},
-R(x^\al,T^\al)^{-1})$ is in fact in a strong $\delta$-neck, which is a
contradiction. 
\end{proof}

 \begin{remark}
If a given $h$ makes Lemma \ref{hexists} work for a given function $r$ then
one can check that logically, $h$ also works for any $r^\prime$
with $r^\prime \ge r$. Because of this, we may assume that $h$ only
depends on $\min r=r(T)$ and is monotonically nondecreasing as a
function of $r(T)$. Similarly, if a given $h$ makes Lemma \ref{hexists}
work for a given value of $\delta$ then $h$ also works for any
$\delta^\prime$ with $\delta^\prime  \ge\delta$.  Thus we may assume
that $h$ is monotonically nondecreasing as a function of $\delta$.
\end{remark} 

\section{Surgery and the pinching condition} \label{surgery}

This section describes how one can take a $\de$-neck 
satisfying the time-$t$ Hamilton-Ivey pinching condition,
and perform surgery so as to obtain a new manifold which also satisfies
the time-$t$ pinching condition, and which is $\de'$-close to the
standard solution modulo rescaling. 
Here $\de'$ is a nonexplicit function of $\de$ but satisifes the important
property that  $\de'(\de)\ra 0$ as $\de\ra 0$.

The main geometric idea which handles the delicate part of the surgery
procedure is contained in the following lemma.  It says that one can ``round off''
the boundary of an approximate round half-cylinder so as to simultaneously 
increase the scalar curvature  and  the minimum of sectional curvature at each
point.

As the statement of the following lemma involves the curvature
operator, we state our conventions. If $M$ has constant sectional
curvature $k$ then the curvature operator acts on $2$-forms as
multiplication by $2k$.  This is consistent with the usual
Ricci flow literature, e.g. \cite{Chow-Knopf}.

Recall that $\epsilon$ is our global parameter, which is taken
sufficiently small.

\begin{lemma}
\label{lemroundoff}
Let $g_{cyl}$ denote the round cylindrical metric  of scalar
curvature $1$ on $\R \times S^2$. Let $z$ denote the coordinate 
in the $\R$-direction.
Given $A>0$, suppose that $f \: : \: (-A,0]\ra \R$ is a 
smooth function such that 

$\bullet$ $f^{(k)}(0)=0$ for all $k\geq 0$.

$\bullet$  On $(-A,0)$, 
\begin{equation}
f(z) <0\,,\quad f^\prime(z) >0\,,\quad f^{\prime \prime}(z)<0.
\end{equation}

$\bullet$ $\|f\|_{C^2}<\eps$.

$\bullet$  For every $z \in (-A,0)$, 
\begin{equation}
\label{eqnderivcomparision}
\max\left(|f(z)|,\left|f^\prime(z) \right|\right)
\leq \eps\:\left| f^{\prime \prime}(z) \right|.
\end{equation}

Then if  $h_0$ is a smooth metric on $(-A,0] \times S^2$
with $\|h_0-g_{cyl}\|_{C^2}<\eps$  and we set
 $h_1= e^{2f(z)}h_0$, it follows that for all $p \in (-A,0) \times S^2$ we have
$R_{h_1}(p) \: > \: R_{h_0}(p) \: - \: f^{\prime\prime}(z(p))$.
Also, if $\la_1(p)$ denotes the lowest eigenvalue of the curvature operator at $p$ then $\la_1^{h_1}(p)>\la_1^{h_0}(p) \: - \: f^{\prime\prime}(z(p))$.
\end{lemma}
\begin{proof}
We will use the variational characterization of $\lambda_1(p)$ :
\begin{equation} \label{Rvar}
\lambda_1(p) \: = \:
\inf_{\omega \neq 0} \frac{\omega_{ij} \: R^{ij}_{\: \: \: kl} \: \omega^{kl}}{\omega_{ij} \: \omega^{ij}}
\end{equation}
where $\omega \in \Lambda^2(T_pM)$.
We will also use the following formulas about curvature quantities for
conformally related metrics in dimension $3$ :
\begin{equation} \label{confscalar}
R_{h_1} \: = \: e^{-2f} \: \left( R_{h_0} \:  - \: 4 \: \triangle f \: - \: 2 \: |\nabla f|^2 \right)
\end{equation}
and
\begin{equation} \label{confcurv}
R^{ij}_{\: \: \: kl}(h_1) \:  = \: e^{-2f} \left( R^{ij}_{\: \: \: kl}(h_0) \: - \: \widetilde{f}^i_{\: \: k} \: 
\delta^j_{\: \: l}
 + \: \widetilde{f}^i_{\: \: l} \: 
\delta^j_{\: \: k}
 + \: \widetilde{f}^j_{\: \: k} \: 
\delta^i_{\: \: l}
 - \: \widetilde{f}^j_{\: \: l} \: 
\delta^i_{\: \: k} \: - \: |\nabla f|^2 
\left( \delta^i_{\: \: k} \: \delta^j_{\: \: l} \: - \:
\delta^i_{\: \: l} \: \delta^j_{\: \: k} \right)
\right),
\end{equation}
where $\widetilde{f}_{ij} \: = \: f_{;ij} \: - \: f_{;i} \: f_{;j}$. That is,
$\widetilde{f} \: = \: \Hess(f) \: - \: df \: \otimes \: df$. The right-hand sides of these
expressions are
computed using the metric $h_0$.

To motivate the proof, let us first consider the linearization of these expressions around
$h_0$. Keeping only the linear terms in $f$ gives 
to leading order,
\begin{equation}
R_{h_1} \: \sim \: R_{h_0} \:  {-2f} \:  R_{h_0} \:  - \: 4 \: \triangle f
\end{equation}
and
\begin{equation}
R^{ij}_{\: \: \: kl}(h_1) \:  \sim \: R^{ij}_{\: \: \: kl}(h_0) \: - \: {2f} \:  R^{ij}_{\: \: \: kl}(h_0) 
\: - \: {f}^{ \: i}_{; \: \: k} \:  \delta^j_{\: \: l}
\: + \: {f}^{ \: i}_{; \: \: l} \:  \delta^j_{\: \: k}
\: + \: {f}^{ \: j}_{; \: \: k} \:  \delta^i_{\: \: l}
\: - \: {f}^{ \: j}_{; \: \: l} \:  \delta^i_{\: \: k}.
\end{equation}
From the assumptions, $R_{h_0} \sim 1$ and $f < 0$ on $(-A,0) \times S^2$.
As $h_0$ is close to $g_{cyl}$, we have $\triangle f \sim f^{\prime \prime}(z)$,
so ${-2f} \:  R_{h_0} \:  - \: 4 \: \triangle f \: \ge \: - \: f^{\prime \prime}(z)$.
Similarly, in the case of $g_{cyl}$ a minimizer
$\omega$ in (\ref{Rvar}) is of the  form $\omega \: = \: X \wedge \partial_z$,
where $X$ is a unit vector in the $S^2$-direction.  As $h_0$ is close to
$g_{cyl}$, a minimizing $\omega$ for $h_0$ will be close to something of
the form $X \wedge \partial_z$. Then
\begin{equation}
\lambda_1^{h_1} \sim \lambda_1^{h_0} \: - \: 2f \: \lambda_1^{h_0}\: - \: 2 \:
f^{\prime \prime}(z) \: \ge \: \lambda_1^{h_0} \: + \: 2 f(z) \: |\lambda_1^{h_0}|\: - \: 2 \:
f^{\prime \prime}(z).
\end{equation}
As $h_0$ is close to $g_{cyl}$, $\lambda_1^{h_0}$ is close to
$\lambda_1^{g_{cyl}} = 0$. Then we can
use (\ref{eqnderivcomparision}) to say that $2 f(z) \: |\lambda_1^{h_0}|\: - \: 2 \:
f^{\prime \prime}(z) \: \ge \: - \: f^{\prime \prime}(z)$.

The remaining issue is to show that the increase in $R$ and $\lambda_1$
coming from the linear approximation is still approximately valid in the
nonlinear case, provided that $\epsilon$ is sufficiently small. For this,
we have to show that the increase from the linear approximation dominates
the error terms that we have neglected.

To deal with the scalar curvature first, from (\ref{confscalar}) we have
\begin{align}
R_{h_1} \: & = \: e^{-2f} \left( R_{h_0} \: - \: 4 \: \triangle_{g_{cyl}} f \right) \: + \:
4 e^{-2f}  \: (\triangle_{g_{cyl}} f \: - \: \triangle f) \:- \: 2 \: e^{-2f} \: |\nabla f|^2 \\
& \ge \: R_{h_0} \: - \: 4 \: f^{\prime \prime}(z)  \: + \:
4 e^{-2f}  \: (\triangle_{g_{cyl}} f \: - \: \triangle f) \:- \: 2 \: e^{-2f} \: |\nabla f|^2. \notag
\end{align}
Next, there is an estimate of the form
\begin{align}
|\triangle_{g_{cyl}} f \: - \: \triangle f| \: & \le \: \const \parallel h_0 - g_{cyl}
\parallel_{C^2} \: \left( |f(z)| \: + \: |f^\prime(z)| \: + \:
|f^{\prime \prime}(z)| \right) \\
& \le \:
\const \: \epsilon \: \left( |f(z)| \: + \: |f^\prime(z)| \: + \:
|f^{\prime \prime}(z)| \right). \notag
\end{align}
As $e^{-2f} \: \le \: e^{2 \epsilon}$, if $\epsilon$ is small then 
\begin{equation}
\left| e^{-2f}  \: (\triangle_{g_{cyl}} f \: - \: \triangle f) \right| \: \le \: \const \:
\epsilon \: \left( |f(z)| \: + \: |f^\prime(z)| \: + \:
|f^{\prime \prime}(z)| \right)
\end{equation}
Similarly, 
\begin{equation}
e^{-2f} \: |\nabla f|^2 \: \le \: \const  \: 
|f^\prime(z)|^2 \: \le \: \const \: \epsilon \: |f^\prime(z)|.
\end{equation}
When combined with (\ref{eqnderivcomparision}), if $\epsilon$ is taken
sufficiently small
then
\begin{equation}
- \: 4 \: f^{\prime \prime}(z)  \: + \:
4 e^{-2f}  \: (\triangle_{g_{cyl}} f \: - \: \triangle f) \:- \: 2 \: e^{-2f} \: |\nabla f|^2
\: \ge \: - \: f^{\prime \prime}(z).
\end{equation}
This shows the desired estimate for $R_{h_1}(p)$.

To estimate $\lambda_1^{h_1}$ we use
(\ref{Rvar}) and (\ref{confcurv}) to write
\begin{equation} \label{min1}
\lambda_1^{h_1}(p) \: = \: e^{-2f(z)} \: \left(
\inf_{\omega \neq 0} \frac{\omega_{ij} \: R^{ij}_{\: \: \: kl}(h_0) \: \omega^{kl}
\: - \: 4 \: \omega_{ij} \widetilde{f}^i_{\: \: k} \omega^{kj}}{\omega_{ij} \: \omega^{ij}} \: - \: 2 \: |\nabla f|^2(z) \right).
\end{equation}
Comparing with
\begin{equation}
\lambda_1^{h_0}(p) \: = \: 
\inf_{\omega \neq 0} \frac{\omega_{ij} \: R^{ij}_{\: \: \: kl}(h_0) \: \omega^{kl}}{
\omega_{ij} \: \omega^{ij}}
\end{equation}
gives
\begin{equation}
\lambda_1^{h_0}(p) \: \le \: e^{2f(z)} \: \lambda_1^{h_1}(p) \: + \: 
\frac{4 \: \omega_{ij} \widetilde{f}^i_{\: \: k} \omega^{kj}}{\omega_{ij} \: \omega^{ij}} \: + \: 2 \: |\nabla f|^2(z),
\end{equation}
where $\omega$ is a minimizer in (\ref{min1}), or
\begin{equation}
\lambda_1^{h_1}(p) \: \ge \: e^{- 2f(z)} \: \lambda_1^{h_0}(p) \: - \:
4 \: e^{- 2f(z)} \: 
\frac{\omega_{ij} \widetilde{f}^i_{\: \: k} \omega^{kj}}{\omega_{ij} \: \omega^{ij}}
\: - \: 2 \: e^{- 2f(z)} \: |\nabla f|^2(z).
\end{equation}
Using the variational formula (\ref{Rvar}), one can show that
$|\lambda_1^{h_0}(p)| \: \le \: \const \epsilon$.
From eigenvalue perturbation theory
\cite[Chapter 12]{Reed-Simon},
$\omega$ will be of the
form $X \wedge \partial_z \: + \: O(\epsilon)$ for some unit vector $X$ tangential
to $S^2$. Then we get an estimate
\begin{equation}
\lambda_1^{h_1}(p) \: \ge \: \lambda_1^{h_0}(p) \: - \:
2 f^{\prime \prime}(z) \: - \:
\const \: \epsilon \: \left( |f(z)| \: + \: |f^\prime(z)| \: + \:
|f^{\prime \prime}(z)| \right).
\end{equation}
From (\ref{eqnderivcomparision}), if $\epsilon$ is taken sufficiently small
then $\lambda_1^{h_1}(p) \: - \: \lambda_1^{h_0}(p) \: \ge \:  - \:
f^{\prime \prime}(z)$.
\end{proof}

Recall that the initial condition $\s_0$ for the standard solution  is
an $O(3)$-symmetric metric 
$g_0$
 on $\R^3$ with nonnegative curvature operator, whose
end is isometric to a round half-cylinder of scalar curvature $1$.
To facilitate the surgery procedure, we will assume that some metric 
ball around
the $O(3)$-fixed point has constant positive curvature. Outside of
this ball we use radial coordinates 
$(z, \theta) \in (-B,\infty) \times S^2$, with
$g_{0} \: = \: e^{2F(z)} g_{cyl}$. 
Here $g_{cyl}$ is the round cylindrical metric of scalar curvature one
and $F \in  C^\infty(-B,\infty)$.

\begin{lemma} \label{capoff}
Given $A > 0$, we can choose $B > A$ and
$F  \in  C^\infty(-B,\infty)$ so that \\
1. $F\equiv 0$ on $[0,\infty) \times S^2$. \\
2. The restriction of $F$ to $(-A,0] \times S^2$ 
satisfies the hypotheses of Lemma \ref{lemroundoff}.\\
3. The metric $e^{2F(z)} g_{cyl}$ on $(-B,\infty) \times S^2$ 
has nonnegative sectional curvature and extends
smoothly to a metric on $\R^3$
by adding a ball of constant positive curvature at
$\{-B\} \times S^2$.
\end{lemma}
\begin{proof}
For a metric of the form $e^{2F(z)} g_{cyl}$, one computes that
the sectional curvatures are
$ - \: e^{- \: 2F} \: F^{\prime \prime}$ and $e^{- \: 2F}
\left( \frac12 \: - \: (F^\prime)^2 \right)$.
In particular, the conditions for positive sectional curvature
are $F^{\prime \prime} < 0$ and $|F^\prime| < 
\frac{1}{\sqrt{2}}$.

The $3$-sphere of constant sectional curvature $k^2$,
with two points removed, has a metric given by
\begin{equation}
F_k(z) \: = \:  \log \left( \frac{\sqrt{2}}{k} \right) \: + \: 
\frac{1}{\sqrt{2}} \: z \: - \: \log \left( 1 + e^{\sqrt{2} z} \right).
\end{equation}
(Shifting $z$ gives other metrics of constant curvature $k^2$.
We have normalized so that the $z=0$ slice is the slice of maximal area.)
Note that the derivative
\begin{equation}
D(z) \: = \: \frac{1}{\sqrt{2}} \: - \:
\sqrt{2} \: \frac{e^{\sqrt{2} z}}{1+e^{\sqrt{2} z}}
\end{equation}
is independent of $k$.

Given $A > 0$, we take $F$ to be $0$ on $[0, \infty)$ and of the form
$c_1 \: e^{c_2/z}$ on $(-A, 0]$. We can take the constant
$c_1 > 0$ sufficiently small and the constant $c_2 < \infty$ sufficiently large
so that the hypotheses of Lemma \ref{lemroundoff} are satisfied. It
remains to smoothly cap off 
$([-A, \infty) \times S^2, e^{2F(z)} g_{cyl})$ with something of positive
sectional curvature.

With our given choice of $F \big|_{(-A, 0]}$, we have
$F^\prime(-A) \in (0, \epsilon)$. As 
$\lim_{z \rightarrow -\infty} 
D(z) \: = \: \frac{1}{\sqrt{2}}$, 
we can choose $B > A$
so that $D(-B) > F^\prime(-A)$. As $F^{\prime \prime}(-A) < 0$ and
$D^{\prime}(-B) < 0$, we can extend $F^\prime$ to a smooth function
$\widetilde{D} \: : \: (-B, \infty) \rightarrow \left( 0, 
\frac{1}{\sqrt{2}} \right)$ 
which has $\widetilde{D}^\prime < 0$ and which
coincides with $D$ on a small interval $(-B, -B + \delta)$.
Putting 
\begin{equation}
F(z) \: = \: F(0) \: + \: \int_0^z \widetilde{D}(w) \: dw,
\end{equation}
we obtain 
$F \in C^\infty(B, \infty)$ which coincides with $F_k$ on
$(-B, -B + \delta)$, for some $k > 0$. Then we can glue on a round metric
ball of constant curvature $k^2$ to $\{-B\} \times S^2$, in order to
obtain the desired metric.
\end{proof}

In the statement of the next lemma we continue with the metric
constructed in Lemma \ref{capoff}.
\begin{lemma}
\label{preservingphipinching}
There exists $\de'=\de'(\de)$ with $\lim_{\de\ra 0}\de'(\de)=0$ and a constant
$\de_0>0$
 such that the following holds.  Suppose that $\de<\de_0$, $x\in \{0\} \times S^2$
and $h_0$ is a Riemannian metric 
on $(-A,\frac{1}{\de}) \times S^2$ with $R(x) > 0$ such that:

\medskip
$\bullet$ $h_0$ satisfies the time-$t$ Hamilton-Ivey pinching condition of
Definition \ref{HamIvey}.

$\bullet$  $R(x)h_0$ is 
$\de$-close to $g_{cyl}$ in the $C^{[\frac{1}{\de}]+1}$-topology.

\bigskip
Then there is a smooth metric 
$h$ on $\R^3 = D^3 \: \cup \: \left( (-B,\frac{1}{\de}) \times S^2 \right)$
such that 

\medskip
$\bullet$ $h$ satisfies the time-$t$ pinching condition.

$\bullet$ The restriction of $h$ to $[0,\frac{1}{\de}) \times S^2$ is $h_0$.

$\bullet$  The restriction of $R(x)h$ to $(-B,-A) \times S^2$ is 
$g_{0}$, the initial metric of a standard solution.

$\bullet$  The restriction of $R(x)h$ to $D^3$ has constant curvature $k^2$.

$\bullet$ $R(x)h$ is $\de'$-close to $e^{2F}g_{cyl}$ in the $C^{[\frac{1}{\de'}]+1}$-topology on 
$\left( -B, \frac{1}{\delta} \right) \times S^2$.

\end{lemma}
\begin{proof}
Put
\begin{equation}
U_1 \: = \:  (-B,-\frac{A}{2}) \times S^2,\quad U_2 \: = \:  (-A,\frac{1}{\de}) \times S^2
\end{equation}
and let $\{\al_1,\al_2\}$ be a $C^{\infty}$ partition of unity subordinate to the open
cover $\{U_1,U_2\}$ of $(-B,\frac{1}{\de}) \times S^2$.  We set
\begin{equation}
h \: = \:  \al_1 \: R(x)^{-1} \: g_{0} \: + \: \al_2 \: e^{2F} \: h_0
\end{equation}
on $\left( -B, \frac{1}{\delta} \right) \times S^2$ and cap it off with a
$3$-ball of constant curvature $k$, as in Lemma \ref{capoff}.

Given $\de'$, we claim that if $\de$ is sufficiently small then the conclusion of
the lemma holds.
The only part of the lemma that is not obvious is the pinching
condition.  Note that on $\left( -\frac{A}{2}, \frac{1}{\delta} \right) \times S^2$ the metric
$h$ agrees with $e^{2F}h_0$ and hence, when $\de$ is sufficiently
small, the pinching
condition will hold on $\left( -\frac{A}{2}, \frac{1}{\delta} \right) \times S^2$
by Lemmas \ref{lemroundoff} and \ref{lempinchingmonotonic}.
On the other hand, when $\de$ is sufficiently small, the 
restrictions of the metrics 
$g_{0} \: =\: e^{2F}g_{cyl}$
and $R(x)e^{2F}h_0$ to $(-A,- \frac{A}{2}) \times S^2$
will be very close and will have strictly positive
curvature. (The positive curvature for $e^{2F}h_0$ also follows from Lemma
\ref{lemroundoff}; if $\delta$ is small enough then $\lambda_1^{h_0}$ will be close
to zero, while $- f^{\prime \prime}(z)$ is strictly positive for
$z \in (-A,- \frac{A}{2})$.) Thus $h$ will have positive curvature on $(-B,- \frac{A}{2}) \times S^2$ 
and the pinching condition will hold there.
\end{proof}

We have now fixed the initial condition 
$g_{0}$
for a
standard solution, along with the procedure to meld 
$g_{0}$
to
an approximate cylinder.

\section{II.4.4. Performing surgery and continuing flows}
\label{II.4.4}
This section discusses the surgery procedure  and shows how to prolong a Ricci 
flow with surgery, provided that the {\it a priori} assumptions hold. 

\begin{definition}(Ricci flow with cutoff)
\label{rdeltacutoffflowdef}
Suppose that $a\geq 0$ and let $\M$ be a Ricci flow with surgery 
defined on $[a,b]$ that satisfies the {\it a priori}
assumptions of Definition \ref{apriori}.
Let $\de:[a,b]\ra (0,\delta_0)$ be a nonincreasing function, where
$\delta_0$ is the parameter of Lemma \ref{preservingphipinching}.
Then $\M$ is a 
{\em  Ricci flow with $(r, \de)$-cutoff} if 
at each singular time $t$, the forward time slice
$\M_{t}^+$ is obtained from the backward time slice $\Omega = 
\M_{t}^-$
by applying the following procedure:

A.  Discard each component of $\Omega$ that does not
intersect
\begin{equation}
\Om_{\rho} = \{(x,t)\in \Omega \mid R(x,t)\leq\rho^{-2}\},
\end{equation}
where $\rho= \de(t)r(t)$.
 
B. In each $\eps$-horn 
$\h_{ij}$ of each of the remaining components
$\Om_i$, find a point $(x_{ij},t)$
such that $R(x_{ij},t)=h^{-2}$, where 
$h=h(t)$
is as in Lemma \ref{hexists}.
  
C. Find a strong $\de$-neck $U_{ij}\times [t-h^2,t]$ 
containing 
$P(x_{ij},t,\de^{-1}R(x_{ij},t)^{-\frac12},-R(x_{ij},t)^{-1})$;  this 
is guaranteed to exist by Lemma \ref{hexists}.
  
D.  For each $ij$, let $S_{ij}\subset U_{ij}$ be
a cross-sectional $2$-sphere containing $(x_{ij},t)$.
Cut $\bigcup_i\,\Om_i$ 
along the $S_{ij}$'s and throw away the tips of the horns
$\h_{ij}$,
to obtain a compact manifold-with-boundary $X$ having a 
spherical boundary component for each $ij$.

E. Glue caps onto $X$, using Lemma \ref{preservingphipinching},  
to obtain the closed manifold $\M_{t}^+$.
\end{definition}

For concreteness, we take the parameter
$A$ of Lemma \ref{preservingphipinching} to be $10$.
The neighborhood of a boundary
component of $X$ is parametrized as $[-A, \delta^{-1}) \times S^2$,
with $S_{ij} \: = \: \{-A \} \times S^2$.
The metric on $[0, \delta^{-1}) \times S^2$ is unaltered by the surgery procedure.
The corresponding region in the new manifold $\M_{t}^+$, minus 
a metric ball of constant curvature, is
parametrized by $(-B, \delta^{-1}) \times S^2$.
Put $S_{ij}^\prime \: = \: \{0\} \times S^2 \subset \M_{t}^-$.
We will consider the part added by surgery on $\h_{ij}$
to be the $3$-disk in $\M_{t}^+$ bounded by $S_{ij}^\prime$. 
In terms of Definition \ref{flowwithsurgery}, if $t = t_k^+$ then
the subset $X_k^+$ of $\Omega_k = \M_{t}^-$ has
boundary $\bigcup_{ij} S_{ij}^\prime$. The added part
$\M_{t}^+ - X_{k+1}^-$ is a union of $3$-balls.

\begin{remark}
Our definition of surgery differs slightly from that in
\cite{Perelman2}. The paper \cite{Perelman2} has two
extra steps involving throwing away certain components
of the postsurgery manifold.
We omit these steps in order to simplify the definition of
surgery, but there is no real loss either way.

First, in the setup of \cite[Section 4.4]{Perelman2},
any component of $\M_{t}^+$ that is $\epsilon$-close
to a metric quotient of the round $S^3$ is thrown away. The
motivation of \cite{Perelman2} was to not have to include 
these in the list of canonical neighborhoods. Such
components
are topologically standard.  We do include such manifolds
in the list of canonical neighborhoods and do not
throw them away in the surgery procedure. 

Second,
when considering the long-time behavior of Ricci flow
in \cite[Section 7]{Perelman2}, 
any component of $\M_{t}^+$ which admits a
metric of nonnegative scalar curvature is thrown away.
The motivation for this extra step is that any such
component admits a metric that is either flat or has
finite extinction time.  In either case one concludes that
the component is a graph manifold and, for the purposes
of the geometrization conjecture, is standard.
(Recall the definition of graph manifolds from
Appendix \ref{geom}.)
Again, we do not throw away such components.
\end{remark}

Note that the definition of Ricci flow with $(r, \de)$-cutoff also depends
on the function $r(t)$ through the {\it a priori} assumption. We now state
how the topology of the time slice changes when going backward through the 
singular time  $t$.  Recall that for $t'<t$ close to $t$, the time slices
$\M_{t'}$ are all diffeomorphic; we refer to this diffeomorphism type
as the {\em presurgery manifold}, and the forward time slice
$\M_t^+$ as the {\em postsurgery manifold}.

\begin{lemma} \label{reconstruct2}
The presurgery manifold may be obtained
from the postsurgery manifold by applying the following
operations finitely many times:

\begin{itemize}
\item Replacing two connected components with their connected sum.
\item Taking the connected sum of a connected component 
with  $S^1\times S^2$ or $\R P^3$. 
\item Taking the disjoint union with an additional 
$S^1 \times S^2$ or an isometric quotient of the round $S^3$.
\end{itemize} 
\end{lemma}
\begin{proof}
The proof is basically the same as that of Lemma \ref{reconstruct}. The only
difference is that we must take
into account the compact components of $\Omega$ that do not intersect
$\Omega_\rho$; these are thrown away in Step A. 
(Such components did not occur in Lemma \ref{reconstruct} because
in Lemma \ref{reconstruct} we were dealing with the first surgery for
the Ricci flow on the initial connected manifold; 
see Lemma \ref{omeganoncompact}, which is valid for the first surgery time.)
Any such component is
diffeomorphic to 
$S^1 \times S^2$, $\R P^3 \# \R P^3$ or
a quotient of the round $S^3$, in view of the canonical
neighborhood assumption; see the proof of Lemma \ref{standard}. 
\end{proof}

\begin{remark}
\label{volumelossremark}
When $\de>0$ is sufficiently small, we will have 
$\vol(\M_t^+)<\vol(\M_t^-)-h(t)^3$
for each surgery time $t\in (a,b)$.
This is because each component that is discarded in step D contains
at least
``half'' of the $\de$-neck $U_{ij}$, which has volume at least $\const \de^{-1}
h(t)^3$,
while the cap added has volume at most $\const h(t)^3$.
\end{remark}

\begin{remark}
For a Ricci flow with surgery whose original manifold is nonaspherical and
irreducible,
one wants to know that the Ricci flow goes extinct within a finite time
\cite{Colding-Minicozzi,Colding-Minicozzi2,Perelman3}. 
Consider the effect of a first surgery, say at time $t$. Among the connected
components of the postsurgery
manifold $\M_{t}^+$, one will
be diffeomorphic to the presurgery manifold
 and the others will be $3$-spheres.
Let ${\mathcal N}_{t}^+$ be a component  of $\M_{t}^+$ that is diffeomorphic to 
the presurgery manifold.
By the nature of the surgery procedure, there is a 
 function
$\xi$ defined on a small interval $(t-\alpha, t)$ so that
$\lim_{t' \rightarrow t} \xi(t') = 1$ and for $t' \in (t - \alpha, t)$,
there is a homotopy-equivalence from $(\M_{t'}, g(t'))$ to ${\mathcal N}_{t}^+$ that
expands distances by at most $\xi(t')$. Following the 
subsequent evolution of
${\mathcal N}_{t}^+$, 
there is a similar statement for
the later singular times. This fact is needed in
\cite{Colding-Minicozzi,Colding-Minicozzi2,Perelman3} 
in order to control the decay of
a certain area functional as one goes through a surgery.
\end{remark}

We discuss how to continue Ricci flows after surgery.
We recall that in Definition \ref{flowwithsurgery}  of a Ricci flow with surgery defined
on an interval $[a,c]$, the final time slice $\M_c$ consists of a single
manifold $\M_c^- = \Omega$ that may or may not be singular.

\begin{lemma}(Prolongation of Ricci flows with cutoff)
\label{extendit}
Take the function $\Phi$  to be the time-dependent pinching
function  associated to 
Definition \ref{HamIvey}
in Appendix \ref{phiappendix}.
Suppose that $r$ and $\de$ are nonincreasing positive functions
defined on $[a,b]$. Let $\M$ be a Ricci flow with
$(r, \de)$-cutoff defined on an interval $[a,c]\subset [a,b]$.
Provided $\sup\de$ is sufficiently small, either 
\begin{enumerate}
\item $\M$ can be prolonged to a Ricci flow with
$(r, \de)$-cutoff defined on $[a,b]$, or
\item There is an extension of $\M$ to a Ricci flow with 
surgery
defined  on an interval $[a,T]$ with $T \in (c, b]$, where \\
a. The restriction of  the flow to any subinterval $[a, T^\prime]$,
$T^\prime < T$,
is a Ricci flow with $(r, \de)$-cutoff, but\\
b. The $r$-canonical neighborhood assumption fails at some
point $(x,T)\in \M_T^-$.
\end{enumerate}
In particular, the only obstacle to prolongation of Ricci
flows with $(r,\de)$-cutoff is the potential breakdown of the
$r$-canonical neighborhood assumption.
\end{lemma}
\begin{proof}
Consider the time slice of $\M_c = \M_c^-$ at time $c$.
If it is
singular then we perform steps A-E of Definition \ref{rdeltacutoffflowdef}
to produce $\M_c^+$; otherwise we set $\M_c^+=\M_c^-$.
Since the surgery is done using Lemma \ref{preservingphipinching},
provided $\de>0$ is sufficiently small, 
the forward time slice $\M_c^+$ will satisfy the $\Phi$-pinching 
assumption.  

We claim that the $r$-canonical neighborhood
assumption holds in $\M_c^+$.
More precisely, if a point
$(x,c) \in \M_c^+$ lies within a distance of $10 \epsilon^{-1} h$ from
the added part $\M_c^+ - \M_c^-$ then it lies in an
$\epsilon$-cap, while if $(x,c)$ lies at distance greater than
$10 \epsilon^{-1} h$ from $\M_c^+ - \M_c^-$ and has
scalar curvature greater than $r(c)^{-2}$ then it lies in a 
canonical neighborhood that was present in the presurgery manifold
$\M_c^-$.
(We are assuming that $\epsilon < \frac{1}{100}$.)
In view of  Lemma \ref{claim5},
the only point to observe is that points at distance roughly
$10 \epsilon^{-1} h$ lie in $\epsilon$-necks, as they are unaltered
by the surgery and they were in $\delta$-necks before the surgery.
This gives the $\epsilon$-neck needed to define an
$\epsilon$-cap.

We now prolong $\M$  by Ricci flow with initial
condition $\M_c^+$. If the flow extends smoothly
up to time $b$ then we are done because either the canonical neighborhood
assumption holds up to time $b$ yielding (a), or it fails at some
time in the interval $(c,b]$, and we have (b).
Otherwise, there is some time $t_{\sing} \le b$ at which it goes singular.
We add the singular limit $\Omega$ at time $t_{\sing}$ 
to obtain a Ricci flow with surgery defined on $[a, t_{\sing}]$.
From Lemma \ref{preservingphipinching},
$\M$ satisfies
the Hamilton-Ivey pinching condition of
Definition \ref{HamIvey} on $[a,t_{\sing}]$.  
As the function $r$ is nonincreasing in $t$, it follows from
Definition \ref{canonicalnbhddef} that the set of times
$t\in [c,t_{\sing}]$ for which the $r$-canonical neighborhood assumption
holds is relatively open to the right (i.e. if the $r$-canonical
neighborhood assumption holds at time $t \in [c, t_{\sing})$ then
it also holds within some interval
$[t, t^\prime)$). Thus the set of times 
$t\in [c,t_{\sing}]$ for which the $r$-canonical neighborhood assumption
holds is either an interval $[c, T)$, with $T \le t_{\sing}$, or
$[c, t_{\sing}]$.
If the set of such times $t$ is
$(c, T)$ for some $T \le t_{\sing}$ then the lemma holds.
Otherwise, the $r$-canonical neighborhood assumption
holds at $t_{\sing}$.  In this case we repeat the construction
with $c$ replaced by $t_{\sing}$, and iterate if necessary. Either we
will reach time $b$ after a finite number  of iterations, or we 
will reach a time $T$ satisfying (2), 
or we will hit an infinite number of singular times
before time $b$. However, the last possibility cannot
occur.  A singular time corresponds to a component going extinct or
to a surgery. The number of components going extinct before time $b$
can be bounded in terms of the number of surgeries before time $b$,
so it suffices to show that the latter is finite.
Each surgery removes a volume of at least $h^3$, but
the lower bound on the scalar curvature during the flow, 
coming from the maximum principle,
gives a finite upper bound on the total volume growth during the complement
of the singular times.
\end{proof}

\begin{remark} \label{30}
The condition $C_1 \ge 30 \epsilon^{-1}$ in 
Definition \ref{canonicalnbhddef}
was in order to ensure that the $\epsilon$-cap coming from a 
surgery satisfies the requirements to be a canonical neighborhood.
\end{remark}

\section{II.4.5. Evolution of a surgery cap}
\label{secII.4.5}

Let $\M$ be a Ricci flow with $(r,\de)$-cutoff.
The next result says that provided $\de$ is small,
after a surgery at scale $h$
there is a ball $B$  of radius $Ah\gg h$ centered in the 
surgery cap whose evolution is close to that
of a standard solution for an elapsed time close to $h^2$, 
unless another surgery occurs  during which
the entire ball is thrown away.   Note that the elapsed time $h^2$
corresponds, modulo parabolic rescaling, to the duration
of the standard solution.

\begin{lemma}(cf. Lemma II.4.5)
\label{II.4.5}\\
For any $A<\infty$, $\th\in (0,1)$ and $\hat r>0$, 
one can find $\hat\de=\hat\de(A,\th,\hat r)>0$
with the following property.  Suppose that we
have a Ricci flow with $(r, \de)$-cutoff
defined on a time interval $[a,b]$ with  $\min r=r(b)\geq \hat r$. 
Suppose that there is a surgery time 
$T_0\in (a,b)$, with $\de(T_0)\leq \hat\de$. Consider a given
surgery at the surgery time and let
$(p,T_0)\in\M_{T_0}^+$ 
be the center of the surgery cap.
Let $\hat h= h(\de(T_0),\eps,r(T_0),\Phi)$ be the surgery 
scale given by Lemma \ref{hexists} and put 
$T_1=\min(b,T_0+\th \hat h^2)$.
Then one of the two following possibilities occurs : \\
(1) The solution is unscathed on $P(p,T_0,A\hat h,T_1-T_0)$.
The pointed solution there (with respect to the basepoint
$(p, T_0)$) is, modulo parabolic
rescaling, $A^{-1}$-close to the pointed flow on
$U_0\times [0,(T_1-T_0)\hat h^{-2}]$, where
$U_0$ is an open subset of the initial time slice $\s_0$ 
of a standard solution
$\s$ and the basepoint is
the center $c$ of the cap in $\s_0$. \\
(2) Assertion (1) holds with $T_1$ replaced by some
$t^+ \in [T_0, T_1)$, where $t^+$ is a surgery time. Moreover,
the entire
ball $B(p,T_0,A\hat h)$ becomes extinct at time $t^+$,
i.e. 
$P(p,T_0,A\hat h,t_+-T_0)\cap \M_{t_+}\subset \M_{t_+}^- - \M_{t_+}^+$.
\end{lemma}
\begin{proof}
We give a proof with the same ingredients as the proof in 
\cite{Perelman2}, but which
is slightly rearranged.  
We first show the following result, which is almost the same as
Lemma \ref{II.4.5}.
\begin{lemma} \label{intermediate}
For any $A<\infty$, $\th\in (0,1)$ and $\hat r>0$, 
one can find $\hat\de=\hat\de(A,\th,\hat r)>0$
with the following property.  Suppose that we
have a Ricci flow with $(r, \de)$-cutoff
defined on a time interval $[a,b]$ with  $\min r=r(b)\geq \hat r$. 
Suppose that there is a surgery time 
$T_0\in (a,b)$, with $\de(T_0)\leq \hat\de$.  Consider a given
surgery at the surgery time and let
$(p,T_0)\in\M_{T_0}^+$
be the center of the surgery cap.
Let $\hat h= h(\de(T_0),\eps,r(T_0),\Phi)$ be the surgery 
scale given by Lemma \ref{hexists} and put 
$T_1=\min(b,T_0+\th \hat h^2)$.
Suppose that the solution is unscathed on $P(p,T_0,A\hat h,T_1-T_0)$.
Then the pointed solution there (with respect to the basepoint
$(p, T_0)$) is, modulo parabolic
rescaling, $A^{-1}$-close to the pointed flow on
$U_0\times [0,(T_1-T_0)\hat h^{-2}]$, where
$U_0$ is an open subset of the initial time slice $\s_0$ 
of a standard solution
$\s$ and the basepoint is
the center $c$ of the cap in $\s_0$.
\end{lemma}
\begin{proof}
Fix $\theta$ and $\hat{r}$.
Suppose that the lemma is not true. Then for some $A > 0$,
there is a sequence $\{\M^\al, (p^\al,T_0^\al) \}_{\al = 1}^\infty$ of 
pointed Ricci flows with $(r^\al, \delta^\al)$-cutoff
that together provide a 
counterexample. In particular,\\
1. $\lim_{\al \rightarrow \infty} 
\delta^\al (T_0^\al) \: = \: 0$. \\
2. $\M^\al$ is unscathed on
$P(p^\al, T_0^\al, A {\hat h}^\al, T_1^\al - T_0^\al)$.\\
3. If
$(\widehat{\M}^\al, (\hat{p}^\al,0))$ is the pointed Ricci
flow arising from $(\M^\al, (p^\al,T_0^\al))$ by
a time shift of $T_0^\al$ and a 
parabolic rescaling by ${\hat h}^\al$ then
$P(\hat{p}^\al, 0, A, (T_1^\al - T_0^\al)(\hat{h}^\al)^{-2})$
is not $A^{-1}$-close to a pointed
subset of a standard solution.

Put $T_2 \: = \: \liminf_{\al \rightarrow \infty} \: (T_1^\al - T_0^\al)(\hat{h}^\al)^{-2}$. (We do
not exclude that $T_2 = 0$.) Clearly $T_2 \le \: \theta$.
After passing to a subsequence, we can assume that
$T_2 \: = \: \lim_{\al \rightarrow \infty} \: (T_1^\al - T_0^\al)(\hat{h}^\al)^{-2}$.
Let $T_3$ be the supremum of the set of times 
$\tau \in [0, T_2]$ with the property that we can apply 
Appendix \ref{subsequence}, if we want,
to take a convergent subsequence of the pointed solutions 
$(\widehat{\M}^\al, (\hat{p}^\al,0))$ on the time interval $[0, \tau]$ to get
a limit solution with bounded curvature.
(In applying Appendix \ref{subsequence}, we use the case $l >0$ of 
Appendix \ref{applocalder} to get bounds on the
curvature derivatives near time $0$. In particular, if
$T_2 > 0$ then $T_3 > 0$.)
From the nature of the
surgery gluing in Lemma \ref{preservingphipinching},
since $\lim_{\al \rightarrow \infty} 
\delta^\al (T_0^\al) \: = \: 0$ 
we know that we can at least
take a limit of the pointed solutions 
$(\widehat{\M}^\al, (\hat{p}^\al,0))$ on the time interval $[0,0]$, 
so $T_3$ is well-defined.

\begin{sublemma}
$T_3 \: = \: T_2$.
\end{sublemma}
\begin{proof}
Suppose not. Consider the interval $[0, T_3)$ (where we define
$[0, 0)$ to be $\{0\}$).
Given $\sigma \in (0, T_2 - T_3]$, for any subsequence
of $\{\M^\al, (\hat{p}^\al,0)\}_{\al=1}^\infty$ (which we relabel as 
$\{\M^\al, (\hat{p}^\al,0)\}_{\al=1}^\infty$)
either \\
1. There is some $\lambda > 0$ and an
infinite number of $\al$
for which the set $B(\hat{p}^{\al}, 0, \lambda)$ becomes 
scathed on $[0, T_3 + \sigma]$, or \\
2. For each $\lambda > 0$ the set
$P(\hat{p}^{\al}, 0, \lambda, T_3 + \sigma)$
is unscathed for large $\al$, but for each
$\Lambda > 0$ there is some
$\lambda_\Lambda > 0$ such that 
$\limsup_{\al \rightarrow \infty} 
\sup_{P(\hat{p}^{\al}, 0, \lambda_\Lambda, T_3 + \sigma)} |\Rm| \: \ge \: \Lambda$.

By Appendix \ref{subsequence}, after
passing to a subsequence, there is a complete limit solution
$(\widehat{\M}^\infty, (\hat{p}^\infty,0))$ defined on the time
interval $[0, T_3)$ with bounded curvature on compact time intervals.
Relabel the subsequence by $\al$. By Lemma \ref{lemstandard0}, 
$(\widehat{\M}^\infty, (\hat{p}^\infty,0))$ must be the same as the
restriction of some standard solution to $[0, T_3)$.
From Lemma \ref{lemclaim4}, the  curvature of $\widehat{\M}^\infty$ is uniformly
bounded on $[0, T_3)$; therefore by the canonical neighborhood assumption and equation
(\ref{II.(1.3)estimates}), 
we can choose $\sigma \in (0, T_2 - T_3]$ and $\Lambda^\prime > 0$ 
so that for any $\lambda > 0$,  we have
$\limsup_{\al \rightarrow \infty} 
\sup_{P(\hat{p}^{\al}, 0, \lambda, T_3 + \sigma)} |\Rm| \: \le \: \Lambda^\prime$.
However, $\lim_{\al \rightarrow \infty} 
\delta^\al (T_0^\al) \: = \: 0$ and surgeries only occur near the
centers of $\delta$-necks. From the curvature bound on the time
interval $[0, T_3 + \sigma]$ and the length distortion
estimates of Lemma \ref{distancedist}, for a given $\lambda$ the 
balls $B(\hat{p}^{\al}, 0, \lambda)$ will stay within 
a uniformly bounded distance from $\hat{p}^{\al}$ on the time interval
$[0, T_3 + \sigma]$. Hence they
cannot be scathed on $[0, T_3 + \sigma]$ for an infinite number of $\al$,
as the collar length of the $\delta$-neck around the supposed surgery locus
would be large enough to prohibit the cap point $\hat{p}^{\al}$ 
from being within a bounded
distance from the surgery locus.
This, along with the fact that $\limsup_{\al \rightarrow \infty} 
\sup_{P(\hat{p}^{\al}, 0, \lambda, T_3 + \sigma)} |\Rm| \: \le \: \Lambda^\prime$
for all $\lambda > 0$,
gives a contradiction.
\end{proof}

Returning to the original sequence 
$\{\M^\al, (p^\al,T_0^\al) \}_{\al = 1}^\infty$ and its rescaling
$\{ \widehat{\M}^\al, (\hat{p}^\al,0) \}_{\al = 1}^\infty$, we can now
take a subsequence that converges on the time interval $[0, T_2)$,
again necessarily to a standard solution.  Then there will be
an infinite subsequence 
$\{ \widehat{\M}^{\al_\beta}, (\hat{p}^{\al_\beta},0) \}_{\beta = 1}^\infty$
of 
$\{ \widehat{\M}^{\al}, (\hat{p}^{\al},0) \}_{\al = 1}^\infty$, 
with
$\lim_{\beta \rightarrow \infty} 
(T_1^{\al_\beta} - T_0^{\al_\beta})(\hat{h}^{\al_\beta})^{-2} \: = \: T_2$,
so that $P(\hat{p}^{\al_\beta}, 0, A, (T_1^{\al_\beta} - T_0^{\al_\beta})
(\hat{h}^{\al_\beta})^{-2})$
is $A^{-1}$-close to a pointed
subset of a standard solution (by the canonical
neighborhood assumption, equation
(\ref{II.(1.3)estimates}) and Appendix \ref{applocalder}). 
This is a contradiction.
\end{proof}

We now finish the proof of Lemma \ref{II.4.5}. If the solution is
unscathed on $P(p,T_0,A\hat h,T_1-T_0)$ then we can apply
Lemma \ref{intermediate} to see
that we are in case (1) of the conclusion of Lemma \ref{II.4.5}. 
Suppose, on the other hand,
that the solution is scathed on $P(p,T_0,A\hat h,T_1-T_0)$. Let
$t^+$ be the largest $t$ so that the solution is unscathed on
$P(p,T_0,A\hat h,t-T_0)$.
We can apply
Lemma \ref{intermediate} to
see that conclusion (1) of Lemma \ref{II.4.5}
holds with $T_1$ replaced by $t^+$.
As surgery is always performed near the middle of a $\delta$-neck,
if $\hat \delta << A^{-1}$ then the final time slice in the parabolic
neighborhood $P(p,T_0,A\hat h,t^+-T_0)$ cannot intersect a
$2$-sphere where a surgery is going to be performed. 
The only other possibility is that
the entire ball $B(p,T_0,A\hat h)$ becomes extinct at time $t^+$.
\end{proof}

\section{II.4.6. Curves that penetrate the surgery region}
\label{II.4.6}

Let $\M$ be a Ricci flow with $(r,\de)$-cutoff.
The next result, Corollary \ref{corollaryII.4.6}, says that if $\de$ is sufficiently small
then an admissible curve $\ga$ which comes close
to a surgery cap at a surgery time will have a large value
of $\int_\gamma \left( R(\gamma(t))+|\dot{\ga}(t)|^2 \right) \: dt$.  
Note that the latter quantity is not quite the same as
${\mathcal L}(\gamma)$, and is invariant under parabolic rescaling.

Corollary \ref{corollaryII.4.6} is used in the extension of
Theorem \ref{nolocalcollapse} to Ricci flows with 
surgery.   The idea is that if $\de$ is small 
and $L(x,t)$ isn't too large then any $\L$-minimizing sequence of
admissible curves joining the basepoint $(x_0,t_0)$ to $(x,t)$
must avoid surgery regions, and will therefore 
accumulate on a minimizing $\L$-geodesic.

\begin{corollary} (cf. Corollary II.4.6)
\label{corollaryII.4.6}
For any $l<\infty$ and $\hat r>0$, we can find $A=A(l,\hat r)<\infty$ and
$\th=\th(l,\hat r)$ with the following property.
Suppose that we are in the situation of Lemma \ref{II.4.5}, with
$\delta(T_0) < \hat{\delta}(A, \th, \hat r)$.
As usual,
$\hat h$ will be the surgery scale coming from Lemma \ref{hexists}.
Let $\ga:[T_0,T_\ga]\ra \M$
be an admissible curve, with
$T_\ga\in (T_0,T_1]$.
Suppose that
$\ga(T_0)\in B(p,T_0,\frac{A\hat h}{2})$,
$\ga([T_0,T_\ga))\subset P(p,T_0,A\hat h, T_\ga-T_0)$, and
either 

a. $T_\ga=T_1=T_0+\th(\hat h)^2$,

or

b. $\ga(T_\ga)\in\D B(p,T_0,A\hat h)\times [T_0,T_\ga]$.

Then 
\begin{equation}
\int_{T_0}^{T_\ga} \left( R(\gamma(t),t)+|\dot{\ga}(t)|^2 \right) \: dt
\: > \: l.
\end{equation}
\end{corollary}
\begin{proof}
For the moment, fix $A<\infty$ and $\th\in (0,1)$. 
Choose $\hat\de=\hat\de(A,\th,\hat r)$
so as to satisfy
Lemma \ref{II.4.5}.   Let $\M$, $(p,T_0)$, etc.,
be as in the hypotheses of Lemma \ref{II.4.5}. Let $\ga:[T_0,T_\ga]\ra\M$
be a curve as in the hypotheses of the Corollary.
From Lemma \ref{II.4.5}, we know that there is a standard solution
$\s$ such that  the parabolic region
$P(p,T_0,A\hat h,T_\ga-T_0) \subset \M$, with basepoint $(p,T_0)$, is
(after parabolic rescaling by $\hat h^{-2}$) $A^{-1}$-close
to a  pointed flow $U_0\times [0,\hat T_\ga] \subset \s$, the latter
having basepoint $(c,0)$. Here
$U_0\subset \s_0$ and $\hat T_\ga= (T_\ga-T_0)\hat h^{-2}$.
Then the image of $\ga$, under the diffeomorphism implicit in the definition
of $A^{-1}$-closeness, gives rise to a smooth curve 
$\ga_0:[0,\hat T_\ga]\ra U_0\times [0,\hat T_\ga]$ so that
(if $A$ is sufficiently large) :

\begin{gather}
\ga_0(0)\in B(c,0,\frac35 A),\\
\int_{T_0}^{T_\ga}|\dot{\ga}|^2dt
\geq \frac12 \int_0^{\hat T_\ga}|\dot{\ga}_0|^2dt,\\
\int_{T_0}^{T_\ga}R({\ga}(t),t)dt
\geq \frac12 \int_0^{\hat T_\ga}R({\ga}_0(t),t)dt,
\end{gather}
and

(a) $\hat T_\ga=\th$,

or

(b) $\ga_0(\hat T_\ga)\not\in P(c,0,\frac45 A,\hat T_\ga)$.

\noindent
In case (a) we have, by Lemma \ref{claim5},
\begin{equation}
\int_{T_0}^{T_\ga} R(\ga(t), t)dt\geq 
\frac12 \int_0^\th R(\ga_0(t),t)dt
\geq \frac12\int_0^\th \const(1-t)^{-1}dt=\const \log(1-\th).
\end{equation}
If we choose $\th$ sufficiently close to $1$ then in this case,
we can ensure that
\begin{equation}
\int_{T_0}^{T_\ga} \left( R(\gamma(t), t)+|\dot{\ga}(t)|^2 \right) \: dt
\: \ge \: 
\int_{T_0}^{T_\ga} R(\gamma(t), t) \: dt \: > \: l.
\end{equation}

In case (b), we may use the fact that the Ricci curvature 
of the standard solution is everywhere 
nonnegative,
and hence the metric tensor is nonincreasing with time.  So if
$\pi:\s=\s_0\times [0,1)\ra\s_\th$ is projection
to the time-$\th$ slice and we put $\eta= \pi\circ\ga_0$
then
\begin{align}
\int_{T_0}^{T_\ga}|\dot{\ga}(t)|^2dt \: 
& \geq \: \frac12 \int_0^{\hat T_\ga}|\dot{\ga_0}(t)|^2dt \: 
\geq \: \frac12 \int_0^{\hat T_\ga}|\dot{\eta}(t)|^2dt
\: \geq \: 
\frac{1}{2 \hat T_\ga} \left(d(\eta(0),\eta(\hat T_\ga))\right)^2 \\
& \ge \:
\frac{1}{2} \left(d(\eta(0),\eta(\hat T_\ga))\right)^2. \notag
\end{align}
With our given value of $\th$, in view of (b), if we take $A$ large enough 
then we can ensure that $\frac12 \left(d(\eta(0),\eta(\hat T_\ga))\right)^2
\: > \: l$.  This proves the lemma.
\end{proof}

\section{II.4.7. A technical estimate}

The next result is a technical result that will not be used in the sequel.

\begin{corollary}(cf. Corollary II.4.7)
For any $Q<\infty$ and $\hat r>0$, there is a $\th=\th(Q,\hat r)\in (0,1)$
with the following property.  Suppose that we are in the situation of
Lemma \ref{II.4.5}, with
$\delta(T_0) < \hat{\delta}(A, \th, \hat r)$
and $A>\eps^{-1}$.  If $\ga:[T_0,T_x]\ra \M$ is a static curve
starting in $B(p,T_0,A\hat h)$, and 
\begin{equation}
\label{Rbound}
Q^{-1}R(\ga(t))\leq R(\ga(T_x))\leq Q(T_x-T_0)^{-1}
\end{equation}
for all $t\in [T_0,T_x]$,  then $T_x\leq T_0+\th\hat h^2$.
\end{corollary}
\begin{remark}
The hypothesis (\ref{Rbound}) in the corollary means
that in the scale of the scalar curvature $R(\gamma(T_x))$ at the
endpoint $\gamma(T_x)$, the scalar curvature on $\gamma$ is bounded
and the elapsed time of $\gamma$ is bounded. The conclusion says that
given these bounds, the elapsed time is strictly less than that of
the corresponding rescaled standard solution.
\end{remark}
\begin{proof}
If $T_x>T_0+\th\hat h^2$ then by Lemma \ref{claim5} and \ref{II.4.5},
\begin{equation}
R(\ga(T_0+\th\hat h^2))\ge \const(1-\th)^{-1}\hat h^{-2}.
\end{equation}
Thus by (\ref{Rbound}) we get
\begin{equation}
Q^{-1}\const (1-\th)^{-1}\hat h^{-2}\leq R(\ga(T_x))\leq Q(T_x-T_0)^{-1},
\end{equation}
or
\begin{equation}
T_x-T_0\leq \const Q^2(1-\th)\hat h^2,
\end{equation}
If we choose $\th$ close enough to $1$ then $\const Q^2(1-\th)\hat h^2$ is less than
$\th\hat h^2$, which gives a contradiction.
\end{proof}

\section{II.5. Statement of the 
the existence theorem for
 Ricci flow with surgery} \label{II.5}

 Our presentation of this material follows Perelman's, except for
some shuffling of the material.  We will be using some terminology
introduced in Section \ref{Ricci flow with surgery}, as well
as results about the $L$-function and noncollapsing from 
Sections \ref{lfunctionandsurgerysection} and \ref{II.5.2section}.

\begin{definition}
A compact Riemannian $3$-manifold is {\em normalized} if
$|\Rm|\leq 1$ everywhere, and the volume of every unit 
ball is at least half the volume of the Euclidean unit
ball.  
\end{definition}

We will use the fact that a smooth normalized
Ricci flow, with bounded curvature on compact time intervals,
satisfies the Hamilton-Ivey pinching condition of
Definition \ref{HamIvey}.

The main result of the surgery procedure is Proposition \ref{surgeryflowexists}
(cf. II.5.1), which
implies that one can choose positive nonincreasing
functions $r:\R_+\ra (0,\infty)$,
$\de:\R_+\ra (0,\infty)$ such that the Ricci 
flow with $(r,\de)$-surgery
flow starting with any normalized initial condition will be defined 
for all time.  

The actual statement is structured to facilitate
a proof by induction:

\begin{proposition}(cf. Proposition II.5.1)
\label{surgeryflowexists}
There exist decreasing sequences $0<r_j<\eps^2$,
$\kappa_j>0$, $0<\bar\de_j<\eps^2$ for $1\leq j<\infty$,
such that for any normalized initial data and any nonincreasing function 
$\de:[0,\infty)\ra(0,\infty)$ such that $\de<\bar\de_j$
on $[2^{j-1}\eps,2^j\eps]$, the Ricci flow with $(r,\de)$-cutoff
is defined for all time and is $\kappa$-noncollapsed
at scales below $\eps$.  
\end{proposition}
\noindent
Here, and in the rest of this section, $r$ and $\kappa$
will always denote functions defined on an interval
$[0,T]\subseteq [0,\infty)$ with the property that 
$r(t)=r_j$ and $\kappa(t)=\kappa_j$
for all $t\in [0,T]\cap [2^{j-1}\eps,2^j\eps)$.
By ``$\kappa$-noncollapsed at scales below $\epsilon$'', we mean
that for each $\rho < \epsilon$ 
and all $(x, t) \in \M$ 
with $t \ge \rho^2$, 
whenever $P(x, t, \rho, - \rho^2)$ is unscathed and
$|\Rm| \: \le \: \rho^{-2}$ on $P(x, t, \rho, - \rho^2)$,
then we also have 
$\vol(B(x,t,\rho)) \: \ge \: \kappa(t) \rho^3$.

Recall that
$\eps$ is a ``global'' parameter which is assumed to be 
small, i.e. all statements involving $\eps$ (explicitly or
otherwise) are true provided $\eps$ is sufficiently small.
Proposition \ref{surgeryflowexists} does not impose any
serious new constraints on $\epsilon$.
For example, instead of using the time intervals
$\{[2^{j-1}\eps,2^j\eps]\}_{j=1}^\infty$,
we could have taken any collection of adjoining time intervals
starting at a small positive time. Also, we just need some
fixed upper bound on $r_j$ and $\overline{\delta}_j$.
We will follow \cite{Perelman2} and write these somewhat
arbitrary constants in terms of the single global parameter $\epsilon$.
Note also that having normalized initial data sets a length
scale for the Ricci flow.

The phrase ``the Ricci flow with $(r,\de)$-cutoff
is defined for all time'' allows for the possibility that the entire
manifold goes extinct, i.e. that after some time we are talking
about the flow on the empty set.

In the rest of this section we give a sketch of the proof.  The details
are in the subsequent sections.

Given positive nonincreasing functions $r$ and $\de$,
if one has a normalized initial condition $(M,g(0))$ then there
will be a maximal time interval on which the Ricci flow with 
$(r,\de)$-cutoff is defined.   This interval can be finite
only if it is of the form $[0,T)$ for some $T<\infty$, and the 
Ricci flow with $(r,\de)$-cutoff on $[0,T)$ extends to a Ricci flow with surgery
on $[0,T]$ for which the $r$-canonical neighborhood
assumption fails at time $T$; see Lemma \ref{extendit}. The main point
here is that the $r$-canonical neighborhood assumption
allows one to run the flow forward up to the singular time, and then perform
surgery, while volume considerations rule out an accumulation of 
surgery times.   Thus the crux of the proof is showing that the 
functions $r$
and $\de$ can be chosen so that the $r$-canonical neighborhood assumption
will continue to hold, and the Ricci flow with surgery 
satisfies a noncollapsing
condition.

The strategy is to  argue
by induction on $i$ that $r_i$, $\bar\de_i$, and $\kappa_i$
can be chosen (and $\bar\de_{i-1}$ can be adjusted)
so that the statement of the proposition holds on the
the finite time interval $[0,2^i\eps]$.  
In the induction step, one establishes the canonical
neighborhood assumption using an argument by contradiction similar
to the proof of Theorem \ref{thmI.12.1}. (We recommend that the reader review this
before proceeding).   The main difference between the proof of
Theorem \ref{thmI.12.1} and that of Proposition \ref{surgeryflowexists} is that the 
non-collapsing assumption, the key ingredient
that allows one to implement the blowup
argument, is no longer available as a direct consequence of Theorem \ref{nolocalcollapse}, 
due to the presence of surgeries.   

We now discuss the augmentations to the non-collapsing argument
of Theorem \ref{nolocalcollapse} necessitated by surgery; this is treated in detail in sections 
\ref{lfunctionandsurgerysection} and \ref{II.5.2section}.
We first recall Theorem \ref{nolocalcollapse} and its proof:
  if a parabolic ball $P(x_0,t_0,r_0,-r_0^2)$ 
in Ricci flow (without surgery) is
sufficiently collapsed then one
uses the $L$-function
with basepoint $(x_0,t_0)$, and the $\L$-exponential map based
at $(x_0,t_0)$, to get a contradiction.   One considers the reduced
volume of a suitably chosen time slice $\M_t$. There is a 
positive lower bound on the reduced volume coming
from the selection of a point where the reduced distance
is at most $\frac{3}{2}$, which in turn comes from an application of the 
maximum principle to the $L$-function. On the other hand, there
is an upper bound on the reduced volume, which the collapsing
forces to be small, thereby giving the contradiction. The upper bound comes
from the monotonicity of the weighted Jacobian of the $\L$-exponential
map.   In fact,  this upper bound works without
significant modification  in the presence of surgery,
provided one considers only the reduced volume contributed
by those points in the time $t$ slice
which may be joined to $(x_0,t_0)$ by minimizing $\L$-geodesics
lying in the regular part of spacetime (see Lemma \ref{localreducedvolume}).

To salvage the lower bound on the reduced volume, the basic idea is that
by making the surgery parameter $\de$ small, one can force the 
$\L$-length of any curve passing close to the surgery locus to 
be large (Lemma \ref{lplusbig}).   This implies that if $(x,t)$ is a point where $L$
isn't too large, then there will  necessarily be an
$\L$-geodesic from $(x_0,t_0)$ to $(x,t)$.  To construct
the minimizer, one takes  a sequence of admissible curves 
from $(x,t)$ to $(x_0,t_0)$ with
$\L$-length tending to the infimum, and argues that they must stay away from the 
surgeries; hence
they remain in a compact part of 
spacetime, and subconverge to a minimizer. 
Therefore the calculations
from Sections \ref{I.7}-\ref{I.7.3}
will be valid near such a point $(x,t)$. The 
maximum principle can then be applied as before 
to show that the minimum of the reduced
length is $\leq \frac{3}{2}$ on each time slice (see Lemma \ref{stillleq3/2}).  

To be more precise, if one makes the surgery parameter $\de(t^\prime)$ small
for a surgery at a given time $t^\prime$ then one can force the $\L$-length of
any curve passing close to the time-$t^\prime$ surgery locus to be large, provided that the
endtime $t_0$ of the curve is not too large compared to $t^\prime$. 
(If $t_0$ is much larger than $t^\prime$ then the curve may spend a long time
in regions of negative scalar curvature after time $t^\prime$. The ensuing negative 
effect on $\L$ could overcome the
positive effect of the small surgery parameter.)  
In the proof of Theorem
\ref{nolocalcollapse}, in order to show noncollapsing at time $t_0$, one went all 
the way back to a time slice near the initial time and found a point there where
$l$ was at most $\frac32$. There would be a problem in using this method
for Ricci flows with surgery - we would have to constantly redefine $\delta(t^\prime)$
to handle the case of larger and larger $t_0$. The resolution is to not go back to
a time slice
near the initial time slice. Instead, in order to
show $\kappa$-noncollapsing in the time slice
$[2^i \eps, 2^{i+1} \eps]$, we will want to get a lower bound on the
reduced volume for
a time $t$-slice with $t$ lying in the preceding time interval
$[2^{i-1}\eps,2^i\eps]$.   As we inductively have control over the geometry
in the time slice $[2^{i-1}\eps,2^i\eps]$, the argument works equally well.

Finally, as mentioned, after obtaining the {\it a priori} $\kappa$-noncollapsing
estimate on the interval $[2^{i}\eps,2^{i+1}\eps]$, one proves that the
$r$-canonical neighborhood assumption holds at time
$T \in [2^{i}\eps,2^{i+1}\eps]$. One difference here is
that because of possible nearby surgeries, there are two ways
to obtain the canonical neighborhood : either from closeness to
a $\kappa$-solution, as in the proof of
Theorem \ref{thmI.12.1}, or from closeness to a standard
solution.

\section{The $L$-function  of I.7 and Ricci flows with surgery}
\label{lfunctionandsurgerysection}

In this section we examine several points which arise when one adapts
the noncollapsing argument of Theorem \ref{nolocalcollapse} to Ricci flows with surgery.
This material is implicit background for 
Lemma \ref{II.5.2} and Proposition \ref{propII.6.3}.  We will use notation 
and terminology introduced in Section \ref{Ricci flow with surgery}.

Let $\M$ be a Ricci flow with surgery, 
and fix a point $(x_0,t_0)\in\M$.   One may
define the $\L$-length of an  admissible curve $\ga$ from 
$(x_0,t_0)$ to some $(x,t)$, for $t<t_0$, using the formula
\begin{equation}
\L(\ga)=\int_t^{t_0}\sqrt{t_0-\bar t}\,\left(R+|\dot{\ga}|^2\right)d\bar t,
\end{equation}
where $\dot{\ga}$ denotes the spatial part of the velocity of $\ga$.
One defines the $L$-function on $\M_{(-\infty, t_0)}$ by setting
$L(x,t)$ to be the infimal $\L$-length of the admissible curves
from $(x_0,t_0)$ to $(x,t)$ if such an admissible curve exists, and 
infinity otherwise. We note that if $(x,t)$ is in a surgery time slice
$\M_t^-$ and is actually removed by the surgery then there will not be an
admissible curve from $(x_0, t_0)$ to $(x,t)$.

If $\ga$ is an admissible curve lying in $\M_{\reg}$ then 
the first variation formula applies.
Hence an admissible curve in $\M_{\reg}$ from $(x_0, t_0)$ to $(x,t)$
whose ${\mathcal L}$-length
equals $L(x,t)$ will satisfy the ${\mathcal L}$-geodesic 
equation.
If $\ga$ is a stable $\L$-geodesic in $\M_{\reg}$ then 
the proof of the monotonicity along $\ga$ of 
the weighted Jacobian
$\tau^{-\frac{3}{2}}\exp(-l(\tau))J(\tau)$ remains valid.
Similarly, if $U\subset \M_{(-\infty,t_0)}$ is an open set such that
every $(x, t)\in U$ is accessible from $(x_0,t_0)$
by a minimizing $\L$-geodesic (i.e. an $\L$-geodesic of $\L$-length
$L(x,t)$) contained in $\M_{\reg}$, 
then the arguments of Section \ref{I.(7.15)} imply that the 
differential inequality 
\begin{equation} \label{newbarrier}
\bar L_\tau+\De \bar L\leq 6
\end{equation}
holds in $U$, in the barrier sense, where $\tau = t_0 - t$,
$\bar L \: = \: 2 \sqrt{\tau} \: L$ and
$\l = \frac{\bar L}{4\tau}$.

\begin{lemma}[Existence of $\L$-minimizers]
\label{minimizersexist}
Let $\M$ be a Ricci flow with surgery
defined on $ [a,b]$. Suppose that $(x_0,t_0)\in \M$ lies in the backward
time slice $\M_{t_0}^-$.

(1) For each $(x,t)\in \M_{[a,t_0)}$ with $L(x,t)<\infty$, there exists
an $\L$-minimizing admissible path $\ga:[t,t_0]\ra\M$ from $(x,t)$ to $(x_0,t_0)$
which satisfies the
$\L$-geodesic equation at every time $\overline{t}\in (t,t_0)$ for which
$\ga(\overline{t})\in\M_{\reg}$.

(2) $L$ is lower semicontinuous on 
$\M_{[a,t_0)}$
and continuous
on 
$\M_{\reg}\cap \M_{[a,t_0)}$. (Note that $\M_a^+ \subset \M_{\reg}$.)

(3) Every sequence $(x_j,t_j)\in \M_{[a,t_0)}$ with $\limsup_j L(x_j,t_j)<\infty$
has  a convergent subsequence.
\end{lemma}
\begin{proof}
(1)   Let 
$\{\ga_j:[t,t_0]\ra\M\}_{j=1}^\infty$ be a sequence of
admissible curves from $(x,t)$ to $(x_0,t_0)$ such that
$\lim_{j\ra\infty}\L(\ga_j)=L(x,t)<\infty$. By restricting the sequence,
we may assume that
$\sup_j \L(\ga_j)<2L(x,t)$.
We claim that there is 
a subsequence of the $\ga_j$'s  that 

(a) converges uniformly to some $\ga_\infty:[t,t_0]\ra\M$,

and

(b) converges weakly to $\ga_\infty$ in $W^{1,2}$ on any  subinterval 
$[t',t'']\subset [t,t_0)$ such that $[t',t'']$ is free of singular times.

\noindent
To see this, note that on any time interval $[c,d]\subset [t,t_0)$
which is free of singular times, one may apply the Schwarz inequality
to the $\L$-length, along with the fact that the metrics on the time
slices $\M_t$, $t \in [c,d]$, are uniformly biLipschitz to each other,
to conclude that the $\ga_j$'s are uniformly
H\"older-continuous on $[c,d]$.  We know
that $\ga_j(t')$ lies in $\M_{t'}^-\cap \M_{t'}^+$ for each surgery
time $t'\in (t,t_0)$, and so one can use similar reasoning to get
H\"older control on a short time interval of the form $[t'',t']$.
Using a change of variable as in  (\ref{changetos}), one obtains
uniform H\"older control near $t_0$ after reparametrizing with $s$.
It follows that the $\ga_j$'s are equicontinuous and map into a
compact part of spacetime,  so 
Arzela-Ascoli applies;  therefore, by passing  to a subsequence we
may assume  that (a) holds. 

To show (b),  we apply weak compactness to the sequence
\begin{equation}
\{\ga_j\restr_{[t',t'']}\};
\end{equation}
\noindent
this is justified by the fact that  the paths $\ga_j\restr_{[t',t'']}$'s remain in a part
of $\M$ with  bounded geometry.  Thus we may assume that our sequence $\{\ga_j\}$
converges uniformly on $[t,t_0]$ and weakly on every subinterval $[t',t'']$ as in (b).
By weak lower semicontinuity of $\L$-length, it follows that the $W^{1,2}$-path
$\ga_\infty$ has $\L$-length $\leq L(x,t)$.  Since any $W^{1,2}$ path  may be approximated
in $W^{1,2}$ by  admissible curves with the same endpoints, it follows that $\ga_\infty$
 minimizes $\L$-length among $W^{1,2}$ paths, and therefore 
it restricts to a smooth solution of the $\L$-geodesic equation on each time
interval $[t',t'']\subset [t,t_0]$ such that $(t',t'')$ is free of singular times.  
Hence $\ga_\infty$ is an $\L$-minimizing admissible curve.

(2) Pick $(x,t)\in\M_{[a,t_0)}$.  To verify lower semicontinuity at 
$(x,t)$ we suppose the sequence $\{(x_j,t_j)\}\subset \M_{[a,t_0)}$ 
converges to $(x,t)$
and $\liminf_{j\ra\infty}L(x_j,t_j)<\infty$.   By (1) there is a sequence 
$\{\ga_j\}$ of $\L$-minimizing admissible curves, where
$\ga_j$ runs from $(x_j,t_j)$ to $(x_0,t_0)$.  By the reasoning above,
a subsequence of $\{\ga_j\}$ converges uniformly and weakly in $W^{1,2}$ to a $W^{1,2}$
curve $\ga_\infty:[t,t_0]\ra\M$ going from $(x,t)$ to $(x_0,t_0)$, with 
\begin{equation}
\L(\ga_\infty)\leq \liminf_{j\ra\infty}\L(\ga_j).
\end{equation}
Therefore $L(x,t)\leq \liminf_{j\ra\infty} L(x_j,t_j)$, and we have established
semicontinuity.  If 
$(x,t)\in \M_{\reg}$,
the opposite inequality obviously
holds, so in this case $(x,t)$ is a point of continuity.

(3) Because $\{L(x_j,t_j)\}$ is uniformly bounded,  any sequence $\{\ga_j\}$ of $\L$-minimizing
paths with $\ga_j(t_j)=(x_j,t_j)$ will be equicontinuous, and hence by Arzela-Ascoli
a subsequence converges uniformly.  Therefore a subsequence of $\{(x_j,t_j)\}$ converges.
\end{proof}

The fact that (\ref{newbarrier}) can hold locally allows one
to appeal -- under appropriate conditions -- to the maximum
principle as in Section \ref{I.(7.15)} to prove
that $\min l\leq \frac{3}{2}$ on 
every time slice. 
Recall that $l \: = \: \frac{L}{2\sqrt{\tau}} \: = \: 
\frac{\overline{L}}{4\tau}$.

\begin{lemma}
\label{stillleq3/2}
Suppose that $\M$ is a Ricci flow with surgery defined on $ [a,b]$.
Take 
$t_0\in (a,b]$
and $(x_0,t_0)\in \M_{t_0}^-$.  Suppose that
for every $t \in [a, t_0)$, every admissible curve
$[t,t_0]\ra\M$ ending at $(x_0,t_0)$ which does not lie
in 
$\M_{\reg}\cup  \M_{t_0}^-$
has reduced length strictly greater than $\frac{3}{2}$. Then
there is a point $(x,a)\in \M_a^+$ where $l(x,a)\leq \frac{3}{2}$.
\end{lemma}
\begin{remark}
In the lemma we consider the Ricci flow with surgery to begin at time $a$.
Hence $\M_{\reg}\cup  \M_{t_0}^- \: = \: \M_a^+ \cup \M_{\reg}\cup  \M_{t_0}^-$
and so the hypothesis of the lemma is a statement about the reduced
lengths of barely admissible curves, in the sense of Section 
\ref{Ricci flow with surgery}.
\end{remark}

\begin{proof}
As in the case when there are no surgeries,
the proof relies on the maximum principle and a continuity argument.

Let $\be:\M_{[a,t_0)}\ra \R\cup \{\infty\}$ be the function
\begin{equation}
\be \: = \: \bar L-6\tau=4\tau\left(l-\frac{3}{2}\right),
\end{equation}
where as usual, $\tau(x,t) \: = \: t_0 - t$. 
Note that for each $\tau\in (0,t_0-a]$, the function
$\be$ attains a minimum $\be_{\min}(\tau)<\infty$ on the 
slice $\M_{t_0 - \tau}$,  because by (2) of Lemma \ref{minimizersexist},
it is continuous on the compact manifold 
$\M_{t_0 - \tau}^+$ (as seen by changing the parameter $a$ of 
Lemma \ref{minimizersexist} to $t_0 - \tau$),
and $\be\equiv\infty$ on 
$\M_{t_0 - \tau}\;-\;\M_{t_0 - \tau}^+$.
Thus it suffices to show that $\be_{\min}(t_0-a)\leq 0$.

From Lemma \ref{neg}, $\beta_{\min}(\tau)<0$  for $\tau>0$ small.
Let $\tau_1\in (0,t_0-a]$ be the supremum of the 
$\bar \tau\in(0,t_0-a]$ such that $\be_{\min}< 0$ on the interval $(0,\bar\tau)$.
 
We claim that 

(a) $\be_{\min}$ is continuous on $(0,\tau_1)$

and

(b) The upper right $\tau$-derivative of $\be_{\min}$ is nonpositive
on $(0,\tau_1)$.

To see (a), pick $\tau\in (0,\tau_1)$, suppose 
that $\{\tau_j\}\subset (0,\tau_1)$ is a sequence
converging to $\tau$ and choose $(x_j,t_0 - \tau_j)\in \M_{t_0 - \tau_j}^+$
such that $\be(x_j,t_0 - \tau_j)=\be_{\min}(\tau_j)< 0$.  
By Lemma \ref{minimizersexist}
part (3), the sequence $\{(x_j,t_0 - \tau_j)\}$ subconverges to some
$(x,t_0 - \tau)\in\M_{t_0 - \tau}^+$
for which $\be(x,t_0 - \tau)\leq \liminf_{j\ra\infty}\be(x_j,t_0 - \tau_j)$.
Thus 
$\be_{\min}$ is lower semicontinuous at $\tau$.  On the other hand,
since $\be_{\min}(\tau)< 0$, the minimum of $\be$ on $\M_{t_0 - \tau}$
will be attained at a point
$(x,t_0 - \tau)\in \M_{t_0 - \tau}^+$ lying 
in the interior of $\M_{t_0 - \tau}^-\cap \M_{t_0 - \tau}^+$,
as $\be>0$ elsewhere on $\M_{t_0 - \tau}$ 
(by Lemma \ref{minimizersexist} and the
hypothesis on admissible curves).  Therefore $\be$ is continuous at 
$(x,{t_0 - \tau})$,
which implies that $\be_{\min}$ is upper semicontinuous at $\tau$. 
This gives (a).

Part (b) of the claim follows from the fact that if $\tau\in (0,\tau_1)$ and
the minimum of $\be$ on $\M_{t_0 - \tau}$ is attained at
$(x,t_0 - \tau)$ then $l(x,\tau)< \frac32$, so there is a neighborhood $U$
of $(x,{t_0 - \tau})$ such that the inequality 
\begin{equation}
\frac{\D \be}{\D \tau}+\De \be \leq 0
\end{equation} 
holds in the barrier sense on $U$ (by Lemma \ref{minimizersexist} and the hypothesis
on admissible curves).
Hence the upper right derivative $\frac{d}{ds} \Big|_{s=0} \beta(x, \tau+s)$ is
nonpositive, so the upper right $\tau$-derivative of 
$\beta_{\min}(\tau)$ is also nonpositive.

The claim implies that $\be_{\min}$ is  nonincreasing on $(0,\tau_1)$,
and so 
$\limsup_{\tau\ra \tau_1^-}\be(\tau)<0$.
By parts (2) and (3) of Lemma \ref{minimizersexist},
we have $\be_{\min}(\tau_1)<0$, and the minimum is attained at some  
$(x,{t_0 - \tau_1})\in
\M_{\reg}$.
(Recall that $\M_a \subset \M_{\reg}$.)
This implies that $\tau_1=t_0-a$, for otherwise $\be_{\min}(\tau)$ 
would be strictly negative for $\tau\geq\tau_1 $ close to $\tau_1$, 
contradicting the definition of $\tau_1$.
\end{proof}

The notion of local collapsing can be adapted to Ricci flows with surgery,
as follows.
\begin{definition}
\label{surgerycollapsedef}
Let $\M$ be a Ricci flow with surgery defined on $[a,b]$. Suppose that
$(x_0,t_0)\in\M$ and $r>0$ are such that $t_0-r^2\geq a$, 
$B(x_0,t_0,r)\subset\M_{t_0}^-$ is a proper ball
and the parabolic ball $P(x_0,t_0,r,-r^2)$ is unscathed.
Then {\em $\M$ is $\kappa$-collapsed at $(x_0,t_0)$ at scale
$r$} if $|\Rm|\leq r^{-2}$ on $P(x_0,t_0,r,-r^2)$ and 
$\vol(B(x_0,t_0,r))<\kappa r^3$; otherwise it is $\kappa$-noncollapsed.
\end{definition}

We make use of the following variant of the noncollapsing
argument from Section \ref{I.7.3}.
 
\begin{lemma}(Local version of reduced volume comparison)
\label{localreducedvolume}
There is a function $\kappa':\R_+\ra \R_+$, satisfying
$\lim_{\kappa\ra 0}\kappa'(\kappa)=0$, with the following property.
Let $\M$ be a Ricci flow with surgery defined on $[a,b]$. Suppose
that we are given $t_0\in (a,b]$,
$(x_0,t_0)\in \M_{t_0}\cap \M_{\reg}$, $t\in [a,t_0)$ and
$r \in (0,\sqrt{t_0-t})$.  Let $Y$ be the set of points $(x,t)\in \M_t$
that are accessible from $(x_0,t_0)$ by means of minimizing $\L$-geodesics
which remain in $\M_{\reg}$.  Assume 
in addition that 
$\M$ is $\kappa$-collapsed at $(x_0,t_0)$ at scale $r$, i.e.
$P(x_0,t_0,r,-r^2) \cap \M_{[t_0 - r^2, t_0)} \subset\M_{\reg}$,
$|\Rm|\leq r^{-2}$ on 
$P(x_0,t_0,r,-r^2)$, and $\vol(B(x_0,t_0,r))<\kappa r^3$. 
Then  the reduced volume of $Y$ is at most 
$\kappa'(\kappa)$.
\end{lemma}
\begin{proof}
Let $\hat Y\subset T_{x_0}\M_{t_0}$ be the set of vectors 
$v\in T_{x_0}\M_{t_0}$
such that there is a minimizing $\L$-geodesic $\ga:[t,t_0]\ra \M_{\reg}$
 running from $(x_0,t_0)$ to some point in $Y$, with
\begin{equation}
\lim_{\bar t\ra t_0}\sqrt{t_0-\bar t}\;\dot{\ga}(\bar t)=-v.
\end{equation}

\noindent
The calculations from Sections \ref{Lremarks}-\ref{monoV} 
apply to $\L$-geodesics sitting in $\M_{\reg}$.
In particular, the monotonicity of the weighted Jacobian
$\tau^{-\frac{n}{2}}\exp(-l(\tau))J(\tau)$ holds.  Now one repeats
 the proof of Theorem \ref{nolocalcollapse}, working with the set $\hat Y$
instead of  the set of initial velocities of {\em all} minimizing
$\L$-geodesics. 
\end{proof}

\section{Establishing noncollapsing in the presence of surgery}
\label{II.5.2section}

The key result of this section, Lemma \ref{II.5.2},
gives conditions under which
one can deduce noncollapsing on a time interval $I_2$, given 
a noncollapsing bound on a preceding interval $I_1$ and
lower bounds on $r$ on $I_1\cup I_2$.

\begin{definition}
The {\em ${\mathcal L}_+$-length} of an admissible curve $\gamma$ is
\begin{equation}
{\mathcal L}_+(\gamma, \tau) \: = \:
\int_{t_0 - \tau}^{t_0} \sqrt{t_0 - t} \: \left(
R_+(\gamma(t), t) \: + \: |\dot{\gamma}(t)|^2 \right) \: dt,
\end{equation}
where $R_+(x,t) = \max(R(x,t), 0)$.
\end{definition}

\begin{lemma}(Forcing $\L_+$ to be large, cf. Lemma II.5.3)
\label{lplusbig}

For all $\La< \infty$, $\bar r>0$ and $\hat r>0$,
there is  a constant  $F_0=F_0(\La,\bar r,\hat r)$
with the following property.  Suppose that

$\bullet$  $\M$ is a Ricci flow  with $(r,\de)$-cutoff
defined on an interval containing $[t,t_0]$, where 
 $r([t,t_0])\subset [\hat r,\eps]$,

$\bullet$ $r_0\geq \bar r$, 
 $B(x_0,t_0,r_0)$ is a proper ball which is unscathed on $[t_0-r_0^2,t_0]$, and
$|\Rm|\leq r_0^{-2}$ on $P(x_0,t_0,r_0,-r_0^2)$,

$\bullet$ $\ga:[t,t_0]\ra\M$ is an admissible curve
ending at $(x_0, t_0)$
whose image 
is not contained in $\M_{\reg}\cup\M_{t_0}$, and

$\bullet$ $\de < F_0(\La,\bar r,\hat r)$
on $[t,t_0]$.  

\noindent
Then $\L_+(\ga)>\La$. 
\end{lemma}
\begin{proof}
The idea is that the hypotheses on  $\ga$  imply that it must touch the
part of the manifold added during surgery at some time $\bar t\in[t,t_0]$.  
Then either $\ga$ has to move very fast at times close to $\bar t$ or
$t_0$, or it
will stay in the surgery region while it develops large
scalar curvature. In the first case $\L_+(\ga)$ will be large because
of the $|\dot{\ga}|^2$ term in the formula for $\L_+$, 
and in the second case it will be large because of 
the $R(\ga)$ term.

First, we can assume that $F_0$ is small enough so that
$F_0 \: < \: \sqrt{\frac{\bar r}{100 \epsilon}}$.
Then since
\begin{equation}
\label{h<<r}
\max_{[t,t_0]}h(t)\leq \left(\max_{[t,t_0]}\de\right)^2
\left(\max_{[t,t_0]}r(t)\right)\leq F_0^2\eps,
\end{equation}
we have $\max_{[t,t_0]}h(t)
<\frac{\bar r}{100}\leq  \frac{r_0}{100}$.

Put 
$\De t=10^{-10} \bar r^4\La^{-2}$.
It always suffices to prove 
the lemma for a larger value of $\La$,
so without loss of generality we can assume that
$\De t\leq \bar r^2\leq r_0^2$. Set 
\begin{equation}
\label{aandthetaconditions}
A=A((\De t)^{-\frac12}\La,\hat r),\quad 
\th=\th((\De t)^{-\frac12}\La,\hat r)
\end{equation}
where $A(\cdot,\cdot)$ and $\th(\cdot,\cdot)$
are  the functions from Corollary \ref{corollaryII.4.6}.
That is, we will eventually be applying Corollary \ref{corollaryII.4.6}
with $l \: = \: (\De t)^{-\frac12}\La$.
We impose the additional constraint on $F_0$ that
\begin{equation}
\label{deltasmall}
 F_0\leq \hat\de(A+2,\th,\hat r)
\end{equation}
on the interval $[t, t_0]$,
where $\hat\de $ is the function from Lemma \ref{II.4.5}.

As $\gamma$ is admissible but is not contained in $\M_{\reg}\cup\M_{t_0}$,
it must pass through the boundary of a surgery cap at some time
in the interval $[t, t_0)$ or it must start in the interior of a
surgery cap at time $t$.
By dropping an initial segment of $\ga$
if necessary, we may assume that $\ga(t)$ lies in a
surgery cap. 

Let $x$ denote the tip of the surgery cap.
Note that 
\begin{equation}
\label{p'sdisjoint}
P(x_0,t_0,r_0,-r_0^2)\cap P(x,t,Ah(t),\th h^2(t))=\emptyset
\end{equation}
since by Lemma \ref{II.4.5} the scalar curvature on 
$P(x,t,Ah(t),\th h^2(t))$ is at least  
$\frac{h^{-2}}{2}>\frac{10^4}{2}r_0^{-2}$, while $|\Rm|\leq r_0^{-2}$
on $P(x_0,t_0,r_0,-r_0^2)$.  Therefore when going backward in time
from $(x_0,t_0)$, $\ga$ must leave the parabolic region
$P(x_0,t_0,r_0,-r_0^2)$
before it arrives at $(x,t)$.
If it exits at a time $\tilde t > t_0-\De t$ then 
applying the Schwarz inequality we get
\begin{equation}
\int_{\tilde{t}}^{t_0} \sqrt{t_0 - s} \: 
|\dot{\gamma}(s)|^2 \: ds \: \ge \:
\left( \int_{\tilde{t}}^{t_0} 
|\dot{\gamma}(s)| \: ds \right)^2 
\left( \int_{\tilde{t}}^{t_0} (t_0 - s)^{-1/2}
\: ds \:\right)^{-1} 
\ge \: \frac{1}{100} \: r_0^2 (\Delta t)^{-1/2}>\La,
\end{equation}
where the factor of $\frac{1}{100}$ comes from 
the length distortion estimate
of Section \ref{secI.8.3},
using the fact that $|\Rm| \le r_0^{-2}$ on $P(x_0, t_0, r_0, -r_0^2)$.
So we can restrict to the case when $\ga$ exits
$P(x_0,t_0,r_0,-\De t)$ through the
initial time slice  at time $t_0-\De t$. In particular,
by (\ref{p'sdisjoint}),
$\gamma$ must exit the parabolic region $P(x,t,Ah(t),\th h^2(t))$
by time $t_0-\De t$.

By Lemma \ref{II.4.5}, the parabolic region $P(x,t,Ah,\theta h^2)$ is 
either unscathed, or it coincides (as a set) with the parabolic
region $P(x,t,Ah,s)$ for some $s\in (0,\th h^2)$ and the entire 
final time slice $P(x,t,Ah,s)\cap \M_{t+s}$
of $P(x,t,Ah,s)$ is thrown away by a surgery 
at time $t+s$. 

One possibility is that  $\gamma$ exits $P(x,t,Ah,\theta h^2)$ 
through the final time slice. If this is the case then $P(x,t,Ah,\theta h^2)$ 
must be unscathed (as otherwise the final face is removed by surgery at 
time $t+s < t+\theta h^2$ and $\gamma$ would have nowhere to go after this time), so
$\gamma$ lies in $P(x,t,Ah,\theta h^2)$ for the entire time interval
$[t, t+\theta h^2]$.

The other possibility is that $\ga$ leaves $P(x,t,Ah,\theta h^2)$ before the
final time slice of $P(x,t,Ah, \theta h^2)$, in which case it exits
the ball $B(x,t,Ah)$ by time $t + \theta h^2$.

Corollary \ref{corollaryII.4.6} applies to either of these
two possibilities. Putting
\begin{equation}
T_\ga= \sup\{\bar t\in [t,t+\theta h^2]\mid \ga([t,\bar t])\subset
P(x,t,Ah,\theta h^2)\}
\end{equation}
and using the fact that $T_\ga \le t_0 - \De t$,
we have
\begin{align}
& \int_{t}^{t_0} \sqrt{t_0 - s} \left( R_+(\gamma(s), s) \: + \:
|\dot{\gamma}(s)|^2 \right) \: ds \: \ge \:
\int_{t}^{T_\gamma} 
\sqrt{t_0 - s} \left( R_+(\gamma(s), s) \: + \:
|\dot{\gamma}(s)|^2 \right) \: ds \: \ge \\
& (\Delta t)^{1/2} \int_{t}^{T_\gamma} 
\left( R_+(\gamma(s), s) \: + \:
|\dot{\gamma}(s)|^2 \right) \: ds \: \ge \:
(\Delta t)^{1/2} \: l \: = \: \La, \notag
\end{align} 
where the last inequality comes from Corollary \ref{corollaryII.4.6} 
and the choice
of $A,\th$, and $\de$ in (\ref{aandthetaconditions}) and (\ref{deltasmall}).
This completes the proof.
\end{proof}

\begin{lemma} \label{Rev}
If $\M$ is a Ricci flow with surgery, with normalized initial
condition at time zero, then for all $t \ge 0$, $R(x,t) \: \ge \:
- \: \frac{3}{2} \: \frac{1}{t+\frac14}$.
\end{lemma}
\begin{proof}
From the initial conditions, $R_{min}(0) \: \ge \: -6$.
If the Ricci flow is smooth then 
(\ref{Rlowerbound}) implies
that $R_{min}(t) \: \ge \: - \: \frac{3}{2} \: \frac{1}{t+\frac14}$.
If there is a surgery
at time $t_0$ then $R_{min}$ on $\M_{t_0}^+$ equals
$R_{min}$ on $\M_{t_0}^-$, as surgery is done in
regions of high scalar curvature.  The lemma follows by applying
(\ref{Rlowerbound}) on the time intervals between
the singular times.
\end{proof}

In the statement of the next lemma, one has successive
time intervals $[a,b)$ and $[b,c)$. As a mnemonic we
use the subscript $-$ for quantities attached to the
earlier interval $[a,b)$, and $+$ for those associated with $[b,c)$.
We will also assume that the global parameter $\epsilon$ is small
enough that the $\Phi$-pinching condition 
implies that whenever $|\Rm(x,t)|\geq \eps^{-2}$, then 
$R(x,t)>\frac{|\Rm(x,t)|}{100}$.
(We  remind the reader of the role of the parameter
$\epsilon$; see Remark \ref{remepsilon}.) 

\begin{lemma}(Noncollapsing estimate) (cf. Lemma II.5.2)
\label{II.5.2}

Suppose $\epsilon \ge r_- \ge r_+ >0$, $\kappa_->0$, $E_->0$ and $E<\infty$.
Then there are constants $\overline{\delta} =
\overline{\delta}(r_-,r_+,\kappa_-,E_-,E)$
and $\kappa_+=\kappa_+(r_-,\kappa_-,E_-,E)$ with
the following property.  Suppose that

$\bullet$  $a<b<c,\quad b-a\geq E_-,\quad c-a\leq E$,

$\bullet$  $\M$ is a Ricci flow with $(r,\de)$-cutoff
with normalized initial condition
defined on a time interval containing $[a,c)$, 

$\bullet$ $r\geq r_-$ on $[a,b)$ and $r\geq r_+$ on $[b,c)$, 

$\bullet$ $r \le \epsilon$ ,

$\bullet$ $\M$ is $\kappa_-$-noncollapsed at scales below
$\eps$ on $[a,b)$ and

$\bullet$ $\de\leq \bar\de$ on $[a,c)$,

Then $\M$ is $\kappa_+$-noncollapsed at scales below $\eps$
on $[b,c)$.
\end{lemma}
\begin{remark}
The important point to notice here is that $\bar\de$ is allowed
to depend on the lower bound $r_+$ on $[b,c)$, but the 
noncollapsing constant $\kappa_+$ does not depend on $r_+$.
\end{remark}
\begin{proof}
In the proof, we can assume that $\frac{r_+}{100} \le \sqrt{E_-/3}$.
If this were not the case then we could prove the lemma with
$r_+$ replaced by $100 \sqrt{E_-/3}$. Then the lemma would 
also hold for the original value of $r_+$.

First, from Lemma \ref{Rev}, $R \ge -6$ on
$\M_{[a,c)}$.

Suppose that $r_0\in (0,\eps)$,  $(x_0,t_0)\in\M_{[b,c)}$,
  $B(x_0,t_0,r_0)$
is a proper ball  unscathed on the interval $[t_0-r_0^2,t_0]$, and
$|\Rm|\leq r_0^{-2}$ on $P(x_0,t_0,r_0,-r_0^2)$.

We first assume that 
$r_0\leq \sqrt{E_-/3}$
and $r_0\geq \frac{r_+}{100}$.

We will consider $\L$-length, $\L_+$-length, etc in 
$\M_{[a,t_0)}$ with basepoint at $(x_0,t_0)$. Suppose that 
$\widehat{t} \in [a,t_0)$. 
Then for any admissible curve 
$\ga:[\widehat{t},t_0]\ra \M_{[a,t_0]}$
ending at $(x_0,t_0)$, we have

\begin{gather}
\label{llplus}
\L(\ga)\leq \L_+(\ga)\leq \L(\ga)+\int_a^c6\sqrt{c-t}\: dt
\leq \L(\ga)+4E^{\frac32} \\
\mbox{and}\quad l(x,t)\geq\frac{L_+-4E^{\frac32}}{2E^{\frac12}}. \notag
\end{gather}

Assume that $\bar\de\leq F_0(4E^{\frac12}+4E^{\frac32},\frac{r_+}{100},r_+)$
where $F_0$ is the function from Lemma \ref{lplusbig}.
Then by (\ref{llplus}) and Lemma \ref{lplusbig},
we conclude that any admissible
curve 
$[\widehat{t},t_0]\ra \M_{[a,t_0]}$
ending at $(x_0,t_0)$
which does not lie in $\M_{\reg}\cup \M_{t_0}$
has reduced length bounded below by $2=\frac32+\frac12$.  By Lemma 
\ref{stillleq3/2} there is an admissible curve $\ga:[a,t_0]\ra\M$
ending at $(x_0,t_0)$ such that
\begin{equation}
\L(\ga)=L(\ga(a))=2\sqrt{t_0-a}\;l(\ga(a))\leq 3\sqrt{t_0-a},
\end{equation} 
so by (\ref{llplus}) it follows that
\begin{equation}
\label{lplusupper}
\L_+(\ga)\leq 3\sqrt{t_0-a}+4E^{\frac32}\leq 3\sqrt{E}+4E^{\frac32}.
\end{equation}

Set 
\begin{equation}
t_1= a+\frac{b-a}{3},\quad t_2= a+\frac{2(b-a)}{3}
\end{equation}
and
\begin{equation}
\rho= \left(3\sqrt{E}+4E^{\frac32}\right)
\left(\frac13 E_-\right)^{-\frac32}.
\end{equation}
By construction, $t_2 \le t_0 - r_0^2$.
Note that there is a $\bar t\in [t_1,t_2]$ such that
$R(\ga(\bar t))\leq \rho$.
Otherwise we would get
\begin{equation}
\L_+(\ga) > \int_{t_1}^{t_2}\sqrt{t_0-t}\;R_+(\ga(t)) \: dt
\geq \sqrt{\frac13 E_-}\int_{t_1}^{t_2}\rho \: dt
\geq \sqrt{\frac13 E_-}\left(\frac13 E_-\right)\rho=3\sqrt{E}+4E^{\frac32},
\end{equation}
contradicting (\ref{lplusupper}).

Put $\overline{x} = \gamma(\overline{t})$.
By Lemma \ref{claim1II.4.2}, there is an estimate of the form
\begin{equation}
\label{Rbound2}
R\leq \const s^{-2}
\end{equation}
on the parabolic ball 
$\hat P = P(\overline{x}, \bar t,s,-s^2)$ with 
$s^{-2} \: = \: \const \: (\rho \: + \: r_-^{-2})$.
Appealing to Hamilton-Ivey curvature pinching as 
usual, we get that $|\Rm|\leq \const s^{-2}$ in $\hat P$.
If $\overline{t} - s^2 < a$ then 
we shrink $s$ (as little as possible)
to ensure that $\hat P\subset\M_{[a,b)}$.
Provided that  $\bar \de$ is less than a small constant 
$c_1=c_1(r_-,E_-,E)$, 
we can guarantee that $\hat P$ is unscathed, by forcing the
curvature in a surgery cap to exceed our bound (\ref{Rbound2})
on $R$.  Put
$U= \M_{\bar t-\frac12 s^2}\cap \hat P$.  
If $s < \epsilon$ then
the $\kappa_-$-noncollapsing assumption on $[a,b)$  gives a lower bound
on $\vol(B(\bar x,\bar t,s))s^{-3}$. If $s \ge \epsilon$ then
the $\kappa_-$-noncollapsing assumption gives a lower bound
on $\vol(B(\bar x,\bar t,\epsilon/2)) (\epsilon/2)^{-3}$.
In either, case, we get a lower bound on $\vol(B(\bar x,\bar t,s))$ and
hence a lower bound
$\vol(U)\geq v=v(r_-,\kappa_-,E_-,E)$.
Now every point
in $U$ can be joined to $(x_0,t_0)$ by  a curve of 
$\L_+$-length at most $\La_+=\La_+(r_-,E_-,E)$,
by concatenating an admissible curve 
$[\bar t-\frac12 s^2,\bar t]\ra \M$ (of controlled
$\L_+$-length) with $\ga\restr_{[\bar t,t_0]}$.
Shrinking $\bar\de$ again, we can apply 
Lemmas \ref{minimizersexist} and \ref{lplusbig} with
(\ref{llplus}) to ensure that every point in 
$U$ can be joined to $(x_0,t_0)$ by a minimizing
$\L$-geodesic lying in $\M_{\reg}\cup \M_{t_0}$.  
Lemma \ref{localreducedvolume}
then implies that 
\begin{equation}
\label{kappa1}
\vol(B(x_0,t_0,r_0))r_0^{-3}\geq \kappa_1=\kappa_1(r_-,\kappa_-,E_-,E).
\end{equation}

(We briefly recall the argument. We have a parabolic ball around
$(\bar x, \bar t)$, of small but controlled size, on which we have
uniform curvature bounds.  The lower volume bound coming from
the $\kappa_-$-noncollapsing assumption
on $[a,b)$ means that we have bounded geometry on the parabolic ball.
As we have a fixed upper bound on $l(\bar x, \bar t)$, we can
estimate from below the reduced volume of the accessible points
$Y \subset \M_{\bar t-\frac12 s^2}^+$. Then we obtain a lower
bound on $\vol(B(x_0,t_0,r_0))r_0^{-3}$ as in Theorem 
\ref{nolocalcollapse}.)

This completes the proof of the lemma when 
$r_0\leq \sqrt{E_-/3}$
and $r_0\geq \frac{r_+}{100}$.

Now suppose that $r_0>\sqrt{E_-/3}$.  Applying our 
noncollapsing estimate (\ref{kappa1}) to the ball of radius
$\sqrt{E_-/3}$ gives 
\begin{equation}
\vol(B(x_0,t_0,r_0))r_0^{-3}\geq \left(\vol(B(x_0,t_0,\sqrt{E_-/3}) (E_-/3)^{-\frac32}
\right)\frac{(E_-/3)^{\frac32}}{r_0^3}
\geq \kappa_1\frac{(E_-/3)^{\frac32}}{\eps^3}=\kappa_2,
\end{equation}
where $\kappa_2=\kappa_2(r_-,\kappa_-,E_-,E)$.

The next sublemma deals with the case when $r_0< \frac{r_+}{100}$.
\begin{sublemma}
\label{r0/100}
If $r_0< \frac{r_+}{100}$ then $\vol(B(x_0,t_0,r_0))r_0^{-3}\geq
\kappa_3=\kappa_3(r_-,\kappa_-,E_-,E)$.
\end{sublemma}
\begin{proof}
Let $s$ be the maximum of the numbers $\bar s\in [r_0,\frac{r_+}{100}]$
such that $B(x_0,t_0,\bar s)$ is unscathed on $[t_0-\bar s^2,t_0]$,
and $|\Rm|\leq \bar s^{-2}$ on $P(x_0,t_0,\bar s,-\bar s^{-2})$.
Then  either

(a) Some point $(x,t)$ on the frontier of $P(x_0,t_0,s,-s^2)$
lies in a surgery cap

or

(b) Some point $(x,t)$ in the closure of $P(x_0,t_0,s,-s^2)$
has $|\Rm|=s^{-2}$

or

(c) $s=\frac{r_+}{100}$.

In case (a) the scalar curvature
at $(x,t)$ will satisfy $R(x,t)\in (\frac{h^{-2}}{2},10h^{-2})$, 
since $(x,t)$ lies in the cap at the surgery time.  Since
$|\Rm(x,t)|\leq s^{-2}$ we conclude that $s \le \const h(t)$.
If $\bar\de$ is small then the pointed time slice
$(\M_t,(x,t))$ will be close, modulo scaling by $h(t)$, to
the initial condition of the standard solution with
the basepoint somewhere in the cap.  Using the
fact that the time slices of $P(x_0,t_0,s,-s^2)$ have  
comparable metrics, along with the fact that 
$r_0 \le s  \le \const h(t)$, we get a lower bound 
$\vol(B(x_0,t_0,r_0))r_0^{-3}\geq \mbox{const}$.

In case (b), we have a static curve $\ga:[t,t_0]\ra\M$
such that $\ga(t)=(x,t)$, $\ga(t_0)\in \ol{B(x_0,t_0,s)}$,
and $|\Rm(x,t)|=s^{-2}\geq 10^4r_+^{-2}$. 
Hence by $\Phi$-pinching,
$R(x,t)\geq \frac{1}{100} s^{-2} \geq 100 r_+^{-2}$
(cf. the remark just before the statement of
Lemma \ref{II.5.2}).
With reference to the constant $\eta$ of
(\ref{II.(1.3)estimates}), 
put $\sigma \: = \: 10^{-6} \: \min \left( 1, 
\frac{1}{\eta} \right)$.
Let 
$\al:[0,\widehat{\rho}]\ra \ol{B(x_0,t_0,s)}\subset \M_{t_0}^-$
be a minimizing geodesic from $(x_0,t_0)$ to 
$\ga(t_0)\in \ol{B(x_0,t_0,s)}$. We can find a point
$z$ along $\alpha$ with $\dist_{t_0}(z, \gamma(t_0)) \le \sigma s$ and
$\dist_{t_0}(z, x_0) \le (1 - \sigma) s$.
Let $\bar\ga:[t,t_0]\ra\M$
be the static curve ending at $z$
and put $(\bar x,t)= \bar\ga(t)$.  In brief,
we get $(\bar x,t)$ by ``pulling $(x,t)$ slightly
inward from the boundary''.

From the distance distortion estimate of Section \ref{secI.8.3},
we can say that $\dist_{t} (\bar x, x) \le 10^6 \sigma s$.
Then applying (\ref{II.(1.3)estimates}) along a minimizing
time-$t$ curve from $\bar x$ to $x$, we conclude that
\begin{equation}
|R^{- \: \frac12}(\bar x, t) \: - \: R^{- \: \frac12}(x, t)| \: \le \:
\frac12 \: \eta \: \dist_{t} (\bar x, x) \: \le \: \frac12 \: s,
\end{equation}
so $R^{- \: \frac12}(\bar x, t) \: \le \:  
R^{- \: \frac12}(x, t) \: + \: \frac12 \: s \: \le 20s$ and hence
$R(\bar x, t) \: \ge \: \frac{1}{400} \: s^{-2} \ge \: 25 \: r_+^{-2}$.
In particular, $(\bar x,t)$ has a canonical neighborhood.
We also know that $R(\bar x, t) \le 6 |\Rm (\bar x, t)| \le 6 s^{-2}$.
It follows from the definition of canonical neighborhoods
(see Definition \ref{canonicalnbhddef}) that
there is some universal constant so that
$\vol(B(\bar x,t,10^{-9}s) ) \ge \const s^3$.
(We recall that from Lemma \ref{lemstandard0},
there is a $\kappa>0$ such that a standard solution
is $\kappa$-noncollapsed as scales $<1$.)
The distance distortion
estimate ensures that $B(\bar x,t,10^{-9}s) \subset
B(z,t_0,10^{-6}s) \subset B(x_0,t_0,s)$. Then the standard volume
distortion estimate implies that 
$\vol(B(x_0,t_0,s)) \: \ge \: \const 
\vol(B(\bar x,t,10^{-9}s) ) \ge \const s^3$,
again for some universal constant.
Finally we use Bishop-Gromov volume comparison to get 
$\vol(B(x_0,t_0,r_0))r_0^{-3}\geq \const \vol(B(x_0,t_0,s))s^{-3}$.
 
In case (c) we apply (\ref{kappa1}), replacing the
$r_0$ parameter there by $s$, and Bishop-Gromov
volume comparison as in case (b). 
\end{proof}

\section{Construction of the Ricci flow with surgery}  \label{surgeryflow}

The proof is by induction on $i$.
To start the induction process, 
we observe that 
the initial normalization
$|\Rm|\leq 1$ at $t=0$
implies that a smooth solution exists for some definite time
\cite[Corollary 7.7]{Chow-Knopf}.
The curvature bound on this time interval, along with the
volume assumption on the initial time balls,
implies that the solution is $\kappa$-noncollapsed below  scale $1$
and satisfies the  $\rho$-canonical neighborhood assumption 
vacuously for small $\rho>0$.

Now assume inductively that $r_j$, $\kappa_j$, and $\bar\de_j$
have been selected for $1\leq j\leq i$, thereby defining
the functions $r$, $\kappa$, and $\bar\delta$ 
on $[0, 2^i\eps]$, such that for 
any nonincreasing function $\de$ on $[0,2^i\eps]$ satisfying 
$0<\de(t)\leq \bar\de(t)$, if one has normalized
initial data then the Ricci flow with $(r,\de)$-cutoff
is defined on $[0,2^i\eps]$ and is $\kappa(t)$-noncollapsed at
scales $<\eps$.

We will determine $\kappa_{i+1}$ using Lemma \ref{II.5.2}.
So suppose it is not possible to choose $r_{i+1}$
and $\bar\de_{i+1}$ (and, if necessary, make $\bar\de_i$ smaller)
so that if we put
$\kappa_{i+1}= \kappa_+(r_i,\kappa_i,2^{i-1}\eps,
3(2^{i-1})\eps)$ (where $\kappa_+$ denotes the function from
Lemma \ref{II.5.2})
then the inductive statement above holds with $i$ replaced by $i+1$.
Then given sequences $r^\al\ra 0$ and $\bar\de^\al\ra 0$, for each
$\al$ there
must be a counterexample, say $(M^{\al},g^{\al}(0))$,
to the statement with $r_{i+1} = r^{\al}$ and
$\bar\de_i = \bar\de_{i+1} = \bar\de^{\al}$. We assume that  
\begin{equation}
  \bar\de^\al < \hat\de(\al,1-\frac{1}{\al},r^\al)
\end{equation}
where $\hat\de$ is the quantity from Lemma \ref{II.4.5}; this
will guarantee that for any $A<\infty$ and $\th\in (0,1)$
we may apply Lemma \ref{II.4.5} with parameters $A$ and $\th$
for sufficiently large $\al$.  We also assume that 
\begin{equation}
\label{F_1condition}
\bar\de^{\al}<\bar\de(r_i,r^\al,\kappa_i,2^{i-1}\eps,3(2^{i-1})\eps)
\end{equation}
where $\bar\de(r_i,r^\al,\kappa_i,2^{i-1}\eps,3(2^{i-1})\eps)$ is from Lemma \ref{II.5.2}.
 By Lemma \ref{extendit}  each initial condition $(M^{\al},g^{\al}(0))$
will prolong to a Ricci flow
with surgery   $\M^{\al}$ defined on a time interval $[0,T^{\al}]$
with 
$T^{\al}\in (2^i\eps,\infty]$,
which restricts to a Ricci flow with $(r,\de)$-cutoff
on any proper subinterval  $[0,\tau]$ of $[0,T^{\al}]$,
 but for which the
$r^\al$-canonical neighborhood assumption fails at some point
$(\bar x^{\al},T^{\al})$ lying in the backward time
slice 
$\M^{\al-}_{T^{\al}}$.
(It is implicit in this statement that $R(\bar x^{\al},T^{\al}) \ge 
\frac{1}{(r^\al)^2}$.)
Since $(M^{\al},g^{\al}(0))$ violates the
theorem, we must have $T^\al\in (2^i\eps,2^{i+1}\eps]$. By (\ref{F_1condition})
and Lemma \ref{II.5.2}, it follows that  $\M^{\al}$ is 
$\kappa_{i+1}$-noncollapsed at scales below $\eps$
on the interval $[2^i \epsilon,T^\al)$, where 
$\kappa_{i+1}= \kappa_+(r_i,\kappa_i,2^{i-1}\eps,
3(2^{i-1})\eps)$ and $\kappa_+$ denotes the function from
Lemma \ref{II.5.2}.  

Let $(\widehat\M^\al,(\bar x^\al,0))$ be the pointed Ricci flow
with surgery obtained from $(\M^\al,(\bar x^\al,T^\al))$
by shifting time by $T^\al$ and parabolically rescaling by
$R(\bar x^\al,T^\al)$.
We also remove the part of
$(\widehat\M^\al,(\bar x^\al,0))$ after time zero and
we take the time-zero slice $\widehat\M^\al_0$ to be
diffeomorphic to $\M^{\al -}_{T^\al}$.
In brief, the rest of the proof goes as follows.  If surgeries
occur further and further away from $(\bar x^\al,0)$
in spacetime as $\al\ra\infty$, then the reasoning of Theorem
\ref{thmI.12.1} applies and we obtain a $\kappa$-solution as a limit.
This would contradict the fact that $(\bar x^\al,T^\al)$ does
not have a canonical neighborhood.  Thus there must be 
surgeries in a parabolic ball of a fixed size centered
at $(\bar x^\al,0)$, for arbitrarily large $\al$. Then one
argues using 
Lemma \ref{II.4.5} that the solution will be close to
the (suitably rescaled and time-shifted) standard solution,
which again leads to a canonical neighborhood and a 
contradiction.

We now return to the proof.
Recall that a metric ball $B$ is
proper if the distance function from the center is a proper
function on $B$. If $T$ is a surgery time for a
Ricci flow with surgery then
a metric ball in $\M_T^-$ need not be proper.

Note that by continuity, 
every  point in $\widehat\M^\al$ whose scalar
curvature is strictly greater than that of $(\bar x^\al,0)$ has a
neighborhood as in Definition \ref{canonicalnbhddef}, except
that the error estimate
is $2 \epsilon$ instead of $\eps$.

\begin{sublemma}
For all $\la<\infty$,
the ball $B(\bar x^\al,0,\la)\subset\widehat\M^\al_0$
is proper for sufficiently
large $\al$.
\end{sublemma}
\begin{proof}
As in Lemma \ref{claim2II.4.2}, for each $\rho < \infty$ 
the scalar curvature on
$B(\bar x^\al, 0, \rho) \subset \widehat\M^{\al}_0$ 
is uniformly bounded in terms of $\al$.
(In carrrying out the proof of Lemma \ref{claim2II.4.2},
we now use the aforementioned property of having
canonical neighborhoods of quality $2\epsilon$.)
Thus $B(\bar x^\al, 0, \rho)$ has compact closure in
$\widehat\M^{\al}_0$
(see Lemma \ref{lemoverlinerisproper}), from which the sublemma follows.
\end{proof}

Let $T_1 \in [-\infty,0]$ be the infimum of the set of numbers
$\tau'\in (-\infty,0]$ such that for all $\la<\infty$,
the ball $B(\bar x^\al,0,\la)\subset\widehat\M^\al_0$
is proper,
and unscathed on $[\tau',0]$ for sufficiently
large $\al$.

\begin{lemma} \label{bdrysmoothness}
After passing to a subsequence
if necessary, the pointed flows  $(\widehat\M^\al,(\bar x^\al,0))$
converge on the time interval $(T_1, 0]$ to a Ricci flow (without 
surgery) $(\M^\infty,(\bar x^\infty,0))$ with 
a smooth complete nonnegatively-curved Riemannian metric
on each time slice, and 
scalar curvature globally bounded above by 
some number $Q<\infty$. 
(We interpret $(0, 0]$ to mean $\{0\}$ rather than the 
empty set.) 
\end{lemma}
\begin{proof}
Suppose first that $T_1 < 0$.
Then the arguments of Theorem \ref{thmI.12.1} apply in the 
time interval $(T_1, 0]$, to give the Ricci flow (without 
surgery) $(\M^\infty,(\bar x^\infty,0))$. 
Since $r^{\al} \rightarrow 0$, Hamilton-Ivey pinching
implies that $\M^\infty$ will have nonnegative curvature.
The fact that the 
canonical neighborhood assumption, 
with $\epsilon$ replaced by $2 \epsilon$,
holds for each $\widehat\M^\al$ allows us
to deduce
that the scalar curvature of $\M^\infty$ is globally bounded above by 
some number $Q<\infty$; compare with Section \ref{pf11.7}
and Step 4 of the proof of Theorem \ref{thmI.12.1}.

Now suppose that $T_1 = 0$.
The argument is similar to Steps 2 and 3 of the proof of Theorem
\ref{thmI.12.1}. 
As in Step 2,
or more precisely as in  Lemma \ref{claim2II.4.2}, 
for each $\rho < \infty$ 
the scalar curvature on
$B(\bar x^\al, 0, \rho) \subset \widehat\M^{\al}_0$ 
is uniformly bounded in terms of $\al$.
Given $\rho < \infty$ and 
$(x^\al,0) \in B(\bar x^\al, 0, \rho) \subset \widehat\M^\al$,
if the parabolic region of
Lemma \ref{claim1II.4.2} (centered around
$(x^\al, 0) \in \widehat\M^\al$) is unscathed then we can
apply 
the $2 \epsilon$-canonical neighborhood assumption
on $\widehat\M^\al$, Lemma \ref{claim1II.4.2} and Appendix \ref{applocalder}
to derive bounds on the curvature derivatives at 
$(x^\al,0) \in \widehat\M^\al_0$
that depend on $\rho$ but are independent of $\al$.
If the parabolic region is scathed then we can apply Lemma
\ref{II.4.5}, along with our scalar curvature bound at 
$(x^\al,0)$,
to again obtain uniform bounds on the curvature derivatives at 
$(x^\al,0)$.
Hence there is a subsequence of
the pointed Riemannian manifolds 
$\{(\widehat\M^\al_0, \bar x^\al)\}_{\al = 1}^\infty$ 
that converges to a smooth complete pointed
Riemannian manifold $(\M^\infty_0, \bar x^\infty)$. 
As in the previous case, it will have 
bounded nonnegative sectional curvature.

This proves the lemma.  Alternatively, in the case $T_1 = 0$
one can argue directly that if the parabolic region of
Lemma \ref{claim1II.4.2} is scathed then 
$(\overline{x}^\al, T^\al)$ has a canonical neighborhood;
see the rest of the proof of Proposition \ref{surgeryflowexists}.
\end{proof}

If we can show that
$T_1 \: = \: -\infty$ then $(\M^\infty,(\bar x^\infty,0))$
will be a $\kappa$-solution, which will contradict the assumption
that $(\bar x^\al,T^\al)$ does not admit a canonical neighborhood.
Suppose that $T_1 \: > \: -\infty$.
We know that for all 
$\tau' \in (T_1, 0]$
and $\la<\infty$, the scalar curvature in $P(\bar x^\al,0,\la,\tau')$
is bounded by $Q+1$ when $\al$ is sufficiently large.  By 
Lemma \ref{claim1II.4.2}, there exists $\sigma< T_1$ such that for all
$\la<\infty$, if (for large $\al$) 
the solution $\widehat\M^\al$ is unscathed on
$P(\bar x^\al,0,\la,t^\al)$ for some $t^\al > \sigma$ then
\begin{equation}
\label{2q+1}
R(x,t)< 8(Q+2)\quad\mbox{for all}\quad (x,t)\in P(\bar x^\al,0,\la,t^\al).
\end{equation}
(In applying Lemma \ref{claim1II.4.2}, we use the fact that in the unscaled
variables,
$r(T^\al)^{-2} \le R(\bar x^\al,T^\al)$ by assumption, along with
the fact that $r(\cdot)$ is a nonincreasing function.)

By the definition of $T_1$, and after passing
to a subsequence if necessary, there exist $\la<\infty$ and a sequence
$\ga^\al:[\si^\al,0]\ra\widehat\M^\al$  of static curves so that  \\
1. $\ga^\al(0)\in B(\bar x^\al,0,\la)$ and \\
2. The point
$\ga^\al(\si^\al)$ is inserted during surgery at time $\si^\al>\si$.

For each $\al$, we may assume that $\si^\al$ is the largest number having
this property.
Put
$\xi^\al= \si^\al+\left( h^\al(\si^\al)\right)^2$.
(In the notation of \cite{Perelman2}, 
$\left( h^\al(\si^\al)\right)^2$ would be written
as $R(\bar x, \bar t) \: h^2(T_0)$. Note that we have no {\it a priori} control
on $h^\al(\si^\al)$.) Then $\xi^\al$
is the blowup time of the rescaled and shifted
standard solution that Lemma \ref{II.4.5}
compares with $(\widehat\M^\al,\ga^\al(\si^\al))$.  We claim
that $\liminf_{\al\ra\infty}\xi^\al>0$. Otherwise, Lemma 
\ref{II.4.5} would imply that after passing to a subsequence,
there are regions of
$\widehat\M^\al$, starting from time $\si^\al$, that are better and better
approximated by rescaled and shifted standard solutions whose blowup times
$\xi^\al$ have a limit that is nonpositive,
thereby contradicting (\ref{2q+1}). Lemma \ref{II.4.5}, along with
the fact that $R({\bar x}^\al, 0) = 1$, also gives a uniform upper bound
on $\xi^\al$.

Now Lemma \ref{II.4.5} implies that for large $\al$,
the restriction of 
$\widehat\M^\al$ to the time interval $[\si^\al,0]$ is well
approximated by the restriction to $[\si^\al,0]$
of a rescaled and shifted standard
solution. Then Lemma \ref{claim5} implies that $(\bar x^\al,T^\al)$
has a canonical neighborhood.
The canonical neighborhood may be either a strong $\epsilon$-neck
or an $\epsilon$-cap. (Note 
a strong $\epsilon$-neck may arise when an 
$\epsilon$-neck around
$(\bar x^\al,T^\al)$ extends smoothly backward in time to form
a strong $\epsilon$-neck that incorporates part of the Ricci flow 
solution that existed before the surgery time $\si^\al$.)

This is a contradiction.
\end{proof}

\section{II.6. Double sided curvature bound in the thick part}
\label{II.6}

Having shown that for a suitable choice of the functions $r$ and $\de$,
the Ricci flow with $(r,\de)$-cutoff exists for all time and for every 
normalized initial condition, one wants to understand
its implications. 
The main results in II.6 are noncollapsing and curvature estimates
which form the basis of the analysis of the large-time behavior given 
in II.7.

\begin{lemma} \label{pscempty}
If $\M$ is a Ricci flow with $(r,\de)$-cutoff on a compact manifold and
$g(0)$ has positive scalar curvature then the solution goes extinct
after a finite time, i.e. $\M_T \: = \: \emptyset$ for some $T > 0$.
\end{lemma}
\begin{proof}
We apply (\ref{Rlowerbound}). This
formula is initially derived for smooth flows but because surgeries are
performed in regions of high scalar curvature, it is also valid for a Ricci 
flow with surgery; cf. the proof of Lemma \ref{Rev}. 
It follows that the flow goes extinct by time
$\frac{3}{2R_{min}(0)}$.
\end{proof}

\begin{lemma} \label{finitetime}
If $\M$ is a Ricci flow with surgery that goes extinct after a finite time,
then the initial (compact connected orientable) $3$-manifold is diffeomorphic to a 
connected sum of $S^1 \times S^2$'s and quotients of the round
$S^3$.
\end{lemma}
\begin{proof}
This follows from Lemma  \ref{reconstruct2}.
\end{proof}

According to \cite{Colding-Minicozzi,Colding-Minicozzi2} and 
\cite{Perelman3}, if 
none of the prime factors in the Kneser-Milnor decomposition of
the initial manifold are aspherical
then the Ricci flow with surgery again goes extinct
after a finite time.
Along with Lemma \ref{finitetime}, this proves the
Poincar\'e Conjecture.

Passing to Ricci flow solutions that may not go extinct after a finite time,
the main result of II.6 is the following :

\begin{corollary} \label{II.6.8}
(cf. Corollary II.6.8) For any $w > 0$ one can find
$\tau = \tau(w) > 0$, $K = K(w) < \infty$,
$\overline{r} = \overline{r}(w) > 0$ and $\theta = \theta(w) > 0$ with
the following property. Suppose we have a solution to the Ricci flow
with $(r, \delta)$-cutoff on the time interval $[0, t_0]$, with
normalized initial data.  
Let $h_{\max}(t_0)$ be the maximal surgery radius on $[t_0/2, t_0]$.
(If there are no surgeries on $[t_0/2, t_0]$ then $h_{\max}(t_0) = 0$.) 
Let $r_0$ satisfy \\
1. $\theta^{-1}(w) h_{\max}(t_0) \: \le \: r_0 \: \le \: 
\overline{r} \sqrt{t_0}$.\\
2. The ball $B(x_0, t_0, r_0)$ has sectional curvatures at least
$- \: r_0^{-2}$ at each point.\\
3. $\vol(B(x_0, t_0, r_0)) \: \ge \: w r_0^3$.\\

Then the solution is unscathed in $P(x_0, t_0, r_0/4, - \tau r_0^2)$
and satisfies $R \: < \: Kr_0^{-2}$ there.
\end{corollary}

Corollary \ref{II.6.8} is an analog of Corollary \ref{corI.12.4}, but there
are some differences.  One minor difference is that Corollary \ref{corI.12.4} is
stated as the contrapositive of Corollary \ref{II.6.8}. Namely,
Corollary \ref{corI.12.4} assumes that $- \: r_0^{-2}$ is achieved as a sectional
curvature in $B(x_0, t_0, r_0)$, and its conclusion is that
$\vol(B(x_0, t_0, r_0)) \: \le \: w r_0^3$. The relation
with Corollary \ref{II.6.8} is the following. Suppose that 
assumptions 1 and 2 of Corollary \ref{II.6.8} hold. If
$- \: r_0^{-2}$ is achieved somewhere as a sectional
curvature in $B(x_0, t_0, r_0)$ then Hamilton-Ivey pinching implies
that the scalar curvature is very large at that point, 
which contradicts the conclusion
of Corollary \ref{II.6.8}. Hence assumption 3 of Corollary \ref{II.6.8}
must not be satisfied.

A more substantial difference is that the smoothness of the flow
in Corollary \ref{corI.12.4}
is guaranteed by the setup, whereas in Corollary \ref{II.6.8} we must
prove that the solution is unscathed in $P(x_0, t_0, r_0/4, - \tau r_0^2)$.

The role of the parameter $\overline{r}$ in Corollary \ref{II.6.8} is
essentially to guarantee that we can use Hamilton-Ivey pinching
effectively.

\section{II.6.5. Earlier scalar curvature bounds on smaller balls from
lower curvature bounds and a later volume bound} \label{II.6.5}

For terminology, with reference to Section \ref{Ricci flow with surgery},
by a time-dependent family $\bigcup_{t \in [c,d]} B(x, t, r)$
of metric balls we mean first that there is a static
curve $\gamma : [c,d] \rightarrow \M$, whose intersection with each
$\M_t$ will be denoted by $x$, and second that there is
a subset ${\mathcal U}$ of $\M$ so that
\begin{enumerate}
\item If $t \notin [c,d]$ then ${\mathcal U} \cap \M_t = \emptyset$.
\item If $t \in [c,d]$ is not
a singularity time then ${\mathcal U} \cap \M_t$ is the $r$-ball around $x$ in
$\M_t$.
\item If $t \in [c,d]$ is a surgery time $t_k^+$ then 
${\mathcal U} \cap \M_t$ is
the image in $\M_t$ of an $r$-ball around $x$ in
$\Omega_k$, that lies entirely in $X_k^+ \subset \Omega_k$.
\end{enumerate}
At the expense of being redundant, if these conditions are
satisfied then we will say that 
we have an {\em unscathed} time-dependent family of metric balls,
to emphasize that the metric balls do not touch the surgery regions.
The restriction of the Ricci-flow-with-surgery to ${\mathcal U}$ 
is a smooth Riemannian metric $g$ on the ``horizontal'' subbundle
of $T{\mathcal U}$.

We first state a consequence of Corollary \ref{I.11.6end}.

\begin{lemma} \label{I.11.6(a)}
Given $w > 0$, there exist $\tau_0 = \tau_0(w) > 0$ and 
$K_0^\prime = K_0^\prime(w) < \infty$
with the following property. Suppose that we have
a Ricci flow with $(r,\de)$-cutoff such that \\
1. $\bigcup_{t \in 
[-\tau r_0^2, 0]} B(x_0, t, r_0)$ is an unscathed time-dependent family of
metric balls, where $\tau \: \le \: \tau_0$. \\
2. The sectional curvatures are bounded below by $- \: r_0^{-2}$ 
on the
above family of balls. \\
3. $\vol(B(x_0, 0, r_0)) \: \ge \: w r_0^3$. \\

Then
\begin{enumerate}
\item 
$R \: \le \: K_0^\prime \: \tau^{-1} \: r_0^{-2}$
on $\bigcup_{ t \in [- \tau r_0^2/2, 0]} B(x_0, t, r_0/2)$, and
\item $\vol(B(x_0, -\tau r_0^2, r_0/2))$ is at least $\frac{1}{10}$
of the volume of the Euclidean ball of the same radius.
\end{enumerate}
\end{lemma}

Here we have changed the conclusion of Corollary \ref{I.11.6end} to obtain an
upper curvature bound on 
$\bigcup_{ t \in [- \tau r_0^2/2, 0]} B(x_0, t, r_0/2)$
instead of $\bigcup_{ t \in [- \tau r_0^2/2, 0]} B(x_0, t, r_0/4)$,
but this clearly follows from the arguments of the proof of Corollary
\ref{I.11.6end}. We have also added a lower volume bound to the conclusion,
which follows from the proof of Corollary \ref{corI.11.6}(b), provided
that $\tau_0$ is sufficiently small.

An analog of Corollary \ref{II.6.8} is the following
Lemma \ref{II.6.5(a)}, which is stated as Lemma II.6.5(a)
in \cite{Perelman2}. The lemma is used there to prove Corollary \ref{II.6.8}.
Our proof of Corollary \ref{II.6.8} will use Lemma \ref{I.11.6(a)} but will
not use Lemma \ref{II.6.5(a)}.
We include the proof of Lemma \ref{II.6.5(a)} for completeness, even
though it will not be used in the sequel.

\begin{lemma} \label{II.6.5(a)}
(cf. Lemma II.6.5(a))
Given $w > 0$, there exist $\tau_0 = \tau_0(w) > 0$ and $K_0 = 
K_0(w) < \infty$ with the following property. Suppose that we have
a Ricci flow with $(r,\de)$-cutoff such that \\
1. The parabolic neighborhood
$P(x_0, 0, r_0, -\tau r_0^2)$ is unscathed,
where $\tau \: \le \: \tau_0$. \\
2. The sectional curvatures are bounded below by $- \: r_0^{-2}$ on 
$P(x_0, 0, r_0, -\tau r_0^2)$. \\
3. $\vol(B(x_0, 0, r_0)) \: \ge \: wr_0^3$. \\

Then $R \: \le \: K_0 \: \tau^{-1} \: r_0^{-2}$ on 
$P(x_0, 0, r_0/4, -\tau r_0^2/2)$.
\end{lemma}
\begin{proof}
If $\tau_0$ is sufficiently small, then for
$t \in [- \tau r_0^2, 0]$ and 
$(x,t) \in B(x_0, t, 9r_0/10)$, the lower
curvature bound 
$\Rm \: \ge \: - \: r_0^{-2}$ on $P(x_0, 0, r_0, -
\tau r_0^2)$ implies that $(x,0) \in B(x_0, 0, r_0)$ (more precisely,
that $(x,t)$ lies on a static curve with one endpoint in $B(x_0,0,r_0)$,
or equivalently, that $(x,t)\in P(x_0,0,r_0,-\tau r_0^2)$).
Thus $\bigcup_{ t \in 
[-\tau r_0^2, 0]} B(x_0, t, 9r_0/10) 
\subset P(x_0, 0, r_0, - \tau r_0^2)$ and so
$\Rm \: \ge \:-r_0^{-2}\:\ge\: - (9r_0/10)^{-2}$
 on 
$\bigcup_{t \in 
[-\tau (9r_0/10)^2, 0]} B(x_0, t, 9r_0/10)$.

Applying 
Lemma
\ref{I.11.6(a)} with $r_0$ replaced by 
$9r_0/10$, and slightly redefining $w$, gives that 
$R \: \le \: K_0^\prime \: \tau^{-1} \: (9r_0/10)^{-2}$
on $\bigcup_{t \in [- \frac34 \tau (9r_0/10)^2, 0]} B(x_0, t, 9r_0/20)$.
Then the length distortion estimate of Lemma \ref{distancedist} implies that for
sufficiently small $\tau_0$, if 
$(x,0) \in B(x_0, 0, r_0/4)$ 
then 
$(x,t) \in B(x_0, t, 9r_0/20)$ 
for
$t \in [-\tau r_0^2/2, 0]$. That is,
\begin{equation}
P(x_0, 0, r_0/4, -\tau r_0^2/2) \subset 
\bigcup_{t \in [- \tau r_0^2/2, 0]} B(x_0, t, 9r_0/20).
\end{equation}
In applying the length distortion estimate we use the fact that
the change in distance is estimated by 
$\Delta d \: \le \:
\const \: \sqrt{K_0^\prime \: \tau^{-1} \: (9r_0/10)^{-2}} \: \cdot \:
\tau r_0^2/2$
which, for small $\tau_0$, is a small fraction
of $r_0$. 

Thus we have shown that 
$R \: \le \: K_0^\prime \: \tau^{-1} \: (9r_0/10)^{-2}$
on $P(x_0, 0, r_0/4, -\tau r_0^2/2)$. 
This proves the lemma.  
\end{proof}

The formulation of \cite[Lemma II.6.5]{Perelman2} specializes Lemma
\ref{II.6.5(a)} to the case
$w \: = \: 1 - \epsilon$. It includes the statement
\cite[Lemma II.6.5(b)]{Perelman2} saying that $\vol(B(x_0, -\tau
r_0^2, r_0/4))$ is at
least $\frac{1}{10}$ of the volume of the Euclidean ball of the same
radius.  This follows from the proof of Corollary \ref{corI.11.6}(b), provided
that $\tau_0$ is sufficiently small.

There is an evident analogy between Lemma \ref{II.6.5(a)} and
Corollary \ref{II.6.8}. However, there is the important difference that
Corollary \ref{II.6.8} (along with Corollary \ref{corI.12.4}) only
assumes a lower sectional curvature bound at the final time slice.

\section{II.6.6. Locating small balls whose subballs have almost
Euclidean volume}

The result of this section is a technical lemma about volumes
of subballs.

\begin{lemma} \label{Alexlemma2} (cf. Lemma II.6.6)
For any $\widehat{\epsilon}, w>0$ 
there exists $\theta_0 = \theta_0(\widehat{\epsilon},w)$
such that if $B(x,1)$ is a metric ball of volume at least $w$, compactly
contained in a manifold without boundary with sectional
curvatures at least $-1$, then there exists a subball
$B(y,\theta_0) \subset B(x,1)$ such that every subball
$B(z,r) \subset B(y,\theta_0)$ of any radius has volume at
least $(1-\widehat{\epsilon})$ times the volume of the Euclidean ball
of the same radius.
\end{lemma}
The proof is similar to that of Lemma \ref{Alexlemma}.
Suppose that the claim is not true.
Then there is a sequence of Riemannian
manifolds $\{M_i\}_{i=1}^\infty$ and balls $B(x_i, 1)
\subset M_i$ with compact closure
so that $\Rm \Big|_{B(x_i, 1)} \: \ge \: - \: 1$ and
$\vol(B(x_i, 1)) \: \ge \: w$, along with
a sequence $r^\prime_i \rightarrow 0$ so
that each subball
$B(x_i^\prime, r_i^\prime) \subset B(x_i, 1)$ has a subball
$B(x_i^{\prime \prime}, r^{\prime \prime}_i) \subset 
B(x_i^\prime, r_i^\prime)$ 
with $\vol(B(x_i^{\prime \prime}, r_i^{\prime \prime})) \: < \: 
(1 - \widehat{\epsilon}) \: \omega_3 
(r^{\prime \prime})^n$. After taking a subsequence, we can assume that
$\lim_{i \rightarrow \infty} (B(x_i, 1), x_i) \: = \: (X, x_\infty)$ in the
pointed Gromov-Hausdorff topology, where $(X, x_\infty)$ is a pointed
Alexandrov space with curvature bounded below by $-1$. From
\cite[Theorem 10.8]{Burago-Gromov-Perelman}, the Riemannian volume
forms $\dvol_{M_i}$ converge weakly to the three-dimensional Hausdorff
measure $\mu$ of $X$. If $x_\infty^\prime$ is a regular point of
$X$ then there is some $\delta > 0$ so that 
$B(x_\infty^\prime, \delta)$ has compact closure in $X$ and
for all $r < \delta$,
$\mu(B(x_\infty^\prime,
r)) \: \ge \: (1-\frac{\widehat{\epsilon}}{10}) \: \omega_3 r^3$. Fixing 
such an $r$ for the moment,
for large $i$ there are balls
$B(x_i^\prime, r) \subset B(x_i, 1)$
with
$\vol(B(x_i^\prime, r)) \: \ge \: (1-\frac{\widehat{\epsilon}}{5}) \: \omega_3 
r^3$.
Recalling the sequence $\{r^\prime_i\}$,
by hypothesis there is a subball
$B(x_i^{\prime \prime}, 
r_i^{\prime \prime}) \subset B(x_i^\prime, r_i^\prime)$ with
$\vol(B(x_i^{\prime \prime}, r_i^{\prime \prime})) \: < \: 
(1 - \widehat{\epsilon}) \: \omega_3 
(r_i^{\prime \prime})^3$. 
Clearly
$B(x_i^\prime,
r) \: \subset \: B(x_i^{\prime \prime}, r + r_i^\prime)$.
From the Bishop-Gromov inequality,
\begin{equation}
\frac{\vol(B(x_i^{\prime \prime}, r + r_i^\prime))}{\vol(B(x_i^{\prime \prime},
r_i^{\prime \prime}))} \: \le \:
\frac{\int_0^{r + r_i^\prime} \sinh^2(s) \: 
ds}{\int_0^{r_i^{\prime \prime}} \sinh^2(s) \: ds}.
\end{equation}
Then
\begin{equation}
\vol(B(x_i^\prime,
r)) \: \le \: 
\vol(B(x_i^{\prime \prime}, r + r_i^\prime)) \: \le \: 
(1-{\widehat{\epsilon}}) \: \omega_3 \: 
\frac{(r_i^{\prime \prime})^3}{\int_0^{r_i^{\prime \prime}} \sinh^2(s) \: ds}
\: \int_0^{r + r_i^\prime} \sinh^2(s) \: 
ds.
\end{equation}
For large $i$ we obtain
\begin{equation}
\vol(B(x_i^\prime,
r)) \: \le \:
(1-\frac{\widehat{\epsilon}}{2}) \: \omega_3 \: \cdot \: 3
\: \int_0^{r} \sinh^2(s) \: 
ds.
\end{equation}
Then if we choose $r$ to be 
sufficiently small, we contradict the fact that
$\vol(B(x_i^\prime, r)) \: \ge \: (1-\frac{\widehat{\epsilon}}{5}) \: \omega_3 
r^3$ for all $i$.

\begin{remark}
By similar reasoning, for every $L>1$ one may find 
$\th_1=\th_1(\widehat{\epsilon},L)$ such that
under the hypotheses of Lemma \ref{Alexlemma2}, there is a subball
$B(y,\th_1)\subset B(x,1)$ which is $L$-biLipschitz to the Euclidean
unit ball.
\end{remark}

\section{II.6.8. Proof of the double sided curvature bound in the
thick part, modulo two propositions}

In this section we explain how Corollary \ref{II.6.8}
follows from Lemma \ref{Alexlemma2} and two other
propositions, which will be proved in subsequent sections. 
We first state the other propositions,
which are Propositions \ref{propII.6.3} and \ref{propII.6.4}.

\begin{proposition} \label{propII.6.3}
(cf. Proposition II.6.3) 
For any 
$A < \infty$
one can find positive constants
$\kappa(A)$, $K_1(A)$, $K_2(A)$, $\overline{r}(A)$,
such that for any $t_0 < \infty$
there exists $\overline{\delta}_{A}(t_0)
> 0$, decreasing in $t_0$, 
with the following property.  Suppose that
we have a Ricci flow with $(r, \delta)$-cutoff on a time interval
$[0, T]$, where $\delta(t) < \overline{\delta}_A(t)$ on $[t_0/2, t_0]$, with
normalized initial data.  Assume that \\
1. The solution is unscathed on a parabolic ball
$P(x_0, t_0, r_0, -r_0^2)$, with $2r_0^2 < t_0$. \\
2. $|\Rm| \le \frac{1}{3r_0^2}$ on $P(x_0, t_0, r_0, -r_0^2)$. \\
3. $\vol(B(x_0, t_0, r_0)) \: \ge \: A^{-1} r_0^3$.\\

Then\\
(a) The solution is $\kappa$-noncollapsed on scales less than $r_0$
in $B(x_0, t_0, Ar_0)$.\\
(b) Every point $x \in B(x_0, t_0, Ar_0)$ with $R(x, t_0) \ge K_1 r_0^{-2}$
has a canonical neighborhood in the sense of
Definition \ref{canonicalnbhddef}.\\
(c) If $r_0 \: \le \: \overline{r} \sqrt{t_0}$ then $R \: \le \: K_2
r_0^{-2}$ in $B(x_0, t_0, Ar_0)$.
\end{proposition}

Proposition \ref{propII.6.3}(a) is an analog of Theorem 
\ref{8.2thm}.

(The reason for the ``$3$'' in the hypothesis 
$|\Rm| \le \frac{1}{3r_0^2}$ comes from Remark \ref{8.2fix}.)
Proposition \ref{propII.6.3}(c) is an analog of Theorem
\ref{thmI.12.2}, but the hypotheses are slightly different.
In Proposition \ref{propII.6.3} one assumes a lower bound on the volume of the
time-$t_0$ ball $B(x_0, t_0, r_0)$, while in Theorem \ref{thmI.12.2}
one assumes assumes a lower bound on the volume of
the time-$(t_0 - r_0^2)$ ball $B(x_0, t_0 - r_0^2, r_0)$. In
view of the curvature assumption on $P(x_0, t_0, r_0, - r_0^2)$,
the hypotheses are 
essentially
equivalent.

Conclusions (a), (b) and (c) of Proposition \ref{propII.6.3}
are similar to the conclusions of Theorem \ref{8.2thm},
Lemma \ref{12.2lemma} and Theorem \ref{thmI.12.2}, respectively.
Conclusions (a) and (b) of Proposition \ref{propII.6.3}
are also related
to what was proved in Proposition \ref{surgeryflowexists} to construct the 
Ricci flow with surgery. The difference
is that the noncollapsing and canonical neighborhood results of Proposition 
\ref{surgeryflowexists}
are statements at or below the scale $r(t)$, whereas Proposition \ref{propII.6.3}
is a statement about much larger scales, comparable to $\sqrt{t_0}$.
We note that the parameter $\overline{\delta}_A$ in
Proposition \ref{propII.6.3} is independent of the function $\delta$
used to define the Ricci flow with $(r, \delta)$-cutoff.

In the proof of the next proposition we will apply Lemma \ref{Alexlemma2} with
$\widehat{\epsilon}$ equal to the global parameter $\epsilon$,
so we will write $\theta(w)$ instead of $\theta(\epsilon, w)$.

\begin{proposition} \label{propII.6.4}
(cf. Proposition II.6.4)
There exist $\tau, \overline{r}, C_1 > 0$ and $K < \infty$ with the
following property. Suppose that we have a Ricci flow with $(r, \delta)$-cutoff
on the time interval $[0, t_0]$, with normalized initial data.
  Let $r_0$ satisfy 
$2C_1 h_{\max}(t_0) \: \le \: r_0 \: \le \: \overline{r} \sqrt{t_0}$,
where $h_{\max}(t_0)$ is the maximal cutoff radius for surgeries in 
$[t_0/2, t_0]$. (If there are no surgeries on $[t_0/2, t_0]$ then 
$h_{\max}(t_0) = 0$.)

Assume \\
1. The ball $B(x_0, t_0, r_0)$ has sectional curvatures at least $- r_0^{-2}$
at each point. \\
2. The volume of any subball $B(x, t_0, r) \subset B(x_0, t_0, r_0)$ with
any radius $r > 0$ is at least $(1-\epsilon)$ times the volume of the
Euclidean ball of the same radius.\\

Then the solution is unscathed on
$P(x_0, t_0, r_0/4, - \tau r_0^2)$
and satisfies $R \: < \: K r_0^{-2}$ there.
\end{proposition}

Proposition \ref{propII.6.4} is an analog of Theorem \ref{thmI.12.3}. However, 
there is the important difference that in Proposition \ref{propII.6.4}
we have to prove that no surgeries occur within
$P(x_0, t_0, r_0/4, - \tau r_0^2)$.

Assuming the validity of Propositions \ref{propII.6.3} and \ref{propII.6.4},
suppose that the hypotheses of Corollary \ref{II.6.8} are satisfied.
We will allow ourselves to shrink the parameter $\overline{r}$ in order
to apply Hamilton-Ivey pinching when needed.
Put $r_0^\prime \: = \: \theta_0(w) \: r_0$, where $\theta_0(w)$ is
from Lemma \ref{Alexlemma2}.
By Lemma \ref{Alexlemma2}, there is a subball $B(x_0^\prime, t_0, r_0^\prime) \subset
B(x_0, t_0, r_0)$ such that every subball of
$B(x_0^\prime,t_0, r_0^\prime)$ has volume at least
$(1-\epsilon)$ times the volume of the Euclidean ball of the same
radius.  As the sectional curvatures are bounded below by $- r_0^{-2}$ on
$B(x_0, t_0, r_0)$, they are bounded below by $- (r_0^\prime)^{-2}$ on
$B(x^\prime_0, t_0, r_0^\prime)$. By an appropriate choice of the
parameters $\theta(w)$ and $\overline{r}$ of Corollary \ref{II.6.8},
in particular taking $\theta(w) \: \le \: \frac{\theta_0(w)}{2C_1}$, we can
ensure that Proposition \ref{propII.6.4} 
applies to $B(x^\prime_0, t_0, r_0^\prime)$.
Then the solution is unscathed on
$P(x_0^\prime, t_0, r_0^\prime/4, - \tau (r_0^\prime)^2)$ and satisfies
$|\Rm| \: \le \: K \: (r_0^\prime)^{-2}$ there, where the lower
bound on $\Rm$ comes from Hamilton-Ivey pinching. 
With $\tau$ being the parameter of Proposition \ref{propII.6.4} and
putting
$r_0^{\prime \prime} \: = \: \min(K^{-1/2}, \tau^{1/2}, \frac14) \: r_0^\prime$,
for all $t_0^{\prime \prime} \in 
[t_0 - (r_0^{\prime \prime})^2, t_0]$ the
solution is unscathed on $P(x_0^\prime, t_0^{\prime \prime}, 
r_0^{\prime \prime}, -(r_0^{\prime \prime})^2)$ and satisfies
$|\Rm| \: \le \: (r_0^{\prime \prime})^{-2}$ there.
From the curvature bound $\Rm \: \ge \: - \: (r_0^\prime)^{-2}$ on 
$P(x_0^\prime, t_0, r_0^\prime/4, - \tau (r_0^\prime)^2)$
(coming from pinching)
and the fact that $B(x^\prime_0, t_0, r_0^{\prime \prime})$
has almost Euclidean volume, we obtain a bound 
$\vol(B(x^\prime_0, t_0^{\prime \prime}, 
r_0^{\prime \prime})) \: \ge \: \const \:
(r_0^{\prime \prime})^{3}$. 
Applying Proposition \ref{propII.6.3} with
$A = \frac{100r_0}{r_0^{\prime \prime}}$ gives 
$R \: \le \: K_2 (r_0^{\prime \prime})^{-2}$ on
$B(x_0, t_0^{\prime \prime}, 10 r_0) \subset
B(x_0^\prime, t_0^{\prime \prime}, 100 r_0)$,
for all $t_0^{\prime \prime} \in 
[t_0 - (r_0^{\prime \prime})^2, t_0]$.
Writing this as
$R \: \le \: \const r_0^{-2}$,
if we further restrict $\theta(w)$ to be sufficiently
small then we can ensure that $R \: \le \: \const \: \theta^2(w) \: h^{-2}
\: \le \: .01 \: h^{-2}$. As surgeries only occur
at spacetime points $(x,t)$ where
$R(x,t) \sim h(t)^{-2}$,
there are no surgeries on 
$\bigcup_{t_0^{\prime \prime} \in 
[t_0 - (r_0^{\prime \prime})^2, t_0]} B(x_0, t_0^{\prime \prime}, 10 r_0)$. 
Using length distortion estimates,
we can find a parabolic neighborhood 
$P(x_0, t_0, r_0/4, - \tau r_0^2) \subset
\bigcup_{t_0^{\prime \prime} \in 
[t_0 - (r_0^{\prime \prime})^2, t_0]} B(x_0, t_0^{\prime \prime}, 10 r_0)$ 
for some fixed $\tau$.
This proves Corollary \ref{II.6.8}.

\section{II.6.3. Canonical neighborhoods and later curvature bounds on bigger balls from
curvature and volume bounds}
We now prove Proposition \ref{propII.6.3}.
We first recall its statement.
\begin{proposition} 
(cf. Proposition II.6.3) 
For any $A > 0$ one can find positive constants $\kappa(A)$, $K_1(A)$, $K_2(A)$, 
$\overline{r}(A)$,
such that for any $t_0 < \infty$
there exists $\overline{\delta}_{A}(t_0)
> 0$, decreasing in $t_0$, with the following property.  Suppose that
we have a Ricci flow with $(r, \delta)$-cutoff on a time interval
$[0, T]$, where $\delta(t) < \overline{\delta}_A(t)$ on $[t_0/2, t_0]$, with
normalized initial data.  Assume that \\
1. The solution is unscathed on a parabolic neighborhood
$P(x_0, t_0, r_0, -r_0^2)$, with $2r_0^2 < t_0$. \\
2. 
$|\Rm| \le \frac13 r_0^{-2}$
on $P(x_0, t_0, r_0, -r_0^2)$. \\
3. $\vol(B(x_0, t_0, r_0)) \: \ge \: A^{-1} r_0^3$.\\

Then\\
(a) The solution is $\kappa$-noncollapsed on scales less than $r_0$
in $B(x_0, t_0, Ar_0)$.\\
(b) Every point $x \in B(x_0, t_0, Ar_0)$ with $R(x, t_0) \ge K_1 r_0^{-2}$
has a canonical neighborhood in the sense of
Definition \ref{canonicalnbhddef}.\\
(c) If $r_0 \: \le \: \overline{r} \sqrt{t_0}$ then $R \: \le \: K_2
r_0^{-2}$ in $B(x_0, t_0, Ar_0)$.
\end{proposition}
\begin{proof}
The proof of part (a) is analogous to the proof of
Theorem \ref{8.2thm}. The proof of part (b) is analogous to
the proof of Lemma \ref{12.2lemma}. The proof of part (c) is analogous
to the proof of Theorem \ref{thmI.12.2}.
We will be brief on the parts of the proof of Proposition \ref{propII.6.3}
that are along the same lines as was done before, and will
concentrate on the differences.

For part (a),
we first remark
 that the $\kappa$-noncollapsing that we want does not follow from
the noncollapsing estimate used in the proof of Proposition
\ref{surgeryflowexists}, which would give a 
time-dependent $\kappa$.
So suppose that $(x,t_0) \in B(x_0, t_0, Ar_0)$, $\rho < r_0$,
the parabolic neighborhood $P(x, t_0, \rho, - \rho^2)$ is
unscathed and $|\Rm| \le \rho^{-2}$ there. We want to get
a lower bound on $\rho^{-3} \: \vol(B(x, t_0, \rho))$.
We first reduce the case $\rho \: < \: \frac{r(t_0)}{100}$ to the
case $\rho \:\geq\: \frac{r(t_0)}{100}$. 

Suppose that $\rho < 
\frac{r(t_0)}{100}$ and let $s$ be the largest
number so that the parabolic neighborhood 
$P(x, t_0, s, - s^2)$ is unscathed and $|\Rm| \le s^{-2}$ there.
Clearly $s \ge \rho$.
If $s \: < \: \frac{r(t_0)}{100}$ then we obtain a lower
bound on $\rho^{-3} \: \vol(B(x, t_0, \rho))$
as in Sublemma \ref{r0/100};
a canonical 
neighborhood of type (d) with small volume cannot occur,
in view of condition 3 of Proposition \ref{propII.6.3}. 

If $s \: \ge \: \frac{r(t_0)}{100}$ and $\frac{r(t_0)}{100} < r_0$
then once we have proved part (a) of the
proposition at scale $\frac{r(t_0)}{100}$, the Bishop-Gromov
inequality will give a lower bound on $\rho^{-3} \: \vol(B(x, t_0, \rho))$
and prove part (a) of the proposition at scale $\rho$.

Suppose that
$s \: \ge \: \frac{r(t_0)}{100} \ge r_0$. Let $D$ be the largest
radius so that $|\Rm| \le r_0^{-2}$ on $B(x,t_0, D)$. 
Clearly $D \ge s \ge \rho$.
If $D \ge
(A+1) r_0$ then assumption 3. of the proposition implies that
$\vol(B(x,t_0,(A+1)r_0) \ge A^{-1} r_0^{-3}$. Since $\rho < r_0$,
the Bishop-Gromov inequality implies a lower bound on 
$\rho^{-3} \: \vol(B(x, t_0, \rho))$. 

Finally, if $D < (A+1)r_0$
then there is some $x_1$ with $d_{t_0}(x,x_1) = D$ so that
$|\Rm(x_1,t_0)| = r_0^{-2} \ge 10000 (r(t_0))^{-2}$. Slightly moving
$x_1$ inward toward $x$, there is a point 
$x_2 \in B(x_1, t_0, D-cr(t_0))$ with
$|\Rm(x_2,t_0)| \ge 5000 (r(t_0))^{-2}$, for a universal constant $c$; 
cf. case (b) of the proof of Sublemma \ref{r0/100}.
If $\overline{r}$ is small then Hamilton-Ivey pinching implies that
$(x_2, t_0)$ is in a canonical neighborhood, which gives an estimate
\begin{equation}
\vol(B(x, t_0, D)) \ge \vol(B(x_2, t_0, cr(t_0))) \ge 
\const (r(t_0))^3 \ge 10^6 \const r_0^3.
\end{equation}
Since $\rho \le D$, the Bishop-Gromov inequality now implies a lower bound on 
$\rho^{-3} \: \vol(B(x, t_0, \rho))$, depending on $A$.

This shows that it suffices to consider scales $\rho$ that are at
least $r(t_0)/100$. 
To continue,
we recall the idea of the proof of Theorem \ref{8.2thm}. 
With the notation of Theorem \ref{8.2thm},
after rescaling so that $r_0 = t_0 = 1$,
we had a point $x \in B(x_0, 1, A)$
around which we wanted to prove noncollapsing.  Defining $l$ using curves
starting at $(x, 1)$, we wanted to find a point
$(y, 1/2) \in
B(x_0, 1/2, 1/2)$ so that
$l(y, 1/2)$ was bounded above by a universal constant.  Given such a point,
we concatenated a minimizing ${\mathcal L}$-geodesic (from $(x, 1)$
to $(y, 1/2)$)
with curves emanating backward in time from
$(y, 1/2)$. Then the bounded geometry near $(y, 1/2)$
allowed us 
to estimate from below the reduced volume at a time slightly less than
$1/2$. 

We knew that there was some point $y \in M$ so that
$l(y, 1/2) \: \le \: \frac{3}{2}$, but the issue in 
Theorem \ref{8.2thm} was to find 
a point $(y, 1/2) \in
B(x_0, 1/2, 1/2)$ with $l(y, 1/2)$ bounded above
by a universal constant. The idea was to take the proof that 
some point $y \in M$ has
$l(y, 1/2) \: \le \: \frac{3}{2}$ and localize it near $x_0$.
The proof of Theorem \ref{8.2thm} used the function
$h(y, t) = \phi(d(y, t) - A (2t-1)) (\overline{L}(y, 1-t) + 7)$.
Here $\phi$ was a certain nondecreasing function that is one on $(- \infty, 
1/20)$ and infinite on $[1/10, \infty)$, and
$\overline{L}(q, \tau) \: = \: 2 \sqrt{\tau} \: L(q, \tau)$.
Clearly $\min h(\cdot, 1) \: \le \: 7$ and $\min h(\cdot, 1/2)$ is
achieved in $B(x_0, 1/2, 1/10)$. The equation
$\square h \: \ge \: - (6+ C(A)) h$ implied that
$\frac{d}{dt} \min h \: \ge \: -(6 + C(A)) \min h$, and so
$(\min h)(t) \: \le \:
7 \: e^{(6+C(A))(1-t)}$. 

In the present case,
if one knew that the possible contribution of a barely admissible curve to
$h(y, t)$ was greater than $7 \: e^{(6+C(A))(1-t)} + \epsilon$
then one could still apply the maximum principle to find a point $(y, 1/2)$
with $h(y, 1/2) \: \le \: 7 \: e^{(6+C(A))/2}$.
For this, it
suffices to know that the possible contribution of a barely admissible curve
to $\overline{L}(q, \tau)$ can be bounded below by a sufficiently large
number.  However, Lemma \ref{lplusbig}  only says that we can make the
contribution of a barely admissible curve to $L$ large
(using the lower scalar curvature bound to pass from
${\mathcal L}_+$ to ${\mathcal L}$).
Because of the factor $2 \sqrt{\tau}$ in the definition of
$\overline{L}(q, \tau)$, we cannot necessarily say that its contribution
to $\overline{L}(q, \tau)$ is large. To salvage the argument,
the idea is to redefine $h$ and
redo the proof of Theorem \ref{8.2thm} in order to get an extra factor
of $\sqrt{\tau}$ in $\min h$. 

(The use of Lemma \ref{lplusbig} is similar
to what was done in the proof of Proposition \ref{surgeryflowexists}.
However, there is a difference in scales.  In Proposition
\ref{surgeryflowexists} one was working at a microscopic scale
in order to construct the Ricci flow with surgery. The function
$\overline{\delta}(t)$ in Proposition \ref{surgeryflowexists} was relevant
to this scale. In the present
case we are working at the macroscopic scale $r_0 \sim \sqrt{t_0}$ in
order to analyze the long-time behavior of the Ricci flow with
surgery.  The function $\delta_A(t)$ of Proposition
\ref{propII.6.3} is relevant to this scale.  Thus we will end up further
reducing the surgery function $\overline{\delta}(t)$ of Proposition \ref{surgeryflowexists}
in order to be able to apply Proposition \ref{propII.6.3}.)

By assumption,
$|\Rm| \: \le \: 1$ at $t = 0$.
From Lemma \ref{Rev}, $R \: \ge \: - \: \frac32 \: \frac{1}{t+1/4}$. 
Then for $t \in [t_0 - r_0^2/2, t_0]$, 
\begin{equation}
R \: r_0^2 \: \ge \: - \: \frac32 \: \frac{r_0^2}{t_0 - r_0^2/2}
\: \ge \: - \: \frac32 \: \frac{r_0^2}{2r_0^2 - r_0^2/2} \: = \: -1.
\end{equation}
After rescaling so that $r_0 = 1$, the time interval
$[t_0 - r_0^2, t_0]$ is shifted to $[0, 1]$.
Then for $t \in [\frac12, 1]$, we certainly have
$R \: \ge \: -3$.

From this, if $0 < \tau \le \frac12$ then 
$\overline{L}(y, \tau) \: \ge \: - 6 \: \sqrt{\tau} \:
\int_0^\tau \sqrt{v} \: dv \: = \: - 4 \tau^2$, so
$\hat{L}(y, \tau) \: \equiv \: \overline{L}(y, \tau) \: + \: 2 \sqrt{\tau}
\: > \: 0$.

Putting
\begin{equation}
h(y, \tau) \: = \: \phi(d_t(x_0, y) - A (2t-1)) \: \hat{L}(y, \tau)
\end{equation} 
and using the fact that $\frac{d}{dt} \sqrt{\tau} \: = \:
- \frac{d}{d\tau} \sqrt{\tau} \: = \: - \: \frac{1}{2 \sqrt{\tau}}$, 
the computations of the proof of Theorem \ref{8.2thm} give
\begin{align}
\square h \: & \ge \: - \: \left( \overline{L} + 2 \sqrt{\tau} \right) 
\: C(A) \: \phi \:
\: - \: 6 \: \phi \: - \: \frac{1}{\sqrt{\tau}} \: \phi \\
& = \: - \: C(A) h \: - \: \left( 6 \: + \: \frac{1}{\sqrt{\tau}} 
\right) \: \phi.
\notag
\end{align}
Then if  $h_0(\tau) \: = \: \min h(\cdot, \tau)$, we have
\begin{align}
\frac{d}{d\tau} \left( \log \left( \frac{h_0(\tau)}{\sqrt{\tau}} \right) 
\right) \: & = \: h_0^{-1} \: \frac{dh_0}{d\tau} \: - \: \frac{1}{2\tau} 
\: \le \: C(A) \: + \: \left( 6 \: + \: \frac{1}{\sqrt{\tau}} \right)
\: \frac{\phi}{h_0} \: - \: \frac{1}{2\tau} \\
& = \: C(A) \: + \: \left( 6 \: + \: \frac{1}{\sqrt{\tau}} \right)
\: \frac{1}{\overline{L} + 2 \sqrt{\tau}} \: - \: \frac{1}{2\tau} \notag \\
& = \: C(A) \: + \: 
\frac{6 \sqrt{\tau} \: + \: 1}{\sqrt{\tau} \overline{L} + 2 {\tau}} \: 
- \: \frac{1}{2\tau}. \notag
\end{align}
As $\overline{L} \: \ge \: - \: 4 \tau^2$,
\begin{equation}
\frac{d}{d\tau} \left( \log \left( \frac{h_0(\tau)}{\sqrt{\tau}} \right) 
\right) \:  \le \:
C(A) \: + \: 
\frac{6 \sqrt{\tau} \: + \: 1}{2 {\tau} - 4 \tau^2 \sqrt{\tau}} \: 
- \: \frac{1}{2\tau} \: \le \: C(A) \: + \: \frac{50}{\sqrt{\tau}}.
\end{equation}

As $\tau \rightarrow 0$, the Euclidean space computation gives
$\overline{L}(q, \tau) \sim |q|^2$, so
$\lim_{\tau \rightarrow 0} \frac{h_0(\tau)}{\sqrt{\tau}} \: = \: 2$.
Then
\begin{equation}
{h_0(\tau)} \: \le \:
2 \: \sqrt{\tau} \: \exp(C(A) \tau \: + \: 100 \sqrt{\tau}). 
\end{equation}
This estimate has the desired extra factor of $\sqrt{\tau}$.

It now suffices to show that for a barely admissible curve $\gamma$
that hits a surgery region at time $1-\tau$,
\begin{equation}
\int_0^\tau \sqrt{v} \: \left( R(\gamma(1-v), v) + 
|\dot{\gamma}(v)|^2
\right) \: dv \: \ge \: 
\exp(C(A) \tau \: + \: 100 \sqrt{\tau}) \: + \:
\epsilon,
\end{equation}
where $0 < \tau \le \frac12$. Choosing $\overline{\delta}_A(t_0)$ small
enough, this follows from Lemma \ref{lplusbig}
along with the lower scalar curvature bound.
Then we can apply the maximum principle and follow the proof of
Theorem \ref{8.2thm}. In our case the bounded geometry near $(y, 1/2) \in
B(x_0, 1/2, 1/2)$ comes from assumptions 2 and 3 of the
Proposition. The function $\overline{\delta}_A$ is now determined.
After reintroducing the scale $r_0$, 
this proves part (a) of the proposition.

The proof of part (b) 
is similar to the proofs of 
Lemma \ref{12.2lemma}
and Proposition \ref{surgeryflowexists}.
Suppose that for some $A > 0$ the claim is not true. Then there
is a sequence of Ricci flows $\M^\al$ which together provide a
counterexample.  In particular, some point $(x^\al, t^\al) \in
B(x_0^\al, t_0^\al, A r_0^\al)$ has $R(x^\al, t^\al) \ge K_1^\al
(r_0^\al)^{-2}$ but does not have a canonical neighborhood, where
$K_1^\al \rightarrow \infty$ as $\al \rightarrow \infty$. 
Because of the canonical neighborhood
assumption, we must have $K_1^\al \: (r_0^{\al})^{-2} \: \le \:
r(t_0^{\al})^{-2}$.
Then $2 K_1^{\al} \: \le \: K_1^{\al} \: t_0^{\al} (r_0^{\al})^{-2} 
\: \le \: t_0^{\al} \: r(t_0^{\al})^{-2}$. Since $K_1^{\al}
\rightarrow \infty$ and
the function $t \rightarrow t r(t)^{-2}$ is bounded on any finite
$t$-interval, it follows that $t_0^\al \rightarrow \infty$.
Applying point selection to each
$\M^\al$ and removing the superscripts, there are points
$\overline{x} \in B(x_0, \overline{t}, 2Ar_0)$ with
$\overline{t} \in [t_0 - r_0^2/2, t_0]$ such that
$\overline{Q} \equiv R(\overline{x}, \overline{t}) \: \ge \:
K_1 r_0^{-2}$ and $(\overline{x}, \overline{t})$ does not have a
canonical neighborhood, but each point $(x, t) \in \overline{P}$
with $R(x, t) \: \ge \: 4 \overline{Q}$ does have a canonical
neighborhood, where $\overline{P} \: = \: 
\{(x, t) \: : \: d_t(x_0, x) \: \le \: d_{\overline{t}}(x_0,
\overline{x}) + K_1^{1/2} \overline{Q}^{-1/2}, t \in
[\overline{t} - \frac14 K_1 \overline{Q}^{-1}, \overline{t}]\}$.
From (a), we have noncollapsing in $\overline{P}$. Rescaling by
$\overline{Q}^{-1}$, we have bounded curvature at bounded distances
from $\overline{x}$; see Lemma \ref{claim2II.4.2}.
Then we can extract a pointed limit $X_\infty$, which we think of as
a time zero slice, that
will have nonnegative sectional curvature.
(The required pinching for the last statement comes from the assumption that
$2r_0^2 < t_0$, along with the fact that $K_1^\al \rightarrow \infty$.)
The fact that
points $(x, t) \in \overline{P}$
with $R(x, t) \: \ge \: 4 \overline{Q}$ have a canonical
neighborhood implies that regions of large scalar curvature in $X_\infty$
have canonical neighborhoods, from which one can
deduce as in Section \ref{pf11.7} that
the sectional curvatures of $X_\infty$ are globally bounded above
by some $Q_0 > 0$.
Then for each $\overline{A}$, 
Lemmas \ref{distancedist} and \ref{claim1II.4.2}
imply that for large $\al$, the parabolic neighborhood
$P(\overline{x}, \overline{t}, \overline{A} \: \overline{Q}^{-1/2},
- \epsilon \eta^{-1} Q_0^{-1} \overline{Q}^{-1})$ is
contained in $\overline{P}$. 
(Here $\epsilon$ is a small parameter, which we absorb in the
global parameter $\epsilon$.)
In applying Lemma \ref{distancedist} we use the curvature
bound near $x_0$ coming from the hypothesis of the proposition along with 
the curvature bound near $\overline{x}$ just derived; cf. the proof
of Lemma \ref{12.2lemma}.
In addition, we claim that 
$P(\overline{x}, \overline{t}, \overline{A} \: \overline{Q}^{-1/2},
- \epsilon \eta^{-1} Q_0^{-1} \overline{Q}^{-1})$ is
unscathed. This is proved as in Section \ref{surgeryflow}.
Recall that the idea is to show that a surgery in
$P(\overline{x}, \overline{t}, \overline{A} \: \overline{Q}^{-1/2},
- \epsilon \eta^{-1} Q_0^{-1} \overline{Q}^{-1})$ implies that
$(\overline{x}, \overline{t})$ lies in a canonical neighborhood,
which contradicts our assumption. In the argument
we use the fact that $t_0^\al \rightarrow \infty$
implies $\delta(t_0^\al) \rightarrow 0$ in order to rule out
surgeries; this is the replacement for the condition
$\overline{\delta}^{\alpha} \rightarrow 0$ that was used in
Section \ref{surgeryflow}.

We extend $X_\infty$ to the maximal backward-time limit and
obtain an ancient $\kappa$-solution, which contradicts the assumption
that the points $(\overline{x}, \overline{t})$ did not have
canonical neighborhoods. This proves part (b) of the proposition.

To prove part (c), 
we can rescale $t_0$ to $1$ and then apply
Lemma \ref{claim2II.4.2}; see the end of the proof of Theorem 
\ref{thmI.12.2}. 
The $\Phi$-pinching that we use comes from the Hamilton-Ivey
estimate of (\ref{estimate}).
We recall that in the proof of Theorem \ref{thmI.12.2}
we need to get nonnegative curvature in the region $W$ near the
blowup point;
this comes from the fact that $\overline{r}^\al \rightarrow 0$ in the
contradiction argument, along with the Hamilton-Ivey pinching.
\end{proof}

This proves the proposition.
In what follows, we will want to apply it freely for arbitrary
$A$, provided that $t_0$ is large enough. To do so,
we reduce the function $\overline{\delta}$ used to define the
Ricci flow with $(r, \delta)$-cutoff, if necessary, in order to ensure
that $\overline{\delta}(t) \: \le \: \overline{\delta}_{2t}(2t)$.
Here $\overline{\delta}_{2t}(2t)$ is the quantity
$\overline{\delta}_{A}(2t)$ from Proposition \ref{propII.6.3} 
evaluated at $A = 2t$.

\section{II.6.4. Earlier scalar curvature bounds on smaller balls from
lower curvature bounds and volume bounds, in the presence of possible 
surgeries}

In this section we prove Proposition \ref{propII.6.4}.
We first recall its statement.
\begin{proposition}
There exist $\tau, \overline{r}, C_1 > 0$ and $K < \infty$ 
with the
following property. Suppose that we have a Ricci flow with $(r, \delta)$-cutoff
on the time interval $[0, t_0]$, with normalized initial data.  Let 
$r_0$ satisfy
$2C_1 h_{\max}(t_0) \: \le \: r_0 \: \le \: \overline{r} \sqrt{t_0}$,
where $h_{\max}(t_0)$ is the maximal cutoff radius for surgeries in 
$[t_0/2, t_0]$. (If there are no surgeries on $[t_0/2, t_0]$ then we put
$h_{\max}(t_0) = 0$.)
Assume \\
1. The ball $B(x_0, t_0, r_0)$ has sectional curvatures at least $- r_0^{-2}$
at each point. \\
2. The volume of any subball $B(x, t_0, r) \subset B(x_0, t_0, r_0)$ with
any radius $r > 0$ is at least $(1-\epsilon)$ times the volume of the
Euclidean ball of the same radius.\\

Then the solution is unscathed on
$P(x_0, t_0, r_0/4, - \tau r_0^2)$ and satisfies $R < K r_0^{-2}$
there.
\end{proposition}

Proposition \ref{propII.6.4}
is an analog of Theorem \ref{thmI.12.3}.  However, the proof of Proposition
\ref{propII.6.4} is more complicated, due to the need to deal with possible
surgeries.

\begin{proof}
The constants $C_1$, $K$ and $\tau$ are fixed numbers, but the
requirements on them will
be specified during the proof. The number $\overline{r}$ will emerge
from the proof, via a contradiction argument.

The first step is to prove an analog of the proposition in which
the parabolic neighborhood of the conclusion is replaced by 
a time-dependent family of metric balls.

\begin{lemma} \label{extra}
There exists $\tau^\prime > 0$
with the
following property. Suppose that we have a Ricci flow with $(r, \delta)$-cutoff
on the time interval $[0, t_0]$, with normalized initial data.  Let 
$r_0$ satisfy
$2C_1 h_{\max}(t_0) \: \le \: r_0 \: \le \: \overline{r} \sqrt{t_0}$,
where $h_{\max}(t_0)$ is the maximal cutoff radius for surgeries in 
$[t_0/2, t_0]$. (If there are no surgeries on $[t_0/2, t_0]$ then we put
$h_{\max}(t_0) = 0$.)
Assume 
\begin{enumerate}
\item The ball $B(x_0, t_0, r_0)$ has sectional curvatures at least $- r_0^{-2}$
at each point.
\item The volume of any subball $B(x, t_0, r) \subset B(x_0, t_0, r_0)$ with
any radius $r > 0$ is at least $(1-\epsilon)$ times the volume of the
Euclidean ball of the same radius.
\end{enumerate}

Then there is an unscathed time-dependent family
of metric balls 
$\bigcup_{t \in [t_0-\tau^\prime r_0^2, t_0]} B(x_0, t, r_0/2)$ 
on which the supremum of $R$ is less than $K r_0^{-2}$.
\end{lemma}
\begin{proof}
Again, $\tau^\prime$ is a fixed number but the
requirements on it will
be specified during the proof. 
The idea of the proof of the lemma
is to put oneself in a setting in
which one can apply Lemma \ref{II.6.5(a)}.

To argue by contradiction,
suppose that we have a sequence of 
Ricci flows with $(r, \delta)$-cutoff $\M^\alpha$ having
balls $B(x_0^\al, t_0^\al, r_0^\al)$ that
satisfy the assumptions of the lemma, with $\overline{r}^\alpha
\rightarrow 0$, but do not satisfy the conclusion.

\begin{sublemma}
$r_0^\al > r(t_0^\al)$ for all but a finite number of $\al$.
\end{sublemma}
\begin{proof}
If not then $r_0^\al \le r(t_0^\al)$ for 
infinitely many
$\al$. After passing to a subsequence, we can assume
that $r_0^\al \le r(t_0^\al)$ for all $\al$. We  claim that
$R \: \le \: (r_0^\al)^{-2}$ on $B(x_0^\al, t_0^\al, \frac{3r_0^\al}{4})$.
If not then $R(x^\al, t_0^\al) \: > \: (r_0^\al)^{-2}$ for some
$x^\al \in B(x_0^\al, t_0^\al, \frac{3r_0^\al}{4})$, and so
$R(x^\al, t_0^\al) \: > \: r(t_0^\al)^{-2}$.
This implies that 
$(x^\al, t_0^\al)$ is in a canonical neighborhood, which
contradicts the almost-Euclidean-volume assumption on subballs of
$B(x_0^\al, t_0^\al, r_0^\al)$.

Thus $R \: \le \: (r_0^\al)^{-2}$ on 
$B(x_0^\al, t_0^\al, \frac{3r_0^\al}{4})$. 
Lemma \ref{claim1II.4.2} 
implies that 
$R \: \le \: 16 \: (r_0^\al)^{-2}$ on
$P(x^\al_0, t^\al_0,  \frac{3r_0^\al}{4}, 
- \: \frac{1}{16} \eta^{-1} (r_0^\al)^2)$.
Furthermore, if \\
(*.1) $C_1 \: \ge \: 100$ \\
then $R \: \le \: \frac{1}{2 h^2}$ on $P(x_0^\al, t_0^\al, \frac{3r_0^\al}{4},
- \frac{1}{16} \eta^{-1} (r_0^\al)^2)$.
As surgeries only occur when $R \: \ge \: h^{-2}$,
there cannot be any surgeries in the region.

If $t \in [t_0^\al - \frac{1}{16} \eta^{-1} (r_0^\al)^2, t_0^\al]$ then
\begin{equation}
\frac{(r_0^\al)^2}{t} \: \le \: 
\frac{(r_0^\al)^2}{t_0^\al - \frac{1}{16} \eta^{-1} (r_0^\al)^2} \: = \: 
\frac{(r_0^\al)^2/t_0^\al}{1 - \frac{1}{16} \eta^{-1} (r_0^\al)^2/t_0^\al},
\end{equation}
The fact that $\lim_{\alpha \rightarrow \infty} \overline{r}^\al = 0$
implies that the Hamilton-Ivey pinching estimate improves with $\al$.
In particular, for large $\al$, we have $|\Rm| \le 
16 \: (r_0^\al)^{-2}$ on
$P(x^\al_0, t^\al_0,  \frac{3r_0^\al}{4}, 
- \: \frac{1}{16} \eta^{-1} (r_0)^\al)^2)$.
By the distance distortion estimates of Section \ref{secI.8.3},
if $\tau \le \min \left( \frac{1}{16} \eta^{-1}, 
\frac{1}{32} \log(\frac32) \right)$ then
\begin{equation}
\bigcup_{t \in [t_0^\al -\tau (r^\al_0)^2, t^\al_0]} 
B(x^\al_0, t, r^\al_0/2) \subset P(x^\al_0, t^\al_0,  \frac{3r_0^\al}{4}, 
- \: \frac{1}{16} \eta^{-1} (r_0)^\al)^2).
\end{equation}
Hence if we have \\
(*.2) $K \: \ge \: 200$ and \\
(*.3) $\tau \le \min \left( \frac{1}{16} \eta^{-1}, 
\frac{1}{32} \log(\frac32) \right)$\\
then the Ricci flows $\M^\al$ and the balls
$B(x_0^\al, t_0^\al, r_0^\al)$ do satisfy the conclusion of
Lemma \ref{extra}, contrary to assumption. This proves the
sublemma.
\end{proof}

Continuing with the proof of Lemma \ref{extra}, we can assume that
for all $\al$, we have $r_0^\al \: > \: r(t_0^\al)$.
For the flow $\M^\alpha$,
let $t_0^\alpha$ be the first time so that for some radius $r_0^\al$,
the ball $B(x_0^\al, t_0^\al, r_0^\al)$ satisfies the hypotheses of the
lemma, but
$\bigcup_{t \in [t_0^\al -\tau^\prime (r^\al_0)^2, t^\al_0]} 
B(x^\al_0, t, r^\al_0/2)$ fails to be unscathed or the
supremum of $R$ on that region is at least
$K (r_0^\al)^{-2}$.
Let $r_0^\al$ be the smallest such radius;
such a radius exists since
$r_0^\alpha \: > \: r(t_0^\alpha)$.

Let 
$\widehat{\tau}\geq 0$ 
be the supremum of the numbers $\widetilde{\tau}$
with the property that for large $\alpha$, 
$\bigcup_{t \in [t_0^\al -\widetilde{\tau} (r^\al_0)^2, t^\al_0]} 
B(x^\al_0, t, r^\al_0)$ is unscathed and
$R \ge - (r_0^\al)^{-2}$ there.

\begin{sublemma}
$\widehat{\tau}$ is bounded below by the parameter
$\tau_0$ of Lemma \ref{I.11.6(a)}, where we take $w = 1-\epsilon$ in
Lemma \ref{I.11.6(a)}.
\end{sublemma}
\begin{proof}
Suppose that $\widehat{\tau} < \tau_0$.
Put $\widehat{t}^\al \: = \: t_0^\al \: - \: (1 - \epsilon^\prime) \: 
\widehat{\tau} \: (r_0^\al)^2$, where $\epsilon^\prime$ will eventually be
taken to be a small positive number. Applying Lemma \ref{I.11.6(a)} to the
solution on $\bigcup_{t \in [t_0^\al - \: (1 - \epsilon^\prime) \: 
\widehat{\tau} \: (r_0^\al)^2, t_0^\al]} B(x_0^\al, t, r_0^\al)$,
the volume
of $B(x_0^\al, \widehat{t}^\al, r_0^\al/2)$ is at least $\frac{1}{10}$ of the
volume of the Euclidean ball of the same radius.
From Lemma \ref{Alexlemma2}, there is a subball 
$B(x_1^\al, \widehat{t}^\al, r^\al) 
\subset
B(x_0^\al, \widehat{t}^\al, r_0^\al/2)$ of radius $r^\al \: = \: 
\theta_0(1/10) \: r_0^\al/2$
with the property that all of its subballs have volume at least
$(1-\epsilon)$ times the volume of the Euclidean ball of the same radius.
The sectional curvature on $B(x_1^\al, \widehat{t}^\al, r^\al)$ 
is bounded below by
$- \: (r^\al)^{-2}$.
Since $B(x_1^\al, \widehat{t}^\al, r^\al)$ has an 
earlier time or smaller radius
than $B(x_0^\al, t_0^\al, r_0^\al)$, it
follows that $\bigcup_{t \in [\widehat{t}^\al - \tau^\prime (r^\al)^2, 
\widehat{t}^\al]} B(x_1^\al, t, r^\al/2)$
is unscathed and has $R \: < \: K (r^\al)^{-2}$ there. 
For $t \in [\widehat{t}^\al - \tau^\prime (r^\al)^2, 
\widehat{t}^\al]$, we have
\begin{equation}
\frac{(r^\al)^2}{t} \: \le \: 
\frac{(r^\al)^2}{t_0^\al - (1-\epsilon) \widehat{\tau} (r_0^\al)^2 -
\tau^\prime (r^\al)^2} =
\frac{(r^\al)^2}{(r_0^\al)^2}
\frac{(r_0^\al)^2}{t_0^\al}
\frac{1}{1 - (1-\epsilon) \widehat{\tau} \frac{(r_0^\al)^2}{t_0^\al} 
- \tau^\prime \frac{(r^\al)^2}{(r_0^\al)^2} \frac{(r_0^\al)^2}{t_0^\al}
}.
\end{equation}
The fact that $\lim_{\al \rightarrow \infty} \overline{r}^\al = 0$
implies that the Hamilton-Ivey pinching improves with $\al$. In
particular, for large $\al$, $|\Rm| \: < \: K (r^\al)^{-2}$ on
$\bigcup_{t \in [\widehat{t}^\al - \tau^\prime (r^\al)^2, 
\widehat{t}^\al]} B(x_1^\al, t, r^\al/2)$. Putting
$\widetilde{r}_0^\al \: = \: K^{-1/2} \: r^\al$,  if \\
(*.4) $K^{-1} \: \le \: \frac{1}{10} \: \tau^\prime$ \\
then $|\Rm| \: \le \: (\widetilde{r}_0^\al)^{-2}$ on
$\bigcup_{t \in [\widetilde{t}^\al - (\widetilde{r}_0^\al)^2, 
\widetilde{t}^\al]} B(x_1^\al, t, \widetilde{r}_0^\al)$ for
$\widetilde{t}^\al \in [\widehat{t}^\al - \frac12 \tau^\prime (r^\al)^2, 
\widehat{t}^\al]$.
Taking $A \: = \: 100 r_0^\al/\widetilde{r}_0^\al$, Proposition 
\ref{propII.6.3}(c)
now implies that
for large $\al$, we have
$R \: \le \: K_2(A) \: (\widetilde{r}^\al_0)^{-2}$ on
$B(x_1^\al, \widetilde{t}^\al, 100 r_0^\al)$. 
Provided that \\
(*.5) $K_2(A) \: K 
\: (\theta_0(1/10))^{-2} \: \le \: \frac{1}{1000} \: C_1^{2}$ \\
we will have $K_2(A) \: (\widetilde{r}^\al_0)^{-2} \: < \: \frac12 \: h^{-2}$
and so such balls will avoid the surgery regions.
Then the length distortion estimates
of Lemma \ref{distancedist} imply that there is some explicit
$c \in (0, \tau^\prime/2)$ so that
$\bigcup_{t \in [\widehat{t}^\al - c (r_0^\al)^2, \widehat{t}^\al]}
B(x_0^\al, t, r_0^\al) \subset 
\bigcup_{t \in [\widehat{t}^\al - c (r_0^\al)^2, \widehat{t}^\al]}
B(x_1^\al, t, 100 r_0^\al)$,
and hence
$R \: \le \: K_2(A) \: (\widetilde{r}_0^\al)^{-2}$ on
$\bigcup_{t \in [\widehat{t}^\al - c (r_0^\al)^2, \widehat{t}^\al]}
B(x_0^\al, t, r_0^\al)$.
Hamilton-Ivey pinching now
implies that for large $\al$, 
$\Rm \: \ge \: - \: (r_0^\al)^{-2}$ on
$\bigcup_{t \in [\widehat{t}^\al - c (r_0^\al)^2, \widehat{t}^\al]}
B(x_0^\al, t, r_0^\al)$.
As $c$ can be taken independent
of the small number $\epsilon^\prime$, taking $\epsilon^\prime \rightarrow 0$
we contradict the maximality of 
$\widehat{\tau}$.
\end{proof}

Continuing with the proof of the lemma,
we can now apply Lemma \ref{I.11.6(a)} to obtain 
$R \: \le \: K_0^\prime \: \tau_0^{-1}
\: (r_0^\al)^{-2}$ on 
$\bigcup_{t \in [t_0^\al - \: \tau_0 (r_0^\al)^2/2, t_0^\al]}
B(x_0^\al, t, r_0^\al/2)$. This will
give a contradiction provided that \\
(*.6) $K_0^\prime \: \tau_0^{-1} \: < \: K/2$, \\
(*.7) $\tau^\prime \: < \: \tau_0/2$ and \\
(*.8) $K_0^\prime \: \tau_0^{-1} \: \le \: \frac{1}{1000} C_1^2$, \\
where the last condition ensures that
$\bigcup_{t \in [t_0^\al-\tau^\prime (r_0^\al)^2, t_0^\al]} 
B(x_0^\al, t, r_0^\al/2)$
does not hit the surgery regions.

We first choose $\tau^\prime$ to satisfy (*.7).
We then choose $K$ to satisfy (*.4) and (*.6).
Finally, we choose $C_1$ to
satisfy (*.5) and (*.8).
This proves the lemma.
\end{proof}

We now finish the proof of Proposition \ref{propII.6.4}.
By the distance distortion estimates of Section \ref{secI.8.3}, if\\
(*.9) $\tau \le \min \left( \tau^\prime, \frac{\log(2)}{2K} \right)$\\
then $P(x_0, t_0, r_0/4, -\tau r_0^2) \subset 
\bigcup_{t_0 - \tau^\prime r_0^2} B(x_0, t, r_0/2)$.
In addition to the conditions on the parameters coming from
Lemma \ref{extra}, we can assume that $C_1$ satisfies (*.1),
$K$ satisfies (*.2), and
$\tau$ satisfies (*.3) and (*.9). The proposition follows. 
\end{proof}

\begin{remark}
The proof of Proposition \ref{propII.6.4} outlined in
\cite[Pf. of II.6.4]{Perelman2}
uses parabolic balls throughout.  Richard Bamler pointed out
to us that there is an apparent problem with this approach,
due to the issue of length distortion.
He also indicated that the problem can be circumvented
by using time-dependent families of metric balls.
\end{remark}

\begin{remark} \label{hmax}
In subsequent sections we will want to know that for any
$w > 0$, with the notation
of Corollary \ref{II.6.8}, we have
$\theta^{-1}(w) \: h_{\max}(t_0) \: \le \: r(t_0)$ if $t_0$ is
sufficiently large (as a function of $w$).
We can always achieve this by lowering the function $\overline{\delta}(\cdot)$
used to define the Ricci flow with $(r, \delta)$-cutoff so that
$\lim_{t_0 \ra \infty} \frac{h_{\max}(t_0)}{r(t_0)} \: = \: 0$.
We will assume hereafter that this is the case.
\end{remark}

\section{II.7.1. Noncollapsed pointed limits are hyperbolic} \label{II.7.1}

In this section we start the analysis of the long-time decomposition
into hyperbolic and graph manifold pieces. In the section,
$\M$ will denote a Ricci flow with $(r,\delta)$-cutoff whose initial time slice
$(\M_0, g(0))$ is compact and has normalized metric.

From Lemma \ref{pscempty}, if 
$g(0)$ has positive scalar curvature then the solution goes extinct
in a finite time.  From Lemma \ref{finitetime}, these manifolds are understood
topologically. If $g(0)$ has nonnegative scalar curvature then
either it acquires positive scalar curvature or it is flat, so again
the topological type is understood.
Hereafter we assume that the flow
does not become extinct, and 
that $R_{\min}<0$ for all $t$. 

\begin{lemma}
$V(t) \: \left( t + \frac14 \right)^{- \: \frac32}$ is nonincreasing in $t$.
\end{lemma}
\begin{proof}
Suppose first that the flow is nonsingular. In the case the lemma follows
from Lemma \ref{Rev} and the equation 
\begin{equation} \label{dVdt}
\frac{dV}{dt} \: = \: - \: \int_M R \: dV \: \le \: - \: R_{min} \: V.
\end{equation}

If there are surgeries then it only has the effect of causing further
decrease in $V$.
\end{proof}

\begin{definition} \label{VR}
Put $\overline{V} \: = \: \lim_{t \rightarrow \infty} V(t) \: \left( t + \frac14 \right)^{- \: \frac32}$
and $\hat{R}(t) \: = \: R_{min}(t) \: V(t)^{\frac23}$.
\end{definition}

\begin{lemma} \label{Rhatev}
On any time interval which is free of singular times,
and on which   $R_{min}(t) \: \le \: 0$ for all $t$
(which we are assuming), we have
\begin{equation} \label{Rhat}
\frac{d\hat{R}}{dt} \: \ge \: \frac23 \: \hat{R} \: V^{-1} \:
\int_M (R_{min} - R) \: dV.
\end{equation}
\end{lemma}
\begin{proof}
From (\ref{3dR}), $\frac{dR_{min}}{dt} \: \ge \: \frac23 \: R_{min}^2$.
Then
\begin{equation} 
\frac{d\hat{R}}{dt} \: =\: \frac{dR_{min}}{dt} \: V^{\frac23} \: + \:
\frac23 \: R_{min} \: V^{- \: \frac13} \: \frac{dV}{dt} \: \ge \:
\frac23 \: R_{min}^2 \: V^{\frac23} \: - \:
\frac23 \: R_{min} \: V^{- \: \frac13} \: \int_M R \: dV,
\end{equation}
from which the lemma follows.
\end{proof}

\begin{corollary} \label{Rmono}
If $R_{min}(t) \: \le \: 0$ for all $t$
(which we are assuming) then $\hat{R}(t)$ is nondecreasing.
\end{corollary}
\begin{proof}
If $\M$ is a nonsingular flow then the corollary follows from Lemma \ref{Rhatev}.
If there are surgeries then it only has the effect of decreasing $V(t)$, and so
possibly increasing $\hat{R}(t)$ (since $R_{min}(t) \le 0$).
\end{proof}

Put $\overline{R} \: = \: \lim_{t \rightarrow \infty} \hat{R}(t)$.

\begin{lemma} \label{RVlemma}
If $\overline{V} \: > \: 0$ then
$\overline{R} \: \overline{V}^{-2/3} \: = \: - \: \frac32$.
\end{lemma}
\begin{proof}
Suppose that $\overline{V} \: > \: 0$. 
Using Lemma \ref{Rev},
\begin{equation}
\overline{R} \: \overline{V}^{-2/3}  \: = \:
\lim_{t \rightarrow \infty} R_{min}(t) V(t)^{\frac23} \: 
\left( V(t) \: \left(t + \frac14 \right)^{- \: \frac32} \right)^{- \: \frac23} \: = \:
\lim_{t \rightarrow \infty} \left( t + \frac14 \right) \: R_{min}(t) \: \ge \: - \: \frac32.
\end{equation}
In particular, there is a limit as $t \rightarrow \infty$ of
$\left( t + \frac14 \right) \: R_{min}(t)$.
Suppose that 
\begin{equation}
\overline{R} \: \overline{V}^{-2/3} \: = \: 
\lim_{t \rightarrow \infty} \left( t + \frac14 \right) \: R_{min}(t) \: = \:
c \: > 
\: - \: \frac32.
\end{equation}
Combining this with (\ref{dVdt}) gives that for any $\mu > 0$,
$V(t) \: \le \: \const \:
t^{\mu-c}$ whenever $t$ is sufficiently large. Then
$V(t) t^{-3/2} \: \le \: \const \: t^{-(c+3/2 - \mu)}$.
Taking $\mu \: = \: \frac12 \: (c \: + \: \frac32)$, we
contradict the assumption that $\overline{V} \: > \: 0$.
\end{proof}

From the proof of Lemma \ref{RVlemma},
if $\overline{V} > 0$ then $R_{min}(t) \: \sim \: - \: \frac{3}{2t}$. 

The next proposition shows that a long-time limit will
necessarily be hyperbolic.
\begin{proposition} \label{propII.7.1}
Given the flow $\M$,
suppose that we have a sequence of parabolic neighborhoods
$P(x^\al, t^\al, r \sqrt{t^\al}, - r^2 t^\al)$, for
$t^\al \rightarrow \infty$ and some fixed $r \in (0,1)$, such that
the scalings of the parabolic neighborhoods with factor
$t^\al$ smoothly converge to some limit
solution $(\M_\infty, (\overline{x}, 1), g_\infty(\cdot))$
defined in a parabolic neighborhood
$P(\overline{x}, 1, r, -r^2)$. Then $g_\infty(t)$ has
constant sectional curvature $- \: \frac{1}{4t}$.
\end{proposition}
\begin{proof}
Suppose first that the flow is surgery-free. Because of the assumed existence of the limit
$(\M_\infty, (\overline{x}, 1), g_\infty(\cdot))$, the
original solution $\M$ has $\overline{V} > 0$.
We claim that the scalar curvature
on $P(\overline{x}, 1, r, -r^2)$ is spatially constant.
If not then there are numbers $c<0$ and $s_0, \mu > 0$ so that
\begin{equation}
\int_{B(\overline{x}, {s}, r)} (R^\prime_{min}({s}) 
\: - \: R(x, {s})) \: dV \: \le \: c
\end{equation}
whenever ${s} \in (s_0 - \mu, s_0 + \mu) \subset [1-r^2, 1]$, 
where $R^\prime_{min}({s})$ 
is the minimum of $R$ over
$B(\overline{x}, {s}, r)$. Then for large $\alpha$,
\begin{equation}
\int_{B(x^\alpha, {s} t^\alpha, r \sqrt{t^\alpha})} 
(R^\prime_{min}({s} t^\alpha) 
\: - \: R(x, {s} t^\alpha)) \: dV \: < \: \frac{c}{2}
\: \sqrt{t^\alpha}, 
\end{equation}
where $R^\prime_{min}({s} t^\alpha)$ is now the minimum of $R$ over
$B(x^\alpha, {s} t^\alpha, r \sqrt{t^\alpha})$.
Thus
\begin{equation} \label{poscont}
\int_{B(x^\alpha, {s} t^\alpha, r \sqrt{t^\alpha})} 
(R_{min}({s} t^\alpha) \: - \: R(x, {s} t^\alpha)) \: 
dV \: < \: \frac{c}{2}
\: \sqrt{t^\alpha}. 
\end{equation}

After passing to a subsequence, we can assume that
$\frac{t^{\alpha + 1}}{t^\al} \: > \: \frac{s_0 + \mu}{s_0 - \mu}$ for
all $\al$. From (\ref{Rhat}),
\begin{align} \label{Ralign}
\overline{R} \: - \: \hat{R}(0) \: & \ge \: 
\frac23 \: \int_0^\infty \hat{R}(t) \: V(t)^{-1} \:
\int_M (R_{min}(t) - R(x,t)) \: dV(x) \: dt \\
& \ge \:
\frac23 \: \sum_{\al} t^\al
\int_{s_0- \mu}^{s_0 + \mu} \hat{R}(s t^\al) \: V(s t^\al)^{-1} \:
\int_M (R_{min}(s t^\al) - R(x,s t^\al)) \: dV(x) \: ds \notag \\
& \ge \:
\frac23 \: \sum_{\al} t^\al
\int_{s_0- \mu}^{s_0 + \mu} \hat{R}(s t^\al) \: V(s t^\al)^{-1} \:
\int_{B(x^\alpha, {s} t^\alpha, r \sqrt{t^\alpha})}
(R_{min}(s t^\al) - R(x,s t^\al)) \: dV(x) \: ds. \notag
\end{align}
Using the definitions of $\overline{V}$ and $\overline{R} \: = \: 
- \: \frac32 \: \overline{V}^{\frac23}$, along with 
(\ref{poscont}), it follows that
the right-hand side of (\ref{Ralign}) is infinite.
This contradicts the fact that $\overline{R} < \infty$.

Thus $R$ is spatially constant on 
$P(\overline{x}, 1, r, -r^2)$.
As $R_{min}(t) \: \sim \: - \: \frac{3}{2t}$ on $M$, we know that the
scalar curvature $R$ at $(x,t) \in P(\overline{x}, 1, r, -r^2)$, which
only depends on $t$, satisfies
$R(t) \: \ge \: - \: \frac{3}{2t}$. It does not immediately follow that
the scalar curvature on $P(\overline{x}, 1, r, -r^2)$
equals $- \: \frac{3}{2t}$,  as $R_{min}(t)$ is the minimum of
the scalar curvature on all of $M$. However, if the scalar curvature is not identically
$- \: \frac{3}{2t}$ on $P(\overline{x}, 1, r, -r^2)$ then
again we can find $c<0$ and $s_0, \mu > 0$ so that
for large $\alpha$, (\ref{poscont}) holds for 
$s \in (s_0-\mu, s_0 +\mu) \subset [1-r^2, 1]$.
Again we get a contradiction using (\ref{Rhat}).  Thus 
$R(t) \: = \: - \: \frac{3}{2t}$ on $P(\overline{x}, 1, r, -r^2)$.
Then from (\ref{3dR}),
each time-slice of $P(\overline{x}, 1, r, -r^2)$
has an Einstein metric.  Thus the sectional curvature on
$P(\overline{x}, 1, r, -r^2)$ is $- \: \frac{1}{4t}$.

The argument goes through if one allows surgeries.
The main ingredient was the monotonicity formulas, which
still hold if there are surgeries. Note that for large 
$\al$ there are no surgeries in $P(x^\al, t^\al, r \sqrt{t^\al}, - r^2 t^\al)$
by assumption.
\end{proof}

\section{II.7.2. Noncollapsed regions with a lower curvature bound are
almost hyperbolic on a large scale}

In this section it is shown that for fixed $A, r,w > 0$ and large time $t_0$, if 
$B(x_0, t_0, r \sqrt{t_0}) \subset \M_{t_0}^+$ has volume at least 
$wr^3 t_0^{\frac32}$ and sectional curvatures at least
$- \: r^{-2} t_0^{-1}$ then the Ricci flow on the parabolic neighborhood 
$P(x_0, t_0, Ar \sqrt{t_0}, A r^2 t_0)$ is close to the flow on a
hyperbolic manifold.

We retain the assumptions of the previous section.

\begin{lemma} \label{lemII.7.2} (cf. Lemma II.7.2) \\
(a) Given $w, r, \xi > 0$ one can find $T = T(w,r,\xi) < \infty$ such that
if the ball $B(x_0, t_0, r \sqrt{t_0}) \subset \M_{t_0}^+$ at some time $t_0 \ge T$ has
volume at least $wr^3 t_0^{\frac32}$ and sectional curvatures at least
$- \: r^{-2} t_0^{-1}$ then the curvature at $(x_0, t_0)$
satisfies
\begin{equation} \label{closeness}
|2tR_{ij}(x_0, t_0) + g_{ij}|^2 \: = \: 
(2tR_{ij}(x_0, t_0) + g_{ij})(2tR^{ij}(x_0, t_0) + g^{ij}) \: < \: \xi^2.
\end{equation}
(b) Given in addition $A < \infty$ and allowing $T$ to depend on $A$, we
can ensure (\ref{closeness}) for all points in 
$B(x_0, t_0, Ar \sqrt{t_0})$. \\
(c) The same is true for $P(x_0, t_0, Ar \sqrt{t_0}, A r^2 t_0)$.

Note that the time $T$ will depend on the initial metric.
\end{lemma}
\begin{proof}
To prove (a), suppose that there is a sequence of points
$(x_0^\al, t_0^\al)$ with $t_0^\al \rightarrow \infty$
that provide a counterexample.  We wish to apply
Corollary \ref{II.6.8} with the parameter $r_0$ of the
corollary equal to $r \sqrt{t_0^\al}$.
Putting 
\begin{equation}
\widehat{w} \: = \: \left( \min_{x \in [0,1]} \frac{\int_0^x \sinh^2(s) \: 
ds}{x^3 \int_0^1 \sinh^2(s) \: ds} \right) \: w,
\end{equation}
if $r > \overline{r}(\widehat{w})$,
where $\overline{r}$ is from Corollary \ref{II.6.8},
then the hypotheses of the lemma will still be satisfied upon replacing
$w$ by $\widehat{w}$ and  
$r$ by $\overline{r}(\widehat{w})$.
Thus 
after redefining   $w$, if necessary,
we may assume that $r \le \overline{r}(w)$.
As the function
$h_{\max}(t)$ 
is nonincreasing, if $t_0^\al$ is sufficiently
large then 
$\theta^{-1}(w)h_{\max}(t_0^\al) \: \le \: r \sqrt{t_0^\al}$.
Using Corollary \ref{II.6.8} with a redefinition of $w$,
we can take a convergent pointed subsequence
as $\al \rightarrow \infty$ of the $t_0^\al$-rescalings,
whose limit is defined in an abstract parabolic neighborhood.
From Proposition \ref{propII.7.1} 
the limit will be hyperbolic, which is a contradiction.

For part (b), Corollary \ref{II.6.8} gives a bound $R \: \le \: K r_0^{-2}$ in
the unscathed parabolic neighborhood $P(x_0, t_0, r_0/4, -\tau r_0^2)$,
where 
$r_0 = \min(r, \overline{r}(w^\prime)) \sqrt{t_0}$.
We apply Proposition \ref{propII.6.3} to the parabolic
neighborhood $P(x_0, t_0, r_0^\prime, - \: (r_0^\prime)^2)$ where
$K \: r_0^{-2} \: = \: (r_0^\prime)^{-2}$. By Proposition \ref{propII.6.3}(b), each
point $ y\in B(x_0, t_0, Ar\sqrt{t_0})$ with scalar curvature at
least $Q = K^\prime(A) r_0^{-2}$ has a canonical neighborhood. 
Suppose that there is such a point.  From part (a) we have $R(x_0, t_0) < 0$, so
along a geodesic from $x_0$ to $y$ there will be some point
$x^\prime_0 \in B(x_0, t_0, Ar\sqrt{t_0})$ 
with scalar curvature $Q$. It also has a canonical neighborhood,
necessarily of type (a) or (b).
We can apply part (a) to a ball around $x^\prime_0$ with a radius
on the order of $(K^\prime(A))^{-1/2} r_0$, and with a value 
of $w$ coming from the canonical neighborhood condition, 
to get a contradiction 
for large $t_0$.
(Note that $(K^\prime(A))^{-1/2} r_0$ is proportionate to
$\sqrt{t_0}$.)
Thus $R \: \le \: K^\prime(A) r_0^{-2}$
on $B(x_0, t_0, Ar\sqrt{t_0})$. If $T$ is large enough then the
$\Phi$-almost nonnegative curvature implies that $|\Rm| \: \le \:
K^\prime(A) r_0^{-2}$. Then the noncollapsing in  
Proposition \ref{propII.6.3}(a) gives a lower local
volume bound. Hence we can apply part (a) of the lemma 
to appropriate-sized balls in
$B(x_0, t_0, \frac{A}{2} r\sqrt{t_0})$. As $A$ is arbitrary, this
proves (b) of the lemma.

For part (c), without loss of generality we can take $\xi$ small.
Suppose that the claim is not true. Then there is a point $(x_0, t_0)$
that satisfies the hypotheses of the lemma but for which there is a point
$\left(x_1, t_1 \right) \in P(x_0, t_0, Ar\sqrt{t_0},
Ar^2 t_0)$ with $|2 t_1 R_{ij}(x_1, t_1) + g_{ij}| \ge \xi$.
Without loss of generality, we can take $(x_1, t_1)$ to be a first
such point in $P(x_0, t_0, Ar\sqrt{t_0}, Ar^2 t_0)$. By part (b),
$t_1 > t_0$. Then $|2 t R_{ij} + g_{ij}| \le \xi$ on
$P(x_0, t_0, Ar\sqrt{t_0}, t_1 -  t_0)$. If $\xi$ is small then this region has
negative sectional curvature and there are no
surgeries in the region. Using the length distortion estimates of
Section \ref{secI.8.3}, we can find
$r^\prime = r^\prime(r, A)>0$ so that the sectional curvature is bounded below
by $- \: (r^\prime)^{-2} \: t_1^{-1}$ on
$B(x_0, t_1, r^\prime \sqrt{t_1})$. Also, by the evolution of volume under
Ricci flow, there will be a $w^\prime = 
w^\prime(r, w, \xi, A)$ so that
the volume of   $B(x_0, t_1, r^\prime \sqrt{t_1})$ is bounded below by
$w^\prime (r^\prime)^3 (t_1)^{\frac32}$. Thus for large $t_0$ we can apply
(b) to $B(x_0, t_1, A^\prime r^\prime \sqrt{t_1})$ with an appropriate choice of
$A^\prime$
to obtain a contradiction.
\end{proof}

\section{II.7.3. Thick-thin decomposition}

This section is concerned with the large-time decomposition of the
manifold in ``thick'' and ``thin'' parts.

\begin{definition}
For $x \in \M_t^+$,
let $\rho(x,t)$ be the unique number $\rho\in (0,\infty)$ such that 
$\inf_{B^+(x,t, \rho)} \Rm \: = \: - \: \rho^{-2}$, if such a $\rho$ exists,
and put $\rho(x,t)=\infty$ otherwise.   
\end{definition}
The function $\rho(x,t)$ is well-defined because $\M_t^+$ is a compact smooth
Riemannian manifold, so for fixed 
$(x,t)\in\M_t^+$ the quantity
$\inf_{B^+(x,t,\rho)} \Rm$ 
is a continuous nonincreasing function of $\rho$ which is negative
for sufficiently large $\rho$ if and only if $(x,t)$ lies in a connected component
with negative sectional curvature somewhere; on the other hand the function $-\rho^{-2}$ is 
continuous and  strictly increasing.  We note that when it is finite, 
the quantity $\rho(x,t)$ may be larger than the diameter of the component
of $\M_t$ containing $(x,t)$. 

As an example, if $\M$ is the flow on a manifold $M$ with spatially constant negative
curvature then for large $t$, $\rho(x,t) \sim 2 \sqrt{t}$ uniformly on $M$.
The ``thin'' part of $\M_t$, in the sense of hyperbolic geometry, can then be
characterized as the points $x$ so that $\vol(B(x,t,\rho(x,t)) \: < \: w \:
\rho^3(x,t)$, for an appropriate constant $w$.

\begin{lemma} \label{bigenough}
For any $w > 0$ we can find $\overline{\rho} = \overline{\rho}(w) > 0$
and $\overline{T}=\overline{T}(w)$
such that if $t\geq \overline{T}$ and $\rho(x,t) < \overline{\rho} \sqrt{t}$ then
\begin{equation}
\vol(B(x,t,\rho(x,t))) \: < \: w \rho^3(x,t).
\end{equation}
\end{lemma}
\begin{proof}
If the lemma is not true then there is a sequence
$(x^\al, t^\al)$ with $t^\al \rightarrow \infty$, 
$\rho(x^\al, t^\al) \: (t^\al)^{-1/2} \rightarrow 0$ and
$\vol(B(x^\al, t^\al, \rho(x^\al, t^\al))) \: \ge \: w \rho^3(x^\al, t^\al)$.
The first step is to apply Corollary \ref{II.6.8}, 
but we need to know that for large $\al$
we have 
$\rho(x^\al, t^\al) \: \ge \: \theta^{-1}(w) \: h_{\max}(t^\al)$, where
$\theta(w)$ and $h_{\max}$ are
from Corollary \ref{II.6.8}. Suppose that 
this is not the case. Then after passing to a subsequence we have
$\rho(x^\al, t^\al) \: < \: \theta^{-1}(w) \: 
h_{\max}(t^\al) \le r(t^{\al})$ for all $\al$,
where we used Remark \ref{hmax} in the last inequality.
There are points
$x^{\al,\prime} \in \overline{B(x^\al, t^\al, \rho(x^\al, t^\al))}$ 
with a sectional curvature
equal to  $- \: \rho^{-2}(x^\al, t^\al)$.
Applying the Hamilton-Ivey pinching estimate of
(\ref{estimate}) with $X^\al = \rho^{-2}(x^\al, t^\al)$, and using the fact 
that
$\lim_{\al \rightarrow \infty} t^\al X^\al \: = \:  \infty$, gives
\begin{equation} \label{Rblowup}
\lim_{\al \rightarrow \infty} 
R(x^{\al,\prime}, t^\al) \: \rho^2(x^\al, t^\al) \: = \: \infty.
\end{equation}

We claim that
the curvatures at the centers of the balls satisfy
\begin{equation}
\lim_{\al \rightarrow \infty} R(x^{\al}, t^\al) \: \rho^2(x^\al, t^\al) \: = \: \infty.
\end{equation} 
Suppose not. Then there is some number $C \in (0, \infty)$ so that
after passing to a subsequence, 
$R(x^{\al}, t^\al) \: \rho^2(x^\al, t^\al) \: \le \: C$ for all $\al$.
By continuity and (\ref{Rblowup}), 
after passing to another subsequence we can assume that there is a point
$x^{\alpha\prime \prime}$ on a time-$t^\al$ geodesic segment between 
$x^{\alpha}$ and $x^{\alpha \prime}$ so that
$R(x^{\al \prime \prime}, t^\al) \: \rho^2(x^\al, t^\al) \: = \: 2C$,
for all $\al$. We now apply Lemma \ref{claim2II.4.2} around
$(x^{\al \prime \prime}, t^\al)$ to get a contradiction to
(\ref{Rblowup}).
(More precisely, we apply a version of 
Lemma \ref{claim2II.4.2} that applies along geodesics, as in Claim 2 of 
II.4.2.) In applying Lemma \ref{claim2II.4.2}, we use the fact that
$\lim_{\al \rightarrow \infty} t^\al \: \rho^{-2}(x^\al, t^\al) \: = \:  
\infty$ in order to say that
$\lim_{\al \rightarrow \infty} 
\frac{\Phi(2C\rho^{-2}(x^\al, t^\al))}{2C\rho^{-2}(x^\al, t^\al)} \: = \: 0$,
with the notation of Lemma \ref{claim2II.4.2}.

Making a similar argument centered at other points in
$B(x^\al, t^\al, \rho(x^\al, t^\al))$, we deduce that
\begin{equation}
\lim_{\al \rightarrow \infty} \left( \inf R \big|_{B(x^\al, t^\al, \rho(x^\al, t^\al))} \right) \: 
\rho^2(x^\al, t^\al) \: = \: \infty.
\end{equation}
In particular, since $\rho(x^\al, t^\al) \: \le \: r(t^\al)$, if $\al$ is
large then each point
$y^\al \in B(x^\al, t^\al, \rho(x^\al, t^\al))$ is the center of a canonical
neighborhood.
As $\al \rightarrow \infty$, the intrinsic scales
$R(y^{\al}, t^\al)^{-1/2}$ become arbitrarily small compared to $\rho(x^\al, t^\al)$.
However, by Lemma \ref{Alexlemma2} there is a subball $B^{\al, \prime}$ of
$B(x^\al,t^\al,\rho(x^\al,t^\al))$ with radius $\theta_0(w) \rho(x^\al,t^\al)$ so that every subball of
$B^{\al,\prime}$ has almost Euclidean volume. This contradicts the existence of a small
canonical neighborhood around each point of $B(x^\al,t^\al,\rho(x^\al,t^\al))$.

We now know that for large $\al$,
$\rho(x^\al, t^\al) \: \ge \: \theta^{-1}(w) \: h_{\max}(t^\al)$.
Thus we can apply Corollary \ref{II.6.8} to get an unscathed solution on the
parabolic neighborhood
$P(x^\al, t^\al, \rho(x^\al, t^\al)/4, - \tau \rho^2(x^\al, t^\al))$, with
$R \: < \: K_0 \: \rho(x^\al, t^\al)^{-2}$ there. Applying 
Proposition \ref{propII.6.3}(c),
along with the fact that
$\lim_{\al \rightarrow \infty} t^\al \: \rho^{-2}(x^\al, t^\al) \: = \:  
\infty$,
gives an estimate $R \: \le \: K_2 \: \rho(x^\al, t^\al)^{-2}$ on
$B(x^\al, t^\al, \rho(x^\al, t^\al))$.
But then for large $\al$, the Hamilton-Ivey
pinching gives $\Rm \: > \: - \: \frac12 \: \rho(x^\al, t^\al)^{-2}$ on 
$B(x^\al, t^\al, \rho(x^\al, t^\al))$, which is a contradiction.
\end{proof}

\begin{remark}
Another approach to the above proof  would be to use the
canonical neighborhood at the center $(x^\al, t^\al)$ of the ball,
along with the Bishop-Gromov
inequality, to contradict the fact that
$\vol(B(x^\al, t^\al, \rho(x^\al, t^\al))) \: \ge \: w \rho^3(x^\al, t^\al)$. For this
to work we would have to know that the relative volume
$(\epsilon^{2} \: R(x^\al, t^\al))^{\frac32} \:  \vol(B(x^\al, t^\al, \epsilon^{-1}
R(x^\al, t^\al)^{- \: \frac12})$ of the canonical neighborhood (of
type (a) or (b)) around $(x^\al, t^\al)$ is small compared to $w$. 
This will be the case if we take the constant $\epsilon$ to
be small enough, but in a $w$-dependent way.  Although $\epsilon$ is supposed
to be a universal constant, this approach will work because when 
characterizing the graph manifold part in Section \ref{II.7.4}, $w$ can be taken to be a 
small but fixed constant.
\end{remark}

\begin{definition}
The {\em $w$-thin part} $M^-(w, t) \subset \M_t^+$ is the set of points
$x \in M$ so that either $\rho(x,t)=\infty$ or
\begin{equation}
\vol(B(x, t, \rho(x, t))) \: < \: w \:
(\rho(x, t))^3.
\end{equation}
The {\em $w$-thick part} is $M^+(w, t) = \M_t^+  - M^-(w, t)$.
\end{definition}

\begin{lemma} \label{applies}
Given $w > 0$, there are
$w^\prime = w^\prime(w) > 0$ and $T^\prime = T^\prime(w) < \infty$ so that
taking $r = \overline{\rho}(w)$
(with reference to Lemma \ref{bigenough}),
if $x_0 \in M^+(w,t)$ and $t_0 \ge T^\prime$ then
$B(x_0, t_0, r\sqrt{t_0})$ has volume at least $w^\prime r^3 t_0^{\frac32}$ and
sectional curvature at least $- \: r^{-2} \: t_0^{-1}$.
\end{lemma}
\begin{proof}
Suppose that $x_0 \in M^+(w, t_0)$. From Lemma \ref{bigenough},
if $t_0$ is big enough (as a function of $w$) then $\rho(x_0, t_0)
\: \ge \: r \: \sqrt{t_0}$.  As $\Rm \: \ge \: 
- \: \rho(x_0, t_0)^{-2}$ on $B(x_0, t_0, \rho(x_0, t_0))$, we have
$\Rm \: \ge \: 
- \: r^{-2} \: t_0^{-1}$ on $B(x_0, t_0, r \sqrt{t_0})$.
As $\vol(B(x_0, t_0, \rho(x_0, t_0))) \: \ge \: w \:
(\rho(x_0, t_0))^3$, the Bishop-Gromov inequality gives a
lower bound on $\left( r \sqrt{t_0} \right)^{-3} \:
\vol(B(x_0, t_0, r \sqrt{t_0}))$ in terms of $w$.
\end{proof}

\section{Hyperbolic rigidity and stabilization of the thick part} 
\label{hyprig}

Lemma \ref{applies} implies that Lemma \ref{lemII.7.2} applies to
$M^+(w,t)$ if $t$ is sufficiently large (as a function of $w$).
That is, if one takes a sequence of points
in the $w$-thick parts at a sequence of times tending 
to infinity then the pointed time slices subconverge,
modulo rescaling the metrics by $t^{-1}$, to complete finite volume
hyperbolic manifolds with sectional curvatures equal to
$- \: \frac14$. (The $- \: \frac14$ comes from the Ricci flow
equation along with the equation $g(t) \: = \: t \: g(1)$ for the rescaled
limit, which implies that $g(1)$ has Einstein constant $- \: \frac12$.) 
In what follows we will take the word ``hyperbolic'' for a $3$-manifold
to mean ``constant sectional curvature $- \: \frac14$''.
The next step, following Hamilton \cite{Hamiltonn},
is to show that
for large time the picture stabilizes, i.e. the limits
are unique in a strong sense. 
\begin{proposition}
\label{slowisotopy}
There exist a number $T_0<\infty$, a nonincreasing function 
$\alpha:[T_0,\infty)\ra (0,\infty)$ with
$\lim_{t \rightarrow \infty} \alpha(t) \: = \: 0$, 
a (possibly empty) collection 
$\{(H_1,x_1),\ldots,(H_k,x_k)\}$
of complete connected pointed finite-volume hyperbolic $3$-manifolds
and a family of smooth maps 
\begin{equation}
f(t):B_t= \bigcup_{i=1}^k \;B \left( x_i, \frac{1}{\alpha(t)} \right)\lra \M_t,
\end{equation}
 defined for $t\in [T_0,\infty)$, such that

1.  $f(t)$ is close to an isometry:
\begin{equation}
\|t^{-1}f(t)^*g_{\M_t}-g_{B_t}\|_{C^{[\frac{1}{\alpha(t)}]}}<\alpha(t),
\end{equation}

2. $f(t)$ defines a smooth family of maps which changes slowly
with time:
\begin{equation}
|\dot{f}(p,t)|<\alpha(t)t^{-\frac12} 
\end{equation}
for all $p\in B_t$, 
where $\dot{f}$ refers to the time derivative (as defined
with admissible curves),

and

3. $f(t)$ parametrizes more and more of the thick part:
 $M^+(\alpha(t),t)\subset \im(f(t))$ for all $t\geq T_0$. 
\end{proposition}

\begin{remark}
The analogous statement in \cite[Section 7.3]{Perelman2} is in terms of a fixed $w$.
That is, for a given $w$ one considers pointed limits of
$\{(\M^+_{t_j}, (x_j, t_j), t_j^{-1} g(t_j))\}_{j=1}^\infty$
with $\lim_{j \rightarrow \infty} t_j = \infty$ and the basepoint
satisfying 
$x_j \in M^+(w, t_j) \subset \M^+_{t_j}$
for all $j$. Considering
the possible limit spaces in a certain order, as described below, 
one extracts complete pointed finite-volume hyperbolic manifolds $\{(H_i, x_i)\}_{i=1}^k$
with $x_i$ in the $w$-thick part of $H_i$. 
There is a number
$w_0 > 0$ so that as long as $w \le w_0$, the hyperbolic manifolds
$H_i$ are independent of $w$.
For any $w^\prime > 0$,
as time goes on the $w^\prime$-thick part of  $\bigcup_{i=1}^k H_i$
better approximates $M^+(w^\prime, t)$. 
Hence the formulation of
Proposition \ref{slowisotopy} is equivalent to that of
\cite[Section 7.3]{Perelman2}. 
\end{remark}

Rather than proving Proposition \ref{slowisotopy}
using harmonic maps 
as in \cite{Hamiltonn},
we give a simple proof using
smooth compactness and a smoothing argument.  Roughly speaking, the idea is
to exploit a variant of Mostow rigidity to show that for large $t$, the components
of the $w$-thick part change slowly with time, and
are close to hyperbolic manifolds which are isolated (due to a refinement
of Mostow-Prasad rigidity). This forces them to eventually stabilize.

\begin{definition}
If $(X,x)$ and $(Y,y)$ are pointed smooth Riemannian manifolds
and $\eps>0$ then $(X,x)$ is $\eps$-close to $(Y,y)$
if there is a pointed map $f:(X,x)\ra (Y,y)$ such that
\begin{equation}
f\restr_{\ol{B(x,\eps^{-1})}}:\ol{B(x,\eps^{-1})}\ra Y
\end{equation}
 is a diffeomorphism onto its image
and 
\begin{equation}
\label{fapprox}
\|f^*g_Y-g_X\|_{C^{\eps^{-1}}}<\eps,
\end{equation}
where the norm is taken on $\ol{B(x,\eps^{-1})}$.
Note that nothing is required of $f$ on the complement of $\ol{B(x,\eps^{-1})}$.
Such a map $f$ is called an {\em $\eps$-approximation}.  
\end{definition}

We will sometimes refer
to a partially defined map $f:(X,x)\supset (W,x)\ra (Y,y)$ as an
$\eps$-approximation provided that $W$ contains $\ol{B(x,\eps^{-1})}$ and
the conditions above are satisfied. By convention we will permit $X$ and $Y$ 
to be disconnnected, in which case $\eps$-closeness only says something about the
components containing the basepoints.
We say that two maps $f_1,f_2:(X,x)\ra Y$ (not necessarily basepoint-preserving) 
are $\eps$-close if
\begin{equation}
\sup_{p\in \ol{B(x,\eps^{-1})}}d_Y(f_1(p),f_2(p))<\eps.
\end{equation}

We  recall some facts about hyperbolic manifolds. 
There is a constant $\mu_0>0$, the Margulis constant,
such that if $X$ is a complete connected finite-volume
hyperbolic $3$-manifold (orientable, as usual), $\mu\leq \mu_0$, and  
\begin{equation}
X_\mu=\{x\in X\mid \injrad(X,x)\geq \mu\}
\end{equation}
is the $\mu$-thick part of $X$, then $X_\mu$ is a nonempty 
compact manifold-with-boundary whose complement $U$
is a finite union of components $U_1,\ldots,U_k$,
where each $U_i$ is isometric either to a geodesic tube
around a closed geodesic or to  a cusp. 
In particular, $X_\mu$ is connected and there is a one-to-one
correspondence between the boundary components of $X_\mu$
and the ``thin'' components $U_i$.   For each $i$,
let $\rho_i:\ol{U_i}\ra\R$ denote either the distance function
from the core geodesic or a Busemann function, in the
tube and cusp cases respectively. In the latter case 
we normalize $\rho_i$ so that $\rho_i^{-1}(0)=\D U_i$.  (The
Busemann function goes to $- \: \infty$ as one goes down
the cusp.)
The {\em radial direction}
is the direction field on $U_i -  \core(U_i)$
defined by $\nabla\rho_i$, where $\core(U_i)$ is the core
geodesic when $U_i$ is a geodesic tube and the empty set otherwise.

\begin{lemma}
\label{prasad}
Let $(X,x)$ be a pointed complete connected finite-volume hyperbolic
$3$-manifold.  Then for each $\zeta>0$ there exists $\xi>0$
such that if $X'$ is a complete finite-volume hyperbolic
manifold with at least as many cusps as $X$,
and $f:(X,x)\ra X'$ is a $\xi$-approximation, then there is an 
isometry $\hat f: (X,x)\ra X'$ which 
is $\zeta$-close to $f$.
\end{lemma}
This was stated as Theorem 8.1 in \cite{Hamiltonn} as going 
back to the work of Mostow. We give a proof here.
The hypothesis about cusps is essential because every pointed noncompact
finite-volume hyperbolic $3$-manifold $(X,x)$ is a pointed limit of a sequence
$\{(X_i,x_i)\}_{i=1}^\infty$ 
of compact hyperbolic manifolds. 
Hence for every $\xi > 0$, if $i$ is sufficiently large then there is a
$\xi$-approximation $f \: : \: (X, x) \rightarrow (X_i, x_i)$, but there is
no isometry from $X$ to $X_i$. 

\begin{proof}
The main step is to show that for the fixed $(X, x)$, if
$\xi$ is sufficiently small then for any $\xi$-approximation
$f \: : \: (X, x) \rightarrow X^\prime$ satisfying the hypotheses of the lemma, 
the manifolds $X$ and $X^\prime$ are
diffeomorphic. The proof of this will use the
Margulis thick-thin decomposition. The rest of the assertion then
follows readily from Mostow-Prasad rigidity
\cite{Mostow,Prasad}.

Pick $\mu_1\in (0,\mu_0)$ so that $X -  X_{\mu_1}$
consists only of cusps $U_1,\ldots, U_k$. The thick part
$X_{\mu_1}$ is compact and connected. Given $\xi > 0$, let 
$f:(B(x,\xi^{-1}),x)\ra (X',x')$ be a $\xi$-approximation
as in (\ref{fapprox}).
The intuitive idea is that because of the compactness
of $X_{\mu_1}$, if $\xi$ is sufficiently small then
$f \big|_{X_{\mu_1}}$ is close to being an isometry
from $X_{\mu_1}$ to its image. 
Then $f(X_{\mu_1})$ is close to a connected component of
the thick part $X^\prime_{\mu_1}$ of $X^\prime$. As
$f(X_{\mu_1})$ and $X^\prime_{\mu_1}$ are connected,
this means that $f(X_{\mu_1})$ is close to $X^\prime_{\mu_1}$.
We will show that in fact $X_{\mu_1}$ is diffeomorphic to 
$X^\prime_{\mu_1}$. The boundary components of
$X_{\mu_1}$ correspond to the cusps of $X$ and the
boundary components of 
$X^\prime_{\mu_1}$ correspond to the connected components of
$X^\prime - X^\prime_{\mu_1}$. As $X^\prime$ has at least
as many cusps as $X$ by assumption, it follows that the
connected components of $X^\prime - X^\prime_{\mu_1}$ are
all cusps. Hence $X$ and
$X^\prime$ are diffeomorphic.

In order to show that $X_{\mu_1}$ is diffeomorphic to 
$X^\prime_{\mu_1}$, we will take a larger region
$W \supset X_{\mu_1}$ that is diffeomorphic to $X_{\mu_1}$
and show that $f(W)$ can be isotoped to $X^\prime_{\mu_1}$
by sliding it inward along the radial direction. More precisely,
in each cusp
$U_i$ put $V_i= \rho_i^{-1}([-3L,-L])$,
where $L\gg 1$ is large enough that  every cuspidal torus
$\rho_i^{-1}(s)$ with $s\in [-3L,-L]$ has diameter much less than one. 
Let $W\subset X$ be the complement of the open horoballs
at height $-2L$, i.e. 
\begin{equation}
W= X -  \bigcup_{i=1}^k  \rho_i^{-1}(-\infty, -2L).
\end{equation}
When $\xi$ 
is sufficiently small, $f$ will preserve injectivity radius
to within a factor close to $1$ for points $p\in B(x,\xi^{-1})$
with $d(p,\D B(x,\xi^{-1}))>2\injrad(X,p)$.  Therefore
when $\xi$ is small, $f$ will map each $V_i$ into 
 $X' -  X'_{\mu_1}$, and hence into one of the
connected components $U'_{k_i}$ of $X' -  X'_{\mu_1}$.
Let $Z_i'$ be the image of $Z_i=\rho_i^{-1}(-2L)$ under $f$.   Note that
$d(\core(U'_{k_i}),Z_i')\gtrsim L$ (if $\core(U'_{k_i}) \neq
\emptyset$), for otherwise 
$f^{-1}(\core(U'_{k_i}))$ would be a closed curve with 
small diameter and curvature lying in $U_i$, which contradicts
the fact that the horospheres have principal curvatures
$- \: \frac12$ (because of our normalization that the sectional curvatures are
$- \: \frac14$). Thus
for each point $z'\in Z_i'$ there is a minimizing radial 
geodesic segment $\ga'$ passing through $z'$ with
$d(z',\D\ga')\gtrsim L$.  The preimage of $\gamma^\prime$ under $f$ is
a curve $\ga\subset X$ with small curvature and length $\gtrsim L$
passing through $Z_i$. This forces the direction of $\ga$
to be nearly radial
and so transverse to $Z_i$.
Hence $Z_i'$ is transverse to the radial direction in
$U'_{k_i}$.  Combining this with the  fact that $Z_i'$ is
embedded  implies that $Z_i'$ is isotopic in $\ol{U'_{k_i}}$ to 
$\D U'_{k_i}$.  It follows that $f(W)$
is isotopic to $X'_{\mu_1}$. Then by the preceding argument
involving counting the number of cusps,
$X$ and $X'$ are diffeomorphic.
We apply Mostow-Prasad rigidity \cite{Mostow,Prasad} to deduce that
$X$ is isometric to $X'$.

We now claim that for any $\zeta > 0$, if $\xi$
is sufficiently small then the map $f$ is $\zeta$-close to an 
isometry from $X$ to $X^\prime$. Suppose not.  Then there are a number $\zeta > 0$ and
a sequence of $\frac{1}{i}$-approximations $f_i \: : \: (X, x) \rightarrow (X^\prime_i, x^\prime_i)$
so that none of the $f_i$'s are $\zeta$-close to  
any 
isometry from
$(X, x)$ to $(X^\prime_i, x^\prime_i)$. Taking a convergent subsequence of
the maps $f_i$ gives a limit isometry
$f_\infty \: : \: (X, x) \rightarrow (X^\prime_\infty, x^\prime_\infty)$.
From what has already been proven,
for large $i$ we know that $X^\prime_i$ is isometric to $X$, and so isometric
to $X^\prime_\infty$.
This is a contradiction.
\end{proof}

Recall the statement of Lemma \ref{lemII.7.2}.

\begin{definition}
Given $w >0$, let $\Lambda_w$ be the space 
of complete pointed finite-volume hyperbolic
$3$-manifolds that arise as pointed limits of
sequences $\{(\M^+_{t_i}, (x_i, t_i), t_i^{-1} g(t_i))\}_{i=1}^\infty$
with $\lim_{i \rightarrow \infty} t_i = \infty$ and the basepoint
$(x_i,t_i)$ satisfying
$(x_i,t_i) \in M^+(w, t_i) \subset \M^+_{t_i}$
for all $i$.
\end{definition}

The space $\La_w$ is compact in the smooth pointed topology.
Any element of $\La_w$ has volume at most $\overline{V}$,
the latter being defined in Definition \ref{VR}.

The next lemma summarizes the content of Lemma \ref{lemII.7.2}.

\bigskip
\begin{lemma}
\label{II.7.2again}
Given $w>0$, there is a decreasing function 
$\beta:[0,\infty)\ra (0,\infty]$
with $\lim_{s\ra\infty}\beta(s)=0$ such that 
if  
$(x,t)\in M^+(w,t)\subset\M_t^+$,
and  $Z_t$ denotes
the forward time slice
$\M_t^+$ rescaled by $t^{-1}$, then

1.  Some $(X,x)\in \La_w$
is $\beta(t)$-close to $(Z_t,(x,t))$.

2.  
$B(x,t,\beta(t)^{-1}\sqrt{t})\subset \M_t^+$
is unscathed
on the interval $[t,2t]$ and if $\ga:[t,2t]\ra \M$ is a static
curve starting at $(x,t)$, $\bar t\in [t,2t]$, then the map 
\begin{equation}
B(x,t,\beta(t)^{-1}\sqrt{t})\ra 
P(x,t,\beta(t)^{-1}\sqrt{t},t)\cap\M_{\bar t}
\end{equation}
defined by following static curves induces a map
\begin{equation}
i_{t,\bar t}:(Z_t,(x,t))\supset (B(x,t,\beta(t)^{-1}),(x,t))
\ra (Z_{\bar t},\ga(\bar t))
\end{equation}
satisfying 
\begin{equation}
\|\left(i_{t,\bar t}\right)^*g_{Z_{\bar
t}}-g_{Z_t}\|_{C^{\beta(t)^{-1}}}
<\beta(t).
\end{equation}
\end{lemma} 
\begin{proof}
This follows immediately from Lemma \ref{lemII.7.2}.
\end{proof}

\bigskip
\noindent
{\em Proof of Proposition \ref{slowisotopy}.}
If for some $w>0$ we have $\La_w=\emptyset$ then
$M^+(w,t)=\emptyset$ for large $t$.  Thus if 
$\La_w=\emptyset$ for all $w>0$,  we can take
the empty collection of pointed hyperbolic manifolds
and then 1 and 2 will be satisfied vacuously, 
and $\alpha(t)$ may be chosen so that 3 holds.
So we assume that $\La_w\neq\emptyset$ for some $w>0$. 

Since every complete finite-volume hyperbolic $3$-manifold
has a point with injectivity radius $\geq \mu_0$, 
there is a $w_0>0$ such that the collections $\{\La_w\}_{w\leq w_0}$
contain the same sets of underlying hyperbolic manifolds
(although the basepoints have more freedom when $w$ is small).
We let $H_1$ be a hyperbolic manifold from this
collection with the fewest
cusps and we choose a basepoint $x_1\in H_1$ so that
$(H_1,x_1)\in \La_{w_0}$.   Put $w_1= \frac{w_0}{2}$. 
Note that $x_1$ lies in the $w_1$-thick part of $H_1$.
In what follows we will use the fact that if $f$ is a
$\epsilon$-approximation from $H_1$, for sufficiently small
$\epsilon$, then $f(x_1)$ will lie in the $.9 w_0$-thick part of
the image.

The idea of the first step of the proof is to define a family $\{f_0(t)\}$
of $\de$-approximations $(H_1,x_1)\ra Z_t$, 
for all $t$ sufficiently large, by taking 
  a $\de$-approximation $(H_1,x_1)\ra Z_t$,
pushing it along static curves, and arguing using Lemma \ref{prasad}
that one can make small adjustments from time to time
to keep it a $\de$-approximation. The family $\{f_0(t)\}$
will not vary continuously with time, but it will have controlled
``jumps''.

More precisely, pick $T_0<\infty$ and let $\xi_1,\ldots,\xi_4>0$ be
parameters to be specified later.
We assume that $T_0$ is large enough so that 
$
2\beta(T_0)<\xi_1,
$
where $\beta$ is from Lemma \ref{II.7.2again}.
By the definition of $\La_{w_1}$, we may pick $T_0$ so that there
is a point 
$(\overline{x}_0,T_0)\in M^+(w_1,T_0)\subset\M_{T_0}^+$ 
and a $\xi_1$-approximation $f_0(T_0):(H_1,x_1)\ra (Z_{T_0},
\overline{x}_0)$.

To do the induction step, for a given $j \ge 0$ suppose that at time $2^j T_0$ 
there is a point
$(\overline{x}_j, 2^j T_0)\in M^+(w_1,2^j T_0)\subset\M_{2^j T_0}^+$ 
and a $\xi_1$-approximation $f_0(2^j T_0):(H_1,x_1)\ra 
(Z_{2^j T_0},\overline{x}_j)$. As mentioned above, if
$\xi_1$ is small then in fact $\overline{x}_j \in M^+(.9 w_0,2^j T_0)$.

By part 2 of Lemma \ref{II.7.2again}, provided that $T_0$ is 
sufficiently large we may define, for all $t\in [2^j T_0, 2^{j+1} T_0]$,
a  $2\xi_1$-approximation $f_0(t):(H_1,x_1)\ra Z_t$ by 
moving $f_0(2^j T_0)$ along static curves.  Provided that
$\xi_1$ is sufficiently small we will have 
$f_0(2^{j+1} T_0)(x_1)\in M^+(w_1,2^{j+1} T_0)$ and then part 1 of 
Lemma \ref{II.7.2again} says there is some $(H',x')\in \La_{w_1}$
with a $\beta(2^{j+1} T_0)$-approximation 
$\phi:(H',x')\ra (Z_{2^{j+1} T_0},f_0(2^{j+1} T_0)(x_1))$.  
Provided that 
$\beta(2^{j+1} T_0)$
and $\xi_1$ are sufficiently small, the partially defined map
$\phi^{-1}\circ f_0(2^{j+1} T_0)$ will define a $\xi_2$-approximation
from $(H_1,x_1)$ to $(H',x')$. Hence provided that $\xi_2$
is sufficiently small, by Lemma \ref{prasad}
the map will be $\xi_3$-close to an isometry 
$\psi:(H_1,x_1)\ra H'$. 
(In applying Lemma \ref{prasad} we use the fact that $H_1$ is
also a manifold with the fewest number of cusps in $\Lambda_{w_1}$.)
 Put $\phi_1= \phi\circ\psi$.
Provided that $\xi_3$ is sufficiently small, 
$f_0(2^{j+1} T_0)$ and $\phi_1$ will be $\xi_4$-close as maps from 
$(H_1,x_1)$ to $Z_{2^{j+1} T_0}$.  
Since $\phi_1$ is a $\beta(2^{j+1} T_0)$-approximation
precomposed with an isometry which shifts basepoints a 
distance at most $\xi_3$, it will be a $2\beta(2^{j+1} T_0)$-approximation
provided that $\xi_3<1$ and $\beta(2^{j+1} T_0)<\frac12$.  We now redefine 
$f_0(2^{j+1} T_0)$ to be $\phi_1$ and let $\overline{x}_{j+1}$ be the
image of $x_1$ under $\phi_1$. This completes the induction step.

In this way we define a family of partially defined maps
$\{f_0(t):(H_1,x_1)\ra Z_t\}_{t\in [T_0,\infty)}$.  
From the construction,
$f_0(2^jT_0)$ is a $2\beta(2^jT_0)$-approximation
for all $j\geq 0$.   Lemma \ref{II.7.2again} then implies
that there is a  function 
$\alpha_1:[T_0,\infty)\ra(0,\infty)$ decreasing to 
zero at infinity such that for all $t\in[T_0,\infty)$,
$f_0(t)$ is an $\alpha_1(t)$-approximation, and for every
$\bar t\in [t,2t]$ we may slide 
$f_0(t)$
along
static curves to define an $\alpha_1(t)$-approximation
 $h(\bar t):(H_1,x_1)\ra Z_t$ which is $\alpha_1(t)$-close
to 
$f_0(\bar t)$. 

One may now employ a standard smoothing argument to convert the
family $\{f_0(t)\}_{t\in[T_0,\infty)}$ into  a family 
$\{f_1(t)\}_{t\in [T_0,\infty)}$ which satisfies the first two 
conditions of the proposition.   If condition 3 fails to hold then
we redefine the $\La_{w}$'s by considering  limits of only those
$\{(\M^+_{t_i}, (x_i, t_i), t_i^{-1} g(t_i))\}_{i=1}^\infty$
with $t_i \rightarrow \infty$ and 
$x_i \in M^+(w,t_i)\subset \M_{t_i}^+$
not in
the image of $f_1(t_i)$.  Repeating the construction we
obtain a pointed hyperbolic manifold $(H_2,x_2)$ and a
family $\{f_2(t)\}$ defined for large $t$ satisfying conditions
1 and 2, where $\im(f_2(t))$ is disjoint from $\im(f_1(t))$ for large
$t$.  Iteration of this procedure must stop after  $k$ 
steps for some finite number $k$, in view of the fact that $\overline{V} < \infty$ and
the fact that there is a positive lower bound on the volumes of
complete hyperbolic $3$-manifolds. We get the
desired family $\{f(t)\}$ by taking the union of the maps
$f_1(t),\ldots,f_k(t)$.

\section{Incompressibility of cuspidal tori} \label{inctori}

By Proposition \ref{slowisotopy}, we know that for large times the thick 
part of the manifold can
be parametrized by a collection of (truncated) finite volume hyperbolic
manifolds.  In this section we show that each cuspidal torus
maps to an embedded incompressible torus in $\M_t$.  
(An alternative argument is given in Section \ref{II.8}.)
The strategy, due
to Hamilton, is to argue by contradiction. If such a torus were compressible
then there would be an embedded compressing disk of least area at each
time. By estimating the rate of change of the area of such disks one concludes
that the area must go to zero in finite time, which is absurd.

Let $T_0$, $\alpha$, $\{(H_1,x_1),\ldots,(H_k,x_k)\}$, $B_t$, and $f(t)$ be as
in Proposition \ref{slowisotopy}.  We will consider a fixed $H_i$, with $1\leq i\leq k$,  which
is noncompact.   Choose a  number
$a>0$ much smaller than the Margulis constant and let
$\{V_1,\ldots,V_l\}\subset H_i$
be  the cusp regions bounded by tori of diameter $a$.  Each
$V_j$ is an embedded $3$-dimensional submanifold (with boundary) of $H_i$
and is isometric to the quotient of a horoball in hyperbolic $3$-space $\mathbb{H}^3$
by the action of a copy of $\Z^2$ sitting in the stabilizer of the horoball.  The 
boundary $\D V_j$ is a totally umbilic torus whose principal curvatures are
equal to $\frac12$ everywhere.   We let $Y\subset H_i$ be the closure of the complement
of $\bigcup_{j=1}^l V_j$ in $H_i$.

Let $T_a<\infty$ be large enough  that  $B_{T_a}$  (defined as in Proposition \ref{slowisotopy}) 
contains $Y$. 

In order to focus on a given cusp, we now fix an integer $1\leq j \leq l$  and put
\begin{equation}
Z =  \partial V_{j},\quad\hat Z_t =  f(t)(Z)\\
\quad \hat Y_t =  f(t)(Y),\quad
\hat W_t =  \M_t^+ -  \Int(\hat Y_t)
\end{equation}
for every $t\geq T_a$.  The objective of this section is:

\begin{proposition}
\label{propincompressible}
The homomorphism 
\begin{equation} \label{injeqn}
\pi_1(f(t)):\pi_1(Z,\star)\ra \pi_1(\M_t^+,f(t)(\star))
\end{equation}
is a monomorphism for all $t\geq T_a$.
\end{proposition}
\begin{proof}
The proof will occupy the remainder of this section.
The first step is:
\begin{lemma} \label{stablelemma}
The kernels of the homomorphisms
\begin{equation}
\label{eqnpi1inclusion}
\pi_1(f(t)):\pi_1(Z,\star)\ra \pi_1(\M_t^+,f(t)(\star)),\quad 
\pi_1(f(t)):\pi_1(Z,\star)\ra \pi_1(\hat W_t,f(t)(\star))
\end{equation}
are independent of $t$, for all $t\geq T_a$.
\end{lemma}
\begin{proof}
We prove the assertion for the first homomorphism. The argument for the 
second one is similar.

The kernel obviously remains constant on any time interval which is free of
singular times. Suppose that $t_0\geq T_a$ is a singular time.  Then the intersection
$
\M_{t_0}^+\cap\M_{t_0}^-
$
includes into $\M_{t_0}^+$ and, by using static curves, into $\M_t$ for $t \neq t_0$
close to $t_0$.  By Van Kampen's theorem, these inclusions induce monomorphisms of the
fundamental groups.  Therefore for $t$ close to $t_0$, 
the kernel of (\ref{eqnpi1inclusion}) is the same as the kernel of 
\begin{equation}
\pi_1(f(t)):\pi_1(Z,\star)\ra \pi_1(\M_{t_0}^+\cap\M_{t_0}^-),
\end{equation}
which is independent of $t$ for times $t$ close to $t_0$.  
\end{proof}

\bigskip
We now assume that the kernel of
\begin{equation}
\pi_1(f(t)):\pi_1(Z,\star)\ra \pi_1(\M_t^+,f(t)(\star))
\end{equation}
is nontrivial for some, and hence every, $t\geq T_a$.  By Van Kampen's 
theorem and the fact that the cuspidal torus
$Z\subset Y$ is incompressible in $Y$, it follows that the kernel $K$
of 
\begin{equation}
\pi_1(f(t)):\pi_1(Z,\star)\ra \pi_1(\hat W_t,f(t)(\star))
\end{equation}
in nontrivial for all $t\geq T_a$.   By Poincar\'e duality, 
$\Image \left( \HH^1(\hat{W}_t; \R) \rightarrow \HH^1(\partial \hat{W}_t; \R) \right)$
is a Lagrangian subspace of $\HH^1(\partial \hat{W}_t; \R)$. In particular,
$\Image \left( \HH^1(\hat{W}_t; \R) \rightarrow \HH^1(Z; \R) \right)$
has
rank one. Dually,
$\Ker \left( \HH_1(Z; \R) \rightarrow \HH_1(\widehat{W}_t; \R) \right)$
has rank one and so $K$,
a subgroup of a rank-two free abelian group,
has rank one. We note that for all large $t$, 
$\hat{Z}_t$ is a convex 
boundary component of $\hat W_t$.
The main theorem of \cite{meeksyau}  implies that for 
every such $t$, there is
is a least-area compressing disk  
\begin{equation}
(N^2_t,\D N^2_t)\subset (\hat W_t,\hat Z_t).
\end{equation}
We  recall that a compressing disk is an embedded disk whose boundary curve is essential
in $\hat Z_t$. We note that by definition, $\hat W_t$ is a compact manifold
even when $t$ is a singular time.   
The embedded curve $f(t)^{-1}(\D N_t)\subset Z$ 
represents a primitive element of $\pi_1(Z)$ which,
since $K$ has  rank one, must therefore generate $K$.
It follows that
modulo taking inverses, the homotopy class of $f(t)^{-1}(\D N_t)\subset Z$ is independent of $t$. 

We define a function $A:[T_a,\infty)\ra (0,\infty)$
by letting $A(t)$ be the infimum of the areas of such embedded compressing disks.
We now show that the least-area compressing disks avoid the surgery regions.

\begin{lemma}
\label{lemavoids}
Let $\de(t)$ be the surgery parameter from 
Section \ref{II.4.4}. There is a $T = T(a) < \infty$ so that whenever $t \ge T$, 
no point in any area-minimizing compressing disk  $N_t\subset \hat W_t$ is in the
center of a $10 \delta(t)$-neck.
\end{lemma}
\begin{proof}
If the lemma were not true then
there would be a sequence of times $t_k\ra\infty$
and  for each $k$ an area-minimizing compressing disk $(N_{t_k},\D N_{t_k})
\subset (\hat W_{t_k},\hat Z_{t_k})$,
along with a point $x_k \in N_{t_k}$ that is in the center of a
$10 \delta(t_k)$-neck.
Note that the scalar curvature
near $\D N_{t_k}$ is comparable to $- \: \frac{3}{2t_k}$.
We now rescale by  $R(x_k, t_k)$, and consider the map of pointed
manifolds $f_k:(N_{t_k},\D N_{t_k},x_k)\hookrightarrow (\hat W_{t_k},\hat Z_{t_k},x_k)$ 
where the domain is equipped with the pullback Riemannian metric.
By \cite{schoen} and standard elliptic regularity, for all $\rho<\infty$
and every integer $j$, the $j^{th}$ covariant derivative of the 
second fundamental form of $f_k$ is uniformly bounded on the ball $B(x_k,\rho)\subset N_{t_k}$,
for sufficiently large $k$.  Therefore the pointed Riemannian 
manifolds $(N_{t_k},\D N_{t_k},x_k)$ 
subconverge in the smooth topology to a pointed, complete, connected, smooth manifold 
$(N_\infty,x_\infty)$.  Using the same bounds on the derivatives
of the second fundamental form, we may extract a limit mapping 
$\phi_{\infty} \: : \: N_{\infty} \ra \R \times S^2$ which 
is a $2$-sided  isometric stable minimal immersion.
By  \cite[Theorem 2]{schoenyau}, $\phi_{\infty}$ is a totally geodesic
immersion whose normal vector field in $M$ has vanishing Ricci curvature.
It follows that $\phi_{\infty}$ is a cover of
a fiber $\{ \pt \}\times S^2$.  This contradicts the fact that $N_\infty$ is noncompact.
\end{proof}

We redefine $T_a$ if necessary so that $T_a$ is greater than the $T$ of
Lemma \ref{lemavoids}.

We can isotope the surface $Z$ by moving it down the cusp $V_j$. In doing so
we do not change the group $K$
but we can make the diameter of $Z$ as small as desired. The next lemma
refers to this isotopy freedom.

\begin{lemma} \label{Hamlemma}
Given $D > 0$, there is a number $a_0 > 0$
so that for any $a \in (0,a_0)$, if $\diam(Z) = a$ and $t$ is
sufficiently large then
$\int_{\D N_t}\kappa_{\D N_t}ds \: \le \: \frac{D}{2}$ and
$\length(\D N_t) \: \le \: \frac{D}{2} \: \sqrt{t}$,
where $\kappa_{\D N_t}$ is the geodesic curvature of
$\D N_t \subset N_t$.
\end{lemma}
\begin{proof}
This is proved in \cite[Sections 11 and 12]{Hamiltonn}. We just state the
main idea.  For the purposes of this proof, we give $\hat W_t$ the
metric $t^{-1} g(t)$. First,  $\D N_t$ is the intersection of $N_t$ with
$\hat Z_t$. Because $N_t$ is minimal with respect to free boundary
conditions (i.e. the only constraint is that $\D N_t$ is in the
right homotopy class in $\hat Z_t$), it follows that $N_t$ meets
$\hat Z_t$ orthogonally. Then $\kappa_{\D N_t} =
\Pi(v,v)$, where $\Pi$ is the second fundamental form of
$\hat Z_t$ in $\hat W_t$ and $v$ is the unit tangent vector of $\D N_t$.
Given $a>0$, let $Z$ be the horospherical torus
in $V_j$ of diameter
$a$.
By Proposition \ref{slowisotopy},
for large $t$ the map $f(t)$ is close to being an isometry of pairs
$(Y, Z) \rightarrow (\hat Y_t, \hat Z_t)$. 
As $Z$ has principal
curvatures $\frac12$, we may assume that
$\Pi(v,v)$ is close to $\frac12$. This reduces the problem to showing that
with an appropriate choice of $a_0$, if $a \in (0, a_0)$ then
for large values of $t$ the length of
$\D N_t$ is guaranteed to be small. The intuition is that since
$\hat W_t$ is close to being the standard cusp $V_j$, a large piece of
the minimal disk
$N_t$ should be like a minimal surface $N_\infty$  in $V_j$
that intersects $Z$ in the given homotopy class.  Such a minimal
surface in $V_j$ essentially consists of a
geodesic curve in $Z$ going all the way down the cusp. The 
length of the intersection of $N_\infty$ with the horospherical
torus of diameter $a$ is proportionate to $a$.  Hence if $a_0$ is small enough, one
would expect that if $a < a_0$ and if $t$ is large then
the length of $\D N_t$ is small.
In particular, the length of $\D N_t$ is uniformly bounded with 
respect to $a$.
 A detailed proof appears in
\cite[Section 12]{Hamiltonn}.

Rescaling from the metric $t^{-1} g(t)$ to the original metric $g(t)$,
$\int_{\D N_t}\kappa_{\D N_t}ds$ is unchanged and
$\length(\D N_t)$ is multiplied by $\sqrt{t}$.
\end{proof}

\begin{lemma} \label{lemsupport}
For every $D>0$ there is a number $a_0>0$ with the following property.
Given $a \in (0, a_0)$, suppose that we take $Z$ to be the torus
cross-section
in $V_j$
of diameter $a$. Then
there is a number $T^\prime_a < \infty$ so that as long as
$t_0\geq T^\prime_a$, there is a smooth function $\bar A$ defined
on a neighborhood of $t_0$ such that $\bar  A(t_0)=A(t_0)$, $\bar A\geq A$
everywhere, and 
\begin{equation} \label{dereqn}
\bar A'(t_0)< \frac{3}{4}\left(\frac{1}{t_0+\frac14}\right) A(t_0)-2\pi+D.
\end{equation}
\end{lemma}
\begin{proof}
Take $a_0$ as in Lemma \ref{Hamlemma}.

For $t_0 > T_a$,
we begin with the minimizing compressing disk $N_{t_0}\subset\M_{t_0}^+$.
If $t_0$ is a surgery time and $N_{t_0}$ intersected the surgery region $\M_{t_0}^+ -
(\M_{t_0}^+ \cap \M_{t_0}^-)$ then $N_{t_0}$ would have to pass through a
$10\delta(t_0)$-neck, which is impossible by Lemma \ref{lemavoids}.
Thus $N_{t_0}$ avoids any parts added by surgery.

For $t$ close to $t_0$ we define an embedded compressing disk $S_t\subset\M_t^+$
as follows. We take $N_{t_0}$ 
and extend it slightly to a smooth surface $N_{t_0}'\subset \M_{t_0}^+$ which contains 
$N_{t_0}$ in its interior.  The surface $N_{t_0}'$ will
be unscathed on some open time interval containing $t_0$. If we let
$S_t'\subset \M_t^+$ be the  surface obtained by moving $N_{t_0}'$ along
static curves then for some $b>0$, the surface $S_t'$ will
intersect $\D \hat W_t$ transversely for all $t\in (t_0-b,t_0+b)$.
Putting
\begin{equation}
S_t =  S_t'\cap \hat W_t 
\end{equation}
defines a compressing disk for $\hat Z_t\subset \hat W_t$.

Define $\bar A:(t_0-b,t_0+b)\ra \R$ by 
\begin{equation}
\bar A(t) =  \area(S_t).
\end{equation}
Clearly $\bar A(t_0)=A(t_0)$ and $\bar A\geq A$.

For the rest of the calculation, we will view $S_t$ as a surface sitting
in a fixed manifold $M$ (a fattening of $\hat W_t$)
 with a varying metric $g(\cdot)$, and
put $S =  S_{t_0}=N_{t_0}$.

By the first variation formula for area, 
\begin{equation}
\label{eqninitial}
\bar A'(t_0)=\int_{\D S}\langle X,\nu_{\D S}\rangle\:ds
+\int_S\frac{d}{dt} \: \Big|_{t=t_0} \dvol_S,
\end{equation}
where $X$ denotes the variation vector field for $\hat Z_t$, viewed
as a surface moving in $M$, and
$\nu_{\D S}$ is the outward normal vector along $\D S$.
By Proposition \ref{slowisotopy}, there is an estimate $|X|\leq \alpha(t_0)t_0^{-\frac12}$,    
where $\alpha(t_0)\ra 0$ as $t_0 \ra \infty$.
Therefore
\begin{equation}
\label{eqnbdyterm}
\left|\int_{\D S_{t_0}}\langle X,\nu_{\D S_{t_0}}\rangle\:ds\right| \:
\leq \: \alpha(t_0) \: t_0^{-\frac12} \: \length(\D N_{t_0}).
\end{equation}
By Lemma \ref{Hamlemma}, the right-hand side of (\ref{eqnbdyterm})
is bounded above by $\frac{D}{2}$ if $t_0$ is large.

We turn to the second term in (\ref{eqninitial}).  Pick $p\in S$
and let $e_1,e_2,e_3$ be an orthonormal basis for $T_pM$ with
$e_1$ and $e_2$ tangent to $S$.  Then 
\begin{equation}
\frac{1}{\dvol_S} \: \frac{d}{dt} \: \Big|_{t=t_0} \dvol_S \: = \: \frac12 
\sum_{i=1}^2 \frac{dg}{dt} \Big|_{t=t_0} (e_i, e_i)
= -\Ric(e_1,e_1)-\Ric(e_2,e_2).
\end{equation}
Now
\begin{align}
-\Ric(e_1,e_1)-\Ric(e_2,e_2)) & =-R+\Ric(e_3,e_3) 
=-R+K(e_3,e_1)+K(e_3,e_2) \\
& =-\frac{R}{2}-K(e_1,e_2)
=-\frac{R}{2}-K_S+\gk_S, \notag
\end{align}
where $K_S$ denotes the Gauss curvature of $S$ and $\gk_S$
denotes the product of the principal curvatures.  Applying the Gauss-Bonnet formula
\begin{equation}
\int_{\D S}\kappa_{\D S}\:ds=2\pi-\int_S K_S\vol_S,
\end{equation}
the fact that $\gk_S\leq 0$
(since $S$ is time-$t_0$ minimal)
and the inequality 
\begin{equation}
R_{\min}(t)\geq -\frac{3}{2}\left(\frac{1}{t+\frac14}\right)
\end{equation}
from Lemma \ref{Rev}, we obtain
\begin{equation} \label{evolve}
\int_S\:\frac{d}{dt} \: \Big|_{t=t_0} \dvol_S \: \le \:
\int_S\:\frac34\left(\frac{1}{t_0+\frac14}\right) \: \dvol_S+\int_{\D S}\kappa_{\D S}ds-2\pi.
\end{equation}
By Lemma \ref{Hamlemma}, if $a \in (0,a_0)$, $\diam(Z) = a$ and $t_0$ is
sufficiently large then
$\int_{\D S}\kappa_{\D S}ds \: \le \: \frac{D}{2}$.
Using (\ref{eqninitial}), (\ref{eqnbdyterm}) and (\ref{evolve}),
if $t_0$ is large then
\begin{equation}
\bar A'(t_0)< \frac34\left(\frac{1}{t_0+\frac14}\right)A(t_0)-2\pi+D.
\end{equation}
This proves the lemma.
\end{proof}
\noindent
{\em Proof of Proposition \ref{propincompressible}.}
Pick $D<2\pi$. Let $a<a_0$ and $T^\prime_a$ be as in Lemma 
\ref{lemsupport}.

By Lemma \ref{lemsupport}, $A$ is bounded on compact subsets of 
$[T^\prime_a, \infty)$.
By Lemma \ref{lemavoids}, for any
$t \in [T^\prime_a, \infty)$ we can find a compact set $K_t \subset \M_t^+$
so that for all $t^\prime \in [T^\prime_a, \infty)$ sufficiently close to $t$, 
the compressing disk $N_{t^\prime}$ lies in $K_t$. If
$(K_t, g(t))$ and $(K_t, g(t^\prime))$ are $e^\sigma$-biLipschitz
equivalent then $A(t) \: \le \: e^{2\sigma} \area \left( N_{t^\prime} \right)
\: = \: e^{2 \sigma} \: A(t^\prime)$ and 
$A(t^\prime) \: \le \: e^{2\sigma} \area \left( N_{t} \right)
\: = \: e^{2 \sigma} \: A(t)$. It follows that $A$ is continuous
on $[T^\prime_a, \infty)$.

For $t \ge T_a^\prime$, put
\begin{equation}
F(t) = 
\left(t+\frac14\right)^{-\frac34} A(t) + 4 (2\pi-D) \left(t+\frac14\right)^{\frac14}.
\end{equation}
We claim that 
$F(t) \le F(T_a^\prime)$ for all $t \ge T_a^\prime$. Suppose not.  Put
$t_0 = \inf \{t \ge T_a^\prime \: : \: F(t) > F(T_a^\prime) \}$.
By continuity, $F(t_0) = F(T_a^\prime)$.
Consider the function $\bar A$ of Lemma \ref{lemsupport}.
Put
\begin{equation}
\bar F(t) =  \left(t+\frac14\right)^{-\frac34} \bar A(t) + 4 (2\pi-D) \left(t+\frac14\right)^{\frac14}.
\end{equation}
Then $\bar F(t_0) \: = \: F(t_0)$ and
in a small interval around $t_0$, we have
$\bar F \ge F$. However, (\ref{dereqn}) implies that
$\bar F^\prime(t_0)
\: < \: 0$.
There is some $\sigma > 0$ so that
for $t \in (t_0, t_0 + \sigma)$, we have
\begin{equation}
F(t) \: \le \: \bar F (t) \: \le \: \bar F(t_0) \: + \: \frac12 \: 
\bar F^\prime(t_0) \: (t - t_0) \: 
< \: \bar F(t_0) \: = \: F(t_0) \: = \: F(T_a^\prime),
\end{equation}
which contradicts the definition of $t_0$.

Thus if $t \ge T_a^\prime$ then
$F(t) \le F(T_a^\prime)$. This implies that
$A(t)$ is negative for large $t$, which contradicts the fact
that an area is nonnegative.

We have shown that the homomorphism (\ref{injeqn}) is injective if $t$ is sufficiently large.
In view of Lemma \ref{stablelemma}, the same statement holds for all $t \ge T_a$.
\end{proof}

\section{II.7.4. The thin part is a graph manifold} \label{II.7.4}

This section is concerned with showing that the thin part
$M^-(w,t)$ is a graph manifold.  We refer to Appendix \ref{geom}
for the definition of a graph manifold.
We remind the reader that this completes the proof the geometrization 
conjecture.

The next two theorems are purely Riemannian. They say that if a $3$-manifold
is locally volume-collapsed, with sectional curvature bounded below, then it
is a graph manifold. They differ slightly in their hypotheses.

\begin{theorem} (cf. Theorem II.7.4) \label{thmII.7.4}
Suppose that $(M^\al, g^\al)$ is a sequence of compact oriented Riemannian
$3$-manifolds, closed or with convex boundary, and $w^\al \rightarrow 0$.
Assume that \\
(1) for each point $x \in M^\al$ there exists a radius $\rho = \rho^\al(x)$
not exceeding the diameter of $M^\al$
such that
the ball $B(x, \rho)$ in the metric $g^\al$ has volume at most
$w^\al \rho^3$ and sectional curvatures at least $- \: \rho^{-2}$.\\
(2) each component of the boundary of $M^\al$ has diameter at most
$w^\al$, and has a (topologically trivial) collar of length one, where
the sectional curvatures are between $- \: \frac14 \: - \: \epsilon$ and
$- \: \frac14 \: - \: \epsilon$. 

Then for large $\alpha$, $M^\al$ is diffeomorphic to a graph manifold.
\end{theorem}

\begin{remark}
A proof of Theorem \ref{thmII.7.4} appears in \cite[Section 8]{Shioya-Yamaguchi}.
The proof in \cite{Shioya-Yamaguchi} is for closed manifolds, but in
view of condition (2) the method of proof clearly goes through to 
manifolds-with-boundary as considered in Theorem \ref{thmII.7.4}.
The statement of the theorem in \cite[Theorem II.7.4]{Perelman2} also has a condition
$\rho < 1$, which seems to be unnecessary. 
\end{remark}

\begin{theorem} (cf. Theorem II.7.4) \label{thmII.7.4second}
Suppose that $(M^\al, g^\al)$ is a sequence of compact oriented Riemannian
$3$-manifolds, closed or with convex boundary, and $w^\al \rightarrow 0$.
Assume that \\
(1) for each point $x \in M^\al$ there exists a radius $\rho = \rho^\al(x)$
such that
the ball $B(x, \rho)$ in the metric $g^\al$ has volume at most
$w^\al \rho^3$ and sectional curvatures at least $- \: \rho^{-2}$.\\
(2) each component of the boundary of $M^\al$ has diameter at most
$w^\al$, and has a (topologically trivial) collar of length one, where
the sectional curvatures are between $- \: \frac14 \: - \: \epsilon$ and
$- \: \frac14 \: - \: \epsilon$. \\
(3) for every $w^\prime > 0$ and $m  \in \{0, 1, \ldots, [\epsilon^{-1}]\}$, 
there exist $\overline{r} = \overline{r}(w^\prime) > 0$
and $K_m = K_m(w^\prime) <\infty$ such that for sufficiently
large $\alpha$, if $r \in (0, \overline{r}]$
and a ball $B(x,r)$ in the metric $g^\al$
has volume at least $w^\prime r^3$ and sectional curvatures at least
$- \: r^{-2}$ then $|\nabla^m \Rm|(x) \: \le \: K_m \: r^{-m-2}$.

Then for large $\alpha$, $M^\al$ is diffeomorphic to a graph manifold.
\end{theorem}

\begin{remark}
The statement of this theorem in \cite[Theorem II.7.4]{Perelman2} has the 
stronger assumption that (3) holds for all $m \ge 0$. In the
application to the locally collapsing part of the Ricci flow,
it is not clear that this stronger condition holds.  However, one does get
a bound on a large number of derivatives, which is good enough.
\end{remark}

\begin{remark}
As pointed out in \cite[Section 7.4]{Perelman2}, adding condition (3) simplifies 
the proof
and allows one to avoid both Alexandrov spaces and Perelman's 
stability theorem. 
(A proof of Perelman's stability theorem appears in \cite{Kapovitch}).
Proofs of Theorem 
\ref{thmII.7.4second} are in \cite{Besson},
\cite{Kleiner-Lott} and \cite{Morgan-Tian2}.
\end{remark}

\begin{remark}
Comparing Theorems \ref{thmII.7.4} and \ref{thmII.7.4second}, 
Theorem \ref{thmII.7.4} has the extra assumption that 
$\rho^\al(x)$ does not exceed the diameter of $M^\al$.
Without this extra assumption, the Alexandrov space arguments could give
that for large $\al$, $M^\al$ is homeomorphic to a nonnegatively curved 
Alexandrov space \cite[Theorem 1.1(2)]{Shioya-Yamaguchi}. This does not
immediately imply that $M^\al$ is a graph manifold.
\end{remark}

\begin{remark}
We give some simple examples where the collapsing theorems apply.
Let $(\Sigma, g_{hyp})$ be a closed surface with the hyperbolic
metric.  Let $S^1(\mu)$ be a circle of length $\mu$ and consider
the Ricci flow on $S^1 \times \Sigma$ with the initial metric
$S^1(\mu) \times (\Sigma, c_0 \: g_{hyp})$. The Ricci flow solution at
time $t$ is $S^1(\mu) \times (\Sigma, (c_0 + 2t) g_{hyp})$.
In the rest of this example we consider 
the rescaled metric $t^{-1} g(t)$.
Its diameter goes like
$O(t^0)$ and its sectional curvatures  go like $O(t^0)$. If we
take $\rho$ to be a small constant then for large $t$,
$B(x, t, \rho)$ is approximately a circle bundle over a ball in a hyperbolic surface
of constant sectional curvature $- \: \frac12$, with circle lengths
that go like $t^{-1/2}$. The sectional curvature on 
$B(x, t, \rho)$ is bounded below by $- \: \rho^{-2}$, and
$\rho^{-3} \: \vol(B(x, t, \rho)) \: \sim \: t^{-1/2}$.
Theorems \ref{thmII.7.4}
and \ref{thmII.7.4second} both apply.

Next, consider a compact
$3$-dimensional nilmanifold that evolves under the Ricci flow.
Let $\theta_1, \theta_2, \theta_3$ be affine-parallel $1$-forms on $M$
which lift to Maurer-Cartan forms on the Heisenberg group, with
$d\theta_1 \: = \: d\theta_2 \: = \: 0$ and $d\theta_3 \: = \:
\theta_1 \: \wedge \: \theta_2$.
Consider the metric
\begin{equation}
g(t) \: = \: \alpha^2(t) \: \theta_1^2 \: + \: 
\beta^2(t) \: \theta_2^2 \: + \: 
\gamma^2(t) \: \theta_3^2.
\end{equation}
Its sectional curvatures are
$R_{1212} \: = \: - \: \frac34 \: \frac{\gamma^2}{\alpha^2 \beta^2}$ and
$R_{1313} \: = \: R_{2323}
\: = \: \frac14 \: \frac{\gamma^2}{\alpha^2 \beta^2}$.
The Ricci tensor is
\begin{equation}
\Ric \: = \: \frac12 \: \frac{\gamma^2}{\alpha^2 \beta^2} \:
\left( - \: \alpha^2 \: \theta_1^2 \: - \: \beta^2 \: \theta_2^2 \: 
+ \: \gamma^2 \: \theta_3^2 \right).
\end{equation}
The general solution to the Ricci flow equation is of the form
\begin{align}
\alpha^2(t) \: & = \: A_0 (t + t_0)^{1/3}, \\
\beta^2(t) \: & = \: B_0 (t + t_0)^{1/3}, \notag \\
\gamma^2(t) \: & = \: \frac{A_0 B_0}{3} (t + t_0)^{-1/3}. \notag
\end{align}
In the rest of this example we consider 
the rescaled metric $t^{-1} g(t)$. 
Its diameter goes like
$t^{-1/3}$, its volume goes like $t^{-4/3}$
and its sectional curvatures
go like $t^0$. If we take
$\rho(x) \: = \: \diam$ then 
$\rho^{-3} \vol(B(x, \rho)) \sim t^{-1/3}$, so both Theorem \ref{thmII.7.4}
and Theorem \ref{thmII.7.4second} apply.
We could also take $\rho(x)$ to be a small constant
$c > 0$, in which case Theorem \ref{thmII.7.4second} applies.

In general, among the eight maximal homogeneous geometries, the
rescaled solution for a compact $3$-manifold with geometry 
$H^2 \times \R$ or $\widetilde{\SL_2(\R)}$
will collapse to a hyperbolic surface of constant 
sectional curvature $- \: \frac12$.
The rescaled solution for a $\Sol$ geometry will collapse to a circle.
The rescaled solution for an $\R^3$ or $\Nil$ geometry will collapse to
a point.

We remark that although these homogeneous solutions
are collapsing in the sense of Theorem
\ref{thmII.7.4}, there is no contradiction with the
no local collapsing result of Theorem \ref{nolocalcollapse},
which only rules out
local collapsing on a finite time interval.
\end{remark}

Returning to our Ricci flow with surgery, recall the statement of
Proposition \ref{slowisotopy}. If the collection
$\{H_1,\ldots, H_k\}$ of Proposition \ref{slowisotopy} is nonempty then
for large $t$, let $\widehat{H}_i(t)$ be the result of removing from $H_i$ the horoballs whose
boundaries are at 
distance approximately $\frac12 \: \alpha(t)$ from the basepoint $x_i$.
(If there are no such horoballs then $\widehat{H}_i(t) = H_i$.)
Put
\begin{equation}
M_{thin}(t) \: = \: \M_t^+ - f(t) (\widehat{H}_1(t) \cup \ldots \widehat{H}_k(t)).
\end{equation}

\begin{proposition} \label{thingraph}
For large $t$, $M_{thin}(t)$ is a graph manifold.
\end{proposition}
\begin{proof}
We give two closely related proofs, one using Theorem \ref{thmII.7.4second} and
one using Theorem \ref{thmII.7.4}.

If the proposition is not true then there is a sequence $t^\al \rightarrow \infty$
so that for each $\al$, $M_{thin}(t^\al)$ is not a graph manifold.  
Let $M^\al$ be the manifold obtained from  $M_{thin}(t^\al)$ 
by throwing away connected components which are closed and admit metrics
of nonnegative sectional curvature, and put $g^\al \: = \: (t^\al)^{-1} g(t^\al)$.
Since any closed manifold of nonnegative sectional curvature is a graph 
manifold by \cite{Hamiltonnnnn}, for each $\al$ the manifold $M^\al$ is not a graph
manifold.

We first show that the assumptions of Theorem \ref{thmII.7.4second} are verified.

\begin{lemma} \label{condition}
Condition (3) in Theorem \ref{thmII.7.4second} holds for the $M^\al$'s.
\end{lemma}
\begin{proof}
With $w^\prime$ being a parameter as in  
Condition (3) in Theorem \ref{thmII.7.4second},
let $\overline{r} = \overline{r}(w^\prime)$ be the parameter of Corollary \ref{II.6.8}. 
It is enough to show that for large $\al$, if $r \in (0, \overline{r} \sqrt{t^\al}]$,
$x^\al \in \M_{t^\al}^+$, and $B(x^\al, t^\al, r)$
has volume at least $w^\prime r^3$ and sectional curvatures
bounded below by $- \: r^{-2}$, then $|\nabla^m \Rm|(x^\al, t^\al) \: \le \:
K_m \: r^{-m-2}$ for an appropriate choice of constants $K_m$.

To prove this by contradiction, we  assume that after passing
to a subsequence if necessary, there are  $r^\al\in (0,\overline{r} \sqrt{t^\al}]$
and $x^\al \in \M_{t^\al}^+$ such that $B(x^\al, t^\al,r^\al)$ has 
volume
at least $w^\prime (r^\al)^3$ and sectional curvature at least 
$-(r^\al)^{-2}$, but $\lim_{\al \rightarrow \infty} \:
(r^\al)^{m+2} \: |\nabla^m \Rm|(x^\al, t^\al) \: = \:
\infty$
for some $m \le [\epsilon^{-1}]$.

In the notation of Corollary \ref{II.6.8}, if
$r^\al \ge \theta^{-1}(w^\prime) h_{\max}(t^\al)$ for infinitely many $\al$
then for these $\al$, Corollary \ref{II.6.8} gives a curvature bound on an unscathed
parabolic neighborhood 
$P \left( x^\al, t^\al, r^\al/4, - \tau (r^\al)^2 \right)$ and hence,
by Appendix \ref{applocalder}, derivative bounds at $(x^\al, t^\al)$.
This is a contradiction.  Therefore we may assume that 
$r^\al < \theta^{-1}(w^\prime) h_{\max}(t^\al) \le r(t^\al)$ for all $\al$,
where we used Remark \ref{hmax} for the last inequality.

Suppose first that 
$R(x^\al, t^\al) \le \left( r(t^\al) \right)^{-2}$.
By Lemma \ref{claim1II.4.2}, there is an estimate
$R \: \le \: 16 \: \left( r(t^\al) \right)^{-2}$ on the parabolic neighborhood
$P \left( x^\al, t^\al, \frac14 \eta^{-1} r(t^\al),
- \: \frac{1}{16} \eta^{-1} \left( r(t^\al) \right)^2 \right)$. 
A surgery in this neighborhood could only occur where 
$R \ge h_{\max}(t^\al)^{-2}$. For large $\al$,
$h_{\max}(t^\al)^{-2} >> r(t^\al)^{-2}$ by Remark \ref{hmax}. 
Hence this neighborhood is unscathed.
Appendix \ref{applocalder} now gives bounds of the form
$|\nabla^m \Rm|(x^\al, t^\al) \: \le \: \const \: 
\left( r(t^\al) \right)^{-m-2} \: \le \: \const \:
(r^\al)^{-m-2}$, which is a contradiction..

Suppose now that $R(x^\al, t^\al) > \left( r(t^\al) \right)^{-2}$.
Then $(x^\al, t^\al)$ is in the center of a canonical neighborhood and there
are universal estimates $|\nabla^m \Rm|(x^\al, t^\al) \: \le \: \const(m) \: 
R(x^\al, t^\al)^{\frac{m+2}{2}}$ for all $m \le [\epsilon^{-1}]$. Hence
in this case, it suffices to show that $R(x^\al, t^\al)$ is bounded above by
a constant times $(r^\al)^{-2}$, i.e. it suffices to get a contradiction 
just to the assumption that
$\lim_{\al \rightarrow \infty} \:
(r^\al)^{2} \: R(x^\al, t^\al) \: = \: \infty$.

So suppose that $\lim_{\al \rightarrow \infty} \:
(r^\al)^{2} \: R(x^\al, t^\al) \: = \: \infty$.
We claim that $\lim_{\al \rightarrow \infty} \: 
(r^\al)^{2} \: \inf R \Big|_{B \left( x^\al, t^\al,
r^\al \right)} \: = \: \infty$. 
Suppose not.  Then there is some $C \in (0, \infty)$ so that
after passing to a subsequence, there are points
$x^{\al \prime} \in B \left( x^\al, t^\al, r^\al
 \right)$ with
$(r^\al)^{2} \: R(x^{\al \prime}, t^\al) < C$. Considering points
along the time-$t^\al$ geodesic segment from
$x^{\al \prime}$ to $x^\al$,
for large $\al$ we can find points 
$x^{\al \prime \prime} \in B \left( x^\al, t^\al,
r^\al \right)$ with
$(r^\al)^{2} \: R(x^{\al \prime\prime}, t^\al) \: = \: 2C$.
Applying Lemma \ref{claim2II.4.2} at $(x^{\al \prime \prime}, t^\al)$,
or more precisely a version that applies along geodesics as in
Claim 2 of II.4.2, we obtain a contradiction to the assumption
that $\lim_{\al \rightarrow \infty} \:
(r^\al)^{2} \: R(x^\al, t^\al) \: = \: \infty$. In applying
Lemma \ref{claim2II.4.2} we use that
$r^\al < r(t^\al)$ and $\lim_{\al \rightarrow \infty} r(t^\al) \: = \: 0$,
giving
$\lim_{\al \rightarrow \infty} \left( r^\al \right)^{-2} \: t^\al \: = \:
\infty$, in order to say that
$\lim_{\al \rightarrow \infty} 
\frac{\Phi \left( 2C (r^\al)^{-2} \right)}{ 2C (r^\al)^{-2}} \: = \: 0$.

Hence for large $\al$, every point $x \in B \left( x^\al, t^\al,
r^\al \right)$ is in the
center of a canonical neighborhood of size comparable to
$R(x, t^\al)^{- \: \frac12}$, which is small compared to $r^\al$.
On the other hand,
from Lemma \ref{Alexlemma2}, 
there is a ball $B^\prime$ of radius $\theta_0(w^{\prime}) \:
r^\al$ in 
$B \left( x^\al, t^\al, r^\al \right)$ so that every subball of $B^\prime$ has 
almost-Euclidean volume. This is a contradiction.

This proves the lemma.
\end{proof}

We continue with the proof of Proposition \ref{thingraph}.
By construction there is a sequence $w^\al \rightarrow 0$ so that
conditions (1) and (2) of Theorem \ref{thmII.7.4second} hold
with $\rho^\al(x^\al) = (t^\al)^{- \frac12} \: \rho(x^\al, t^\al)$.
Hence for large $\al$, $M^\al$ is
diffeomorphic to a graph manifold.  This is a contradiction to the choice of
the $M^\al$'s and proves the theorem.

We now give a proof that instead uses Theorem \ref{thmII.7.4}.
Let $d^\al$ denote the diameter of $M^\al$.
If we take $\rho^\al(x^\al) = (t^\al)^{- \frac12} \: \rho(x^\al, t^\al)$
then we can apply Theorem \ref{thmII.7.4} as long as
that the diameter statement in condition (1) of Theorem \ref{thmII.7.4} is
satisfied. 
If it is not satisfied then there is some point $x^\al \in M^\al$ with
$\rho^\al(x^\al) \: > \:  d^\al$. 
The sectional curvatures of $M^\al$ are bounded below by 
$- \: \frac{1}{\rho^\al(x^\al)^2}$, and so are bounded below by
$- \: \frac{1}{(d^\al)^2}$.
If there is a subsequence  with
$\frac{\rho^\al(x^\al)}{d^\al} \le C < \infty$  then
\begin{equation}
\vol(M^\al) \: = \: \vol(B(x^\al, \rho^\al(x^\al))) \: \le \: w^\al \:
\rho^\al(x^\al)^3 \: \le \: w^\al \: C^3 \: (d^\al)^3
\end{equation}
and we can apply Theorem \ref{thmII.7.4} with $\rho^\al \: = \: d^\al$,
after redefining $w^\al$.
Thus we may assume that 
$\lim_{\al \rightarrow \infty} \frac{\rho^\al(x^\al)}{d^\al} = \infty$.
If there is a subsequence with
$\frac{\vol(M^\al)}{(d^\al)^3}  \rightarrow 0$ then we
can apply Theorem \ref{thmII.7.4} with $\rho^\al \: = \: d^\al$.
Thus we may assume that $\frac{\vol(M^\al)}{(d^\al)^3}$ is bounded away from zero.
After rescaling the metric to make the diameter one, we are in a
noncollapsing situation 
with the lower sectional curvature bound going to zero.
By the argument in the proof of Lemma \ref{condition} (which used
Corollary \ref{II.6.8}) there are uniform $L^\infty$-bounds on $\Rm(M^\al)$ and
its covariant derivatives.  After passing to a subsequence, there is a
limit $(M_\infty, g_\infty)$ in the smooth topology which is diffeomorphic
to $M^\al$ for large $\al$, and carries a metric of nonnegative sectional curvature. 
As any boundary component of $M_\infty$ would have to have a neighborhood
of negative sectional curvature (see 
the definition of $M_{thin}$ and 
condition (2) of Theorem \ref{thmII.7.4}),
$M_\infty$ is closed. However, by construction $M^\al$ has no connected components
which are closed and admit metrics of nonnegative sectional curvature.
 This contradiction shows that the diameter statement in condition
(1) of Theorem \ref{thmII.7.4} is satisfied. The other conditions of
Theorem \ref{thmII.7.4} are satisfied as before.
\end{proof}

Thus for large $t$,  $\M_t^+$ has a decomposition into
a piece  $f(t) (\widehat{H}_1(t) \cup \ldots \widehat{H}_k(t))$, whose interior admits
a complete finite-volume hyperbolic metric, and the complement, which is a
graph manifold.

In addition, by Section \ref{inctori} 
the cuspidal tori are incompressible in $\M_t^+$.
By Lemma \ref{reconstruct2},
the initial (connected) manifold $\M_0$ is diffeomorphic to a connected sum of
the connected components of $\M_t$, along with some possible additional connected sums
with a finite number of $S^1 \times S^2$'s and quotients of the round $S^3$.
This proves the geometrization conjecture of Appendix \ref{geom}. 

\section{II.8. Alternative proof of cusp incompressibility} \label{II.8}

The goal of this section is to prove Perelman's Proposition II.8.2, which
gives a numerical characterization of
the geometric type of a compact $3$-manifold. It also contains an
independent proof of the incompressibility of the cuspidal ends
of the hyperbolic piece in the geometric decomposition.

We recall from Section \ref{I.2.2}
that $\lambda(g)$ is the first eigenvalue of
$-4 \triangle \: + \: R$, and can also be expressed as
\begin{equation} \label{Rayleigh}
\lambda(g) \: = \: \inf_{\Phi \in C^\infty(M) \: : \:
\Phi \neq 0} \frac{\int_M \left( 4 |\nabla \Phi|^2 \: + \:
R \: \Phi^2 \right) \: dV}{\int_M \Phi^2 \: dV}.
\end{equation}

From Lemma \ref{2/3}, if $g(\cdot)$ is a Ricci flow and
$\lambda(t) \: = \: \lambda(g(t))$ then
\begin{equation} \label{ll}
\frac{d}{dt} \lambda(t) \: \ge \: \frac23 \: 
\lambda^2(t).
\end{equation}
From Lemma \ref{expandingprop},
$\lambda(t) \: V(t)^{\frac23}$ is nondecreasing
when it is nonpositive.

For any metric $g$, there are inequalities
\begin{equation} \label{test}
\min R \: \le \: \lambda(g) \: \le \: \frac{\int_M R \: dV}{\vol(M, g)},
\end{equation}
where the first inequality follows directly from (\ref{Rayleigh})
and the second inequality comes from using $1$ as a
test function in (\ref{Rayleigh}).

Perelman's proof of his Proposition II.8.2 uses the functional
$\lambda(g) V(g)^{\frac23}$. The functional
$R_{min} V(g)^{\frac23}$ plays a similar role.  For example, from
Corollary \ref{Rmono}, 
$\widehat{R}(t) = R_{min}(t) V(t)^{\frac23}$ is nondecreasing
when it is nonpositive. We first give a proof of an analog of
Proposition II.8.2 that uses $R_{min} V(g)^{\frac23}$ instead 
of $\lambda(g) V(g)^{\frac23}$. The technical simplification is that 
when
$R_{min} V(g)^{\frac23}$ is nonpositive, it is nondecreasing under a surgery,
as surgeries are only done in regions of large positive scalar curvature,
so $R_{min}$ doesn't change,
and a surgery reduces volume. (A possible extinction of a component
clearly doesn't change $R_{min} V(g)^{\frac23}$.)
We show that a minimal-volume hyperbolic submanifold of $M$
has incompressible tori, which gives a different approach to
Section \ref{inctori}. 

Perelman's alternative approach to
Section \ref{inctori} uses the functional $\lambda(g) V(g)^{\frac23}$
instead of $R_{min} V(g)^{\frac23}$. Our use of $R_{min} V(g)^{\frac23}$ and
the sigma-invariant
$\sigma(M)$, instead of $\lambda(g) V(g)^{\frac23}$ and $\overline{\lambda}$,
is inspired by \cite{Anderson2}.

We then give the arguments using $\lambda(g) V(g)^{\frac23}$,
thereby proving Perelman's Proposition II.8.2.
The main technical difficulty is to
control how $\lambda(g) V(g)^{\frac23}$ changes under a surgery.

\subsection{The approach using the $\sigma$-invariant}

We first give some well-known results about the sigma-invariant.
We recall that the sigma-invariant of a closed connected
manifold $M$ of dimension $n \ge 3$ is given by
\begin{equation} \label{yamabe}
\sigma(M) \: = \: \sup_{{\mathcal C}} \inf_{g \in {\mathcal C}}
\frac{\int_M R(g) \: \dvol(g)}{\vol(M, g)^{\frac{n-2}{n}}},
\end{equation}
where ${\mathcal C}$ runs over the conformal classes of Riemannian
metrics on $M$. 
From the solution to the Yamabe problem,
the infimum in (\ref{yamabe}) is realized by a metric of constant
scalar curvature in the given conformal class.
It follows that if $\sigma(M) > 0$ then $M$ admits a metric with
positive scalar curvature. Conversely, suppose that $M$ admits
a metric $g_0$ with positive scalar curvature. Let ${\mathcal C}$
be the conformal class containing $g_0$.  Then
\begin{equation}
\inf_{g \in {\mathcal C}}
\frac{\int_M R(g) \: \dvol(g)}{\vol(M, g)^{\frac{n-2}{n}}} \: = \: 
\inf_{u > 0} \frac{
\int_M \left(
\frac{4(n-1)}{n-2} \:  |\nabla u|^2 \: + \:
R(g_0) \: u^2 \right) \dvol_M(g_0)
}{
\left( \int_M u^{\frac{2n}{n-2}} \: \dvol_M(g_0) \right)^{\frac{n-2}{n}}}
\end{equation}
is positive, in view of the Sobolev embedding theorem, and so $\sigma(M) > 0$.

We claim that if $\sigma(M) \le 0$ then
\begin{equation} \label{testt}
\sigma(M) \: = \: \sup_{g} R_{min}(g) \: V(g)^{\frac2n}.
\end{equation}
To see this, as
the infimum in (\ref{yamabe}) is realized by a metric of constant
scalar curvature in the given conformal class, it follows that 
$\sigma(M) \: \le \: \sup_{g} R_{min}(g) \: V(g)^{\frac2n}$.
Now given a Riemannian metric $g$, the infimum in (\ref{yamabe})
within the corresponding conformal class ${\mathcal C}$ 
equals $\widetilde{R} \: V(\widetilde{g})^{\frac{2}{n}}$
for a metric $\widetilde{g} \: = \: u^{\frac{4}{n-2}} g$ with
constant scalar curvature $\widetilde{R}$. Then
\begin{equation}
\widetilde{R} \: V(\widetilde{g})^{\frac{2}{n}} \: = \: 
\inf_{u > 0} \frac{
\int_M \left(
\frac{4(n-1)}{n-2} \:  |\nabla u|^2 \: + \:
R \: u^2 \right) \dvol_M
}{
\left( \int_M u^{\frac{2n}{n-2}} \: \dvol_M \right)^{\frac{n-2}{n}}}.
\end{equation}
As
\begin{equation}
\frac{ \int_M \left(
\frac{4(n-1)}{n-2} \:  |\nabla u|^2 \: + \:
R \: u^2 \right) \dvol_M
}{
\left( \int_M u^{\frac{2n}{n-2}} \: \dvol_M \right)^{\frac{n-2}{n}}}
\: \ge \: R_{min}(g) \: \frac{
\int_M u^2 \: \dvol_M
}{
\left( \int_M u^{\frac{2n}{n-2}} \: \dvol_M \right)^{\frac{n-2}{n}}}
\end{equation}
and $R_{min}(g) \le 0$,
Holder's inequality implies that
$\widetilde{R} \: V(\widetilde{g})^{\frac{2}{n}} \: \ge \:
{R}_{min}(g) \: V({g})^{\frac{2}{n}}$.
It follows that 
$\sigma(M) \: \ge \: \sup_{g} R_{min}(g) \: V(g)^{\frac2n}$.

The next proposition
answers conjectures of Anderson \cite{Anderson3}.

\begin{proposition} \label{propYamabe}
Let $M$ be a closed connected oriented $3$-manifold. \\
(a) If $\sigma(M) > 0$ then $M$ is
diffeomorphic to a connected sum of a finite number of
$S^1 \times S^2$'s and metric quotients of the round $S^3$.
Conversely, each such manifold has $\sigma(M) > 0$. \\
(b) $M$ is a graph manifold if and only if $\sigma(M)
\ge 0$. \\
(c) If $\sigma(M) < 0$
 then
$\left( - \: \frac23 \: \sigma(M)
\right)^{\frac32}$ is the minimum
of the numbers
$V$ with the following property : $M$ can be decomposed as a connected sum
of a finite collection of $S^1 \times S^2$'s, metric quotients of the round
$S^3$ and some other components, the union of which is denoted by
$M^\prime$, and there exists a (possibly disconnected) complete
finite-volume
manifold $N$ with constant sectional curvature $- \: \frac14$ and
volume $V$ which can be embedded in $M^\prime$ so that the
complement $M^\prime - N$ (if nonempty) is a graph manifold.  

Moreover, if
$\vol(N) \: = \: \left( - \: \frac23 \: \sigma(M)
\right)^{\frac32}$
then the
cusps of $N$ (if any) are incompressible in $M^\prime$. 
\end{proposition}

\begin{proof}
If $\sigma(M) > 0$ then $M$ has a metric $g$ of positive scalar curvature.
From Lemmas \ref{pscempty} and \ref{finitetime}, $M$ is a connected sum of $S^1 \times S^2$'s and
metric quotients of the round $S^3$.
Conversely, if $M$ is a connected sum of $S^1 \times S^2$'s and
metric quotients of the round $S^3$ then $M$ admits a metric $g$ of
positive scalar curvature and so $\sigma(M) > 0$.

Now suppose that $\sigma(M) \le 0$.
If $M$ is a graph manifold then $M$ volume-collapses with bounded
curvature, so (\ref{testt}) implies that $\sigma(M) = 0$.

Suppose that $M$ is not a graph manifold. Suppose that we have
a given decomposition of $M$ 
as a connected sum
of a finite collection of $S^1 \times S^2$'s, metric quotients of the round
$S^3$ and some other components, the union of which is denoted by
$M^\prime$, and there exists a (possibly disconnected) finite-volume complete
manifold $N$ with constant sectional curvature $- \: \frac14$ 
which can be embedded in $M^\prime$ so that the
complement (if nonempty) is a graph manifold.
Let
$V_{hyp}$ denote the hyperbolic volume of $N$.  We do not
assume that the cusps of $N$ are incompressible in $M^\prime$.
For any $\epsilon > 0$, we claim that
there is a metric $g_\epsilon$ on $M$ with 
$R \: \ge \: - \: 6 \cdot \frac14 \: - \: \epsilon$ and volume
$V(g_\epsilon) \: \le \: V_{hyp} + \epsilon$. This comes from
collapsing the graph manifold pieces, along with the
fact that the connected sum operation can be performed
while  decreasing the scalar curvature arbitrarily little and increasing
the volume arbitrarily little.
Then
$R_{min}(g_\epsilon) \: V(g_\epsilon)^{\frac23}
 \: \ge \: - \: \frac32 \: V_{hyp}^{2/3} \: - \: \const \: 
\epsilon$. Thus
$\sigma(M) \: \ge \:
- \: \frac32 \: V_{hyp}^{2/3}$.

Let 
$\widehat{V}$
denote the minimum of $V_{hyp}$ over
all such decompositions of $M$. (As the set of
volumes of complete finite-volume $3$-manifolds with constant curvature
$- \: \frac14$ is well-ordered, there is a minimum.)
Then
$\sigma(M) \: \ge \:
- \: \frac32 \: \widehat{V}^{2/3}$.

Next, take an arbitrary metric
$g_0$ on $M$ and consider the Ricci flow $g(t)$ with initial metric $g_0$.
From Sections \ref{hyprig} and  \ref{II.7.4}, there is a nonempty
manifold $N$ with a complete finite-volume metric
of constant curvature $- \: \frac14$
so that for large $t$, there is a decomposition $\M^+_t = M_1(t) \cup M_2(t)$ of the time-$t$ manifold,
where $M_1(t)$ is a graph manifold and $(M_2(t), \frac{1}{t} \: g(t) \big|_{M_2(t)})$
is close to a large
piece of $N$. In terms of condition (c) of Proposition \ref{propYamabe}, we will think
of $M^\prime$ as being $\M^+_t$.
Because of the presence of $N$, we know that $t \: R_{min}(t) \: \le \: - \:
\frac32 \: + \: \epsilon(t)$ and
$V(t) \: \ge \:
t^{2/3} \: V_{hyp}(N) \: - \: \epsilon(t)$ for a function $\epsilon(t)$ 
with $\lim_{t \rightarrow \infty}
\epsilon(t) \: = \: 0$. The monotonicity of
$R_{min}(t) \: V^{2/3}(t)$, even through surgeries, implies that
\begin{equation}
R_{min}(g_0) \: V^{2/3}(g_0) \: \le \: - \: \frac32 \: V_{hyp}(N)^{2/3}
 \: \le \: - \: \frac32 \: \widehat{V}^{2/3}.
 \end{equation}
Thus $\sigma(M) \: \le \: - \: \frac32 \: \widehat{V}^{2/3}$.

This shows that 
$\sigma(M) \: = \: - \: \frac32 \: \widehat{V}^{2/3}$.
Now take a decomposition of $M$ as in
condition (c) of Proposition \ref{propYamabe}, with
$V_{hyp}(N) \: = \: \widehat{V}$. We claim that the cuspidal $2$-tori of $N$
are incompressible in $M^\prime$. If not 
then there would be a metric  $g$ on $M$ with $R(g) \: \ge \: - \: \frac32$ and
$\vol(g) \: < \: V_{hyp}(N)$ \cite[Pf. of Theorem 2.9]{Anderson}. 
This would contradict the fact that  
$\sigma(M) \: = \: - \: \frac32 \: V_{hyp}(N)^{2/3}$.
\end{proof}

\subsection{The approach using the $\overline{\lambda}$-invariant}

\begin{proposition} (cf. II.8.2) \label{propII.8.2}
Let $M$ be a closed connected oriented $3$-manifold. \\
(a) If $M$ admits a metric $g$ with $\lambda(g) > 0$ then it is
diffeomorphic to a connected sum of a finite number of
$S^1 \times S^2$'s and metric quotients of the round $S^3$.
Conversely, each such manifold admits a metric $g$ with 
$\lambda(g) > 0$. \\
(b) Suppose that $M$ does not admit any metric $g$ with
$\lambda(g) > 0$. Let $\overline{\lambda}$ denote the
supremum of $\lambda(g) V(g)^{\frac23}$ over all metrics $g$ on $M$.
Then $M$ is a graph manifold if and only if $\overline{\lambda}
= 0$. \\
(c) Suppose that $M$ does not admit any metric $g$ with
$\lambda(g) > 0$, and $\overline{\lambda} < 0$.
Then
$\left( - \: \frac23 \: \overline{\lambda}
\right)^{\frac32}$ is the minimum
of the numbers
$V$ with the following property : $M$ can be decomposed as a connected sum
of a finite collection of $S^1 \times S^2$'s, metric quotients of the round
$S^3$ and some other components, the union of which is denoted by
$M^\prime$, and there exists a (possibly disconnected) complete
manifold $N$ with constant sectional curvature $- \: \frac14$ and
volume $V$ which can be embedded in $M^\prime$ so that the
complement $M^\prime - N$ (if nonempty) is a graph manifold.  

Moreover, if
$\vol(N) \: = \: \left( - \: \frac23 \: \overline{\lambda}
\right)^{\frac32}$
then the cusps of $N$ (if any) are incompressible in $M^\prime$. 
\end{proposition}

\begin{proof}
We first give the argument for Proposition \ref{propII.8.2} under the 
pretense that all Ricci flows are smooth, except for possible
extinction of components. (Of course this is not the
case, but it will allow us to present the main idea of the proof.)

If $\lambda(g) > 0$ for some metric $g$ then from (\ref{ll}), the Ricci flow
starting from $g$ will become extinct within time
$\frac{3}{2\lambda(g)}$. Hence Lemma \ref{finitetime} applies.
Conversely, if $M$ is a connected sum of $S^1 \times S^2$'s and
metric quotients of the round $S^3$ then $M$ admits a metric $g$ of
positive scalar curvature.  From (\ref{test}), $\lambda(g) > 0$. 

Now suppose that $M$ does
not admit any metric $g$ with $\lambda(g) > 0$.
If $M$ is a graph manifold then $M$ volume-collapses with bounded
curvature, so (\ref{test}) implies that $\overline{\lambda} = 0$.

Suppose that $M$ is not a graph manifold. Suppose that we have
a given decomposition of $M$ 
as a connected sum
of a finite collection of $S^1 \times S^2$'s, metric quotients of the round
$S^3$ and some other components, the union of which is denoted by
$M^\prime$, and there exists a (possibly disconnected) complete
manifold $N$ with constant sectional curvature $- \: \frac14$ 
which can be embedded in $M^\prime$ so that the
complement (if nonempty) is a graph manifold.
Let
$V_{hyp}$ denote the hyperbolic volume of $N$.  We do not
assume that the cusps of $N$ are incompressible in $M^\prime$.
For any $\epsilon > 0$, we claim that
there is a metric $g_\epsilon$ on $M$ with 
$R \: \ge \: - \: 6 \cdot \frac14 \: - \: \epsilon$ and volume
$V(g_\epsilon) \: \le \: V_{hyp} + \epsilon$. This comes from
collapsing the graph manifold pieces, along with the
fact that the connected sum operation can be performed
while  decreasing the scalar curvature arbitrarily little and increasing
the volume arbitrarily little.
Then (\ref{test}) implies that
$\lambda(g_\epsilon) \: V(g_\epsilon)^{\frac23}
 \: \ge \: - \: \frac32 \: V_{hyp}^{2/3} \: - \: \const \: 
\epsilon$. Thus
$\overline{\lambda} \: \ge \:
- \: \frac32 \: V_{hyp}^{2/3}$.

Let $\widehat{V}$ denote the minimum of $V_{hyp}$ over
all such decompositions of $M$. (As the set of
volumes of complete finite-volume $3$-manifolds with constant curvature
$- \: \frac14$ is well-ordered, there is a minimum.)
Then
$\overline{\lambda} \: \ge \:
- \: \frac32 \: \widehat{V}^{2/3}$.

Next, take an arbitrary metric
$g_0$ on $M$ and consider the Ricci flow $g(t)$ with initial metric $g_0$.
From Sections \ref{hyprig} and \ref{II.7.4}, there is a
nonempty manifold $N$ with a finite-volume complete metric
of constant curvature $- \: \frac14$
so that for large $t$, there is a decomposition of the time-$t$ manifold
$\M^+_t = M_1(t) \cup M_2(t)$
where $M_1(t)$ is a graph manifold and $(M_2(t), \frac{1}{t} \: g(t) \big|_{M_2(t)})$
is close to a large
piece of $N$.
As $N$ has finite volume, a constant function on $N$ is square-integrable
and so  $\inf \spec(- \triangle_N) \: = \: 0$.
Equivalently,
\begin{equation}
\inf_{f \in C^\infty_c(N), f \neq 0}
\frac{\int_N |\nabla f|^2 \: \dvol_N}{\int_N f^2 \: \dvol_N} \: = \: 0.
\end{equation}
Taking an appropriate test function $\Phi$ on $M$ with compact
support in $M_2(t)$
gives $t \: \lambda(t) \: \le \: - \:
\frac32 \: + \: \epsilon_1(t)$, with $\lim_{t \rightarrow \infty}
\epsilon_1(t) \: = \: 0$. In terms of condition (c) of Proposition \ref{propII.8.2},
we will think of $M^\prime$ as being $\M^+_t$.
From the presence of $N$, we know that $V(t) \: \ge \:
t^{2/3} \: V_{hyp}(N) \: - \: \epsilon_2(t)$, with $\lim_{t \rightarrow \infty}
\epsilon_2(t) \: = \: 0$. As we are assuming that the Ricci flow is nonsingular,
the monotonicity of
$\lambda(t) \: V^{2/3}(t)$ implies that
\begin{equation}
\lambda(g_0) \: V^{2/3}(g_0) \: \le \: - \: \frac32 \: V_{hyp}(N)^{2/3}
 \: \le \: - \: \frac32 \: \widehat{V}^{2/3}.
 \end{equation}
Thus $\overline{\lambda} \: \le \: - \: \frac32 \: \widehat{V}^{2/3}$.

This shows that 
$\overline{\lambda} \: = \: - \: \frac32 \: \widehat{V}^{2/3}$.
Now take a decomposition of $M$ as in condition (c) of Proposition \ref{propII.8.2},
with
$V_{hyp}(N) \: = \: \widehat{V}$. We claim that the cuspidal $2$-tori of $N$ are
incompressible in $M^\prime$. If not 
then there would be a metric $g$ on $M$ with $R(g) \: \ge \: - \: \frac32$ and
$\vol(g) \: < \: V_{hyp}(N)$ \cite[Pf. of Theorem 2.9]{Anderson}. 
Using (\ref{test}), one
would obtain a contradiction to the fact that  
$\overline{\lambda} \: = \: - \: \frac32 \: V_{hyp}(N)^{2/3}$.
 
To handle the behaviour of $\lambda(t) \: V^{\frac23}(t)$ under
Ricci flows with surgery, we 
first state a couple of general facts about
Schr\"odinger operators.

\begin{lemma} \label{newprop1}
Given a closed Riemannian manifold $M$, let $X$ be a codimension-$0$
submanifold-with-boundary of $M$. Given $R \in C^\infty(M)$, let
$\lambda_M$ be the lowest eigenvalue of $-4\triangle \: + \: R$
on $M$, with corresponding eigenfunction $\psi$. 
Let $\lambda_X$ be the lowest eigenvalue of the corresponding
operator on $X$, with Dirichlet boundary conditions, and similarly
for $\lambda_{M - \Int(X)}$. Then for all
$\eta \in C^\infty_c(\Int(X))$, we have
\begin{equation} \label{eig1}
\lambda_M \: \le \: \min(\lambda_X, \lambda_{M - \Int(X)})
\end{equation}
and
\begin{equation} \label{eig2}
\lambda_X \: \le \: \lambda_M \: + \: 4 \:
\frac{\int_M |\nabla \eta|^2 \: \psi^2 \: dV}{\int_M \eta^2 \: 
\psi^2 \: dV}.
\end{equation}
\end{lemma}
\begin{proof}
Equation (\ref{eig1}) follows from Dirichlet-Neumann bracketing
\cite[Chapter XIII.15]{Reed-Simon}. To prove (\ref{eig2}),
$\eta \psi$ is supported in $\Int(X)$ and so
\begin{align}
\lambda_X \: & \le \: \frac{\int_M \left( 4 |\nabla(\eta \psi)|^2 \: + \:
R \: \eta^2 \psi^2 \right) \: dV}{\int_M \eta^2 \psi^2 \: dV} \\
& = \: \frac{\int_M \left( 4 |\nabla \eta|^2 \psi^2 \: + \:
8 \langle \nabla\eta, \nabla \psi \rangle \eta \psi \: + \:
4 \eta^2 |\nabla \psi|^2 \: + \:
R \: \eta^2 \psi^2 \right) \: dV}{\int_M \eta^2 \psi^2 \: dV}. \notag
\end{align}
As $- 4 \triangle \psi \: + \: R \psi \: = \: \lambda_M \psi$, we have

\begin{align}
\lambda_M \int_M \eta^2 \psi^2 dV \: & = \:
- \: 4 \: \int_M \eta^2 \psi \triangle \psi \: dV \: + \: 
\int_M R \eta^2 \psi^2 \: dV \\
& = \:
4 \: \int_M \langle \nabla(\eta^2 \psi), \nabla \psi \rangle \: dV \: 
+ \: \int_M R \eta^2 \psi^2 \: dV \notag \\
& = \: 4 \: \int_M \eta^2 \: |\nabla \psi|^2 \: dV \: \: + \:  
8 \int_M \langle \nabla\eta, \nabla \psi \rangle \eta \psi \: dV
\: + \: \int_M R \eta^2 \psi^2 \: dV. \notag
\end{align}
Equation (\ref{eig2}) follows.
\end{proof}

The next result is an Agmon-type estimate. 

\begin{lemma} \label{newprop2}
With the notation of Lemma \ref{newprop1}, 
given a nonnegative function $\phi \in C^\infty(M)$,
suppose that $f \in C^\infty(M)$ satisfies
\begin{equation}
4 |\nabla f|^2 \: \le \: R \: - \: \lambda_M \: - \: c
\end{equation}
on $\supp(\phi)$, for some $c > 0$. Then
\begin{equation} \label{eig3}
\|e^f \phi \psi\|_2 \: \le \: 4 \: c^{-1} \: \left( \|e^f \triangle \phi\|_\infty
\: + \: \|e^f \nabla \phi\|_\infty (\lambda_M \: - \: \min R)^{1/2} \right)
\: \|\psi\|_2.
\end{equation}
\end{lemma}
\begin{proof}
Put $H \: = \: - \: 4 \triangle \: + \: R$. 
By assumption, 
\begin{equation}
\phi (R \: - \: 4 |\nabla f|^2 \: - \: \lambda_M) \: \phi \: \ge \: c \: \phi^2
\end{equation}
and so there is an
 inequality of operators on $L^2(M)$ :
\begin{equation}
\phi (H \: - \: 4 |\nabla f|^2 \: - \: \lambda_M) \: \phi \: = \:
4 \phi d^* d \phi \: + \: 
\phi (R \: - \: 4 |\nabla f|^2 \: - \: \lambda_M) \: \phi \:
\ge \: c \: \phi^2.
\end{equation}
In particular,
\begin{equation} \label{part1}
\int_M e^f \psi\phi \: (H \: - \: 4 |\nabla f|^2 \: - \: \lambda_M) \: \phi
e^f \psi dV \: \ge \: c \int_M \phi^2 e^{2f} \psi^2 \: dV.
\end{equation}

For $\rho \in C^\infty(M)$,
\begin{equation} \label{Hf}
e^f H (e^{-f} \rho) \: = \: H \rho \: + \: 4 \: \nabla \cdot 
((\nabla f) \rho) \: + \: 4 \: \langle \nabla f, \nabla \rho \rangle
\: - \: 4 \: |\nabla f|^2 \rho
\end{equation}
and so
\begin{equation}
\int_M \rho e^f H(e^{-f} \rho) \: dV \: = \:
\int_M \rho (H \: - \: 4 \: |\nabla f|^2) \rho \: dV.
\end{equation}
Taking $\rho \: = \: e^f \phi \psi$ gives
\begin{equation} \label{part2}
\int_M e^{2f} \phi \psi H(\phi \psi) \: dV \: = \:
\int_M e^f \phi \psi (H \: - \: 4 \: 
|\nabla f|^2)  e^{f} \phi \psi \: dV.
\end{equation}
From (\ref{part1}) and (\ref{part2}),
\begin{equation}
c \|e^f \phi \psi \|_2^2 \: \le \:
\int_M e^{2f} \phi \psi (H \: - \: \lambda_M) (\phi \psi) \: dV \: = \:
\int_M e^{2f} \phi \psi [H, \phi] \psi \: dV \: = \:
\langle e^{f} \phi \psi, e^{f} [H, \phi] \psi \rangle_2.
\end{equation}
Thus
\begin{equation}
c \|e^f \phi \psi \|_2 \: \le \:
\|e^{f} [H, \phi] \psi\|_2.
\end{equation}
Now
\begin{equation}
e^{f} [H, \phi] \psi \: = \: - \: 4 \: e^f (\triangle \phi) \psi \: - \:
8 \: e^f \langle \nabla\phi, \nabla \psi \rangle.
\end{equation}
Then 
\begin{equation}
\|e^{f} [H, \phi] \psi\|_2 \: \le \:
4 \: \|e^f \triangle \phi\|_\infty \: \|\psi\|_2 \: + \: 8 \:
\|e^f \nabla\phi\|_\infty \: \|\nabla \psi\|_2.
\end{equation}
Finally,
\begin{equation}
4 \: \|\nabla \psi\|_2^2 \: = \: \int_M (\lambda_M \: - \: R) \: \psi^2 \:
dV
\end{equation}
and so
\begin{equation}
2 \|\nabla \psi\|_2 \: \le \: (\lambda_M \: - \: \min R)^{1/2} \:
\|\psi\|_2.
\end{equation}
This proves the lemma.
\end{proof}

Clearly Lemma \ref{newprop2} is also true if 
$f$ is just assumed to be Lipschitz-regular.

We now apply Lemmas \ref{newprop1} and \ref{newprop2} 
to a Ricci flow with surgery.
A singularity caused by extinction of a component will not be a
problem, so
let $T_0$ be a surgery time  and let
$M_+ = \M_{T_0}^+$ be the postsurgery manifold.  We will write
$\lambda^+$ instead of $\lambda_{M_+}$. Let $M_{cap} = 
\M_{T_0}^+ - (\M_{T_0}^+ \cap \M_{T_0}^-)$ be the added
caps and put $X \: = \: \overline{M_+ - M_{cap}} = 
\overline{\M_{T_0}^+ \cap \M_{T_0}^-}$.
For simplicity,
let us assume that $M_{cap}$ has a single component; the argument in
the general case is similar.
From the nature of the surgery procedure, the surgery is done in an
$\epsilon$-horn extending from $\Omega_\rho$, where
$\rho \: = \: \delta(T_0) r(T_0)$. 
In fact, because of the canonical
neighborhood assumption, we can extend the $\epsilon$-horn inward
until $R \sim r(T_0)^{-2}$. Appplying (\ref{Rayleigh}) with a test function
supported in an $\epsilon$-tube near this inner boundary, it follows that
$\lambda^+ \: \le \: c^\prime \: r(T_0)^{-2}$ for some universal constant 
$c^\prime >> 1$.

In what follows we take $\delta(T_0)$ to be small.
As $R$ is much greater than $r(T_0)^{-2}$ on $M_{cap}$, it follows that
$\lambda_{M - \Int(X)}$ is much greater than $r(T_0)^{-2}$. Then from
(\ref{eig1}), $\lambda^+ \: \le \: \lambda_X$. We can apply
(\ref{eig2}) to get an inequality the other way.
We take the function $\eta$ to interpolate from being $1$ outside of the
$h(T_0)$-neighborhood 
$N_h M_{cap}$ of $M_{cap}$, to being $0$ on $M_{cap}$. 
In terms of the normalized eigenfunction $\psi$ on $M_+$, 
this gives a bound of the form
\begin{equation} \label{bound11}
\lambda_X \: \le \: \lambda^+ \: + \: \const \: h(T_0)^{-2} \:
\frac{\int_{N_h M_{cap}} \psi^2 \: dV}{1 - \int_{N_h M_{cap}} \psi^2 \: dV}.
\end{equation}

We now wish to show that $\int_{N_h M_{cap}} \psi^2 \: dV$ is small.
For this we apply Lemma \ref{newprop2} with
$c \: = \: c^\prime \: r(T_0)^{-2}$. Take an $\epsilon$-tube $U$,
in the $\epsilon$-horn,
whose center has scalar curvature roughly $200 \: c^\prime \: r(T_0)^{-2}$
and which is the closest tube to the cap with this property.
Let $x \: : \: U \rightarrow
(- \: \epsilon^{-1}, \epsilon^{-1})$ 
be the longitudinal parametrization of the tube,
which we take to be increasing in the direction of the surgery cap. 
Let $\Phi \: : \: (-1, 1) \rightarrow [0,1]$
be a fixed nondecreasing smooth function which is zero on $(-1, 1/4)$ and one on
$(1/2, 1)$. Put $\phi \: = \: \Phi \circ x$ on $U$.
Extend $\phi$ to $M_+$ by making it zero to the left of $U$ and one to the right of $U$,
where ``right of $U$'' means the connected component of $M_+ - U$ containing the
surgery cap.
Dimensionally, $|\nabla \phi|_\infty \: \le \: \const \: 
r(T_0)^{-1}$ and $|\triangle \phi|_\infty \: \le \: \const \:  r(T_0)^{-2}$.
Define a function $f$ to the right of $x^{-1}(0)$ by setting it to be the distance from $x^{-1}(0)$
with respect to the metric $\frac14 \: (R - \lambda^+ - c) \: g_{M_+}$. 
(Note that to the right of $x^{-1}(0)$, we have
$R \: \ge \: 200 c^\prime r(T_0)^{-2} \: \ge \: \lambda^+ + c$.) 
Then
equation (\ref{eig3}) gives
$|e^f \phi \psi|_2 \: \le \: \const (T_0)$. The point is that $\const (T_0)$ is independent
of the (small) surgery parameter $\delta(T_0)$. 

Hence
\begin{equation} \label{bound12}
\int_{N_h M_{cap}} \psi^2 \: dV \: \le \:
\left( \sup_{N_h M_{cap}} e^{-2f} \right) \: 
\int_{N_h M_{cap}} \: e^{2f} \psi^2 \: dV \: \le \:
\const \: \sup_{N_h M_{cap}} e^{-2f}.
\end{equation}
To estimate $\sup_{N_h M_{cap}} e^{-2f}$, we use the fact that
the $\epsilon$-horn consists of a
sequence of $\epsilon$-tubes stacked together.  In the 
region of $M_+$ from $x^{-1}(0)$ to the surgery cap, the scalar curvature
ranges from roughly $200 \: c^\prime \: r(T_0)^{-2}$
to $h(T_0)^{-2}$.
On a given
$\epsilon$-tube, if $\epsilon$ is sufficiently small then
the ratio of the scalar curvatures between the
two ends is bounded by $e$. Hence in going 
from $x^{-1}(0)$ to the surgery cap,
one must cross at least $N$ disjoint $\epsilon$-tubes,
with $e^N \: = \: \frac{1}{200 c^\prime} \: r(T_0)^2 \: h(T_0)^{-2}$. 
Traversing a given $\epsilon$-tube (say of
radius $r^\prime$) in going towards the surgery cap, $f$ increases by
roughly $\const \: \int_{- \epsilon^{-1}r^\prime}^{\epsilon^{-1}r^\prime}
(r^\prime)^{-1} \: ds$, which is $\const \epsilon^{-1}$. Hence near
the surgery cap, we have
\begin{equation}
\sup_{N_h M_{cap}} e^{-2f} \: \le \: \const \: e^{-\const N \epsilon^{-1}}
\: = \: \const \: (r(T_0)^2 \: h(T_0)^{-2})^{-\const \epsilon^{-1}}.
\end{equation}
Combining this with (\ref{bound11}) and (\ref{bound12}), we obtain
\begin{equation}
\lambda_X \: \le \: \lambda^+ \: + \: \const \: h(T_0)^{-2} \: 
(r(T_0)^2 \: h(T_0)^{-2})^{- \const \epsilon^{-1}}.  
\end{equation}

By making a single redefinition of $\epsilon$,
we can ensure that
$\lambda_X \: \le \: \lambda^+ \: + \: \const \: h(T_0)^{4}$. The last
constant will depend on $r(T_0)$ but is independent of $\delta(T_0)$.
Thus if $\delta(T_0)$ is small enough, we can ensure that
$|\lambda_X \: - \: \lambda^+|$ is small in comparison to the
volume change $V^-(T_0) - V^+(T_0)$, which is comparable to $h(T_0)^3$.

If $\lambda^-$ is the smallest eigenvalue of $- \:4 \: \triangle \: + \: R$ 
on the presurgery manifold $\M_t$, for $t$ slightly less than $T_0$,
then we can estimate $|\lambda_X \: - \: \lambda^-|$ in a similar way.
Hence for an arbitrary positive continuous function $\xi(t)$, 
we can make the parameters $\overline{\delta}_j$ of
Proposition \ref{surgeryflowexists} small enough to ensure that
\begin{equation}
|\lambda^+(T_0) - \lambda^-(T_0)| \: \le \: \xi(T_0) \: 
(V^-(T_0) - V^+(T_0))
\end{equation}
for a surgery at time $T_0$.

We now redo the argument for the proposition, as given above in the
surgery-free case, in the presence of surgeries.
Suppose first that $\lambda(g_0) > 0$ for some metric $g_0$ on $M$. 
After possible rescaling, we can assume that $g_0$ is the initial
condition for a Ricci flow with surgery $(\M, g(\cdot))$, with normalized initial condition.
Using the lower scalar curvature bound of Lemma \ref{Rev} and the Ricci
flow equation,
the volume on the time interval $[0, \frac{3}{\lambda(0)}]$ has an {\it a priori} upper
bound of the form $\const V(0)$. As a surgery at time $T_0$ removes a 
volume comparable to $h(T_0)^3$, we have
$\sum h(T_0)^3 \: \le \: \const V(0)$, where the sum is over the
surgeries and $T_0$ denotes the
surgery time.  From the above discussion,
the change in $\lambda$ due to the 
surgeries is bounded below by $-\: \const \sum_{T_0} h(T_0)^4$.
Then the decrease in $\lambda$ due to surgeries on the time interval
$[0, \frac{3}{\lambda(0)}]$ is bounded above by
\begin{equation}
\const \sum_{T_0 \in [0, \frac{3}{\lambda(0)}]} h(T_0)^4 \: \le \:  \const \:
\left( \sup_{t \in [0, \frac{3}{\lambda(0)}]} h(t) \right) \: V(0).
\end{equation}
By choosing the function $\delta(t)$ to be sufficiently small,
the decrease in $\lambda$ due to surgeries
is not enough to prevent the blowup of $\lambda$ on the time
interval $[0, \frac{3}{\lambda(0)}]$ coming from the increase of
$\lambda$ between the surgeries. Hence the solution goes extinct.

Now suppose that $M$ does not admit a metric $g$ with $\lambda(g) > 0$.
Again, if $M$ is a graph manifold then $\overline{\lambda} \: = \: 0$.

Suppose that $M$ is not a graph manifold.
As before, $\overline{\lambda} \: \ge \: - \: \frac32 \: \widehat{V}^{2/3}$.
Given an initial metric $g_0$,
we wish to show that by choosing the function $\delta(t)$ 
small enough we can make the function $\lambda(t) V^{2/3}(t)$
arbitrarily close to being nondecreasing. To see this, we consider
the effect of a surgery on $\lambda(t) V^{2/3}(t)$. Upon performing a surgery
the volume decreases, which in itself cannot decrease
$\lambda(t) V^{2/3}(t)$. (We are using the fact that $\lambda(t)$ is
nonpositive.) 
Then from the above discussion, the
change in $\lambda(t) V^{2/3}(t)$ due from a surgery at time $T_0$,
is bounded below by $- \: \xi(T_0) \: (V^-(T_0) - V^+(T_0)) \: 
V^-(T_0)^{2/3}$.
With normalized initial conditions,
we have an {\it a priori} upper bound on $V(t)$ in terms of $V(0)$ and $t$.
Over any time interval $[T_1, T_2]$, we must have
\begin{equation}
\sum_{T_0 \in [T_1, T_2]} \left( V^-(T_0) - V^+(T_0) \right) \: \le \: 
\sup_{t \in [T_1,T_2]} V(t),
\end{equation}
where the sum is over the surgeries in the
interval $[T_1,T_2]$.
Then along with the monotonicity of $\lambda(t) V^{2/3}(t)$ in between the
 surgery times,  by choosing the function
$\delta(t)$ appropriately we can ensure that 
for any 
$\sigma > 0$
there is a Ricci flow with $(r, \delta)$-cutoff
starting from $g_0$ so that
$\lambda(g_0) \: V^{2/3}(g_0) \: \le \: \lambda(t) \: V^{2/3}(t) \: + \:
\sigma$ for all $t$.
It follows that
$\overline{\lambda} \: \le \: - \: \frac32 \: \widehat{V}^{2/3}$.

This shows that $\overline{\lambda} \: = \: - \: \frac32 \: \widehat{V}^{2/3}$.
The same argument as before shows that if we have a decomposition
with the hyperbolic volume of $N$ equal to  $\widehat{V}$ then the cusps of $N$ 
(if any) are incompressible in $M^\prime$. 
\end{proof}

\begin{remark}
It follows that if the three-manifold $M$ does not admit a metric of positive scalar curvature then
$\sigma(M) \: = \: \overline{\lambda}$. In fact, this is true in any dimension
$n \ge 3$ \cite{Akutagawa-Ishida-Lebrun}.
\end{remark}

\appendix
\section{Maximum principles} \label{maxprin}

In this appendix we list some maximum principles and their
consequences.  Our main source is \cite{Chow-Lu-Ni}, where
references to the original literature can be found.

The first type of maximum principle is a {\em weak} maximum
principle which says that under certain conditions, a spatial inequality
on the initial condition implies a time-dependent inequality at later times.

\begin{theorem}
Let $M$ be a closed manifold. Let $\{g(t)\}_{t \in [0,T]}$ be a smooth
one-parameter family of Riemannian metrics on $M$ and let 
$\{X(t)\}_{t \in [0, T]}$ be
a smooth one-parameter family of vector fields on $M$. Let $F \: : \:
\R \times [0, T] \rightarrow \R$ be a Lipschitz function.  Suppose that
$u=u(x,t)$ is $C^2$-regular in $x$, $C^1$-regular in $t$ and
\begin{equation}
\frac{\partial u}{\partial t} \: \le \: \triangle_{g(t)} u \: + \: X(t) u \: + \: F(u, t).
\end{equation}
Let $\phi \: : \: [0,T] \rightarrow \R$ be the solution of 
$\frac{d\phi}{dt} \: = \: F(\phi(t), t)$ with $\phi(0) = \alpha$. If
$u(\cdot, 0) \le \alpha$ then $u(\cdot, t) \le \phi(t)$ for all $t \in [0, T]$.
\end{theorem}

There are various noncompact versions of the weak maximum principle.
We state one here.

\begin{theorem}
Let $(M, g(\cdot))$ be a complete Ricci flow solution on the interval
$[0,T]$ with uniformly bounded curvature.  If $u=u(x,t)$ is a
$W^{1,2}_{loc}$ function that weakly satisfies
$\frac{\partial u}{\partial t} \: \le \: \triangle_{g(t)} u$, with
$u(\cdot, 0) \le 0$ and
\begin{equation}
\int_0^T \int_M e^{- \: c \: d_t^2(x, x_0)} \: u^2(x,s) \: dV(x) \: ds < \infty
\end{equation}
for some $c > 0$, then $u(\cdot,t) \le 0$ for all $t \in [0,T]$.
\end{theorem}

A {\em strong} maximum principle says that under certain conditions,
a strict inequality at a given time implies strict inequality at later
times and also slightly earlier times.  It does not require complete metrics.

\begin{theorem}
Let $M$ be a connected manifold. Let $\{g(t)\}_{t \in [0,T]}$ be a smooth
one-parameter family of Riemannian metrics on $M$ and let 
$\{X(t)\}_{t \in [0, T]}$ be
a smooth one-parameter family of vector fields on $M$. Let $F \: : \:
\R \times [0, T] \rightarrow \R$ be a Lipschitz function.  Suppose that
$u=u(x,t)$ is $C^2$-regular in $x$, $C^1$-regular in $t$ and
\begin{equation}
\frac{\partial u}{\partial t} \: \le \: \triangle_{g(t)} u \: + \: X(t) u \: + \: F(u, t).
\end{equation}
Let $\phi \: : \: [0,T] \rightarrow \R$ be a solution of 
$\frac{d\phi}{dt} \: = \: F(\phi(t), t)$. If
$u(\cdot, t) \le \phi(t)$ for all $t \in [0, T]$ and
$u(x_0, t_0) < \phi(t)$ for some $x_0 \in M$ and $t_0 \in (0,T]$ then
there is some $\epsilon > 0$ so that $u(\cdot, t) < \phi(t)$ for
$t \in (t_0 - \epsilon, T]$.
\end{theorem}

A consequence of the strong maximum principle is a statement about
restricted holonomy for Ricci flow solutions with nonnegative curvature
operator $\Rm$.

\begin{theorem} \label{splitting}
Let $M$ be a connected manifold. Let $\{g(t)\}_{t \in [0,T]}$ be a smooth
one-parameter family of Riemannian metrics on $M$ with nonnegative
curvature operator that satisfy the Ricci flow equation.  Then for each $t \in (0, T]$, the
image $\Image(\Rm_{g(t)})$ of the curvature operator
is a smooth subbundle of $\Lambda^2(T^* M)$ which is invariant under
spatial parallel translation. There is a sequence of times
$0 = t_0 < t_1 < \ldots < t_k = T$ such that for each
$1 \le i \le k$,
$\Image(\Rm_{g(t)})$ is a Lie subalgebra of $\Lambda^2(T^*_m M)
\cong o(n)$ that is independent of $t$ for $t \in (t_{i-1}, t_i]$.
Furthermore, $\Image(\Rm_{g(t_i)}) \subset \Image(\Rm_{g(t_{i+1})})$.
\end{theorem}

In particular, under the hypotheses of Theorem \ref{splitting}, a
local isometric splitting at a given time implies a local
isometric splitting at earlier times.

\section{$\phi$-almost nonnegative curvature}
\label{phiappendix}

In three dimensions,
the Ricci flow equation implies that
\begin{equation} \label{3dR}
\frac{dR}{dt} \: = \:
\triangle R \: + \: \frac{2}{3} R^2 \: + 
2 \: |R_{ij} - \frac{R}{3} g_{ij}|^2 .  
\end{equation}
The maximum principle of Appendix \ref{maxprin} implies that
if $(M, g(\cdot))$ is a Ricci flow solution defined for
$t \in [0, T)$, with complete time slices and bounded curvature
on compact time intervals, then
\begin{equation} \label{Rlowerbound}
(\inf R)(t) \: \ge \: 
\frac{(\inf R)(0)}
{1 - \frac23 t (\inf R)(0)}.
\end{equation}
In particular,
$tR(\cdot, t) \: > \: - \: \frac{3}{2}$ for all $t \ge 0$
(compare \cite[Section 2]{Hamiltonn}).

Recall that the curvature operator is an operator on $2$-forms.
We follow the usual Ricci flow convention that if a manifold
has constant sectional curvature $k$ then its curvature operator
is multiplication by $2k$. In general, the trace of the
curvature operator equals the scalar curvature.

In three dimensions, having nonnegative
curvature operator is equivalent to having nonnegative
sectional curvature.  Each eigenvalue of the curvature operator
is twice a sectional curvature.

Hamilton-Ivey pinching, as given
in \cite[Theorem 4.1]{Hamiltonn}, says the following.

Assume that at $t = 0$ the eigenvalues $\lambda_1 \le \lambda_2 \le
\lambda_3$ of the
curvature operator at each point satisfy $\lambda_1 \ge -1$.
(One can
always achieve this by rescaling. Note that it implies
$R(\cdot, t) \ge - \frac32 \frac{1}{t+\frac12}$.)
Given a point $(x, t)$, put $X \: = \: - \lambda_1$.
If $X > 0$ then
\begin{equation}
R(x, t) \: \ge \: X \left( \ln X + \ln(1+t) - 3 \right),
\end{equation}
or equivalently,
\begin{equation} \label{estimate}
tR(x, t) \: \ge \: tX \left( \ln(tX) + \ln(\frac{1+t}{t}) - 3 \right).
\end{equation}

\begin{definition} \label{HamIvey}
Given $t\geq 0$, a Riemannian $3$-manifold $(M,g)$ satisfies the 
{\em  time-$t$ Hamilton-Ivey pinching condition} if for every $x\in M$,
if $\la_1\leq \la_2\leq \la_3$ are the eigenvalues of the curvature operator at $x$,
then either 

$\bullet$ $\la_1\geq 0$, i.e. the curvature is nonnegative at $x$,

or

$\bullet$ If $\la_1 < 0$ and $X \: = \: - \lambda_1$ then
$tR(x)\geq tX \left(\ln (tX) + \ln \left( \frac{1+t}{t} \right)-3
\right)$.
\end{definition}
This condition has the following monotonicity property:

\begin{lemma}
\label{lempinchingmonotonic}
Suppose that $\Rm$ and $\Rm'$  
are $3$-dimensional  curvature operators whose scalar curvatures
and first eigenvalues satisfy
$R'\geq R\geq 0$ and $\la_1'\geq \la_1$.
If $\Rm$ satisfies the time-$t$ Hamilton-Ivey pinching condition then so does 
$\Rm'$.
\end{lemma}
\begin{proof}
We may assume that  $\la_1'<0$ and 
$\log (tX') \: + \: \ln \left( \frac{1+t}{t} \right) \: - \: 3 \: > \: 0$,
since otherwise the condition will be satisfied
(because $R'>0$ by hypothesis).  The function 
\begin{equation}
Y\mapsto tY \left( \ln(t Y) \: + \: \ln \left( \frac{1+t}{t} \right) \: - \: 3 \right)
\end{equation}
is monotone increasing on the interval on which  
$\ln(t Y) \: + \: \ln \left( \frac{1+t}{t} \right) \: - \: 3$
is nonnegative, so 
$tX \left( \ln(t X) \: + \: \ln \left( \frac{1+t}{t} \right) \: - \: 3 \right)
\: \ge \:
tX^\prime \left( \ln(t X^\prime) \: + \: \ln \left( \frac{1+t}{t} \right) \: - \: 3 \right)$.
Hence $\Rm'$ satisfies the pinching condition too.
\end{proof}

The content of the pinching equation is that for any $s \in \R$,
if $t R(\cdot, t) \: \le \:
s$ then there is a lower bound $t \Rm (\cdot, t) \: \ge \: \const
(s, t)$. Of course, this is a vacuous statement
if $s \: \le \: - \: \frac{3}{2}$.

Using equation (\ref{estimate}), we can find a positive
function $\Phi \in C^\infty(\R)$ such that \\
1. $\Phi$ is nondecreasing. \\
2. For $s > 0$, $\frac{\Phi(s)}{s}$ is decreasing. \\
3. For large $s$,
$\Phi(s) \sim \frac{s}{\ln s}$. \\
4. For all $t$,
\begin{equation} \label{phi}
\Rm(\cdot, t) \: \ge \: - \: \Phi(R(\cdot, t)).
\end{equation}

This bound has the most consequence when $s$ is large.

We note that for the original unscaled Ricci flow
solution, the precise
bound that we obtain 
depends on $t_0$ and
the time-zero metric, through its lower curvature bound.

\section{Ricci solitons} \label{Riccisolitons}

Let $\{V(t)\}$ be a time-dependent family of vector fields on a manifold $M$.
The solution to the equation
\begin{equation}
\frac{dg}{dt} \: = \: {\mathcal L}_{V(t)} g
\end{equation}
is 
\begin{equation}
g(t) \: = \: \phi^{-1}(t)^* g(t_0)
\end{equation}
where $\{ \phi(t) \}$ is the $1$-parameter group of diffeomorphisms generated 
by $- V$, normalized by $\phi(t_0) \: = \: \Id$.
(If $M$ is noncompact then we assume that $V$ can be integrated.
The reason for the funny signs is that if a $1$-parameter family of
diffeomorphisms $\eta(t)$ is generated by vector fields $W(t)$ then
${\mathcal L}_{W(t)} \: = \: \eta^{-1}(t)^* \:
\frac{d\eta(t+\epsilon)^*}{d\epsilon} \Big|_{\epsilon = 0} \: = \: - \:
\frac{d\eta^{-1}(t+\epsilon)^*}{d\epsilon} \Big|_{\epsilon = 0} \: \eta(t)^*$.)

The equation for a {\em steady soliton} is
\begin{equation}
2 \Ric \: + \: {\mathcal L}_V g \: = \: 0,
\end{equation}
where $V$ is a time-independent
vector field.
The corresponding Ricci flow is given by
\begin{equation}
g(t) \: = \: \phi^{-1}(t)^* g(t_0),
\end{equation}
where $\{ \phi(t) \}$ is the $1$-parameter group of diffeomorphisms generated 
by $- V$. (Of course, in this case
$\{ \phi^{-1}(t) \}$ is the $1$-parameter group of diffeomorphisms generated 
by $V$.)

A {\em gradient steady soliton}
satisfies the equations
\begin{align} \label{steady}
\frac{\partial g_{ij}}{\partial t} \: & = \: - \: 2 R_{ij} \: = \:
2 \nabla_i \nabla_j f, \\
\frac{\partial f}{\partial t} \: & = \: |\nabla f|^2. \notag
\end{align}
It follows from (\ref{steady}) that
\begin{equation}
\frac{\partial}{\partial t} \left( g^{ij} \: \partial_j f \right) \: = \:
2 \: R^{ij} \: \nabla_j f \: + \: g^{ij} \nabla_j |\nabla f|^2 \: = \:
- 2 \: (\nabla^i \nabla^j f) \: \nabla_j f \: + \: \nabla^i
|\nabla f|^2 \: = \: 0,
\end{equation}
showing that $V \: = \: \nabla f$ is indeed constant in $t$.
The solution to (\ref{steady}) is
\begin{align} \label{steadysoln}
g(t) \: & = \: \phi^{-1}(t)^* g(t_0), \\
f(t) \: & = \: \phi^{-1}(t)^* f(t_0). \notag
\end{align}

Conversely, given a metric $\widehat{g}$ and a function $\widehat{f}$
satisfying 
\begin{equation}
\widehat{R}_{ij} \: + \: \widehat{\nabla}_i \widehat{\nabla}_j \widehat{f}
\: = \: 0,
\end{equation}
put $V \: = \: \widehat{\nabla} \widehat{f}$. 
If we define $g(t)$ and $f(t)$ by 
\begin{align}
g(t) \: & = \: \phi^{-1}(t)^* \widehat{g}, \\
f(t) \: & = \: \phi^{-1}(t)^* \widehat{f} \notag
\end{align}
then they satisfy (\ref{steady}).

A solution to (\ref{steady}) satisfies
\begin{equation}
\frac{\partial f}{\partial t} \: = \: |\nabla f|^2 \: - \: \triangle f
\: - \: R,
\end{equation}
or
\begin{equation}
\frac{\partial}{\partial t} \: e^{-f} \: = \: - \: \triangle e^{-f}
\: + \: R \: e^{-f}.
\end{equation}
This perhaps motivates Perelman's use of the backward heat
equation (\ref{backwards}).

A {\em shrinking soliton} lives on a time interval $(-\infty, T)$.
For convenience, we take $T \: = \: 0$. Then the equation is
\begin{equation}
2 \Ric \: + \: {\mathcal L}_V g \: + \: \frac{g}{t} \: = \: 0.
\end{equation}
The vector field $V = V(t)$ satisfies
$V(t) \: = \: - \: \frac{1}{t} \: V(-1)$.
The corresponding Ricci flow is given by
\begin{equation} \label{shrinkingflow}
g(t) \: = \: - \: t \: \phi^{-1}(t)^* g(-1),
\end{equation}
where $\{ \phi(t) \}$ is the $1$-parameter group of diffeomorphisms generated 
by $- V$, normalized by $\phi(-1) \: = \: \Id$.

A {\em gradient shrinking soliton} satisfies the equations
\begin{align} \label{shrinking}
\frac{\partial g_{ij}}{\partial t} \: & = \: - \: 2 R_{ij} \: = \:
2 \nabla_i \nabla_j f \: + \: \frac{g_{ij}}{t}, \\
\frac{\partial f}{\partial t} \: & = \: |\nabla f|^2. \notag
\end{align}
It follows from (\ref{shrinking}) that
$V \: = \: \nabla f$ satisfies $V(t) \: = \: - \: \frac{1}{t} \: V(-1)$.
The solution to (\ref{shrinking}) is
\begin{align} \label{shrinkingsoln}
g(t) \: & = \: - \: t \: \phi^{-1}(t)^* g(-1), \\
f(t) \: & = \: \phi^{-1}(t)^* f(-1). \notag
\end{align}

Conversely, given a metric $\widehat{g}$ and a function $\widehat{f}$
satisfying 
\begin{equation}
\widehat{R}_{ij} \: + \: \widehat{\nabla}_i \widehat{\nabla}_j \widehat{f} 
\: - \: \frac{1}{2} \:
\widehat{g}
\: = \: 0,
\end{equation}
put $V(t) \: = \: - \: \frac{1}{t} \: \widehat{\nabla} \widehat{f}$. 
If we define $g(t)$ and $f(t)$ by 
\begin{align}
g(t) \: & = \: - \: t \: \phi^{-1}(t)^* \widehat{g}, \\
f(t) \: & = \: \phi^{-1}(t)^* \widehat{f} \notag
\end{align}
then they satisfy (\ref{shrinking}). 

An {\em expanding soliton} lives on a time interval $(T, \infty)$.
For convenience, we take $T \: = \: 0$. Then the equation is
\begin{equation}
2\Ric \: + \: {\mathcal L}_V g \: + \: \frac{g}{t} \: = \: 0.
\end{equation}
The vector field $V = V(t)$ satisfies
$V(t) \: = \: \frac{1}{t} \: V(1)$.
The corresponding Ricci flow is given by
\begin{equation}
g(t) \: = \: t \: \phi^{-1}(t)^* g(1),
\end{equation}
where $\{ \phi(t) \}$ is the $1$-parameter group of diffeomorphisms generated 
by $- V$, normalized by $\phi(1) \: = \: \Id$.

A {\em gradient expanding soliton} satisfies the equations
\begin{align} \label{expanding}
\frac{\partial g_{ij}}{\partial t} \: & = \: - \: 2 R_{ij} \: = \:
2 \nabla_i \nabla_j f \: + \: \frac{g_{ij}}{t}, \\
\frac{\partial f}{\partial t} \: & = \: |\nabla f|^2. \notag
\end{align}
It follows from (\ref{expanding}) that
$V \: = \: \nabla f$ satisfies $V(t) \: = \: \frac{1}{t} \: V(1)$.
The solution to (\ref{expanding}) is
\begin{align} \label{expandingsoln}
g(t) \: & = \: t \: \phi^{-1}(t)^* g(1), \\
f(t) \: & = \: \phi^{-1}(t)^* f(1). \notag
\end{align}

Conversely, given a metric $\widehat{g}$ and a function $\widehat{f}$
satisfying 
\begin{equation}
\widehat{R}_{ij} \: + \: \widehat{\nabla}_i \widehat{\nabla}_j \widehat{f} 
\: + \: \frac{1}{2} \:
\widehat{g}
\: = \: 0,
\end{equation}
put $V(t) \: = \: \frac{1}{t} \: \widehat{\nabla} \widehat{f}$. 
If we define $g(t)$ and $f(t)$ by 
\begin{align}
g(t) \: & = \: t \: \phi^{-1}(t)^* \widehat{g}, \\
f(t) \: & = \: \phi^{-1}(t)^* \widehat{f} \notag
\end{align}
then they satisfy (\ref{expanding}). 

Obvious examples of solitons are given by Einstein metrics, with 
$V = 0$. Any steady or expanding soliton on a closed
manifold comes from an Einstein metric.  Other examples of solitons
(see \cite[Chapter 2]{Chow-Knopf})
are : \\ \\
1. (Gradient steady soliton) The cigar soliton on $\R^2$ and the Bryant soliton
on $\R^3$. \\
2. (Gradient shrinking soliton) Flat $\R^n$ with $f \: = \: - \:
\frac{|x|^2}{4t}$. \\
3. (Gradient shrinking soliton) The shrinking cylinder 
$\R \times S^{n-1}$ with 
$f \: = \: - \: \frac{x^2}{4t}$, where $x$ is the coordinate on $\R$. \\
4. (Gradient shrinking soliton) The Koiso soliton on $\C P^2 \#
\overline{\C P^2}$.

\section{Local derivative estimates} \label{applocalder}

\begin{theorem} \label{localder}
For any $\alpha, K, K^\prime, l \: \ge \: 0$ and $m, n \in \Z^+$,
there is some $C \: = \: C(\alpha, K, K^\prime, l,m, n)$ with the
following property. 
Given $r > 0$, 
suppose that $g(t)$ is a Ricci flow solution 
for $t \in [0, \overline{t}]$, where $0 \: < \: 
\overline{t} \: \le \:
\frac{\alpha r^2}{K}$, defined on an open neighborhood $U$ of a point
$p \in M^n$. Suppose that $\overline{B(p, r, 0)}$ is a compact
subset of $U$, that
\begin{equation}
|\Rm(x, t)| \: \le \: \frac{K}{r^{2}}
\end{equation}
for all $x \in U$ and $t \in [0, \overline{t}]$, and that
\begin{equation}
|\nabla^\beta \Rm(x, 0)| \: \le \: \frac{K^\prime}{r^{|\beta|+2}}
\end{equation}
for all $x \in U$ and $|\beta| \: \le \: l$.
Then
\begin{equation}
|\nabla^\beta \Rm(x, t)| \: \le \: 
\frac{C}{
r^{|\beta|+2} \: \left( \frac{t}{r^2} \right)^{\frac{\max(m-l,0)}{2}}
}
\end{equation}
for all $x \in B \left(p, \frac{r}{2}, 0 \right)$,
$t \in (0, \overline{t}]$ and $|\beta| \: \le \: m$. 

In particular, 
$|\nabla^\beta \Rm(x, t)| \: \le \: 
\frac{C}{r^{|\beta|+2}}$ whenever 
$|\beta| \: \le \: l$.
\end{theorem}

The main case $l = 0$ of Theorem \ref{localder} is due to Shi
\cite{Shi2 (1989)}. The extension to $l \: \ge \: 0$ appears in 
\cite[Appendix B]{Lu-Tian}.

\section{Convergent subsequences of Ricci flow solutions} 
\label{subsequence}

\begin{theorem} 
\label{accumulate}
Given $r_0 \in (0, \infty]$, 
let $\{g_i(t)\}_{i=1}^\infty$ be a sequence of Ricci flow solutions
on connected pointed manifolds $(M_i, m_i)$, 
defined for $t \in (A, B)$ with $-\infty \le A < 0 < B \le \infty$.
We assume that for all $i$, $M_i$ equals the time-zero ball $B_0(m_i, r_0)$
and for all $r \in (0, r_0)$, $\overline{B_0(m_i, r)}$ is compact. 
Suppose that the following
two conditions are satisfied : \\
1. For each $r \in (0, r_0)$ and each compact interval $I \subset (A,B)$, 
there is an $N_{r,I} < \infty$ so that for all $t \in I$ and all $i$,
$\sup_{B_0(m_i, r) \times I} |\Rm(g_i)| \le N_{r,I}$, and \\
2. The
time-$0$ injectivity radii 
$\{\inj(g_i(0))(m_i)\}_{i=1}^\infty$ 
are uniformly bounded below by a positive number.

Then after passing to a subsequence, the solutions converge smoothly
to a
Ricci flow solution $g_\infty(t)$ on a connected pointed manifold
$(M_\infty, m_\infty)$, defined for
$t \in (A, B)$, for which $M_\infty = B_0(m_\infty, r_0)$ and
$\overline{B_0(m_\infty, r)}$ is compact for all $r \in (0, r_0)$. 
That is, for any compact interval $I \subset (A, B)$
and any $r < r_0$,
there are pointed
time-independent diffeomorphisms $\phi_{r,i} \: : \: 
B_0(m_\infty, r) \rightarrow B_0(m_i,r)$ so that
$\{(\phi_{r,i} \times \Id)^* g_i\}_{i=1}^\infty$  converges smoothly
to $g_\infty$ on $B_0(m_\infty, r) \times I$. 
\end{theorem}

Given the sectional curvature bounds, the lower
bound on the injectivity radii is equivalent to a lower bound
on the volumes of balls around $m_i$ 
\cite[Theorem 4.7]{Cheeger-Gromov-Taylor}.
Theorem \ref{accumulate} is a slight generalization of
\cite[Main Theorem]{Hamiltonnn}, in which $r_0 = \infty$ and
$N_{r,I}$ is independent of $r$; see Corollary \ref{accumulate3}
below.

There are many variants of the theorem with alternative hypotheses.
One can replace the interval $(A,B)$ with an interval $(A, B]$,
$-\infty \le A < 0 \le B < \infty$.
One can also replace the interval $(A, B)$ with an interval
$[A, B)$, $-\infty < A< 0 < B
\le \infty$, if in addition one has uniform time-$A$ bounds 
$\sup_{B_A(m_i, r)} |\nabla^j \Rm(g_i(A))| \: \le \: C_{r, j}$.
Then using Appendix \ref{applocalder},
one gets smooth convergence to a limit solution $g_\infty$ 
on the time interval $[A, B)$. 
(Without the time-$A$ bounds one would only get
$C^0$-convergence on $[A, B)$ and $C^\infty$-convergence
on $(A, B)$.)
There is a similar statement if
one replaces the interval $(A, B)$ with an interval
$[A, B]$, $-\infty < A< 0 \le B
< \infty$. There is a version in which balls are replaced
by annuli. One can generalize the hypotheses to allow for
an $r$-dependent time interval.

In the setting of Theorem \ref{accumulate}, suppose that
$r_0 = \infty$. Then the time-zero slice $(M_\infty, m_\infty, g_\infty(0))$
is complete, but it does not immediately follow that the
other time slices are complete.  We now give a 
condition which will guarantee completeness, and which will be
sufficient for our purposes.

\begin{corollary}
\label{accumulate2}
Let $\{g_i(t)\}_{i=1}^\infty$ be a sequence of Ricci flow solutions
on connected pointed manifolds $(M_i, m_i)$, 
defined for $t \in (A, 0]$ with $-\infty \le A < 0$.
Suppose that each time-zero slice $(M_i, m_i, g_i(0))$ is complete.
Suppose that the following
three conditions are satisfied : \\
1. For each $r \in (0, \infty)$ and each compact interval $I \subset (A,0]$, 
there is an $N_{r,I} < \infty$ so that for all $i$,
$\sup_{B_0(m_i, r) \times I} |\Rm(g_i)| \le N_{r,I}$,  \\
2. For each compact interval $I \subset (A, 0]$, there  is some
$N^\prime_I \in (0, \infty)$ with the following property : for each 
$r \in (0, \infty)$, there is some $J_{r,I} \in \Z^+$ so that whenever
$i \ge J_{r,I}$, we have
$\Ric(g_i) \ge - N^\prime_I g_i$ on $B_0(m_i, r) \times I$, and \\
3. The
time-$0$ injectivity radii 
$\{\inj(g_i(0))(m_i)\}_{i=1}^\infty$ 
are uniformly bounded below by a positive number.

Then after passing to a subsequence, the solutions converge smoothly
to a
Ricci flow solution $g_\infty(t)$ on a connected pointed manifold
$(M_\infty, m_\infty)$, defined for
$t \in (A, 0]$, with complete time slices.
\end{corollary}
\begin{proof}
Let $(M_\infty, m_\infty, g_\infty(\cdot))$ be constructed as in Theorem
\ref{accumulate}. Then on each compact time interval $I \subset (A, 0]$, 
we have
$\Ric(g_\infty) \ge - N^\prime_I g_\infty$. By (\ref{similar}),
if $t \in I$ then for any $m_0, m_1 \in M_\infty$, we have 
\begin{equation} \label{similar2}
\dist_t(m_0, m_1) \ge e^{- N^\prime_I t} \dist_0(m_0, m_1).
\end{equation}

Suppose that $\{m_j\}_{j=1}^\infty$ is a Cauchy sequence
in $(M_\infty, g_\infty(t))$. From (\ref{similar2}), it is also
a Cauchy sequence in $(M_\infty, g_\infty(0))$, and so has a
subsequence, which we relabel as $\{m_j\}_{j=1}^\infty$, that
converges to some $m^\prime$ in $(M_\infty, g_\infty(0))$.
However, for small $\epsilon > 0$, the restriction of the
identity map $(M_\infty, g_\infty(0)) \rightarrow 
(M_\infty, g_\infty(t))$ to
$B_0(m^\prime, \epsilon)$ is biLipschitz.  It follows that
$\lim_{j \rightarrow \infty} m_j = m^\prime$ in
$(M_\infty, g_\infty(t))$.
\end{proof}

There is an obvious analog to Corollary \ref{accumulate2} in
which we assume that the Ricci flows are defined for $[0, B)$ and 
for any compact interval $I \subset [0, B)$, 
there is an upper bound $\Ric(g_i) \le N^\prime_I g_i$
on $B_0(m_i, r) \times I$ for large $i$.
If we assume double-sided curvature bounds 
then the statement is as follows.

\begin{corollary}
\label{accumulate3}
Let $\{g_i(t)\}_{i=1}^\infty$ be a sequence of Ricci flow solutions
on connected pointed manifolds $(M_i, m_i)$, 
defined for $t \in (A, B)$ with $-\infty \le A < 0 < B \le \infty$.
We assume that for all $i$, the time slice
$(M_i, m_i, g_i(0))$ is complete.
Suppose that the following
two conditions are satisfied : \\
1. For each compact interval $I \subset (A,B)$, 
there is an $N_{I} < \infty$ with the following property :
for each $r \in (0, \infty)$, there is some $J_{r,I} \in \Z^+$
so that whenever $i \ge J_{r,I}$, we have
$\sup_{B_0(m_i, r) \times I} |\Rm(g_i)| \le N_{I}$, and \\
2. The
time-$0$ injectivity radii 
$\{\inj(g_i(0))(m_i)\}_{i=1}^\infty$ 
are uniformly bounded below by a positive number.

Then after passing to a subsequence, the solutions converge smoothly
to a
Ricci flow solution $g_\infty(t)$ on a connected pointed manifold
$(M_\infty, m_\infty)$, defined for
$t \in (A, B)$, with complete time slices.
\end{corollary}

In the case when $J_{r,I} = 1$ for all $r$ and $I$,
i.e. $\sup_{M_i \times I} |\Rm(g_i)| \le N_{I}$,
Corollary \ref{accumulate3} is the same as
\cite[Main Theorem]{Hamiltonnn}.

\section{Harnack inequalities for Ricci flow} \label{appharnack}

We first recall the statement of the matrix Harnack inequality.
Put 
\begin{align}
P_{abc} \: & = \: \nabla_a R_{bc} \: - \: \nabla_b R_{ac}, \\
M_{ab} \: & = \: \triangle R_{ab} \: - \: \frac12 \: \nabla_a \nabla_b
R \: + \: 2 \: R_{acbd} R_{cd} \: - \: R_{ac} R_{bc} \: + \:
\frac{R_{ab}}{2t}. \notag
\end{align}
Given a $2$-form $U$ and a $1$-form $W$, put
\begin{equation}
Z(U, W) \: = \: M_{ab} W_a W_b \: + \: 2 \: P_{abc} U_{ab} W_c
\: + \: R_{abcd} U_{ab} U_{cd}.
\end{equation}
Suppose that we have a Ricci flow for $t > 0$ 
on a complete manifold with bounded
curvature on each 
compact time interval
and nonnegative curvature operator.
Hamilton's matrix Harnack inequality says that for all $t > 0$ and all
$U$ and $W$, $Z(U, W) \: \ge \: 0$ \cite[Theorem 14.1]{Hamilton}.

Taking
$W_a \: = \: Y_a$ and $U_{ab} \: = \: (X_a Y_b \: - \:
Y_a X_b)/2$ and 
using the fact that
\begin{equation} \label{Rict}
\Ric_t(Y, Y) \: = \: (\triangle R_{ab}) Y^a Y^b \: + \:
2 \: R_{acbd} R_{cd} Y^a Y^b \: - \: 2 \: R_{ac} R_{bc} Y^a Y^b, 
\end{equation}
we can write $2Z(U, W) \: = \: H(X, Y)$, where
\begin{align} \label{Hexpression}
H(X, Y) \: = \: & - \: \Hess_R (Y,Y) \: - \: 2 \langle R(Y, X)Y, X \rangle
\: +\: 4 \left( \nabla_X \Ric(Y, Y) \: - \: \nabla_Y \Ric(Y, X) \right) \\
&  + \: 2 \Ric_t(Y, Y) \: + \: 2 \big| \Ric(Y, \cdot) \big|^2
\: + \: \frac{1}{t} \: \Ric(Y, Y). \notag
\end{align}
Substituting the elements of an orthonormal basis $\{e_i\}_{i=1}^n$ for
$Y$ and summing over $i$ gives
\begin{align}
\sum_i H(X, e_i) \: = & \: - \: \triangle R \: + \: 2 \Ric(X, X) \: +  \\
& 4 (\langle \nabla R, X \rangle
 \: - \: \sum_i \nabla_{e_i} \Ric(e_i, X)) \: + \:
2 \sum_i \Ric_t(e_i, e_i) \: + \: 2 |\Ric|^2 \: + \: \frac{1}{t} \: R.
\notag
\end{align}
Tracing the second Bianchi identity gives
\begin{equation}
\sum_i \nabla_{e_i} \Ric(e_i, X) \: = \: 
\frac{1}{2} \: \langle \nabla R, X \rangle.
\end{equation}
From (\ref{Rict}),
\begin{equation}
\sum_i \Ric_t(e_i, e_i) \: = \: \triangle R.
\end{equation}
Putting this together gives
\begin{equation} \label{together}
\sum_i H(X, e_i) \: = \: H(X),
\end{equation}
where 
\begin{equation} \label{Heqn}
H(X) \: = \: R_t \: + \: \frac{1}{t} \: R \: + \:
2 \langle \nabla R, X \rangle \: 
+ \: 2 \Ric(X, X).
\end{equation}
We obtain Hamilton's trace Harnack 
inequality, saying that $H(X) \: \ge \: 0$ for all $X$. 

In the rest of this section we assume that
the solution is defined for all $t \in (- \infty, 0)$. 
Changing the origin point of time, we have
\begin{equation}
R_t \: + \: \frac{1}{t- t_0} \: R \: + \:
2 \langle \nabla R, X \rangle \: 
+ \: 2 \Ric(X, X) \: \ge \: 0
\end{equation} 
whenever $t_0 \: \le \: t$. Taking $t_0 \rightarrow - \infty$ gives
\begin{equation}
R_t \: + \:
2 \langle \nabla R, X \rangle \: 
+ \: 2 \Ric(X, X) \: \ge \: 0
\end{equation} 
In particular, taking $X = 0$ shows that the scalar curvature is
nondecreasing in $t$ for any ancient solution with nonnegative
curvature operator, assuming again that the metric is complete
on each time slice with bounded curvature
on each compact time interval.
More generally,
\begin{equation}
0 \: \le \: R_t \: + \: 
2 \langle \nabla R, X \rangle \: 
+ \: 2 \Ric(X, X) \: \le \: R_t \: + \: 
2 \langle \nabla R, X \rangle \: 
+ \: 2 R \langle X, X \rangle.
\end{equation} 
If $\gamma \: : \: [t_1, t_2] 
\rightarrow M$ is a curve parametrized by $s$ then taking
$X \: = \: \frac12 \: \frac{d\gamma}{ds}$ gives
\begin{equation}
\frac{d R(\gamma(s), s)}{ds}  \: = \:
R_t(\gamma(s), s) \: + \: \left\langle \frac{d\gamma}{ds}, \nabla R 
\right\rangle
\: \ge \: - \: \frac12 \: R \: \left\langle \frac{d\gamma}{ds}, \frac{d\gamma}{ds}
\right\rangle.
\end{equation}
Integrating $\frac{d \ln R(\gamma(s),s)}{ds}$ with respect to $s$ and using the fact that
$g(t)$ is nonincreasing in $t$ gives
\begin{equation}
\label{harnackinequality}
R(x_2,t_2)\geq 
\exp\left(-\frac{d_{t_1}^2(x_1,x_2)}{2(t_2-t_1)}\right)R(x_1,t_1).
\end{equation}
whenever $t_1 \: < \: t_2$ and $x_1, x_2 \in M$.
(If $n = 2$ then one can replace $\frac{d_{t_1}^2(x_1,x_2)}{2(t_2-t_1)}$ by
$\frac{d_{t_1}^2(x_1,x_2)}{4(t_2-t_1)}$.) In particular, if
$R(x_2, t_2) \: = \: 0$ for some $(x_2, t_2)$ then $g(t)$ must be flat
for all $t$.

\section{Alexandrov spaces} \label{Alexandrov}

We recall some facts about Alexandrov spaces
(see \cite[Chapter 10]{Burago-Burago-Ivanov},
\cite{Burago-Gromov-Perelman}).
Given 
points $p,x,y$ in a nonnegatively curved Alexandrov space,
we let $\cangle_p(x,y)$ denote the comparison angle at $p$,
i.e. the angle of the Euclidean
comparison triangle at the vertex corresponding to $p$. 

The Toponogov splitting theorem says that if $X$ is a
proper nonnegatively curved Alexandrov space which contains a line, then $X$
splits isometrically as a product $X=\R\times Y$, where
$Y$ is a proper, nonnegatively curved Alexandrov space
\cite[Theorem 10.5.1]{Burago-Burago-Ivanov}.

Let $M$ be an $n$-dimensional Alexandrov space with nonnegative curvature,
$p\in M$,
and $\la_k\ra 0$. Then the sequence $(\la_kM,p)$ of pointed spaces
Gromov-Hausdorff converges to the Tits cone $\ctits M$
(the Euclidean cone over the Tits boundary $\tits M$)
which is a nonnegatively curved Alexandrov
space of dimension $\leq n$
\cite[p. 58-59]{Ballmann-Gromov-Schroeder}.  
If the Tits cone splits isometrically
as a product $\R^k\times Y$, then $M$ itself splits off
a factor of $\R^k$; using triangle comparison, one finds $k$
orthogonal lines passing through a basepoint, and applies the
Toponogov splitting theorem.  

Now suppose that $x_k\in M$ is a sequence with $d(x_k,p)\ra\infty$
and $r_k \in \R^+$  is a sequence with  
$\frac{r_k}{d(x_k,p)}\ra 0$.  Then the sequence 
$(\frac{1}{r_k} M,x_k)$ 
subconverges to a pointed Alexandrov
space $(N_\infty,x_\infty)$ which splits off a line.  To see this,
observe that since 
$(\frac{1}{d(x_k,p)}M,p)$ 
converges to a cone,
we can find a 
sequence $y_k\in M$ such that $\frac{d(y_k,x_k)}{d(x_k,p)}\ra 1$,
and $\cangle_{x_k}(p,y_k)\ra\pi$.   
Observe that for any $\rho<\infty$,
we can find sequences $p_k\in \ol{px_k}$, $z_k\in\ol{x_ky_k}$
such that $\frac{d(x_k,p_k)}{r_k}\ra \rho$, 
$\frac{d(x_k,z_k)}{r_k}\ra \rho$, and by monotonicity of comparison
angles 
\cite[Chapter 4.3]{Burago-Burago-Ivanov}
we will have $\cangle_{x_k}(p_k,z_k)\ra \pi$.  Passing
to the Gromov-Hausdorff limit, we find $p_\infty, z_\infty\in N_\infty$
such that $d(p_\infty,x_\infty)=d(z_\infty,x_\infty)=\rho$ and 
$\cangle_{x_\infty}(p_\infty,z_\infty)=\pi$.  Since this
construction applies for all $\rho$, it follows that $N_\infty$
contains a line passing through $x_\infty$.
Hence, by the Toponogov splitting theorem, it is isometric
to a metric product $\R\times N'$ for some Alexandrov space
$N'$.

If $M$ is a complete Riemannian manifold of nonnegative sectional 
curvature and $C\subset M$ is a 
compact connected domain with weakly convex boundary then 
the subsets $C_t= \{x\in C\mid d(x,\partial C)\geq t\}$ 
are convex in $C$ 
\cite[Chapter 8]{Cheeger-Ebin}. 
If the second fundamental form of 
$\partial C$ is $\geq \frac{1}{r}$ at each point
of $\partial C$, then for all $x\in C$ we have $d(x,\partial C)\leq r$,
since the first focal point of $\partial C$ along any inward pointing
normal geodesic occurs at distance $\leq r$. 

Finally, we recall the statement of the Bishop-Gromov volume
comparison inequality.  Suppose that $M$ is an
$n$-dimensional
Riemannian manifold.
Given $p \in M$ and $r_2 \ge r_1 > 0$,
suppose that $B(p, r_2)$ has compact closure in $M$ and that the
sectional curvatures of $B(p, r_2)$ are bounded below by
$K \in \R$. Then
\begin{equation} \label{BG}
\frac{\vol(B(p, r_2))}{\vol(B(p, r_1))} \: \le \:
\begin{cases}
\frac{\int_0^{r_2} \sin^{n-1}(kr) \: dr}{\int_0^{r_1} \sin^{n-1}(kr) \: dr} 
& \text{ if } K = k^2, \\ \\
\frac{r_2^n}{r_1^n} 
& \text{ if } K = 0, \\ \\
\frac{\int_0^{r_2} \sinh^{n-1}(kr) \: dr}{\int_0^{r_1} \sinh^{n-1}(kr) \: dr} 
& \text{ if } K =  -  k^2.
\end{cases}
\end{equation}
(If $K = k^2$ with $k>0$ then we restrict to $r_2 \le \frac{\pi}{k}$.)
The same inequality holds if we just assume that
$B(p, r_2)$ has Ricci curvature bounded below by
$(n-1) \: K$. Equation (\ref{BG}) also holds if
$M$ is an Alexandrov space and $B(p, r_2)$ has Alexandrov curvature
bounded below by $K$.

\section{Finding balls with controlled curvature} \label{pointpicking}

\begin{lemma}
\label{selection}
Let $X$ be a Riemannian manifold
with $R \ge 0$
and suppose $\ol{B(x,5r)}$
is a compact subset of $X$.  Then there is a ball 
$B(y,\bar r)\subset B(x,5r)$, $\overline{r} \le r$,
such that $R(z)\leq 2R(y)$
for all $z\in B(y,\bar r)$ and $R(y)\bar r^2\geq R(x)r^2$.
\end{lemma}
\begin{proof}
Define sequences $x_i\in B(x,5r)$, $r_i>0$
inductively as follows.  Let $x_1=x$, $r_1=r$.
For $i>1$, let $x_{i+1}=x_i$, $r_{i+1}=r_i$ if
$R(z)\leq 2R(x_i)$ for all $z\in B(x_i,r_i)$;
otherwise let $r_{i+1}=\frac{r_i}{\sqrt{2}}$, and let 
$x_{i+1}\in B(x_i,r_i)$ be a point such that
$R(x_{i+1})>2R(x_i)$.  The sequence of balls
$B(x_i,r_i)$ is contained in $B(x,5r)$, so the
sequences $x_i,\,r_i$ are eventually constant,
and we can take $y=x_i$, $\bar r=r_i$ for large
$i$.
\end{proof}

There is an evident spacetime version of the lemma.

\section{Statement of the geometrization conjecture} \label{geom}

Let $M$ be a connected orientable closed (= compact boundaryless) $3$-manifold.
One formulation of the geometrization conjecture says that
$M$ is the connected sum of closed $3$-manifolds $\{M_i\}_{i=1}^n$, each of
which admits a codimension-$0$ compact submanifold-with-boundary $G_i$
so that 
\begin{itemize}
\item $G_i$ is a graph manifold
\item $M_i - G_i$ is hyperbolic, i.e. admits a
complete Riemannian metric with constant negative sectional curvature
and finite volume
\item Each component $T$ of $\partial G_i$ is an incompressible torus in $M_i$,
i.e. with respect to a basepoint $t \in T$, the induced map
$\pi_1(T, t) \rightarrow \pi_1(M_i, t)$ is injective.
\end{itemize}
We allow $G_i = \emptyset$ or $G_i = M_i$. 

A reference for graph manifolds is
\cite[Chapter 2.4]{Matveev}. The definition is as follows.
One takes a collection $\{P_i\}_{i=1}^N$ of pairs of pants
(i.e. closed $2$-disks with two balls removed)
and a collection of closed $2$-disks $\{D^2_j\}_{j=1}^{N^\prime}$. The
$3$-manifolds $\{S^1 \times P_i\}_{i=1}^N
\cup \{S^1 \times D^2_j\}_{j=1}^{N^\prime}$ have toral boundary components.
One takes an even number of these tori, matches them in pairs by
homeomorphisms, and glues $\{S^1 \times P_i\}_{i=1}^N
\cup \{S^1 \times D^2_j\}_{j=1}^{N^\prime}$ by these homeomorphisms.
The
resulting $3$-manifold $G$ is a graph manifold, and all 
graph manifolds
arise in this way.  
We will assume that the gluing homeomorphisms are such that
$G$ is orientable.  
Clearly the boundary of $G$, if nonempty, is a
disjoint union of tori. It is also clear that the result of
gluing two graph manifolds along some collection of boundary tori is
a graph manifold. The connected sum of two $3$-manifolds is a 
graph manifold if and only if each factor is a graph manifold 
\cite[Proposition 2.4.3]{Matveev}.

The reason to require incompressibility of the tori $T$ in the
statement of the geometrization conjecture is to
exclude phony decompositions, such as writing $S^3$ as the
union of a solid torus and a hyperbolic knot complement.

A more standard version of the geometrization conjecture uses
some facts from $3$-manifold theory
\cite{Scott}. First $M$ has a Kneser-Milnor decomposition as
a connected sum of uniquely defined prime factors. Each prime factor is
$S^1 \times S^2$ or is irreducible, i.e. any embedded $S^2$ bounds
a $3$-ball. If $M$ is irreducible then it has a JSJ decomposition, i.e.
there is a minimal collection
of disjoint incompressible 
embedded tori $\{T_k\}_{k=1}^K$ in $M$, unique up
to isotopy, with the property that
if $M^\prime$ is the metric completion of a component of 
$M - \bigcup_{k=1}^K T_k$ (with respect to an induced 
Riemannian metric from $M$) then
\begin{itemize}
\item $M^\prime$ is a Seifert $3$-manifold or
\item $M^\prime$ is non-Seifert and any embedded incompressible torus in
$M^\prime$ can be isotoped into $\partial M^\prime$.
\end{itemize}
The second version of the geometrization conjecture reduces to
the conjecture that in the latter case, the interior
of $M^\prime$ is hyperbolic. Thurston proved that this is true when
$\partial M^\prime \neq \emptyset$.
The reason for the word ``geometrization''
is explained in \cite{Scott,Thurston}.

An orientable Seifert $3$-manifold is a graph manifold
\cite[Proposition 2.4.2]{Matveev}. It follows that
the second version of the geometrization conjecture implies the
first version. One can show directly that any graph manifold
is a connected sum of prime graph manifolds, each of which can be split
along incompressible tori to
obtain a union of Seifert manifolds \cite[Proposition 2.4.7]{Matveev},
thereby showing the equivalence of the two versions.


\begin{thebibliography}{20}
\bibitem{Akutagawa-Ishida-Lebrun} K. Akutagawa, M. Ishida and C. Lebrun,
``Perelman's Invariant, Ricci Flow, and the Yamabe Invariants of Smooth Manifolds'',
Arch. Math. 88, p. 71-76 (2007)

\bibitem{Anderson3} M. Anderson,
``Scalar Curvature and Geometrization Conjectures for $3$-Manifolds'',
in \underline{Comparison Geometry}, MSRI Publ. 30,
Cambridge University Press, Cambridge, p. 49-82 (1997)

\bibitem{Anderson} M. Anderson, 
``Scalar Curvature and the Existence of Geometric Structures on 
$3$-Manifolds. I'', J. Reine Angew. Math. 553, p. 125-182 (2002)

\bibitem{Anderson2} M. Anderson,
``Canonical Metrics on 3-Manifolds and 4-Manifolds'',
Asian Journal of Math. 10, p. 127-164 (2006)

\bibitem{Anderson-Chow} G. Anderson and B. Chow,
``A Pointwise Bound for Solutions of the Linearized Ricci Flow
Equation on 3-Manifolds'', Calculus of Variations 23, p. 1-12 (2005)

\bibitem{Angenent-Knopf}  S. Angenent and D. Knopf,
``Precise Asymptotics of the Ricci Flow
Neckpinch'', Comm. Anal. Geom. 15, p. 773-844 (2007)

\bibitem{Ballmann-Gromov-Schroeder} 
W. Ballmann, M. Gromov and V. Schroeder,
\underline{Manifolds of Nonpositive Curvature},
Progress in Mathematics 61, 
Birkh\"auser, Boston (1985)

\bibitem{Beckner} 
W. Beckner,
``Geometric Asymptotics and the Logarithmic Sobolev Inequality'',
Forum Math. 11, p. 105-137 (1999)

\bibitem{Besson} L. Bessi\`eres, G. Besson, M. Boileau, S. Maillot and
J. Porti, ``Weak Collapsing and Geometrisation of Aspherical 3-Manifolds'',
preprint, http://arxiv.org/abs/0706.2065 (2007)

\bibitem{Bourguignon} J.-P. Bourguignon, ``Une Stratification de l'Espace 
des Structures Riemanniennes'',
Compositio Math. 30, p. 1-41 (1975)

\bibitem{Burago-Burago-Ivanov}
D. Burago, Y. Burago and S. Ivanov, 
\underline{A Course in Metric Geometry},
Graduate Studies in Mathematics 33,  
American Mathematical Society, Providence (2001)

\bibitem{Burago-Gromov-Perelman}
Y. Burago, M. Gromov and G. Perelman,
``A. D. Aleksandrov Spaces with Curvatures Bounded Below'',
Russian Math. Surveys 47, p. 1-58 (1992)

\bibitem{Calabi} E. Calabi,
``An Extension of E. Hopf's Maximum Principle with an Application to 
Riemannian Geometry'',
Duke Math. J. 25, p. 45-56 (1957)

\bibitem{Cao-Chow} H-D. Cao and B. Chow,
``Recent Developments on the Ricci Flow'', Bull. of the AMS 36, p. 59-74 (1999)

\bibitem{Cao-Zhu} H.-D. Cao and X.-P. Zhu,
``A Complete Proof of the Poincar\'e and Geometrization Conjectures - Application of the Hamilton-Perelman theory of the Ricci Flow'',
Asian Journal of Math. 10, p. 165-492 (2006) \\
Erratum to
``A Complete Proof of the Poincar\'e and Geometrization Conjectures - Application of the Hamilton-Perelman theory of the Ricci Flow'',  Asian Journal of Math. 10, p. 663-664 (2006)            

\bibitem{Chau-Tam-Yu (2007)} A. Chau, L.-F. Tam and C. Yu,
``Pseudolocality for the Ricci Flow and Applications'',
preprint, http://arxiv.org/abs/math/0701153 (2007)

\bibitem{Chavel} I. Chavel,
\underline{Riemannian Geometry - a Modern Introduction}, 
Cambridge Tracts in Mathematics 108, 
Cambridge University Press, Cambridge (1993)

\bibitem{Cheeger-Ebin} J. Cheeger and D. Ebin,
\underline{Comparison Theorems in Riemannian Geometry}. 
North-Holland Publishing Co., Amsterdam-Oxford (1975)

\bibitem{Cheeger-Gromoll} J. Cheeger and D. Gromoll,
``On the Structure of Complete Manifolds of Nonnegative Curvature'', 
Ann. of Math. 96, p. 413-443 (1972)

\bibitem{Cheeger-Gromov-Taylor} J. Cheeger, M. Gromov and M. Taylor,
``Finite Propagation Speed, Kernel Estimates for Functions of the
Laplace Operator, and the Geometry of Complete Riemannian Manifolds'',
J. Diff. Geom. 17, p. 15-53 (1982)

\bibitem{Chen-Zhu} B. Chen and X. Zhu,
``Uniqueness of the Ricci Flow on Complete
Noncompact Manifolds'', J. of Diff. Geom. 74, p. 119-154
(2006)

\bibitem{Chow-Knopf} B. Chow and D. Knopf,
\underline{The Ricci Flow : An Introduction}, AMS Mathematical Surveys and
Monographs, Amer. Math. Soc., Providence (2004)

\bibitem{Chow-Lu-Ni} B. Chow, P. Lu and L. Ni, \underline{Hamilton's Ricci Flow},
Graduate Studies in Mathematics 77, AMS, Providence (2006)

\bibitem{Colding-Minicozzi} T. Colding and W. Minicozzi,
``Estimates for the Extinction Time for the Ricci Flow on Certain 3-Manifolds and a 
Question of Perelman'',
J. Amer. Math. Soc. 18, p.  561-569 (2005)

\bibitem{Colding-Minicozzi2} T. Colding and W. Minicozzi,
``Width and Finite Extinction Time of Ricci Flow'',
preprint, http://arxiv.org/abs/0707.0108 (2007)

\bibitem{Dodziuk} J. Dodziuk,  ``Maximum Principle for Parabolic 
Inequalities and the Heat Flow on Open Manifolds'',
Indiana Univ. Math. J. 32, p. 703-716 (1983)

\bibitem{Eells-Sampson} J. Eells and J. Sampson,
``Harmonic Mappings of Riemannian Manifolds'',
Amer. J. Math. 86, p. 109-160 (1964)

\bibitem{Eschenburg} J.-H. Eschenburg, ``Local Convexity and Nonnegative 
Curvature---Gromov's Proof of the Sphere Theorem'', 
Inv. Math. 84, p. 507-522 (1986)

\bibitem{Hamiltonnnn} R. Hamilton,
``Four-Manifolds with Positive Curvature Operator'', 
J. Diff. Geom. 24, p. 153-179 (1986)

\bibitem{Hamsurface} R. Hamilton,
``The Ricci Flow on Surfaces'', in
\underline{Mathematics and General Relativity},
Contemp. Math. 71, 
AMS, Providence, RI, p. 237-262 (1988).

\bibitem{harnack} R. Hamilton,
``The Harnack Estimate for the Ricci Flow'', 
J. Diff. Geom. 37, p. 225-243 (1993)

\bibitem{Hamiltonnn} R. Hamilton,
``A Compactness Property for Solutions of the Ricci Flow''
Amer. J. Math. 117, p. 545-572 (1995)

\bibitem{Hamilton} R. Hamilton,
``The Formation of Singularities in the Ricci Flow'', in
\underline{Surveys in Differential Geometry}, Vol. II,
Internat. Press, Cambridge, p. 7-136 (1995) 

\bibitem{Hamiltonn} R. Hamilton,
``Nonsingular Solutions of the Ricci Flow on Three-Manifolds'',
Comm. Anal. Geom. 7, p. 695-729 (1999)

\bibitem{Hamiltonnnnn} R. Hamilton,
``Three-Manifolds with Positive Ricci Curvature'', 
J. Diff. Geom. 17, p. 255-306 (1982)

\bibitem{Hamiltonnnnnn} R. Hamilton, 
``Four-Manifolds with Positive Isotropic Curvature'', 
Comm. Anal. Geom. 5, p. 1-92 (1997)

\bibitem{Ishii} H. Ishii,
``On the Equivalence of Two Notions of Weak Solutions, Viscosity Solutions 
and Distribution Solutions'',
Funkcial. Ekvac. 38, p. 101-120 (1995)

\bibitem{Kapovitch} V. Kapovitch, ``Perelman's Stability
Theorem,'' in \underline{Surveys in Differential Geometry},
vol. XI, Metric and Comparison Geometry, eds. J. Cheeger and
K. Grove, International Press, Somerville, MA, p. 103-136

\bibitem{KleinerLottwebsite} 
B. Kleiner and J. Lott, webpage at
http://www.math.lsa.umich.edu/\~{}lott/ricciflow/perelman.html

\bibitem{Kleiner-Lott} B. Kleiner and J. Lott, ``Locally Collapsed
$3$-Manifolds'', to appear

\bibitem{Lu-Tian} P. Lu and G. Tian,
``Uniqueness of Standard Solutions in the Work of Perelman'',
preprint, http://www.math.lsa.umich.edu/\~{}lott/ricciflow/perelman.html  (2005)

\bibitem{Matveev} 
S. Matveev,
\underline{Algorithmic Topology and Classification of 3-Manifolds},
Springer-Verlag, Berlin (2003)

\bibitem{meeksyau}
W. Meeks III and S.-T. Yau,
``Topology of Three-Dimensional Manifolds and the Embedding Problems in
Minimal Surface Theory'',
Ann. of Math. 112, p. 441-484 (1980)

\bibitem{Morgan} J. Morgan, 
``Recent Progress on the Poincar\'e Conjecture and the Classification of $3$-Manifolds'',
Bull. Amer. Math. Soc. 42, p. 57-78 (2005)

\bibitem{Morgan-Tian} J. Morgan and G. Tian,
\underline{Ricci Flow and the Poincar\'e Conjecture},
Clay Mathematics Monographs 3, Amer. Math. Soc., Providence (2007)

\bibitem{Morgan-Tian2} J. Morgan and G. Tian,
``Completion of Perelman's Proof of the Geometrization Conjecture'',
to appear

\bibitem{Mostow} G. Mostow, 
\underline{Strong Rigidity of Locally Symmetric Spaces},
Annals of Math. Studies No. 78,
Princeton University Press, Princeton (1973)

\bibitem{Ni3} L. Ni, ``The Entropy Formula for Linear Heat Equation'',
J. Geom. Anal. 14, p. 87-100 (2004)

\bibitem{Ni2} L. Ni, ``Ricci Flow and Nonnegativity of Sectional
Curvature'', Math. Res. Let. 11, p. 883-904 (2004)

\bibitem{Ni} L. Ni, ``A Note on Perelman's LYH Inequality'', 
Comm. Anal. Geom. 14, p. 883-905 (2006) 

\bibitem{Perelman}
G. Perelman, 
``The Entropy Formula for the Ricci Flow and its Geometric
Applications'', http://arXiv.org/abs/math.DG/0211159

\bibitem{Perelman2}
G. Perelman, 
``Ricci Flow with Surgery on Three-Manifolds'',
http://arxiv.org/abs/math.DG/0303109

\bibitem{Perelman3}
G. Perelman, 
``Finite Extinction Time for the Solutions to the Ricci Flow
on Certain Three-Manifolds'',
http://arxiv.org/abs/math.DG/0307245

  
\bibitem{Prasad} G. Prasad, ``Strong Rigidity of $\Q$-Rank $1$ Lattices'',
Invent. Math. 21, p. 255-286 (1973)

\bibitem{Reed-Simon} M. Reed and B. Simon, 
\underline{Methods of Modern Mathematical Physics IV, Analysis of 
Operators}, Academic Press, New York (1978)

\bibitem{Rothaus (1981)} O. Rothaus,
``Logarithmic Sobolev Inequalities and the Spectrum of Schr\"odinger 
Operators'', J. Funct. Anal. 42, p. 110-120 (1981)

\bibitem{schoen} R. Schoen,
``Estimates for Stable Minimal Surfaces in Three-Dimensional Manifolds'',
in \underline{Seminar on Minimal Submanifolds}, 
Ann. of Math. Stud. 103, Princeton Univ. Press, Princeton, N.J., p. 111-126 (1983).

\bibitem{schoenyau}
R. Schoen and S.-T. Yau,
``Complete Three-Dimensional Manifolds with Positive Ricci Curvature
  and Scalar Curvature'', in
\underline{Seminar on Differential Geometry}, Ann. of
  Math. Stud. 102,  Princeton Univ. Press, Princeton, N.J., p.  209-228 (1982)

\bibitem{Scott} P. Scott,
``The Geometries of $3$-Manifolds'',
Bull. London Math. Soc. 15, p. 401-487 (1983) 

\bibitem{Sharafutdinov} V. Sharafutdinov,
``The Pogorelov-Klingenberg Theorem for Manifolds that are Homeomorphic to 
$\R^n$'',
Siberian Math. J. 18, p. 649-657 (1977)

\bibitem{Shi (1989)} W.-X. Shi, ``Complete Noncompact Three-Manifolds
with Nonnegative Ricci Curvature'',
J. Diff. Geom. 29, p. 353-360 (1989)

\bibitem{Shi2 (1989)} W.-X. Shi, ``Deforming the Metric on Complete 
Riemannian Manifolds'',
J. Diff. Geom. 30, p. 223-301 (1989)

\bibitem{Shi3 (1989)} W.-X. Shi, `` Ricci Deformation of the Metric on 
Complete Noncompact Riemannian Manifolds'',
J. Diff Geom. 30, p. 303-394 (1989)

\bibitem{Shioya-Yamaguchi} T. Shioya and T. Yamaguchi,
``Volume-Collapsed Three-Manifolds with a Lower Curvature Bound'',
Math. Ann. 333, p. 131-155 (2005)

\bibitem{Thurston} W. Thurston,
``Three-Dimensional Manifolds, Kleinian Groups and Hyperbolic Geometry'',
Bull. Amer. Math. Soc. 6, p. 357-381 (1982)

\bibitem{Topping} P. Topping, \underline{Lectures on the Ricci Flow},
LMS Lecture Notes 325, London Mathematical Society and Cambridge University
Press, 
http://www.maths.warwick.ac.uk/~topping/RFnotes.html
 
\bibitem{Ye1} R. Ye, ``On the $l$-Function and the Reduced Volume of
Perelman I,II'', Trans. Amer. Math. Soc. 360, p. 507-531, 533-544 (2008) 

\bibitem{Ye2} R. Ye, ``On the Uniqueness of 2-Dimensional
$\kappa$-Solutions'',
http://www.math.ucsb.edu/\~{}yer/ricciflow.html
\end{thebibliography}
\end{document}